\documentclass[11pt]{amsart}
\usepackage[square]{natbib}
\newif\ifamsartclass
\makeatletter\@ifclassloaded{amsart}{\amsartclasstrue}{}\makeatother

\usepackage[utf8]{inputenc}
\usepackage[T1]{fontenc}

\DeclareFontShape{T1}{cmr}{m}{scit}{<->ssub*cmr/m/sc}{}

\usepackage[margin=1in]{geometry}
\usepackage{amsmath,amsthm,amssymb,amsfonts}
\usepackage{mathrsfs}
\usepackage{mathtools}
\usepackage{esint}
\usepackage{xfrac}
\usepackage{nicefrac}
\usepackage{dsfont}
\usepackage{upgreek}
\usepackage{accents}

\let\olddiamond\diamond
\let\oldsquare\square
\usepackage{mathabx}
\renewcommand{\square}{\oldsquare}
\renewcommand{\diamond}{\olddiamond}

\usepackage{booktabs}
\usepackage{array}
\usepackage{paralist}
\usepackage{enumitem}
\usepackage{verbatim}
\usepackage[normalem]{ulem}
\usepackage{cancel}

\usepackage[dvipsnames]{xcolor}
\usepackage{graphicx}
\usepackage{subcaption}

\usepackage{tikz}
\usepackage{pgf}
\usetikzlibrary{calc}
\usetikzlibrary{external}
\usepackage[
colorlinks=true,
pdfstartview=FitV,
linkcolor=blue,
citecolor=blue,
urlcolor=blue
]{hyperref}
\ifamsartclass
\hypersetup{citecolor=PineGreen,linkcolor=RoyalBlue,urlcolor=BrickRed}
\fi

\providecommand*\showkeyslabelformat[1]{{\normalfont\tiny#1}}

\makeatletter
\newcommand{\runinsubsubsections}{%
	\renewcommand{\subsubsection}{\@startsection{subsubsection}{3}%
		\z@{.5\linespacing\@plus.7\linespacing}{-.5em}%
		{\normalfont\bfseries}}}
\newcommand{\runinparagraphs}[1][3.25ex\@plus1ex\@minus.2ex]{%
	\renewcommand{\paragraph}{\@startsection{paragraph}{4}%
		\z@{#1}{-1em}%
		{\normalfont\bfseries}}}
\makeatother

\ifamsartclass\else
\usepackage{titlesec}

\newcommand{\addperiod}[1]{#1.}
\titleformat{\section}{\centering\normalfont\Large}{\thesection.}{0.5em}{}
\titleformat*{\subsection}{\bfseries}
\titleformat{\subsubsection}[runin]{\normalfont\bfseries}{\thesubsubsection.}{0.5em}{\addperiod}
\titleformat*{\paragraph}{\bfseries}
\titleformat*{\subparagraph}{\large\bfseries}

\usepackage[nottoc,notlot,notlof]{tocbibind}
\usepackage{fancyhdr}
\usepackage{extramarks}

\fi

\numberwithin{equation}{section}
\numberwithin{figure}{section}

\newtheorem{theorem}{Theorem}[section]

\newtheorem{corollary}[theorem]{Corollary}
\newtheorem{proposition}[theorem]{Proposition}
\newtheorem{lemma}[theorem]{Lemma}

\newtheorem{theoremA}{Theorem}

\theoremstyle{definition}
\newtheorem{definition}[theorem]{Definition}

\newtheorem{remark}[theorem]{Remark}

\DeclareMathOperator{\dist}{dist}

\DeclareMathOperator*{\osc}{osc}

\DeclareMathOperator*{\argmin}{argmin}
\DeclareMathOperator{\var}{var}

\DeclareMathSymbol{\shortminus}{\mathbin}{AMSa}{"39}

\let\originalleft\left
\let\originalright\right
\renewcommand{\left}{\mathopen{}\mathclose\bgroup\originalleft}
\renewcommand{\right}{\aftergroup\egroup\originalright}

\newcommand\thickbar[1]{\accentset{\rule{.45em}{.6pt}}{#1}}
\renewcommand{\bar}{\thickbar}
\renewcommand{\hat}{\widehat}
\renewcommand*{\tilde}{\widetilde}

\makeatletter
\let\oldoint\oint
\renewcommand{\oint}{\oldoint\nolimits}
\makeatother

\newcommand*{\N}{\ensuremath{\mathbb{N}}}
\newcommand*{\Z}{\ensuremath{\mathbb{Z}}}

\newcommand*{\R}{\ensuremath{\mathbb{R}}}

\newcommand*{\Zd}{\ensuremath{\mathbb{Z}^d}}
\newcommand*{\Rd}{\ensuremath{\mathbb{R}^d}}
\newcommand{\E}{\mathbb{E}}
\renewcommand{\P}{\ensuremath{\mathbb{P}}}

\renewcommand{\O}{\ensuremath{\mathcal{O}}}

\newcommand{\X}{\ensuremath{\mathcal{X}}}

\newcommand{\A}{\mathcal{A}}
\renewcommand{\S}{\mathcal{S}}

\renewcommand{\a}{\mathbf{a}}
\renewcommand{\b}{\ensuremath{\mathbf{b}}}
\newcommand{\f}{\mathbf{f}}
\newcommand{\g}{\mathbf{g}}
\newcommand{\h}{\mathbf{h}}

\newcommand{\ahom}{\bar{\a}}

\newcommand{\sol}{\mathrm{sol}}
\newcommand{\sym}{\mathrm{sym}}

\newcommand*{\Id}{\ensuremath{\mathrm{I}_d}}

\newcommand{\Rsymp}{\R^{d\times d}_{\mathrm{sym},+}}

\newcommand{\Lsol}{L^2_{\mathrm{sol}}}
\newcommand{\Lsolo}{L^2_{\mathrm{sol,0}}}

\newcommand{\eps}{\varepsilon}
\newcommand{\ep}{\eps}

\newcommand{\nf}{\nicefrac}
\newcommand{\qand}{\quad \mbox{and} \quad}

\renewcommand{\subset}{\subseteq}
\renewcommand{\supset}{\supseteq}

\newcommand{\gmcdelta}{\textcolor{blue}{\overline \delta}}

\makeatletter
\newcommand{\negphantom}{\v@true\h@true\negph@nt}
\newcommand{\neghphantom}{\v@false\h@true\negph@nt}
\newcommand{\negph@nt}{\ifmmode\expandafter\mathpalette
	\expandafter\mathnegph@nt\else\expandafter\makenegph@nt\fi}
\newcommand{\makenegph@nt}[1]{%
	\setbox\z@\hbox{\color@begingroup#1\color@endgroup}\finnegph@nt}
\newcommand{\finnegph@nt}{%
	\setbox\tw@\null
	\ifv@ \ht\tw@\ht\z@\dp\tw@\dp\z@\fi
	\ifh@\wd\tw@-\wd\z@\fi\box\tw@}
\newcommand{\mathnegph@nt}[2]{%
	\setbox\z@\hbox{$\m@th #1{#2}$}\finnegph@nt}
\makeatother

\newcommand{\Hminusul}{\hat{\phantom{H}}\negphantom{H}\underline{H}^{-1}}
\newcommand{\Hminusuls}[1]{\hat{\phantom{H}}\negphantom{H}\underline{H}^{#1}}

\newcommand{\Wminusul}[2]{\hat{\phantom{W}}\negphantom{W}\underline{W}^{#1,#2}}

\newcommand{\Besov}[3]{\mathring{\phantom{B}}\negphantom{B}\underline{B}^{#1}_{#2,#3}}

\def\Xint#1{\mathchoice
	{\XXint\displaystyle\textstyle{#1}}%
	{\XXint\textstyle\scriptstyle{#1}}%
	{\XXint\scriptstyle\scriptscriptstyle{#1}}%
	{\XXint\scriptscriptstyle\scriptscriptstyle{#1}}%
	\!\int}
\def\XXint#1#2#3{{\setbox0=\hbox{$#1{#2#3}{\int}$}
		\vcenter{\hbox{$#2#3$}}\kern-.5\wd0}}

\def\fint{\Xint-}

\newcommand{\avsum}{\mathop{\mathpalette\avsuminner\relax}\displaylimits}
\makeatletter
\newcommand\avsuminner[2]{%
	{\sbox0{$\m@th#1\sum$}%
		\vphantom{\usebox0}%
		\ooalign{%
			\hidewidth
			\smash{\,\rule[.23em]{8.8pt}{1.1pt} \relax}%
			\hidewidth\cr
			$\m@th#1\sum$\cr
		}%
	}%
}
\makeatother

\makeatletter
\newcommand\avsuminnerr[2]{%
	{\sbox0{$\m@th#1\sum$}%
		\vphantom{\usebox0}%
		\ooalign{%
			\hidewidth
			\smash{\,\rule[.23em]{6pt}{0.7pt} \relax}%
			\hidewidth\cr
			$\m@th#1\sum$\cr
		}%
	}%
}
\makeatother

\newcommand{\cu}{\square}

\newcommand{\avsumcube}[3]{\avsum_{#1 \in 3^{#2}\Zd\cap \cu_{#3}}}
\newcommand{\maxcube}[3]{\max_{#1 \in 3^{#2}\Zd\cap \cu_{#3}}}

\newcommand{\spx}{\triangle}

\DeclareMathAlphabet{\mathmybb}{U}{bbold}{m}{n}
\newcommand{\indc}{{\mathmybb{1}}}

\renewcommand*\showkeyslabelformat[1]{}
\title[Anomalous scaling limit]{Anomalous scaling limit of a Brownian particle in a log-correlated potential}

\author[S. Armstrong]{Scott Armstrong}
\address[S. Armstrong]{CNRS \& Laboratoire Jacques-Louis Lions, Sorbonne Universit\'e, and Courant Institute of Mathematical Sciences, New York University}
\email{scottnarmstrong@gmail.com}

\author[A. Bou-Rabee]{Ahmed Bou-Rabee}
\address[A. Bou-Rabee]{Department of Mathematics, University of Pennsylvania}
\email{ahmedmb@sas.upenn.edu}

\author[T. Kuusi]{Tuomo Kuusi}
\address[T. Kuusi]{Department of Mathematics and Statistics, University of Helsinki}
\email{tuomo.kuusi@helsinki.fi}

\date{}

\begin{document}

\begin{abstract}
We study a Brownian particle in~$\Rd$,~$d\geq2$, with drift given by the gradient of a log-correlated Gaussian potential. At weak disorder, we prove convergence to a scaling limit whose law is singular with respect to Brownian motion and whose invariant measure is given by Gaussian multiplicative chaos. The proof uses a renormalization group induction: at each scale, the elliptic operator is well approximated by a Laplacian with effective diffusivity decaying as a power of scale. We compute this exponent exactly in dimension two.
\end{abstract}

\maketitle

\begingroup\renewcommand{\thefigure}{1}
\begin{figure}[h]
\centering
\includegraphics[width=\textwidth]{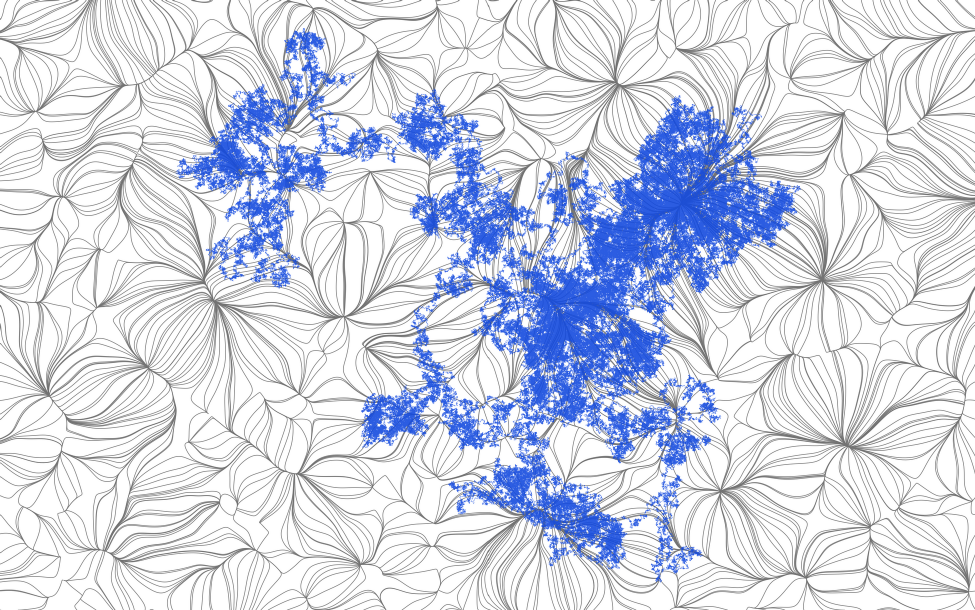}
\caption{A Brownian particle advected by the gradient of a log-correlated random potential.}
\end{figure}
\endgroup

\smallskip

\tableofcontents

\section{Introduction}

\subsection{Brownian motion in a log-correlated potential}

We study the long-time behavior of a Brownian particle in~$\Rd$,~$d\geq2$, whose drift is the gradient of a log-correlated random potential~$\log\a$, defined in~\eqref{e.intro.anchored.coefficient} below. The particle position~$X_t$ is described by the stochastic differential equation
\begin{equation}
	dX_t
	=
	\nabla\log\a(X_t)\,dt+\sqrt2\,dW_t\,,
	\qquad
	X_0=x\in\Rd\,,
	\label{e.sde.intro}
\end{equation}
where~$W$ is a standard Brownian motion on~$\Rd$. The diffusion~$X$ has generator~$\a^{-1}\nabla\cdot(\a\nabla)=\Delta+\nabla\log\a\cdot\nabla$ and reversible measure~$\a(x)\,dx$.

\smallskip

In the language of statistical physics,~$X$ is an overdamped particle at unit temperature in the energy landscape~$-\log\a$. The drift pushes it toward the large values of~$\a$, where the reversible measure puts its mass, and to travel it must cross the barriers of this landscape. For a bounded stationary potential with rapidly decaying correlations, homogenization theory shows that these barriers average out: on large scales the particle behaves like a Brownian motion with a deterministic, positive effective diffusivity. In the model considered here, the potential has independent fluctuations of equal strength at every length scale. As more scales contribute, their effects accumulate and the effective diffusivity tends to zero. The particle is then \textbf{subdiffusive}: the time it needs to travel a large distance~$r$ grows faster than~$r^2$.

\smallskip

Diffusion in random force fields was studied in the physics literature by~\citet{Fisher},~\citet*{KLY},~\citet*{FFQSS} and~\citet*{BCGLD}; see~\citet*{BG} for a review. For gradient drifts,~\citet*{BCGLD} found that the strength of the logarithmic correlations of the potential is not renormalized at any order of perturbation theory. Consequently, the time the particle needs to travel a distance~$r$ grows like a power of~$r$ whose exponent depends continuously on this strength. When the strength is small, their renormalization group computation predicts that the exponent exceeds~$2$ by an amount proportional to it, and they suggested that this excess has no corrections of higher order~\citep[Section~4.3.2.3]{BG}. The same formula appears in the later work of~\citet*[equation~(3)]{CastilloLeDoussal}, who studied how the exponent changes when the strength is large.

\smallskip

The existence of a scaling limit for a diffusion driven by the gradient of a smoothed planar Gaussian free field was also posed as an open problem by~\citet*[Problem~16.13]{BiskupSurvey}.

\smallskip

In this paper, we prove that, when the fluctuations of the potential are small, the diffusion~$X$ has a unique scaling limit. We also show that the effective diffusivity at length scale~$r$, characterized by the Brownian scaling limit~\eqref{e.intro.effective.diffusivity} below, is bounded above and below by negative powers of~$r$. In dimension two it is exactly a power of~$r$ with the exponent predicted by~\citet*{BCGLD}; in every dimension its exponent agrees with that prediction up to an error of higher order in the size of the fluctuations. The scaling limit is a continuous, conservative, strong Markov process. Its path law is mutually singular with every mixture of Brownian laws with arbitrary covariance matrices, and its reversible measure is a multiplicative chaos measure.

\subsection{The model}
\label{ss.assumptions}

We write~$\P$ and~$\E$ for the law and expectation of the environment~$\omega$, a realization of the fields. We call the law of the trajectories in a fixed environment   \textbf{quenched} and its average over the environment \textbf{annealed}.

\smallskip

Fix~$d\geq2$ and the \textbf{disorder strength}~$\gmcdelta\in(0,\nf12]$, and let~$\g_0,\g_1,\g_2,\ldots$ be independent random fields such that, for every~$k\in\N_0$,
\begin{equation}
	\g_k
	\quad\text{and}\quad
	\g_0(3^{-k}\,\cdot)
	\quad\text{have the same law}\,.
	\label{e.same.law.shift}
\end{equation}
For~$m\in\Z$, let~$\cu_m\coloneqq(-\frac12\,3^m,\frac12\,3^m)^d$ be the cube of side~$3^m$. We assume the following.
\begin{enumerate}[label=(\textbf{g\arabic*})]
\setcounter{enumi}{0}

\item
\label{a.g1}
The field~$\g_0$ is stationary, has mean zero, and has range of dependence~$\sqrt d$: the sigma-fields generated by its restrictions to measurable sets~$U,V\subseteq\Rd$ are independent whenever~$|x-y|_2\geq\sqrt d$ for every~$x\in U$ and~$y\in V$.

\item
\label{a.g2}
Almost surely,~$\g_0\in C^{1,1}_{\mathrm{loc}}(\Rd)$, and
\begin{equation}
	\E\left[\exp\left(\gmcdelta^{-2}\left(\|\g_0\|_{L^\infty(\cu_0)}+\|\nabla\g_0\|_{L^\infty(\cu_0)}+\|\nabla^2\g_0\|_{L^\infty(\cu_0)}\right)^2\right)\right]\leq2\,.
	\label{e.coefficient.field.regularity}
\end{equation}
\item
\label{a.g3}
The law of~$\g_0$ is invariant under signed permutations of the coordinates and under~$\g_0\mapsto-\g_0$.

\item
\label{a.g4}
The number~$\tau^2=\log\E[\exp(\g_0(0))]$ is positive.
\end{enumerate}
We say that a result holds \textbf{at weak disorder} if it holds whenever~$\gmcdelta\leq\gmcdelta_0$, for some~$\gmcdelta_0>0$ depending only on~$d$ and the indicated parameters.

\smallskip

The coefficient is
\begin{equation}
	\a(x)
	\coloneqq
	\exp\left(
	\sum_{k=0}^\infty(\g_k(x)-\g_k(0))
	\right)\,.
	\label{e.intro.anchored.coefficient}
\end{equation}
The field~$\log\a$ is log-correlated: the variance of~$\log\a(x)-\log\a(y)$ lies between two positive multiples of~$\log|x-y|$ for all sufficiently large~$|x-y|$.

\smallskip

The effective diffusivity at a given scale is defined through an \textbf{infrared cutoff}. For~$L\in\N_0$ we call~$L$ the \textbf{cutoff},~$3^L$ the \textbf{cutoff scale}, and the \textbf{cutoff coefficient} is
\begin{equation}
	\a_L(x)
	\coloneqq
	\exp\left(
	\sum_{k=0}^L(\g_k(x)-\tau^2)
	\right)\,.
	\label{e.intro.cutoff.coefficients}
\end{equation}
This coefficient is stationary, has mean one and has range of dependence~$\sqrt d\,3^L$. Let~$X^{(L)}$ be the diffusion defined by~\eqref{e.sde.intro} with~$\a_L$ in place of~$\a$; its generator is
\begin{equation}
	\a_L^{-1}\nabla\cdot(\a_L\nabla)
	=
	\Delta+\nabla\log\a_L\cdot\nabla\,.
	\label{e.intro.cutoff.generator}
\end{equation}
For fixed~$L$, with~$X^{(L)}$ started at the origin, the \textbf{effective diffusivity}~$\ahom_L$ is characterized by the annealed scaling limit
\begin{equation}
	\bigl(\varepsilon X^{(L)}_{t/\varepsilon^2}\bigr)_{t\geq0}
	\longrightarrow
	\bigl(\sqrt{2\ahom_L}\,W_t\bigr)_{t\geq0}
	\quad\text{in law as }\varepsilon\downarrow0\,.
	\label{e.intro.effective.diffusivity}
\end{equation}
Its precise definition is given in~\eqref{e.ahom.ahomstar.defs}. Theorem~\ref{t.renormalized.diffusivities} gives power-law bounds for~$\ahom_L$ and the exact formula~$\ahom_L=e^{-(L+1)\tau^2}$ in dimension two.

\subsection{The scaling limit}

For the limiting process, we take~$(\g_j)_{j\in\Z}$ independent, with~$\g_j$ having the law of~$\g_0(3^{-j}\,\cdot)$ for every~$j\in\Z$; see Section~\ref{s.scaling.limit}.

\begin{theoremA}[Scaling limit]
\label{t.A}
Assume~\ref{a.g1}--\ref{a.g4}. There exists~$\gmcdelta_0(d)>0$ such that, if~$0<\gmcdelta\leq\gmcdelta_0$, then there is a family of quenched path laws~$(\mathbf P_x)_{x\in\Rd}$, measurably determined by~$(\g_j)_{j\in\Z}$, which almost surely defines a continuous strong Markov process~$Z$ on~$\Rd$, started at~$x$ under~$\mathbf P_x$, with corresponding expectations~$\mathbf E_x$, such that the following statements hold.
\begin{enumerate}[label=(\roman*)]
\item \underline{\emph{Convergence}}: For every~$x\in\Rd$, as the integer~$N$ tends to infinity, the annealed law of the process~$(3^{-N}X_{3^{2N}t/\ahom_N})_{t\geq0}$ with~$X_0=3^Nx$ converges to the annealed law of~$Z$ started at~$x$.
\item \underline{\emph{Reversible measure}}: Almost surely,~$Z$ is reversible with respect to a measure~$\mu$ that is singular with respect to Lebesgue measure.
\item \underline{\emph{Anomalous scaling}}: There exist~$C<\infty$ and~$\eta>0$, depending only on~$d$ and the law of~$\g_0$, such that, if~$\sigma_r$ denotes the first time at which~$Z$ leaves the cube of side~$r$ centered at its starting point, then, for every~$x\in\Rd$ and~$r\in(0,1]$,
\begin{equation*}
	\E\bigl[\mathbf E_x[\sigma_r]\bigr]\leq Cr^{2+\eta}\,.
\end{equation*}
For every~$x\in\Rd$, almost surely,~$\mathbf P_x$-almost every path satisfies
\begin{equation*}
	\limsup_{t\downarrow0}\,t^{-\gamma}|Z_t-Z_0|=\infty
	\qquad\text{for every }\gamma>\nf1{(2+\eta)}\,.
\end{equation*}
In dimension two, one can take~$\eta=\tau^2/\log3$.
\item \underline{\emph{Non-Gaussian marginals}}: Almost surely, for every~$x\in\Rd$ and~$t>0$, the law of~$Z_t$ under~$\mathbf P_x$ is absolutely continuous with respect to~$\mu$ and is not Gaussian.
\end{enumerate}
\end{theoremA}

Theorem~\ref{t.scaling.limit}(i) gives a quenched version of~(i): the rescaled processes can be defined in one environment so that, almost surely, their quenched laws converge to those of~$Z$, uniformly for starting points in compact sets.

\smallskip

Standard Brownian motion needs mean time of order~$r^2$ to leave a cube of side~$r$, whereas~(iii) bounds the annealed mean exit time of~$Z$ by~$Cr^{2+\eta}$ for~$r\leq1$. Consequently, the path law of~$Z$ is mutually singular with every mixture of Brownian laws with arbitrary covariance matrices (Theorem~\ref{t.scaling.limit}(iii)).

\subsection{Homogenization at the cutoff scale}

We compare the Dirichlet problem for the normalized cutoff coefficient~$\ahom_L^{-1}\a_L$, on cubes of side at least~$3^L$, with the Dirichlet problem for the Laplacian. For~$L,M\in\N_0$ with~$L\leq M$, define~$A_{L,M}(x)\coloneqq\ahom_L^{-1}\a_L(3^Mx)$ for~$x\in\cu_0$. The coefficient~$A_{L,M}$ is stationary, has homogenized matrix equal to the identity, and has range of dependence~$\sqrt d\,3^{L-M}$. Given~$f\in L^2(\cu_0)$ and~$h\in H^2(\cu_0)$, let~$u_{L,M},u_{\mathrm{hom}}\in H^1(\cu_0)$ solve
\begin{equation*}
	\left\{
	\begin{aligned}
		-\nabla\cdot(A_{L,M}\nabla u_{L,M})&=f
		&&\text{in }\cu_0\,,\\
		u_{L,M}&=h
		&&\text{on }\partial \cu_0\,,
	\end{aligned}
	\right.
	\qquad
	\left\{
	\begin{aligned}
		-\Delta u_{\mathrm{hom}}&=f
		&&\text{in }\cu_0\,,\\
		u_{\mathrm{hom}}&=h
		&&\text{on }\partial \cu_0\,,
	\end{aligned}
	\right.
\end{equation*}
and define the \textbf{homogenization error}
\begin{equation*}
	\mathcal D_{L,M}
	\coloneqq
	\|u_{L,M}-u_{\mathrm{hom}}\|_{L^2(\cu_0)}
	+
	\|\nabla u_{L,M}-\nabla u_{\mathrm{hom}}\|_{H^{-1}(\cu_0)}
	+
	\|A_{L,M}\nabla u_{L,M}-\nabla u_{\mathrm{hom}}\|_{H^{-1}(\cu_0)}\,.
\end{equation*}
The negative norms in~$\mathcal D_{L,M}$ are defined in~\eqref{e.ordinary.negative.norm.convention}.

\begin{theoremA}[Quantitative homogenization of the Dirichlet problem]
\label{t.B}
Assume~\ref{a.g1}--\ref{a.g4}. The following statements hold.
\begin{enumerate}[label=(\roman*)]
\item \underline{\emph{Estimate uniform in the cutoff}}: For every~$\vartheta\in(0,1)$ and~$q\in[1,\infty)$, there exist~$\gmcdelta_0(\vartheta,q,d)>0$ and~$C(\vartheta,q,d)<\infty$ such that, if~$\gmcdelta\leq\gmcdelta_0$, then there exist random variables~$\mathcal Z_{L,M}\geq1$ satisfying
\begin{equation*}
	\sup_{0\leq L\leq M}
	\|\mathcal Z_{L,M}\|_{L^q(\P)}
	\leq C
\end{equation*}
such that, almost surely, for every~$L,M\in\N_0$ with~$L\leq M$,~$f\in L^2(\cu_0)$ and~$h\in H^2(\cu_0)$,
\begin{equation*}
	\mathcal D_{L,M}
	\leq
	\mathcal Z_{L,M}\,\gmcdelta^\vartheta
	\left(\|f\|_{L^2(\cu_0)}+\|h\|_{H^2(\cu_0)}\right)\,.
\end{equation*}
If~$f=0$, then the factor~$\gmcdelta^\vartheta$ may be replaced by~$C\gmcdelta$.

\item \underline{\emph{Algebraic rate at a fixed cutoff}}: There exists~$\alpha_{\mathrm{hom}}(d)>0$ such that, for every~$L\in\N_0$ and~$q\in[1,\infty)$, there exists a random variable~$\mathcal Y_{L,q}\geq1$ with~$\|\mathcal Y_{L,q}\|_{L^q(\P)}\leq C(L,q,d,\gmcdelta)$ such that, almost surely, for every~$M\in\N_0$ with~$M\geq L$,~$f\in L^2(\cu_0)$ and~$h\in H^2(\cu_0)$,
\begin{equation*}
	\mathcal D_{L,M}
	\leq
	\mathcal Y_{L,q}\,3^{-\alpha_{\mathrm{hom}}(M-L)}
	\left(\|f\|_{L^2(\cu_0)}+\|h\|_{H^2(\cu_0)}\right)\,.
\end{equation*}
\end{enumerate}
\end{theoremA}

The estimate in~(i) includes the case~$M=L$, and its constants do not depend on the cutoff. In~(ii), the cutoff~$L$ is fixed, the rate is algebraic in~$M-L$, and no smallness of~$\gmcdelta$ is assumed.

\smallskip

For~$f=1$ and~$h=0$, the solution~$u_{L,M}$ is the mean exit time from~$\cu_0$ of the diffusion with generator~$\nabla\cdot(A_{L,M}\nabla)$, as a function of its starting point; for~$f=0$, it is the expected value of~$h$ at the exit position. Theorem~\ref{t.B} compares these quantities with the corresponding quantities for the diffusion with generator~$\Delta$.

\smallskip

The law of~$X$ is determined by its Dirichlet form~$\int\a|\nabla u|^2\,dx$ on~$L^2(\a\,dx)$, and constructing the scaling limit requires the convergence of the Dirichlet forms of the rescaled processes. The bound in Theorem~\ref{t.B}(i) is of order~$\gmcdelta^\vartheta$ and does not tend to zero as the cutoff tends to infinity, so it does not give this convergence. For the convergence argument, see Section~\ref{sec:mfd-uniqueness}.

\subsection{Large-scale regularity and a Liouville theorem}

At weak disorder, we obtain large-scale H\"older estimates with constants independent of the cutoff, despite the absence of uniform ellipticity bounds.

\smallskip

We call~$u$~$\a_L$-harmonic if~$\nabla\cdot(\a_L\nabla u)=0$. We use the notation~$(u)_U\coloneqq\fint_Uu$ and~$\|u\|_{\underline L^p(U)}\coloneqq(\fint_U|u|^p)^{\nf1p}$ for a measurable set~$U$ of finite positive volume and~$p\in[1,\infty)$. We read~$\a_\infty$ as~$\a$.

\begin{theoremA}[Large-scale H\"older regularity and a Liouville theorem]
\label{t.C}
Assume~\ref{a.g1}--\ref{a.g4}. There exist~$\gmcdelta_0(d)>0$,~$C_0(d)<\infty$ and~$C(d)<\infty$ such that, if~$\gmcdelta\leq\gmcdelta_0$, then the following statements hold with
\begin{equation*}
	\gamma_{\mathrm{reg}}
	\coloneqq
	1-C_0\gmcdelta|\log\gmcdelta|^{\nf12}\,.
\end{equation*}
\begin{enumerate}[label=(\roman*)]
\item \underline{\emph{Large-scale H\"older estimate}}: Let~$\gamma\in[\nf12,\gamma_{\mathrm{reg}}]$. For every~$L\in\N_0\cup\{\infty\}$ and~$m\in\N$, there exists an~$\N$-valued random variable~$\mathcal L_L(\gamma,m)$ such that, for every~$k\in\N$,
\begin{equation}
	\P[\mathcal L_L(\gamma,m)>k]
	\leq
	C\exp\left(
	-\frac{(1-\gamma)^2(k-C)_+}
	{C\gmcdelta^2|\log\gmcdelta|}
	\right)\,.
	\label{e.intro.regularity.minimal.scale}
\end{equation}
Moreover, almost surely, for every~$L$ and~$m$, if~$u\in H^1(\cu_m)$ is~$\a_L$-harmonic in~$\cu_m$, then, for every~$n\in\N_0$ satisfying~$n\leq m-\mathcal L_L(\gamma,m)$ and every~$z\in3^n\Zd$ with~$z+\cu_n\subseteq\cu_{m-1}$,
\begin{equation}
	\left\|
	u-(u)_{z+\cu_n}
	\right\|_{\underline L^2(z+\cu_n)}
	\leq
	C3^{-\gamma(m-n)}
	\left\|
	u-(u)_{\cu_m}
	\right\|_{\underline L^2(\cu_m)}
	\label{e.intro.large.scale.Holder}
\end{equation}
and
\begin{equation}
	\left\|\a_L^{\nf12}\nabla u
	\right\|_{\underline L^2(z+\cu_n)}
	\leq
	C3^{(1-\gamma)(m-n)}
	\left\|\a_L^{\nf12}\nabla u
	\right\|_{\underline L^2(\cu_m)}\,.
	\label{e.intro.large.scale.energy}
\end{equation}

\item \underline{\emph{Liouville theorem}}: Almost surely, for every~$L\in\N_0\cup\{\infty\}$, if~$u\in H^1_{\mathrm{loc}}(\Rd)$ satisfies
\begin{equation*}
	-\nabla\cdot(\a_L\nabla u)=0
	\quad\text{in }\Rd
	\quad\text{and}\quad
	\liminf_{R\to\infty}
	R^{-\gamma_{\mathrm{reg}}}
	\inf_{c\in\mathbb R}
	\|u-c\|_{\underline L^2(B_R)}=0\,,
\end{equation*}
then~$u$ is constant.
\end{enumerate}
\end{theoremA}

The estimate~\eqref{e.intro.large.scale.energy} controls concentration of the energy measure~$\a_L|\nabla u|^2\,dx$. At the scales covered by the theorem, the energy in~$z+\cu_n$ is bounded by~$C^2 3^{-[d-2(1-\gamma)](m-n)}$ times the energy in~$\cu_m$. In this sense, the energy behaves as though its dimension were at least~$d-2(1-\gamma)$. Taking~$\gamma=\gamma_{\mathrm{reg}}$ bounds the deficit from~$d$ by~$2C_0\gmcdelta|\log\gmcdelta|^{\nf12}$, so weak disorder permits only a small loss of dimension at these scales. This control is used to prove uniqueness of the limiting Dirichlet form in Section~\ref{sec:mfd-uniqueness}.

\subsection{Overview of the proofs}

\begingroup\runinsubsubsections\runinparagraphs[0pt]

We give an outline of the arguments leading to Theorems~\ref{t.A},~\ref{t.B} and~\ref{t.C}. The proof is divided into five main steps.

\begin{enumerate}
\item \textbf{Deterministic coarse-graining (Section~\ref{s.coarse.graining}).} Following~\citet{AK.HC}, we develop elliptic estimates in terms of variational problems on cubes. These estimates describe the operator at each spatial scale through averaged energies and give quantitative harmonic approximation without requiring pointwise ellipticity bounds uniform in the cutoff.

\item \textbf{The renormalization group argument (Sections~\ref{s.infrared.cutoffs}--\ref{s.homogenized.diffusivities}).} We prove by induction that, on cubes of side~$3^m$, the operator with coefficient~$\a_m$ is well approximated by~$\ahom_m\Delta$. The induction alternates between averaging over smaller cubes, which improves the approximation, and adding the missing fields, which restores the coefficient at the new cutoff. We control the relative error uniformly in~$m$ and then determine the decay of~$\ahom_m$.

\item \textbf{Large-scale regularity (Section~\ref{s.Holder.regularity}).} We combine harmonic approximation with an excess-decay iteration to prove the H\"older and energy estimates~\eqref{e.intro.large.scale.Holder}--\eqref{e.intro.large.scale.energy} of Theorem~\ref{t.C}. Approximation improves on good scales, and control of the frequency of bad scales allows us to iterate, uniformly in the cutoff. The Dirichlet homogenization estimates of Theorem~\ref{t.B} follow directly by applying coarse-graining on the whole cube.

\item \textbf{Uniqueness of the limiting Dirichlet form (Section~\ref{sec:mfd-uniqueness}).} We prove that the Dirichlet forms of the rescaled processes on cubes have a unique limit, which we call the \textbf{limit form}. Its precise definition is given in~\eqref{eq:mfd-18}. The regularity estimates control the effect of resampling fine-scale fields and give compactness of the minimal energies on cubes with prescribed boundary values. We identify their limits by adapting the comparison and local replacement arguments of~\citet*{GM,DGUniqueness,DGComparison} for Liouville quantum gravity metrics, with energy minimizers in place of geodesics.

\item \textbf{Construction of the scaling limit (Sections~\ref{sec:mfd-process-convergence}--\ref{lim:sec-process-measure}).} We combine the limit form with the limit~$\mu$ of the reversible measures to construct the limiting diffusion. Convergence of the resolvents first gives convergence in probability of the quenched path laws. To obtain almost-sure convergence, we adapt the comparison argument of~\citet[Section~3.4]{DGComparison}, strengthened by~\citet*[Lemma~3.11]{Devlin2024}. We then prove the properties of the limit stated in Theorem~\ref{t.A}.
\end{enumerate}

\noindent Figure~\ref{f.dependence} shows the dependence of the main results. Below we explain each step in more detail. The moment estimates are collected in Appendix~\ref{s.orlicz}, the concentration inequalities in Appendix~\ref{app.concentration}, and the Besov norm comparisons in Appendix~\ref{s.Besov.norm.comparisons}.

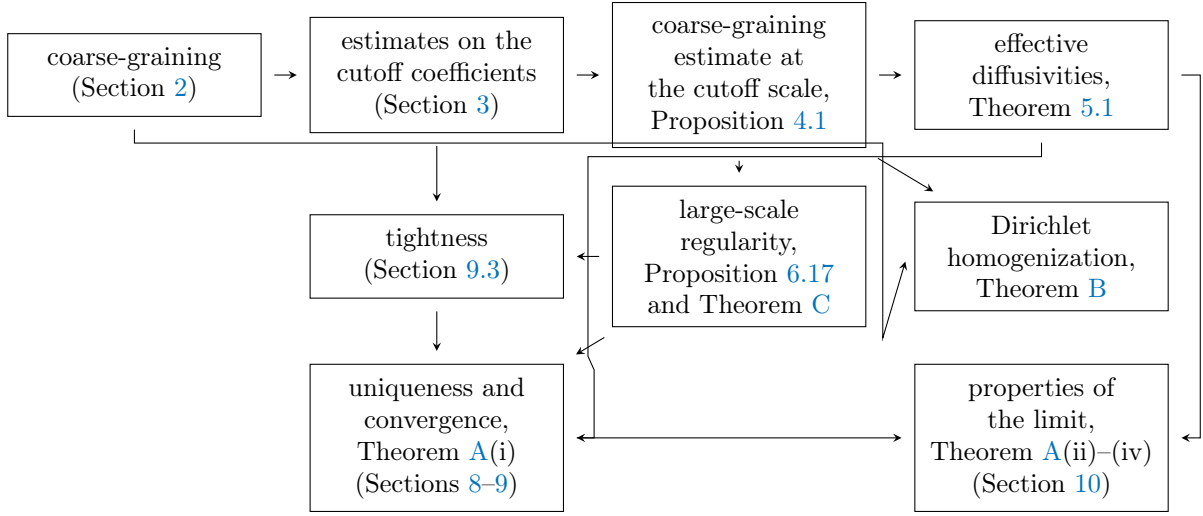
\begin{figure}[ht]
\centering
\tikzset{external/export next=false}
\begin{tikzpicture}[x=1cm,y=1.2cm,every node/.style={draw,rectangle,align=center,font=\small,text width=3cm,inner sep=5pt,outer sep=3pt},>=stealth,shorten >=2pt,shorten <=2pt]
\node (s2) at (0,0) {coarse-graining\\ (Section~\ref{s.coarse.graining})};
\node (s3) at (4,0) {estimates on the\\ cutoff coefficients\\ (Section~\ref{s.infrared.cutoffs})};
\node (s4) at (8,0) {coarse-graining estimate at\\ the cutoff scale,\\ Proposition~\ref{p.coarse.grained.bound}};
\node (s5) at (12,0) {effective\\ diffusivities,\\ Theorem~\ref{t.renormalized.diffusivities}};
\node (s6) at (8,-2) {large-scale regularity,\\ Proposition~\ref{p.cutoff.Holder.regularity}\\ and Theorem~\ref{t.C}};
\node (s6b) at (12,-2) {Dirichlet\\ homogenization,\\ Theorem~\ref{t.B}};
\node (s8) at (4,-2) {tightness\\ (Section~\ref{s.tightness})};
\node (s9) at (4,-4) {uniqueness and convergence,\\ Theorem~\ref{t.A}(i)\\ (Sections~\ref{sec:mfd-uniqueness}--\ref{sec:mfd-process-convergence})};
\node (s10) at (12,-4) {properties of the limit,\\ Theorem~\ref{t.A}(ii)--(iv)\\ (Section~\ref{lim:sec-process-measure})};
\draw[->] (s2) -- (s3);
\draw[->] (s3) -- (s4);
\draw[->] (s4) -- (s5);
\draw[->] (s4) -- (s6);
\draw[->] (s6) -- (s8);
\draw[->] (s3) -- (s8);
\draw[->] (s2.south) -- (0,-0.75) -- (9.9,-0.75) -- (9.9,-2.9) -- (s6b.west);
\draw[->] (s4) -- (s6b);
\draw[->] (s5.south) -- (12,-0.9) -- (6,-0.9) -- (6,-3.1) -- ([xshift=0.3cm,yshift=0.9cm]s9.east) -- ([xshift=0.3cm]s9.east) -- (s9.east);
\draw[->] (s8) -- (s9);
\draw[->] (s6) -- (s9);
\draw[->] (s9) -- (s10);
\draw[->] (s5.east) -- (14.1,0) -- (14.1,-4) -- (s10.east);
\end{tikzpicture}
\caption{Logical dependence of the main results. An arrow means that the first result is used in the proof of the second.}
\label{f.dependence}
\end{figure}

\subsubsection{Deterministic coarse-graining (Section~\ref{s.coarse.graining})}

The first step is a deterministic theory for the elliptic operator~$\nabla\cdot\a_L\nabla$. Since the pointwise ellipticity bounds for~$\a_L$ deteriorate as~$L$ increases, we need estimates that describe its behavior through averages at the scale under consideration. Following~\citet{AK.HC}, we associate two \textbf{coarse-grained matrices} with each cube, defined in~\eqref{e.coarse.grained.matrices}. The first records the minimal energy for affine boundary values; the inverse of the second records the energy of the Neumann solution with the boundary flux of a constant vector field. For a constant coefficient, both matrices equal that coefficient. Their agreement measures how well the operator is approximated by a constant-coefficient operator on the cube.

The coarse-grained matrices on subcubes enter the constants in Poincar\'e and Caccioppoli inequalities. When they are close to a common constant matrix, the same estimates give quantitative harmonic approximation. Thus both ellipticity and homogenization are measured by averaged energies, rather than by the pointwise bounds on~$\a_L$. The coefficient field in Section~\ref{s.coarse.graining} is arbitrary and fixed; randomness enters only when we apply this theory to the cutoff coefficients.

\paragraph{Coarse-grained sensitivity estimates.}
To pass from one cutoff to the next, we must compare the coarse-grained matrices before and after multiplying the coefficient by a slowly varying field. Lemma~\ref{l.sensitivity.general} uses the variational formulas to give sensitivity estimates in terms of the size of this multiplicative change, rather than the pointwise ellipticity ratio of the original coefficient. This is the comparison needed to add successive fields while measuring errors relative to the decreasing effective diffusivity.

\paragraph{Equations with nonzero right-hand side.}
The coarse-grained matrices give estimates for harmonic functions. To obtain homogenization estimates for general Dirichlet problems, we also need to treat solutions with a nonzero right-hand side. For this, we prove the deterministic comparison estimate~\eqref{e.general.coarse.graining.estimate.ASD}: solutions of the heterogeneous and constant-coefficient equations with the same forcing and boundary values have close gradients and fluxes, in negative norms, when the coarse-grained matrices are close to the constant coefficient and the forcing varies little on the subcubes used for coarse-graining. This is the symmetric-coefficient case of~\citet*[Proposition~2.21]{ASD}. Applying this estimate with the renormalization bounds on the coarse-grained matrices gives Theorem~\ref{t.B}.

\subsubsection{The renormalization group argument (Sections~\ref{s.infrared.cutoffs}--\ref{s.homogenized.diffusivities})}

As the cutoff increases, the coefficient acquires new fields at larger scales. Homogenizing a fixed cutoff on much larger cubes therefore does not by itself describe the coefficient at the scale of the new cutoff. We instead follow the coarse-grained matrices together with the effective diffusivity~$\ahom_m$, measuring their errors relative to~$\ahom_m$ at each scale. The main estimate says that, on~$\cu_m$, the coarse-grained matrices of~$\a_m$ remain close to~$\ahom_m\Id$ in relative terms, uniformly in~$m$, even as~$\ahom_m$ tends to zero.

A key point is that the cube and the cutoff are at the same scale. At a fixed cutoff~$L$, finite range of dependence allows us to average over many independent cubes once we consider a cube much larger than~$\cu_L$. But the coefficient at this larger scale contains fields absent from~$\a_L$. To return to an estimate at the cutoff scale, we must add these fields and change the normalization from~$\ahom_L$ to the effective diffusivity at the new cutoff. The induction must control both changes while keeping the same bound on the relative error.

\paragraph{Contraction by averaging.}
Fix the cutoff~$L$ and consider~$\a_L$ on a much larger cube. Its finite range of dependence makes the coarse-grained quantities on well-separated subcubes independent. Subadditivity and concentration then bring the averages of these quantities close to deterministic matrices, and the coarse-graining estimates convert this into an improvement of the relative error (Proposition~\ref{p.homogenization.step}). Passing through enough scales makes this improvement strong enough to absorb the constants introduced elsewhere in the induction.

This contraction is what allows the induction to close with the same constant at each scale. The deterministic estimates introduce constants, and adding the new fields enlarges the bound on the relative error. An estimate that merely preserves this bound under averaging would therefore deteriorate at every step. We use the scale separation to obtain a fixed contraction before adding the fields, leaving enough room for these losses.

\paragraph{Adding the coarser fields.}
After averaging, we are on a cube at scale~$3^m$ but still working with~$\a_L$. To obtain the coefficient~$\a_m$, we add the fields~$\g_{L+1},\ldots,\g_m$. Their slow variation on the smaller cubes allows us to apply the sensitivity estimates of Lemma~\ref{lemma:infrared.approx.cutoffs} and bound the change relative to the effective diffusivity. This step enlarges the bound on the relative error, but at weak disorder its cost is smaller than the improvement obtained by averaging. Independence and the scaling identity~\eqref{e.same.law.shift} put us in the same setting at the new cutoff.

On a small cube, a field of much larger scale is nearly constant, so multiplying the coefficient by its exponential approximately multiplies both coarse-grained matrices by the same factor. The sensitivity estimates control the error in this approximation. Across the larger cube, moment bounds control the multiplicative change of the coefficient, and the sensitivity estimate combines these bounds with the contraction obtained by averaging. Thus the number of fields added in one step cannot be chosen independently of the contraction obtained by averaging.

\paragraph{Closing the induction.}
The number of scales between cutoffs must satisfy two competing requirements: it must be large enough for averaging to contract the error, yet small enough for the added fields to remain a controlled perturbation. In Section~\ref{ss.homogenization.step}, we first prove moment bounds by induction on the cutoff scale, with the parameter choices in~\eqref{e.params.first.induction.step}. These bounds supply the low-moment estimates needed for concentration. We then repeat the induction on the cutoff scale using these estimates to improve the dependence on the moment order, with the revised parameter choices in~\eqref{e.params.second.induction.step}. The result is Proposition~\ref{p.coarse.grained.bound}: for every fixed~$p\in[1,\infty)$ and sufficiently small~$\gmcdelta$, both coarse-grained matrices of~$\a_m$ on~$\cu_m$ differ from~$\ahom_m\Id$ by a relative error of order~$\gmcdelta$ in~$L^p(\P)$, uniformly in~$m$.

At each step, concentration uses moment estimates already established on smaller cubes. Both inductions start from the first finitely many scales, where weak disorder gives the required smallness directly.

\paragraph{Bounds on the effective diffusivity.}
We next determine how~$\ahom_m$ decreases as~$m$ increases. When new fields are added to the coefficient, an affine function must be corrected to become harmonic. We approximate this correction by the solution of a Poisson equation. The energy of the approximation, computed in Lemma~\ref{l.laplacian.corrector.energy}, gives the leading decrease~$2\tau^2/d$ in~$\log\ahom_m$ per scale. Summing this calculation over blocks of scales gives~\eqref{e.main.asymptotic}.

For a general field law, the error in this calculation can be larger than the leading term. To prove power-law decay for all the fields allowed here, we use Proposition~\ref{p.homogenized.coefficient.strict.decay}: a nonconstant field strictly reduces the expected minimal energy, and this reduction can be repeated at separated scales. In dimension two, duality gives the exact formula~$\ahom_m=e^{-(m+1)\tau^2}$ (Proposition~\ref{p.special.two.d.exact.formula}).

\subsubsection{Large-scale regularity (Section~\ref{s.Holder.regularity})}

We use the coarse-graining estimates to prove H\"older estimates for~$\a_L$-harmonic functions, uniformly in~$L$. The moment bounds from renormalization do not give a small homogenization error on every cube. We therefore need an iteration that permits some of the scales to have a large error.

The quantity we iterate is the \textbf{excess}: the averaged~$L^2$ distance of the solution from its best affine approximation, divided by the side length of the cube. A scale is good when the coarse-graining error is small enough, and the coarse-grained ellipticity bounds hold, for Lemma~\ref{l.harmonic.approximation.good.scales.GMC} to compare the solution with a Laplace-harmonic function having the same boundary values on a smaller cube. The excess of this comparison function decreases on still smaller cubes. Transferring this decay to the original solution gives~\eqref{e.excess.decay.one.step}, with an additional error proportional to the slope of its affine approximation.

We must control this slope as well as the excess. Equation~\eqref{e.iteration.gradient.stability.GMC} bounds the change in slope between successive cubes by their excesses. Lemma~\ref{l.iteration.lemma.GMC} combines the two estimates: the loss in the resulting oscillation bound grows exponentially with the number of bad scales and with the sum of the approximation errors on the good scales. The proof thus requires bounds on both quantities, not just on the probability that an individual scale is bad.

We obtain these bounds from the independence of the fields~$\g_j$. The good-scale events themselves are not independent, since each involves several field indices. In the sums defining these events, we separate the contributions from nearby field indices from the tails. Independence controls the retained contributions, and concentration controls the tails. Proposition~\ref{p.density.of.good.scales} bounds the number of bad scales, and Lemma~\ref{l.sum.the.errors} bounds the sum of the approximation errors.

For a target H\"older exponent below one, Lemmas~\ref{l.sum.the.errors} and~\ref{l.min.scale.good.scale} show that, after a random number of scales, both losses are small enough to be absorbed in the desired oscillation estimate. This number is the scale separation in Theorem~\ref{t.C}, and has the exponential tail~\eqref{e.intro.regularity.minimal.scale}. The cutoff estimates of Proposition~\ref{p.cutoff.Holder.regularity} also cover scales larger than~$3^L$ and are uniform in~$L$, so they pass to~$\a$. The same iteration gives the energy bound~\eqref{e.energy.density.estimate}, used in Section~\ref{sec:mfd-uniqueness} to control the energy of small cubes and thin strips along their faces. Theorem~\ref{t.B} follows directly from coarse-graining on the whole cube.

\subsubsection{Uniqueness of the limit form (Section~\ref{sec:mfd-uniqueness})}

To prove uniqueness of the scaling limit, we must show that every subsequence gives the same Dirichlet form with the same normalization. The reversible measures already converge by the multiplicative chaos construction of~\citet{Kahane}. We first prove compactness of the minimal energies, then compare the forms obtained along any two convergent subsequences.

\paragraph{Compactness by resampling.}
For a field whose scale is small, partition a fixed cube into boxes separated by thin strips. The strips separate the box interiors by more than the field's range of dependence, so its restrictions to those interiors are independent. The regularity estimates bound the energy in the strips and in each box (Lemma~\ref{mfd:lem-strips}), and hence the cost of removing the field and of resampling each restriction. Lemma~\ref{mfd:lem-15} applies the Efron--Stein inequality to these independent restrictions. Summing over the omitted fields gives~\eqref{eq:mfd-16}: the minimal energies are approximated by their conditional expectations given finitely many fields, uniformly in the cutoff. This proves compactness in~$L^p(\P)$ and gives subsequential limits of the Green operators. The associated forms are local, consistent under restriction to subcubes, and satisfy the rule that multiplying the coefficient by a positive continuous function multiplies the energy measure by that function.

\paragraph{Comparison of two limits.}
We adapt the uniqueness arguments for Liouville quantum gravity metrics of~\citet{GM,DGUniqueness} and the comparison argument of~\citet{DGComparison}. Two subsequential limit forms~$\mathscr E$ and~$\mathscr F$ have the same domain and bound each other by deterministic multiples (Theorem~\ref{mfd:thm-C0}). We optimize the lower and upper multiples over all cubes and all functions in this domain; these are the \textbf{optimal comparison constants}. Resampling a field multiplies both energy measures by the same weight and preserves the comparison constants. Independence of the fields then makes these constants deterministic (Lemma~\ref{mfd:lem-endpoints}).

Suppose the two constants were different. On a positive proportion of scales, concentration would give a comparison of the minimal energies with affine boundary values that improves one of the optimal constants by a fixed fraction of their gap. On good cubes, harmonic approximation makes the solution's boundary values close enough to affine for this improvement to apply. We replace the solution inside disjoint such cubes by the harmonic extensions for the other form, keeping its boundary values fixed. If these cubes carry a fixed fraction of the energy, the sum of the local improvements gives a better comparison for the forms themselves, contradicting optimality (Theorem~\ref{mfd:thm-prop}). This adapts the path replacement argument of~\citet*[Theorem~1.9]{GM}: energy minimizers replace geodesics, and harmonic replacement replaces the modification of a path segment. The optimal constants therefore agree, so one limit form is a deterministic positive multiple of the other.

\paragraph{Errors relative to the gap.}
The gap between the constants may be arbitrarily small. Estimating the two minimal energies separately would give errors larger than this gap. Proposition~\ref{mfd:prop-conc} instead applies Efron--Stein to their difference; the resulting error is proportional to the gap. Affine approximation then transfers the improved comparison to the boundary values of the solution. We also need the replacement cubes to carry a fixed fraction of its energy. This plays the role of geodesic confluence in~\citet*[Section~4]{GM} and the field modifications in annuli in~\citet*[Section~4.2]{DGUniqueness}. The energy argument applies in every dimension~$d\geq2$.

The solution and its energy measure depend on the fields, so a probability bound for a fixed cube would not prevent most of the energy from lying on bad cubes. Proposition~\ref{mfd:prop-allchain} gives one event on which every chain of nested cubes has few bad cells. Integrating this count against the solution's energy measure, and using the energy bound from regularity, gives good cubes carrying enough energy for harmonic replacement (Lemma~\ref{mfd:lem-mass}).

\paragraph{The proportionality constant.}
Proportionality still allows different speeds for the limiting process: with the reversible measure~$\mu$ fixed, multiplying the Dirichlet form by a constant multiplies the generator by that constant. To obtain a unique limit with our prescribed rescaling of time, we must show that this factor is one. The coefficient~$\ahom_N^{-1}\a_N(3^N\,\cdot)$ is normalized to have homogenized matrix equal to the identity. Its coarse-grained matrices converge to the identity in expectation on large cubes, uniformly in~$N$ (Theorem~\ref{mfd:thm-eta}). Both limit forms inherit this normalization, which forces the factor to be one (Theorem~\ref{mfd:thm-c1}).

\subsubsection{Construction of the scaling limit (Sections~\ref{sec:mfd-process-convergence}--\ref{lim:sec-process-measure})}

We now combine the unique limit form with the limiting reversible measure~$\mu$ to construct the diffusion. We first identify the resolvents of the processes killed on leaving a cube. Their equations involve integration against the approximating reversible measures, whose densities have no uniform bound. The mass bounds of Proposition~\ref{mfd:prop-chaos-growth} make this integration continuous in a fractional Sobolev norm, and the energy controls that norm. Mollifying the measure-valued source then allows us to use the H\"older estimate for bounded sources. This gives uniform convergence of the killed resolvents, including their continuous zero-boundary representatives (Proposition~\ref{mfd:prop-uniform-resolvent}).

\paragraph{Convergence in probability.}
To pass from resolvents to path laws, we must prevent the paths from leaving a large cube or oscillating too much in a short time. A resolvent equation expresses the Laplace transform of the exit time. The local bound for its solution, proved by Moser iteration with respect to the reversible measure, bounds the probability of an early exit. Applying this bound at stopping times controls the oscillation of the paths and gives tightness (Section~\ref{s.tightness}). We identify every limiting path law from the killed resolvents: the Markov property expresses discounted integrals of finite-dimensional observations through compositions of resolvents, and these integrals determine a continuous path law. Since the limit form is unique, the quenched path laws converge in probability (Proposition~\ref{mfd:prop-quenched-convergence}).

\paragraph{Almost-sure convergence.}
Convergence in probability does not control the exceptional cutoffs in a single realization of the fields. To obtain almost-sure convergence, we compare minimal energies at two finite cutoffs and prove summable bounds on the probability of a failed comparison. On a small cube, the two energies have different deterministic normalizations. Dividing by these normalizations and rescaling to the unit cube gives two energies that converge in probability to the same limit (Lemma~\ref{mfd:lem-local-normalizations}). We must show that the normalizations themselves are close on enough scales to combine these local comparisons.

Suppose their ratios were less than a fixed number below one on a positive proportion of scales. Joining harmonic extensions on cubes of those scales would bound the minimal energy with a fixed nonconstant affine boundary value at a finite cutoff by a fixed factor less than one times its positive limit. This contradicts convergence in probability. Interchanging the two cutoffs rules out ratios above a fixed number greater than one. This adapts the comparison argument of~\citet{DGComparison}, in the form strengthened by~\citet*[Lemma~3.11]{Devlin2024}.

We then combine the local energy comparisons on scales where the normalization ratios are close to one. Lemma~\ref{mfd:lem-finite-good-levels} controls the failed comparisons along every chain of nested cubes, and the stopping argument of Lemma~\ref{mfd:lem-finite-stopping} joins harmonic extensions while keeping the total error small. The resulting error probabilities are summable, so the Borel--Cantelli lemma gives almost-sure convergence of the minimal energies (Proposition~\ref{mfd:prop-as-dirichlet}). The same comparison for Dirichlet and Neumann sources gives convergence of the forms and killed resolvents, and hence of the quenched path laws (Propositions~\ref{mfd:prop-as-forms} and~\ref{mfd:prop-as-quenched}).

\paragraph{Properties of the limit.}
We prove anomalous scaling by estimating exit times. Lemma~\ref{lim:lem-strict-decay} bounds the ratio of effective diffusivities at two arbitrary scales. Lemma~\ref{lim:lem-mean-exit} uses this ratio to bound the annealed mean exit time of the rescaled processes from a cube of side~$3^{-k}$ by~$C3^{-(2+\eta)k}$, uniformly for cutoffs at least~$k$, and passes this bound to the limit. Markov's inequality and the Borel--Cantelli lemma show that the limiting path leaves small cubes too quickly to have pointwise H\"older exponent greater than~$1/(2+\eta)$ at time zero (Corollary~\ref{lim:cor-paths}). This also proves singularity with respect to Brownian path laws.

To prove singularity of~$\mu$, we remove the positive smooth factor contributed by fields at scales larger than one. For the remaining stationary measure, we use the multiplicative chaos change of measure of~\citet[equation~(4.3)]{RVsurvey}. Weighting the joint law of the fields and a point in a bounded cube by this measure exponentially tilts the independent field values at that point. The strong law then gives opposite behavior of the approximating densities: they diverge at almost every point sampled from the limiting measure, and vanish at Lebesgue-almost every point. Restoring the positive factor preserves this distinction (Theorem~\ref{lim:thm-measure}(ii)).

Finally, a Nash inequality bounds the killed transition laws by a multiple of~$\mu$. Passing this bound to the limit and exhausting~$\Rd$ by cubes proves that every transition law at positive time is absolutely continuous with respect to~$\mu$. Since~$\mu$ is carried by a Lebesgue-null set and gives no mass to affine hyperplanes, the transition laws are singular with respect to every Gaussian law (Theorem~\ref{lim:thm-nongaussian}).

\endgroup

\subsection{Comparison with previous work}
\label{ss.lit.review}

The planar lattice analogue of our model is a random walk whose edge conductances are exponentials of the discrete Gaussian free field. In these walks, the probability of choosing a neighbor is proportional to the conductance of the joining edge, and the reversible measure at a vertex is the sum of the incident conductances. Thus the field affects both the transition probabilities and the reversible measure, as it does in~\eqref{e.sde.intro}. For this model,~\citet*{BDG} proved recurrence, return-probability estimates and exponents for expected exit times. At small disorder,~\citet*{DW} proved tightness of the random quenched path laws for a version using a smooth approximation to the field.

\smallskip

As mentioned above,~\citet*[Problem~16.13]{BiskupSurvey} asked whether the field and the diffusion converge jointly, after a suitable rescaling of time, with continuous, nonconstant limiting paths. We answer this question for the continuum model~\eqref{e.sde.intro}: under~\ref{a.g1}--\ref{a.g4}, at weak disorder and in every dimension~$d\geq2$, we prove joint convergence of the field and the diffusion to a unique scaling limit with continuous, nonconstant paths (Theorem~\ref{t.scaling.limit}). We expect the methods developed here also to yield scaling limits for lattice models.

\smallskip

When~$\g_0$ is Gaussian, the limiting reversible measure~$\mu$ is Gaussian multiplicative chaos. This construction goes back to~\citet{Kahane}; see~\citet*{RVsurvey} and~\citet{Berestycki1} for accounts of the theory, and~\citet*{DRSV} for log-correlated fields.

Gaussian multiplicative chaos also appears as the reversible measure of Liouville Brownian motion, constructed by~\citet*{GRV} and~\citet{Berestycki2}. That process is obtained by changing the time parametrization of Brownian motion using the chaos measure of the planar Gaussian free field. Its Dirichlet energy remains~$\int|\nabla u|^2\,dx$; for a smooth positive density~$\a$, the corresponding generator is~$\a^{-1}\Delta$. In our model, the field also weights the energy:~$X$ has Dirichlet form~$\int\a|\nabla u|^2\,dx$ on~$L^2(\a\,dx)$ and generator~$\a^{-1}\nabla\cdot(\a\nabla)$. Convergence of the reversible measures therefore leaves a further problem: proving that the rescaled energies converge to a unique Dirichlet form, which determines the limiting dynamics together with~$\mu$.

\smallskip

For a critically correlated incompressible random drift,~\citet*{ABK.SD} proved a quenched invariance principle with logarithmic superdiffusive scaling. Our renormalization argument is more closely related to~\citet*{ASD}, who proved power-law superdiffusion in a self-similar incompressible random drift with sufficiently small positive Hurst exponent. As in that paper, we follow an induction across scales for the elliptic operator, using variationally defined coarse-grained matrices in place of pointwise ellipticity bounds. In the incompressible models the effective diffusivity increases with scale; here it decreases.

The principal new arguments are those needed to construct and identify the non-Gaussian scaling limit. The renormalization estimates give an error controlled by the disorder strength, rather than an error that vanishes as the scale increases. We use the regularity estimates to bound the effect of resampling fine-scale fields on minimal energies and obtain compactness of the Dirichlet forms. To prove uniqueness, we adapt the comparison and replacement arguments for Liouville quantum gravity metrics of~\citet{GM,DGUniqueness,DGComparison}. Harmonic replacement on cubes carrying a fixed fraction of the energy makes two subsequential limit forms proportional, and the effective-diffusivity normalization makes them equal. Comparing minimal energies at two finite cutoffs then gives almost-sure convergence, adapting the argument of~\citet{Devlin2024}. Passing from the forms to the path laws also requires bounds on integration against the singular limiting reversible measure.

\begingroup\runinsubsubsections

\subsection{Notation}
\label{s.notation}

\subsubsection*{Basic notation}
The Euclidean norm on~$\R^m$ is denoted by~$|x|_2=(\sum_{i=1}^m x_i^2)^{\nf12}$. In energies, variational formulas and function-space norms, we also write~$|x|$ for this norm, and~$p\cdot q=\sum_{i=1}^d p_iq_i$ denotes the Euclidean pairing. Unit vectors in these formulas are Euclidean unit vectors. We write~$r\wedge s\coloneqq\min\{r,s\}$,~$r\vee s\coloneqq\max\{r,s\}$ and~$r_+\coloneqq r\vee0$. The H\"older conjugate exponent of~$p\in[1,\infty]$ is denoted by~$p'$, with~$p^{-1}+(p')^{-1}=1$; in particular,~$1'=\infty$ and~$\infty'=1$. Indicator functions for events and subsets of~$\Rd$ are denoted by~$\indc$.

\subsubsection*{Linear algebra}
The set of positive definite symmetric~$d$-by-$d$ matrices is denoted by~$\Rsymp$, and~$\Id$ is the identity matrix. The Loewner ordering of symmetric matrices is denoted by~$\leq$:~$A\leq B$ means that~$B-A$ is positive semidefinite. The norm~$|A|$ is the Euclidean operator norm,~$|A|=\sup_{p\neq0}|Ap|_2/|p|_2$. We identify a positive scalar coefficient with the corresponding scalar multiple of~$\Id$.

\subsubsection*{Measure and integration}
The Lebesgue measure of a measurable subset~$U\subseteq\Rd$ is denoted by~$|U|$. Volume-normalized integrals and~$L^p$ norms are denoted by~$(f)_U$ and~$\|f\|_{\underline L^p(U)}$, as defined before Theorem~\ref{t.C}. Similarly, for a nonempty finite set~$A$, we write~$\avsum_{z\in A}c_z\coloneqq\frac1{\#A}\sum_{z\in A}c_z$. If a cube~$U$ is the disjoint union of cubes~$(z+r\cu_0)_{z\in A}$ of equal side~$r$, then
\begin{equation}
	\avsum_{z\in A}
	\|f\|_{\underline L^p(z+r\cu_0)}^p
	=
	\|f\|_{\underline L^p(U)}^p\,.
	\label{e.average.normalized.norm.partition}
\end{equation}

\subsubsection*{Function spaces}
We use the standard H\"older spaces~$C^{k,\alpha}(U)$ and Sobolev spaces~$W^{s,p}(U)$. When~$p=2$, we write~$H^s(U)=W^{s,2}(U)$. The space~$H^1_0(U)$ is the closure of~$C_c^\infty(U)$ in~$H^1(U)$. In~\ref{a.g2}, the gradient magnitude is~$\sum_{i=1}^d|\partial_i\g_0|$, and~$\|\nabla^2\g_0\|_{L^\infty(\cu_0)}$ is the Lipschitz constant of~$\nabla\g_0$ on~$\cu_0$, using this gradient norm and the spatial distance~$|x-y|_\infty=\max_{1\leq i\leq d}|x_i-y_i|$. For~$s\in (0,1)$ and~$p\in [1,\infty)$, the volume-normalized fractional Sobolev seminorm on~$\cu_m$ is
\begin{equation*}
	[  u ]_{\underline{W}^{s,p}(\cu_m)}
	\coloneqq
	\biggl(
	s
	\fint_{\cu_m} \int_{\cu_m}
	\frac{|u(x) - u(y)|^p}{|x-y|^{d+sp}}
	\, dx \, dy
	\biggr)^{\! \nf1p}
	\,.
\end{equation*}
In this fractional Sobolev seminorm, distance and vector magnitude are Euclidean. At~$s=1$, the seminorm is the classical first-order one,~$\|\nabla u\|_{\underline L^p(\cu_m)}$, with Euclidean gradient magnitude. For a bounded Lipschitz set~$U$, the same formulas define~$[\,\cdot\,]_{\underline W^{s,p}(U)}$. Let
\begin{equation*}
	\|u\|_{\underline W^{s,p}(U)}
	\coloneqq
	[u]_{\underline W^{s,p}(U)}+
	|U|^{-\nf{s}{d}}\|u\|_{\underline L^p(U)}\,.
\end{equation*}
For~$p=\infty$ we write~$[u]_{\underline W^{s,\infty}(U)}\coloneqq\sup_{x\neq y\in U}|u(x)-u(y)|\,|x-y|^{-s}$ for the H\"older seminorm of order~$s$, and define the norm by the same formula with~$\|u\|_{L^\infty(U)}$. For~$s\in(0,1]$,~$q\in(1,\infty)$,~$N\in\N$, and~$F:U\to\R^N$, define the inhomogeneous volume-normalized negative Sobolev norm by
\begin{equation*}
	\|F\|_{\Wminusul{-s}{q}(U)}
	\coloneqq
	\sup_{\varphi\in C^\infty(\overline U;\R^N),\,\varphi\neq0}
	\frac{\left|\fint_U F\cdot\varphi\right|}
	{[\varphi]_{\underline W^{s,q'}(U)}+
	|U|^{-\nf{s}{d}}\|\varphi\|_{\underline L^{q'}(U)}}\,.
\end{equation*}
When~$q=2$ and~$N=d$, we write
\begin{equation}
	\|F\|_{\Hminusuls{-s}(U)}
	\coloneqq
	\sup_{\varphi\in C^\infty(\overline U;\Rd),\,\varphi\neq0}
	\frac{\left|\fint_U F\cdot\varphi\right|}
	{[\varphi]_{\underline H^s(U)}+
	|U|^{-\nf{s}{d}}\|\varphi\|_{\underline L^2(U)}}\,.
	\label{e.hatted.negative.norm.convention}
\end{equation}
Here and below,~$\underline H^s=\underline W^{s,2}$.  If~$F\in(H^s(U;\Rd))^*$, then the numerator in~\eqref{e.hatted.negative.norm.convention} is interpreted as~$|U|^{-1}|\langle F,\varphi\rangle|$; by density of smooth functions, the supremum is then the corresponding dual norm. For the unadorned negative norm in Theorem~\ref{t.B} we use zero-boundary tests:
\begin{equation}
	\|F\|_{\underline H^{-1}(U)}
	\coloneqq
	\sum_{i=1}^d
	\sup_{\substack{\varphi\in C_c^\infty(U)\\ \|\nabla\varphi\|_{\underline L^2(U;\ell^2)}\leq1}}
	\left|\fint_U F_i\varphi\right|\,.
	\label{e.ordinary.negative.norm.convention}
\end{equation}
Each summand is the volume-normalized norm of~$F_i$ in the dual of~$H^1_0(U)$; the hatted norm uses unrestricted smooth tests. On~$\cu_0$, the volume-normalized and unnormalized conventions coincide.

\smallskip

For a bounded convex domain~$U\subseteq\Rd$, the spaces of solenoidal vector fields are
\begin{equation}
	L^2_{\sol}(U)\coloneqq\bigl(\nabla H^1_0(U)\bigr)^\perp
	\qquad\text{and}\qquad
	\Lsolo(U)\coloneqq\bigl(\nabla H^1(U)\bigr)^\perp\,.
	\label{e.sols.and.pots}
\end{equation}
All orthogonal complements in~\eqref{e.sols.and.pots} are taken in~$L^2(U;\Rd)$.

\subsubsection*{Besov-type seminorms}
The volume-normalized Besov-type seminorm on~$\cu_m$ for exponents~$s\in(0,1)$ and~$p,q,r\in[1,\infty)$ is
\begin{equation}
	[  u ]_{\underline{B}^s_{p,q,r}(\cu_m)}
	\coloneqq
	\biggl(
	s
	\sum_{n=-\infty}^m
	3^{-nsq}
	\biggl(
	\avsum_{z\in 3^{n-1} \Zd \cap \cu_m, \, z+\cu_n \subseteq \cu_m}
	\bigl\| u - (u)_{z+\cu_n} \bigr\|_{\underline{L}^p(z+\cu_n)}^r
	\biggr)^{\!\nf qr} \biggr)^{\!\nf1q} \,.
	\label{e.Besov.Morrey.seminorm}
\end{equation}
For~$r=p$, the seminorm reads
\begin{equation*}
	[  u ]_{\underline{B}^s_{p,q,p}(\cu_m)}
	=
	\biggl(
	s \sum_{n=-\infty}^m
	3^{-nsq}
	\biggl(
	\avsum_{z\in 3^{n-1} \Zd \cap \cu_m, \, z+\cu_n \subseteq \cu_m}
	\bigl\| u - (u)_{z+\cu_n} \bigr\|_{\underline{L}^p(z+\cu_n)}^p
	\biggr)^{\!\nf qp} \biggr)^{\!\nf1q} \,.
\end{equation*}
When the third index is omitted, it equals the first: we write~$\underline B^s_{p,q}\coloneqq\underline B^s_{p,q,p}$. The corresponding norm is~$\|f\|_{\underline B^s_{p,q,r}(\cu_m)} \coloneqq[f]_{\underline B^s_{p,q,r}(\cu_m)} +3^{-sm}|(f)_{\cu_m}|$. By Lemma~\ref{l.Wsp.vs.Bspp}, there exists~$C(d)<\infty$ such that, for every~$s\in(0,1)$ and~$p\in[1,\infty)$,
\begin{equation}
	C^{-1} [  u ]_{\underline{W}^{s,p}(\cu_m)}
	\leq
	[ u ]_{\underline{B}^s_{p,p,p}(\cu_m)}
	\leq
	C [  u ]_{\underline{W}^{s,p}(\cu_m)}
	\,.
	\label{e.Wsp.vs.Bspp.intro}
\end{equation}
The corresponding norms are equivalent for each fixed~$s,p$; unlike the seminorm comparison, their constants need not be uniform in these parameters.

\smallskip

We allow each of~$p,q,r$ to equal~$\infty$.  If~$p=\infty$, the volume-normalized~$L^p$ norm is replaced by the~$L^\infty$ norm; if~$r=\infty$, the average over cube centers is replaced by a supremum; and if~$q=\infty$, the sum over scales is replaced by a supremum.  We also allow~$s=1$ when~$q=\infty$.

\smallskip

We use Euclidean distance and vector magnitude in the H\"older seminorms. When~$q=r=\infty$, the Besov seminorm is equivalent to the H\"older seminorm of order~$s$ for every~$p$. More precisely,
\begin{align}
	C^{-1}s \sup_{x,y\in \cu_m, \, x\neq y}
	\frac{|u(x) - u(y)| }{|x-y|^s }
	&\leq
	[  u ]_{\underline{B}^s_{1,\infty,\infty}(\cu_m)}
	\leq
	[  u ]_{\underline{B}^s_{\infty,\infty,\infty}(\cu_m)}
	\notag \\ &
	\leq
	C \sup_{x,y\in \cu_m, \, x\neq y}
	\frac{|u(x) - u(y)| }{|x-y|^s }\,.
	\label{e.Besov.to.Holder}
\end{align}
The first inequality of~\eqref{e.Besov.to.Holder} is the only nontrivial one, and it is the Campanato characterization of H\"older spaces.

\smallskip

For~$s\in(0,1]$ and~$q,r\in[1,\infty]$, define
\begin{equation}
	\| f \|_{\Besov{-s}{r}{q}(\cu_m)}
	\coloneqq
	\biggl(
	s\sum_{n=-\infty}^m
	3^{nsq}
	\biggl(
	\avsum_{z\in 3^{n} \Zd \cap \cu_m}
	\bigl| (f)_{z+\cu_n} \bigr|^r
	\biggr)^{\!\nf qr} \biggr)^{\!\nf1q} \,.
	\label{e.Besov.weaknorm.def}
\end{equation}
If~$q^{-1}+(q')^{-1}=r^{-1}+(r')^{-1}=1$ and either~$s<1$ or~$(s,q)=(1,1)$, then Lemma~\ref{l.circ.norm.dominates} gives, for every~$f,g$ for which the right side is finite and~$fg\in L^1(\cu_m)$,
\begin{equation}
	\biggl| \fint_{\cu_m} f g \biggr|
	\leq
	C s^{-1}  \| g \|_{\underline{B}^{s}_{1,q',r'}(\cu_m)}
	\| f \|_{\Besov{-s}{r}{q}(\cu_m)}\,.
	\label{e.duality.splitting}
\end{equation}

\subsubsection*{Triadic cubes and simplices}
Spatial balls and distances use~$|x-y|_\infty$:~$B_r(x)=x+2r\cu_0$, and~$\dist(U,V)=\inf_{x\in U,y\in V}|x-y|_\infty$. The axis-aligned triadic cube~$\cu_n$, defined in Section~\ref{ss.assumptions}, has side length~$3^n$ and is centered at the origin. More generally,~$z+\cu_n$ denotes its translate centered at~$z\in\Rd$. For a permutation~$\pi$ of~$\{1,\ldots,d\}$,~$n\in\Z$, and~$z\in\Rd$, the simplex of size~$3^n$ is
\begin{equation*}
	\spx_n^\pi(z)  \coloneqq  z + 3^n  \Bigl\{ (x_1,\ldots,x_d) \in\Rd \,:\, -\nf12 < x_{\pi(1)} < x_{\pi(2)} < \cdots < x_{\pi(d)} < \nf12  \Bigr\}
	\,.
\end{equation*}
These~$d!$ \textbf{Kuhn simplices} partition~$z+\cu_n$ up to null sets and each has volume~$3^{nd}/d!$. They form the \textbf{Kuhn triangulation} of the cube. For~$n\leq m$, the triangulations of the scale-$n$ triadic subcubes of~$z+\cu_m$ refine its Kuhn triangulation. Each fine simplex belongs to a coarse simplex, even when its supporting cube crosses a diagonal face of that coarse simplex.

\endgroup

\subsection{Formalization in the Lean~4 proof assistant}
\label{ss.formalization}

Theorems~\ref{t.A}--\ref{t.C}, together with their supporting results, have been formalized and machine-checked in the Lean~4 proof assistant for the diffusions constructed from elliptic resolvents. The identification of these diffusions with the solutions of~\eqref{e.sde.intro} is proved in the paper and is not part of the formalization. The formal development and its paper-to-Lean correspondence are maintained in the private repository
\begin{center}
\href{https://github.com/nitromannitol/subdiffusive-process}{\texttt{github.com/nitromannitol/subdiffusive-process}}\,.
\end{center}

\subsection*{Use of artificial intelligence}

The authors developed all the main ideas of this paper in 2024. Versions of GPT were used to check the arguments, most recently GPT-6.1 Sol, and versions of Claude were used to help with the writing, most recently Claude Opus~5.5 and Claude Sonnet~5.5. The authors take responsibility for the mathematical content and presentation of the paper.

\subsection*{Acknowledgments}

S.A. was supported by NSF grant DMS-2350340. A.B. was supported by NSF grant DMS-2202940. T.K. was supported by the Academy of Finland and by the European Research Council (ERC) under the European Union's Horizon 2020 research and innovation programme (grant agreement No 818437).

\section{Coarse-graining theory}\label{s.coarse.graining}

In this section we develop the deterministic estimates used to compare a heterogeneous elliptic equation with a constant-coefficient equation. Throughout,~$\a:\Rd\to\Rsymp$ is a fixed uniformly elliptic coefficient field; no randomness is used. The estimates depend on coarse-grained matrices, allowing us to use the multiscale bounds proved later in place of cutoff-dependent pointwise ellipticity bounds.

\smallskip

The definitions are from~\citet{AK.HC}. Section~\ref{ss.cg.mats} recalls the coarse-grained matrices and their variational properties. Section~\ref{ss.cg.ellipticity} uses them to measure ellipticity across scales and to obtain Poincar\'e inequalities. The sensitivity estimates of Section~\ref{ss.sensitivity} control the effect of multiplying the coefficient by a slowly varying field; they are used when adding fields in the renormalization argument. Finally, Section~\ref{ss.general.coarse.graining.estimates} gives Caccioppoli inequalities and combines the multiscale bounds into a comparison of Dirichlet problems, including equations with a nonzero right-hand side. We indicate where an estimate is adapted to the present setting.

\begingroup\runinsubsubsections

\subsection{Coarse-grained matrices}
\label{ss.cg.mats}

The two coarse-grained matrices record the energies of affine Dirichlet data and constant Neumann fluxes. We recall their definitions from~\citet{AK.HC}, using one quadratic quantity~$J$ to track both matrices. Its subadditivity is the starting point for averaging across scales.

\smallskip

We denote the linear space of solutions~$u$ of~$-\nabla \cdot \a \nabla u = 0$ in an open set~$U \subseteq \Rd$ by
\begin{equation*}
	\mathcal{A}(U; \a)  \coloneqq  \{ u \in H^1_{\mathrm{loc}}(U) : -\nabla \cdot \a \nabla u = 0 \quad \text{in~$U$}\}
	\,.
\end{equation*}
For every~$p, q \in \Rd$ and bounded convex domain~$U\subseteq\Rd$, we define the quantity
\begin{equation}
	J(U,p,q; \a)  \coloneqq  \sup_{v \in \mathcal{A}(U;\a)\cap H^1(U)} \fint_{U} \Bigl( -\frac{1}{2} \nabla v \cdot \a \nabla v - p \cdot \a \nabla v + q \cdot \nabla v \Bigr)
	\,.
	\label{e.J.general}
\end{equation}
The supremum is attained in~$H^1(U)$, and the maximizer is unique up to an additive constant. We denote its mean-zero representative by~$v(\cdot,U,p,q;\a)$, so that~$(v(\cdot,U,p,q;\a))_U=0$. For a bounded convex domain~$U \subseteq \Rd$, define the~\textbf{coarse-grained matrices}~$\a(U)$ and~$\a_*(U)$ by the identity, valid for every~$p,q\in\Rd$,
\begin{equation}
	J(U, p,q; \a) \eqqcolon \frac{1}{2} p \cdot \a(U) p + \frac{1}{2} q \cdot \a_{*}^{-1}(U) q - p \cdot q \,.
	\label{e.coarse.grained.matrices}
\end{equation}
For the coarse-grained matrices, superscripts denote matrix powers; for example,~$\a_*^{-1}(U)=(\a_*(U))^{-1}$. When the coefficient needs to be specified, it is written after a semicolon, as in~$\a(U;\a^{-1})$. The variational characterizations below express the matrices as minimal energies. They let us compare two coefficient fields by evaluating their energies on the same competitors, as used in Section~\ref{ss.sensitivity}. For~$p, q \in \Rd$,
\begin{equation}
	\frac12 p \cdot \a(U) p  = \min_{u \in p\cdot x + H_0^1(U)} \fint_{U}\frac12 \nabla u \cdot \a \nabla u
	=  \max_{\g \in \Lsol(U)} \fint_U \Bigl(-\frac12 \g \cdot \a^{-1} \g + p \cdot \g  \Bigr) \,,
	\label{e.variational.a}
\end{equation}
and
\begin{equation}
	\frac12 q \cdot \a^{-1}_*(U) q
	=
	\max_{w \in H^1(U)} \fint_{U} \Bigl(-\frac12 \nabla w \cdot \a \nabla w + q \cdot \nabla w \Bigr)
	= \min_{\g \in q + \Lsolo(U)} \fint_{U} \frac12 \g \cdot \a^{-1} \g \,.
	\label{e.variational.astar}
\end{equation}

\subsubsection{Properties of the coarse-grained matrices}

We recall the variational properties from~\citet{AK.HC}. Subadditivity bounds the energy on a domain by the energies on its parts, while the identities for the maximizer relate~$J$ to averaged gradients and fluxes. Throughout,~$U\subseteq\Rd$ is a bounded convex domain.

\smallskip

The coarse-grained matrices are bounded by averages of the field, in the Loewner order:
\begin{equation}
	\Bigl( \fint_{U} \a^{-1} \Bigr)^{\!-1} \leq \a_{*}(U) \leq \a(U)
	\leq
	\fint_{U} \a  \,.
	\label{e.ord.bounds.truncated}
\end{equation}

\smallskip

The coarse-grained matrices are~\textbf{subadditive}: for every partition~$\{ U_i \}_{i=1}^N$ of~$U$ (up to a Lebesgue null set)
\begin{equation}
	\a(U) \leq  \sum_{i=1}^N \frac{|U_i|}{|U|} \a(U_i)  \qand \a_*(U) \geq \biggl( \sum_{i=1}^N \frac{|U_i|}{|U|}\a^{-1}_*(U_i) \biggr)^{-1} \,,
	\label{e.subadditivity.a}
\end{equation}
which generalizes the first and third inequalities of~\eqref{e.ord.bounds.truncated}. Equivalently, for every~$p, q \in \Rd$,
\begin{equation}
	J(U, p, q; \a) \leq
	\sum_{i=1}^N \frac{|U_i|}{|U|} J(U_i, p , q; \a) \,.
	\label{e.subadditivity.J}
\end{equation}
The spatial averages of the gradient and the flux of the maximizer can be written explicitly in terms of the coarse-grained matrices: for~$p, q \in \Rd$,
\begin{equation}
	\left\{
	\begin{aligned}
		& \fint_U
		\nabla v(\cdot,U,p,q; \a)
		=
		- p + \a_{*}^{-1}(U) q
		\\ &
		\fint_U
		\a \nabla v(\cdot,U,p,q; \a)
		=
		q - \a(U) p\,.
	\end{aligned}
	\right.
	\label{e.v.spatial.averages}
\end{equation}
For every~$e\in\Rd$, the maximizer~$v(\cdot,U,0,\a_{*}(U)e;\a)$ satisfies, by~\eqref{e.v.spatial.averages},
\begin{equation*}
	\begin{pmatrix}
		\bigl(\nabla v(\cdot,U,0,\a_{*}(U)e; \a) \bigr)_{U}  \\
		\bigl( \a   \nabla v(\cdot,U,0,\a_{*}(U)e; \a) \bigr)_{U}
	\end{pmatrix}
	=
	\begin{pmatrix}
		e  \\  \a_{*}(U)  e
	\end{pmatrix}
	\,.
\end{equation*}
The identity~\eqref{e.J.plug.identity} below gives~$J(U,e,\a_{*}(U)e\,;\a)=\frac12 e\cdot(\a(U)-\a_{*}(U))e$. For every~$w\in\A(U;\a)\cap H^1(U)$ and~$p,q\in\Rd$, the first variation of~\eqref{e.J.general} gives
\begin{align}
	q\cdot \fint_U \nabla w - p \cdot \fint_U \a \nabla w
	=
	\fint_U \nabla w \cdot \a \nabla v(\cdot,U,p,q; \a)\,.
	\label{e.firstvar}
\end{align}
The second variation gives
\begin{align}
	& J(U,p,q; \a) - \fint_U \Bigl  ( -\frac12 \nabla w \cdot \a \nabla w -p\cdot \a \nabla w+ q\cdot \nabla w   \Bigr  )
	\notag \\ & \qquad \qquad\qquad
	=
	\fint_U \frac12 \bigl ( \nabla v(\cdot,U,p,q;\a) - \nabla w \bigr )\cdot \a \bigl ( \nabla v(\cdot,U,p,q;\a) - \nabla w \bigr )\,.
	\label{e.secondvar}
\end{align}
It follows that~$J$ can be written as the energy of the maximizer: for every~$p,q\in\Rd$,
\begin{equation}
	J(U,p,q;\a) = \fint_U \frac12 \nabla v(\cdot,U,p,q;\a) \cdot \a \nabla v(\cdot,U,p,q;\a)
	\,.
	\label{e.Jenergy.v}
\end{equation}
The Euler equation also gives
\begin{equation*}
	J(U,p,q;\a) = \frac12 \Bigl( q\cdot \fint_U \nabla v(\cdot,U,p,q;\a) - p \cdot \fint_U \a \nabla v(\cdot,U,p,q;\a) \Bigr) \,.
\end{equation*}
The coarse-grained matrix~$\a_{*}(U)$ bounds the energy of an arbitrary~$H^1(U)$ function from below in terms of the spatial average of its gradient: for every~$u\in H^1(U)$,
\begin{align}
	\frac12 \bigl(  \nabla u \bigr)_U \cdot \a_{*}(U) \bigl(  \nabla u \bigr)_U
	\leq
	\fint_U \frac12 \nabla u \cdot \a \nabla u
	\,.
	\label{e.energymaps}
\end{align}
The dual inequality is, for every~$\h\in L^2_{\sol}(U)$,
\begin{align*}
	\frac12\bigl( \h \bigr)_U \cdot \a^{-1} (U) \bigl( \h \bigr)_U
	\leq
	\fint_U \frac12 \h \cdot \a^{-1} \h
	\,.
\end{align*}
For every constant positive-definite matrix~$\b\in\Rsymp$ and every~$u\in\mathcal A(U;\a)\cap H^1(U)$,
\begin{equation}
	\bigl| \bigl(  ( \a - \b) \nabla u \bigr)_U
	\bigr|
	\leq
	2\biggl(\max_{|e|=1}J(U,e,\b e\,;\a)\biggr)^{\nf12}
	\biggl( \fint_U \frac12 \nabla u \cdot \a \nabla u  \biggr)^{\!\nf12}
	\,,
	\label{e.fluxmaps}
\end{equation}
in particular, with the coarse-graining error~$\mathcal{E}_{s,\infty,1}$ of Definition~\ref{d.mathcal.E} below, for every~$m \in \Z$ and every~$u\in\mathcal A(\cu_m;\a)\cap H^1(\cu_m)$,
\begin{equation*}
	3^{-s m} \| ( \a-\b) \nabla u \|_{\Besov{-s}{2}{1}(\cu_m)}
	\leq
	4s^{-1}
	|\b|^{\nf12}
	\|\a^{\nf12} \nabla u \|_{\underline{L}^2(\cu_m)}
	\mathcal{E}_ {s,\infty,1} (\cu_m;\a,\b)
	\,.
\end{equation*}
The bound~\eqref{e.fluxmaps} follows from~\eqref{e.firstvar},~\eqref{e.Jenergy.v}, and the Cauchy--Schwarz inequality.

\smallskip

The identity~\eqref{e.coarse.grained.matrices} can be rewritten for every symmetric matrix~$\b \in \R^{d \times d}_{\sym}$ as
\begin{align}
	J(U,e,\b e; \a)
	&=
	\frac12 e \cdot \bigl( \a(U) - \a_{*}(U) \bigr) e
	+
	\frac12 e \cdot \bigl(\b - \a_{*}(U)\bigr) \a_{*}^{-1}(U) \bigl(\b - \a_{*}(U)\bigr) e
	\,.
	\label{e.J.plug.identity}
\end{align}
If~$\b\in\Rsymp$, then completing the square in the first argument of~$J$ instead of the second gives
\begin{align}
	J(U,\b^{-1} e, e; \a)
	&
	=
	\frac12 e \cdot \bigl(\a_{*}^{-1}(U) - \a^{-1}(U)\bigr)e
	+
	\frac12 e \cdot \bigl(\b^{-1} - \a^{-1} (U)\bigr) \a(U) \bigl(\b^{-1} - \a^{-1} (U)\bigr) e
	\,.
	\label{e.J.plug.identity.inverses}
\end{align}
Applying~\eqref{e.J.plug.identity} and~\eqref{e.J.plug.identity.inverses} with~$\b^{-\nf12} e$ and~$\b^{\nf12} e$, respectively, in place of~$e$ yields
\begin{align}
	\label{e.J.half.half.identity}
	J(U,\b^{-\nf12}e,\b^{\nf12} e;\a)
	&
	=
	\frac12 e \cdot \b^{-\nf12}\bigl( \a(U) - \a_{*}(U) \bigr) \b^{-\nf12}e
	\notag \\ & \qquad
	+
	\frac12 e \cdot \bigl(\Id - \b^{-\nf12}\a_{*}(U)\b^{-\nf12}\bigr) \b^{\nf12}\a_{*}^{-1}(U) \b^{\nf12} \bigl(\Id - \b^{-\nf12}\a_{*}(U)\b^{-\nf12}\bigr) e
	\\ &
	\label{e.J.half.half.identity.numbertwo}
	=
	\frac12 e \cdot \b^{\nf12}\bigl( \a_{*}^{-1}(U)- \a^{-1}(U)  \bigr) \b^{\nf12}e
	\notag \\ & \qquad
	+
	\frac12 e \cdot \bigl(\Id - \b^{\nf12}\a^{-1}(U)\b^{\nf12}\bigr) \b^{-\nf12}\a(U) \b^{-\nf12} \bigl(\Id - \b^{\nf12}\a^{-1}(U)\b^{\nf12}\bigr) e
	\,.
\end{align}

\begin{lemma}
\label{l.sharp.compare.J}
For every bounded convex domain~$U\subseteq\Rd$ and~$\b\in\Rsymp$,
\begin{align}
	\max_{|e| = 1} J(U, \b^{-\nf12}   e, \b^{\nf12}  e; \a)
	&
	\leq
	\frac12 \bigl|\b^{-\nf12} \bigl( \a(U) - \a_{*}(U) \bigr) \b^{-\nf12}  \bigr|
	+
	\frac12 \bigl| \b^{-\nf12} \a^{\nf12}_*(U) -  \b^{\nf12} \a^{-\nf12}_*(U)\bigr|^2
	\notag \\ &
	\leq
	2 \max_{|e| = 1} J(U, \b^{-\nf12}   e, \b^{\nf12}  e; \a)  \,,
	\label{e.bound.for.J.general.preplugin}
\end{align}
and
\begin{align}
	\max_{|e| = 1} J(U, \b^{-\nf12}   e, \b^{\nf12}  e; \a)
	&
	\leq
	\frac12 \bigl|\b^{\nf12}   \bigl( \a^{-1}(U) - \a^{-1}_{*}(U) \bigr) \b^{\nf12}  \bigr|
	+
	\frac12 \bigl| \b^{\nf12} \a^{-\nf12}(U) -  \b^{-\nf12} \a^{\nf12}(U)\bigr|^2
	\notag \\ &
	\leq
	2 \max_{|e| = 1} J(U, \b^{-\nf12}   e, \b^{\nf12}  e; \a)  \,.
	\label{e.bound.for.J.general.preplugin2}
\end{align}
In particular, there exists a universal constant~$C<\infty$ such that
\begin{multline}
	\bigl|\b^{-\nf12} \bigl( \a(U) - \a_{*}(U) \bigr) \b^{-\nf12} \bigr|
	+
	\bigl|\b^{\nf12} \a^{-1}_*(U) \b^{\nf12} - \Id \bigr|^2
	+
	\bigl|\b^{-\nf12} \a(U) \b^{-\nf12} - \Id\bigr|^2
	\\
	\leq
	C\max_{|e| = 1} J(U, \b^{-\nf12}   e, \b^{\nf12}  e; \a) \bigl(1 + \max_{|e| = 1} J(U, \b^{-\nf12}   e, \b^{\nf12}  e; \a) \bigr)  \,.
	\label{e.coarse.grained.matrix.deviation.by.J}
\end{multline}

\end{lemma}
\begin{proof}
For positive-semidefinite matrices~$M$ and~$N$, we have~$\max\{|M|,|N|\}\leq|M+N|\leq|M|+|N|\leq2|M+N|$. Applying this inequality to the two summands in~\eqref{e.J.half.half.identity} proves~\eqref{e.bound.for.J.general.preplugin}; applying it to the two summands in~\eqref{e.J.half.half.identity.numbertwo} proves~\eqref{e.bound.for.J.general.preplugin2}.

\smallskip

For every positive-definite matrix~$A$ and~$J\geq0$,
\begin{equation}
	\frac 12 \bigl| (\Id-A)A^{-1} (\Id-A) \bigr|
	\leq J
	\quad \implies \quad
	\bigl|\Id-A\bigr|^2 \leq 4J(1+J)\,.
	\label{e.from.a.ratios.to.J}
\end{equation}
Indeed, the hypothesis implies
\begin{equation*}
	\bigl| \Id - A \bigr|^2
	\leq
	2J|A|
	\leq
	2J|\Id - A | + 2J
	\leq
	\frac12 |\Id - A |^2 + 2J^2 + 2J\,.
\end{equation*}
Since~$(\Id-A)A^{-1}(\Id-A)=(\Id-A^{-1})A(\Id-A^{-1})$, apply~\eqref{e.from.a.ratios.to.J} with~$A=(\b^{-\nf12}\a_*(U)\b^{-\nf12})^{-1}$ in~\eqref{e.J.half.half.identity} and with~$A=(\b^{\nf12}\a^{-1}(U)\b^{\nf12})^{-1}$ in~\eqref{e.J.half.half.identity.numbertwo}. This gives
\begin{align}
	\lefteqn{
	\max\Bigl\{
	\bigl|\b^{-\nf12}\a(U)\b^{-\nf12} - \Id \bigr|^2
	\, ,\,
	\bigl|\b^{\nf12}\a_{*}^{-1}(U; \a)\b^{\nf12} - \Id \bigr|^2
	\Bigr\}
	} \qquad &
	\notag \\ &
	\leq
	10\sup_{|e|=1}
	J(U,\b^{-\nf12}e,\b^{\nf12} e; \a)
	\bigl( 1 + \sup_{|e|=1}
	J(U,\b^{-\nf12}e,\b^{\nf12} e; \a)  \bigr)
	\,.
	\label{e.aratios.to.J}
\end{align}
Combining~\eqref{e.aratios.to.J} with~\eqref{e.bound.for.J.general.preplugin} gives~\eqref{e.coarse.grained.matrix.deviation.by.J}.
\end{proof}

For every~$\lambda>0$, bounded convex domain~$U\subseteq\Rd$ and~$p,q\in\Rd$,
\begin{equation}
	J(U,p,q\,;\a) = J(U,\lambda^{-\nf12} p, \lambda^{\nf12}q \,; \lambda \a) \,.
	\label{e.homogeneous.J}
\end{equation}

\subsection{Coarse-grained elliptic and functional inequalities}
\label{ss.cg.ellipticity}

Elliptic estimates on a cube require information about the coefficient on its subcubes. Following~\citet[Definition~2.1]{AK.HC}, we collect the coarse-grained matrices over all smaller scales into the ellipticity constants~$\lambda_{s,q}$ and~$\Lambda_{s,q}$. The weights discount fine scales, allowing the resulting Poincar\'e and Caccioppoli inequalities to depend on multiscale energy averages rather than the pointwise ellipticity ratio.

\smallskip

The weights~$(1-3^{-s})3^{-ns}$ sum to one over~$n\geq0$. For~$0<s\leq2$, their normalizing factor satisfies~$\frac49 s\leq1-3^{-s}\leq(\log3)s$.

\begin{definition}[Coarse-grained ellipticity constants]
For~$s>0$,~$q\in[1,\infty)$,~$m\in\Z$, and a uniformly elliptic field~$\a:\cu_m\to\Rsymp$, define the \textbf{coarse-grained ellipticity constants} by
\begin{equation}
	\left\{
	\begin{aligned}
		& {\Lambda}_{s,q}(\cu_m\,;\a)
		\coloneqq
		\biggl(
		(1-3^{-sq}) \sum_{k=-\infty}^{m}
		3^{-sq(m-k)}
		\max_{z\in 3^k\Zd \cap \cu_m}
		\bigl| \a(z+\cu_k; \a) \bigr|^{\nf q2}
		\biggr)^{\!\nf{2}{q}}
		\,, \\  &
		{\lambda}_{s,q}(\cu_m;\a)
		\coloneqq
		\biggl((1-3^{-sq}) \sum_{k=-\infty}^{m}
		3^{-sq(m-k)}
		\max_{z\in 3^k\Zd \cap \cu_m}
		\bigl| \a_{*}^{-1}(z+\cu_k; \a) \bigr|^{\nf q2}
		\biggr)^{\!- \nf{2}{q}}
		\,,
	\end{aligned}
	\right.
	\label{e.coarse.grained.ellipticity}
\end{equation}
and, for~$q=\infty$,
\begin{equation}
	\left\{
	\begin{aligned}
		& {\Lambda}_{s,\infty}(\cu_m\,;\a)
		\coloneqq
		\sup_{k\in(-\infty,m]\cap \Z}
		3^{-2s(m-k)}
		\max_{z\in 3^k\Zd \cap \cu_m}
		\bigl| \a(z+\cu_k; \a) \bigr|
		\,, \\  &
		{\lambda}_{s,\infty}(\cu_m;\a)
		\coloneqq
		\biggl(
		\sup_{k\in(-\infty,m]\cap \Z}
		3^{-2s(m-k)}
		\max_{z\in 3^k\Zd \cap \cu_m}
		\bigl| \a_{*}^{-1}(z+\cu_k; \a) \bigr|
		\biggr)^{\!-1}
		\,.
	\end{aligned}
	\right.
	\label{e.coarse.grained.ellipticity.infty}
\end{equation}
When the exponent~$q$ is omitted, it is understood that~$q=1$, that is,~$\Lambda_{s}(\cu_m;\a)=\Lambda_{s,1}(\cu_m;\a)$ and~$\lambda_{s}(\cu_m;\a)=\lambda_{s,1}(\cu_m;\a)$. Equation~\eqref{e.coarse.grained.ellipticity.infty} is the limit of~\eqref{e.coarse.grained.ellipticity} as~$q\to\infty$.  On~$y+\cu_m$, every subcube~$z+\cu_k$ in these definitions is replaced by~$y+z+\cu_k$, with~$z\in3^k\Zd\cap\cu_m$.
\end{definition}
For every~$m\in\Z$,~$q\in [1,\infty]$ and~$t,s\in (0,\infty)$ with~$t<s$, we have
\begin{equation}
	\lambda_{t,q}(\cu_m; \a)
	\leq
	\lambda_{s,q} (\cu_m; \a)
	\leq |\a_*^{-1}(\cu_m;\a) |^{-1}
	\leq
	|\a(\cu_m;\a)|
	\leq
	\Lambda_{s,q}(\cu_m; \a)
	\leq
	\Lambda_{t,q}(\cu_m; \a)
	\,.
	\label{e.ellipticities.monotone.ordered}
\end{equation}
There exists~$C(d)<\infty$ such that, for every~$s\in(0,1]$ and~$p,q\in[1,\infty]$ with~$p\leq q$,
\begin{equation*}
	\Lambda_{s,q} (\cu_m;\a)
	\leq
	C s^{\nf2q - \nf2p}
	\Lambda_{s,p} (\cu_m;\a)
	\leq
	C s^{\nf2q - 2}
	\Lambda_{s} (\cu_m;\a)
	\,,
\end{equation*}
and
\begin{equation*}
	\lambda^{-1}_{s,q} (\cu_m;\a)
	\leq
	C s^{\nf2q - \nf2p}
	\lambda^{-1}_{s,p} (\cu_m;\a)
	\leq
	C s^{\nf2q -2}
	\lambda^{-1}_{s} (\cu_m;\a)
	\,.
\end{equation*}
For every~$k\in(-\infty,m]\cap\Z$, applying~\eqref{e.ellipticities.monotone.ordered} on the subcubes and comparing the definitions gives
\begin{equation}
	\left\{
	\begin{aligned}
		&
		\max_{z\in 3^k\Zd\cap \cu_m}
		\bigl| \a(z+\cu_k; \a) \bigr|
		\leq
		\max_{z\in 3^k\Zd\cap \cu_m}
		\Lambda_{s,q} (z+\cu_k; \a)
		\leq 3^{2s(m-k)} \Lambda_{s,q} (\cu_m;\a)
		\,,
		\\ &
		\max_{z\in 3^k\Zd\cap \cu_m}
		\bigl| \a_{*}^{-1}(z+\cu_k; \a) \bigr|
		\leq
		\max_{z\in 3^k\Zd\cap \cu_m}
		\lambda_{s,q}^{-1} (z+\cu_k; \a)
		\leq 3^{2s(m-k)} \lambda_{s,q}^{-1} (\cu_m;\a)
		\,.
	\end{aligned}
	\right.
	\label{e.bound.one.cube.by.lambdas}
\end{equation}

\smallskip

To compare the equation with a constant-coefficient equation, we also need the two coarse-grained matrices to be close to the same matrix~$\a_0$. The quantities~$J(z+\cu_l,\a_0^{-\nf12}e,\a_0^{\nf12}e;\a)$ measure this error on each subcube. The following definition combines them across scales; its different indices allow either averaging or maximization over the subcubes.

\begin{definition}[Coarse-graining error]
\label{d.mathcal.E}
For~$s\in(0,1]$,~$p,q\in[1,\infty]$,~$m,n\in\Z$ with~$n\leq m$, a uniformly elliptic field~$\a:\cu_m\to\Rsymp$, and a matrix~$\a_0\in\Rsymp$, define the \textbf{coarse-graining error} by
\begin{align*}
	\lefteqn{
	\mathcal{E}_{s,p,q} (\cu_m , n; \a, \a_0 )
	} \ \notag \\ &
	=
	\left\{
	\begin{aligned}
		&
		\biggl(
		(1-3^{-sq})
		\sum_{l=-\infty}^n
		3^{-sq(n-l)}
		\biggl(
		\avsum_{z \in 3^l\Zd \cap \cu_m}
		\max_{|e| =1} \bigl( J\bigl(z {+} \cu_l,\a_0^{-\nf 12} e, \a_0^{\nf 12} e \,; \a\bigr) \bigr)^{\nf p2}
		\biggr)^{\!\nf qp}
		\biggr)^{\!\nf1q}
		\,,
		\quad  p,q<\infty\,,
		\\ &
		\biggl(
		(1-3^{-sq})
		\sum_{l=-\infty}^n
		3^{-sq(n-l)} \!\!
		\max_{z \in 3^l\Zd \cap \cu_m}
		\max_{|e| =1} \bigl( J\bigl(z {+} \cu_l,\a_0^{-\nf 12} e, \a_0^{\nf 12} e \,; \a\bigr) \bigr)^{\nf q2}
		\biggr)^{\!\nf1q}
		\,,
		\quad p=\infty,\ q<\infty\,,
		\\
		&
		\sup_{l \in (-\infty,n]\cap \Z}
		3^{-s(n-l)}
		\biggl(
		\avsum_{z \in 3^l\Zd \cap \cu_m}
		\max_{|e| =1} \bigl( J\bigl(z {+} \cu_l,\a_0^{-\nf 12} e, \a_0^{\nf 12} e \,; \a\bigr) \bigr)^{\nf p2}
		\biggr)^{\!\nf 1p}
		\,,
		\quad p<\infty,\ q=\infty \,,
		\\ &
		\sup_{l \in (-\infty,n]\cap \Z}
		3^{-s(n-l)}
		\max_{z \in 3^l\Zd \cap \cu_m}
		\max_{|e| =1} \bigl( J\bigl(z {+} \cu_l,\a_0^{-\nf 12} e, \a_0^{\nf 12} e \,; \a\bigr) \bigr)^{\nf 12}
		\,,
		\quad p = q = \infty \,.
	\end{aligned}
	\right.
\end{align*}
On~$y+\cu_m$, every subcube~$z+\cu_l$ in this definition is replaced by~$y+z+\cu_l$, with~$z\in3^l\Zd\cap\cu_m$.  We use the abbreviation
\begin{equation*}
	\mathcal{E}_{s,p,q} (\cu_m ; \a, \a_0 )
	\coloneqq
	\mathcal{E}_{s,p,q} (\cu_m , m; \a, \a_0 )
	\,.
\end{equation*}
When the exponent~$q$ is omitted, it is understood that~$q=1$.
\end{definition}

The definition gives, for~$m\in\Z$,~$p\in[2,\infty]$,~$q\in[1,\infty]$, and~$0<t<s\leq1$,
\begin{equation*}
	\max_{|e| =1}  J\bigl(\cu_m,\a_0^{-\nf 12} e, \a_0^{\nf 12} e \,; \a\bigr)^{\nf12}
	\leq
	\mathcal{E}_{s,p,q} (\cu_m ; \a, \a_0 )
	\leq
	\mathcal{E}_{t,p,q} (\cu_m ; \a, \a_0 )
	\,,
\end{equation*}
and, for every~$k\in (-\infty,m]\cap \Z$,
\begin{equation}
	\left\{
	\begin{aligned}
		&
		\max_{z\in 3^k\Zd\cap \cu_m}
		\mathcal{E}_{s,\infty,q} (z+\cu_k;\a, \a_0)
		\leq 3^{s(m-k)} \mathcal{E}_{s,\infty,q}  (\cu_m;\a,\a_0)
		\,,
		\\ &
		\avsum_{z\in 3^k\Zd\cap \cu_m}
		\mathcal{E}_{s,p,q} (z+\cu_k;\a, \a_0)
		\leq 3^{(s+d(\frac1p-\frac1q)_+)(m-k)}
		\mathcal{E}_{s,p,q}  (\cu_m;\a,\a_0)
		\,.
	\end{aligned}
	\right.
	\label{e.bound.one.cube.by.mathcalE}
\end{equation}

\smallskip

For a positive scalar~$\a_0$,~$s\in(0,1]$ and~$q\in[1,\infty]$,
\begin{align}
	\frac12 \mathcal{E}_{s,\infty,q}
	(\cu_m;\a,\a_0)^2
	&
	\leq
	\max\bigl\{
	|\a_0|^{-1} \Lambda_{s,q}(\cu_m;\a)
	\,,
	|\a_0| \lambda_{s,q}^{-1} (\cu_m;\a)
	\bigr\}
	\notag \\ &
	\leq
	\bigl(1+\sqrt2\,\mathcal{E}_{s,\infty,q}(\cu_m;\a,\a_0)\bigr)^2
	\,.
	\label{e.bound.Lambdas.by.Es}
\end{align}
To see the upper bound, let~$J_U=\max_{|e|=1}J(U,\a_0^{-\nf12}e,\a_0^{\nf12}e;\a)$. The ordering~$\a_*(U)\leq\a(U)$ and~\eqref{e.coarse.grained.matrices} imply that each eigenvalue~$t$ of~$\a_0^{-1}\a(U)$ or~$\a_0\a_*^{-1}(U)$ satisfies~$t+t^{-1}-2\leq2J_U$. Thus~$\sqrt t\leq1+\sqrt{2J_U}$ and~$t\leq1+2J_U+\sqrt{2J_U}$. The first inequality, followed by the weighted~$\ell^q$ triangle inequality, gives the upper bound in~\eqref{e.bound.Lambdas.by.Es}. For the lower bound, use~$\sqrt{J_U}\leq2^{-\nf12} (|\a_0^{-1}\a(U)|^{\nf12}+|\a_0\a_*^{-1}(U)|^{\nf12})$ and the same triangle inequality. The second eigenvalue bound also gives
\begin{equation*}
	\max\{\a_0^{-1}\Lambda_{s,2}(\cu_m;\a),\a_0\lambda_{s,2}^{-1}(\cu_m;\a)\}
	\leq1+2\mathcal E_{s,\infty,2}^2(\cu_m;\a,\a_0)+2\mathcal E_{s,\infty,2}(\cu_m;\a,\a_0)\,.
\end{equation*}
Moreover, in the case~$q = 2$,
\begin{equation*}
	1+\frac12 \mathcal{E}_{s,\infty,2}
	(\cu_m;\a,\a_0)^2
	\leq
	\max\bigl\{
	|\a_0|^{-1} \Lambda_{s,2}(\cu_m;\a)
	\,,
	|\a_0| \lambda_{s,2}^{-1} (\cu_m;\a)
	\bigr\} \,.
\end{equation*}
By Definition~\ref{d.mathcal.E}, for every~$p,q \in [1,\infty]$ and~$s \in (\nf{d}{p},1)$,
\begin{equation*}
	\mathcal{E}_{s, \infty, q} (\cu_m ; \a, \a_0 )
	\leq
	\Bigl(
	\frac{(1-3^{-sq})}{(1-3^{-q(s-\frac{d}{p})})}
	\Bigr)^{\nf1q}
	\mathcal{E}_{s - \nf{d}{p}, p, q}(\cu_m;\a,\a_0) \,.
\end{equation*}
For~$1\leq q\leq p<\infty$ and~$s>0$, we have~$\mathcal{E}_{s,\infty,q}(\cu_m;\a,\a_0)\leq\mathcal{E}_{\nf{sq}{p},\infty,p}(\cu_m;\a,\a_0)$. The same definitions give
\begin{equation}
	\Lambda_{s, q}(\cu_m; \a)
	\leq
	\Lambda_{\nf{sq}{p}, p}(\cu_m; \a)
	\qand
	\lambda_{s, q}^{-1}(\cu_m; \a)
	\leq
	\lambda_{\nf{sq}{p}, p}^{-1}(\cu_m; \a) \,.
	\label{e.compareLambdaqs}
\end{equation}
By H\"older's inequality, there exists a universal constant~$C<\infty$ such that, for~$0<s<t\leq1$,
\begin{equation}
	\mathcal{E}_{t,\infty,1}(\cu_m , n; \a, \a_0)
	\leq
	\frac{Ct}{(t-s)^{\nf12}\, s^{\nf12}} \, \mathcal{E}_{s,\infty,2}(\cu_m , n; \a, \a_0)
	\,.
	\label{e.mathcalE.one.vs.mathcalE.two}
\end{equation}

\smallskip

The coarse-grained ellipticity constants satisfy the following stability estimate, which is~\citet[Lemma~2.5, p.~33]{ASD}.

\begin{lemma}[Stability of the coarse-grained ellipticity constants]
\label{l.lambdas.stability}
There exists a constant~$C(d)<\infty$ such that, for every~$q\in[1,\infty]$,~$s,t\in(0,\nf12]$ with~$s<t$, and~$x\in\cu_0$ such that~$x+\cu_{-1}\subseteq\cu_0$, we have the estimates
\begin{equation}
	\lambda_{t,q}^{-1}(x+\cu_{-1}\,;\a)
	\leq
	\frac{C}{1-2s} \Bigl( \frac{t}{t-s}\Bigr)^{\! \nf 2q} \lambda_{s,q}^{-1}(\cu_{0}\,;\a) \,,
	\label{e.lambda.stability}
\end{equation}
\begin{equation}
	\Lambda_{t,q}(x+\cu_{-1}\,;\a)
	\leq
	\frac{C}{1-2s} \Bigl( \frac{t}{t-s}\Bigr)^{\! \nf 2q} \Lambda_{s,q}(\cu_{0}\,;\a)\,,
	\label{e.big.Lambda.stability}
\end{equation}
and, for every~$\a_0\in\Rsymp$ (see~\citet[(2.62), p.~33]{ASD}),
\begin{equation*}
	\mathcal{E}_{t,\infty,2}(x+\cu_{-1}\,;\a,\a_0)
	\leq
	\frac{C}{(1-2s)^{\nf12}} \Bigl( \frac{t}{t-s}\Bigr)^{\! \nf 12}
	\mathcal{E}_{s,\infty,2}(\cu_0\,;\a,\a_0)
	\,.
\end{equation*}
\end{lemma}

The next proposition is the weak-norm estimate of~\citet[Lemma~2.2, p.~27]{AK.HC}, stated there for the gradient and flux of a solution. The same argument uses~\eqref{e.energymaps} and its dual on each subcube, so it applies to every~$u\in H^1(\cu_m)$ and every~$\h\in L^2_{\mathrm{sol}}(\cu_m)$. In the flux case, a zero-trace test function on a subcube extends by zero to one on~$\cu_m$, so the restriction of~$\h$ is still solenoidal. The positive-order estimate~\eqref{e.CG.Poincare.mean.zero} below corresponds to~\citet[Lemma~2.3, p.~27]{AK.HC}.
\begin{proposition}[Coarse-grained Poincar\'e inequality]
\label{p.coarse.grained.poincare}
For every~$s\in (0,1]$,~$q \in [1,\infty]$,~$m \in \Z$, every~$u\in H^1(\cu_m)$ and~$\h \in L^2_{\mathrm{sol}}(\cu_m)$,
\begin{equation}
	3^{-s m} [ \nabla u ]_{\Besov{-s}{2}{q}(\cu_m)}
	\leq
	(1-3^{-sq})^{-\nf1q}
	\lambda_{s,q}^{-\nf12}  (\cu_m; \a) \| \a^{\nf 12} \nabla u \|_{\underline{L}^2(\cu_m)}
	\,,
	\label{e.besov.grad.poincare}
\end{equation}
and
\begin{equation}
	3^{-s m} [ \h ]_{\Besov{-s}{2}{q}(\cu_m)}
	\leq
	(1-3^{-sq})^{-\nf1q}
	\Lambda^{\nf12 }_{s,q} (\cu_m; \a)  \| \a^{-\nf 12} \h \|_{\underline{L}^2(\cu_m)}
	\,.
	\label{e.besov.flux.poincare}
\end{equation}
\end{proposition}

We combine the functional inequality
\begin{equation}
	[ u ]_{\underline{B}^{1-s}_{2,\infty}(\cu_0)}
	\leq
	C(d) [ \nabla u ]_{\Besov{-s}{2}{1}(\cu_0)}
	\,,
	\label{e.nabla.u.detach}
\end{equation}
with~\eqref{e.besov.grad.poincare} and the bound~$\| u - (u)_{\cu_0} \|_{\underline{L}^2(\cu_0)} \leq C \liminf_{s\uparrow 1}  [ u ]_{\underline{B}^{1-s}_{2,\infty}(\cu_0)}$, and let~$s\uparrow1$. This gives, for some~$C(d)<\infty$ and every~$u\in H^1(\cu_0)$,
\begin{equation}
	\| u - (u)_{\cu_0} \|_{\underline{L}^2(\cu_0)}
	\leq
	C
	\lambda_{1,1}^{-\nf12}  (\cu_0; \a) \| \a^{\nf 12} \nabla u \|_{\underline{L}^2(\cu_0)}
	\,.
	\label{e.CG.Poincare.mean.zero}
\end{equation}
For zero boundary data we also have the following estimate. Let~$w\in H^1_0(\cu_0)$ be the \textbf{torsion function} of~$\cu_0$, that is, the solution of
\begin{equation*}
	\left\{
	\begin{aligned}
		& -\Delta w = 1 & \text{in} & \ \cu_0\,, \\
		& w = 0 & \text{on} & \ \partial \cu_0\,.
	\end{aligned}
	\right.
\end{equation*}
Then, for every~$u\in H^1_0(\cu_0)$,
\begin{equation*}
	\bigl| ( u )_{\cu_0} \bigr|
	=
	\biggl| \fint_{\cu_0}
	\nabla w \cdot \nabla u
	\biggr|
	\leq
	C[ \nabla w ]_{\underline{B}^{1}_{2,\infty}(\cu_0)}
	[ \nabla u ]_{\Besov{-1}{2}{1}(\cu_0)}
	\leq
	C [ \nabla u ]_{\Besov{-1}{2}{1}(\cu_0)}\,.
\end{equation*}
Combining this with~\eqref{e.CG.Poincare.mean.zero} yields, for some~$C(d)<\infty$ and every~$u\in H^1_0(\cu_0)$,
\begin{equation}
	\| u \|_{\underline{L}^2(\cu_0)}
	\leq
	C
	\lambda_{1,1}^{-\nf12}  (\cu_0; \a) \| \a^{\nf 12} \nabla u \|_{\underline{L}^2(\cu_0)}
	\,.
	\label{e.CG.Poincare.trace.zero}
\end{equation}

\subsection{Sensitivity estimates}
\label{ss.sensitivity}

Adding fields changes the scalar coefficient by multiplication. To repeat the renormalization argument, we must control this change relative to the coarse-grained matrices, even when the original coefficient has a large ellipticity ratio. The first lemma compares the matrices through the multiplicative ratio of two coefficients. The second estimates~$J$ itself, allowing a constant rescaling before comparing the fields. Both follow by using the minimizers for one coefficient as competitors for the other.

\begin{lemma}
\label{l.sensitivity.general}
Let~$U$ be a bounded convex domain and let~$\a,\b:U\to(0,\infty)$ be scalar, uniformly elliptic coefficient fields. Then
\begin{equation*}
	\bigl|\b_*^{-\nf12}(U)\a_*(U)\b_*^{-\nf12}(U)-\Id\bigr|
	\leq\|\a\b^{-1}-1\|_{L^\infty(U)}\vee\|\a^{-1}\b-1\|_{L^\infty(U)}
	\,,
\end{equation*}
and
\begin{equation*}
	\bigl|\b^{-\nf12}(U)\a(U)\b^{-\nf12}(U)-\Id\bigr|
	\leq\|\a\b^{-1}-1\|_{L^\infty(U)}\vee\|\a^{-1}\b-1\|_{L^\infty(U)}
	\,.
\end{equation*}
\end{lemma}
\begin{proof}
Both inequalities follow from the more precise bounds
\begin{equation}
	\left\{
	\begin{aligned}
		&
		\bigl( \inf_U (\a \b^{-1})  \bigr) \b(U) \leq\a(U) \leq \bigl( \sup_U (\a \b^{-1}) \bigr) \b(U)\,,
		\\
		&
		\bigl( \inf_U (\a^{-1} \b)  \bigr) \b^{-1}_*(U) \leq
		\a^{-1}_*(U) \leq
		\bigl( \sup_U (\a^{-1} \b )  \bigr)
		\b^{-1}_*(U)
		\,.
	\end{aligned}
	\right.
	\label{e.coarse-grained-matrix-inequalities.general}
\end{equation}
For the upper bound in the first line of~\eqref{e.coarse-grained-matrix-inequalities.general}, the variational formula~\eqref{e.variational.a} gives, for every~$e \in \Rd$,
\begin{align*}
	\frac{1}{2} e \cdot \a(U) e
	=
	\inf_{v \in e\cdot x + H_0^1(U)} \fint_{U} \frac{1}{2} \nabla v \cdot \a \nabla v
	&
	\leq \sup_U (\a \b^{-1})  \inf_{v \in e\cdot x + H_0^1(U)} \fint_{U} \frac{1}{2} \nabla v \cdot \b \nabla v
	\notag \\ &
	=
	\sup_U (\a \b^{-1})
	\frac{1}{2} e \cdot \b(U) e
	\,.
\end{align*}
The dual variational formula~\eqref{e.variational.astar} gives, for every~$e\in\Rd$,
\begin{align*}
	\frac{1}{2} e \cdot \a_*^{-1}(U) e
	=
	\min_{\h \in e + \Lsolo(U)} \fint_{U} \frac{1}{2} \h \cdot \a^{-1} \h
	&
	\leq
	\sup_U(  \a^{-1} \b )
	\min_{\h \in e + \Lsolo(U)} \fint_{U} \frac{1}{2} \h \cdot \b^{-1} \h
	\notag \\ &
	=
	\sup_U(  \a^{-1} \b )
	\frac{1}{2} e \cdot \b^{-1}_*(U) e
	\,.
\end{align*}
Interchanging~$\a$ and~$\b$ gives the lower bounds in~\eqref{e.coarse-grained-matrix-inequalities.general}.
\end{proof}

\begin{lemma}
\label{l.J.sensitivity}
Let~$U\subseteq\Rd$ be a bounded convex domain and~$\a,\b:U\to(0,\infty)$ be scalar, uniformly elliptic coefficient fields. Then, for every~$\lambda \in (0,\infty)$,~$p,q \in \Rd$ and~$\delta \in (0,1]$,
\begin{align}
	\lefteqn{
	J(U, \lambda^{-\nf12}  p, \lambda ^{\nf12} q \, ; \b)
	}  \  &
	\notag \\ &
	\leq
	(1+\delta) J(U,p,q\,;\a)
	+
	\frac3{\delta} \Bigl(
	\| \lambda \b^{-1} \a - 1\|_{L^\infty(U)}^2  +
	\| \lambda^{-1} \b \a^{-1} - 1\|_{L^\infty(U)}^2
	\Bigr) \bigl( J(U,p,q\,; \a)  + p\cdot q \bigr)\,.
	\label{e.J.sensitivity}
\end{align}
\end{lemma}
\begin{proof}
By~\eqref{e.homogeneous.J}, we may reduce to the case~$\lambda=1$. Let~$v$ and~$u$ be the maximizers of~$J(U,p,0\,; \a)$ and~$J(U,0,q\,; \a)$, respectively. Then~$-v+c \in p\cdot x + H_0^1(U)$ for an additive constant~$c$, and~$\a \nabla u \in q + \Lsolo(U)$. Therefore, by~\eqref{e.variational.a} and~\eqref{e.variational.astar},
\begin{align*}
	J(U,p, 0 \, ; \b)
	=
	\min_{w \in p\cdot x + H_0^1(U) } \fint_{U} \frac12 \nabla w \cdot \b \nabla w
	&
	\leq
	\fint_{U} \frac12 \nabla v \cdot \b \nabla v
	\notag \\ &
	= J(U,p,0\,;\a)
	+
	\fint_{U} \frac12 \nabla v \cdot (\b \a^{-1} - 1) \a \nabla v
	\,,
\end{align*}
and
\begin{align*}
	J(U,0,q\, ; \b)
	= \min_{\h \in q + \Lsolo(U) } \fint_{U} \frac12 \h \cdot \b^{-1} \h
	&
	\leq
	\fint_{U} \frac12 \a \nabla u \cdot \b^{-1} \a \nabla u
	\notag \\ &
	=
	J(U,0,q\, ; \a) + \fint_{U} \frac12 \nabla u \cdot (\a \b^{-1} - 1) \a \nabla u
	\,.
\end{align*}
Hence,
\begin{equation*}
	J(U,p,q\, ; \b)
	\leq
	J(U,p,q\, ; \a)
	+
	\fint_{U} \frac12 \nabla v \cdot (\b \a^{-1} - 1) \a \nabla v
	+
	\fint_{U} \frac12 \nabla u \cdot (\a \b^{-1} - 1) \a \nabla u
	\,.
\end{equation*}
Since
\begin{equation*}
	\a \b^{-1} - 1
	= - \a \b^{-1} (\b \a^{-1} - 1) = - (\b \a^{-1} - 1) + \a \b^{-1} (\b \a^{-1} - 1)^2
	\,,
\end{equation*}
we have
\begin{equation*}
	\fint_{U} \frac12 \nabla u (\a \b^{-1} - 1) \cdot \a \nabla u
	\leq
	- \fint_{U} \frac12 \nabla u (\b \a^{-1} - 1) \cdot \a \nabla u
	+ \| \a \b^{-1} (\b \a^{-1} - 1)^2 \|_{L^\infty(U)}
	J(U,0,q\, ; \a)
	\,.
\end{equation*}
Pointwise,~$\a\b^{-1}(\b\a^{-1}-1)^2\leq2(\a\b^{-1}-1)^2+2(\b\a^{-1}-1)^2$, and, hence,
\begin{multline*}
	\| \a \b^{-1} (\b \a^{-1} - 1)^2 \|_{L^\infty(U)}
	J(U,0,q\, ; \a)
	\\
	\leq
	2 \bigl( \| \b \a^{-1} - 1 \|_{L^\infty(U)}^2 + \| \b^{-1} \a - 1 \|_{L^\infty(U)}^2 \bigr) \bigl(J(U,p,q\, ; \a)  + p\cdot q\bigr) \,.
\end{multline*}
By H\"older's and Young's inequalities, for every~$\delta>0$,
\begin{align*}
	\lefteqn{
	\fint_{U} \frac12 \nabla v (\b \a^{-1} - 1) \cdot \a \nabla v
	- \fint_{U} \frac12 \nabla u (\b \a^{-1} - 1) \cdot \a \nabla u
	} \qquad &
	\notag \\ &
	=
	\fint_{U} \frac12 \a^{\nf12} \nabla (v + u)
	(\b \a^{-1} - 1) \cdot \a^{\nf12}\nabla (v-u)
	\notag \\ &
	\leq
	\| \b \a^{-1} - 1\|_{L^\infty(U)}  \bigl( J(U,p,q\,;\a) \bigr)^{\nf12}
	2^{-\nf12} \| \a^{\nf12}\nabla (v-u) \|_{\underline{L}^2(U)}
	\notag \\ &
	\leq
	\delta J(U,p,q\,;\a) + \frac1{8 \delta}
	\| \b \a^{-1} - 1\|_{L^\infty(U)}^2
	\| \a^{\nf12}\nabla (v-u) \|_{\underline{L}^2(U)}^2
	\notag \\ &
	\leq
	\delta J(U,p,q\,;\a) +
	\frac1{2 \delta} \| \b \a^{-1} - 1\|_{L^\infty(U)}^2 \bigl( J(U,p,0\,; \a) + J(U,0,q\,; \a)  \bigr)
	\,.
\end{align*}
Combining these displays gives~\eqref{e.J.sensitivity}.
\end{proof}

\subsection{General coarse-graining estimates}
\label{ss.general.coarse.graining.estimates}

We now turn multiscale ellipticity and coarse-graining bounds into a comparison of solutions. Proposition~\ref{p.general.coarse.graining} controls the gradient and flux differences between heterogeneous and constant-coefficient Dirichlet problems with the same data. It separates the error due to the coefficient from the oscillation of the forcing, so it applies both to harmonic approximation and to inhomogeneous equations.

\smallskip

We first record the gradient, flux and boundary energy estimates needed for a nonzero right-hand side. These are the symmetric-coefficient cases of~\citet[Lemmas~2.11--2.12 and Proposition~2.21]{ASD}; the boundary estimate has the direct proof given below.

\smallskip

The first lemma is taken from~\citet[Lemma~2.11, p.~39]{ASD}.
\begin{lemma}[Coarse-grained estimates with a divergence-form right-hand side]
\label{l.coarse.graining.RHS}
There exists~$C(d)<\infty$ such that the following holds. Fix~$m\in\N$ and~$s\in(0,1)$. Let~$\a:\cu_m\to\Rsymp$ be uniformly elliptic and~$\g\in H^s(\cu_m;\Rd)$. If~$u\in H^1(\cu_m)$ solves
\begin{equation*}
	-\nabla \cdot \a \nabla u = \nabla \cdot \g \quad \text{in} \ \cu_m\,,
\end{equation*}
then we have the estimates
\begin{equation*}
	3^{-sm}\| \nabla u \|_{\Besov{-s}{2}{2}(\cu_m)}
	\leq
	C s^{-\nf32}  \lambda_{\nf s2,2}^{-\nf12}  (\cu_m;\a)
	\| \a^{\nf12} \nabla u\|_{\underline{L}^2(\cu_m)}
	+
	C s^{-3}
	\lambda_{\nf s2,2}^{-1} (\cu_m;\a)
	3^{sm} [ \g ]_{\underline{H}^s(\cu_m) }  \,,
\end{equation*}
\begin{equation*}
	3^{-sm} \| \a \nabla u \|_{\Besov{-s}{2}{2}(\cu_m)}
	\leq
	C s^{-\nf32} \Lambda_{\nf s2,2}^{\nf12}  (\cu_m;\a)
	\| \a^{\nf12} \nabla u\|_{\underline{L}^2(\cu_m)}
	+
	C s^{-\nf 92}
	\frac{\Lambda_{\nf s 2,2}^{\nf12}  (\cu_m;\a)  }{\lambda_{\nf s 2,2}^{\nf12}  (\cu_m;\a)  }
	3^{sm} [ \g ]_{\underline{H}^s(\cu_m) }  \,.
\end{equation*}
Moreover, for every~$h\in H^{1+s}(\cu_m)$, the solutions~$v,w \in H^1(\cu_m)$  of the Dirichlet and Neumann problems
\begin{equation*}
	\left\{
	\begin{aligned}
		& -\nabla \cdot \a\nabla v = \nabla \cdot \g
		& \text{in} & \ \cu_m\,, \\
		& v = h & \text{on} & \ \partial \cu_m\,,
	\end{aligned}
	\right.
	\qand
	\left\{
	\begin{aligned}
		& -\nabla \cdot \a\nabla w = \nabla \cdot \g
		& \text{in} & \ \cu_m\,, \\
		& \mathbf{n} \cdot (\a\nabla w +\g - (\g)_{\cu_m}) = 0 & \text{on} & \ \partial \cu_m
	\end{aligned}
	\right.
\end{equation*}
satisfy the estimate
\begin{align}
	\lefteqn{
	\| \a^{\nf12} \nabla v\|_{\underline{L}^2(\cu_m)}
	+
	\| \a^{\nf12} \nabla w\|_{\underline{L}^2(\cu_m)}
	} \qquad &
	\notag \\ &
	\leq
	Cs^{-3} \lambda_{\nf s 2,2}^{-\nf12} (\cu_m\,;\a)
	3^{sm} [ \g ]_{\underline{H}^s(\cu_m)}
	+ Cs^{-\nf32} \Lambda_{\nf s2,2}^{\nf12}  (\cu_m;\a)
	3^{sm}
	\| \nabla h \|_{\underline{H}^s(\cu_m)}
	\,.
	\label{e.cg.RHS}
\end{align}
\end{lemma}
The same estimates may be used at scale zero without changing the positive-scale statement of the lemma. Given data on~$\cu_0$, apply the lemma on~$\cu_1$ to
\begin{equation*}
	\a'(y)=\a(y/3),\qquad u'(y)=3u(y/3),\qquad
	\g'(y)=\g(y/3),\qquad h'(y)=3h(y/3)\,.
\end{equation*}
The weak equations and boundary conditions are preserved. The gradients and fluxes are unchanged under this dilation, the coarse-grained ellipticity constants have the same values, and each order-$s$ seminorm or negative norm gains precisely the factor canceled by~$3^{sm}$ or~$3^{-sm}$ in the estimates. Thus changing variables gives the estimates on~$\cu_0$ used below.

The next lemma is a symmetric-coefficient adaptation of the boundary Caccioppoli estimate in~\citet[Lemma~2.12, p.~42]{ASD}. Its direct proof below gives the boundary factor~$t^{-3}$ and does not require the additional factor~$(1-2t)^{-1}$ in that source.

\begin{lemma}[Coarse-grained Caccioppoli inequality with right-hand side]
\label{l.coarse.grained.Caccioppoli.RHS.ASD}
There exists~$C(d)<\infty$ such that, for every~$s\in (0,1)$,~$t \in (0,\nf 12)$ with~$1-s-t>0$, every~$x \in \cu_0$, every~$\a \in L^{\infty}(\cu_0 ;\R_{\sym,+}^{d \times d})$,~$u \in H^1(\cu_0)$,~$\g \in H^{2t}(\cu_0;\Rd)$ and~$h \in H^{1+2t}(\cu_0)$ with~$(h)_{\cu_0}=0$ satisfying
\begin{equation*}
	\left\{
	\begin{aligned}
		& -\nabla \cdot \a\nabla u = \nabla \cdot \g
		&  \text{in} & \ \cu_0 \,,
		\\ &
		u = h &  \text{on} & \ (\partial \cu_0) \cap (x + \cu_{-1})\,,
	\end{aligned}
	\right.
\end{equation*}
we have the estimate
\begin{align}
	\lefteqn{
	\| \a^{\nf12} \nabla u \|_{\underline{L}^2((x + \cu_{-2}) \cap \cu_0)}^2
	} \qquad &
	\notag \\ &
	\leq
	\biggl( \frac{C}{1-s-t} \biggr)^{\! 2 + \nf{4s}{(1-s-t)}}
	\biggl( \frac{\Lambda_{s}(\cu_0; \a)}{\lambda_{t}(\cu_0; \a)} \biggr)^{\!\nf{(1-t)}{(1-s-t)}}
	\Bigl(
	{\lambda}_{t}(\cu_0; \a)
	\| u \|_{\underline{L}^2(\cu_{0})}^2
	\notag \\ & \qquad
	+
	t^{-11}
	\lambda_{t}^{-1} (\cu_0\,;\a)
	\| \g \|_{\underline{H}^{2t} (\cu_0)}^2
	+
	t^{-3}
	\Lambda_t(\cu_0; \a)
	\| \nabla h \|_{\underline{H}^{2t}(\cu_0)}^2
	\Bigr)
	\,.
	\label{e.cg.Caccioppoli.with.RHS.ASD}
\end{align}
\end{lemma}
\begin{proof}
In this proof, all norms and coarse-grained ellipticity constants are taken on~$\cu_0$ unless otherwise indicated, and~$\|\cdot\|_2$ denotes~$\|\cdot\|_{\underline{L}^2(\cu_0)}$. Let~$v$ solve the Dirichlet problem in~$\cu_0$ with forcing~$\nabla\cdot\g$ and boundary values~$h$ on all of~$\partial\cu_0$. By subadditivity,~$\Lambda_{t,2}\leq\Lambda_t$ and~$\lambda_{t,2}^{-1}\leq\lambda_t^{-1}$. Consequently,~\eqref{e.cg.RHS}, applied with~$2t<1$ in place of~$s$, gives
\begin{equation*}
	\|\a^{\nf12}\nabla v\|_2\leq Ct^{-3}\lambda_t^{-\nf12}\|\g\|_{\underline H^{2t}}
	+Ct^{-\nf32}\Lambda_t^{\nf12}
	\|\nabla h\|_{\underline H^{2t}}\,.
\end{equation*}
Since~$v-h\in H_0^1(\cu_0)$ and~$(h)_{\cu_0}=0$, testing~$v-h$ against the torsion function of~$\cu_0$, as in the argument preceding~\eqref{e.CG.Poincare.trace.zero}, controls~$(v)_{\cu_0}$ by the negative Besov seminorms of~$\nabla v$ and~$\nabla h$. Combining this bound with~\eqref{e.CG.Poincare.mean.zero} and with~\eqref{e.besov.grad.poincare} for~$s=q=1$, we obtain~$\|v\|_2\leq C\lambda_t^{-\nf12}\|\a^{\nf12}\nabla v\|_2+C\|\nabla h\|_2$. Apply the homogeneous boundary Caccioppoli inequality~\citep[Proposition~2.9, p.~36]{ASD} to~$u-v$, which has zero trace on the part of~$\partial\cu_0$ where the boundary condition is imposed. Since~$x\in\cu_0$, each side of~$(x+\cu_{-2})\cap\cu_0$ has length at least~$1/18$, so~$|\cu_0|/|(x+\cu_{-2})\cap\cu_0|\leq18^d$. Hence~$\|\a^{1/2}\nabla v\|_{\underline L^2((x+\cu_{-2})\cap\cu_0)}^2\leq18^d\|\a^{\nf12}\nabla v\|_2^2$. The triangle inequality now bounds the left side of~\eqref{e.cg.Caccioppoli.with.RHS.ASD} by
\begin{equation*}
	C\left(\frac{C}{1-s-t}\right)^{2+\nf{4s}{(1-s-t)}}
	\left(\frac{\Lambda_s}{\lambda_t}\right)^{\nf{s}{(1-s-t)}}\Lambda_s
	(\|u\|_2^2+\|v\|_2^2)+C(d)\|\a^{\nf12}\nabla v\|_2^2\,.
\end{equation*}
Insert the preceding estimates and absorb this dimensional multiple of~$\|\a^{\nf12}\nabla v\|_2^2$ to obtain~\eqref{e.cg.Caccioppoli.with.RHS.ASD}.
\end{proof}
The homogeneous case of the next proposition is~\citet[Proposition~2.5, p.~29]{AK.HC}; the forcing term follows from Lemma~\ref{l.coarse.grained.Caccioppoli.RHS.ASD}.
\begin{proposition}[Coarse-grained Caccioppoli inequality]
\label{p.coarse.grained.caccioppoli}
There exists a constant~$C(d)<\infty$ such that, for every~$s\in (0,1)$,~$t\in (0,\nf12)$ with~$s+t < 1$, every coefficient field~$\a \in L^{\infty}(\cu_0 ;\R_{\sym,+}^{d \times d})$ and every~$u \in H^1(\cu_0)$ and~$\f \in H^{2t}(\cu_0;\Rd)$ satisfying
\begin{equation*}
	-\nabla \cdot \a\nabla u = \nabla \cdot \f \quad \text{in} \ \cu_0
	\,,
\end{equation*}
we have the estimate
\begin{align*}
	\| \a^{\nf12} \nabla u \|_{\underline{L}^2(\cu_{-1})}^2
	&
	\leq
	\biggl( \frac{C}{1-s-t} \biggr)^{\! 2 + \nf{4s}{(1-s-t)}}
	\biggl( \frac{\Lambda_{s,1}(\cu_0; \a)}{\lambda_{t,1}(\cu_0; \a)} \biggr)^{\!\nf{s}{(1-s-t)}}
	{\Lambda}_{s,1}(\cu_0; \a)
	\| u \|_{\underline{L}^2(\cu_{0})}^2
	\notag \\ & \qquad
	+
	t^{-11}
	\biggl( \frac{C}{1-s-t} \biggr)^{\! 2 + \nf{4s}{(1-s-t)}}
	\biggl( \frac{\Lambda_{s,1}(\cu_0; \a)}{\lambda_{t,1}(\cu_0; \a)} \biggr)^{\!\nf{(1-t)}{(1-s-t)}}
	\lambda_{t,1}^{-1} (\cu_0\,;\a)
	\| \f \|_{\underline{H}^{2t} (\cu_0)}^2
	\,.
\end{align*}
\end{proposition}
\begin{proof}
Apply Lemma~\ref{l.coarse.grained.Caccioppoli.RHS.ASD} with zero boundary datum on each cube of a fixed finite cover of~$\cu_{-1}$ by cubes of side~$3^{-2}$ compactly contained in~$\cu_0$, and sum the resulting estimates. The forcing term has the displayed exponent since~$(1-t)/(1-s-t)=1+s/(1-s-t)$.
\end{proof}

The preceding estimates allow us to compare gradients and fluxes on subcubes and sum these comparisons in a negative norm. The resulting estimate is the symmetric-coefficient case of~\citet[Proposition~2.21]{ASD}.

\begin{proposition}[General coarse-graining estimate]
\label{p.general.coarse.graining.ASD}
\label{p.general.coarse.graining}
Let~$m,n\in\Z$ with~$n<m$, let~$\a \in L^{\infty}(\cu_m ;\allowbreak\R_{\sym,+}^{d \times d})$ and~$\a_0 \in (0,\infty)$. Fix~$p \in [2,\infty)$,~$s, s_1, s_2 \in (0,1)$ with~$s_1 < s < s_2 < 1$ and~$sp^{\prime}<1$, where~$p^{\prime}=p/(p-1)$. There exists~$C(p,d)<\infty$ such that, for every~$\g \in W^{s_2,p}(\cu_m;\Rd)$ and~$u,v\in H^1(\cu_m)$ satisfying
\begin{equation*}
	\left\{
	\begin{aligned}
		& -\nabla \cdot \a \nabla u  = \nabla \cdot \g  & \text{in} & \ \cu_m\,,
		\\
		& -\nabla \cdot \a_0 \nabla v = \nabla \cdot \g  & \text{in} & \ \cu_m\,,
		\\
		&
		u-v \in H_0^1(\cu_m)\,,
	\end{aligned}
	\right.
\end{equation*}
we have the estimate
\begin{align}
	\lefteqn{
	3^{-ms} \a_0 \|  \nabla u - \nabla v \|_{\Wminusul{-s}{p} (\cu_m)}
	+
	3^{-ms} \| \a \nabla u - \a_0 \nabla v \|_{\Wminusul{-s}{p} (\cu_m)}
	}
	\qquad &
	\notag \\ &
	\leq
	C s^{-1-\nf{1}{p^{\prime}}} \a_0^{\nf12}
	\mathcal{E}_{s_1, \infty, 1} (\cu_m, n;\a,\a_0)
	\biggl(
	\sum_{k=-\infty}^n
	3^{-(s-s_1) p (n-k)}
	\avsum_{z\in 3^k\Zd\cap\cu_m}
	\| \a^{\nf12} \nabla u \|_{\underline{L}^2(z+\cu_k)}^{p}
	\biggr)^{\!\nf 1p}
	\notag \\ & \qquad
	+
	Cs^{-\nf{11}{2}}
	(s_2 - s)^{-1}
	\bigl(1 + \mathcal{E}^{2}_{\nf{s_1}{2}, \infty, 2}(\cu_m, n;\a,\a_0) \bigr)
	3^{s_2 n}
	[ \g ]_{\underline{W}^{s_2,p}(\cu_m)}
	\,.
	\label{e.general.coarse.graining.estimate.ASD}
\end{align}
\end{proposition}
The fractional seminorms in~\citet{ASD} omit the leading factor~$s$ used here. Converting the negative norms contributes~$s^{-1/p'}$, and converting the forcing seminorm contributes~$s_2^{-1/p}$. The energy coefficient is therefore~$s^{-1-1/p'}$, while the forcing coefficient is bounded by~$s^{-11/2}(s_2-s)^{-1}$, since~$s<s_2$ and~$1/p+1/p'=1$. The gradient difference has zero mean because~$u-v\in H^1_0(\cu_m)$, so its pairing with a test function is unchanged when the mean of the test function is subtracted. The mean-zero dual norm used for this difference in that reference thus gives our full dual norm, up to a factor of two. The flux uses the full dual norm in both conventions. For symmetric coefficients, the doubled variational quantity in that reference is the sum of its two directional energies. Maximizing over the two directions with total squared length one gives the directional maximum defining~$\mathcal E$.

\smallskip

By Sobolev embedding (Appendix~\ref{s.Besov.norm.comparisons}), if~$p-2$ is small in terms of~$s$, then the left side of~\eqref{e.general.coarse.graining.estimate.ASD} controls the quantity~$\a_0 3^{-m}\|u-v\|_{\underline{L}^2(\cu_m)}$.

\smallskip

\endgroup

\section{Estimates for the cutoff coefficients}
\label{s.infrared.cutoffs}

In this section we prove the probabilistic estimates on the cutoff coefficients used in the rest of the paper.  They are of two kinds. Sensitivity bounds (Lemmas~\ref{l.sensitivity.for.cutoffs} and~\ref{lemma:infrared.approx.cutoffs}), the first of which underlies the induction of Section~\ref{s.self.similar}, measure how much the coefficient and its coarse-grained matrices change when the fields of a range of scales are added to the coefficient. Negative Besov bounds (Lemma~\ref{l.aL.negative.Besov}) control the ratio~$\a_m\a_n^{-1}$ of two cutoff coefficients. We also prove, in Proposition~\ref{p.special.two.d.exact.formula}, the planar identity~$\ahom_m=e^{-\tau^2(m+1)}$, which gives the exact two-dimensional exponents.

\smallskip

Throughout this section,~$\P$ is a probability measure satisfying Assumptions~\ref{a.g1}--\ref{a.g4}.  On a given cube, a single field has little effect in two opposite situations. A field whose scale is much smaller than the cube oscillates on it, and is therefore small in a negative Besov norm on the cube. A field whose scale is much larger than the cube is nearly constant on it, and multiplying the coefficient by a constant does not change the solutions of the equation.

\smallskip

Stochastic bounds are written in the Orlicz notation: for~$s,A>0$ we write~$X\leq\O_{\Gamma_s}(A)$ if~$\E[\exp((A^{-1}\max\{X,0\})^s)]\leq2$, and~$X=\O_{\Gamma_s}(A)$ if this holds for both~$X$ and~$-X$; Appendix~\ref{s.orlicz} collects the properties we use.

\smallskip

Recall from~\eqref{e.intro.cutoff.coefficients} that, for~$L\in\N_0$, the cutoff coefficient is~$\a_L(x)=\prod_{k=0}^L\exp(\g_k(x)-\tau^2)$. Assumption~\ref{a.g2}, together with the symmetry~$\g_0\mapsto-\g_0$ in~\ref{a.g3}, gives~$C(d)<\infty$ such that
\begin{equation}
	\tau^2\leq C\gmcdelta^2\,,
	\label{e.tau.bound}
\end{equation}
and~$\bigl|\tau^2-\frac12\E[|\g_0(0)|^2]\bigr|\leq C\gmcdelta^4$. We also let~$\a_{-1}\equiv1$. For every~$m\in\N_0$ and~$n\in\Z$ with~$-1\leq n<m$, we have~$\E[(\a_m\a_n^{-1})(0)]=1$ and~$\E[(\a_m^{-1}\a_n)(0)]=\exp(2\tau^2(m-n))$. Moreover,
\begin{equation}
	\text{the fields~$\a_m \a_n^{-1}$ and~$\exp \bigl( -2 \tau^2 (m-n) \bigr) \a_m^{-1} \a_n$ have the same law}
	\,.
	\label{e.inverse.symmetry.in.law}
\end{equation}

\subsection{Besov space estimates}

We first estimate the oscillation of the ratios~$\a_m\a_n^{-1}$.

\begin{lemma}
\label{l.sensitivity.for.cutoffs}
\label{l.am.Lp.moments}
There exists a constant~$C(d) < \infty$ such that, for every~$n\in\N_0$ and~$k\in\Z$,
\begin{equation}
	\sup_{m \geq n}\sup_{x,y \in \cu_k}  \log \biggl(  \frac{(\a_m \a_n^{-1})(x)}{(\a_m \a_n^{-1})(y)} \biggr)
	\leq
	\O_{\Gamma_2} ( C \gmcdelta 3^{k-n} ) \,,
	\label{e.sensitivity.field}
\end{equation}
and, for every~$\xi \in [1,\infty)$,
\begin{equation}
	\E \Biggl[ \biggl| \sup_{m \geq n} \sup_{x,y \in \cu_k}  \frac{ (\a_m \a_n^{-1})(x)}{ (\a_m \a_n^{-1})(y)}  - 1 \biggr|^{\xi}\Biggr]^{\nf1\xi}
	\leq
	C \xi^{\nf12}  \gmcdelta 3^{k - n }\exp(C \xi \gmcdelta^2 3^{2 (k-n)} ) \,.
	\label{e.infraredhom.approx.cutoffs.large.waves.field}
\end{equation}
Moreover, for every~$m\in\N_0$ with~$m>n$,
\begin{align}
	&
	\bigl\| \log (\a_m \a_n^{-1})  + (m-n) \tau^2 \bigr\|_{L^\infty(\cu_k)}
	\leq
	\O_{\Gamma_2} \bigl( C \gmcdelta \bigl( |m-n|^{\nf12} +  (k- n)_+ \bigr) \bigr) \,,
	\label{e.infrared.approx.cutoffs.field}
\end{align}
and, for every~$\xi \in [1,\infty)$,
\begin{align}
	\lefteqn{
	\E \Bigl[ \bigl\| \a_m \a_n^{-1} - 1 \bigr\|_{L^\infty(\cu_k)}^\xi  \Bigr]^{\nf 1\xi}
	+  \E \Bigl[  \bigl\|  \a_m^{-1}  \a_n -1 \bigr\|_{L^\infty(\cu_k)}^\xi
	\Bigr]^{\nf 1\xi}
	} \qquad &
	\notag \\ &
	\leq
	\bigl( C \xi^{\nf12} \gmcdelta \bigl( (m-n)^{\nf12} +  (k- n)_+  \bigr)  \bigr)
	\exp\bigl( C \xi \gmcdelta^2 \bigl( m-n + (k - n)^2_+ \bigr) \bigr)  \,.
	\label{e.aman.Linfty.moments}
\end{align}
Finally, for every bounded domain~$U\subseteq\Rd$, every~$p\in[1,\infty)$, and every~$m\in\N_0$ and~$n\in\Z$ with~$-1\leq n<m$,
\begin{equation}
	\E \Bigl[
	\bigl\| \a_m\a_{n}^{-1} - \Id \bigr\|_{\underline{L}^p(U) }^p \Bigr]^{\nf1p}
	\leq
	\bigl(Cp^{\nf12} \gmcdelta (m-n)^{\nf12}  \bigr) \exp\bigl( C p \gmcdelta^2 (m-n)  \bigr)
	\,,
	\label{e.am.Lp.moments}
\end{equation}
and
\begin{equation}
	\E \Bigl[
	\bigl\| \a_m^{-1} \a_{n}- \Id \bigr\|_{\underline{L}^p(U) }^p \Bigr]^{\nf1p}
	\leq
	\bigl(Cp^{\nf12} \gmcdelta (m-n)^{\nf12}  \bigr) \exp\bigl( C p \gmcdelta^2 (m-n)  \bigr)
	\,.
	\label{e.am.Lp.moments.inv}
\end{equation}
\end{lemma}
\begin{proof}
We bound the logarithm~$\log(\a_m\a_n^{-1})=\sum_{j=n+1}^m(\g_j-\tau^2)$ on~$\cu_k$ and then exponentiate. Its oscillation gives~\eqref{e.sensitivity.field} and~\eqref{e.infraredhom.approx.cutoffs.large.waves.field}. The main step is the bound~\eqref{e.infrared.approx.cutoffs.field} on its size: we split the sum at the scale of~$\cu_k$, use that the fields of larger scale are nearly constant on~$\cu_k$, and control the remaining fields on the subcubes of side~$3^n$ by a bound on the maximum of many random variables. The moment bounds~\eqref{e.aman.Linfty.moments} and~\eqref{e.am.Lp.moments} then follow by exponentiating, and~\eqref{e.am.Lp.moments.inv} follows from~\eqref{e.inverse.symmetry.in.law}.

\smallskip

We first prove~\eqref{e.sensitivity.field}. By~\eqref{e.coefficient.field.regularity},
\begin{equation*}
	\sup_{m\geq n}\sup_{x,y\in\cu_k}
	\log\frac{(\a_m\a_n^{-1})(x)}{(\a_m\a_n^{-1})(y)}
	\leq
	\sum_{j=n+1}^\infty\sup_{x,y\in\cu_k}|\g_j(x)-\g_j(y)|
	\leq\O_{\Gamma_2}(C3^{k-n}\gmcdelta)\,.
\end{equation*}
Indeed, for~$j\geq k$ the mean value theorem bounds the oscillation by~$C3^k\|\nabla\g_j\|_{L^\infty(\cu_k)} \leq\O_{\Gamma_2}(C3^{k-j}\gmcdelta)$. For~$j<k$, cover~$\cu_k$ by the~$3^{d(k-j)}$ cubes~$z+\cu_j$ with~$z\in3^j\Zd\cap\cu_k$ and use~\eqref{e.max.scaling.O.Gamma.Two} to bound the oscillation by~$2\|\g_j\|_{L^\infty(\cu_k)} \leq\O_{\Gamma_2}(C\gmcdelta\sqrt{1+k-j}) \leq\O_{\Gamma_2}(C\gmcdelta3^{k-j})$. Then~\eqref{e.infraredhom.approx.cutoffs.large.waves.field} follows from~\eqref{e.moments.general.lognormal}, applied to the supremum in~\eqref{e.sensitivity.field}. Since~$\log(\a_m\a_n^{-1})+(m-n)\tau^2=\sum_{j=n+1}^m\g_j$, the bound~\eqref{e.infrared.approx.cutoffs.field} follows once we prove, for a fixed~$m\in\N_0$ with~$m>n$, that
\begin{equation}
	\sup_{x \in \cu_k} \biggl| \sum_{j=n+1}^m
	\g_j (x) \biggr| \leq \O_{\Gamma_2}\Bigl( C \gmcdelta \bigl( (m-n)^{\nf12} +  (k-n)_+ \bigr) \Bigr) \,.
	\label{e.bound-b-inside}
\end{equation}
By~\eqref{e.coefficient.field.regularity},~\eqref{e.CLT.scaling.O.Gamma.Two}, and~\eqref{e.orlicz.triangle.inequality},
\begin{align}
	\sup_{x \in \cu_k} \biggl| \sum_{j= k \vee (n+1)}^m
	\g_j (x)
	\biggr|
	\leq
	\underbrace{\biggl| \sum_{j= k \vee (n+1)}^m
	\g_j (0)
	\biggr|}_{ \leq \O_{\Gamma_2}( \gmcdelta (m-n)^{\nf12})}
	+
	\sum_{j=k \vee (n+1)}^m \underbrace{\| \g_j - \g_j(0) \|_{L^{\infty}(\cu_k)}}_{\leq \O_{\Gamma_2}(C 3^{k - j} \gmcdelta)}
	\leq \O_{\Gamma_2}(C \gmcdelta (m-n)^{\nf12})  \,.
	\label{e.bound.sum.large.scales}
\end{align}
This proves~\eqref{e.bound-b-inside} when~$k\leq n$. Suppose that~$k>n$.  The same three estimates, now on~$\cu_n$ and with the sums stopped at~$m\wedge k$, give
\begin{align*}
	\sup_{x \in \cu_n} \biggl| \sum_{j= n+1}^{m\wedge k}
	\g_j (x)
	\biggr|
	\leq
	\underbrace{\biggl| \sum_{j= n+1}^{m\wedge k}
	\g_j (0)
	\biggr|}_{\O_{\Gamma_2}( \gmcdelta ((m\wedge k)-n)^{\nf12})}
	+
	\sum_{j= n+1}^{m\wedge k} \underbrace{\| \g_j - \g_j(0) \|_{L^{\infty}(\cu_n)}}_{\O_{\Gamma_2}(C 3^{n- j} \gmcdelta)}
	\leq \O_{\Gamma_2}(C  \gmcdelta ((m\wedge k)-n)^{\nf12})  \,.
\end{align*}
The maximum of~$N\geq2$ random variables, each of which is~$\O_{\Gamma_2}(A)$, is~$\O_{\Gamma_2}((3\log N)^{\nf12}A)$; this is~\citet[Lemma~2.3]{ABK.SD}. Applying it with~$N=3^{d(k-n)}$ gives
\begin{equation*}
	\sup_{x \in \cu_k} \biggl| \sum_{j=n+1}^{m\wedge k} \g_j (x) \biggr|
	= \max_{z \in 3^n \Zd \cap \cu_k} \biggl \| \sum_{j=n+1}^{m\wedge k} \g_j \biggr\|_{L^{\infty}(z + \cu_n)}
	\leq  \O_{\Gamma_2} \bigl( C \gmcdelta  (k-n) \bigr) \,.
\end{equation*}
If~$m>k$, then, by~\eqref{e.bound.sum.large.scales}, the remaining sum~$\sum_{j=k+1}^m\g_j$ contributes at most~$\O_{\Gamma_2}(C\gmcdelta(m-n)^{\nf12})$ on~$\cu_k$.  This proves~\eqref{e.bound-b-inside}.

\smallskip

To prove~\eqref{e.aman.Linfty.moments}, we use
\begin{equation*}
	\|\a_m\a_n^{-1}-1\|_{L^\infty(\cu_k)}
	\vee
	\|\a_m^{-1}\a_n-1\|_{L^\infty(\cu_k)}
	\leq \bigl\|\log(\a_m\a_n^{-1})\bigr\|_{L^\infty(\cu_k)}
	\exp\bigl(\bigl\|\log(\a_m\a_n^{-1})\bigr\|_{L^\infty(\cu_k)}\bigr)\,.
\end{equation*}
Equations~\eqref{e.tau.bound} and~\eqref{e.infrared.approx.cutoffs.field}, together with~\eqref{e.moments.general.lognormal} and~\eqref{e.moments.OGamma2} at exponent~$2\xi$, bound the~$L^{2\xi}(\P)$ norms of the random variable~$\|\log(\a_m\a_n^{-1})\|_{L^\infty(\cu_k)}$ and of its exponential. H\"older's inequality then gives~\eqref{e.aman.Linfty.moments}.

\smallskip

At a fixed point, the pointwise version of~\eqref{e.infrared.approx.cutoffs.field} and the exponential moment bound used above give
\begin{equation*}
	\|\a_m(x)\a_n^{-1}(x)-1\|_{L^p(\P)}
	\leq
	Cp^{\nf12}\gmcdelta(m-n)^{\nf12}
	\exp\bigl(Cp\gmcdelta^2(m-n)\bigr)\,.
\end{equation*}
For~$n=-1$, separate the independent factor~$e^{\g_0-\tau^2}$ and apply the estimate with~$n=0$. Stationarity now gives~\eqref{e.am.Lp.moments}. Finally,~\eqref{e.inverse.symmetry.in.law} and~\eqref{e.tau.bound} give~\eqref{e.am.Lp.moments.inv}.
\end{proof}

\begin{lemma}
\label{l.aL.negative.Besov}
There exist~$C(d)<\infty$ and~$c(d)>0$ such that, for every~$s\in (0,1]$ and~$p\in [2,\infty)$ satisfying
\begin{equation}
	s \geq C p \gmcdelta^{2}
	\qquad\text{and}\qquad
	p\gmcdelta^2\log(2+p)\leq c
	\,,
	\label{e.aL.Besov.p.condition}
\end{equation}
and every~$m\in\N_0$ and~$n\in\Z$ with~$-1\leq n<m$, we have the estimates
\begin{equation}
	3^{-sm} \E \Bigl[ \| \a_m\a_n^{-1} - 1 \|_{\Besov{-s}{p}{p} (\cu_m) }^p \Bigr]^{\nf1p}
	\leq
	Cs^{-1}\gmcdelta p^{\nf12} (\log p)
	\,,
	\label{e.aL.Besov}
\end{equation}
and
\begin{equation}
	3^{-sm} \E \Bigl[ \| \exp( - 2 \tau^2 (m-n) ) \a_m^{-1} \a_n - 1  \|_{\Besov{-s}{p}{p} (\cu_m) }^p \Bigr]^{\nf1p}
	\leq
	Cs^{-1} \gmcdelta p^{\nf12} (\log p)
	\,.
	\label{e.aL.Besov.invese}
\end{equation}
\end{lemma}
\begin{proof}
The estimates~\eqref{e.aL.Besov} and~\eqref{e.aL.Besov.invese} are equivalent by~\eqref{e.inverse.symmetry.in.law}; we prove~\eqref{e.aL.Besov}. By~\eqref{e.Besov.weaknorm.def}, its left side is controlled by the~$p$th moments of the averages of~$\a_m\a_n^{-1}-1$ on the cubes of size~$3^j$,~$j\leq m$, and~\eqref{e.aL.negative.Besov.iden1} below splits each such average into the contributions of the fields~$\g_{k+1}$ with~$n\leq k<m$. The proof proceeds in three steps. Step~1 bounds the contribution of one field~$\g_{k+1}$ on a cube of size~$3^j$ with~$k\leq j\leq m$, Step~2 sums these bounds over~$k<j$ and bounds the terms with~$k\geq j$ directly, to bound the average on a single cube, and Step~3 sums over the scales~$j$. The heart of the matter is the decay in~$j-k$ obtained in Step~1, which comes from the independence of~$\g_{k+1}$ from~$\a_m\a_{k+1}^{-1}$ and keeps the sum over~$k<j$ in Step~2 of order~$\log(2+p)$.

\smallskip

When~$n=-1$, the summand with~$k=-1$ in the decomposition~\eqref{e.aL.negative.Besov.iden1} below contains~$\g_0$ and~$\a_m\a_0^{-1}$. Assumptions~\ref{a.g1}--\ref{a.g2}, independence, and Lemma~\ref{l.am.Lp.moments} with~$n=0$ give the bound required for this summand. Every remaining summand has~$k\geq0$. Here and below, Lemma~\ref{l.am.Lp.moments} is applied to a ratio~$\a_m\a_{k+1}^{-1}$ only when~$k+1<m$; when~$k+1=m$ (for instance, when~$m=0$ in the summand just treated), the ratio is identically~$1$ and the bounds used for it hold trivially.

\smallskip

By~\eqref{e.Besov.weaknorm.def},
\begin{align}
	3^{-sm}
	\E \Bigl[ \| \a_m\a_n^{-1} - 1 \|_{\Besov{-s}{p}{p} (\cu_m) }^p \Bigr]^{\nf1p}
	=
	\E \Biggl[
	\sum_{j=-\infty}^m
	s 3^{sp(j-m)}
	\avsum_{z\in 3^j\Zd \cap \cu_m}
	\bigl| \bigl( \a_m\a_n^{-1} \bigr) _{z+\cu_j} - 1 \bigr|^p
	\Biggr]^{\nf 1p}
	\,.
	\label{e.aL.negative.Besov.iden0}
\end{align}
For each~$j \leq m$, we write
\begin{align}
	\fint_{\cu_j}
	\a_m \a_n^{-1}  - 1
	=
	\sum_{k=n}^{m-1}
	\fint_{\cu_j}
	(\a_{m} \a_{k}^{-1} -\a_{m} \a_{k+1}^{-1} )
	=
	\sum_{k=n}^{m-1}
	\fint_{\cu_j}
	\a_m \a_{k+1}^{-1}
	(\exp( \g_{k+1} - \tau^2) - 1 )
	\,.
	\label{e.aL.negative.Besov.iden1}
\end{align}

\smallskip

\emph{Step 1.} We show the existence of~$C(d)<\infty$ such that, for every~$j\in\{k,\ldots,m\}$,
\begin{align}
	\lefteqn{
	\E \Biggl[ \biggl|
	\fint_{\cu_j}
	\a_m \a_{k+1}^{-1}(\exp( \g_{k+1} - \tau^2) - 1 )
	\biggr|^p
	\Biggr]^{\nf1p}
	} \qquad &
	\notag \\ &
	\leq
	C \gmcdelta p^{\nf12}
	\Bigl( 3^{-\frac d2(j-k)}
	+
	\min \bigl\{ p 3^{-d(1 - \frac 1p) (j-k)} \, , 1 \bigr\}  \Bigr)
	\exp\bigl( C p \gmcdelta^2 (m-k)  \bigr)
	\,.
	\label{e.condy.moment.bound}
\end{align}
For each~$k\leq j$, we have
\begin{equation*}
	\fint_{\cu_j}
	\a_m \a_{k+1}^{-1}(\exp( \g_{k+1} - \tau^2) - 1 )
	=
	\avsum_{z \in 3^k\Zd \cap \cu_j}
	\fint_{z+\cu_k}
	\a_m \a_{k+1}^{-1}(\exp( \g_{k+1} - \tau^2) - 1 )
	\,.
\end{equation*}
Condition on~$\a_m\a_{k+1}^{-1}$, which is independent of the field~$\g_{k+1}$. Since~$\E[\exp(\g_{k+1}(x)-\tau^2)]=1$ for every~$x$, each summand then has conditional mean zero. We bound the conditional moments of one summand and then apply Lemma~\ref{l.Rosenthal} to the sum. By the Cauchy--Schwarz inequality,~\eqref{e.aman.Linfty.moments},~\eqref{e.aL.Besov.p.condition}, and~\eqref{e.tau.bound},
\begin{align}
	\lefteqn{
	\E \biggl[
	\biggl|
	\fint_{\cu_k}
	\a_m \a_{k+1}^{-1}(\exp( \g_{k+1} - \tau^2) - 1 )
	\biggr|^p
	\,\bigg\vert\,
	\a_m\a_{k+1}^{-1}
	\biggr] ^{\nf1p}
	} \qquad &
	\notag \\ &
	\leq
	\| \a_m \a_{k+1}^{-1} \|_{\underline{L}^2(\cu_k)}
	\E \biggl[
	\biggl(
	\fint_{\cu_k}
	\bigl| \exp( \g_{k+1} - \tau^2) - 1 \bigr|^2
	\biggr)^{\!\nf p2}
	\biggr]^{\nf1p}
	\leq
	C p^{\nf12}  \gmcdelta
	\| \a_m \a_{k+1}^{-1} \|_{\underline{L}^2(\cu_k)}
	\,.
	\label{e.condiepees}
\end{align}
In particular,
\begin{align*}
	\lefteqn{ \biggl(
	\avsum_{z\in 3^k\Zd \cap \cu_j}
	\E \biggl[
	\biggl|
	\fint_{z+\cu_k}
	\a_m \a_{k+1}^{-1}(\exp( \g_{k+1} - \tau^2) - 1 )
	\biggr|^2
	\,\bigg\vert\,
	\a_m\a_{k+1}^{-1}
	\biggr]
	\biggr)^{\!\nf12}
	} \qquad &
	\notag \\ &
	\leq
	C \gmcdelta
	\biggl(
	\avsum_{z\in 3^k\Zd \cap \cu_j}
	\| \a_m \a_{k+1}^{-1} \|_{\underline{L}^2(z+\cu_k)}^2
	\biggr)^{\!\nf12}
	=
	C \gmcdelta \| \a_m \a_{k+1}^{-1} \|_{\underline{L}^2(\cu_j)}
	\,,
\end{align*}
and
\begin{align*}
	\lefteqn{
	\E \biggl[
	\max_{z\in 3^k\Zd \cap \cu_j}
	\biggl|
	\fint_{z+\cu_k}
	\a_m \a_{k+1}^{-1}(\exp( \g_{k+1} - \tau^2) - 1 )
	\biggr|^p
	\,\bigg\vert\,
	\a_m\a_{k+1}^{-1}
	\biggr] ^{\nf1p}
	} \qquad &
	\notag \\ &
	\leq
	\E \biggl[
	\sum_{z\in 3^k\Zd \cap \cu_j}
	\biggl|
	\fint_{z+\cu_k}
	\a_m \a_{k+1}^{-1}(\exp( \g_{k+1} - \tau^2) - 1 )
	\biggr|^p
	\,\bigg\vert\,
	\a_m\a_{k+1}^{-1}
	\biggr] ^{\nf1p}
	\notag \\ &
	\leq
	Cp^{\nf12} \gmcdelta
	\biggl(
	\sum_{z\in 3^k\Zd \cap \cu_j}
	\| \a_m \a_{k+1}^{-1} \|_{\underline{L}^2(z+\cu_k)}^p
	\biggr)^{\!\nf1p}
	\leq
	C 3^{\frac dp(j-k)} p^{\nf12} \gmcdelta
	\| \a_m \a_{k+1}^{-1} \|_{\underline{L}^p(\cu_j)}
	\,.
\end{align*}
Conditionally on~$\a_m\a_{k+1}^{-1}$, the field~$\a_m \a_{k+1}^{-1}(\exp( \g_{k+1} - \tau^2) - 1 )$ has range of dependence at most~$\sqrt{d}\,3^{k+1}$. The summands indexed by~$z \in 3^k\Zd \cap \cu_j$ can therefore be partitioned into~$C(d)$ classes within each of which they are independent. Applying Lemma~\ref{l.Rosenthal} to each class, we obtain
\begin{align*}
	\lefteqn{
	\E \Biggl[ \biggl|
	\sum_{z\in 3^k\Zd\cap \cu_j}
	\fint_{z+\cu_k} \a_m \a_{k+1}^{-1}(\exp( \g_{k+1} - \tau^2) - 1 )
	\biggr|^p
	\,\Bigg\vert\,
	\a_m\a_{k+1}^{-1}
	\Biggr]^{\nf1p}
	} \qquad &
	\notag \\ &
	\leq
	Cp^{\nf12}
	\biggl(
	\sum_{z\in 3^k\Zd \cap \cu_j}
	\E \biggl[
	\biggl|
	\fint_{z+\cu_k}
	\a_m \a_{k+1}^{-1}(\exp( \g_{k+1}-\tau^2) - 1 )
	\biggr|^2
	\,\bigg\vert\,
	\a_m\a_{k+1}^{-1}
	\biggr]
	\biggr)^{\!\nf12}
	\notag \\ & \qquad
	+
	Cp\E \biggl[
	\max_{z\in 3^k\Zd \cap \cu_j}
	\biggl|
	\fint_{z+\cu_k}
	\a_m \a_{k+1}^{-1}(\exp( \g_{k+1}-\tau^2) - 1 )
	\biggr|^p
	\,\bigg\vert\,
	\a_m\a_{k+1}^{-1}
	\biggr] ^{\nf1p}
	\notag \\ &
	\leq
	C 3^{\frac d2(j-k)} \gmcdelta p^{\nf12}\| \a_m \a_{k+1}^{-1} \|_{\underline{L}^2(\cu_j)}
	+
	Cp^{\nf32} \gmcdelta
	\biggl(
	\sum_{z\in 3^k\Zd \cap \cu_j}
	\| \a_m \a_{k+1}^{-1} \|_{\underline{L}^2(z+\cu_k)}^p
	\biggr)^{\!\nf1p}
	\,.
\end{align*}
Dividing by~$3^{d(j-k)}$ gives
\begin{align*}
	\lefteqn{
	\E \Biggl[ \biggl|
	\avsum_{z\in 3^k\Zd\cap \cu_j}
	\fint_{z+\cu_k} \a_m \a_{k+1}^{-1}(\exp( \g_{k+1}-\tau^2) - 1 )
	\biggr|^p
	\,\Bigg\vert\,
	\a_m\a_{k+1}^{-1}
	\Biggr]^{\nf1p}
	} \qquad &
	\notag \\ &
	\leq
	C 3^{-\frac d2(j-k)} \gmcdelta p^{\nf12}\| \a_m \a_{k+1}^{-1} \|_{\underline{L}^2(\cu_j)}
	+
	C 3^{-d(1 - \frac 1p) (j-k)} p^{\nf32} \gmcdelta
	\| \a_m \a_{k+1}^{-1} \|_{\underline{L}^p(\cu_j)}
	\,.
\end{align*}
Separately, the triangle inequality and~\eqref{e.condiepees} give the cruder bound
\begin{align*}
	\E \Biggl[
	\biggl|
	\avsum_{z\in 3^k\Zd \cap \cu_j}
	\fint_{z+\cu_k}
	\a_m \a_{k+1}^{-1}(\exp( \g_{k+1}-\tau^2) - 1 )
	\biggr|^p
	\,\Bigg\vert\,
	\a_m\a_{k+1}^{-1}
	\Biggr] ^{\nf1p}
	\leq
	C p^{\nf12}  \gmcdelta
	\| \a_m \a_{k+1}^{-1} \|_{\underline{L}^p(\cu_j)}
	\,.
\end{align*}
Combining the previous two displays yields
\begin{align*}
	\lefteqn{
	\E \Biggl[ \biggl|
	\fint_{\cu_j}
	\a_m \a_{k+1}^{-1}(\exp( \g_{k+1}-\tau^2) - 1 )
	\biggr|^p
	\,\Bigg\vert\,
	\a_m\a_{k+1}^{-1}
	\Biggr]^{\nf1p}
	} \qquad &
	\notag \\ &
	\leq
	C \gmcdelta p^{\nf12} 3^{-\frac d2(j-k)} \| \a_m \a_{k+1}^{-1} \|_{\underline{L}^2(\cu_j)}
	+
	C \gmcdelta p^{\nf12}
	\min \bigl\{ p 3^{-d(1 - \frac 1p) (j-k)} \, , 1 \bigr\}
	\| \a_m \a_{k+1}^{-1} \|_{\underline{L}^p(\cu_j)}
	\,.
\end{align*}
Taking~$L^p(\P)$ norms and using~\eqref{e.am.Lp.moments} gives
\begin{align*}
	\lefteqn{
	\E \Biggl[ \biggl|
	\fint_{\cu_j}
	\a_m \a_{k+1}^{-1}(\exp( \g_{k+1}-\tau^2) - 1 )
	\biggr|^p
	\Biggr]^{\nf1p}
	} \qquad &
	\notag \\ &
	\leq
	C \gmcdelta p^{\nf12}
	\Bigl( 3^{-\frac d2(j-k)}
	+
	\min \bigl\{ p 3^{-d(1 - \frac 1p) (j-k)} \, , 1 \bigr\}  \Bigr)
	\E \Bigl[ \| \a_m \a_{k+1}^{-1} \|_{\underline{L}^p(\cu_j)}^p \Bigr]^{\nf1p}
	\notag \\ &
	\leq
	C \gmcdelta p^{\nf12}
	\Bigl( 3^{-\frac d2(j-k)}
	+
	\min \bigl\{ p 3^{-d(1 - \frac 1p) (j-k)} \, , 1 \bigr\}  \Bigr)
	\exp\bigl( C p \gmcdelta^2 (m-k)  \bigr)
	\,.
\end{align*}
This is~\eqref{e.condy.moment.bound}.

\smallskip

\emph{Step 2.} We show that, for every~$j\leq m$,
\begin{equation}
	\E \Biggl[ \biggl| \fint_{\cu_j}  \a_m \a_n^{-1}  - 1 \biggr|^p \Biggr]^{\nf1p}
	\leq
	C \gmcdelta p^{\nf12} (\log p) (1+m-j) \exp\bigl( C p \gmcdelta^2 (m-j)  \bigr)
	\,.
	\label{e.am.Besov.onecube}
\end{equation}
By the triangle inequality and~\eqref{e.aL.negative.Besov.iden1} we have
\begin{align*}
	&\biggl| \fint_{\cu_j}  \a_m \a_n^{-1}  - 1 \biggr|  \\
	&\quad \leq
	\sum_{k=n}^{j-1}
	\biggl|  \fint_{\cu_j} \a_m \a_{k+1}^{-1}(\exp( \g_{k+1}-\tau^2) - 1 )
	\biggr|
	+
	\sum_{k=j}^{m-1}
	\biggl|  \fint_{\cu_j} \a_m \a_{k+1}^{-1}(\exp( \g_{k+1}-\tau^2) - 1 )
	\biggr|
	\,.
\end{align*}
We bound the first term on the right side by summing~\eqref{e.condy.moment.bound} over~$k\in \{ n,\ldots,j-1\}$. The resulting sum, after factoring out~$\exp(Cp\gmcdelta^2(m-j))$, is at most~$C\log(2+p)$ by~\eqref{e.aL.Besov.p.condition}. Thus
\begin{align*}
	\E \Biggl[
	\biggl|
	\sum_{k=n}^{j-1}
	\fint_{\cu_j} \a_m \a_{k+1}^{-1}(\exp( \g_{k+1}-\tau^2) - 1 )
	\biggr|^p \Biggr]^{\nf1p}
	&
	\leq
	\sum_{k=n}^{j-1}
	\E \Biggl[
	\biggl|
	\fint_{\cu_j} \a_m \a_{k+1}^{-1}(\exp( \g_{k+1}-\tau^2) - 1 )
	\biggr|^p \Biggr]^{\nf1p}
	\notag \\ &
	\leq
	C \gmcdelta p^{\nf12} (\log p)
	\exp\bigl( C p \gmcdelta^2 (m-j)  \bigr)
	\,.
\end{align*}
For the second sum over~$k\geq j$ we use the cruder bound coming from~\eqref{e.aman.Linfty.moments} and~$\tau \leq C \gmcdelta$, to obtain
\begin{align*}
	\lefteqn{
	\E \Biggl[ \biggl( \sum_{k=j}^{m-1}
	\biggl|  \fint_{\cu_j} \a_m \a_{k+1}^{-1}(\exp( \g_{k+1}-\tau^2) - 1 )
	\biggr| \biggr)^p \Biggr]^{\nf1p}
	} \qquad &
	\notag \\ &
	\leq
	\sum_{k=j}^{m-1}
	\E \Bigl[ \|  \a_m \a_{k+1}^{-1} \|_{L^\infty(\cu_j)}^{2p}\Bigr]^{\nf1{2p}}  \E \Bigl[ \| (\exp( \g_{k+1}-\tau^2) - 1 ) \|_{L^\infty(\cu_j)}^{2p}  \Bigr]^{\nf1{2p}}
	\notag \\ &
	\leq
	\sum_{k=j}^{m-1}
	C\gmcdelta p^{\nf12}
	\exp\bigl( C p \gmcdelta^2 (m-k)  \bigr)
	\leq
	C\gmcdelta p^{\nf12}
	(m-j)
	\exp\bigl( C p \gmcdelta^2 (m-j)  \bigr)
	\,.
\end{align*}
Combining these yields~\eqref{e.am.Besov.onecube} when~$j\geq n$.

\smallskip

Suppose now that~$j<n$. By~\eqref{e.aL.negative.Besov.iden1}, the triangle inequality, independence of~$\g_{k+1}$ from~$\a_m\a_{k+1}^{-1}$, Assumption~\ref{a.g2}, and Lemma~\ref{l.am.Lp.moments},
\begin{align*}
	\E\left[
	\left|
	\fint_{\cu_j}\a_m\a_n^{-1}-1
	\right|^p
	\right]^{\nf1p}
	&\leq
	\sum_{k=n}^{m-1}
	\E\left[
	\left|
	\fint_{\cu_j}
	\a_m\a_{k+1}^{-1}
	\bigl(\exp(\g_{k+1}-\tau^2)-1\bigr)
	\right|^p
	\right]^{\nf1p}
	\\
	&\leq
	C\gmcdelta p^{\nf12}(m-n)
	\exp\bigl(Cp\gmcdelta^2(m-n)\bigr)
	\\
	&\leq
	C\gmcdelta p^{\nf12}(1+m-j)
	\exp\bigl(Cp\gmcdelta^2(m-j)\bigr)\,.
\end{align*}
This proves~\eqref{e.am.Besov.onecube} for every~$j\leq m$.

\smallskip

\emph{Step 3.} We sum over the scales to prove~\eqref{e.aL.Besov}. Returning to~\eqref{e.aL.negative.Besov.iden0} and applying~\eqref{e.am.Besov.onecube} and the triangle inequality, we obtain
\begin{align}
	\lefteqn{
	3^{-spm}
	\E \Bigl[ \| \a_m\a_n^{-1} - 1 \|_{\Besov{-s}{p}{p} (\cu_m) }^p \Bigr]
	}
	\qquad &
	\notag \\ &
	\leq
	C
	\sum_{j=-\infty}^m
	s 3^{sp(j-m)}
	\E \Biggl[
	\avsum_{z\in 3^j\Zd \cap \cu_m}
	\bigl| \bigl( \a_m\a_n^{-1} \bigr) _{z+\cu_j} - 1 \bigr|^p
	\Biggr]
	\notag \\ &
	\leq
	s \bigl( C \gmcdelta p^{\nf12} (\log p) \bigr)^p
	\sum_{j=-\infty}^m
	(1+m-j)^p
	\exp \bigl( \bigl(  (\log 3)sp - C p^2 \gmcdelta^2  \bigr)  (j-m) \bigr)
	\,.
	\label{e.aL.Besov.sumscales}
\end{align}
By~\eqref{e.aL.Besov.p.condition}, after increasing the constant there, the exponent in the last line of~\eqref{e.aL.Besov.sumscales} satisfies~$(\log 3)sp - Cp^2\gmcdelta^2\geq\frac12(\log 3)sp$, so we obtain
\begin{align*}
	3^{-sm}
	\E \Bigl[ \| \a_m\a_n^{-1} - 1 \|_{\Besov{-s}{p}{p} (\cu_m) }^p \Bigr]^{\nf1p}
	&
	\leq
	C s^{\nf1p} \gmcdelta p^{\nf12} (\log p)
	\biggl(
	\sum_{j=-\infty}^m
	(1+m-j)^p
	\exp \Bigl( \frac12 sp (j-m) \Bigr)
	\biggr)^{\!\nf1p}
	\notag \\ &
	\leq
	Cs^{-1} \gmcdelta p^{\nf12} (\log p)
	\,.
\end{align*}
This is~\eqref{e.aL.Besov}, which is equivalent to~\eqref{e.aL.Besov.invese}, as noted at the beginning of the proof.
\end{proof}

\subsection{Sensitivity estimates for coarse-grained matrices}

\begin{lemma}
\label{lemma:infrared.approx.cutoffs}
There exists a constant~$C(d) < \infty$ such that, for every~$n\in\N_0$ and~$k\in\Z$, every~$\xi \in [1,\infty)$, and every convex domain~$U \subseteq \cu_k$,
\begin{multline}
	\E \Biggl[
	\sup_{L \geq n }
	\biggl(
	\biggl |
	\frac{\a_{L,*}(U)^{-\nf12}\a_{n,*}(U)\a_{L,*}(U)^{-\nf12}}
	{(\a_n\a_L^{-1})(0)}-\Id
	\biggr|
	+
	\biggl |
	\frac{\a_L(U)^{-\nf12}\a_n(U)\a_L(U)^{-\nf12}}
	{(\a_n\a_L^{-1})(0)}-\Id
	\biggr|
	\biggr)^{\xi}
	\Biggr]^{\nf 1\xi}
	\\
	+
	\E \Biggl[
	\sup_{L \geq n }
	\biggl(
	\biggl |
	\frac{\a_{n,*}(U)^{-\nf12}\a_{L,*}(U)\a_{n,*}(U)^{-\nf12}}
	{(\a_L\a_n^{-1})(0)}-\Id
	\biggr|
	+
	\biggl |
	\frac{\a_n(U)^{-\nf12}\a_L(U)\a_n(U)^{-\nf12}}
	{(\a_L\a_n^{-1})(0)}-\Id
	\biggr|
	\biggr)^{\xi}
	\Biggr]^{\nf 1\xi}
	\\
	\leq
	C \xi^{\nf12} \gmcdelta 3^{k- n}
	\exp(C \xi \gmcdelta^2 3^{2 (k-n)} )
	\,,
	\label{e.sensitivity.large.waves}
\end{multline}
and, for every~$m\in\N_0$ with~$m>n$,
\begin{multline}
	\E \Bigl[ \bigl |
	\a_{m,*}(U)^{-\nf12}\a_{n,*}(U)\a_{m,*}(U)^{-\nf12}-\Id
	\bigr|^{\xi}
	+
	\bigl |
	\a_m(U)^{-\nf12}\a_n(U)\a_m(U)^{-\nf12}-\Id
	\bigr|^{\xi}  \Bigr]^{\nf1 \xi}
	\\
	+\E \Bigl[ \bigl |
	\a_{n,*}(U)^{-\nf12}\a_{m,*}(U)\a_{n,*}(U)^{-\nf12}-\Id
	\bigr|^{\xi}
	+
	\bigl |
	\a_n(U)^{-\nf12}\a_m(U)\a_n(U)^{-\nf12}-\Id
	\bigr|^{\xi}  \Bigr]^{\nf1 \xi}
	\\
	\leq
	\bigl( C \xi^{\nf12} \gmcdelta \bigl( (m-n)^{\nf12} +  (k- n)_+  \bigr)  \bigr)
	\exp\bigl( C \xi \gmcdelta^2 \bigl( m-n + (k - n)^2_+ \bigr) \bigr)
	\,.
	\label{e.infrared.approx.cutoffs}
\end{multline}
\end{lemma}
\begin{proof}
Apply Lemmas~\ref{l.sensitivity.for.cutoffs} and~\ref{l.sensitivity.general} on~$U\subset\cu_k$, first to~$(\a_n,(\a_n\a_L^{-1})(0)\a_L)$ and then with the two coefficients interchanged. Taking the supremum over~$L\geq n$ and using~\eqref{e.infraredhom.approx.cutoffs.large.waves.field} gives~\eqref{e.sensitivity.large.waves}. For~\eqref{e.infrared.approx.cutoffs}, apply Lemma~\ref{l.sensitivity.general} directly to the unnormalized pair~$(\a_n,\a_m)$, in both orders, and use~\eqref{e.aman.Linfty.moments}.
\end{proof}

For~$m\in\N_0$, take the expectation of~$\a_m(U)$ and the inverse of the expectation of~$\a_{m,*}^{-1}(U)$:
\begin{equation}
	\ahom_{m}(U)  \coloneqq  \E[\a_m(U)]
	\qand
	\ahom_{m,*}(U) \coloneqq \E[\a_{m,*}^{-1}(U)]^{-1} \,.
	\label{e.annealed.matrices.def}
\end{equation}
By subadditivity, the matrices~$\ahom_m(\cu_n)$ and~$\ahom_{m,*}^{-1}(\cu_n)$ are nonincreasing in~$n$ in the Loewner order. Using~\eqref{e.ord.bounds.truncated} and~\eqref{e.subadditivity.J}, the limits of~$\ahom_m(\cu_n)$ and~$\ahom_{m,*}(\cu_n)$ as~$n\to\infty$ exist. We define
\begin{equation}
	\ahom_m  \coloneqq  \lim_{n \to \infty} \ahom_{m}(\cu_n)
	\,.
	\label{e.ahom.ahomstar.defs}
\end{equation}
For each fixed~$m$, these two limits agree: the ratio~$b_n/c_n$ below is the renormalized ellipticity ratio of~$\a_m$ at scale~$3^n$ in the sense of~\citet[(2.106), p.~41]{AK.HC}, and tends to one by~\citet[Theorem~3.1, pp.~52--53]{AK.HC}. We verify the hypotheses of that theorem for the fixed-cutoff field. After a fixed spatial rescaling the field is stationary with finite range of dependence, and the supremum of~$|\log\a_m|$ on each fixed cube has a sub-Gaussian bound by Assumption~\ref{a.g2} and the finiteness of the sum defining~$\a_m$. Thus Proposition~D.3 of~\citet[pp.~145--146]{AK.HC}, with exponent~$\gamma=1/4$, supplies the coarse-grained ellipticity assumption~(P2\textdagger) of that theorem \citep[Section~2.6, p.~33]{AK.HC}. The theorem's stationarity and mixing assumptions~(P1) and~(P3) \citep[Section~1.2, pp.~7--9]{AK.HC} hold with~$\beta=0$ and~$\nu=d/2>\gamma$. All constants in this application may depend on the fixed cutoff~$m$. The same verification, for the normalized coefficient, is used in the proof of Theorem~\ref{t.fixed.cutoff.Dirichlet.algebraic} in Section~\ref{s.Holder.regularity}.

\smallskip

By the symmetries of Assumption~\ref{a.g3}, write~$\ahom_m(\cu_n)=b_n\Id$ and~$\ahom_{m,*}(\cu_n)=c_n\Id$. Theorem~3.1 of~\citet{AK.HC} gives~$b_n/c_n\to1$ after undoing the fixed spatial rescaling. Subadditivity gives decreasing~$b_n$ and increasing~$c_n$, and the constant competitors and Cauchy--Schwarz give~$\E[\a_m^{-1}(0)]^{-1}\leq c_n\leq b_n\leq\E[\a_m(0)]$. Their limits are therefore positive and equal, so~$\ahom_{m,*}(\cu_n)\to\ahom_m$. This argument concerns the annealed matrices themselves: the dual quantity throughout is~$\E[\a_{m,*}^{-1}(\cu_n)]^{-1}$, with the inverse taken after the expectation. Its admissible fields are~$q+L^2_{\mathrm{sol},0}(\cu_n)$, as in the dual variational formula. No interchange of a samplewise limit and expectation is needed.

\begin{lemma}
\label{l.annealed.matrix.bounds}
For every bounded convex domain~$U\subseteq\Rd$ and~$m,n\in\N_0$ with~$n<m$, we have
\begin{equation}
	\left\{
	\begin{aligned}
		& \ahom_{m}(U) \leq \ahom_n(U) \leq \E \Bigl [\sup_U (\a_m^{-1} \a_n) \Bigr] \ahom_{m}(U) \,,
		\\
		&
		\exp\bigl( -2\tau^2(m-n) \bigr) \ahom_{m,*}^{-1} (U) \leq \ahom_{n,*}^{-1} (U)  \leq \E \Bigl [\sup_U (\a_m \a_n^{-1}) \Bigr] \ahom_{m,*}^{-1} (U)
		\,.
	\end{aligned}
	\right.
	\label{e.annealed.matrix.bounds}
\end{equation}
In particular, for every~$m,n\in\N_0$ with~$n<m$, we have
\begin{equation}
	\ahom_m \leq \ahom_n \leq \exp( 2 \tau^2 (m-n)) \ahom_m
	\,.
	\label{e.annealed.ordering}
\end{equation}
\end{lemma}
\begin{proof}
Since~$\a_m\a_n^{-1}$ is independent of~$\a_n$, it is independent of~$v(\cdot,U,e,0;\a_n)$. Moreover, for a constant~$c$, the function~$-v(\cdot,U,e,0;\a_n)+c$ belongs to~$e\cdot x+H_0^1(U)$. Hence
\begin{align*}
	e \cdot \ahom_{m}(U) e
	&
	=
	\E \biggl[
	\min_{w \in e\cdot x + H^1_0(U) }
	\fint_U
	\nabla w \cdot \a_m \nabla w
	\biggr]
	\notag \\ &
	\leq
	\E \biggl[
	\fint_U
	\nabla v(\cdot,U,e,0;\a_n) \cdot \a_m \nabla v(\cdot,U,e,0;\a_n)
	\biggr]
	\notag \\ &
	=
	\fint_U
	\E \bigl[
	\nabla v(x,U,e,0\,;\a_n) \cdot \a_m(x) \nabla v(x,U,e,0\,;\a_n)
	\bigr]
	\,dx
	\notag \\ &
	=
	\fint_U
	\underbrace{\E \bigl[ \a_m(x) \a_n^{-1}(x) \bigr]}_{=1}
	\E \bigl[
	\nabla v(x,U,e,0\,;\a_n) \cdot \a_n(x) \nabla v(x,U,e,0\,;\a_n)
	\bigr]
	\,dx
	=
	e \cdot \ahom_n(U) e\,.
\end{align*}
For the other direction, we have
\begin{align*}
	e \cdot \ahom_{m}(U) e
	&
	=
	\E \biggl[
	\min_{w \in e\cdot x + H^1_0(U) }
	\fint_U
	\nabla w \cdot \a_m \nabla w
	\biggr]
	\notag \\ &
	\geq
	\E \biggl[
	\inf_{U} \bigl( \a_m \a_n^{-1}\bigr)
	\min_{w \in e\cdot x + H^1_0(U) }
	\fint_U
	\nabla w \cdot \a_n \nabla w
	\biggr]
	\notag \\ &
	=
	\E \Bigl[ \inf_{U} \bigl( \a_m \a_n^{-1}\bigr) \Bigr]
	\E \biggl[ \min_{w \in e\cdot x + H^1_0(U) }
	\fint_U
	\nabla w\cdot \a_n \nabla w
	\biggr]
	=
	\E \Bigl[ \inf_{U} \bigl( \a_m \a_n^{-1}\bigr) \Bigr] \bigl( e\cdot \ahom_n(U) e\bigr)\,.
\end{align*}
Since~$t\mapsto t^{-1}$ is convex on~$(0,\infty)$ and~$(\inf_U(\a_m\a_n^{-1}))^{-1}=\sup_U(\a_m^{-1}\a_n)$, Jensen's inequality gives~$\E[\inf_U(\a_m\a_n^{-1})]^{-1}\leq\E[\sup_U(\a_m^{-1}\a_n)]$; together with the previous two displays, this gives the first line of~\eqref{e.annealed.matrix.bounds}. For the second line of~\eqref{e.annealed.matrix.bounds}, use the dual variational formula with the minimizing flux for~$\a_n$:
\begin{align*}
	e \cdot \ahom_{m,*}^{-1} (U) e
	&
	=
	\E \biggl[
	\min_{\h \in e + L^2_{\sol,0}(U) }
	\fint_U
	\h\cdot \a_m^{-1} \h
	\biggr]
	\notag \\ &
	\leq
	\E \biggl[
	\fint_U
	(\a_n \nabla v(\cdot,U,0,e;\a_n)) \cdot \a_m^{-1} (\a_n \nabla v(\cdot,U,0,e;\a_n) )
	\biggr]
	\notag \\ &
	=
	\fint_U
	\E \bigl[
	(\a_n(x) \nabla v(x,U,0,e\,;\a_n)) \cdot \a_m^{-1}(x) (\a_n(x) \nabla v(x,U,0,e\,;\a_n) )
	\bigr]
	\,dx
	\notag \\ &
	=
	\fint_U
	\underbrace{\E \bigl[ \a_m^{-1}(x) \a_n(x) \bigr]}_{= \exp(2\tau^2(m-n))}
	\E \bigl[
	\nabla v(x,U,0,e\,;\a_n) \cdot \a_n(x) \nabla v(x,U,0,e\,;\a_n)
	\bigr]
	\,dx
	\notag \\ &
	=
	\exp(2\tau^2(m-n))
	\bigl( e \cdot \ahom_{n,*}^{-1} (U) e\bigr)
	\,,
\end{align*}
and
\begin{align*}
	e \cdot \ahom_{m,*}^{-1} (U) e
	&
	=
	\E \biggl[
	\min_{\h \in e + L^2_{\sol,0}(U) }
	\fint_U
	\h \cdot \a_m^{-1} \h
	\biggr]
	\notag \\ &
	\geq
	\E \biggl[
	\inf_{U}
	(\a_m^{-1} \a_n)
	\min_{\h \in e + L^2_{\sol,0}(U) }
	\fint_U
	\h \cdot \a_n^{-1} \h
	\biggr]
	\notag \\ &
	=
	\E \Bigl[
	\inf_{U}
	(\a_m^{-1} \a_n)
	\Bigr]
	\E \biggl[
	\min_{\h \in e + L^2_{\sol,0}(U) }
	\fint_U
	\h \cdot \a_n^{-1} \h
	\biggr]
	=
	\E \Bigl[
	\inf_{U}
	(\a_m^{-1} \a_n)
	\Bigr]
	\bigl( e \cdot \ahom_{n,*}^{-1}(U) e \bigr)
	\,.
\end{align*}
In the same way,~$\E[\inf_U(\a_m^{-1}\a_n)]^{-1}\leq\E[\sup_U(\a_m\a_n^{-1})]$, and the last two displays give the second line of~\eqref{e.annealed.matrix.bounds}.
\end{proof}

By~\eqref{e.annealed.ordering} and~\eqref{e.tau.bound}, there exists~$C(d)<\infty$ such that, for every~$m,n\in\N_0$ with~$n\leq m$,
\begin{equation}
	\bigl| \ahom_{m}\ahom^{-1}_{n} - \Id \bigr|
	+
	\bigl| \ahom_{m}^{-1} \ahom_{n} - \Id \bigr|
	\leq
	C \gmcdelta^2 (m-n) \exp(C \gmcdelta^2 (m-n))
	\,.
	\label{e.annealed.ordering.ratio}
\end{equation}
In dimension two, planar duality determines~$\ahom_m$ exactly.

\begin{proposition}
\label{p.special.two.d.exact.formula}
In dimension~$d=2$, for every~$m\in\N_0$, we have~$\ahom_m=\exp(-\tau^2(m+1))$.
\end{proposition}
\begin{proof}
Let~$R$ be the matrix with rows~$(0,1)$ and~$(-1,0)$. On the simply connected square~$\cu_n$, the divergence-free fields with zero normal component on the boundary are exactly the rotated gradients of functions vanishing on the boundary:~$\Lsolo(\cu_n)=R\nabla H_0^1(\cu_n)$. Substituting this into the dual minimum in~\eqref{e.variational.astar} therefore gives, for every scalar coefficient~$a$ bounded above and away from zero on~$\cu_n$, the identity~$a_*^{-1}(\cu_n;a)=R\,a(\cu_n;a^{-1})R^{\mathrm T}$. By the invariance of the law of~$\g_0$ under~$\g_0\mapsto-\g_0$ in~\ref{a.g3},~$\a_m^{-1}$ has the same law as~$e^{2(m+1)\tau^2}\a_m$. Taking expectations in the preceding identity and using homogeneity and the invariance of the law under signed permutations, we obtain~$\ahom_{m,*}^{-1}(\cu_n)=e^{2(m+1)\tau^2}\ahom_m(\cu_n)$. Let~$n\to\infty$ and use that~$\ahom_{m,*}(\cu_n)\to\ahom_m$, as shown after~\eqref{e.ahom.ahomstar.defs}. The positive scalar~$\ahom_m$ satisfies~$\ahom_m^{-1}=e^{2(m+1)\tau^2}\ahom_m$, hence~$\ahom_m=e^{-(m+1)\tau^2}$. This argument is a finite-volume variational form of planar duality~\citep{Dykhne}.
\end{proof}

\begin{lemma}
\label{lemma.crude.J.bound}
There exists a constant~$C(d)<\infty$ such that, for every~$m \in \N_0$, every~$\xi \in [1, \infty)$, and every convex domain~$U \subseteq \cu_{m}$,
\begin{equation}
	\E \biggl[ \max_{|e| = 1} J(U, \ahom^{-\nf12}_{m} e, \ahom_{m}^{\nf12} e \, ; \a_{m})^{\xi} \biggr]^{\nf1\xi} \leq
	C \xi (m+1)  \gmcdelta^2
	\exp(C \xi (m+1) \gmcdelta^2)
	\,.
	\label{e.crude.J.upperbound}
\end{equation}
\end{lemma}
\begin{proof}
By~\eqref{e.coarse.grained.matrices} and~\eqref{e.ord.bounds.truncated},
\begin{align*}
	\max_{|e|=1}J(U,\ahom_m^{-\nf12}e,\ahom_m^{\nf12}e;\a_m)
	&
	\leq\frac12\fint_U\Bigl(\frac{\a_m}{\ahom_m}+\frac{\ahom_m}{\a_m}-2\Bigr)
	=\frac12\fint_U\Bigl(\frac{\a_m}{\ahom_m}-1\Bigr)\Bigl(1-\frac{\ahom_m}{\a_m}\Bigr)
	\notag \\ &
	\leq\frac14\fint_U\Bigl(\Bigl(\frac{\a_m}{\ahom_m}-1\Bigr)^2+\Bigl(\frac{\ahom_m}{\a_m}-1\Bigr)^2\Bigr)\,.
\end{align*}
Taking expectations in~\eqref{e.ord.bounds.truncated} also gives~$e^{-2(m+1)\tau^2}\leq\ahom_m\leq1$. Consequently, Lemma~\ref{l.am.Lp.moments}, with~$n=-1$ and exponent~$2\xi$, and~\eqref{e.tau.bound} bound the~$L^\xi(\P)$ norm of the right side of the last display by the right side of~\eqref{e.crude.J.upperbound}.
\end{proof}

\section{Inductive coarse-graining up the scales}
\label{s.self.similar}

In this section we prove the coarse-graining estimate at the cutoff scale, stated below in Proposition~\ref{p.coarse.grained.bound}. For a bounded convex domain~$U$ and~$m\in\N_0$, abbreviate
\begin{equation}
	\mathcal J_m(U)
	\coloneqq
	\max_{|e|=1}J(U,\ahom_m^{-\nf12}e,\ahom_m^{\nf12}e;\a_m)\,.
	\label{e.normalized.J.maximum.definition}
\end{equation}
On the cube~$\cu_m$, whose size matches the cutoff scale, the quantity~$\mathcal J_m(\cu_m)$ is bounded in~$L^\xi(\P)$ by~$C\xi\log(2+\xi)\gmcdelta^2$ whenever~$\xi\leq C^{-1}\gmcdelta^{-2}|\log\gmcdelta|^{-1}$. In particular, the bound does not deteriorate as further scales are added.

\smallskip

The quantity~$\mathcal J_m(\cu_m)$ is the gap between the two variational principles that define the coarse-grained matrices.  It measures the error in replacing the coefficient on a cube by a constant matrix: by Lemma~\ref{l.sharp.compare.J}, it controls~$\a_m(\cu_m)-\a_{m,*}(\cu_m)$ and the distance of these matrices from~$\ahom_m$, and it vanishes if and only if~$\a_m(\cu_m)=\a_{m,*}(\cu_m)=\ahom_m\Id$.  We argue by induction up the scales: the inductive hypothesis of Definition~\ref{d.mathcalS.def} bounds~$\mathcal J_m(\cu_m)$ in~$L^\xi(\P)$ for every~$m$ up to some~$m_0$, and the induction step extends this bound to the next scale.  The field added at each step has a scale larger than the cubes already treated, so, by Section~\ref{s.infrared.cutoffs}, it acts on them almost as a constant factor; the gain comes from estimating its effect on the larger cube in negative norms.

\begin{proposition}[Coarse-graining estimate at the cutoff scale]
\label{p.coarse.grained.bound}
There are constants~$c(d)>0$ and~$C(d)<\infty$ such that, for every~$\xi\in[1,\infty)$, if the parameter~$\gmcdelta$ in Assumption~\ref{a.g2} satisfies
\begin{equation}
	\xi  \leq C^{-1} \gmcdelta^{-2} \left| \log \gmcdelta\right|^{-1} \,,
	\label{e.xi.upper.bound.gmc}
\end{equation}
then, for every~$m\in\N_0$, we have
\begin{equation}
	\E \biggl[ \max_{|e| = 1}
	\bigl(
	J(\cu_m,\ahom_m^{-\nf12} e , \ahom_m^{\nf12} e\,;\a_m)
	\bigr) ^{\xi}
	\biggr]^{\nf1\xi}
	\leq
	C\xi \log(2+\xi) \gmcdelta^2
	\,.
	\label{e.homogenize.at.scale.m}
\end{equation}
Moreover, for every~$s \in (0,1]$ satisfying
\begin{equation*}
	4ds^{-1} \leq \xi\,,
	\qquad
	C\xi \log(2+\xi) \gmcdelta^2 \leq cs
	\quad\text{and}\quad
	s\log(2+\xi)\geq c\,,
\end{equation*}
we have
\begin{equation}
	\E \bigl[ \mathcal{E}_{s,\infty}(\cu_m;\a_m,\ahom_m)^{\xi} \bigr]^{\nf1\xi}
	\leq
	Cs^{-1} \bigl( \xi \log(2+\xi) \gmcdelta^2 \bigr)^{\nf12}
	\,.
	\label{e.homogenize.at.scale.m.E}
\end{equation}
\end{proposition}

After reducing the constant in~\eqref{e.xi.upper.bound.gmc}, apply~\eqref{e.homogenize.at.scale.m} with moment~$2\xi$ and use Lemma~\ref{l.sharp.compare.J} to obtain
\begin{equation*}
	\E \bigl[ \bigl| \ahom_m^{-1} \a_m(\cu_m) - \Id \bigr|^{2\xi} \bigr] ^{\nf 1\xi}
	+
	\E \bigl[ \bigl| \ahom_m \a_{m,*}^{-1} (\cu_m) - \Id \bigr|^{2\xi} \bigr]^{\nf 1\xi}
	\leq
	C \xi \log(2+\xi) \gmcdelta^2
	\,.
\end{equation*}
The dependence on~$\gmcdelta$ is optimal: the scale-$m$ field alone produces fluctuations of order~$\gmcdelta$ in the coarse-grained matrices.

\smallskip

\begin{definition}[Inductive hypothesis {$\mathcal{S}(m_0,\xi,\delta_1)$}]
\label{d.mathcalS.def}
Given~$m_0\in\N_0$, an exponent~$\xi\in[1,\infty)$, and~$\delta_1\in(0,1)$, let~$\S(m_0,\xi,\delta_1)$ denote the statement
\begin{equation}
	\max_{m \in \{0,\ldots,m_0\}}
	\left\|
	\max_{|e|=1}
	J(\cu_m,\ahom_m^{-\nf12}e,\ahom_m^{\nf12}e;\a_m)
	\right\|_{L^\xi(\P)}
	\leq \delta_1
	\,.
	\label{e.induction.hypothesis.h=0}
\end{equation}
Let~$\mathcal S(\infty,\xi,\delta_1)$ denote the statement
\begin{equation*}
	\sup_{m\in\N_0}
	\left\|
	\max_{|e|=1}
	J(\cu_m,\ahom_m^{-\nf12}e,\ahom_m^{\nf12}e;\a_m)
	\right\|_{L^\xi(\P)}
	\leq \delta_1 \,.
\end{equation*}
\end{definition}

The bound~\eqref{e.homogenize.at.scale.m} is the assertion that~$\mathcal S(\infty,\xi,C\xi\log(2+\xi)\gmcdelta^2)$ holds throughout the stated range of~$\xi$. The bound~\eqref{e.homogenize.at.scale.m.E} on~$\mathcal{E}_{s,\infty}$ then follows from Lemma~\ref{l.ellipticity.bound} with~$\delta_1 = C\xi\log(2+\xi)\gmcdelta^2$.

\smallskip

Proposition~\ref{p.coarse.grained.bound} is used in every later section: in Section~\ref{s.homogenized.diffusivities} to derive the decay of~$\ahom_m$, in Section~\ref{s.Holder.regularity} to show that the scales with a large homogenization error are rare, and in Sections~\ref{sec:mfd-uniqueness}--\ref{lim:sec-process-measure} to bound the coarse-grained ellipticity of~$\a_L$. The proof has three parts and uses two inductions, the first of which supplies the low-moment estimate used in the second.
\begin{itemize}
\item In Section~\ref{ss.multiscale.estimate} we use~$\mathcal S(m_0,\xi,\delta_1)$ to bound the coarse-graining error~$\mathcal E_{s,\infty}$ of~$\a_L$ on cubes below the cutoff scale (Lemma~\ref{l.ellipticity.bound}) and above it (Lemma~\ref{l.multiscale.response.large.cubes}).
\item In Section~\ref{ss.homogenization.step} we prove the contraction of Proposition~\ref{p.homogenization.step}: under~$\mathcal S(m_0,\xi,\delta_1)$, the quantity~$\mathcal J_L(\cu_m)$ with~$m\geq L+C\log(2+\xi)$ is at most~$\frac{1}{4}\delta_1$ in~$L^\xi(\P)$.
\item At the end of Section~\ref{ss.homogenization.step} we add the omitted fields~$\g_{L+1},\ldots,\g_m$ using Lemma~\ref{l.J.sensitivity} and close the inductions.
\end{itemize}

\subsection{A multiscale estimate}
\label{ss.multiscale.estimate}

The statement~$\mathcal S(m_0,\xi,\delta_1)$ bounds, for every~$m\leq m_0$, the~$L^\xi(\P)$ norm of the maximum in~\eqref{e.induction.hypothesis.h=0}. We now deduce multiscale bounds from this estimate.

\begin{definition}[{The coefficient~$\b_{L,m}$}]
\label{d.bLm}
For~$L,m\in\N_0$, define
\begin{equation*}
	\b_{L,m}
	\coloneqq
	\ahom_{m\wedge L}\frac{\a_L}{\a_{m\wedge L}}\,.
\end{equation*}
Thus, if~$m\leq L$, then~$\b_{L,m}=\ahom_m\exp\bigl(\sum_{j=m+1}^L(\g_j-\tau^2)\bigr)$, whereas~$\b_{L,m}=\ahom_L$ if~$m\geq L$.
\end{definition}

\begin{lemma}
\label{l.ellipticity.bound}
Let~$m_0 \in \N_0$,~$\xi \in [6,\infty)$ and~$\delta_1 \in [\gmcdelta^2, 1)$, and assume that~$\mathcal{S}(m_0, \xi, \delta_1)$ holds. There exist constants~$c(d)\in (0,1]$ and~$C(d) \in [1, \infty)$ such that, for every~$s\in (0,1]$ satisfying
\begin{equation}
	4d s^{-1} \leq \xi \leq  c s \gmcdelta^{-2}  \delta_1
	\,,
	\label{e.gammazeta.cond}
\end{equation}
and every~$m \in \N_0$ and~$n\in\Z \cap (-\infty,m_0]$ with~$m \geq n$, we have
\begin{align}
	\E \biggl[ \sup_{L\geq m} \maxcube{z}{n}{m}  \bigl(  \mathcal{E}_{s,  \infty} (z+\cu_n; \a_L,  (\b_{L,m})_{\cu_m} ) \bigr) ^{\xi} \biggr]^{\nf1{\xi}}
	\leq
	C s^{-1} \delta_1^{\nf 12} 3^{(\frac{d}{\xi}+\frac{s}{2})(m-n)}\,,
	\label{e.mathcal.E.Lms.moments}
\end{align}
and, for every~$L \in \N_0$ with~$L\leq m_0$,
\begin{equation}
	\E \bigl[   \mathcal{E}_{s,  \infty} (\cu_L; \a_L, \ahom_L )^{\xi}   \bigr]^{\nf1{\xi}}
	\leq
	C s^{-1} \delta_1^{\nf 12}
	\,.
	\label{e.mathcal.E.Lms.moments.small.L}
\end{equation}
\end{lemma}
\begin{proof}

Let~$s\in (0,1]$ satisfy~\eqref{e.gammazeta.cond}. By Definition~\ref{d.mathcal.E}, the coarse-graining error in~\eqref{e.mathcal.E.Lms.moments} is controlled by the values of~$J$ on the subcubes of~$z+\cu_n$ of every size~$3^k$ with~$k\leq n$. The proof proceeds in three steps. Step~1 bounds the moments of these values, Step~2 sums these bounds over the scales, and Step~3 specializes to~$m=n=L$. The heart of the matter is the case~$k\geq0$ of Step~1, where Lemma~\ref{l.J.sensitivity} transfers the induction hypothesis from the coefficient~$\a_k$ to~$\a_L$; the case~$k<0$ is estimated directly from the coefficient.

\smallskip

\emph{Step 1.} We first treat the scales at which the induction hypothesis applies. We show the existence of~$C(d)<\infty$ such that for every~$m,k \in \N_0$ with~$k\leq m\wedge m_0$ and every~$z \in 3^k \Zd \cap \cu_m$,
\begin{align}
	\E\biggl[   \sup_{L \geq m}\max_{|e| = 1}
	\Bigl( J \bigl(z + \cu_{k}, (\b_{L,m})_{\cu_m} ^{-\nf12}   e, (\b_{L,m})_{\cu_m} ^{\nf12}  e \,;\a_L \bigr) \Bigr)^{\xi} \biggr]^{\nf1\xi}
	\leq
	C \delta_1 3^{s(m-k)}
	\,.
	\label{e.moment.bound.for.J}
\end{align}
We apply Lemma~\ref{l.J.sensitivity} with~$\delta=1$ and~$\lambda=(\a_L\a_m^{-1})_{\cu_m} \ahom_m \ahom_k^{-1}$, and with~$\b =\a_L$ and~$\a=  \a_k$. We obtain, for every~$L,m,k\in\N_0$ and~$z \in 3^k \Zd \cap \cu_m$ with~$L\geq m \geq k$, and~$e\in \Rd$ with~$|e|=1$,
\begin{align*}
	\lefteqn{
	J(z+\cu_k,(\b_{L,m})_{\cu_m} ^{-\nf12} e,(\b_{L,m})_{\cu_m} ^{\nf12}e\,;\a_L)
	} \qquad &
	\notag \\ &
	\leq
	2 J(z+\cu_k,\ahom_k^{-\nf12}  e, \ahom_k^{\nf12} e \,; \a_k)
	\notag \\ & \qquad
	+
	3 \Bigl(
	\sup_{L\geq m} \sup_{x \in \cu_m,y \in z+\cu_k}
	\bigl(
	\bigl|\ahom_m \ahom_k^{-1}  (\a_L\a_m^{-1})(x)  (\a_L^{-1} \a_k)(y) - 1 \bigr|
	\notag \\ & \qquad \qquad \qquad \qquad
	+\bigl|\ahom^{-1}_m \ahom_k  (\a^{-1}_L\a_m)(x)  (\a_L \a^{-1}_k)(y) - 1 \bigr|
	\bigr)
	\Bigr)^{2}
	\bigl( J(z+\cu_k,\ahom_k^{-\nf12}  e, \ahom_k^{\nf12} e \,; \a_k)  + 1 \bigr)
	\,.
\end{align*}
We first bound the moments of the supremum in this display:
\begin{multline}
	\E\Biggl[
	\Bigl(
	\sup_{L\geq m} \sup_{x \in \cu_m,y \in z+\cu_k}
	\bigl(
	\bigl|\ahom_m \ahom_k^{-1}  (\a_L\a_m^{-1})(x)  (\a_L^{-1} \a_k)(y) - 1 \bigr|
	\\
	+\bigl|\ahom^{-1}_m \ahom_k  (\a^{-1}_L\a_m)(x)  (\a_L \a^{-1}_k)(y) - 1 \bigr|
	\bigr)
	\Bigr)^{2 \xi}
	\Biggr]^{\nf 1 \xi}
	\\
	\leq C \delta_1 3^{s(m-k)} \,.
	\label{e.XmLk.bound}
\end{multline}
Given~\eqref{e.XmLk.bound}, the induction hypothesis~$\mathcal{S}(m_0, \xi, \delta_1)$ and the independence of this supremum from the field~$\a_k$ give~\eqref{e.moment.bound.for.J}.

\smallskip

To prove~\eqref{e.XmLk.bound}, we write
\begin{equation*}
	\bigl|\ahom_m \ahom_k^{-1}  (\a_L\a_m^{-1})(x)  (\a_L^{-1} \a_k)(y) - 1 \bigr|
	=
	\biggl|\ahom_m \ahom_k^{-1}  \frac{(\a_L\a_m^{-1})(x)}{(\a_L\a_m^{-1})(y)}  (\a_m^{-1} \a_k)(y) - 1 \biggr| \,,
\end{equation*}
and use the general bound
\begin{multline*}
	|uvw-1| \leq
	|u-1|+|v-1|+|w-1|+|(u-1)(v-1)| \\
	+|(u-1)(w-1)|+|(v-1)(w-1)|+|(u-1)(v-1)(w-1)| \,,
\end{multline*}
	together with the independence of the fields~$\a_L \a^{-1}_m$ and~$\a^{-1}_m \a_k$ and the following three inequalities, invoked at moment~$2\xi$, which follow from~\eqref{e.annealed.ordering.ratio},~\eqref{e.infraredhom.approx.cutoffs.large.waves.field},~\eqref{e.aman.Linfty.moments} and~\eqref{e.gammazeta.cond}. For the deterministic factor,
\begin{equation*}
	|\ahom_m\ahom_k^{-1}-1|
	\leq C\gmcdelta^2(m-k)\exp(C\gmcdelta^2(m-k))
	\leq C\delta_1^{\nf12}3^{\frac{s}{6}(m-k)}\,,
\end{equation*}
for the ratio of the values of~$\a_L\a_m^{-1}$ at two points,
\begin{equation*}
	\E\biggl[
	\biggl(\sup_{L\geq m}\sup_{x,y\in\cu_m}
	\frac{(\a_L\a_m^{-1})(x)}{(\a_L\a_m^{-1})(y)}-1\biggr)^{2\xi}
	\biggr]^{\nf1\xi}
	\leq C\xi\gmcdelta^2\exp(C\xi\gmcdelta^2(m-k))
	\leq C\delta_1 3^{\frac{s}{3}(m-k)}\,,
\end{equation*}
and for the factor~$\a_m^{-1}\a_k$,
\begin{equation*}
	\E\left[
	\|\a_m^{-1}\a_k-1\|_{L^\infty(\cu_k)}^{2\xi}
	\right]^{\nf1\xi}
	\leq C\xi\gmcdelta^2(m-k)\exp(C\xi\gmcdelta^2(m-k))
	\leq C\delta_1 3^{\frac{s}{3}(m-k)}\,.
\end{equation*}
Each of the three deviations thus has~$L^{2\xi}$ norm at most~$C\delta_1^{\nf12}3^{s(m-k)/6}$. The two random factors depend on disjoint sets of fields, and the remaining factor is deterministic, so the~$L^{2\xi}$ norm of a product of deviations is the product of their~$L^{2\xi}$ norms. Since~$\delta_1<1$, even the product of all three is bounded by~$C\delta_1^{\nf12}3^{s(m-k)/2}$. The bound on~$|XYZ-1|$ above therefore proves~\eqref{e.XmLk.bound}, and hence~\eqref{e.moment.bound.for.J}.

\smallskip

We next estimate the terms in Definition~\ref{d.mathcal.E} with cube size smaller than one.  Equations~\eqref{e.coarse.grained.matrices} and~\eqref{e.ord.bounds.truncated} imply that, for every positive scalar coefficient~$\a$, positive scalar~$b$, bounded convex domain~$V$, and~$|e|=1$,
\begin{equation}
	J(V,b^{-\nf12}e,b^{\nf12}e\,;\a)
	\leq
	\frac12\fint_V\left(\frac{\a}{b}+\frac{b}{\a}-2\right)\,.
	\label{e.direct.J.bound.subunit}
\end{equation}
Indeed,~\eqref{e.ord.bounds.truncated} bounds the two quadratic terms in~\eqref{e.coarse.grained.matrices} by the corresponding averages of~$\a$ and~$\a^{-1}$. To estimate the right side with~$\a=\a_L$ and~$b=(\b_{L,m})_{\cu_m}$, we factor
\begin{equation*}
	\frac{\a_L}{(\b_{L,m})_{\cu_m}}
	=
	\frac{\a_m}{\ahom_m}
	\cdot
	\frac{\a_L\a_m^{-1}}
	{(\a_L\a_m^{-1})_{\cu_m}}\,.
\end{equation*}
The elementary inequality~$r+r^{-1}-2\leq(r-1)^2+(r^{-1}-1)^2$, valid for every~$r>0$, allows us to estimate the two factors separately.  We first estimate~$\a_m/\ahom_m$. For each unit cube~$z+\cu_0$ meeting~$\cu_m$, independence and Assumptions~\ref{a.g1}--\ref{a.g2} give, for every~$p\geq1$,
\begin{equation*}
	\left\|\sum_{i=0}^m\g_i(z)\right\|_{L^p(\P)}
	\leq C\sqrt{p}\,\gmcdelta\sqrt{m+1}
	\,,
\end{equation*}
and, for the oscillations of the fields on the same cube,
\begin{equation*}
	\left\|\sum_{i=0}^m
	\osc_{z+\cu_0}\g_i\right\|_{L^p(\P)}
	\leq C\sqrt{p}\,\gmcdelta\,.
\end{equation*}
There are at most~$C3^{dm}$ such unit cubes.  We apply~\eqref{e.moments.general.lognormal} on each cube and then take a union bound over these cubes.  The variational bounds~\eqref{e.ord.bounds.truncated}, together with~$\E[\a_0(0)]=1$ and~$\E[\a_0(0)^{-1}]=e^{2\tau^2}$, give~$e^{-2\tau^2}\leq\ahom_0\leq1$.  Combining this with~\eqref{e.annealed.ordering.ratio} controls the deterministic factor~$\ahom_m$.  This yields
\begin{equation*}
	\E\!\left[
	\left\|
	\frac{\a_m}{\ahom_m}
	+
	\frac{\ahom_m}{\a_m}
	-2
	\right\|_{L^\infty(\cu_m)}^\xi
	\right]^{\nf1\xi}
	\leq
	C\xi\gmcdelta^2(m+1)
	\exp\!\left(
	\frac{d(\log3)m+C}{\xi}
	+C\xi\gmcdelta^2(m+1)
	+C\gmcdelta^2m
	\right)\,.
\end{equation*}
For the second factor, the average in the denominator lies between the minimum and the maximum of~$\a_L\a_m^{-1}$ on~$\cu_m$.  Hence~\eqref{e.infraredhom.approx.cutoffs.large.waves.field}, applied with~$n=k=m$, gives the same bound with~$m+1$ replaced by a dimensional constant.  The two factors are independent.  H\"older's inequality and~$|uv-1|\leq|u-1||v|+|v-1|$, applied also to the inverses, therefore give
\begin{align*}
	\lefteqn{
	\E\!\left[
	\sup_{L\geq m}
	\left\|
	\frac{\a_L}{(\b_{L,m})_{\cu_m}}
	+
	\frac{(\b_{L,m})_{\cu_m}}{\a_L}
	-2
	\right\|_{L^\infty(\cu_m)}^\xi
	\right]^{\nf1\xi}
	}\qquad &
	\notag \\ &
	\leq
	C\xi\gmcdelta^2(m+1)
	\exp\!\left(
	\frac{d(\log3)m+C}{\xi}
	+C\xi\gmcdelta^2(m+1)
	+C\gmcdelta^2m
	\right)
	\leq C\delta_1 3^{sm}\,.
\end{align*}
The last inequality follows from~\eqref{e.gammazeta.cond}. Thus, by~\eqref{e.direct.J.bound.subunit}, for every~$k<0$,
\begin{equation}
	\E\!\left[
	\sup_{L\geq m}\max_{z\in3^k\Zd\cap\cu_m}\max_{|e|=1}
	J\bigl(z+\cu_k,(\b_{L,m})_{\cu_m}^{-\nf12}e,(\b_{L,m})_{\cu_m}^{\nf12}e\,;\a_L\bigr)^\xi
	\right]^{\nf1\xi}
	\leq C\delta_1 3^{sm}\,.
	\label{e.moment.bound.for.J.subunit}
\end{equation}

\smallskip

\emph{Step 2.} We prove~\eqref{e.mathcal.E.Lms.moments}. Split the sum in Definition~\ref{d.mathcal.E} into~$k\geq0$ and~$k<0$. Minkowski's inequality,~\eqref{e.moment.bound.for.J},~\eqref{e.moment.bound.for.J.subunit}, and the union bound over the~$O(3^{d(m-k)})$ cubes of scale~$k$ give
\begin{align*}
	\left\|\sup_{L\geq m} \maxcube{z}{n}{m}\!
	\mathcal{E}_{s,\infty}
	\bigl(z+\cu_n;\a_L,(\b_{L,m})_{\cu_m}\bigr)
	\right\|_{L^\xi(\P)}
	&\leq
	C \delta_1^{\nf 12}
	3^{\frac{d}{\xi}(m-n)}
	\sum_{k=0}^n 3^{-(\frac s2-\frac d\xi)(n-k)+\frac s2(m-n)}
	\indc_{\{n\geq0\}}
	\notag \\ & \qquad
	+
	C\delta_1^{\nf12}3^{\frac s2m}
	\sum_{k=-\infty}^{n\wedge(-1)}3^{-s(n-k)}
	\\
	&\leq
	Cs^{-1}\delta_1^{\nf12}3^{(\frac{d}{\xi}+\frac{s}{2})(m-n)}
	\,.
\end{align*}
Here we used~$\xi\geq4ds^{-1}$ from~\eqref{e.gammazeta.cond}.

\smallskip

\emph{Step 3.} We prove~\eqref{e.mathcal.E.Lms.moments.small.L}. Specializing~\eqref{e.mathcal.E.Lms.moments} to~$m=n=L$ and using~$\b_{L,L}=\ahom_L$ gives~\eqref{e.mathcal.E.Lms.moments.small.L}.
\end{proof}

\begin{lemma}[Coarse-graining error on cubes above the cutoff scale]
\label{l.multiscale.response.large.cubes}
There exist~$c(d)>0$ and~$C(d)<\infty$ such that the following holds. Let~$L,m_0\in\N_0$ satisfy~$L\leq m_0$, let~$s\in(0,1]$ and~$\delta_1\in[\gmcdelta^2,1)$, and let~$\xi$ satisfy~$4ds^{-1}\leq\xi\leq cs\gmcdelta^{-2}\delta_1$. If~$\mathcal S(m_0,\xi,\delta_1)$ is valid, then the following estimate holds for every~$M\geq L$, every~$w\in\Rd$, and~$r\in\{1,2\}$:
\begin{equation}
	\bigl\|
	\mathcal E_{s,\infty,r}(w+\cu_M;\a_L,\ahom_L)
	\bigr\|_{L^\xi(\P)}
	\leq
	Cs^{-\nf{1}{r}}\delta_1^{\nf12}\,.
	\label{e.multiscale.response.large.cubes}
\end{equation}
\end{lemma}
\begin{proof}
It is enough, by stationarity, to take~$w=0$.  Recall the notation~$\mathcal J_L(\cdot)$ of~\eqref{e.normalized.J.maximum.definition}. Definition~\ref{d.mathcal.E} and Minkowski's inequality give
\begin{align}
	\bigl\|\mathcal E_{s,\infty,r}(\cu_M;\a_L,\ahom_L)\bigr\|_{L^\xi}^r
	\leq
	C\sum_{k=-\infty}^M
	3^{-sr(M-k)}
	\left\|
	\max_{z\in3^k\Zd\cap\cu_M}\mathcal J_L(z+\cu_k)^{\nf12}
	\right\|_{L^\xi}^r\,.
	\label{e.multiscale.response.large.cubes.sum}
\end{align}
We split the sum into three ranges.

\smallskip

If~$L\leq k\leq M$, then subadditivity~\eqref{e.subadditivity.J} on the partition into cubes of scale~$L$, the induction hypothesis, and Jensen's inequality yield~$\|\mathcal J_L(z+\cu_k)\|_{L^\xi}\leq C\delta_1$. Taking the maximum over at most~$C3^{d(M-k)}$ cubes gives a contribution to~\eqref{e.multiscale.response.large.cubes.sum} bounded by
\begin{equation*}
	C\delta_1^{\nf{r}{2}}
	\sum_{j=0}^{M-L}3^{-r(s-\frac{d}{\xi})j}
	\leq Cs^{-1}\delta_1^{\nf{r}{2}}\,.
\end{equation*}

\smallskip

If~$0\leq k<L$, then~\eqref{e.moment.bound.for.J} with~$m=L$ and cutoff~$L$ gives~$\|\mathcal J_L(z+\cu_k)\|_{L^\xi}\leq C\delta_1 3^{s(L-k)}$. After taking the maximum over the cubes of scale~$k$ in~$\cu_M$, we bound the corresponding summand by
\begin{equation*}
	C\delta_1^{\nf{r}{2}}
	3^{-r(s-\frac{d}{\xi})(M-L)}
	3^{-r(\frac{s}{2}-\frac{d}{\xi})(L-k)}\,.
\end{equation*}
The sum of these bounds over~$k$ is at most~$Cs^{-1}\delta_1^{\nf{r}{2}}$.

\smallskip

Finally, if~$k<0$, apply~\eqref{e.moment.bound.for.J.subunit} with~$m=L$ on each scale-$L$ parent cube.  The maximum over the at most~$C3^{d(M-L)}$ such cubes gives
\begin{align*}
	&\sum_{k=-\infty}^{-1}
	3^{-sr(M-k)}
	\left\|
	\max_{z\in3^k\Zd\cap\cu_M}\mathcal J_L(z+\cu_k)^{\nf12}
	\right\|_{L^\xi}^r
	\\
	&\qquad\leq
	C\delta_1^{\nf{r}{2}}
	3^{-r(s-\frac{d}{\xi})(M-L)}
	3^{-r(\frac{s}{2}-\frac{d}{\xi})L}
	\sum_{j=1}^{\infty}3^{-srj}
	\leq Cs^{-1}\delta_1^{\nf{r}{2}}\,.
\end{align*}
Inserting the estimates into~\eqref{e.multiscale.response.large.cubes.sum} and taking the~$r$th root proves the lemma.
\end{proof}

\subsection{The contraction step}
\label{ss.homogenization.step}

The coefficient fields in Lemmas~\ref{l.J.bound.by.Besov},~\ref{l.Besov.norms} and~\ref{l.var.bounds} are understood in the class used to define the coarse-grained matrices in Section~\ref{s.coarse.graining}: they are measurable and symmetric, and on every bounded cube in their domain there are constants~$0<\lambda\leq\Lambda<\infty$ such that~$\lambda\Id\leq\a\leq\Lambda\Id$ almost everywhere. These constants may depend on the cube. For a random field this condition holds almost surely, and the field and its coarse-grained matrices are measurable random variables. Variances and nonnegative expectations are allowed to be infinite, with the variance convention stated before Lemma~\ref{l.var.bounds}; no additional square-integrability assumption is imposed.

\smallskip

We now prove the contraction step: if~$\mathcal J_m(\cu_m)$ is at most~$\delta_1$ in~$L^\xi(\P)$ for every scale~$m$ up to~$m_0$, then, for~$L\leq m_0$ and a cube~$\cu_m$ of scale~$m\geq L+C\log(2+\xi)$, the quantity~$\mathcal J_L(\cu_m)$, which measures~$\a_L$ against~$\ahom_L$, is at most~$\frac{1}{4}\delta_1$ in~$L^\xi(\P)$.  This is Proposition~\ref{p.homogenization.step}, and the factor~$\frac{1}{4}$ compensates for the losses incurred elsewhere in the induction.  We prove it by the multiscale variational argument of~\citet{AKMBook}, using the estimates of~\citet{AK.HC}, whose constants depend on coarse-grained ellipticity constants rather than on pointwise bounds on the coefficient.  We first bound~$J$ by a negative Sobolev norm of the gradient of its maximizer and by the amount by which the values of~$J$ on the subcubes exceed its value on the cube (Lemma~\ref{l.J.bound.by.Besov}).  Throughout this subsection, for~$p,q\in\Rd$, a coefficient field~$\a$,~$n\in\Z$ and~$z\in3^n\Zd$, we write~$v_{n,z}\coloneqq v(\cdot,z+\cu_n,p,q\,;\a)$ for the maximizer of~$J(z+\cu_n,p,q\,;\a)$.

\begin{lemma}[{Estimate of~$J$ by a negative Sobolev norm}]
\label{l.J.bound.by.Besov}
There exists~$C(d)<\infty$ such that, for every~$m\in\N$, coefficient field~$\a:\cu_m \to \R^{d\times d}_{\mathrm{sym},+}$ and~$p,q\in\Rd$,
\begin{align}
	J(\cu_{m-1},p,q\,;\a)
	&
	\leq
	C {\Lambda}_{\nf14}(\cu_m;\a)
	\biggl(\frac{{\Lambda}_{\nf14}(\cu_m\,;\a)}{{\lambda}_{\nf14}(\cu_m\,;\a)} \biggr)^{\!\nf12}
	3^{-2m} \| \nabla v(\cdot,\cu_m,p,q\,;\a) \|_{\Hminusul(\cu_m)}^2
	\notag \\ & \qquad
	+
	\sum_{z\in 3^{m-1}\Zd\cap \cu_m}
	2( J(z+\cu_{m-1},p,q\,;\a) - J (\cu_m,p,q\,;\a) )
	\,.
	\label{e.J.L.bound}
\end{align}
\end{lemma}
\begin{proof}
We compute, using the triangle inequality,~\eqref{e.secondvar} and~\eqref{e.Jenergy.v},
\begin{align*}
	J(\cu_{m-1},p,q\,;\a)
	&
	=
	\fint_{\cu_{m-1}}
	\frac12
	\nabla v_{m-1,0}
	\cdot \a\nabla v_{m-1,0}
	\notag \\ &
	\leq
	\fint_{\cu_{m-1}}
	\nabla v_{m,0}
	\cdot \a\nabla v_{m,0}
	+
	\fint_{\cu_{m-1}}
	( \nabla v_{m,0} - \nabla v_{m-1,0})
	\cdot
	\a( \nabla v_{m,0} - \nabla v_{m-1,0})
	\notag \\ &
	\leq
	\fint_{\cu_{m-1}}
	\nabla v_{m,0}
	\cdot \a\nabla v_{m,0}
	\notag \\ & \qquad
	+
	\sum_{z\in 3^{m-1}\Zd\cap \cu_m}
	\fint_{z+\cu_{m-1}}
	( \nabla v_{m,0}{-}\nabla v_{m-1,z} )
	\cdot
	\a( \nabla v_{m,0} {-} \nabla v_{m-1,z})
	\notag \\ &
	=
	\fint_{\cu_{m-1}}
	\nabla v_{m,0}
	\cdot \a\nabla v_{m,0}
	+
	\sum_{z\in 3^{m-1}\Zd\cap \cu_m}
	\!\!
	2( J(z+\cu_{m-1},p,q\,;\a) - J (\cu_m,p,q\,;\a) )
	\,.
\end{align*}
The coarse-grained Caccioppoli inequality (Proposition~\ref{p.coarse.grained.caccioppoli}) yields
\begin{equation*}
	\fint_{\cu_{m-1}}
	\nabla v_{m,0}
	\cdot \a \nabla v_{m,0}
	\leq
	C {\Lambda}_{\nf14}(\cu_m;\a)
	\biggl(\frac{{\Lambda}_{\nf14}(\cu_m\,;\a)}{{\lambda}_{\nf14}(\cu_m\,;\a)} \biggr)^{\!\nf12}
	3^{-2m}  \fint_{\cu_m}  | v_{m,0}-(v_{m,0})_{\cu_m}|^2
	\,.
\end{equation*}
Combining these and using~$\| v_{m,0}-(v_{m,0})_{\cu_m} \|_{\underline{L}^2(\cu_m)} \leq C \| \nabla v_{m,0} \|_{\Hminusul(\cu_m)}$, we obtain~\eqref{e.J.L.bound}.
\end{proof}

\begin{lemma}[Negative-norm estimates]
\label{l.Besov.norms}
There exists~$C(d)<\infty$ such that, for every~$m\in\N$, coefficient field~$\a:\cu_m \to \R^{d\times d}_{\mathrm{sym},+}$, exponent~$s\in (0,1)$, integer~$L_0\leq m$,  and~$p,q\in\Rd$, we have
\begin{align}
	\lefteqn{
	3^{-m} \bigl\| \nabla v(\cdot,\cu_m,p,q\,;\a) - ( \a_{*}^{-1}(\cu_m) q -p ) \bigr\|_{\Hminusul(\cu_m)}
	}
	\qquad &
	\notag \\ &
	\leq
	C \lambda_s^{-\nf12}(\cu_m;\a)
	\sum_{k=L_0}^m
	3^{-(1-s)(m-k) }
	\biggl( \avsum_{z\in3^k\Zd \cap \cu_m}
	\!\!
	J(z+\cu_k,p,q\,;\a) - J(\cu_m,p,q\,;\a) \biggr)^{\!\nf12}
	\notag \\ & \qquad
	+
	C\sum_{n=L_0}^m
	3^{-(m-n)}
	\biggl(\avsumcube{z}{n}{m}
	\bigl| \a_{*}^{-1}(\cu_m) q - \a_{*}^{-1}(z+\cu_n ) q \bigr|^2
	\biggr)^{\!\nf12}
	\notag \\ & \qquad
	+\frac{C{\lambda}_{s}^{-\nicefrac12}(\cu_m;\a) }{1-s} 3^{-(1-s) (m-L_0)}
	\| \a^{\nicefrac12} \nabla v(\cdot,\cu_m,p,q\,;\a)  \|_{\underline{L}^2(\cu_m)}
	\notag \\ & \qquad
	+
	C 3^{-(m-L_0)}  \bigl| \a_{*}^{-1}(\cu_m) q -p\bigr|
	\,.
	\label{e.Besov.norms.gradient}
\end{align}
\end{lemma}
\begin{proof}
The duality inequality~\eqref{e.duality.splitting}, applied to~$\nabla v_{m,0}-(\nabla v_{m,0})_{\cu_m}$ and an~$H^1$ test function, together with the Poincar\'e inequality on cubes, gives the first inequality below. For the second, we split the sum at~$L_0$ and use the triangle inequality:
\begin{align}
	\lefteqn{
	3^{-m} \bigl\| \nabla v_{m,0} - ( \a_{*}^{-1}(\cu_m) q -p ) \bigr\|_{\Hminusul(\cu_m)}
	} \qquad &
	\notag \\ &
	\leq
	C\sum_{n=-\infty}^m
	3^{-(m-n)}
	\biggl(\avsumcube{z}{n}{m}
	\bigl| (\nabla v_{m,0})_{z+\cu_n} - (\a_{*}^{-1}(\cu_m) q -p)\bigr|^2
	\biggr)^{\!\nf12}
	\notag \\ &
	\leq
	C\sum_{n=L_0}^m
	3^{-(m-n)}
	\biggl(\avsumcube{z}{n}{m}
	\bigl| (\nabla v_{m,0})_{z+\cu_n} - ( \a_{*}^{-1}(z+\cu_n ) q - p)\bigr|^2
	\biggr)^{\!\nf12}
	\notag \\ & \qquad
	+
	C\sum_{n=L_0}^m
	3^{-(m-n)}
	\biggl(\avsumcube{z}{n}{m}
	\bigl| \a_{*}^{-1}(\cu_m) q - \a_{*}^{-1}(z+\cu_n ) q \bigr|^2
	\biggr)^{\!\nf12}
	\notag \\ & \qquad
	+
	C\sum_{n=-\infty}^{L_0-1}
	3^{-(m-n)}
	\biggl(\avsumcube{z}{n}{m}
	\Bigl(
	\bigl| (\nabla v_{m,0})_{z+\cu_n} \bigr|^2
	+ \bigl| \a_{*}^{-1}(\cu_m) q -p  \bigr|^2
	\Bigr)
	\biggr)^{\!\nf12}
	\,.
	\label{e.msp.plus.split}
\end{align}
For~$n \in \Z$ with~$n \leq L_0-1$, we have, by~\eqref{e.energymaps} and~\eqref{e.bound.one.cube.by.lambdas},
\begin{align*}
	\avsumcube{z}{n}{m}
	\bigl| (\nabla v_{m,0})_{z+\cu_n} \bigr|^2
	&
	\leq
	\max_{z\in 3^n\Zd \cap \cu_m}
	\bigl| \a_{*}^{-1}(z+\cu_n) \bigr| \| \a^{\nicefrac12} \nabla v_{m,0} \|_{\underline{L}^2(\cu_m)}^2
	\notag \\ &
	\leq
	3^{2s(m-n)} {\lambda}_{s}^{-1}(\cu_m;\a)
	\| \a^{\nicefrac12} \nabla v_{m,0} \|_{\underline{L}^2(\cu_m)}^2
	\,.
\end{align*}
This gives the following bound for the third term on the right side of~\eqref{e.msp.plus.split}:
\begin{align}
	\lefteqn{
	\sum_{n=-\infty}^{L_0-1}
	3^{-(m-n)}
	\biggl(\avsumcube{z}{n}{m}
	\Bigl(
	\bigl| (\nabla v_{m,0})_{z+\cu_n} \bigr|^2
	+ \bigl| \a_{*}^{-1}(\cu_m) q -p  \bigr|^2
	\Bigr)
	\biggr)^{\!\nicefrac12}
	} \qquad &
	\notag \\ &
	\leq
	\frac{C}{1-s} 3^{-(1-s) (m-L_0)}  {\lambda}_{s}^{-\nicefrac12}(\cu_m;\a)
	\| \a^{\nicefrac12} \nabla v_{m,0} \|_{\underline{L}^2(\cu_m)} +
	C 3^{-(m-L_0)}  \bigl| \a_{*}^{-1}(\cu_m) q -p\bigr|
	\,.
	\label{e.msp.plus.split.term.1}
\end{align}
Next we use~\eqref{e.v.spatial.averages} and~\eqref{e.energymaps} to obtain, for each~$n\leq m$,
\begin{align*}
	\bigl|
	\a_*^{\nf12} (z+\cu_n)
	\bigl(
	( \nabla v_{m,0} )_{z+\cu_n}   - \bigl(\a_{*}^{-1}(z+\cu_n) q - p \bigr) \bigr) \bigr|^{2}
	&
	=
	\bigl|
	\a_*^{\nf12} (z+\cu_n)
	( \nabla v_{m,0} - \nabla v_{n,z} )_{z+\cu_n}
	\bigl|^{2}
	\notag \\ &
	\leq
	\fint_{z+\cu_n}
	( \nabla v_{m,0} - \nabla v_{n,z} )\cdot \a ( \nabla v_{m,0} - \nabla v_{n,z} )\,.
\end{align*}
Summing this over~$z\in 3^n\Zd \cap \cu_m$ and applying~\eqref{e.secondvar} yields, for each~$n\leq m$,
\begin{align*}
	\lefteqn{
	\avsumcube{z}{n}{m}
	\Bigl|
	\a_*^{\nf12} (z{+}\cu_n)
	\bigl(
	( \nabla v_{m,0} )_{z+\cu_n}   - \bigl(\a_{*}^{-1}(z+\cu_n) q - p \bigr) \bigr) \Bigr|^{2}
	}   &
	\notag \\ &
	\leq
	\avsumcube{z}{n}{m}
	\fint_{z+\cu_n}
	( \nabla v_{m,0} {-} \nabla v_{n,z} )\cdot \a ( \nabla v_{m,0} {-} \nabla v_{n,z} )
	\notag \\ &
	=
	2 \!\!
	\avsum_{z\in3^n\Zd \cap \cu_m}
	\!\!
	\bigl( J(z{+}\cu_n,p,q\,;\a) - J(\cu_m,p,q\,;\a) \bigr)
	\,.
\end{align*}
Using that, by~\eqref{e.bound.one.cube.by.lambdas}, for every~$P \in \Rd$, every~$n\leq m$ and every~$z\in 3^n\Zd \cap \cu_m$,
\begin{equation*}
	| P |^2
	\leq
	\bigl| \a_*^{-1} (z+\cu_n)
	\bigr|
	\bigl|
	\a_*^{\nf12} (z+\cu_n)
	P \bigr|^{2}
	\leq
	3^{2s(m-n)} {\lambda}_{s}^{-1}(\cu_m;\a)
	\bigl|
	\a_*^{\nf12} (z+\cu_n)
	P \bigr|^{2}\,,
\end{equation*}
we obtain the following bound for the first term on the right side of~\eqref{e.msp.plus.split}:
\begin{align}
	\lefteqn{
	\sum_{n=L_0}^m
	3^{-(m-n)}
	\biggl(\avsumcube{z}{n}{m}
	\bigl| (\nabla v_{m,0})_{z+\cu_n} - ( \a_{*}^{-1}(z+\cu_n ) q - p)\bigr|^2
	\biggr)^{\!\nf12}
	} \qquad &
	\notag \\ &
	\leq
	C \lambda_s^{-\nf12}(\cu_m;\a)
	\sum_{k=L_0}^m
	3^{-(1-s)(m-k) }
	\biggl( \avsum_{z\in3^k\Zd \cap \cu_m}
	\!\!
	J(z+\cu_k,p,q\,;\a) - J(\cu_m,p,q\,;\a) \biggr)^{\!\nf12}\,.
	\label{e.msp.plus.split.term.3}
\end{align}
Combining~\eqref{e.msp.plus.split},~\eqref{e.msp.plus.split.term.1} and~\eqref{e.msp.plus.split.term.3} yields~\eqref{e.Besov.norms.gradient}.
\end{proof}

For a square-integrable random matrix~$M$, we write~$\var[M]\coloneqq\E[|M-\E M|^2]$, with the matrix norm fixed in Section~\ref{s.notation}. If~$\E[|M|^2]=+\infty$, we let~$\var[M]\coloneqq+\infty$.

\begin{lemma}[Variance estimate]
\label{l.var.bounds}
For every~$\Zd$-stationary random coefficient field~$\a:\Rd \to \R^{d\times d}_{\mathrm{sym},+}$, every~$n,k\in\Z$ with~$k<n$, and every positive-definite matrix~$\b\in \R^{d\times d}_{\mathrm{sym},+}$, we have
\begin{align}
	\var\big[ \a_{*}^{-1}(\cu_n)  \bigr]
	&
	\leq
	2\var\Biggl[
	\avsum_{z\in 3^k\Zd \cap \cu_n} \!\!\!\!
	\a_{*}^{-1} (z+\cu_k) \Biggr]
	\notag \\ & \qquad
	+
	32 \E \Biggl[
	\biggl( \avsum_{z\in 3^k\Zd\cap \cu_n}
	\sum_{i=1}^d
	J(z+\cu_k,\b^{-1} e_i, e_i\,;\a)
	\biggr)^{\!2}
	\Biggr]
	\,.
	\label{e.var.a.star}
\end{align}
\end{lemma}
\begin{proof}
If the right side of~\eqref{e.var.a.star} is infinite, the assertion is immediate. Otherwise the average of the positive-definite matrices~$\a_*^{-1}(z+\cu_k)$ is in~$L^2(\P)$ by the variance convention. Subadditivity gives
\begin{equation*}
	0\leq\a_*^{-1}(\cu_n)
	\leq\avsum_{z\in3^k\Zd\cap\cu_n}\a_*^{-1}(z+\cu_k)\,.
\end{equation*}
Hence~$\a_*^{-1}(\cu_n)$ is also in~$L^2(\P)$. Each subcube matrix is bounded in the Loewner order by the number of subcubes times their average, so it too is in~$L^2(\P)$. The centered matrices and their differences below are therefore well defined and square-integrable. For the operator norm, the triangle inequality and Jensen's inequality give~$\var[M]\leq4\E[|M|^2]$ for every square-integrable random matrix~$M$. We apply this to the difference between~$\a_*^{-1}(\cu_n)$ and the average of~$\a_*^{-1}(z+\cu_k)$ over~$z\in3^k\Zd\cap\cu_n$, and obtain, for every~$k<n$,
\begin{align}
	\var\big[ \a_{*}^{-1}(\cu_n)  \bigr]
	\leq
	2\var\Biggl[ \avsum_{z\in 3^k\Zd \cap \cu_n}
	\!\!\!
	\a_{*}^{-1} (z+\cu_k) \Biggr]
	+
	8\E \Biggl[
	\biggl|
	\avsum_{z\in 3^k\Zd \cap \cu_n}
	\!\!\!
	\a_{*}^{-1} (z+\cu_k)
	-
	\a_{*}^{-1} (\cu_n)
	\biggr|^2 \Biggr]\,.
	\label{e.var.splitting}
\end{align}
We compare the arithmetic and harmonic means. Let~$\b_1,\ldots,\b_N$ be symmetric positive-definite matrices, and write~$\mathbf h\coloneqq(\avsum_{i=1}^N\b_i^{-1})^{-1}$ for their harmonic mean. For every positive-definite matrix~$\b$,
\begin{align}
	\avsum_{i=1}^N\b_i
	&
	=
	\mathbf{h}+\avsum_{i=1}^N  (\b_i - \b ) \b_i^{-1} (\b_i - \b)
	-
	\underbrace{(\mathbf{h} - \b) \mathbf{h}^{-1} (\mathbf{h} - \b)}_{\geq 0}
	\leq
	\mathbf{h}+\avsum_{i=1}^N  (\b_i - \b ) \b_i^{-1} (\b_i - \b)
	\,.
	\label{e.arithmetic.harmonic.mean}
\end{align}
Apply~\eqref{e.arithmetic.harmonic.mean} with~$\{ \a^{-1} (z+\cu_k) \,:\, z \in 3^k\Zd \cap \cu_n\}$ in place of the~$\b_i$ and combine this with~\eqref{e.J.plug.identity.inverses} and~\eqref{e.ord.bounds.truncated} to obtain, for every positive-definite matrix~$\b$,
\begin{align*}
	\avsum_{z\in 3^k\Zd \cap \cu_n} \a_{*}^{-1}  (z+\cu_k)
	-
	\biggl( \avsum_{z\in 3^k\Zd \cap \cu_n} \a (z+\cu_k) \biggr)^{\!-1}
	\leq
	2 \avsum_{z\in 3^k\Zd \cap \cu_n}
	\sum_{i=1}^d
	J(z+\cu_k, \b^{-1} e_i,e_i\,;\a) \Id
	\,.
\end{align*}
Using the previous display and the subadditivity of the coarse-grained matrices~$\a_{*}^{-1}$ and~$\a$, we find
\begin{align}
	\a^{-1}(\cu_n)
	\leq
	\a_{*}^{-1} (\cu_n)
	&
	\leq
	\avsum_{z\in 3^k\Zd \cap \cu_n}
	\a_{*}^{-1} (z+\cu_k)
	\notag \\ &
	\leq
	\biggl( \avsum_{z\in 3^k\Zd \cap \cu_n} \a (z+\cu_k) \biggr)^{\!-1}
	+ 2 \biggl (\avsum_{z\in 3^k\Zd\cap \cu_n}
	\sum_{i=1}^d
	J(z+\cu_k,\b^{-1} e_i, e_i\,;\a) \biggr )\Id
	\notag \\ &
	\leq
	\a^{-1}(\cu_n) +
	2\biggl (
	\avsum_{z\in 3^k\Zd\cap \cu_n}
	\sum_{i=1}^d
	J(z+\cu_k,\b^{-1} e_i, e_i \,;\a)\biggr )\Id
	\,.
	\label{e.astar.inverse.comparison}
\end{align}
We deduce from~\eqref{e.astar.inverse.comparison} that
\begin{equation*}
	0 \leq
	\avsum_{z\in 3^k\Zd \cap \cu_n}
	\a_{*}^{-1}(z+\cu_k)
	-
	\a_{*}^{-1} (\cu_n)
	\leq
	2\biggl (
	\avsum_{z\in 3^k\Zd\cap \cu_n}
	\sum_{i=1}^d
	J(z+\cu_k,\b^{-1} e_i, e_i\,;\a)\biggr )\Id\,,
\end{equation*}
and therefore
\begin{equation*}
	\E \Biggl[
	\biggl|
	\avsum_{z\in 3^k\Zd \cap \cu_n} \a_{*}^{-1} (z+\cu_k)
	-
	\a_{*}^{-1} (\cu_n)
	\biggr|^2 \Biggr]
	\leq
	4 \E \Biggl[
	\biggl( \avsum_{z\in 3^k\Zd\cap \cu_n}
	\sum_{i=1}^d
	J(z+\cu_k,\b^{-1} e_i, e_i\,;\a)
	\biggr)^{\!2}
	\Biggr]
	\,.
\end{equation*}
Combining this estimate with~\eqref{e.var.splitting} yields~\eqref{e.var.a.star}.
\end{proof}

We next combine the previous four lemmas into an upper bound on~$\E \bigl[ J(\cu_m,p,q\,;\a_L) \bigr]$.

\begin{proposition}
\label{p.combine.under.S}
There exist constants~$c(d)>0$ and~$C(d)<\infty$ such that the following holds.  Let~$L\in\N_0$,~$\xi\in[16d,\infty)$ and~$\delta_1\in[C\xi\gmcdelta^2,c]$ be such that~$\S(L,\xi,\delta_1)$ is valid. For every~$m\in\N$ with~$m\geq L$ and~$C3^{-(m-L)}\leq\nf14$, and every~$p,q\in\Rd$ satisfying
\begin{equation}
	q = \ahom_L p
	\quad \text{and} \quad
	|q|^2 \leq \ahom_L \,,
	\label{e.pq.normed}
\end{equation}
we have the estimate
\begin{align}
	\lefteqn{
	\E \bigl[ J(\cu_m,p,q\,;\a_L) \bigr]
	} \qquad &
	\notag \\ &
	\leq
	C
	\delta_1 \bigl( \delta_1 + 3^{-(m-L)} \bigr)
	+
	C
	\sum_{n=L}^m
	3^{-(m-n)}
	\E \bigl[  J(\cu_n,p,q\,;\a_L) - J(\cu_m,p,q\,;\a_L) \bigr]
	\,.
	\label{e.J.bound.by.indyhyp}
\end{align}
\end{proposition}
\begin{proof}
The proof proceeds in five steps. Step~1 bounds~$\E[J(\cu_m,p,q\,;\a_L)]$ by the sum on the right side of~\eqref{e.J.bound.by.indyhyp} and three further terms: the distances of the expectations~$\ahom_{L,*}^{-1}(\cu_n)$ of~$\a_{L,*}^{-1}(\cu_n)$ from~$\ahom_L^{-1}$, the variances of~$\a_{L,*}^{-1}(\cu_n)$, and the contribution of the complement of the event~$\mathsf G$ defined below, on which the coarse-grained ellipticity constants of~$\a_L$ on~$\cu_m$ lie within a factor~$2$ of~$\ahom_L$. Steps~2,~3 and~4 bound these three terms using~$\S(L,\xi,\delta_1)$, and Step~5 combines the bounds. The heart of the matter is Step~3: there the independence of~$\a_L$ on distant cubes of scale~$L$ makes the variance of the average of~$\a_{L,*}^{-1}(z+\cu_L)$ over~$z\in3^L\Zd\cap\cu_n$ decay like~$3^{-d(n-L)}$, which gives the term~$\delta_13^{-(m-L)}$ in~\eqref{e.J.bound.by.indyhyp}.

\smallskip

Define the event~$\mathsf{G}$ by
\begin{equation}
	\mathsf{G}  \coloneqq
	\Bigl\{
	\frac 12 \ahom_{L}
	\leq
	{\lambda}_{\nf14}(\cu_m;\a_L)
	\leq
	{\Lambda}_{\nf14}(\cu_m;\a_L)
	\leq
	2 \ahom_{L}
	\Bigr\}\,.
	\label{e.goodevent.combine}
\end{equation}

\smallskip

\emph{Step 1.} We combine the previous lemmas to obtain an upper bound for~$\E \bigl[ J(\cu_{m},p,q\,;\a_L) \bigr]$. We show that there exists~$C(d)<\infty$ such that, if~$3^{-\frac{3}{2}(m-L)}C \leq \nf14$, then
\begin{align}
	\E \bigl[ J(\cu_{m},p,q\,;\a_L) \bigr]
	&
	\leq
	C
	\sum_{n=L}^m 3^{-(m-n)}
	\E \bigl[  J (\cu_n,p,q\,;\a_L) - J (\cu_m,p,q\,;\a_L) \bigr]
	\notag \\ & \qquad
	+
	C \ahom^2_{L}
	\sum_{n=L}^m 3^{-(m-n)}
	\Bigl(
	\var
	\bigl[  \a_{L,*}^{-1}(\cu_n) \bigr]
	+ \bigl| \ahom_{L,*}^{-1}(\cu_n) -  \ahom_L^{-1} \bigr|^2 \Bigr)
	\notag \\ & \qquad
	+
	C \E \bigl[ J (\cu_{m-1},p,q\,;\a_L) \indc_{\mathsf{G}^{\mathrm c}}  \bigr]
	\,.
	\label{e.initial.JL.bound}
\end{align}
Taking expectations in~\eqref{e.J.L.bound} and using the subadditivity of~$J$ and stationarity, we obtain
\begin{align*}
	\E \bigl[ J(\cu_{m},p,q\,;\a_L) \bigr]
	&
	\leq
	\E \bigl[ J (\cu_{m-1},p,q\,;\a_L) \bigr]
	\notag \\ &
	=
	\E \bigl[ J (\cu_{m-1},p,q\,;\a_L) \indc_{\mathsf{G}}  \bigr]
	+
	\E \bigl[ J (\cu_{m-1},p,q\,;\a_L) \indc_{\mathsf{G}^{\mathrm c}}  \bigr]
	\notag \\ &
	\leq
	C \E \biggl[  {\Lambda}_{\nf14}(\cu_m;\a_L)
	\biggl( \frac{{\Lambda}_{\nf14}(\cu_m;\a_L)}{{\lambda}_{\nf14}(\cu_m;\a_L)} \biggr)^{\!\nf12}
	3^{-2m}\| \nabla v(\cdot,\cu_m,p,q\,;\a_L) \|_{\Hminusul(\cu_m)}^2 \indc_{\mathsf{G}}  \biggr]
	\notag \\ &  \qquad
	+
	C\E \bigl[  J (\cu_{m-1},p,q\,;\a_L) - J (\cu_m,p,q\,;\a_L) \bigr]
	\notag \\ & \qquad
	+
	C \E \bigl[ J (\cu_{m-1},p,q\,;\a_L) \indc_{\mathsf{G}^{\mathrm c}}  \bigr]
	\,.
\end{align*}
Applying Lemma~\ref{l.Besov.norms}, squaring both sides of~\eqref{e.Besov.norms.gradient}, taking expectations and using that
\begin{equation*}
	\biggl( \frac{{\Lambda}_{\nf14}(\cu_m;\a_L)}{{\lambda}_{\nf14}(\cu_m;\a_L)} \biggr)^{\!\nf32}
	\indc_{\mathsf{G}}
	\leq
	8
	\,,
\end{equation*}
and that~${\Lambda}_{\nf14}(\cu_m;\a_L)\indc_{\mathsf{G}}\leq2\ahom_{L}$, we obtain
\begin{align*}
	\lefteqn{
	\E \biggl[  {\Lambda}_{\nf14}(\cu_m;\a_L)
	\biggl( \frac{{\Lambda}_{\nf14}(\cu_m;\a_L)}{{\lambda}_{\nf14}(\cu_m;\a_L)} \biggr)^{\!\nf12}
	3^{-2m}\| \nabla v(\cdot,\cu_m,p,q\,;\a_L) \|_{\Hminusul(\cu_m)}^2 \indc_{\mathsf{G}}  \biggr]
	}
	\qquad &
	\notag \\ &
	\leq
	C
	\sum_{n=L}^m
	3^{-(m-n)}
	\E \bigl[  J (\cu_n,p,q\,;\a_L) - J (\cu_m,p,q\,;\a_L) \bigr]
	\notag \\ & \qquad
	+
	C\ahom_L \sum_{n=L}^m
	3^{-(m-n)}
	\E \biggl[ \avsumcube{z}{n}{m}
	\bigl| \a_{L,*}^{-1}(\cu_m) q - \a_{L,*}^{-1}(z+\cu_n ) q \bigr|^2 \biggr]
	\notag \\ & \qquad
	+
	C\ahom_L \E\Bigl[ \bigl| \a_{L,*}^{-1} (\cu_m)q  - p  \bigr|^2 \Bigr]
	\notag \\ & \qquad
	+
	C
	3^{-\frac{3}{2} (m-L)}
	\E \Bigl[ \| \a_L^{\nicefrac12} \nabla v \|_{\underline{L}^2(\cu_m)}^2 \Bigr]
	\,.
\end{align*}
Under the condition~$C3^{-3(m-L)/2}\leq\nf14$, we use~\eqref{e.Jenergy.v} to bound the last term by~$\frac{1}{2}\E[J(\cu_m,p,q\,;\a_L)]$, which we absorb into the left side. By stationarity and the triangle inequality,
\begin{align*}
	\lefteqn{
	\E \biggl[ \avsumcube{z}{n}{m}
	\bigl| \a_{L,*}^{-1}(\cu_m) - \a_{L,*}^{-1}(z+\cu_n ) \bigr|^2 \biggr]
	} \   &
	\notag \\ &
	\qquad \leq
	3 \var\bigl[ \a_{L,*}^{-1}(\cu_m ) \bigr]
	+
	3 \var\bigl[ \a_{L,*}^{-1}(\cu_n ) \bigr]
	+
	3\bigl| \ahom_{L,*}^{-1} (\cu_m) - \ahom_{L,*}^{-1}(\cu_n) \bigr|^2
	\,.
\end{align*}
Combining the previous two displays and using~\eqref{e.pq.normed}, we obtain~\eqref{e.initial.JL.bound}.

\smallskip

\emph{Step 2.} Assuming~$\S(L,\xi,\delta_1)$, we prove that, for every~$n\geq L$,
\begin{equation}
	\bigl| \ahom_{L,*}^{-1}(\cu_n)-\ahom_L^{-1}\bigr|^2
	\leq
	C \ahom_{L}^{-2} \delta_1^2
	\,.
	\label{e.means.downscale.by.defect}
\end{equation}
Indeed, subadditivity on the partition into cubes of scale~$L$ and the induction hypothesis give, for every~$|e|=1$,
\begin{equation*}
	\E\bigl[
	J(\cu_n,\ahom_L^{-\nf12}e,\ahom_L^{\nf12}e;\a_L)
	\bigr]
	\leq\delta_1\,.
\end{equation*}
On the other hand,~\eqref{e.coarse.grained.matrices}, the definition of the expected coarse-grained matrices, and their monotonicity imply
\begin{equation*}
	2\E\bigl[
	J(\cu_n,\ahom_L^{-\nf12}e,\ahom_L^{\nf12}e;\a_L)
	\bigr]
	\geq
	\ahom_L\ahom_{L,*}^{-1}(\cu_n)-1
	\geq0\,.
\end{equation*}
This proves~\eqref{e.means.downscale.by.defect}.

\smallskip

\emph{Step 3.} Using the variance estimates in Lemma~\ref{l.var.bounds}, the independence assumption, and the hypothesis~$\S(L, \xi, \delta_1)$, we show that, for every~$n\geq L$,
\begin{equation}
	\var\big[ \a_{L,*}^{-1}(\cu_n)  \bigr]
	\leq
	C  \ahom_{L}^{-2} \delta_1 \bigl( \delta_1 + 3^{-d(n-L)} \bigr)
	\,.
	\label{e.var.bound.astar.withS}
\end{equation}
For~$n=L$,~\eqref{e.var.bound.astar.withS} follows directly from~$\S(L,\xi,\delta_1)$ and~\eqref{e.coarse.grained.matrix.deviation.by.J}. By Lemma~\ref{l.var.bounds}, for every~$n>L$ we have
\begin{align*}
	\var\bigl[ \a_{L,*}^{-1}(\cu_n)  \bigr]
	\leq
	2\var\Biggl[
	\avsum_{z\in 3^L\Zd \cap \cu_n} \!\!\!\!
	\a_{L,*}^{-1} (z+\cu_L) \Biggr]
	+
	32 \E \Biggl[
	\biggl( \avsum_{z\in 3^L\Zd\cap \cu_n}
	\sum_{i=1}^d
	J(z+\cu_L,\ahom_L^{-1} e_i, e_i\,;\a_L)
	\biggr)^{\!2}
	\Biggr]
	\,.
\end{align*}
By assumption~\ref{a.g1} and the definition of~$\a_L$, for every~$y,z\in 3^L\Zd$ with~$|y-z|> d 3^{L+1}$, the random variables~$ \a_{L,*}^{-1}(y+\cu_L)$ and~$ \a_{L,*}^{-1}(z+\cu_L)$ are independent. Consequently, by~$\S(L,\xi, \delta_1)$ and~\eqref{e.coarse.grained.matrix.deviation.by.J},
\begin{align*}
	\var\Biggl[
	\avsum_{z\in 3^{L}\Zd \cap \cu_n} \!\!\!
	\a_{L,*}^{-1} (z+\cu_{L}) \Biggr]
	\leq
	C 3^{-d(n-L)}
	\var\bigl[\a_{L,*}^{-1} (\cu_{L}) \bigr]
	&
	\leq
	C 3^{-d(n-L)}
	\E \Bigl[ \bigl| \a_{L,*}^{-1}(\cu_{L}) - \ahom_{L}^{-1} \bigr|^2 \Bigr]
	\notag \\ &
	\leq
	C  \ahom_{L}^{-2}  3^{-d(n-L)} \delta_1
	\,.
\end{align*}
The identical calculation for the averages over~$z$ of~$\sum_{i=1}^dJ(z+\cu_L,\ahom_L^{-1}e_i,e_i\,;\a_L)$ gives
\begin{equation*}
	\var \Biggl[
	\avsum_{z\in 3^{L}\Zd\cap \cu_n}
	\sum_{i=1}^d
	J(z+\cu_{L},\ahom_{L}^{-1} e_i, e_i\,;\a_L)
	\Biggr]
	\leq
	C  \ahom_{L}^{-2}  3^{-d(n-L)} \delta_1^2
	\,.
\end{equation*}
Finally, using~$\S(L, \xi, \delta_1)$ again, we have
\begin{align*}
	\E \bigl[ J(z+\cu_{L}, \ahom_{L}^{-1} e_i, e_i\,;\a_L)\bigr]^2
	&
	=
	\ahom_{L}^{-2}  \E \bigl[ J(\cu_{L},\ahom_{L}^{-\nf12}e_i, \ahom_{L}^{\nf12}  e_i\,;\a_L)\bigr]^2
	\leq
	\ahom_{L}^{-2}  \delta_1^2 \,.
\end{align*}
Combining the previous three displays yields~\eqref{e.var.bound.astar.withS}.

\smallskip

\emph{Step 4.} We prove that
\begin{equation}
	\E \bigl[ J (\cu_{m-1},p,q\,;\a_L) \indc_{\mathsf{G}^{\mathrm c}}  \bigr]
	\leq
	C\delta_1^2\,.
	\label{e.badeventbound.sect42}
\end{equation}
By~\eqref{e.goodevent.combine} and~\eqref{e.bound.Lambdas.by.Es}, we have~$\mathsf G^{\mathrm c}\subseteq\{\mathcal E_{\nf14,\infty}(\cu_m;\a_L,\ahom_L)\geq\nf15\}$. Apply the bound~\eqref{e.multiscale.response.large.cubes} of Lemma~\ref{l.multiscale.response.large.cubes} with~$s=\nf14$ and the same exponent~$\xi$.  This gives, for a constant~$C_B(d)<\infty$,
\begin{equation}
	\P \bigl[ \mathsf{G}^{\mathrm c} \bigr]
	\leq
	\P \bigl[ \mathcal{E}_{\nf14,\infty} (\cu_m; \a_L, \ahom_{L} ) \geq \nf 1{5}  \bigr]
	\leq
	(C_B \delta_1)^{\nf \xi 2}
	\,.
	\label{e.bad.B.bound}
\end{equation}
We have~$\E \bigl[ J (\cu_{m-1},p,q\,;\a_L)^2 \bigr]\leq C\delta_1^2$ by subadditivity down to scale~$L$, the induction hypothesis, and~\eqref{e.pq.normed}. After decreasing~$c(d)$, we may suppose that~$C_B\delta_1\leq1$. The Cauchy--Schwarz inequality,~\eqref{e.bad.B.bound}, and~$\xi\geq4$ give
\begin{equation*}
	\E \bigl[ J (\cu_{m-1},p,q\,;\a_L) \indc_{\mathsf{G}^{\mathrm c}}  \bigr]
	\leq
	C\delta_1(C_B\delta_1)^{\nf\xi4}
	\leq
	C \delta_1^2
	\,.
\end{equation*}

\smallskip

\emph{Step 5.} We conclude. Combining~\eqref{e.initial.JL.bound},~\eqref{e.means.downscale.by.defect},~\eqref{e.var.bound.astar.withS}, and~\eqref{e.badeventbound.sect42} gives~\eqref{e.J.bound.by.indyhyp}.
\end{proof}

\begin{proposition}[Contraction above the cutoff scale]
\label{p.homogenization.step}
There exist constants~$c(d)>0$ and~$C(d)<\infty$ such that the following holds.  Let~$m_0\in\N$,~$\xi\in[16d,\infty)$ and~$\delta_1\in[C\xi\gmcdelta^2,c]$ be such that~$\S(m_0,\xi,\delta_1)$ is valid. Then, for every~$m\in\N$ and~$L\in\N_0$ with~$L\leq m_0$ and~$m\geq L+C\log(2+\xi)$, we have
\begin{equation}
	\E \biggl[
	\max_{|e| = 1}
	\bigl(
	J(\cu_m,\ahom_L^{-\nf12}e,\ahom_L^{\nf12}e;\a_L)
	\bigr) ^{\xi}
	\biggr]^{\nf1\xi}
	\leq
	\frac14 \delta_1
	\,.
	\label{e.homogenization.step.bound}
\end{equation}
\end{proposition}
\begin{proof}
Fix~$p,q\in\Rd$ satisfying~\eqref{e.pq.normed}. By subadditivity on the partition into cubes of scale~$L$, the induction hypothesis and~\eqref{e.pq.normed}, for every~$n\geq L$,
\begin{equation}
	\E \bigl[ J (\cu_n,p,q\,;\a_L)^\xi \bigr]^{\nf 1\xi}  \leq \delta_1\,.
	\label{e.iter.init}
\end{equation}
Proposition~\ref{p.combine.under.S} gives~\eqref{e.J.bound.by.indyhyp} at the scales~$m\geq L$ with~$C3^{-(m-L)}\leq\nf14$; at each of the at most~$C(d)$ smaller scales~$n\geq L$, the bound~\eqref{e.iter.post} below follows from~\eqref{e.iter.init} and Jensen's inequality once its constant is large enough. Iterating~\eqref{e.J.bound.by.indyhyp} over the larger scales, starting from~\eqref{e.iter.init}, we obtain an exponent~$\alpha(d)>0$ and a constant~$C(d)<\infty$ such that, for every~$n\in\N$ with~$n\geq L$,
\begin{equation}
	\E \bigl[ J (\cu_n,p,q\,;\a_L) \bigr]
	\leq
	C\delta_1 \bigl( \delta_1 + 3^{-\alpha (n-L)} \bigr)
	\,.
	\label{e.iter.post}
\end{equation}
Using subadditivity and the triangle inequality, we have, for every~$m\geq n\geq L$,
\begin{align*}
	\lefteqn{
	\E \bigl[ J (\cu_m,p,q\,;\a_L)^\xi \bigr] ^{\nf1\xi}
	} \quad &
	\notag \\ &
	\leq
	\E \biggl[
	\biggl(
	\avsumcube{z}{n}{m}   J (z+\cu_n,p,q\,;\a_L) \biggr)^\xi
	\biggr]^{\nf1\xi}
	\notag \\ &
	\leq
	\E \Biggl[
	\biggl|
	\avsumcube{z}{n}{m}
	\bigl( J (z+\cu_n,p,q\,;\a_L) - \E \bigl[ J (\cu_n,p,q\,;\a_L) \bigr] \bigr)  \biggr|^\xi
	\Biggr]^{\nf1\xi}
	+
	\E \bigl[ J (\cu_n,p,q\,;\a_L) \bigr]
	\,.
\end{align*}
Since~$J(z+\cu_n,p,q\,;\a_L)$ depends only on~$\a_L|_{z+\cu_n}$ and~$\a_L$ has range of dependence at most~$\sqrt{d}3^n$, we may group the summands into~$C(d)$ families of independent random variables. Applying Rosenthal's inequality (Lemma~\ref{l.Rosenthal}) to each family and using~\eqref{e.iter.init}, we find that
\begin{align}
	\lefteqn{
	\E \Biggl[
	\biggl|
	\avsumcube{z}{n}{m}
	\bigl( J (z+\cu_n,p,q\,;\a_L) - \E \bigl[ J (\cu_n,p,q\,;\a_L) \bigr] \bigr)  \biggr|^\xi
	\Biggr]^{\nf1\xi}
	} \quad &
	\notag \\ &
	\leq
	C 3^{-\frac{d}{2}(m-n)}
	\xi^{\nf12}
	\E \bigl[ J (\cu_n,p,q\,;\a_L)^2 \bigr] ^{\nf12}
	+ C 3^{-\frac{d(\xi-1)}{\xi} (m-n)} \xi
	\E \bigl[
	J (\cu_n,p,q\,;\a_L)^\xi
	\bigr]^{\nf1\xi}
	\notag \\ &
	\leq
	C 3^{-\frac{d}{2}(m-n)} \xi^{\nf12}
	\delta_1
	+
	C 3^{-\frac{d(\xi-1)}{\xi} (m-n)} \xi \delta_1
	\notag \\ &
	\leq
	C 3^{-\frac{d}{2}(m-n)} \xi \delta_1
	\,.
	\label{e.appl.Rosenthal}
\end{align}
Combining the previous two displays with~\eqref{e.iter.post}, we obtain
\begin{equation*}
	\E \bigl[ J (\cu_m,p,q\,;\a_L)^\xi \bigr] ^{\nf1\xi}
	\leq
	C 3^{-\frac{d}{2}(m-n)} \xi \delta_1 +
	C\delta_1 \bigl( \delta_1 + 3^{-\alpha (n-L)} \bigr) \,.
\end{equation*}
Selecting~$n \coloneqq  L + \lceil \frac{1}{2} (m-L) \rceil$ yields
\begin{equation*}
	\E \bigl[ J (\cu_m,p,q\,;\a_L)^\xi \bigr] ^{\nf1\xi}
	\leq
	C 3^{-\frac{d}{4} (m-L)} \xi \delta_1 +
	C\delta_1 \bigl( \delta_1 + 3^{-\frac{1}{2} \alpha (m-L)} \bigr) \,.
\end{equation*}
If~$m\geq L+C\log(2+\xi)$, then the right side is at most~$c(d)\delta_1$.  To pass from a fixed direction to the supremum in~\eqref{e.homogenization.step.bound}, we use a fixed~$\nf14$-net of the unit sphere.  Since the map~$e\longmapsto J(\cu_m,\ahom_L^{-\nf12}e,\ahom_L^{\nf12}e;\a_L)$ is a nonnegative quadratic form, its supremum on the sphere is at most twice its maximum on the net.  Bounding the maximum over the net by the sum over its points, whose number depends only on~$d$, we obtain
\begin{equation*}
	\E \biggl[
	\max_{|e|=1}
	J(\cu_m,\ahom_L^{-\nf12}e,\ahom_L^{\nf12}e;\a_L)^\xi
	\biggr] ^{\nf1\xi}
	\leq
	\frac 14 \delta_1
	\,.
\end{equation*}
This completes the proof.
\end{proof}

To complete the proof of Proposition~\ref{p.coarse.grained.bound}, we compare the coefficient at each scale with the coefficient~$h$ scales below it using Lemma~\ref{l.J.sensitivity}, and choose~$h$ so that the resulting error is absorbed by the gain in Proposition~\ref{p.homogenization.step}.

\begin{proof}[{Proof of~Proposition~\ref{p.coarse.grained.bound}}]
We use two inductions. The first proves~$\S(\infty,\xi,C\xi\log^2(2+\xi)\gmcdelta^2)$; the second removes one factor of~$\log(2+\xi)$.

\smallskip

\emph{Step 1.} We prove the base case. For both inductions, it is the assertion that, for some constant~$C(d)<\infty$ and every~$\xi\in [1,\infty)$ and~$h \in \N$ satisfying~$C \xi \gmcdelta^2 h \leq 1$,
\begin{equation}
	\S\bigl( 3 h, \xi, C \xi \gmcdelta^2 h\bigr) \quad \text{is valid}\,.
	\label{e.seed}
\end{equation}
Lemma~\ref{lemma.crude.J.bound} implies that, for every~$n\in\{0,\ldots,3h\}$,
\begin{equation*}
	\|\mathcal J_n(\cu_n)\|_{L^\xi(\P)}
	\leq C\xi\gmcdelta^2h\exp(C\xi\gmcdelta^2h)\,.
\end{equation*}
This implies~\eqref{e.seed}.

\smallskip

\emph{Step 2.} We prove that there exists~$C(d)<\infty$ such that~$\S(\infty,\xi,C\xi\log^2(2+\xi)\gmcdelta^2)$ is valid for every~$\xi\in[16d,C^{-1}\gmcdelta^{-2}|\log\gmcdelta|^{-2}]$. We assume that
\begin{equation}
	h \geq K \log(2+\xi)
	\quad \text{and} \quad
	\xi \gmcdelta^2  h^2 \leq K^{-1} \delta_1 \leq K^{-2} \,,
	\label{e.params.first.induction.step}
\end{equation}
where~$K(d)$ is a large constant to be determined below. We iterate in blocks of~$h$ scales. For parameters satisfying~\eqref{e.params.first.induction.step} and every~$m\in\N \cap [3 h,\infty)$, we prove
\begin{equation}
	\S(m, \xi, \delta_1)
	\implies
	\S(m+ h, \xi,\delta_1)
	\,.
	\label{e.imp}
\end{equation}
Equations~\eqref{e.seed} and~\eqref{e.imp} imply by induction that~$\mathcal S(\infty,\xi,\delta_1)$ holds for every choice of parameters satisfying~\eqref{e.params.first.induction.step}, after enlarging~$K$.

\smallskip

To prove~\eqref{e.imp}, we fix~$m\in\N$ with~$m\geq 3h$ such that~$\S(m, \xi, \delta_1)$ is valid. By Proposition~\ref{p.homogenization.step}, if~$K$ is taken sufficiently large, then
\begin{equation}
	\E \biggl[ \max_{|e| = 1} J\bigl(\cu_{m+h}, \ahom_m^{-\nf12}  e, \ahom_m^{\nf12}  e \,;\a_m \bigr)^{\xi}\biggr]^{\nf{1}{\xi}} \leq \frac14 \delta_1 \,.
	\label{e.gain-in-delta}
\end{equation}
We apply Lemma~\ref{l.J.sensitivity} with~$\delta=1$ and~$\lambda=\ahom_{m+h}\ahom_m^{-1}$, and with~$\b=\a_{m+h}$ and~$\a=\a_m$. This gives, for every bounded convex domain~$U\subseteq\Rd$ and~$|e|=1$,
\begin{align}
	J(U, \ahom_{m+h}^{-\nf12} e, \ahom_{m+h}^{\nf12}e\,;\a_{m+h})
	\leq
	2 J(U, \ahom_m^{-\nf12} e, \ahom_m^{\nf12}e\,;\a_m)
	+
	3 \X_{m,h}(U)^2 \bigl( J(U, \ahom_m^{-\nf12} e, \ahom_m^{\nf12}e\,;\a_m)    + 1 \bigr)
	\,,
	\label{e.J.localization.block}
\end{align}
where we have defined
\begin{equation*}
	\X_{m,h}(U)  \coloneqq
	\bigl\| \ahom_{m+h} \ahom^{-1}_m  (\a^{-1}_{m+h} \a_{m}) - 1 \bigr\|_{L^{\infty}(U)}
	+
	\bigl\| \ahom^{-1}_{m+h} \ahom_m  (\a_{m+h} \a^{-1}_{m}) - 1 \bigr\|_{L^{\infty}(U)}
	\,.
\end{equation*}
Taking the maximum over~$|e|=1$ in~\eqref{e.J.localization.block} gives
\begin{equation}
	\mathcal J_{m+h}(U)
	\leq
	2\mathcal J_m(U)
	+3\X_{m,h}(U)^2\bigl(\mathcal J_m(U)+1\bigr)\,.
	\label{e.J.localization.block.maximum}
\end{equation}
The fields~$\a_{m+h}^{-1}\a_m$ and~$\a_m$ are independent; hence so are~$\X_{m,h}(\cu_{m+h})$ and~$\mathcal J_m(\cu_{m+h})$. The triangle inequality therefore gives
\begin{align}
	\|\mathcal J_{m+h}(\cu_{m+h})\|_{L^\xi(\P)}
	&
	\leq
	2\|\mathcal J_m(\cu_{m+h})\|_{L^\xi(\P)}
	\notag \\ & \qquad
	+
	3\|\X_{m,h}(\cu_{m+h})^2\|_{L^\xi(\P)}
	\bigl(1+\|\mathcal J_m(\cu_{m+h})\|_{L^\xi(\P)}\bigr)
	\,.
	\label{e.J.localization.block.moment}
\end{align}
By~\eqref{e.annealed.ordering.ratio} and~\eqref{e.aman.Linfty.moments}, applied with moment~$2\xi$, together with~\eqref{e.params.first.induction.step}, we have the estimate
\begin{equation*}
	\E\bigl[ \X_{m,h}(\cu_{m+h})^{2 \xi} \bigr]^{\nf 1 \xi} \leq C \xi \gmcdelta^2 h^2 \leq \frac{1}{12} \delta_1 \,.
\end{equation*}
Inserting this and~\eqref{e.gain-in-delta} into~\eqref{e.J.localization.block.moment} gives
\begin{equation*}
	\|\mathcal J_{m+h}(\cu_{m+h})\|_{L^\xi(\P)}
	\leq
	\frac12 \delta_1 + C\xi \gmcdelta^2 h^2
	\leq \delta_1
	\,.
\end{equation*}
The same argument, with~$(m,m+h)$ replaced by~$(k-h,k)$, gives the same bound on~$\cu_k$ for each~$k\in\{m+1,\ldots,m+h\}$; together with~$\S(m,\xi,\delta_1)$, this proves~\eqref{e.imp}. Hence~$\S(\infty,\xi,C\xi\log^2(2+\xi)\gmcdelta^2)$ holds for all~$\xi\geq16d$ in the stated range. The remaining range~$1\leq\xi<16d$ follows from monotonicity of~$L^\xi(\P)$ norms.

\smallskip

\emph{Step 3.} We prove~\eqref{e.homogenize.at.scale.m}. We now consider parameters~$\xi,\delta_1$ and~$h$ satisfying
\begin{equation}
	h \geq K \log(2+\xi)
	\quad \text{and} \quad
	\xi \gmcdelta^2  h \leq K^{-1} \delta_1 \leq K^{-2} \,,
	\label{e.params.second.induction.step}
\end{equation}
where~$K(d)$ is a large constant to be determined below. We also assume~$\xi\geq 4$.

\smallskip

For every~$m\in\N \cap [3 h,\infty)$, we prove the implication
\begin{equation}
	\S(m, \xi, \delta_1)
	\implies
	\S(m+ h, \xi,\delta_1)
	\,.
	\label{e.homogenization.induction.closure}
\end{equation}
By subadditivity,~\eqref{e.J.localization.block.maximum}, H\"older's inequality, and the independence of the fields defining~$\X_{m,h}$ and~$\mathcal J_m$, we have
\begin{align*}
	\|\mathcal J_{m+h}(\cu_{m+h})\|_{L^\xi(\P)}
	&
	\leq
	\E \Biggl[
	\biggl(
	\avsum_{z\in 3^m\Zd \cap \cu_{m+h}}
	\mathcal J_{m+h}(z+\cu_m)\biggr)^{\!\xi}
	\Biggr]^{\nf1\xi}
	\notag \\ &
	\leq
	3 \E \Biggl[
	\avsum_{z\in 3^m\Zd \cap \cu_{m+h}}
	\bigl(
	\X_{m,h}(z {+} \cu_m)
	\bigr)^{2\xi}
	\Biggr]^{\nf1\xi}
	\notag \\ & \qquad
	+
	3 \E \Biggl[
	\avsum_{z\in 3^m\Zd \cap \cu_{m+h}}
	\!\!\!
	\bigl( 1 +
	\X_{m,h}(z {+} \cu_m)^2
	\bigr)^{\xi}
	\Biggr]^{\nf1\xi}
	\notag \\ & \qquad\qquad \times
	\E \Biggl[
	\biggl(
	\avsum_{z\in 3^m\Zd \cap \cu_{m+h}}
	\mathcal J_m(z+\cu_m)^2
	\biggr)^{\!\nf \xi2}
	\Biggr]^{\nf1\xi}
	\notag \\ &
	\leq
	C \xi \gmcdelta^2 h
	+
	3
	(1 + C \xi \gmcdelta^2 h )
	\E \Biggl[
	\biggl(
	\avsum_{z\in 3^m\Zd \cap \cu_{m+h}}
	\mathcal J_m(z+\cu_m)^2
	\biggr)^{\!\nf \xi2}
	\Biggr]^{\nf1\xi}
	\,.
\end{align*}
Each summand depends only on the restriction of~$\a_m$ to its own cube, and~$\a_m$ has range of dependence at most~$\sqrt d\,3^m$. We can therefore color the cubes with~$C(d)$ colors so that the summands of each color are independent. Applying Rosenthal's inequality (Lemma~\ref{l.Rosenthal}) to each color as in~\eqref{e.appl.Rosenthal}, with~$\mathcal J_m^2$ in place of~$J$ and~$\nf\xi2$ in place of~$\xi$, and using~$\mathcal S(m,\xi,\delta_1)$ and the estimate~$\mathcal S(\infty,2,C\gmcdelta^2)$ proved in Step~2, we obtain
\begin{align*}
	\lefteqn{
	\E \Biggl[
	\biggl(
	\avsum_{z\in 3^m\Zd \cap \cu_{m+h}}
	\mathcal J_m(z+\cu_m)^2
	\biggr)^{\!\nf \xi2}
	\Biggr]^{\nf2\xi}
	} \qquad &
	\notag \\ &
	\leq
	C 3^{-\frac{d}{2} h}
	\xi^{\nf12}
	\E \bigl[\mathcal J_m(\cu_m)^4\bigr]^{\nf12}
	+ C
	3^{-\frac{d(\xi-2)}{\xi} h} \xi
	\E \bigl[\mathcal J_m(\cu_m)^\xi\bigr]^{\nf2\xi}
	\notag \\ &\qquad
	+
	\E \bigl[\mathcal J_m(\cu_m)^2\bigr]
	\notag \\ &
	\leq
	C 3^{-\frac{d}{2} h}
	\xi^{\nf12} \delta_1^2
	+ C 3^{-\frac{d}{2} h} \xi
	\delta_1^2
	+
	C \gmcdelta^4
	\notag \\ &
	\leq
	\frac1{36} \delta_1^2
	\,.
\end{align*}
Combining the previous two displays and using~$\delta_1\leq1$, we obtain
\begin{equation*}
	\|\mathcal J_{m+h}(\cu_{m+h})\|_{L^\xi(\P)}
	\leq
	C\xi\gmcdelta^2h
	+\frac12\bigl(1+C\xi\gmcdelta^2h\bigr)\delta_1
	\leq
	\frac12\delta_1+C\xi\gmcdelta^2h
	\leq \delta_1
	\,.
\end{equation*}
The last inequality follows from~\eqref{e.params.second.induction.step} after fixing~$K(d)$ sufficiently large. The same calculation, with~$(m,m+h)$ replaced by~$(k-h,k)$, gives the same bound on~$\cu_k$ for each~$k\in\{m+1,\ldots,m+h\}$; together with~$\S(m,\xi,\delta_1)$, this proves~\eqref{e.homogenization.induction.closure}.

\smallskip

The base case~\eqref{e.seed} and~\eqref{e.homogenization.induction.closure} imply that~$\S(\infty,\xi,C\xi\log(2+\xi)\gmcdelta^2)$ holds for~$\xi\geq4$ whenever~$\xi\log(2+\xi)\gmcdelta^2\leq C^{-1}$. The range~$1\leq\xi<4$ follows by monotonicity of~$L^\xi(\P)$ norms. This proves~\eqref{e.homogenize.at.scale.m}.

\smallskip

\emph{Step 4.} We prove~\eqref{e.homogenize.at.scale.m.E}. It follows from~\eqref{e.mathcal.E.Lms.moments.small.L}, applied with~$\delta_1=C\xi\log(2+\xi)\gmcdelta^2$.
\end{proof}

\section{Effective diffusivities}
\label{s.homogenized.diffusivities}

In this section we determine how fast the effective diffusivities~$\ahom_m$ decay.  We prove that they are nonincreasing in~$m$, that~$e^{-(m+1)\tau^2}\leq\ahom_m\leq3^{-\eta_{\mathrm{sub}}(m+1)}$ for every~$m\in\N_0$, with an exponent~$\eta_{\mathrm{sub}}>0$ depending only on~$d$ and the law of~$\g_0$, and that, at weak disorder,
\begin{equation*}
	\left|
	\log\ahom_m+\frac{2\tau^2(m+1)}d
	\right|
	\leq
	C\gmcdelta^2|\log\gmcdelta|(1+m\gmcdelta)\,.
\end{equation*}
In dimension two, planar duality gives the exact identity~$\ahom_m=e^{-\tau^2(m+1)}$.  These statements are collected in Theorem~\ref{t.renormalized.diffusivities}.  The decay of~$\ahom_m$ gives the time scale~$3^{2N}/\ahom_N$ in Theorem~\ref{t.A} and the exponent~$\eta$ in Theorem~\ref{t.A}(iii).

\smallskip

We prove the asymptotics at weak disorder through matching upper and lower bounds on~$\ahom_m/\ahom_{m-h}$ for integers~$h\leq\gmcdelta^{-1}$. We compare the corrector of~$\a_m$ with the solution of the linearized corrector equation. The corrector is the function added to an affine function to make it solve the equation with coefficient~$\a_m$. The linearized corrector solves a Poisson equation forced by the added fields, and we compute its energy exactly in Lemma~\ref{l.laplacian.corrector.energy}.  The error term in these asymptotics exceeds the main term unless~$\tau^2$ is bounded below by a power of~$\gmcdelta$.  We therefore prove the two-sided bound~$e^{-(m+1)\tau^2}\leq\ahom_m\leq3^{-\eta_{\mathrm{sub}}(m+1)}$ separately, by a nonperturbative argument which assumes no smallness of~$\gmcdelta$.  Its upper bound makes the process subdiffusive in every dimension, and the asymptotics identify the exponent at weak disorder.

\begin{theorem}[Bounds on the effective diffusivities]
\label{t.renormalized.diffusivities}
There exists~$\eta_{\mathrm{sub}}=\eta_{\mathrm{sub}}(d,\operatorname{law}(\g_0))>0$ such that, for every~$0\leq n\leq m$,
\begin{equation*}
	e^{-(m+1)\tau^2}
	\leq
	\ahom_m
	\leq
	\ahom_n
	\leq
	\min\{1,e^{2\tau^2(m-n)}\ahom_m\}\,,
	\qquad
	\ahom_m\leq3^{-\eta_{\mathrm{sub}}(m+1)}\,.
\end{equation*}
In dimension two, for every~$m\in\N_0$,
\begin{equation}
	\ahom_m=e^{-(m+1)\tau^2}\,.
	\label{e.renormalized.diffusivities.two.dimensions}
\end{equation}
Moreover, there exist~$\gmcdelta_0(d)>0$ and~$C(d)<\infty$ such that, if~$\gmcdelta\leq\gmcdelta_0(d)$, then, for every~$m\in\N_0$,
\begin{equation*}
	\left|\log\ahom_m+\frac{2\tau^2(m+1)}{d}\right|
	\leq
	C\gmcdelta^2|\log\gmcdelta|(1+m\gmcdelta)\,.
\end{equation*}
\end{theorem}

\begin{proof}
The monotonicity in~$m$ and the comparison between two scales are~\eqref{e.annealed.ordering}, and the affine competitor gives~$\ahom_0\leq1$.  The remaining assertions are Propositions~\ref{p.special.two.d.exact.formula},~\ref{p.sharp.asymptotic},~\ref{p.homogenized.coefficient.strict.decay}, and~\ref{p.homogenized.coefficient.reciprocal.lower}.
\end{proof}

\begin{proposition}[Asymptotics of the effective diffusivities at weak disorder]
	\label{p.sharp.asymptotic}
		There exist constants~$\gmcdelta_0(d)>0$ and~$C(d)<\infty$ such that, if~$\gmcdelta\leq\gmcdelta_0(d)$, then, for every~$m\in\N_0$,
	\begin{equation}
		\biggl|\log\ahom_m + \frac{2\tau^2 (m+1)}{d}\biggr| \leq C\gmcdelta^2\left|\log\gmcdelta\right|\bigl(1+m\gmcdelta\bigr)\,.
		\label{e.main.asymptotic}
	\end{equation}
\end{proposition}

\begin{lemma}[Energy of the linearized corrector]
\label{l.laplacian.corrector.energy}
Let~$m\in\N_0$ and~$h\in\N$ satisfy~$h\leq m+1$ and~$h\leq\gmcdelta^{-1}$, let~$X\coloneqq\prod_{k=m-h+1}^m\exp(\g_k-\tau^2)-1$, and let~$p\in\Rd$ be a unit vector.  If~$w_p^{(K)}\in H_0^1(\cu_K)$ solves
\begin{equation*}
	-\Delta w_p^{(K)}=\nabla\cdot(Xp)
	\quad\text{in }\cu_K\,,
\end{equation*}
then, with an implicit constant that depends only on~$d$,
\begin{equation}
	\lim_{K\to\infty}
	\E\bigl[\|\nabla w_p^{(K)}\|_{\underline L^2(\cu_K)}^2\bigr]
	=
	\frac{2\tau^2h}{d}+O(\gmcdelta^4h^2)\,.
	\label{e.laplacian.corrector.energy}
\end{equation}
The same conclusion holds if~$w_p^{(K)}$ is the mean-zero solution of
\begin{equation*}
	-\Delta w_p^{(K)}=-\nabla\cdot(Xp)
	\quad\text{in }\cu_K\,,
	\qquad
	\partial_{\mathbf n}w_p^{(K)}=Xp\cdot\mathbf n
	\quad\text{on }\partial\cu_K\,.
\end{equation*}
\end{lemma}
\begin{proof}
We show that the limit equals~$\E[|\nabla\phi(0)|^2]$, where~$\nabla\phi$ is the orthogonal projection of~$Xp$ onto the stationary gradients, and then compute this quantity. First, the spectral representation of~$X$ gives the decomposition~$Xp=\nabla\phi+\nabla\cdot S$ into a scalar and a skew potential, and both potentials grow sublinearly. Next, the potentials, cut off near~$\partial\cu_K$, serve as competitors in the variational formulas for the Dirichlet and the Neumann problems; the main point is that, by the sublinearity, the cutoffs change the energy by an amount that vanishes as~$K\to\infty$ and then~$L\to\infty$. Finally, the symmetry of the spectral measure gives~$\E[|\nabla\phi(0)|^2]=\frac1d\E[X(0)^2]$, and an expansion of~$\E[X(0)^2]$ in the last~$h$ fields completes the proof.

\smallskip

The vector products and unit probes in this proof are Euclidean, and every constant~$C$ and every implicit constant depends only on~$d$. We construct the scalar and skew potentials, with row divergence~$(\nabla\cdot S)_i=\sum_{j=1}^d\partial_jS_{ij}$. The field~$X$ is centered, stationary, mean-square continuous and square integrable. Its finite range of dependence makes its covariance compactly supported and integrable. Its spectral measure~$\mu$ is therefore finite and has no atom at zero. Write the real-field spectral representation as~$X(x)=\int e^{ik\cdot x}Z(dk)$, with~$Z(-dk)=\overline{Z(dk)}$, where~$\E|Z(dk)|^2=\mu(dk)$ and~$\mu(\Rd)=\E[X(0)^2]$. Define potentials anchored at the origin by
\begin{equation*}
	\phi(x)=\int(e^{ik\cdot x}-1)\frac{-i\,p\cdot k}{|k|_2^2}Z(dk)\,,
	\qquad
	S_{ij}(x)=\int(e^{ik\cdot x}-1)\frac{-i(p_i k_j-p_j k_i)}{|k|_2^2}Z(dk)\,.
\end{equation*}
The multipliers are understood to be zero at~$k=0$. Their squared magnitudes, after multiplication by~$|e^{ik\cdot x}-1|^2$, are bounded by~$C|p|_2^2\min\{|x|_2^2,|k|_2^{-2}\}$. The derivative multipliers are bounded. Taking jointly measurable versions, we thus have~$\phi,S\in L^2(\P;H^1_{\mathrm{loc}})$, with stationary square-integrable gradients. Differentiation gives the multipliers
\begin{equation*}
	\widehat{\nabla\phi}=\frac{k(p\cdot k)}{|k|_2^2}\widehat X\,,
	\qquad
	\widehat{\nabla\cdot S}=\left(p-\frac{k(p\cdot k)}{|k|_2^2}\right)\widehat X\,.
\end{equation*}
Consequently~$Xp=\nabla\phi+\nabla\cdot S$,~$S=-S^{\mathrm T}$, and the two vector fields are orthogonal in~$L^2(\P;\Rd)$. Cutting the scalar multiplier off near~$k=0$ approximates~$\nabla\phi$ by stationary gradients. The complementary multiplier is divergence free and is orthogonal to stationary gradients by stationary integration by parts. Thus~$\nabla\phi$ is the orthogonal projection of~$Xp$ onto their closure.

\smallskip

For a centered cube~$R\cu_0$ of side~$R$, since
\begin{equation*}
	\fint_{R\cu_0}\left|e^{ik\cdot x}-\fint_{R\cu_0}e^{ik\cdot y}\,dy\right|^2\,dx
	=1-\left|\fint_{R\cu_0}e^{ik\cdot y}\,dy\right|^2
	\leq\min\{1,R^2|k|_2^2/12\}\,,
\end{equation*}
the spectral isometry and~$\sum_{i,j}|p_i k_j-p_j k_i|^2\leq2|p|_2^2|k|_2^2$ imply
\begin{equation*}
	R^{-2}\E\fint_{R\cu_0}\left(|\phi-(\phi)_{R\cu_0}|^2+|S-(S)_{R\cu_0}|^2\right)
	\leq3|p|_2^2\int\min\{R^{-2}|k|_2^{-2},1/12\}\,\mu(dk)
	\longrightarrow0\,.
\end{equation*}
Here the squared spectral matrix norm is bounded by the sum of the squared entries used in the calculation. Dominated convergence applies because~$\mu$ is finite and~$\mu(\{0\})=0$. At~$R=3^K$ this proves
\begin{equation}
	\lim_{K\to\infty}3^{-2K}\E\left[
	\fint_{\cu_K}\left(
	|\phi-(\phi)_{\cu_K}|^2+|S-(S)_{\cu_K}|^2
	\right)\right]=0\,.
	\label{e.stationary.helmholtz.sublinearity}
\end{equation}
The argument applies in dimension two as well: it requires stationary gradients and anchored increments, and makes no finite-variance stationarity assertion about the values of the potentials.

\smallskip

Fix~$L\in\N$.  Choose~$\eta_{K,L}\in C_c^\infty(\cu_K)$ equal to one at distance at least~$3^{K-L}$ from~$\partial\cu_K$ and satisfying~$0\leq\eta_{K,L}\leq1$ and~$|\nabla\eta_{K,L}|\leq C3^{L-K}$.  The support of~$1-\eta_{K,L}$ has relative volume at most~$C3^{-L}$.  Consequently, by stationarity and~\eqref{e.stationary.helmholtz.sublinearity}, the function~$-\eta_{K,L}(\phi-(\phi)_{\cu_K})\in H_0^1(\cu_K)$ satisfies
\begin{equation*}
	\limsup_{K\to\infty}
	\E\bigl[\|\nabla\phi-\nabla\bigl(\eta_{K,L}(\phi-(\phi)_{\cu_K})\bigr)\|_{
	\underline L^2(\cu_K)}^2\bigr]
	\leq C3^{-L}\E[|\nabla\phi(0)|^2]\,,
\end{equation*}
and, for the same reasons, the field
\begin{equation*}
	\mathbf j_{K,L}\coloneqq
	\nabla\cdot\left(
	\eta_{K,L}\bigl(S-(S)_{\cu_K}\bigr)
	\right)
	\in\Lsolo(\cu_K)
\end{equation*}
satisfies
\begin{equation}
	\limsup_{K\to\infty}
	\E\bigl[\|\mathbf j_{K,L}-\nabla\cdot S\|_{
	\underline L^2(\cu_K)}^2\bigr]
	\leq C3^{-L}\E[|\nabla\cdot S(0)|^2]\,.
	\label{e.stationary.solenoidal.approximation}
\end{equation}

\smallskip

For the Dirichlet problem, the variational identity
\begin{equation*}
	\frac12\|\nabla w_p^{(K)}\|_{\underline L^2(\cu_K)}^2
	=
	\sup_{\varphi\in H_0^1(\cu_K)}
	\fint_{\cu_K}
	\left(-Xp\cdot\nabla\varphi-\frac12|\nabla\varphi|^2\right)
\end{equation*}
and the orthogonality of~$\nabla\cdot S$ to compactly supported gradients give
\begin{equation*}
	\E\bigl[\|\nabla w_p^{(K)}\|_{
	\underline L^2(\cu_K)}^2\bigr]
	\leq \E[|\nabla\phi(0)|^2]\,.
\end{equation*}
Using~$\varphi=-\eta_{K,L}(\phi-(\phi)_{\cu_K})$ in the variational identity, then sending~$K\to\infty$ and~$L\to\infty$ in this order, proves the reverse inequality.

\smallskip

For the Neumann problem, let~$w_p^{(K)}$ now denote the mean-zero solution. Its weak formulation says that~$\nabla w_p^{(K)}$ is the orthogonal projection of~$Xp$ onto~$\nabla H^1(\cu_K)$.  Using the restriction of~$\phi$ to~$\cu_K$ as a competitor in the variational formula for the norm of this projection, and using~$\E[(Xp)(0)\cdot\nabla\phi(0)]=\E[|\nabla\phi(0)|^2]$, we obtain
\begin{equation*}
	\E\bigl[\|\nabla w_p^{(K)}\|_{
	\underline L^2(\cu_K)}^2\bigr]
	\geq \E[|\nabla\phi(0)|^2]\,.
\end{equation*}
On the other hand, the field~$\mathbf j_{K,L}\in\Lsolo(\cu_K)$ is orthogonal to~$\nabla w_p^{(K)}$, and hence
\begin{equation*}
	\|\nabla w_p^{(K)}\|_{\underline L^2(\cu_K)}^2
	\leq
	\|Xp-\mathbf j_{K,L}\|_{\underline L^2(\cu_K)}^2\,.
\end{equation*}
Sending first~$K\to\infty$ and then~$L\to\infty$ and using~\eqref{e.stationary.solenoidal.approximation} proves the reverse inequality.  Consequently, in both cases,
\begin{equation*}
	\lim_{K\to\infty}
	\E\bigl[\|\nabla w_p^{(K)}\|_{\underline L^2(\cu_K)}^2\bigr]
	=
	\E\bigl[|\nabla\phi(0)|^2\bigr]\,.
\end{equation*}
The spectral measure of the scalar stationary field~$X$ is invariant under coordinate reflections and permutations.  The off-diagonal entries of the Fourier multiplier of the projection therefore integrate to zero against this measure, while the diagonal entries have equal integrals, which sum to~$\E[X(0)^2]$.  Hence~$\E[|\nabla\phi(0)|^2]=\frac1d\E[X(0)^2]$. Independence of the fields and the definition~$\tau^2=\log\E[\exp(\g_0(0))]$ yield~$\E[X(0)^2]=\exp(h(\log\E[\exp(2\g_0(0))]-2\tau^2))-1$. The symmetry and sub-Gaussian bound in Assumptions~\ref{a.g2}--\ref{a.g3} give~$\log\E[\exp(2\g_0(0))]-2\tau^2=2\tau^2+O(\gmcdelta^4)$. Since~$h\gmcdelta^2\leq\gmcdelta$, expansion of the exponential proves~\eqref{e.laplacian.corrector.energy}.
\end{proof}

\begin{lemma}[Upper bound on~$\ahom_m/\ahom_{m-h}$]
	\label{l.one.step.upper}
	There exists~$C(d) < \infty$ such that, for every~$m,h\in\N$ with
	\begin{equation*}
		h\leq\gmcdelta^{-1}
		\qquad\text{and}\qquad
		m-h\geq16\lceil|\log_3\gmcdelta|\rceil\,,
	\end{equation*}
	we have
	\begin{equation}
		\ahom_m \leq \ahom_{m-h}  \biggl(1 - \frac{2 \tau^2 h}{d} + C\gmcdelta^4 h^{2}  + C \gmcdelta^2 \left|\log \gmcdelta\right| \biggr)  \,.
		\label{e.one.step.upper}
	\end{equation}
\end{lemma}

\begin{proof}
We first record a bound on~$\tau^2$, settle the planar case by the exact formula, fix the moment orders, and treat the range~$\gmcdelta>\gmcdelta_*$, where the ordering~\eqref{e.annealed.ordering} suffices. The main argument, for~$d\geq3$ and~$\gmcdelta\leq\gmcdelta_*$, bounds~$\ahom_m$ from above by the Dirichlet principle, with a competitor built from the linearized corrector; it runs in four steps, outlined before Step~1, and its heart is Step~3, where Lemma~\ref{l.laplacian.corrector.energy} gives the main term~$-2\tau^2h/d$.

\smallskip

Symmetry and Assumption~\ref{a.g2}, followed by Jensen with exponent~$\gmcdelta^2/2\leq1$, give
\begin{equation*}
	e^{\tau^2}=\E\cosh(\g_0(0))
	\leq\E e^{\g_0(0)^2/2}
	\leq\left(\E e^{\g_0(0)^2/\gmcdelta^2}\right)^{\gmcdelta^2/2}
	\leq2^{\gmcdelta^2/2}\,.
\end{equation*}
Thus~$\tau^2\leq(\log2/2)\gmcdelta^2$. First consider~$d=2$. Proposition~\ref{p.special.two.d.exact.formula} gives~$\ahom_m/\ahom_{m-h}=e^{-\tau^2h}$. The bound~$\tau^2\leq(\log2/2)\gmcdelta^2$ and~$h\gmcdelta\leq1$ permit the estimate~$e^{-x}\leq1-x+x^2$ with~$x=\tau^2h$, which proves~\eqref{e.one.step.upper}. Henceforth~$d\geq3$.

\smallskip

We fix all moment orders before restricting the disorder. Let~$c_P,C_P$ be constants for Proposition~\ref{p.coarse.grained.bound}, let~$s=\nf14$ and~$\xi_*\coloneqq\lceil\max\{16d,e^{4c_P}\}\rceil$. Choose a dimensional~$\gmcdelta_*\in(0,e^{-1}]$ so that
\begin{equation*}
	C_P\xi_*\gmcdelta_*\leq1\,,\qquad
	C_P\xi_*\log(2+\xi_*)\gmcdelta_*^2\leq c_P/4\,,\qquad
	16C_P^2\xi_*\log(2+\xi_*)\gmcdelta_*^2\leq1\,.
\end{equation*}
For~$0<\gmcdelta\leq\gmcdelta_*$ these ensure every condition of the proposition at~$(\xi_*,s)$: indeed~$\gmcdelta^2|\log\gmcdelta|\leq\gmcdelta$,~$4d/s\leq\xi_*$ and~$s\log(2+\xi_*)\geq c_P$. Its first part also applies at moment~$1$. The field moments used below have fixed orders at most~$16$ and follow from Assumption~\ref{a.g2}; they impose no scale-dependent restriction.

\smallskip

For the complementary range~$\gmcdelta>\gmcdelta_*$, the ordering~\eqref{e.annealed.ordering} gives~$\ahom_m\leq\ahom_{m-h}$. Since~$h\geq1$,
\begin{equation*}
	\gmcdelta^2h\leq\gmcdelta_*^{-2}\gmcdelta^4h^2\,,
	\qquad
	\frac{2\tau^2h}{d}\leq\frac{\log2}{d}\gmcdelta_*^{-2}\gmcdelta^4h^2\,.
\end{equation*}
Increasing~$C(d)$ proves the conclusion in this range. We may therefore assume~$\gmcdelta\leq\gmcdelta_*$ throughout the remaining proof.

\smallskip

	We bound~$\E[p\cdot\a_m(\cu_K)p]$, for a unit vector~$p$ and large~$K$, from above by the Dirichlet principle. The competitor is~$p\cdot x+w$, where~$w$ solves a Poisson equation forced by the last~$h$ fields, corrected on each interior cell of the scale~$3^n$ fixed in Step~1 by solutions of cell problems for~$\a_m$. Step~1 sets up~$w$ and the cells. Step~2 bounds the energy of the competitor by the sum of the cell energies, up to an error~$C\gmcdelta^{15}\ahom_{m-h}$. Step~3 uses the independence of the last~$h$ fields from~$\a_{m-h}$ to bound this sum, up to a further error of the same size, by~$\ahom_{m-h}(\cu_n)$ times an expectation involving only these fields, and evaluates this expectation. This is the heart of the matter: Lemma~\ref{l.laplacian.corrector.energy} gives the main term~$-2\tau^2h/d$, and reversing the signs of the last~$h$ fields shows that the remaining error is~$O(\gmcdelta^4h^2)$. Step~4 concludes.

\smallskip

	\emph{Step 1.} We set up the Dirichlet problem and record the interior estimates on the \textbf{cells} of scale~$3^n$, the cubes~$z+\cu_n$ with~$z\in3^n\Zd$. Fix~$p\in\Rd$ with~$|p|_2=1$, and set
	\begin{equation}
		\h\coloneqq\sum_{k=m-h+1}^{m}(\g_k-\tau^2)\,.
		\label{e.h.def}
	\end{equation}
	Choose
	\begin{equation}
		n\coloneqq m-h-16\lceil|\log_3\gmcdelta|\rceil\,,
		\label{e.parameter.choice}
	\end{equation}
	so that~$3^{-(m-h-n)}\leq\gmcdelta^{16}$. For~$K>m$, let~$w\in H_0^1(\cu_K)$ solve
	\begin{equation*}
		-\Delta w=\nabla\cdot((e^\h-1)p)
		\quad\text{in }\cu_K\,.
	\end{equation*}
	Testing the equation with~$w$ and applying the gradient Calder\'on--Zygmund estimate at exponent~$8$ on the cube (see the reflection argument below), followed by H\"older's inequality for the lower exponents, give, for~$1\leq\xi\leq8$,
	\begin{equation*}
		\|\nabla w\|_{L^\xi(\P;\underline L^\xi(\cu_K))}
		\leq C\gmcdelta h^{\nf{1}{2}}\,,
		\qquad
		\|\nabla w\|_{\underline L^2(\cu_K)}^2
		=
		-\fint_{\cu_K}(e^\h-1)p\cdot\nabla w\,.
	\end{equation*}
	Lemma~\ref{l.laplacian.corrector.energy} yields
	\begin{equation*}
		\lim_{K\to\infty}
		\E\|\nabla w\|_{\underline L^2(\cu_K)}^2
		=
		\frac{2\tau^2h}{d}+O(\gmcdelta^4h^2)\,.
	\end{equation*}

\smallskip

	For every scale-$3^n$ cell~$z+\cu_n$ satisfying~$z+\cu_{m-h}\subseteq\cu_K$, define
	\begin{equation*}
		p_z=(p+\nabla w)_{z+\cu_n}\,,
		\qquad
		F_z=\nabla w-(\nabla w)_{z+\cu_n}\,,
		A_z=\|F_z\|_{\underline L^2(z+\cu_n)}\,,
		\qquad
		B_z=3^n\|\nabla F_z\|_{\underline L^2(z+\cu_n)}\,.
	\end{equation*}
    Odd reflection of the zero-Dirichlet solution across the cube faces, with the forcing field~$(e^\h-1)p$ reflected so that its normal component is even and its tangential components are odd, reduces the equation to a periodic Poisson equation. Its right-hand side~$\nabla\cdot((e^\h-1)p)=p\cdot\nabla e^\h$ is a scalar function, and the reflected equation has the odd reflection of that function as right-hand side. The periodic Calder\'on--Zygmund estimates, followed by scaling, give the gradient estimate above and the following Hessian estimate with constants independent of~$K$. In particular, the weak Hessian of~$w$ is in~$L^4(\cu_K)$ and~$\|D^2w\|_{\underline L^4(\cu_K)}\leq C\|\nabla e^\h\|_{\underline L^4(\cu_K)}$. At each point, the geometric scaling of the derivatives of the last~$h$ fields gives~$\|\nabla\h\|_{L^8(\P)}\leq C\gmcdelta3^{-(m-h)}$, while~$\E e^{8\h}\leq C$ because~$h\gmcdelta^2\leq\gmcdelta$. H\"older's inequality and stationarity consequently give~$(\E\fint_{\cu_K}|\nabla e^\h|^4)^{\nf14}\leq C\gmcdelta3^{-(m-h)}$. We use this estimate on the finite family of cells, with no spatial maximum of the Hessian. Jensen's inequality on each cell and the partition into scale-$3^n$ cubes give
\begin{equation*}
	\limsup_{K\to\infty}\avsum_z\E[B_z^4]
	\leq C3^{4n}\sup_K\E\fint_{\cu_K}|D^2w|^4
	\leq C\gmcdelta^{68}\,.
\end{equation*}
The passage from all cells to the retained cells contributes their volume ratio, which tends to one. Mean-zero Poincar\'e on~$z+\cu_n$ gives~$A_z\leq CB_z$. In particular,
\begin{equation*}
	\limsup_{K\to\infty}\left(\avsum_z\E[A_z^4]\right)^{\nf14}
	+\limsup_{K\to\infty}\left(\avsum_z\E[B_z^4]\right)^{\nf14}
	\leq C\gmcdelta^{17}\leq C\gmcdelta^{15}\,.
\end{equation*}
The Hessians and the cell norms have measurable versions: the Poisson solution depends continuously on its datum, and difference quotients converge to its weak derivatives. Here and below the averages run over the retained interior cells. The omitted cells lie within distance~$C3^{m-h}$ of~$\partial\cu_K$, and their relative volume tends to zero as~$K\to\infty$.

	\smallskip
	\emph{Step 2.} We solve cell problems on the interior cells, patch their solutions into a competitor for the Dirichlet principle, and bound the energy by a sum over cells. On~$z+\cu_n$, let~$\chi_z,\phi_z\in H_0^1(z+\cu_n)$ solve
	\begin{equation*}
		-\nabla\cdot(\a_m\nabla\chi_z)
		=
		\nabla\cdot(\a_mp_z)\,,
		\qquad
		-\nabla\cdot(\a_m\nabla\phi_z)
		=
		\nabla\cdot(\a_mF_z)\,,
	\end{equation*}
	and set~$S_z=p_z+\nabla\chi_z$ and~$T_z=F_z+\nabla\phi_z$. Then~$\nabla\cdot(\a_mS_z)= \nabla\cdot(\a_mT_z)=0$,~$(T_z)_{z+\cu_n}=0$, and~$\fint_{z+\cu_n}|S_z|^2\a_m=p_z\cdot\a_m(z+\cu_n)p_z$. Testing the Dirichlet principle with the function equal to~$p\cdot x+w+\chi_z+\phi_z$ on each interior cell~$z+\cu_n$ and to~$p\cdot x+w$ on the omitted cells gives
	\begin{equation*}
		\frac12\E[p\cdot\a_m(\cu_K)p]
		\leq
		\avsum_z\E\left[
		\frac12p_z\cdot\a_m(z+\cu_n)p_z
		+
		\fint_{z+\cu_n}S_z\cdot\a_mT_z
		+
		\frac12\fint_{z+\cu_n}|T_z|^2\a_m
		\right]
		+o_K(1)\,.
	\end{equation*}

\smallskip

	Since~$\phi_z$ minimizes the energy of~$F_z+\nabla\phi$ and~$(F_z)_{z+\cu_n}=0$, we have~$\fint_{z+\cu_n}|T_z|^2\a_m=\fint_{z+\cu_n}F_z\cdot\a_mT_z$. Apply~\eqref{e.duality.splitting},~\eqref{e.besov.flux.poincare}, and
	\begin{equation*}
		3^{\frac{1}{4}n}[F_z]_{\underline B^{\nf{1}{4}}_{2,\infty}(z+\cu_n)}
		\leq
		C\bigl(A_z+A_z^{\nf{3}{4}}B_z^{\nf{1}{4}}\bigr)
	\end{equation*}
	to obtain
	\begin{equation*}
		\fint_{z+\cu_n}|T_z|^2\a_m
		\leq
		C\Lambda_{1/4,1}(z+\cu_n;\a_m)
		\bigl(A_z^2+A_z^{\nf{3}{2}}B_z^{\nf{1}{2}}\bigr)\,.
	\end{equation*}
	We verify the coefficient moment used here. The fixed choice~$(\xi_*,s)$ above and stationarity give~$\|\mathcal E_{1/4,\infty}(z+\cu_n;\a_n,\ahom_n)\|_{L^8(\P)}\leq1$. The squared bound~\eqref{e.bound.Lambdas.by.Es}, with~$q=1$, therefore gives
\begin{equation*}
	\|\Lambda_{1/4,1}(z+\cu_n;\a_n)\|_{L^4(\P)}\leq C\ahom_n\,,
	\qquad
	\|\lambda_{1/4,1}^{-1}(z+\cu_n;\a_n)\|_{L^4(\P)}\leq C\ahom_n^{-1}\,.
\end{equation*}
Applying the pathwise comparison~\eqref{e.coarse-grained-matrix-inequalities.general} on every subcube in the definitions of the coarse ellipticities gives
\begin{equation*}
	\Lambda_{1/4,1}(z+\cu_n;\a_m)
	\leq(\sup_{z+\cu_n}\a_m\a_{m-h}^{-1})(\sup_{z+\cu_n}\a_{m-h}\a_n^{-1})\Lambda_{1/4,1}(z+\cu_n;\a_n)\,,
\end{equation*}
and the same comparison for~$\lambda_{1/4,1}^{-1}$ with the reciprocal ratios~$\a_{m-h}\a_m^{-1}$ and~$\a_n\a_{m-h}^{-1}$. The field estimate~\eqref{e.aman.Linfty.moments}, on a cube of scale~$n$, bounds both ratios in~$L^8(\P)$ by~$C(d)$: their block lengths are~$h\leq\gmcdelta^{-1}$ and~$m-h-n\leq C|\log\gmcdelta|$. H\"older with~$\nf18+\nf18+\nf14=\nf12$ and~\eqref{e.annealed.ordering} yield
\begin{equation*}
	\|\Lambda_{1/4,1}(z+\cu_n;\a_m)\|_{L^2(\P)}\leq C\ahom_{m-h}\,,
	\qquad
	\|\lambda_{1/4,1}^{-1}(z+\cu_n;\a_m)\|_{L^2(\P)}\leq C\ahom_{m-h}^{-1}\,.
\end{equation*}
This comparison is simultaneous over the subcubes and requires no maximum of separate moment estimates. To combine the coefficient and cell bounds, take H\"older on the product of~$\P$ with the uniform measure on the finite retained family. In particular,
\begin{multline*}
	\avsum_z\E[\Lambda_{1/4,1}(z+\cu_n;\a_m)(A_z^2+A_z^{\nf32}B_z^{\nf12})]
	\\
	\leq\left(\avsum_z\E[\Lambda_{1/4,1}^2(z+\cu_n;\a_m)]\right)^{\nf12}
	\left(\avsum_z\E[(A_z^2+A_z^{\nf32}B_z^{\nf12})^2]\right)^{\nf12}\,.
\end{multline*}
The last factor is bounded by~$C(\avsum_z\E[A_z^4+B_z^4])^{\nf12}$. Thus the estimates of Step~1 imply
	\begin{equation*}
		\limsup_{K\to\infty}\avsum_z
		\E\fint_{z+\cu_n}|T_z|^2\a_m
		\leq C\gmcdelta^{30}\ahom_{m-h}\,.
	\end{equation*}
	The same coefficient bounds give~$\avsum_z\E\fint_{z+\cu_n}|S_z|^2\a_m\leq C\ahom_{m-h}$. Cauchy--Schwarz therefore yields
	\begin{equation*}
		\limsup_{K\to\infty}\avsum_z
		\E\left|
		\fint_{z+\cu_n}S_z\cdot\a_mT_z
		\right|
		\leq C\gmcdelta^{15}\ahom_{m-h}\,.
	\end{equation*}
	With the bounds on the~$T_z$ terms and on the cross terms, the estimate obtained from the Dirichlet principle above bounds the energy of the competitor by the sum of the cell energies, up to an error~$C\gmcdelta^{15}\ahom_{m-h}$ in the limit~$K\to\infty$.

	\smallskip
	\emph{Step 3.} We estimate the sum of the cell energies~$p_z\cdot\a_m(z+\cu_n)p_z$. Since~$\a_m=e^\h\a_{m-h}$,
	\begin{equation*}
		p_z\cdot\a_m(z+\cu_n)p_z
		\leq
		\left(\sup_{z+\cu_n}e^\h\right)
		p_z\cdot\a_{m-h}(z+\cu_n)p_z\,.
	\end{equation*}
	We have the deterministic estimate
	\begin{equation*}
		\left|
		\left(\sup_{z+\cu_n}e^\h\right)|p_z|^2
		-
		\fint_{z+\cu_n}e^\h|p+\nabla w|^2
		\right|
		\leq
		Ce^{\|\h\|_{L^\infty(z+\cu_n)}}
		\left(
		\osc_{z+\cu_n}\h\,|p_z|^2+(|p_z|+A_z)A_z
		\right)\,.
	\end{equation*}
	The choice of~$n$ gives~$\|\osc_{z+\cu_n}\h\|_{L^8(\P)}\leq C\gmcdelta^{17}$. The exponential factor and~$p_z$ have bounded eighth moments, while Step~1 controls~$A_z$. Combining these bounds with the independence of~$\h$ and~$\a_{m-h}$, we obtain
	\begin{equation*}
		\limsup_{K\to\infty}\avsum_z
		\E[p_z\cdot\a_m(z+\cu_n)p_z]
		\leq
		\ahom_{m-h}(\cu_n)
		\limsup_{K\to\infty}
		\E\fint_{\cu_K}e^\h|p+\nabla w|^2
		+
		C\gmcdelta^{15}\ahom_{m-h}\,.
	\end{equation*}

\smallskip

	It remains to evaluate the expectation in the last display. The equation for~$w$, the zero trace of~$w$, and~$\E[e^\h]=1$ give
	\begin{equation*}
		\E\fint_{\cu_K}e^\h|p+\nabla w|^2
		=
		1-\E\|\nabla w\|_{\underline L^2(\cu_K)}^2
		+
		\E\fint_{\cu_K}|\nabla w|^2(e^\h-1)\,.
	\end{equation*}
	We now show that the last term is~$O(\gmcdelta^4h^2)$, uniformly in~$K$. Let~$v\in H_0^1(\cu_K)$ solve the first equation below; then~$w-v\in H_0^1(\cu_K)$ solves the second:
	\begin{equation*}
		-\Delta v=\nabla\cdot(\h p)\,,
		\qquad
		-\Delta(w-v)
		=
		\nabla\cdot((e^\h-1-\h)p)\,.
	\end{equation*}
	The Calder\'on--Zygmund and exponential-moment estimates give
	\begin{equation*}
		\|\nabla v\|_{L^4(\P;\underline L^4)}
		\leq C\gmcdelta h^{\nf{1}{2}}\,,
		\qquad
		\|\nabla(w-v)\|_{L^4(\P;\underline L^4)}
		\leq C\gmcdelta^2h\,.
	\end{equation*}
	Since constants have zero divergence,~$v$ is linear in the sum~$\h+h\tau^2=\sum_{k=m-h+1}^m\g_k$ of the last~$h$ fields. Reversing the signs of these fields simultaneously therefore gives
	\begin{equation*}
		\E\fint_{\cu_K}|\nabla v|^2(\h+h\tau^2)=0\,,
		\qquad
		\E\fint_{\cu_K}|\nabla v|^2\h
		=
		-h\tau^2\E\|\nabla v\|_{\underline L^2}^2
		=
		O(\gmcdelta^4h^2)\,.
	\end{equation*}
	The preceding moment bounds control the terms containing~$\nabla(w-v)$, and~$|e^s-1-s|\leq\frac12s^2e^{|s|}$ controls the nonlinear remainder. H\"older's inequality and~$\gmcdelta h\leq1$ thus give
	\begin{equation*}
		\sup_K\left|
		\E\fint_{\cu_K}|\nabla w|^2(e^\h-1)
		\right|
		\leq C\gmcdelta^4h^2\,.
	\end{equation*}
	Lemma~\ref{l.laplacian.corrector.energy} now yields
	\begin{equation*}
		\limsup_{K\to\infty}
		\E\fint_{\cu_K}e^\h|p+\nabla w|^2
		\leq
		1-\frac{2\tau^2h}{d}+C\gmcdelta^4h^2\,.
	\end{equation*}

	\smallskip
	\emph{Step 4.} We conclude the proof of~\eqref{e.one.step.upper}. Combining Steps~2 and~3 and passing to the thermodynamic limit gives
	\begin{equation*}
		\frac12\ahom_m
		\leq
		\left(
		\frac12-\frac{\tau^2h}{d}+C\gmcdelta^4h^2
		\right)
		\ahom_{m-h}(\cu_n)
		+
		C\gmcdelta^{15}\ahom_{m-h}\,.
	\end{equation*}
    For the finite-volume readout we use the first moment, rather than a bound on the expected square root of the defect. For a Euclidean unit vector~$e$,~\eqref{e.coarse.grained.matrices} in expectation gives
\begin{equation*}
	2\E J(\cu_n,\ahom_n^{-\nf12}e,\ahom_n^{\nf12}e;\a_n)
	=\left(\frac{\ahom_n(\cu_n)}{\ahom_n}-1\right)
	+\left(\ahom_n\ahom_{n,*}^{-1}(\cu_n)-1\right)\,.
\end{equation*}
Both summands are nonnegative by the monotone limits~\eqref{e.ahom.ahomstar.defs}. Proposition~\ref{p.coarse.grained.bound} at moment~$1$ bounds each by~$C\gmcdelta^2$. The finite-volume ordering~\eqref{e.annealed.matrix.bounds} and~\eqref{e.annealed.ordering} then give
\begin{equation*}
	\ahom_{m-h}\leq\ahom_{m-h}(\cu_n)
	\leq\ahom_n(\cu_n)
	\leq(1+C\gmcdelta^2)\ahom_n
	\leq(1+C\gmcdelta^2)e^{2\tau^2(m-h-n)}\ahom_{m-h}\,.
\end{equation*}
Since~$\tau^2(m-h-n)\leq C\gmcdelta^2|\log\gmcdelta|$, this proves
\begin{equation*}
	\ahom_{m-h}\leq\ahom_{m-h}(\cu_n)
	\leq(1+C\gmcdelta^2|\log\gmcdelta|)\ahom_{m-h}\,.
\end{equation*}
The preceding bracket is bounded in absolute value by~$C(d)$. For substitution, one may use~$Bx\leq Ba+|B|(x-a)$ for~$x\geq a$; no sign assumption on the bracket is needed. Also~$\gmcdelta^{15}\leq\gmcdelta^2|\log\gmcdelta|$ in the present range. Substituting this into the previous display and using~$\tau^2\leq C\gmcdelta^2$ and~$\gmcdelta h\leq1$ proves~\eqref{e.one.step.upper}.
\end{proof}

\begin{lemma}[Lower bound on~$\ahom_m/\ahom_{m-h}$]
	\label{l.one.step.lower}
	There exists~$C(d) < \infty$ such that, for every~$m,h\in\N$ with
	\begin{equation*}
		h\leq\gmcdelta^{-1}
		\qquad\text{and}\qquad
		m-h\geq16\lceil|\log_3\gmcdelta|\rceil\,,
	\end{equation*}
	we have
	\begin{equation}
		\ahom_m^{-1} \leq \ahom_{m-h}^{-1}  \biggl(1 + \frac{2 \tau^2 h}{d} + C\gmcdelta^4 h^{2}  + C \gmcdelta^2 \left|\log \gmcdelta\right| \biggr)  \,.
		\label{e.one.step.lower}
	\end{equation}
\end{lemma}

\begin{proof}
In dimension two, Proposition~\ref{p.special.two.d.exact.formula} gives~$\ahom_{m-h}/\ahom_m=e^{\tau^2h}$. With~$x=\tau^2h\leq(\log2/2)\gmcdelta\leq1/2$, the estimate~$e^x\leq1+x+2x^2$ proves~\eqref{e.one.step.lower}. Henceforth~$d\geq3$. Fix the moment orders and dimensional threshold~$\gmcdelta_*$ as in the upper-bound proof. If~$\gmcdelta>\gmcdelta_*$,~\eqref{e.annealed.ordering} gives~$\ahom_{m-h}/\ahom_m\leq e^{2\tau^2h}$. For~$t=\gmcdelta^2h\in[0,1]$, convexity gives~$e^{(\log2)t}\leq1+t$. Since~$t\leq\gmcdelta_*^{-2}\gmcdelta^4h^2$, increasing~$C(d)$ proves the conclusion in this range. We may therefore assume~$\gmcdelta\leq\gmcdelta_*$ without adding a hypothesis to the lemma.

\smallskip

The proof is dual to that of Lemma~\ref{l.one.step.upper}. We bound~$\E[q\cdot\a_{m,*}^{-1}(\cu_K)q]$, for a unit vector~$q$ and large~$K$, from above by the dual variational principle. The competitor is built from a divergence-free field~$\mathbf z$, defined through the solution~$W$ of a Neumann problem forced by the last~$h$ fields: on each interior cell of the scale~$3^n$ fixed in Step~1, it is the flux of the solution of the Neumann problem for~$\a_m$ whose normal flux on the cell boundary is that of~$\mathbf z$. Step~1 sets up~$W$,~$\mathbf z$, and the cells. Step~2 bounds the dual energy of the competitor by the sum of the cell energies, up to an error~$C\gmcdelta^{15}\ahom_{m-h}^{-1}$. Step~3 bounds this sum, up to a further error of the same size, by~$\ahom_{m-h,*}^{-1}(\cu_n)$ times an expectation involving only the last~$h$ fields, and evaluates this expectation. As in the upper bound, this is the heart of the matter: Lemma~\ref{l.laplacian.corrector.energy} gives the main term~$2\tau^2h/d$, and reversing the signs of the last~$h$ fields shows that the remaining error is~$O(\gmcdelta^4h^2)$. Step~4 concludes.

\smallskip

\emph{Step 1.} We set up the Neumann problem and record the interior estimates on the scale-$3^n$ cells. Fix~$q\in\Rd$ with~$|q|_2=1$, and define~$\h$ and~$n$ by~\eqref{e.h.def} and~\eqref{e.parameter.choice}.  For~$K>m$, let~$W$ be the mean-zero solution of
	\begin{equation*}
		-\Delta W=-\nabla\cdot((e^\h-1)q)
		\quad\text{in }\cu_K\,,
		\qquad
		\partial_{\mathbf n}W=(e^\h-1)q\cdot\mathbf n
		\quad\text{on }\partial\cu_K\,,
	\end{equation*}
	and let~$\mathbf z=e^\h q-\nabla W$.  Then~$\nabla\cdot\mathbf z=0$ and~$\mathbf z\cdot\mathbf n=q\cdot\mathbf n$ on~$\partial\cu_K$. Even reflection of~$W$, with odd reflection of the normal forcing component and even reflection of the tangential components, preserves its weak Neumann equation and reduces it to the periodic divergence-form Poisson equation. Apply the periodic gradient estimate at exponent~$8$ and then H\"older's inequality for the lower exponents. Its constant is independent of~$K$ by scaling, and the same estimate applies on each parent cube below. This estimate and Lemma~\ref{l.laplacian.corrector.energy} give, for~$1\leq\xi\leq8$,
	\begin{equation*}
		\E\bigl[\|\nabla W\|_{\underline L^\xi(\cu_K)}^\xi\bigr]^{\nf{1}{\xi}}
		\leq C\gmcdelta h^{\nf{1}{2}}\,,
		\qquad
		\lim_{K\to\infty}
		\E\bigl[\|\nabla W\|_{\underline L^2(\cu_K)}^2\bigr]
		=
		\frac{2\tau^2h}{d}+O(\gmcdelta^4h^2)\,.
	\end{equation*}

\smallskip

    We retain a finite family of interior cells on which the weak Hessian satisfies a fourth-moment bound. To make the choice explicit, let~$N=1+\max\{2,68\lceil|\log_3\gmcdelta|\rceil\}$. For~$K$ sufficiently large, take parents~$z'+\cu_{n+N}$ with~$z'\in3^{n+N-1}\Zd$ contained in~$\cu_K$. Their concentric inner thirds tile the retained region. Retain all scale-$3^n$ cells in these inner thirds. The parents have overlap bounded by a dimensional constant, and the omitted cells lie within distance~$C3^{n+N}$ of~$\partial\cu_K$. Their relative volume tends to zero with~$m,h,\gmcdelta,N$ fixed.

\smallskip

Here is the two-radius decomposition giving the required fourth-moment bound. On each parent~$z'+\cu_{n+N}$, let~$U_{z'}\in H_0^1(z'+\cu_{n+N})$ solve~$-\Delta U_{z'}=-\nabla\cdot((e^\h-1)q)$, and let~$W_{z'}$ solve the corresponding mean-zero Neumann problem with normal derivative~$(e^\h-1)q\cdot\mathbf n$. In the interior of~$z'+\cu_{n+N}$, we decompose~$W=U_{z'}+(W_{z'}-U_{z'})+(W-W_{z'})$. The last two terms are harmonic, as follows by using compactly supported tests in the Neumann equations. The Dirichlet weak-Hessian estimate used in the upper proof gives a cell-averaged fourth-moment bound~$C\gmcdelta^{68}$ for the first term. For each harmonic term, the interior Hessian estimate on the inner third gives its scale-$3^n$ Hessian norm bounded by~$C3^{-N}$ times its normalized gradient norm on~$z'+\cu_{n+N}$. The local Dirichlet and Neumann gradient estimates bound their fourth moments uniformly, since~$(\E\fint_{z'+\cu_{n+N}}|e^\h-1|^4)^{\nf14}\leq C\gmcdelta h^{\nf12}\leq C$. For the restriction of~$\nabla W$ to the parents, Jensen and their bounded overlap reduce the averaged fourth moments to the global bound already proved. Consequently the averages over the finite family of the fourth moments of the two harmonic contributions are at most~$C3^{-4N}$. In particular they are at most~$C3^{-3(N-1)}\leq C\gmcdelta^{68}$. The inequality~$(a+b+c)^4\leq64(a^4+b^4+c^4)$ now gives the same fourth-moment bound for the weak Hessian of~$W$ on the retained cells. This argument uses the interior scalar equation and does not replace the global Neumann boundary condition by a scalar Dirichlet condition.

\smallskip

For every retained cell put
\begin{equation*}
	q_z=(\mathbf z)_{z+\cu_n}\,,\qquad
	G_z=\mathbf z-q_z\,,\qquad
	A_z=\|G_z\|_{\underline L^2(z+\cu_n)}\,,\qquad
	B_z=3^n\|\nabla G_z\|_{\underline L^2(z+\cu_n)}\,.
\end{equation*}
The additional derivative of~$e^\h q$ in~$\nabla G_z$ satisfies the fourth-moment bound~$C\gmcdelta^{68}$ by the field estimate in the upper proof. Mean-zero Poincar\'e gives~$A_z\leq CB_z$. Thus~$\limsup_{K\to\infty}\avsum_z\E[A_z^4+B_z^4]\leq C\gmcdelta^{68}$, and, in particular, both fourth-moment roots are at most~$C\gmcdelta^{15}$. All averages in this proof run over the retained cells. Measurability follows from continuous dependence of the weak solutions and their weak-derivative carriers, as in the upper proof. The cells are fixed by their geometry and the parameters, independently of the realized fields.

\smallskip

	\emph{Step 2.} We solve Neumann problems on the interior cells, patch the resulting fluxes into an admissible flux, and bound the dual energy by a sum over cells. On every interior cell, let~$\psi_z,\eta_z$ be the mean-zero solutions of
	\begin{equation*}
		-\nabla\cdot(\a_m\nabla\psi_z)=0\,,
		\qquad
		\a_m\nabla\psi_z\cdot\mathbf n=q_z\cdot\mathbf n\,,
	\end{equation*}
	and
	\begin{equation*}
		-\nabla\cdot(\a_m\nabla\eta_z)=0\,,
		\qquad
		\a_m\nabla\eta_z\cdot\mathbf n=G_z\cdot\mathbf n\,,
	\end{equation*}
		respectively. These problems are compatible because~$G_z$ is divergence free. The normal trace of the flux~$\a_m\nabla(\psi_z+\eta_z)$ agrees with that of~$\mathbf z$ on every cell boundary.  Therefore the field equal to~$\a_m\nabla(\psi_z+\eta_z)$ on each interior cell and to~$\mathbf z$ on the omitted cells belongs to~$q+\Lsolo(\cu_K)$, and the contribution of the omitted cells to its dual energy vanishes in expectation as~$K\to\infty$. Indeed, if~$O_K$ is their union, H\"older with exponents~$(4,4,2)$ on probability times normalized volume gives
\begin{equation*}
	\E\fint_{\cu_K}\indc_{O_K}|\mathbf z|^2\a_m^{-1}
	\leq C\left(\frac{|O_K|}{|\cu_K|}\right)^{\nf14}\longrightarrow0\,.
\end{equation*}
Here~$C$ may depend on~$m$ and~$h$ but not on~$K$; it is finite by the eighth-moment gradient estimate and the second moment of~$\a_m^{-1}$. Passing from the volume-weighted sum to the average over retained cells introduces another deterministic remainder tending to zero in expectation. Taking expectations in the dual variational inequality therefore gives
	\begin{multline*}
		\frac12\E[q\cdot\a_{m,*}^{-1}(\cu_K)q]
		\leq
		\avsum_z\E\Bigl[
		\frac12q_z\cdot\a_{m,*}^{-1}(z+\cu_n)q_z
		\\
		+\fint_{z+\cu_n}\a_m\nabla\psi_z\cdot\nabla\eta_z
		+\frac12\fint_{z+\cu_n}\a_m|\nabla\eta_z|^2
		\Bigr]+o_K(1)\,.
	\end{multline*}
	The argument of Step~2 of the proof of Lemma~\ref{l.one.step.upper}, applied to fluxes instead of gradients and with~\eqref{e.besov.grad.poincare} in place of~\eqref{e.besov.flux.poincare}, gives
	\begin{equation*}
		\fint_{z+\cu_n}\a_m|\nabla\eta_z|^2
		\leq
		C\lambda_{1/4,1}^{-1}(z+\cu_n;\a_m)
		\bigl(A_z^2+A_z^{\nf{3}{2}}B_z^{\nf{1}{2}}\bigr)\,.
	\end{equation*}
	The fixed-moment squared-error calculation and the simultaneous pathwise subcube comparison in Step~2 of the upper proof, now with inverse ratios, give the required coefficient bound. In particular, using~\eqref{e.annealed.ordering} to compare~$\ahom_n^{-1}$ with~$\ahom_{m-h}^{-1}$,
	\begin{equation*}
		\E\bigl[\lambda_{1/4,1}^{-2}(\cu_n;\a_m)\bigr]^{\nf{1}{2}}
		\leq C\ahom_{m-h}^{-1}\,.
	\end{equation*}
	The preceding three displays and Cauchy--Schwarz thus yield
	\begin{equation*}
		\limsup_{K\to\infty}\avsum_z
		\E\fint_{z+\cu_n}\a_m|\nabla\eta_z|^2
		\leq C\gmcdelta^{30}\ahom_{m-h}^{-1}\,,
	\end{equation*}
	\begin{equation*}
		\limsup_{K\to\infty}\avsum_z
		\E\left|\fint_{z+\cu_n}
		\a_m\nabla\psi_z\cdot\nabla\eta_z\right|
		\leq C\gmcdelta^{15}\ahom_{m-h}^{-1}\,.
	\end{equation*}
	With these two bounds, the dual variational inequality above bounds the dual energy of the competitor by the sum of the cell energies, up to an error~$C\gmcdelta^{15}\ahom_{m-h}^{-1}$ in the limit~$K\to\infty$.

\smallskip

	\emph{Step 3.} We estimate the sum of the cell energies~$q_z\cdot\a_{m,*}^{-1}(z+\cu_n)q_z$. The pair~$(q_z,G_z)$ is a function of the last~$h$ fields and is independent of~$\a_{m-h}$. Since~$\a_m^{-1}=e^{-\h}\a_{m-h}^{-1}$, the deterministic estimate in Step~3 of the proof of Lemma~\ref{l.one.step.upper}, with~$\h$ replaced by~$-\h$,~$p+\nabla w$ by~$\mathbf z$, and~$p_z$ by~$q_z$, gives
	\begin{equation*}
		\limsup_{K\to\infty}\avsum_z
		\E\bigl[q_z\cdot\a_{m,*}^{-1}(z+\cu_n)q_z\bigr]
		\leq
		\ahom_{m-h,*}^{-1}(\cu_n)
		\limsup_{K\to\infty}
		\E\fint_{\cu_K}|\mathbf z|^2e^{-\h}
		+C\gmcdelta^{15}\ahom_{m-h}^{-1}\,.
	\end{equation*}
	We next show that
	\begin{equation}
		\limsup_{K\to\infty}
		\left|
		\E\fint_{\cu_K}|\nabla W|^2(e^{-\h}-1)
		\right|
		\leq C\gmcdelta^4h^2\,.
		\label{e.lower.quadratic.error}
	\end{equation}
	Indeed, let~$V$ solve the Neumann problem defining~$W$ with~$e^\h-1$ replaced by its linearization~$\h$.  If~$V_0$ is the solution of the same problem with~$\h$ replaced by the sum~$\h+h\tau^2=\sum_{k=m-h+1}^m\g_k$ of the last~$h$ fields, then
	\begin{equation*}
		\nabla V=\nabla V_0-h\tau^2q\,,
		\qquad
		\E\|\nabla V_0\|_{\underline L^2(\cu_K)}^2
		\leq C\gmcdelta^2h\,.
	\end{equation*}
	Reversing the signs of the last~$h$ fields simultaneously gives
	\begin{equation*}
		\E\fint_{\cu_K}|\nabla V_0|^2(\h+h\tau^2)=0\,,
		\qquad
		\E\fint_{\cu_K}q\cdot\nabla V_0=0\,,
	\end{equation*}
	while Cauchy--Schwarz gives~$|\E\fint_{\cu_K}(q\cdot\nabla V_0)(\h+h\tau^2)|\leq C\gmcdelta^2h$. Expanding~$|\nabla V|^2\h$ with the identities for~$\nabla V$ and~$\h$ therefore yields~$\E\fint_{\cu_K}|\nabla V|^2\h=O(\gmcdelta^4h^2)$. The Neumann Calder\'on--Zygmund estimate gives~$\E[\|\nabla(W-V)\|_{\underline L^4}^4]^{\nf{1}{4}} \leq C\gmcdelta^2h$. The H\"older estimates used in Step~3 of the upper-bound proof, together with~$|e^{-s}-1+s|\leq\frac12s^2e^{|s|}$, control every term containing~$\nabla(W-V)$ and the nonlinear remainder by~$C\gmcdelta^4h^2$.  This proves~\eqref{e.lower.quadratic.error}.

\smallskip

	Since~$\mathbf z=e^\h q-\nabla W$,
	\begin{equation*}
		\fint_{\cu_K}|\mathbf z|^2e^{-\h}
		=
		\fint_{\cu_K}e^\h
		-2\fint_{\cu_K}q\cdot\nabla W
		+\fint_{\cu_K}e^{-\h}|\nabla W|^2\,.
	\end{equation*}
	The expectation of the first term is one.  Taking expectations in the mean-zero Neumann problem gives~$\E W=0$, so the expectation of the second term vanishes.  Lemma~\ref{l.laplacian.corrector.energy} and~\eqref{e.lower.quadratic.error} consequently give
	\begin{equation*}
		\limsup_{K\to\infty}
		\E\fint_{\cu_K}|\mathbf z|^2e^{-\h}
		\leq
		1+\frac{2\tau^2h}{d}+C\gmcdelta^4h^2\,.
	\end{equation*}

\smallskip

	\emph{Step 4.} We conclude the proof of~\eqref{e.one.step.lower}. Combining Steps~2 and~3, letting~$K\to\infty$, and recalling that~$\ahom_{m,*}(\cu_K)\to\ahom_m$, we obtain
	\begin{equation*}
		\frac12\ahom_m^{-1}
		\leq
		\left(\frac12+\frac{\tau^2h}{d}+C\gmcdelta^4h^2\right)
		\ahom_{m-h,*}^{-1}(\cu_n)
		+C\gmcdelta^{15}\ahom_{m-h}^{-1}\,.
	\end{equation*}
    The first-moment identity in Step~4 of the upper proof gives~$\ahom_{n,*}^{-1}(\cu_n)\leq(1+C\gmcdelta^2)\ahom_n^{-1}$. The finite-volume dual ordering~\eqref{e.annealed.matrix.bounds} and the monotone limit~\eqref{e.ahom.ahomstar.defs} yield
\begin{equation*}
	\ahom_{m-h}^{-1}\leq\ahom_{m-h,*}^{-1}(\cu_n)
	\leq e^{2\tau^2(m-h-n)}\ahom_{n,*}^{-1}(\cu_n)
	\leq(1+C\gmcdelta^2)e^{2\tau^2(m-h-n)}\ahom_{m-h}^{-1}\,.
\end{equation*}
Since~$\tau^2(m-h-n)\leq C\gmcdelta^2|\log\gmcdelta|$, we obtain
\begin{equation*}
	\ahom_{m-h}^{-1}\leq\ahom_{m-h,*}^{-1}(\cu_n)
	\leq(1+C\gmcdelta^2|\log\gmcdelta|)\ahom_{m-h}^{-1}\,.
\end{equation*}
The preceding bracket is positive and bounded by a dimensional constant. Absorbing~$\gmcdelta^{15}$ into~$\gmcdelta^2|\log\gmcdelta|$ and using~\eqref{e.tau.bound} and~$\gmcdelta h\leq1$ proves~\eqref{e.one.step.lower}.
\end{proof}

\begin{proof}[Proof of Proposition~\ref{p.sharp.asymptotic}]
	It suffices to consider~$\gmcdelta\leq c(d)$.  Let~$J\coloneqq16\lceil|\log_3\gmcdelta|\rceil$. The ordering~\eqref{e.annealed.ordering} gives
	\begin{equation}
		\left|\log\ahom_J+\frac{2\tau^2(J+1)}{d}\right|
		\leq C\gmcdelta^2(J+1)
		\leq C\gmcdelta^2|\log\gmcdelta|\,.
		\label{e.sharp.base.scale}
	\end{equation}
	By Lemmas~\ref{l.one.step.upper} and~\ref{l.one.step.lower} and the inequality~$\log t\leq t-1$ for~$t>0$, applied to~$\ahom_m/\ahom_{m-h}$ and to~$\ahom_{m-h}/\ahom_m$, whenever~$1\leq h\leq\gmcdelta^{-1}$ and~$m-h\geq J$,
	\begin{equation}
		\biggl|\log\ahom_m - \log\ahom_{m-h} + \frac{2\tau^2 h}{d}\biggr| \leq C\gmcdelta^4 h^2 + C\gmcdelta^2|\log\gmcdelta|\,.
		\label{e.log.increment}
	\end{equation}
	We may assume~$m>J$.  Starting at~$m$, use~\eqref{e.log.increment} repeatedly with~$h=\lfloor\gmcdelta^{-1}\rfloor$, each time lowering the scale by~$h$, for as long as the new scale stays at least~$J$, and then once more with the remaining~$h<\gmcdelta^{-1}$, if any, to arrive at~$J$.  This takes at most~$1+(m-J)/\lfloor\gmcdelta^{-1}\rfloor\leq C(1+m\gmcdelta)$ steps, and each step has~$h\leq\gmcdelta^{-1}$ and therefore contributes an error at most~$C\gmcdelta^2+C\gmcdelta^2|\log\gmcdelta|$, because~$\gmcdelta^4h^2\leq\gmcdelta^2$.  Together with~\eqref{e.sharp.base.scale}, this gives
	\begin{equation*}
		\left|\log\ahom_m+\frac{2\tau^2(m+1)}{d}\right|
		\leq
		C\gmcdelta^2|\log\gmcdelta|(1+m\gmcdelta)\,,
	\end{equation*}
	which is~\eqref{e.main.asymptotic}.
\end{proof}

\subsection{Further bounds on the effective diffusivities}

The error term in Proposition~\ref{p.sharp.asymptotic} dominates~$2\tau^2m/d$ unless~$\tau^2$ is bounded below by a multiple of~$\gmcdelta^3|\log\gmcdelta|$, and at~$m=0$ the proposition gives no upper bound. The following decay estimate assumes no such lower bound.

\begin{proposition}[Power-law decay of the effective diffusivities]
	\label{p.homogenized.coefficient.strict.decay}
	There exists~$\eta_{\mathrm{sub}}(d,\operatorname{law}(\g_0))>0$ satisfying
	\begin{equation}
		\eta_{\mathrm{sub}}=\frac{\tau^2}{\log3}\quad\text{if }d=2\,,
		\qquad
		0<\eta_{\mathrm{sub}}\leq\frac{\tau^2}{d\log3}
		\quad\text{if }d\geq3\,,
		\label{e.strict.decay.eta.range}
	\end{equation}
	and, for every~$m\in\N_0$,
	\begin{equation}
		\ahom_m
		\leq
		3^{-\eta_{\mathrm{sub}}(m+1)}\,.
		\label{e.homogenized.coefficient.strict.decay}
	\end{equation}
\end{proposition}

\begin{proof}
	In dimension~$d=2$, Proposition~\ref{p.special.two.d.exact.formula} proves the result with~$\eta_{\mathrm{sub}}=\tau^2/\log3$.  We may therefore suppose that~$d\geq3$. We then compare~$\ahom_m$ with the coarse-grained matrices of the product of the factors~$\exp(\g_k-\tau^2)$ over the multiples~$k$ of a fixed integer~$R$. Step~1 shows that, with the weight~$\exp(\g_0-\tau^2)$, the expected minimal energy of affine boundary data on a simplex is strictly smaller than that of the affine function. Step~2 replaces this minimal energy by a discretized one, which is still smaller than one and can be iterated. Step~3 iterates it over the scales~$0,R,2R,\ldots$, using the independence of the fields; this is the heart of the matter. Step~4 reinserts the omitted factors, which have mean one and are independent of the others, and concludes.

\smallskip

	Put
	\begin{equation*}
		A_N^{(R)}
		\coloneqq
		\prod_{j=0}^N \exp(\g_{jR}-\tau^2)
		\qquad(N\in\N_0)\,,
	\end{equation*}
	where the integer~$R\geq1$ will be fixed below.  Thus~$\E[\exp(\g_k(x)-\tau^2)]=1$ for every~$k$ and~$x$.  Define, with the maximum over all permutations~$\pi$ of~$\{1,\ldots,d\}$,
	\begin{equation}
		q_*
		\coloneqq
		\max_{\pi}\
		\sup_{|p|=1}
		\E\left[
		\min_{u\in p\cdot x+H_0^1(\spx_0^\pi)}
		\fint_{\spx_0^\pi}\exp(\g_0-\tau^2)|\nabla u|^2
		\right]\,.
		\label{e.strict.decay.qstar.def}
	\end{equation}
	\emph{Step 1.} We prove that~$q_*<1$. The affine competitor gives~$q_*\leq1$.  If equality held, there would be a permutation~$\pi$ and~$|p|=1$ for which the expected minimum in~\eqref{e.strict.decay.qstar.def} equals~$\E[\fint_{\spx_0^\pi}\exp(\g_0-\tau^2)|p|^2]=1$. The energy of~$p\cdot x$ minus the minimum is nonnegative and would then have zero expectation, so it would vanish almost surely.  Hence~$p\cdot x$ would be a minimizer, and the Euler--Lagrange equation would give~$p\cdot\nabla \exp(\g_0-\tau^2)=0$ in~$\spx_0^\pi$ almost surely.  By stationarity and the almost-sure~$C^1$ regularity of~$\exp(\g_0-\tau^2)$,~$p\cdot\nabla \exp(\g_0-\tau^2)=0$ on~$\Rd$.  Consequently,~$\exp(\g_0(0)-\tau^2)=\exp(\g_0(tp)-\tau^2)$ for every~$t$.  Choosing~$t$ with~$|tp|_2>\sqrt d$, we see that the random variables~$\exp(\g_0(0)-\tau^2)$ and~$\exp(\g_0(tp)-\tau^2)$, which are independent by Assumption~\ref{a.g1}, coincide almost surely; hence each is deterministic. It follows that~$\g_0(0)$ is deterministic and, since it has mean zero, that~$\tau^2=0$, contrary to Assumption~\ref{a.g4}.  Thus
	\begin{equation}
		q_*<1\,.
		\label{e.strict.decay.qstar.gap}
	\end{equation}

\smallskip

	\emph{Step 2.} We show that~$q_R$ converges to~$q_*$, and fix~$R$ with~$q_R<1$. For a permutation~$\pi$ and~$p\in\Rd$, let~$T$ run over all Kuhn simplices contained in~$\spx_0^\pi$ whose supporting cubes have side~$3^{-R}$, and let~$v$ run over the continuous functions on~$\overline{\spx_0^\pi}$ which are affine on every such~$T$ and equal to~$p\cdot x$ on~$\partial\spx_0^\pi$.  Set
	\begin{equation}
		Q_{R,\pi}(p)
		\coloneqq
		\min_{v}
		\frac1{|\spx_0^\pi|}
		\sum_{T}
		|T|\bigl(\sup_{\overline T}\exp(\g_0-\tau^2)\bigr)
		\bigl|\nabla v\vert_T\bigr|^2\,,
		\qquad
		q_R
		\coloneqq
		\max_{\pi}\
		\sup_{|p|=1}\E[Q_{R,\pi}(p)]\,.
		\label{e.strict.decay.qR.def}
	\end{equation}
	Since~$\exp(\g_0-\tau^2)$ is almost surely continuous and conforming piecewise-affine functions are dense in~$H_0^1(\spx_0^\pi)$, we have, almost surely as~$R\to\infty$,
	\begin{equation*}
		Q_{R,\pi}(p)
		\longrightarrow
		\min_{u\in p\cdot x+H_0^1(\spx_0^\pi)}
		\fint_{\spx_0^\pi}\exp(\g_0-\tau^2)|\nabla u|^2\,.
	\end{equation*}
	The affine competitor bounds~$Q_{R,\pi}(p)$ by~$|p|^2\sup_{\spx_0^\pi}\exp(\g_0-\tau^2)$, whose expectation is finite by Assumption~\ref{a.g2}.  Dominated convergence and~\eqref{e.strict.decay.qstar.gap} show that~$q_R\to q_*<1$. Fix~$R$ such that
	\begin{equation}
		q_R<1\,.
		\label{e.strict.decay.qR.gap}
	\end{equation}

\smallskip

	\emph{Step 3.} We iterate~\eqref{e.strict.decay.qR.gap} over scales separated by~$R$. We show by induction on~$N$ that, for every~$N\in\N_0$, every permutation~$\pi$ and every~$p\in\Rd$,
	\begin{equation}
		\E\left[
		p\cdot A_N^{(R)}(\spx_{NR}^{\pi})p
		\right]
		\leq q_R^{N+1}|p|^2\,.
		\label{e.strict.decay.sparse.induction}
	\end{equation}
	For~$N=0$, this follows from~\eqref{e.strict.decay.qR.def}, because restricting the minimization to piecewise-affine functions and replacing~$\exp(\g_0-\tau^2)$ by its supremum on each simplex can only increase the energy. Suppose that~\eqref{e.strict.decay.sparse.induction} holds with~$N-1$ in place of~$N$. Fix a parent simplex~$\spx_{NR}^\pi$ and a vector~$p$.  The Kuhn simplices~$T$ contained in~$\spx_{NR}^\pi$ and supported on cubes of side~$3^{(N-1)R}$ form, after rescaling, the mesh used in~\eqref{e.strict.decay.qR.def}.  On this mesh, we consider the continuous function with boundary values~$p\cdot x$ that is affine on each simplex and minimizes the energy in~\eqref{e.strict.decay.qR.def} with the weights~$\sup_{\overline T}\exp(\g_{NR}-\tau^2)$ in place of~$\sup_{\overline T}\exp(\g_0-\tau^2)$, and let~$p_T$ be its slope on~$T$.  On every~$T$, we then take the minimizer of~$\fint_TA_{N-1}^{(R)}|\nabla u|^2$ among the functions whose boundary values on~$\partial T$ agree with this piecewise-affine function.  These minimizers glue to a competitor in~$p\cdot x+H_0^1(\spx_{NR}^\pi)$, and the variational formula gives
	\begin{equation*}
		p\cdot A_N^{(R)}(\spx_{NR}^\pi)p
		\leq
		\sum_T\frac{|T|}{|\spx_{NR}^\pi|}
		\bigl(\sup_{\overline T}\exp(\g_{NR}-\tau^2)\bigr)
		p_T\cdot A_{N-1}^{(R)}(T)p_T\,.
	\end{equation*}
	Condition on~$\g_{NR}$.  The slopes~$p_T$ and the displayed weights are then fixed, whereas~$A_{N-1}^{(R)}$ is independent of~$\g_{NR}$.  Stationarity and the induction hypothesis bound the conditional expectation of the right side by
	\begin{equation*}
		q_R^N
		\sum_T\frac{|T|}{|\spx_{NR}^\pi|}
		\bigl(\sup_{\overline T}\exp(\g_{NR}-\tau^2)\bigr)|p_T|^2\,.
	\end{equation*}
	Taking the remaining expectation and using the scaling relation~$\g_{NR}\overset{\mathrm{law}}=\g_0(3^{-NR}\,\cdot)$ and~\eqref{e.strict.decay.qR.def}, we obtain
	\begin{equation*}
		\E\left[
		p\cdot A_N^{(R)}(\spx_{NR}^\pi)p
		\right]
		\leq q_R^{N+1}|p|^2\,,
	\end{equation*}
	which proves~\eqref{e.strict.decay.sparse.induction}.

\smallskip

	\emph{Step 4.} We pass from the products~$A_N^{(R)}$, which contain only every~$R$-th field, to~$\ahom_m$ and conclude. Fix~$m\in\N_0$ and let~$N\coloneqq\lfloor m/R\rfloor$.  The field~$\a_m/A_N^{(R)}$ is the product of the factors~$\exp(\g_k-\tau^2)$ with~$0\leq k\leq m$ and~$R\nmid k$.  It is independent of~$A_N^{(R)}$ and satisfies~$\E[(\a_m/A_N^{(R)})(x)]=1$.  If~$u$ minimizes~$\fint_VA_N^{(R)}|\nabla u|^2$ among the functions on a bounded convex domain~$V$ with boundary values~$p\cdot x$, then~$u$ is a function of~$A_N^{(R)}$, so conditioning on~$A_N^{(R)}$ and using its independence from~$\a_m/A_N^{(R)}$ give
	\begin{equation}
		\E\left[p\cdot\a_m(V)p\right]
		\leq
		\E\left[\fint_V \a_m|\nabla u|^2\right]
		=
		\E\left[p\cdot A_N^{(R)}(V)p\right]\,.
		\label{e.strict.decay.full.sparse.comparison}
	\end{equation}
	For~$K\geq NR$, partition~$\cu_K$ into the Kuhn simplices of the cubes of side~$3^{NR}$.  Using on each simplex the minimizer for~$A_N^{(R)}$ with boundary values~$p\cdot x$ as a competitor, we deduce from~\eqref{e.strict.decay.sparse.induction} that
	\begin{equation*}
		\E\left[
		p\cdot A_N^{(R)}(\cu_K)p
		\right]
		\leq q_R^{N+1}|p|^2\,.
	\end{equation*}
	Letting~$K\to\infty$ in~\eqref{e.strict.decay.full.sparse.comparison} and using~\eqref{e.ahom.ahomstar.defs}, we obtain~$\ahom_m\leq q_R^{\lfloor \nf{m}{R}\rfloor+1}$. Since~$R(\lfloor m/R\rfloor+1)\geq m+1$, the choice~$\eta_{\mathrm{sub}}\coloneqq-\log q_R/(R\log3)>0$ proves~\eqref{e.homogenized.coefficient.strict.decay}.  Replacing this value by its minimum with~$\tau^2/(d\log3)$ proves~\eqref{e.strict.decay.eta.range} as well.
\end{proof}

Proposition~\ref{p.sharp.asymptotic} controls~$\ahom_m$ only up to its error term. The following lower bound has no error term and holds for every~$\gmcdelta$.

\begin{proposition}[Lower bound on the effective diffusivities]
	\label{p.homogenized.coefficient.reciprocal.lower}
	For every~$m\in\N_0$,
	\begin{equation}
		\ahom_m
		\geq
		e^{-(m+1)\tau^2}\,.
		\label{e.homogenized.coefficient.reciprocal.lower}
	\end{equation}
\end{proposition}

\begin{proof}
	Fix~$m,n\in\N_0$. Since~$e^{(m+1)\tau^2}\a_m=\exp(\sum_{k=0}^m\g_k)$, Assumption~\ref{a.g3} gives~$e^{-(m+1)\tau^2}\a_m^{-1}\overset{\mathrm{law}}=e^{(m+1)\tau^2}\a_m$.  Since coarse-grained matrices are homogeneous of degree one in the coefficient field, it follows that
	\begin{equation}
		\E[\a(\cu_n;e^{(m+1)\tau^2}\a_m)]
		=
		\E[\a(\cu_n;e^{-(m+1)\tau^2}\a_m^{-1})]
		=
		e^{(m+1)\tau^2}\ahom_m(\cu_n)\,.
		\label{e.reciprocal.finite.volume.isotropy}
	\end{equation}
	Here the second matrix is the coarse-grained matrix of the reciprocal coefficient field, and not the inverse of the first.

\smallskip

	For~$i\in\{1,\ldots,d\}$, let~$u_i$ and~$v_i$ be the minimizers in~\eqref{e.variational.a}, with~$p=e_i$, for the respective coefficients~$e^{(m+1)\tau^2}\a_m$ and~$e^{-(m+1)\tau^2}\a_m^{-1}$ on~$\cu_n$.  Let~$E$ and~$F$ be the Jacobian matrices of~$x\mapsto(u_i(x))_i$ and~$x\mapsto(v_i(x))_i$.  Since~$u_i$ and~$v_i$ have boundary values~$e_i\cdot x$, both maps have trace~$x$ on~$\partial\cu_n$.  For two~$d$-by-$d$ matrices, set
	\begin{equation*}
		Q(E,F)
		\coloneqq
		\operatorname{tr}(E)\operatorname{tr}(F)
		-
		\operatorname{tr}(EF)\,,
		\qquad
		|M|_{\mathrm F}\coloneqq\operatorname{tr}(M^{\mathrm T}M)^{\nf12}\,,
	\end{equation*}
	where~$|\cdot|_{\mathrm F}$ is the Frobenius norm.  The bilinear form~$Q$ is the polarization of the sum of the second-order principal minors.  These minors are null Lagrangians, and~$(u_i)_i$ and~$(v_i)_i$ have the same trace, so
	\begin{equation}
		\fint_{\cu_n}Q(E,F)
		=
		Q(\Id,\Id)
		=
		d(d-1)\,.
		\label{e.reciprocal.mixed.minors}
	\end{equation}
	Equivalently, this follows by integration by parts from
	\begin{equation*}
		Q(E,F)
		=
		\sum_{i,j=1}^d
		\left(
		\partial_i(u_i\partial_jv_j)
		-
		\partial_j(u_i\partial_iv_j)
		\right)\,.
	\end{equation*}

\smallskip

	For arbitrary matrices~$E,F$, Cauchy--Schwarz and~$\operatorname{tr}(F)^2\leq d|F|_{\mathrm F}^2$ give
	\begin{equation*}
		Q(E,F)
		=
		\operatorname{tr}\bigl(E^{\mathrm T}(\operatorname{tr}(F)\Id-F^{\mathrm T})\bigr)
		\leq
		(d-1)|E|_{\mathrm F}|F|_{\mathrm F}\,,
	\end{equation*}
	and
	\begin{equation*}
		\left|\operatorname{tr}(F)\Id-F^{\mathrm T}\right|_{\mathrm F}^2
		=
		|F|_{\mathrm F}^2+(d-2)\operatorname{tr}(F)^2
		\leq
		(d-1)^2|F|_{\mathrm F}^2\,.
	\end{equation*}
	Consequently, pointwise on~$\cu_n$,
	\begin{equation*}
		e^{(m+1)\tau^2}\a_m|E|_{\mathrm F}^2+e^{-(m+1)\tau^2}\a_m^{-1}|F|_{\mathrm F}^2
		\geq
		2|E|_{\mathrm F}|F|_{\mathrm F}
		\geq
		\frac{2}{d-1}Q(E,F)\,.
	\end{equation*}
	Integrating, taking expectations, and using~\eqref{e.reciprocal.finite.volume.isotropy},~\eqref{e.variational.a}, and~\eqref{e.reciprocal.mixed.minors}, we obtain
	\begin{equation*}
		2e^{(m+1)\tau^2}\operatorname{tr}\ahom_m(\cu_n)
		=
		\E\left[
		\fint_{\cu_n}
		\left(e^{(m+1)\tau^2}\a_m|E|_{\mathrm F}^2+e^{-(m+1)\tau^2}\a_m^{-1}|F|_{\mathrm F}^2\right)
		\right]
		\geq
		2d\,.
	\end{equation*}
	The matrix~$\ahom_m(\cu_n)$ is scalar by Assumption~\ref{a.g3}. Dividing by~$2d$ and letting~$n\to\infty$ proves~\eqref{e.homogenized.coefficient.reciprocal.lower}.
\end{proof}

\section{Large-scale regularity and Dirichlet homogenization}
\label{s.Holder.regularity}

In this section we prove the large-scale H\"older estimate (Theorem~\ref{t.C}) and the homogenization error bound for the Dirichlet problem (Theorem~\ref{t.B}). The regularity estimates are used in Sections~\ref{sec:mfd-uniqueness} and~\ref{sec:mfd-process-convergence}.

\smallskip

Heuristically, a solution of an equation with a rapidly oscillating coefficient cannot be regular on small scales, but it can be regular on large ones: if on every scale it is close to a solution of the homogenized equation, then it inherits the regularity of that solution down to the scale at which the approximation fails.  The quantity we iterate is the excess, the~$L^2$ distance from a solution to the affine functions on a cube.  On a scale where the coarse-graining error is small, the excess contracts by a fixed factor; if such scales are frequent enough, the excess therefore decays geometrically and a H\"older estimate follows.

\smallskip

The scales on which the error fails to be small occur in a proportion that we can make as small as we wish by lowering the H\"older exponent; here randomness and independence between the scales enter through a concentration inequality. The cutoff coefficient~$\a_L$ contains no fields of scale larger than the cutoff scale~$3^L$, so averaging over cubes larger than the cutoff scale improves its comparison with a constant coefficient. For this reason no relation between the size of the domain and the cutoff scale is needed, and the estimates are uniform in~$L$ and pass to the limiting coefficient.

\smallskip

The section is organized as follows.
\begin{itemize}
\item In Section~\ref{ss.density.of.good.scales} we define the good event, on which the coarse-graining error at a given scale is deterministically bounded (Proposition~\ref{p.good.scale.mathcal.E}), and we show that, outside an event of exponentially small probability, the proportion of scales on which the good event fails is small (Proposition~\ref{p.density.of.good.scales}).
\item In Section~\ref{ss.good.scale.estimates} we show that, on a good event, a solution is close to a harmonic function and its excess contracts by a fixed factor (Lemmas~\ref{l.harmonic.approximation.good.scales.GMC} and~\ref{l.excess.decay.good.scales.GMC}).  The iteration of Lemma~\ref{l.iteration.lemma.GMC} combines this contraction on the good scales with a loss of at most a bounded factor on each of the other scales.
\item In Section~\ref{ss.regularity.iteration} we define the minimal scale and prove Proposition~\ref{p.Holder.regularity}, the large-scale H\"older estimate in~$\cu_m$ for~$L\geq m$.
\item In Section~\ref{ss.cutoff.uniform.regularity} we remove the restriction~$L\geq m$ (Proposition~\ref{p.cutoff.Holder.regularity}) and show that the estimates are stable when the coefficient is multiplied by a positive continuous function that is uniformly close to a constant (Proposition~\ref{p.cutoff.Holder.bounded.multiplier}).
\item Theorems~\ref{t.cutoff.Dirichlet.homogenization} and~\ref{t.fixed.cutoff.Dirichlet.algebraic} are the two parts of Theorem~\ref{t.B}: they bound the homogenization error of the Dirichlet problem uniformly in the cutoff and, at a fixed cutoff, with an algebraic rate.  They follow from the general coarse-graining estimate of Proposition~\ref{p.general.coarse.graining}, applied in the whole domain, and do not use the regularity estimates.
\item In the last subsection we extend the regularity estimates to the limiting coefficient~$\a$, which proves Theorem~\ref{t.C} in the form of Theorem~\ref{t.large.scale.Holder.multifractal}.
\end{itemize}

In~\eqref{e.excess.def} and~\eqref{e.excess.affine.def} below,~$\ell$ denotes an affine function on~$\Rd$, and the minimum and the minimizer are taken over all affine functions.  For~$m,n \in \Z$ with~$n \leq m$ and~$x \in \cu_m$, we define the~\textbf{excess} of a function~$u\in L^2((x +\cu_n) \cap \cu_m)$ by
\begin{equation}
	E(u,(x +\cu_n) \cap \cu_m)  \coloneqq
	3^{-n}
	\min_{\ell}
	\| u - \ell \|_{\underline{L}^2((x +\cu_n) \cap \cu_m)}
	\label{e.excess.def}
\end{equation}
and its best affine approximation by
\begin{equation}
	\ell (u,(x +\cu_n) \cap \cu_m)   \coloneqq  \argmin_{\ell}
	\| u - \ell \|_{\underline{L}^2((x +\cu_n) \cap \cu_m)}
	\,.
	\label{e.excess.affine.def}
\end{equation}
Recall that~$\b_{L,m}$ is defined in Definition~\ref{d.bLm}.

\begin{proposition}[Large-scale H\"older estimate for~$L\geq m$]
\label{p.Holder.regularity}
There exists a constant~$C(d)<\infty$ such that, if~$\gmcdelta \leq C^{-1}$, then, for every~$\gamma \in [\nf12, 1 - C\gmcdelta \left| \log \gmcdelta \right|^{\nf 12} ]$ and~$m\in\N_0$, there exists an~$\N$-valued random variable~$\mathcal{L}_L(\gamma,m)$ (the same for every~$L\geq m$) satisfying
\begin{equation}
	\P \bigl[ \mathcal{L}_L(\gamma,m) > k \bigr]
	\leq C\exp \biggl(
	-\frac{(1-\gamma)^2(k-C)_+}
	{C\gmcdelta^2 \left| \log \gmcdelta \right|}
	\biggr)
	\qquad(k\in\N)
	\label{e.mathcal.L.bound}
\end{equation}
such that, for every~$L\in\N_0$ with~$L\geq m$ and every~$\g \in C^{0,\nf12}(\cu_m;\Rd)$,~$h \in C^{1,\nf12}(\cu_m)$ and~$u\in H^1(\cu_m)$ solving
\begin{equation}
	\left\{
	\begin{aligned}
		& -\nabla \cdot \a_L \nabla u = \nabla \cdot \g  && \text{in} \ \cu_m\,, \\
		& u = h && \text{on} \ \partial \cu_m\,,
	\end{aligned}
	\right.
	\label{e.eq.for.u.holder.estimate}
\end{equation}
we have the following estimates:
\begin{itemize}

\item \underline{\emph{Large-scale H\"older estimate}}: for every~$n\in\N_0$ with~$n\leq m-\mathcal{L}_L(\gamma,m)$ and~$x\in\cu_m$,
\begin{align}
	\lefteqn{
	\sup_{\substack{j\in\N_0\\ j\leq n}}
	\max_{y\in3^j\Zd\cap((x+\cu_n)\cap\cu_m)}
	3^{\gamma(n-j)}
	\left\|u-(u)_{(y+\cu_j)\cap\cu_m}\right\|_{\underline{L}^2((y+\cu_j)\cap\cu_m)}
	}\qquad &
	\notag \\ &
	\leq
	C 3^{-\gamma(m-n)}
	\Bigl(
	\| u - (u)_{\cu_{m}} \|_{\underline{L}^2(\cu_{m})}
	+
	(\b_{L,m})_{\cu_m}^{-1} 3^{\frac{3}{2} m} [\g]_{\underline{W}^{\nf12,\infty}(\cu_m)}
	\Bigr)
	\notag \\ & \qquad
	+
	C 3^{-\gamma(m-n)} 3^{\frac{3}{2} m} \|\nabla h\|_{\underline{W}^{\nf12,\infty}(\cu_m)}
	\indc_{\{x\notin\cu_{m-1}\}}\,;
	\label{e.Holder.estimate.boxes.local}
\end{align}

\item \underline{\emph{Energy estimate}}: for every~$n\in\N_0$ with~$n\leq m-\mathcal{L}_L(\gamma,m)$ and~$x\in\cu_m$,
\begin{align}
	\| \a_L^{\nicefrac12} \nabla u \|_{\underline{L}^2((x+\cu_{n}) \cap \cu_m)}
	&\leq
	C 3^{(1-\gamma)(m-n)}
	\Bigl(
	\| \a_L^{\nicefrac12} \nabla u \|_{\underline{L}^2(\cu_{m})}
	+
	(\b_{L,m})_{\cu_m}^{-\nf12} 3^{\frac{1}{2} m} [\g]_{\underline{W}^{\nf12,\infty}(\cu_m)}
	\Bigr)
	\notag \\ & \qquad
	+
	C 3^{(1-\gamma)(m-n)} (\b_{L,m})_{\cu_m}^{\nf12} 3^{\frac{1}{2} m}
	\|\nabla h\|_{\underline{W}^{\nf12,\infty}(\cu_m)} \indc_{\{x\notin\cu_{m-1}\}}\,;
	\label{e.energy.density.estimate}
\end{align}

\item \underline{\emph{Excess decay estimate}}: if~$n\in\N_0$ is such that~$\mathcal{L}_L(\gamma,m) \leq m-n \leq (1-\gamma)^{-1}$, then for every~$l\in\N_0$ satisfying~$n\leq l\leq m-5$ and every~$x\in\cu_m$,
\begin{align}
	\lefteqn{E(u,(x+\cu_{n})\cap\cu_m)} \qquad &
	\notag \\ &
	\leq
	C3^{-\frac{1}{4}(l-n)} E(u,(x+\cu_{l})\cap\cu_m)
	+
	C(1-\gamma)^{\nf12} 3^{-l}
	\| u - (u)_{(x+\cu_{l})\cap\cu_m} \|_{\underline{L}^2((x+\cu_{l})\cap\cu_m)}
	\notag \\ & \qquad
	+
	C(\b_{L,m})_{\cu_m}^{-1} 3^{\frac{1}{2} l} [\g]_{\underline{W}^{\nf12,\infty}(\cu_m)}
	\notag \\ & \qquad
	+
	C\Bigl( (1-\gamma)(m-n)\|\nabla h\|_{L^\infty(\cu_m)}
	+ 3^{\frac{1}{2} l} \|\nabla h\|_{\underline{W}^{\nf12,\infty}(\cu_m)} \Bigr)
	\indc_{\{x\notin\cu_{m-1}\}}\,.
	\label{e.excess.decay.estimate}
\end{align}

\end{itemize}
\end{proposition}

At scale~$3^m$, Proposition~\ref{p.general.coarse.graining} compares an~$\a_L$-harmonic function with a harmonic function when the coarse-graining error~$\mathcal{E}_{\nicefrac14, \infty} (\cu_m\,;\a_L,(\b_{L,m})_{\cu_m})$ is small.

\subsection{Harmonic approximation and independence between scales}
\label{ss.density.of.good.scales}

We now introduce the good event~$\mathcal G_m(\ep,s)$ at scale~$m$, on which the coarse-graining error at that scale is deterministically bounded, and estimate the number of scales at which it fails.  Proposition~\ref{p.good.scale.mathcal.E} bounds the coarse-graining error on a scale where the good event holds, and Proposition~\ref{p.density.of.good.scales} shows that, among consecutive scales, the proportion of those where it fails is small.  These two estimates are the input to the iteration of Lemma~\ref{l.iteration.lemma.GMC} in Section~\ref{ss.regularity.iteration}, and the second uses the independence of the fields across scales.

\begin{definition}[Good event]
\label{d.good.events}
For~$s,\ep\in(0,1]$ and~$m\in\N_0$, let the \textbf{good event}~$\mathcal G_m(\ep,s)$ be the event on which the following inequalities hold. The gradient magnitude in~\eqref{e.good.gj} is~$|\nabla\g_i|_2$, and the unit vectors in~\eqref{e.good.sum.J} are Euclidean:
\begin{align}
	&
	\sup_{j \geq 0}  3^{-\frac{1}{8} s j}
	\sum_{i=(m-j)\vee0}^{m+j} \bigl\|  | \g_i |
	+
	3^{i} | \nabla \g_{i}| \bigr\|_{L^\infty(\cu_{m+1+j})}
	\leq
	\ep
	\,,
	\label{e.good.gj}
	\\ &
	\sup_{j \geq 0}  3^{-\frac{1}{8} s j}
	\biggl\| \, \prod_{i=(m-j)\vee0}^{m+j} \exp\bigl(| \g_i |\bigr)
	+
	\prod_{i=m+j}^{\infty} \exp\bigl(4 | \g_i - \g_i(0) | \bigr)
	\, \biggr\|_{L^\infty(\cu_{m+1+j})}
	\leq
	6 \,,
	\label{e.good.discounted.am.vs.aj}
	\\ &
	\sup_{0\leq j \leq m} \sup_{0\leq n \leq j-2} 3^{-\frac{1}{8} s(m-n)}
	\max_{z \in 3^n \Zd \cap (\cu_j \setminus \cu_{j-1})}
	\max_{|e| = 1} J(z+\cu_n,e,e\,; \ahom_{n}^{-1} \a_{n})
	\leq
	\ep^2
	\,.
	\label{e.good.sum.J}
\end{align}
For~$y \in \Rd$, we let~$\mathcal{G}_{m,y} ( \ep,s)$ denote the event obtained by translating the displayed quantities by~$y$. The supremum in~\eqref{e.good.sum.J} is defined to be zero when there is no pair~$(j,n)$ satisfying the displayed constraints.
\end{definition}

We call a scale~$m$ \textbf{good} if the event~$\mathcal G_m(\ep,s)$ holds: the fields at scales~$m\pm j$ and the quantities~$J$ on subcubes of scale~$m-j$ are small, up to weights~$3^{-sj/8}$ that decay with~$j$. A scale that is not good is called \textbf{bad}. We use the same words for the translated events~$\mathcal G_{m,y}(\ep,s)$, with respect to the center~$y$.

\smallskip

Throughout this section, we use the following notation. For~$m,j \in \N_0$ with~$j\leq m$, we denote
\begin{equation}
	\g_{m,j} \coloneqq \sum_{i=j+1}^m \g_i \,.
	\label{e.g.m.j}
\end{equation}
Moreover, for~$0\leq j\leq m$,
\begin{equation}
	\frac{\a_m}{\a_j}
	=
	\exp\bigl(\g_{m,j}-(m-j)\tau^2\bigr)\,.
	\label{e.a.m.j.ratio}
\end{equation}
For~$s,\ep\in(0,1]$, condition~\eqref{e.good.gj} implies
\begin{equation}
	\biggl( \sup_{0\leq j \leq m}
	3^{-\frac{1}{8} s(m-j)} \| \g_{m,j} \|_{L^\infty(\cu_{m})}
	+
	\sum_{j=m}^\infty
	3^{m}
	\| \nabla \g_j \|_{L^\infty(\cu_m)}
	\biggr) \indc_{\mathcal{G}_m ( \ep,s)}
	\leq
	C \ep
	\,.
	\label{e.what.e.good.gj.gives}
\end{equation}

\begin{proposition}[Coarse-graining on good scales]
\label{p.good.scale.mathcal.E}
There is a constant~$C(d)<\infty$ such that, for every~$s\in[64\gmcdelta^2,\nf12]$ and~$L,m\in\N_0$ with~$L\geq m$, we have
\begin{align}
	\lefteqn{
	\mathcal E_{s,\infty,2}
	\bigl(\cu_m;\a_L,(\b_{L,m})_{\cu_m}\bigr)
	\indc_{\mathcal G_m(1,s)}
	} \qquad &
	\notag \\ &
	\leq
	C \sup_{0\leq j \leq m} \sup_{0\leq n \leq j-2}
	3^{-\frac{1}{2} s(m-n)}
	\max_{\substack{z \in 3^n \Zd \cap (\cu_j \setminus \cu_{j-1})\\ |e|=1}}
	J\bigl(z+\cu_n,e,e\,; \ahom_n^{-1}\a_n\bigr)^{\nf12}
	\notag \\ & \qquad
	+
	C s^{-1} \gmcdelta^2
	\notag \\ & \qquad
	+
	C\sup_{0\leq j \leq m}
	3^{-\frac{1}{8} s(m-j)} \| \g_{m,j} \|_{L^\infty(\cu_{m})}
	\notag \\ & \qquad
	+
	C3^{-\frac{1}{8}sm}
	\left\|\sum_{i=0}^m\g_i\right\|_{L^\infty(\cu_m)}
	\notag \\ & \qquad
	+
	C\sum_{j=m}^\infty
	3^{m}
	\| \nabla \g_j \|_{L^\infty(\cu_m)}
	\,.
	\label{e.mathcalE.bound.for.good}
\end{align}
In particular, if~$\ep\in[s^{-1}\gmcdelta^2,1]$, then~\eqref{e.what.e.good.gj.gives} yields, on~$\mathcal G_m(\ep,s)$,
\begin{equation*}
	\mathcal E_{s,\infty,2}
	\bigl(\cu_m;\a_L,(\b_{L,m})_{\cu_m}\bigr)
	\leq C\ep\,.
\end{equation*}
\end{proposition}

We will also use the following elementary inequality. If~$\ep\in(0,1]$,~$s>0$, and~$D\geq0$, then every~$x\in[0,\ep3^{\frac{s}{16}D}]$ satisfies~$3^{-\frac{s}{2}D}x\leq3^{-\frac{s}{2}D}\min\{x,1\}+\ep^8$. Indeed, the second term on the right is needed only when~$x>1$, and then~$3^{-\frac{s}{16}D}<\ep$, so that~$3^{-\frac{s}{2}D}x\leq\ep3^{-\frac{7s}{16}D}\leq\ep^8$.

\smallskip

\begin{proposition}[Proportion of good scales]
\label{p.density.of.good.scales}
There exists~$C(d)<\infty$ such that, for every~$s,\theta,\ep \in (0,1]$ satisfying
\begin{equation}
	\theta \geq C s^{-6} \ep^{-2} \gmcdelta^{2} \left| \log \gmcdelta\right|
	\,,
	\label{e.theta.vs.gmcdelta.cond.reg}
\end{equation}
and every~$m_0\in \N_0$ and~$M \in \N_0$, we have the estimate
\begin{equation}
	\P \Biggl[
	\avsum_{m=m_0}^{m_0+M}
	\indc{\bigl\{ \mathcal{G}_m (  \ep,s) \bigr\}}
	\leq
	1- \theta \Biggr]
	\leq
	\exp \biggl(  -  \frac{s^6 \ep^2 \theta}{C \gmcdelta^2 \left| \log \gmcdelta\right| } (M+1) \biggr)
	\,.
	\label{e.good.H.one}
\end{equation}
\end{proposition}

\begin{proof}[{Proof of Proposition~\ref{p.good.scale.mathcal.E}}]
Fix~$s \in [64 \gmcdelta^2,\nf 12]$. Throughout the proof, we denote~$\tilde{\a}_{L,m}  \coloneqq (\b_{L,m})_{\cu_m}^{-1}  \a_{L}$. The proof compares each quantity~$J(z+\cu_n,e,e\,;\tilde\a_{L,m})$ entering~$\mathcal E_{s,\infty,2}$ with the same quantity for the coefficient~$\ahom_n^{-1}\a_n$. Step~1 reduces~$\mathcal E_{s,\infty,2}$ to the cubes lying in the annuli~$\cu_j\setminus\cu_{j-1}$. Step~2 is the heart of the matter: it treats the cubes of side at most one directly and, on the larger cubes, bounds the ratio of the two coefficients by means of conditions~\eqref{e.good.gj} and~\eqref{e.good.discounted.am.vs.aj}. Step~3 combines the bounds.

\smallskip

\emph{Step 1.} We bound the error by the contributions of cubes lying in the annuli~$\cu_j\setminus\cu_{j-1}$. There exists a constant~$C(d)<\infty$ such that, for every~$m \in \N_0$ and~$s \in (0,\nf 12]$,
\begin{align}
	\mathcal{E}_{s,\infty,2}(\cu_m ; \a_L,(\b_{L,m})_{\cu_m})^2
	&
	=
	(1-3^{-2s}) \sum_{n=-\infty}^m 3^{-2s(m-n)} \max_{z \in 3^n \Zd \cap \cu_m} \max_{|e| = 1}  J(z+\cu_n,e,e\,;\tilde{\a}_{L,m})
	\notag \\ &
	\leq
	C \sup_{j \leq m} \sup_{n \leq j-2} 3^{-\frac{3}{2} s(m-n)}
	\max_{z \in 3^n \Zd \cap (\cu_j \setminus \cu_{j-1})}
	\max_{|e| = 1} J(z+\cu_n,e,e\,;\tilde{\a}_{L,m})
	\,.
	\label{e.mathcalE.annular.decomp.pre}
\end{align}
For a domain~$U$ and~$e \in \Rd$, write~$J(U) \coloneqq J(U,e,e\,;\tilde{\a}_{L,m})$. By decomposing into annular regions and interchanging the order of summation, we obtain
\begin{align*}
	\lefteqn{
	s \sum_{n=-\infty}^m 3^{-2s(m-n)} \max_{z \in 3^n \Zd \cap \cu_m} J(z+\cu_n)
	} \qquad &
	\\ &
	\leq
	s \sum_{n=-\infty}^m 3^{-2s(m-n)} \max_{n+1 \leq j \leq m} \max_{z \in 3^n \Zd \cap (\cu_j \setminus \cu_{j-1})} J(z+\cu_n)
	+ s\sum_{n=-\infty}^m 3^{-2s(m-n)} J(\cu_n)
	\\ &
	\leq C
	\sup_{j \leq m} \sup_{n \leq j-1} 3^{-\frac{3}{2} s(m-n)}   \max_{z \in 3^n \Zd \cap (\cu_j \setminus \cu_{j-1})} J(z+\cu_n)
	+ s \sum_{n=-\infty}^m 3^{-2s(m-n)} J(\cu_n)
	\,.
\end{align*}
For the second term on the right, we use subadditivity. For every~$n \leq m$, we have
\begin{equation*}
	J(\cu_n)
	\leq
	\sum_{k=-\infty}^{n-1} \sum_{z \in 3^{k} \Zd \cap (\cu_{k+1} \setminus \cu_{k})}
	\frac{|\cu_k|}{|\cu_n|}J(z+\cu_k)
	\leq
	\sum_{k=-\infty}^{n-1} 3^{-d(n-1-k)}
	\max_{z \in 3^{k} \Zd \cap (\cu_{k+1} \setminus \cu_{k})} J(z+\cu_k)
	\,.
\end{equation*}
Substituting back and interchanging the order of summation yields
\begin{align*}
	\sum_{n=-\infty}^m 3^{-2s(m-n)} J(\cu_n)
	&
	\leq
	\sum_{n=-\infty}^m 3^{-2s(m-n)} \sum_{k=-\infty}^{n-1} 3^{-d(n-1-k)} \max_{z \in 3^{k}\Zd \cap (\cu_{k+1} \setminus \cu_{k})} J(z + \cu_{k})
	\\ &
	\leq
	3^{2s}\sum_{k=-\infty}^{m-1} 3^{-2s(m-k)} \max_{z \in 3^{k}\Zd \cap (\cu_{k+1} \setminus \cu_{k})} J(z + \cu_{k}) \sum_{n=k+1}^{m} 3^{- (d-2s)(n-1-k)}
	\\ &
	\leq
	\frac32\,3^{2s} \sum_{k=-\infty}^{m-1} 3^{-2s(m-k)} \max_{z \in 3^{k}\Zd \cap (\cu_{k+1} \setminus \cu_{k})} J(z + \cu_{k})
	\,.
\end{align*}
After relabeling~$k \to n$ and~$k+1 \to j$, the sum on the right becomes
\begin{equation*}
	\sum_{j=-\infty}^{m} 3^{-2s(m-j)} \max_{z \in 3^{j-1}\Zd \cap (\cu_{j} \setminus \cu_{j-1})} J(z + \cu_{j-1})
	\leq
	\frac{C}{s} \sup_{j \leq m}3^{-\frac{3}{2} s(m-j)} \max_{z \in 3^{j-1}\Zd \cap (\cu_{j} \setminus \cu_{j-1})} J(z + \cu_{j-1})
	\,,
\end{equation*}
which is of the form appearing on the right side of~\eqref{e.mathcalE.annular.decomp.pre} with~$n = j-1$. We obtain~\eqref{e.mathcalE.annular.decomp.pre} by combining the previous displays with subadditivity in the form
\begin{equation*}
	\max_{z \in 3^{j-1}\Zd \cap (\cu_{j} \setminus \cu_{j-1})} J(z + \cu_{j-1})
	\leq
	\max_{z \in 3^{j-2}\Zd \cap (\cu_{j} \setminus \cu_{j-1})} J(z + \cu_{j-2})
	\,.
\end{equation*}

\smallskip

\emph{Step 2.} We first treat the cubes of side at most one.  By~\eqref{e.direct.J.bound.subunit}, with~$b=(\b_{L,m})_{\cu_m}$, the contribution of the indices~$n\leq0$ to the left side of~\eqref{e.mathcalE.annular.decomp.pre} is at most
\begin{equation}
	C3^{-2sm}
	\left\|
	\frac{\a_L}{(\b_{L,m})_{\cu_m}}
	+
	\frac{(\b_{L,m})_{\cu_m}}{\a_L}
	-2
	\right\|_{L^\infty(\cu_m)}\,.
	\label{e.good.scale.subunit.tail}
\end{equation}
On~$\mathcal G_m(1,s)$, the relations~\eqref{e.combined.coefficient.ratio.reg},~\eqref{e.rho.m.n.bound}, and~\eqref{e.long.ratio.reg} established below, together with~\eqref{e.good.discounted.am.vs.aj}, bound~\eqref{e.good.scale.subunit.tail} by the last four terms on the right side of~\eqref{e.J.localization.reg} below.  It remains to bound the contribution of the cubes of side larger than one.

\smallskip

We estimate only~$n\leq j-2$, the range used in Step~3. The scale-$j-1$ cubes were subdivided at the end of Step~1 with the coefficient~$\tilde\a_{L,m}$ unchanged; sensitivity is applied only after that subdivision. We show that there exists~$C(d)<\infty$ such that
\begin{align}
	\lefteqn{
	\sup_{1\leq j \leq m} \sup_{1\leq n \leq j-2}
	3^{-\frac{3}{2} s(m-n)}
	\max_{z \in 3^n \Zd \cap (\cu_j \setminus \cu_{j-1})}
	\max_{|e| = 1} J(z+\cu_n,e,e\,;\tilde{\a}_{L,m})
	\indc_{\mathcal G_m(1,s)}
	} \qquad &
	\notag \\ &
	\leq
	C \sup_{0\leq j \leq m} \sup_{0\leq n \leq j-2} 3^{-s(m-n)}
	\max_{z \in 3^n \Zd \cap (\cu_j \setminus \cu_{j-1})}
	\max_{|e| = 1} J(z+\cu_n,e,e\,; \ahom_{n}^{-1} \a_{n})
	\notag \\ & \qquad
	+
	C s^{-2} \gmcdelta^4
	\notag \\ & \qquad
	+
	C \biggl(
	\sum_{n=m}^\infty
	3^m \| \nabla \g_n \|_{L^\infty(\cu_m)}
	\biggr)^2
	\notag \\ & \qquad
	+
	C \sup_{0\leq n \leq m}
	3^{-s(m-n)}\|\g_{m,n}\|_{L^\infty(\cu_m)}^2
	\notag \\ & \qquad
	+
	C3^{-sm}\left\|\sum_{i=0}^m\g_i\right\|_{L^\infty(\cu_m)}^2
	\,.
	\label{e.J.localization.reg}
\end{align}
We apply Lemma~\ref{l.J.sensitivity} with~$\a=\ahom_n^{-1}\a_n$ and~$\b=\tilde\a_{L,m}$. For every~$n\in\N$ with~$n<m$, this gives
\begin{align}
	\lefteqn{
	J(z+\cu_n,e, e \, ; \tilde{\a}_{L,m})
	} \qquad &
	\notag \\ &
	\leq
	2J(z+\cu_n,e,e\,;\ahom_{n}^{-1} \a_{n})
	\notag \\ &  \qquad
	+
	3\Biggl(
	\biggl\|\frac{\ahom_{n}^{-1} \a_{n}} {\tilde{\a}_{L,m}} - 1 \biggr\|_{L^\infty(z+\cu_n)}^2
	+
	\biggl\|\frac{\tilde{\a}_{L,m}} {\ahom_{n}^{-1} \a_{n}} - 1 \biggr\|_{L^\infty(z+\cu_n)}^2
	\Biggr)
	\bigl( J(z+\cu_n,e,e\,; \ahom_{n}^{-1} \a_{n})   + 1 \bigr)\,.
	\label{e.J.sensitivity.applied.reg}
\end{align}
By~\eqref{e.a.m.j.ratio},
\begin{equation}
	\frac{\tilde{\a}_{L,m}}{\ahom_n^{-1}\a_n}
	=
	\frac{\a_L\a_m^{-1}}{(\a_L\a_m^{-1})_{\cu_m}}
	\exp\bigl(\g_{m,n}+\log(\ahom_n/\ahom_m)-(m-n)\tau^2\bigr)\,.
	\label{e.combined.coefficient.ratio.reg}
\end{equation}
By~\eqref{e.annealed.ordering},
\begin{equation}
	|\bigl(\log(\ahom_n/\ahom_m)-(m-n)\tau^2\bigr)|\leq\tau^2(m-n)\,.
	\label{e.rho.m.n.bound}
\end{equation}
Let~$S_m\coloneqq\sum_{j=m}^{\infty}3^m\|\nabla\g_j\|_{L^\infty(\cu_m)}$. Conditions~\eqref{e.good.gj} and~\eqref{e.good.discounted.am.vs.aj} in the definition of~$\mathcal G_m(1,s)$, the mean-value theorem, and~$|\exp(r)-1|\leq |r|\exp(|r|)$ give
\begin{equation}
	\left(
	\left\|\frac{\a_L\a_m^{-1}}{(\a_L\a_m^{-1})_{\cu_m}}-1\right\|_{L^\infty(\cu_m)}
	+
	\left\|\frac{(\a_L\a_m^{-1})_{\cu_m}}{\a_L\a_m^{-1}}-1\right\|_{L^\infty(\cu_m)}
	\right)\indc_{\mathcal G_m(1,s)}
	\leq C\min\{1,S_m\}\,,
	\label{e.long.ratio.reg}
\end{equation}
and also~$S_m\indc_{\mathcal G_m(1,s)}\leq C$.  Moreover,~\eqref{e.good.discounted.am.vs.aj},~\eqref{e.rho.m.n.bound},~$\tau^2\leq4\gmcdelta^2$, and~$s\geq64\gmcdelta^2$ imply
\begin{equation*}
	\exp\bigl(\|\g_{m,n}\|_{L^\infty(z+\cu_n)}\bigr)
	\indc_{\mathcal G_m(1,s)}
	\leq C3^{\frac{s}{8}(m-n)}\,,
	\qquad
	\exp(|\bigl(\log(\ahom_n/\ahom_m)-(m-n)\tau^2\bigr)|)\leq3^{\frac{s}{16}(m-n)}\,.
\end{equation*}
Applying these bounds to~\eqref{e.combined.coefficient.ratio.reg} and its reciprocal yields
\begin{multline*}
	\biggl\|\frac{\ahom_{n}^{-1} \a_{n}} {\tilde{\a}_{L,m}} - 1 \biggr\|_{L^\infty(z+\cu_n)}
	+
	\biggl\|\frac{\tilde{\a}_{L,m}} {\ahom_{n}^{-1} \a_{n}} - 1 \biggr\|_{L^\infty(z+\cu_n)}
	\\
	\leq C3^{\frac{3}{16}s(m-n)} \min\biggl\{ 1,S_m
	+
	\| \g_{m,n} \|_{L^\infty(z+\cu_n)}
	+
	\tau^2(m-n)
	\biggr\}
	\,.
\end{multline*}
On~$\mathcal G_m(1,s)$, condition~\eqref{e.good.sum.J} bounds the term~$J(z+\cu_n,e,e\,;\ahom_n^{-1}\a_n)$ in~\eqref{e.J.sensitivity.applied.reg} by~$3^{\frac{s}{8}(m-n)}$. Hence the square of the right side of the last display, multiplied by~$J(z+\cu_n,e,e\,;\ahom_n^{-1}\a_n)+1$ and by~$3^{-\frac{3s}{2}(m-n)}$, is bounded by
\begin{equation*}
	C3^{-s(m-n)}
	\left(
	S_m^2
	+
	\|\g_{m,n}\|_{L^\infty(z+\cu_n)}^2
	+
	\bigl(\tau^2(m-n)\bigr)^2
	\right)\,.
\end{equation*}
Finally,~$\sup_{r\geq0}3^{-sr}(\tau^2r)^2\leq Cs^{-2}\gmcdelta^4$. Combining these estimates with~\eqref{e.J.sensitivity.applied.reg} proves~\eqref{e.J.localization.reg}.

\smallskip

\emph{Step 3.} We conclude. Combining~\eqref{e.mathcalE.annular.decomp.pre}, the estimate~\eqref{e.good.scale.subunit.tail} for cubes of side at most one, and~\eqref{e.J.localization.reg} gives~\eqref{e.mathcalE.bound.for.good}.
\end{proof}

For each of the three conditions~\eqref{e.good.gj},~\eqref{e.good.discounted.am.vs.aj}, and~\eqref{e.good.sum.J} separately, we estimate the probability that it fails on a given proportion of the scales.

\begin{lemma}
\label{l.bad.scales.for.nabla.gj}
There exists a constant~$C(d)<\infty$ such that, for every~$\theta,s,\ep \in (0,1]$ satisfying~$ \theta \geq C s^{-6} \ep^{-2} \gmcdelta^2$, we have, for every~$m_0 \in \N_0$ and~$M \in\N_0$,
\begin{equation}
	\P \Biggl[ \avsum_{m=m_0}^{m_0 + M} \indc{\biggl\{
	\sup_{j \geq 0}  3^{-s j}
	\sum_{i=(m-j)\vee0}^{m+j} \bigl\|  | \g_i |
	+
	3^{i} | \nabla \g_{i}| \bigr\|_{L^\infty(\cu_{m+1+j})}
	>
	\ep
	\biggr\}} \geq \theta \Biggr]
	\leq
	\exp \biggl(
	- \frac
	{  s^6 \ep^2 \theta}
	{C \gmcdelta^2}(M+1)
	\biggr)
	\,.
	\label{e.bad.scales.for.nabla.gj}
\end{equation}
\end{lemma}

\begin{lemma}
\label{l.exp.goodscales}
There exists a constant~$C(d)<\infty$ such that, for every~$m_0,M \in \N_0$ and~$s,\theta\in (0,1]$ satisfying~$\theta \geq C s^{-2} \gmcdelta^{2} $, we have the estimate
\begin{multline}
	\P \Biggl[
	\avsum_{m=m_0}^{m_0+M}
	\indc\biggl\{
	\sup_{j \geq 0}  3^{- s j}
	\biggl\| \, \prod_{i=(m-j)\vee0}^{m+j} \exp\bigl(| \g_i |\bigr)
	+
	\prod_{i=m+j}^{\infty} \exp\bigl(4 | \g_i - \g_i(0) | \bigr)
	\, \biggr\|_{L^\infty(\cu_{m+1+j})}
	>
	6
	\biggr\}
	\geq
	\theta
	\Biggr]
	\\
	\leq
	\exp\biggl( - \frac{s^2\theta}{C \gmcdelta^{2}} (M+1)  \biggr)
	\,.
	\label{e.amaj.good.scales.one}
\end{multline}
\end{lemma}

\begin{lemma}
\label{l.discounted.J.local}
There exists a constant~$C(d)<\infty$ such that, if~$\ep,s,\theta \in (0,1]$ satisfy the conditions
\begin{equation*}
	\gmcdelta^{2} \left| \log \gmcdelta\right|
	\leq
	C^{-1} s^3 \ep^2
	\qand
	\theta \geq C s^{-3} \ep^{-2} \gmcdelta^{2} \left| \log \gmcdelta\right| \,,
\end{equation*}
then we have, for every~$m_0 \in \N_0$ and~$M \in \N_0$, the estimate
\begin{multline}
	\P\Biggl[
	\avsum_{m=m_0}^{m_0+M} \indc\biggl\{ \sup_{0\leq j \leq m} \sup_{0\leq n \leq j-2} 3^{-s(m-n)}
	\max_{z \in 3^n \Zd \cap (\cu_j \setminus \cu_{j-1})}
	\max_{|e| = 1} J(z+\cu_n,e,e\,; \ahom_{n}^{-1} \a_{n})  > \ep^2 \biggr\} > \theta
	\Biggr]
	\\
	\leq
	\exp\biggl(- \frac{s^{3} \ep^2  \theta}{C\left| \log \gmcdelta \right| \gmcdelta^2  }  (M+1)  \biggr)
	\,.
	\label{e.discounted.J.local.lemma}
\end{multline}
\end{lemma}

\begin{proof}[{Proof of Proposition~\ref{p.density.of.good.scales}}]
Fix~$s,\theta,\ep\in(0,1]$ and assume~\eqref{e.theta.vs.gmcdelta.cond.reg} with its constant~$C$ large enough that Lemmas~\ref{l.bad.scales.for.nabla.gj},~\ref{l.exp.goodscales}, and~\ref{l.discounted.J.local} are applicable with parameters~$(\nf{s}{8},\nf{\theta}{4},\ep)$ in place of~$(s,\theta,\ep)$. If the good event fails on at least a proportion~$\theta$ of the scales, then one of the three conditions~\eqref{e.good.gj}--\eqref{e.good.sum.J} fails on at least a proportion~$\theta/3$ of them. Since~$\theta/3>\theta/4$, this case is covered by the strict inequality in Lemma~\ref{l.discounted.J.local}. A union bound over the three conditions gives~\eqref{e.good.H.one}.
\end{proof}

\begin{proof}[{Proof of Lemma~\ref{l.bad.scales.for.nabla.gj}}]
Interchanging the order of summation gives
\begin{align*}
	\lefteqn{
	\sup_{j \geq 0}
	3^{-sj}
	\sum_{i = (m-j)\vee0}^{m+j}
	\bigl\|  | \g_i |
	+
	3^{i} | \nabla \g_{i}| \bigr\|_{L^\infty(\cu_{m+1+j})}
	} \qquad &
	\notag \\ &
	\leq
	\sum_{j = 0}^\infty
	3^{-sj}
	\sum_{i = (m-j)\vee0}^{m+j}
	\bigl\|  | \g_i |
	+
	3^{i} | \nabla \g_{i}| \bigr\|_{L^\infty(\cu_{m+1+j})}
	\notag \\ &
	=
	\sum_{i \in \N_0}
	3^{-s|i-m|}
	\biggl(
	\sum_{j = |i-m|}^\infty
	3^{-s(j - |i-m|)}
	\bigl\|  | \g_i |
	+
	3^{i} | \nabla \g_{i}| \bigr\|_{L^\infty(\cu_{m+1+j})}
	\biggr)
	\,.
\end{align*}
For a universal constant~$C_2$ to be fixed and~$i\in\N_0$, set
\begin{equation*}
	X_{m,i}
	\coloneqq
	C_2 s^{-1} \ep^{-1} 3^{-\frac{1}{2}s|m-i|}
	\sum_{j = |i-m|}^\infty
	3^{-s(j - |i-m|)}
	\bigl\|  | \g_i |
	+
	3^{i} | \nabla \g_{i}| \bigr\|_{L^\infty(\cu_{m+1+j})}
	\,,
\end{equation*}
and let~$X_{m,i}=0$ if~$m<0$ or~$i<0$. We shall apply Proposition~\ref{p.concentration.for.scales} to this extension to all indices in~$\Z$. We have
\begin{equation*}
	\sup_{j \geq 0}
	3^{-sj}
	\sum_{i = (m-j)\vee0}^{m+j}
	\bigl\|  | \g_i |
	+
	3^{i} | \nabla \g_{i}| \bigr\|_{L^\infty(\cu_{m+1+j})}
	\leq
	C_2^{-1} s \ep  \sum_{i \in \N_0}
	3^{-\frac{1}{2} s| m-i|} X_{m,i}
	\,.
\end{equation*}
By~\eqref{e.max.scaling.O.Gamma.Two},~\eqref{e.same.law.shift}, and~\eqref{e.coefficient.field.regularity},
\begin{equation*}
	\bigl\|  | \g_i |
	+
	3^{i} | \nabla \g_{i}| \bigr\|_{L^\infty(\cu_{m+1+j})}
	\leq
	\O_{\Gamma_2}\bigl( C\gmcdelta (|m- i| + j + 1)^{\nf 12} \bigr)
	\,.
\end{equation*}
By~\eqref{e.moments.OGamma2}, there exists a constant~$C_1(d)$ such that, for every~$p \in [1,\infty)$, we have~$\E\bigl[ X_{m,i}^p  \bigr]^{\nf 1p} \leq C_1 C_2 s^{-\nf 52} \ep^{-1} \gmcdelta p^{\nf 12}$. Moreover, for~$m,m',i,i' \in \N_0$ with~$|i-i'| \geq 1$, the variables~$X_{m,i}$ and~$X_{m',i'}$ are independent. We fix~$p \coloneqq C_1^{-2} C_2^{-2} s^5 \ep^2 \gmcdelta^{-2}$.  With this choice of~$p$, if~$C(d)$ is large enough in terms of~$C_1$ and~$C_2$, the assumption~$\theta \geq C s^{-6} \ep^{-2} \gmcdelta^{2}$, which since~$\theta\leq1$ also gives~$\gmcdelta \leq C^{-\nf12} s^{3} \ep$, implies that~$\theta \geq p^{-1}$ and~$s \geq 2 p^{-1}$. By the inequality displayed after the definition of~$X_{m,i}$, we then have
\begin{equation*}
	\biggl\{
	\sup_{j \geq 0}
	3^{-sj}
	\sum_{i = (m-j)\vee0}^{m+j}
	\bigl\|  | \g_i |
	+
	3^{i} | \nabla \g_{i}| \bigr\|_{L^\infty(\cu_{m+1+j})} >  \ep
	\biggr\}
	\subseteq
	\biggl\{
	\sum_{i \in \N_0} 3^{-\frac{1}{2} s|m-i|} X_{m,i}  >  C_2 s^{-1} \biggr\}
	\,.
\end{equation*}
For~$C_2$ large enough, Proposition~\ref{p.concentration.for.scales} gives~\eqref{e.bad.scales.for.nabla.gj}.
\end{proof}

\begin{proof}[Proof of Lemma~\ref{l.exp.goodscales}]
We apply Propositions~\ref{p.concentration.for.scales.exp.sequence} and~\ref{p.concentration.for.scales.exp.field}. For sufficiently large~$C(d)$, let~$X_n\coloneqq C3^n\|\nabla\g_n\|_{L^\infty(\cu_{n+1})}$ for~$n\in\N_0$ and~$X_n\coloneqq0$ for~$n<0$. On~$\cu_{k+1+j}$,
\begin{align*}
	\prod_{i=k+j}^{\infty} \exp\bigl(4 | \g_i - \g_i(0) | \bigr)
	&
	\leq
	\prod_{i=k+j}^{\infty} \exp\Bigl( 3^{k+j-i}
	\bigl( C 3^i \| \nabla \g_i \|_{L^\infty(\cu_{i+1})} \bigr)\Bigr)
	\leq
	\prod_{i=0}^{\infty}
	\exp( 3^{-i} X_{k+i+j})
	\,.
\end{align*}
The sequence~$(X_n)$ satisfies the hypotheses of Proposition~\ref{p.concentration.for.scales.exp.sequence} with~$\mathsf a_{\mathrm{sg}}\coloneqq C\gmcdelta$. If~$\gmcdelta\leq C^{-1}s\theta^{\nf12}$, then that proposition gives
\begin{equation*}
	\P \Biggl[
	\avsum_{m=m_0}^{m_0+M}
	\indc\biggl\{
	\sup_{j \geq 0}  3^{- s j}
	\biggl\|
	\prod_{i=m+j}^{\infty} \exp\bigl(4 | \g_i - \g_i(0) | \bigr)
	\, \biggr\|_{L^\infty(\cu_{m+1+j})}
	>
	3
	\biggr\}
	\geq
	\theta
	\Biggr]
	\leq
	\exp\biggl( - \frac{s^2\theta}{C \gmcdelta^{2}} (M+1)  \biggr)  \,.
\end{equation*}
If~$\gmcdelta \leq C^{-1}s \theta^{\nf12}$ for sufficiently large~$C(d)$, then Proposition~\ref{p.concentration.for.scales.exp.field}, applied with~$\g_i$ in place of~$X_i$ for~$i\in\N_0$ and with~$0$ in place of~$X_i$ for~$i<0$, yields
\begin{equation*}
	\P \Biggl[ \avsum_{k=m_0}^{m_0+M} \! \! \! \indc{\biggl\{ \sup_{j \geq 0}
	3^{- s j }
	\biggl\| \, \prod_{i=(k-j)\vee0}^{k+j} \exp\bigl(| \g_i |\bigr) \, \biggr\|_{L^\infty(\cu_{k+1+j})}  > 3 \biggr\}} \geq \theta \Biggr]
	\leq
	\exp \biggl( - \frac{s^2 \theta}{C \gmcdelta^2} (M+1) \biggr)
	\,.
\end{equation*}
Apply the previous two displays with~$\theta/2$ in place of~$\theta$. If the sum of the two products exceeds six, then at least one exceeds three, so a union bound gives~\eqref{e.amaj.good.scales.one}.
\end{proof}

\begin{proof}[Proof of Lemma~\ref{l.discounted.J.local}]
First,
\begin{align*}
	\lefteqn{
	\sup_{0\leq j \leq m} \sup_{0\leq n \leq j-2} 3^{-s(m-n)}
	\max_{z \in 3^n \Zd \cap (\cu_j \setminus \cu_{j-1})}
	\max_{|e| = 1} J(z+\cu_n,e,e\,; \ahom_{n}^{-1} \a_{n})
	} \qquad &
	\notag \\ &
	\leq
	\sum_{j = 0}^m 3^{-s(m-j)}
	\sum_{n = 0}^{j-2} 3^{-s(j-n)}
	\max_{z \in 3^n \Zd \cap (\cu_j \setminus \cu_{j-1})}
	\max_{|e| = 1} J(z+\cu_n,e,e\,; \ahom_{n}^{-1} \a_{n})
	\,.
\end{align*}
We set, for a universal constant~$K \geq 1$ to be fixed,
\begin{equation*}
	X_{m,j}
	\coloneqq
	K s^{-1} \ep^{-2}
	\sum_{n = 0}^{j-2} 3^{-s(j-n)}
	\max_{z \in 3^n \Zd \cap (\cu_j \setminus \cu_{j-1})}
	\max_{|e| = 1} J(z+\cu_n,e,e\,; \ahom_{n}^{-1} \a_{n})
	\indc_{\{ j \leq m\}}
	\,.
\end{equation*}
Then
\begin{multline*}
	\biggl\{ \sup_{0\leq j \leq m} \sup_{0\leq n \leq j-2} 3^{-s(m-n)}
	\max_{z \in 3^n \Zd \cap (\cu_j \setminus \cu_{j-1})}
	\max_{|e| = 1} J(z+\cu_n,e,e\,; \ahom_{n}^{-1} \a_{n})  \geq \ep^2 \biggr\}
	\\
	\subseteq
	\biggl\{
	\sum_{j = 0}^m 3^{-s(m-j)}
	X_{m,j}
	\geq K s^{-1} \biggr\}
	\,.
\end{multline*}
There is an integer~$r(d)<\infty$ such that the columns~$(X_{m,j})_m$, indexed by~$j$, are~$r(d)$-dependent in the sense of~\eqref{e.independent.columns.twosided}. To verify this, fix one of the fields~$\g_i$. Column~$j$ depends on~$\g_i$ only if~$i\leq j-2$, and then only through its values on the annulus~$\cu_j\setminus\cu_{j-1}$. If~$j'>j$ and~$j'-j\geq r(d)$, then the distance between these annuli is at least~$(3^{j'-1}-3^j)/2>\sqrt d\,3^{j-2}\geq\sqrt d\,3^i$, provided that~$r(d)$ is large enough. Thus, by Assumption~\ref{a.g1}, the restrictions of~$\g_i$ to the regions used by two sets of columns at distance at least~$r(d)$ are independent. Since the fields are mutually independent, so are the two sets of columns. Moreover, by the triangle inequality, stationarity, and Proposition~\ref{p.coarse.grained.bound}, there is a constant~$C_1(d)<\infty$ such that, if~$p$ satisfies~$2ds^{-1} \leq p \leq C_1^{-1} \gmcdelta^{-2} \left| \log \gmcdelta\right|^{-1}$, then
\begin{align*}
	\E\bigl[ X_{m,j}^p \bigr]^{\nf 1p}
	&
	\leq
	K s^{-1} \ep^{-2}
	\sum_{n = 0}^{j-2} 3^{-s(j-n)} \E\biggl[
	\max_{z \in 3^n \Zd \cap (\cu_j \setminus \cu_{j-1})}
	\max_{|e| = 1} J(z+\cu_n,e,e\,; \ahom_{n}^{-1} \a_{n})^p  \biggr]^{\nf 1p}
	\notag \\ &
	\leq
	K s^{-1} \ep^{-2}
	\sum_{n = 0}^{j-2} 3^{-(s- \frac{d}{p})(j-n)}
	\E\Bigl[
	\max_{|e| = 1} J(\cu_n,e,e\,; \ahom_{n}^{-1} \a_{n})^p  \Bigr]^{\nf 1p}
	\leq
	C_1 K s^{-2} \ep^{-2}\gmcdelta^2  \left| \log \gmcdelta\right|  p
	\,.
\end{align*}
Taking~$p \coloneqq C_1^{-1} K^{-1} s^{2} \ep^{2} \gmcdelta^{-2} \left| \log \gmcdelta\right|^{-1} $ and assuming that~$\gmcdelta$ satisfies
\begin{equation*}
	\gmcdelta^{2} \left| \log \gmcdelta\right|
	\leq
	\frac{s^3 \ep^2}{2d C_1 K}
	\quad \implies \quad
	2ds^{-1} \leq p = C_1^{-1} K^{-1} s^{2} \ep^{2} \gmcdelta^{-2} \left| \log \gmcdelta\right|^{-1}
	\,,
\end{equation*}
we deduce that~$\E[ X_{m,j}^p ]^{\nf 1p} \leq 1$. Assuming also~$\theta \geq p^{-1}$ and applying Proposition~\ref{p.concentration.for.scales} with~$r=r(d)$, we obtain~\eqref{e.discounted.J.local.lemma} for a sufficiently large universal constant~$K$.
\end{proof}

\subsection{Excess decay on good scales}
\label{ss.good.scale.estimates}

\begin{lemma}[Centered energy at a good scale]
\label{l.centered.energy.good.scale}
Let~$s,L,m,n,z,x,u,h,\g$ satisfy the hypotheses of Lemma~\ref{l.harmonic.approximation.good.scales.GMC}, and let~$\a_0=(\b_{L,n+2})_{z+\cu_{n+2}}$. On~$\mathcal G_{n+2,z}(1,\nf s8)$, for every~$q\in (x+\cu_{n-1})\cap\cu_m$,
\begin{align}
	\a_0^{-1}\fint_{(q+\cu_{n-4})\cap\cu_m}\a_L|\nabla u|_2^2
	&\leq C3^{-2n}\|u-(u)_{(x+\cu_n)\cap\cu_m}\|_{\underline L^2((x+\cu_n)\cap\cu_m)}^2
	\notag \\ & \qquad
	+C|(\nabla h)_{(x+\cu_n)\cap\cu_m}|_2^2
	\indc_{\{\overline{(x+\cu_n)\cap\cu_m}\cap\partial\cu_m\neq\varnothing\}}
	\notag \\ & \qquad
	+Cs^{-12}\a_0^{-2}3^{2sn}[\g]_{\underline H^s((x+\cu_n)\cap\cu_m)}^2
	\notag \\ & \qquad
	+Cs^{-4}3^{2sn}[\nabla h]_{\underline H^s((x+\cu_n)\cap\cu_m)}^2
	\indc_{\{\overline{(x+\cu_n)\cap\cu_m}\cap\partial\cu_m\neq\varnothing\}}
	\,.
	\label{e.centered.energy.good.scale}
\end{align}
Here~$C$ depends only on~$d$. The same estimate holds with the cutoff modification of the good event and with the hypothesis~$L\geq m$ removed.
\end{lemma}

\begin{proof}
We give the boundary argument with its centering scalar retained. The proof tracks the energy~$F(r)$ of~$u$ on the windows around~$x$ defined below. We begin with ellipticity bounds for every cube contained in~$z+\cu_{n+2}$, including cubes off the triadic grid. Next, for each cube of a cover scale~$k\leq n-4$ whose parent cube has a face on~$\partial\cu_m$, we compare~$u$ on this parent with the solution of the same equation with boundary data~$h$; this gives an energy estimate in which the energy of~$u$ on the parent has the small coefficient~$\theta$. Interior cubes use the coarse-grained Caccioppoli inequality of Proposition~\ref{p.coarse.grained.caccioppoli}. The heart of the matter is a parent with a face on~$\partial\cu_m$: there the Caccioppoli estimate for the difference of~$u$ and this solution contains the mean~$\mu$ of the difference on the parent, which we bound through slabs along that face, at a depth chosen so that the coefficient of the energy of~$u$ on the parent is at most~$\theta$. Finally, we sum over the cubes, whose overlap is bounded, and iterate the resulting inequality for~$F$. All energies in the proof are for~$\a_L$; the comparison scalar is always the fixed~$\a_0$ above.

\smallskip

First we record the ellipticity bounds needed below the scale~$n$. For every translated cube~$y+\cu_k$ contained in~$z+\cu_{n+2}$, including~$k<0$, stability gives
\begin{equation}
	\mathcal E_{s/6,\infty,2}(y+\cu_k;\a_L,\a_0)
	\leq C3^{(s/8)(n+2-k)}
	\mathcal E_{s/8,\infty,2}(z+\cu_{n+2};\a_L,\a_0)\,.
	\label{e.boundary.energy.offgrid.error}
\end{equation}
For completeness, the arbitrary-depth form of stability follows by partitioning an off-grid cube into the maximal cubes of the anchor's triadic grid that it contains. At depth~$j$, their total volume is at most~$C(d)3^{-j}$ times the volume of the off-grid cube: a maximal cube has its parent meeting the boundary, so it lies in a boundary layer of width~$C3^{-j}$ times the side length. Subadditivity of~$J$ and~\eqref{e.bound.one.cube.by.mathcalE} bound the value of~$J$ on an off-grid cube of scale~$k-\ell$ by
\begin{equation*}
	\frac{C(d)}{1-2t'}\,
	3^{2t'(n+2-k+\ell)}
	\mathcal E_{t',\infty,2}(z+\cu_{n+2};\a_L,\a_0)^2\,,
	\qquad 0<t'<\nf12\,.
\end{equation*}
The maximal cubes cover up to a null set. To justify countable subadditivity, restrict each admissible harmonic field in the variational definition of~$J$ to these open subcubes. The integrand of~$J$ for this field is integrable, since the fixed finite cutoff is uniformly elliptic on the anchor. The integral equals the sum over the subcubes, and each summand is bounded by the volume of that subcube times the value of~$J$ on it. The displayed summable majorant permits this comparison before taking the supremum over the admissible fields. Apply this to both scalar parts of the doubled quantity~$J$ and then to its unit-sphere maximum. Summing with the normalized weights~$(1-3^{-2t})3^{-2t\ell}$ gives
\begin{equation*}
	\mathcal E_{t,\infty,2}(y+\cu_k;\a_L,\a_0)^2
	\leq \frac{C(d)t}{(t-t')(1-2t')}
	3^{2t'(n+2-k)}
	\mathcal E_{t',\infty,2}(z+\cu_{n+2};\a_L,\a_0)^2\,.
\end{equation*}
Taking~$(t,t')=(s/6,s/8)$ proves~\eqref{e.boundary.energy.offgrid.error}, with a dimensional constant and no inverse power of~$s$. On the stated event, the anchor error is bounded by~$C$. By~\eqref{e.bound.Lambdas.by.Es}, monotonicity in the order, and~\eqref{e.compareLambdaqs}, this gives
\begin{equation}
	\max\left\{
	\a_0^{-1}\Lambda_{s/6,2}(y+\cu_k),\,
	\a_0\lambda_{s/6,2}^{-1}(y+\cu_k),\,
	\a_0^{-1}\Lambda_{1/2,1}(y+\cu_k),\,
	\a_0\lambda_{s/3,1}^{-1}(y+\cu_k)
	\right\}
	\leq C3^{(s/4)(n+2-k)}\,.
	\label{e.boundary.energy.offgrid.caps}
\end{equation}

\smallskip

We work at a free cover scale~$k\leq n-4$. For a center~$q$, let~$\widehat q+\cu_k$ be the coordinatewise clamping of~$q+\cu_k$ inside~$\cu_m$. The set covered is~$(q+\cu_{k-2})\cap\cu_m$. Its scale-$(k-1)$ patch meets no artificial face of~$\widehat q+\cu_k$. If this patch leaves~$\cu_m$, then~$\widehat q+\cu_k$ has a face on~$\partial\cu_m$; this is the boundary case. Write~$E_{\widehat q}(w)=\fint_{\widehat q+\cu_k}\a_L|\nabla w|_2^2$.

\smallskip

Let~$w_{\widehat q}$ solve the same forced equation as~$u$ in~$\widehat q+\cu_k$, with boundary data~$h$ on all of~$\partial(\widehat q+\cu_k)$. Its energy satisfies
\begin{align}
	E_{\widehat q}(w_{\widehat q})
	&\leq C3^{s(n+2-k)/4}\Bigl(
	\a_0|(\nabla h)_{\widehat q+\cu_k}|_2^2
	+\a_0s^{-4}3^{2sk}[\nabla h]_{\underline H^s(\widehat q+\cu_k)}^2
	\notag \\ & \qquad
	+\a_0^{-1}s^{-12}3^{2sk}[\g]_{\underline H^s(\widehat q+\cu_k)}^2
	\Bigr)\,.
	\label{e.boundary.energy.competitor}
\end{align}
Here is the separation of the three terms. Split~$h$ into an affine function with slope~$(\nabla h)_{\widehat q+\cu_k}$ and a remainder whose gradient has mean zero. The affine Dirichlet minimum is bounded by~$|(\nabla h)_{\widehat q+\cu_k}|_2^2\Lambda_{s/6,2}(\widehat q+\cu_k)$, by the definition of the coarse matrix and its ordering. The centered remainder is controlled by the Dirichlet energy estimate~\citep[Lemma~2.11, (2.82)]{ASD}; its squared datum term is at most~$Cs^{-3}\Lambda_{s/2,2}(\widehat q+\cu_k)$ times the squared scale-normalized positive norm of the remainder's gradient. This norm is bounded by~$C3^{sk}s^{-1/2}[\nabla h]_{\underline H^s(\widehat q+\cu_k)}$, giving the~$s^{-4}$ term. Add the forced zero-boundary solution. The same estimate gives the squared forcing term~$Cs^{-7}\lambda_{s/2,2}^{-1}(\widehat q+\cu_k)3^{2sk}[\g]_{\underline H^s(\widehat q+\cu_k)}^2$, which is bounded by the displayed~$s^{-12}$ term. The factor~$s^{-1/2}$ in both conversions comes from the factor~$s$ in this paper's squared fractional seminorm. All constants in these comparisons are uniform for~$0<s\leq1/4$. The energy triangle inequality and~\eqref{e.boundary.energy.offgrid.caps} prove~\eqref{e.boundary.energy.competitor}. In particular, the affine slope has no inverse power of~$s$.

\smallskip

Apply Lemma~\ref{l.coarse.grained.Caccioppoli.RHS.ASD}, rescaled to~$\widehat q+\cu_k$, to~$u-w_{\widehat q}$, with zero forcing and zero boundary datum, and with orders~$1/2$ and~$s/2$ at the center~$q$. This difference has zero trace on the physical faces needed by the estimate. Equation~\eqref{e.CG.Poincare.mean.zero}, rescaled to~$\widehat q+\cu_k$, bounds the oscillation of~$w_{\widehat q}$ by its own energy. Here the Caccioppoli parameter of the lemma is~$1-1/2-s/2=(1-s)/2$, and its ellipticity factor satisfies
\begin{equation*}
	\left(\frac{\Lambda_{1/2,1}(\widehat q+\cu_k)}{\lambda_{s/2,1}(\widehat q+\cu_k)}\right)^{(2-s)/(1-s)}
	\lambda_{s/2,1}(\widehat q+\cu_k)
	\leq C\a_0 3^{17s(n+2-k)/12}\,.
\end{equation*}
Indeed,~$(2-s)/(1-s)\leq7/3$ for~$s\leq1/4$; the powers of the bound in~\eqref{e.boundary.energy.offgrid.caps} are thus at most~$17/3\leq7$ in the solution oscillation term and~$20/3\leq8$ after the extra inverse ellipticity factor from Poincar\'e in the competitor term. Thus, for every scalar~$c_0$ and with~$\mu\coloneqq(u-w_{\widehat q})_{\widehat q+\cu_k}$,
\begin{multline}
	\fint_{(q+\cu_{k-2})\cap\cu_m}\a_L|\nabla u|_2^2
	\\
	\leq C\left(
	3^{7s(n+2-k)/4}\a_0 3^{-2k}
	\bigl(\|u-c_0\|_{\underline L^2(\widehat q+\cu_k)}^2+\mu^2\bigr)
	+3^{2s(n+2-k)}E_{\widehat q}(w_{\widehat q})\right)
	\,.
	\label{e.boundary.energy.scalar.row}
\end{multline}
Indeed,~$\|u-(w_{\widehat q})_{\widehat q+\cu_k}\|_{\underline L^2(\widehat q+\cu_k)}^2 \leq\|u-c_0\|_{\underline L^2(\widehat q+\cu_k)}^2+\mu^2$, by the identity for variance and the minimizing property of~$(u)_{\widehat q+\cu_k}$.

\smallskip

To estimate~$\mu$, extend~$w_{\widehat q}-h\in H^1_0(\widehat q+\cu_k)$ by zero to~$\cu_m$ and let~$\rho=(u-h)-\operatorname{ext}_{\widehat q}(w_{\widehat q}-h)\in H^1_0(\cu_m)$. On~$\widehat q+\cu_k$ this is~$u-w_{\widehat q}$. Its zero extension across a met face vanishes on the outer half of every cube centered on that face. For such a tile~$y+\cu_\kappa$ of side~$3^\kappa$, the rescaled inequality~\eqref{e.CG.Poincare.mean.zero} and the vanishing half give
\begin{equation*}
	|(\rho)_{y+\cu_\kappa}|^2\leq
	C3^{(s/4)(n+2-\kappa)}\a_0^{-1}3^{2\kappa}
	\fint_{y+\cu_\kappa}\a_L|\nabla\rho|_2^2\,.
\end{equation*}
This applies to either sign of the met face. At a lower face the vanishing half is below the face; at an upper face it is above. The coefficient is kept in its original coordinates.

\smallskip

Tile the doubled slab adjoining one met face at depth~$j$, so~$\kappa=k-j$. Its inner half~$S$ lies in~$\widehat q+\cu_k$, and~$3^{kd}/|S|\leq C(d)3^j$. Jensen on the equal-volume tiles bounds the slab mean by the tile energies. On the inner slab, the triangle inequality between its mean and the parent mean gives
\begin{align}
	\mu^2
	&\leq C_f(d)3^{s(n+2-k+j)/4}\a_0^{-1}3^{2(k-j)}3^j
	\bigl(E_{\widehat q}(u)+E_{\widehat q}(w_{\widehat q})\bigr)
	\notag \\ & \qquad
	+C3^j\left(
	\|u-c_0\|_{\underline L^2(\widehat q+\cu_k)}^2+
	\|w_{\widehat q}-(w_{\widehat q})_{\widehat q+\cu_k}\|_{\underline L^2(\widehat q+\cu_k)}^2\right)
	\,.
	\label{e.boundary.energy.face.mean}
\end{align}
The doubled slab's energy is at most~$C3^j$ times the normalized parent energy, since its exterior half contributes zero. The last competitor oscillation is at most~$C3^{s(n+2-k)/4}\a_0^{-1}3^{2k}E_{\widehat q}(w_{\widehat q})$. This proves~\eqref{e.boundary.energy.face.mean} while keeping the scalar~$\mu$.

\smallskip

After substitution in~\eqref{e.boundary.energy.scalar.row}, the coefficient of~$E_{\widehat q}(u)$ is at most~$C_f'(d)3^{2s(n+2-k)}3^{-(1-s/4)j}$. Fix~$\theta=(16\cdot28^d)^{-1}$ and choose~$j=\left\lceil\frac{16}{15}\left(2s(n+2-k)+\max\{0,\log_3(C_f'/\theta)\}\right)\right\rceil$. Since~$1-s/4\geq15/16$, this coefficient is at most~$\theta$, and~$3^j\leq C(d)3^{(32/15)s(n+2-k)}$. The factor of the solution oscillation term is at most~$3^{(7/4+32/15)s(n+2-k)}=3^{(233/60)s(n+2-k)}\leq3^{4s(n+2-k)}$. The factors of the competitor terms, including the additional~$3^{s(n+2-k)/4}$, are at most~$C3^{12s(n+2-k)}$. Multiplying by the volume of the set covered and using that it is at most~$3^{kd}$ gives the parent-local estimate
\begin{align}
	\lefteqn{
	\int_{(q+\cu_{k-2})\cap\cu_m}\a_L|\nabla u|_2^2
	} \qquad &
	\notag \\ &
	\leq\theta\int_{\widehat q+\cu_k}\a_L|\nabla u|_2^2
	+C3^{4s(n+2-k)}\a_0 3^{-2k}\int_{\widehat q+\cu_k}|u-c_0|^2
	\notag \\ & \qquad
	+C3^{12s(n+2-k)}\Biggl(
	\a_0\int_{\widehat q+\cu_k}|\nabla h|_2^2
	+\a_0s^{-4}3^{2sk}3^{kd}[\nabla h]_{\underline H^s(\widehat q+\cu_k)}^2
	\notag \\ & \qquad\qquad
	+\a_0^{-1}s^{-12}3^{2sk}3^{kd}[\g]_{\underline H^s(\widehat q+\cu_k)}^2
	\Biggr)
	\,.
	\label{e.boundary.energy.parent.local}
\end{align}

\smallskip

We now sum before converting the parent-local terms into window terms. Denote by~$\mathsf B$ the sum of the four terms on the right of~\eqref{e.centered.energy.good.scale}, multiplied by~$\a_0$ and with their constants omitted. For~$3^{n-1}\leq r<R\leq4\cdot3^{n-2}$, choose~$k\in\Z$ so that~$9\cdot3^k\leq R-r\leq27\cdot3^k$. The scale-$(k-2)$ cubes meeting~$\cu_m\cap\overline{x+2r\cu_0}$ have their projected parents inside~$\cu_m\cap\overline{x+2R\cu_0}$ and~$(x+\cu_n)\cap\cu_m$. Their overlap is at most~$28^d$: their centers have spacing~$3^{k-2}$ and clamping moves a center by at most~$3^k/2$ in each coordinate. Cubes whose centers are outside~$\cu_m$ have zero intersection with~$\cu_m$, since the scales are integer scales of the triadic grid.

\smallskip

In the boundary class use~\eqref{e.boundary.energy.parent.local} with~$c_0=(u)_{(x+\cu_n)\cap\cu_m}$; in the interior class use Proposition~\ref{p.coarse.grained.caccioppoli}, rescaled to the patch~$q+\cu_{k-1}$, with the same~$c_0$: apply it to~$u-c_0$ with the orders~$1/2$ and~$s/2$ and the forcing~$\g$ minus its mean on the patch. In the boundary class,~$\overline{(x+\cu_n)\cap\cu_m}$ meets~$\partial\cu_m$, so both datum terms have precisely the indicator in~$\mathsf B$. Bounded overlap controls the integrals of~$|u-c_0|^2$ and~$|\nabla h|_2^2$. For every field~$F$, positivity of the fractional kernel gives~$\sum_q3^{kd}[F]_{\underline H^s(\widehat q+\cu_k)}^2 \leq28^d|(x+\cu_n)\cap\cu_m|[F]_{\underline H^s((x+\cu_n)\cap\cu_m)}^2$. Finally,
\begin{equation*}
	\|\nabla h\|_{\underline L^2((x+\cu_n)\cap\cu_m)}^2
	\leq |(\nabla h)_{(x+\cu_n)\cap\cu_m}|_2^2+
	Cs^{-1}3^{2sn}[\nabla h]_{\underline H^s((x+\cu_n)\cap\cu_m)}^2\,.
\end{equation*}
This is the variance identity followed by the diameter bound in the fractional integral; every side of~$(x+\cu_n)\cap\cu_m$ is at least~$3^n/2$. Its fractional term is bounded by the~$s^{-4}$ datum term.

\smallskip

With~$k$ as chosen above, the oscillation transport is~$3^{4s(n-k+2)+2(n-k)}\leq C3^{3(n-k)}$. The fractional transports are~$3^{12s(n-k+2)-2s(n-k)}\leq C3^{3(n-k)}$, and the datum mean transport is~$3^{12s(n-k+2)}\leq C3^{3(n-k)}$. All three inequalities use~$s\leq1/4$. Therefore, with~$F(r)=\int_{\cu_m\cap\overline{x+2r\cu_0}}\a_L|\nabla u|_2^2$,
\begin{equation}
	F(r)\leq\frac1{16}F(R)
	+\frac{C3^{n(d+3)}\mathsf B}{(R-r)^3}\,.
	\label{e.boundary.energy.window.step}
\end{equation}
The total energy bounds~$F$. Iterating~\eqref{e.boundary.energy.window.step} with~$r$ and~$R$ consecutive radii~$3^{n-1}+3^{n-2}(1-2^{-i})$ gives
\begin{multline*}
	F(3^{n-1})
	\\
	\leq16^{-M}F\bigl(3^{n-1}+3^{n-2}(1-2^{-M})\bigr)
	+C3^{n(d+3)}\mathsf B3^{-3(n-2)}
	\sum_{i=0}^{M-1}16^{-i}8^{i+1}
	\\
	\leq16^{-M}F\bigl(3^{n-1}+3^{n-2}(1-2^{-M})\bigr)+C 3^{nd}\mathsf B\,.
\end{multline*}
Let~$M\to\infty$. Every set~$(q+\cu_{n-4})\cap\cu_m$ in the lemma lies in~$\cu_m\cap\overline{x+2\cdot3^{n-1}\cu_0}$ and has volume at least~$c(d)3^{nd}$. Dividing proves~\eqref{e.centered.energy.good.scale}. For the cutoff modification, the same argument uses Proposition~\ref{p.cutoff.regularity.good.scales} for the anchor error bound. All subsequent steps concern~$\a_L$ and~$\a_0$ deterministically and do not require~$L\geq m$.
\end{proof}

The next lemma adapts~\citet[Lemma~4.10]{ASD}.
\begin{lemma}[{Harmonic approximation on good scales}]
\label{l.harmonic.approximation.good.scales.GMC}
There exists~$C(d)<\infty$ such that, for every~$s\in [512\gmcdelta^2,\nf 14]$, every~$L,m,n\in\N_0$ with~$L\geq m\geq n+5$, every~$z \in \cu_m$ and~$x \in (z+\cu_{n-3})\cap\cu_m$, and every~$u\in H^1(\cu_m)$,~$\g \in H^{s}(\cu_m;\Rd)$, and~$h \in H^{1+s}(\cu_m)$ satisfying
\begin{equation*}
	\left\{
	\begin{aligned}
		& -\nabla \cdot \a_L \nabla u = \nabla \cdot \g  & \text{in} & \ \cu_m\,, \\
		& u = h & \text{on} & \ \partial \cu_m\,,
	\end{aligned}
	\right.
\end{equation*}
the following holds. Let~$y \in \cu_m$ be such that~$(x+\cu_{n-4})\cap\cu_m \subset y+\cu_{n-2} \subset (x+\cu_{n-1})\cap\cu_m$. Let~$v$ be the unique solution of the Dirichlet problem
\begin{equation*}
	\left\{
	\begin{aligned}
		& -\Delta v = 0
		& \text{in} & \ y+\cu_{n-2} \,,
		\\
		& v = u & \text{on} & \ \partial ( y+\cu_{n-2} ) \,.
	\end{aligned}
	\right.
\end{equation*}
Then
\begin{align}
	\lefteqn{
	\| u - v \|_{\underline{L}^2((x+\cu_{n-4})\cap\cu_m)}
	\indc_{\mathcal{G}_{n+2,z}(1,\nf s8)}
	} \qquad &
	\notag \\ &
	\leq
	Cs^{-\nf32}
	\mathcal{E}_{\nf s8,\infty,2}(z+\cu_{n+2} ; \a_L , (\b_{L,n+2})_{z+\cu_{n+2}})
	\Bigl( \| u - (u)_{(x+\cu_{n})\cap\cu_m} \|_{\underline{L}^2((x+\cu_{n})\cap\cu_m)}
	\notag \\ & \qquad
	+ s^{-\nf32}  3^{n}\bigl| (\nabla h) _{(x+\cu_{n})\cap\cu_m} \bigr| \indc_{\{\overline{(x+\cu_{n})\cap\cu_m} \cap \partial \cu_m \neq \emptyset \}}
	\Bigr)
	\notag \\ & \qquad
	+
	Cs^{-\nf{15}2} (\b_{L,n+2})_{z+\cu_{n+2}}^{-1} 3^{(1+s)n} [ \g ]_{\underline{H}^s((x+\cu_{n})\cap\cu_m)}
	\notag \\ & \qquad
	+
	Cs^{-\nf72} 3^{(1+s) n}  [ \nabla h ]_{\underline{H}^{s}((x+\cu_{n})\cap\cu_m)} \indc_{\{\overline{(x+\cu_{n})\cap\cu_m} \cap \partial \cu_m \neq \emptyset \}}
	\,.
	\label{e.good.scale.homog.GMC}
\end{align}
\end{lemma}

\begin{proof}
Let~$\a_0=(\b_{L,n+2})_{z+\cu_{n+2}}$. The proof has three parts. First, the good event bounds the coarse-graining errors on~$y+\cu_{n-2}$. Second, the coarse-graining estimate of~\citet[Proposition~2.21]{ASD} and a duality argument bound the~$\underline L^2$ norm of~$u-v$ by these errors times the energy norm of~$u$ on~$y+\cu_{n-2}$, plus a forcing term, as in~\eqref{e.harmonic.approximation.coarse.graining}; most of the work is in tracking the powers of~$s$ in this part. Third, Lemma~\ref{l.centered.energy.good.scale}, applied on a partition of~$y+\cu_{n-2}$, bounds this energy norm by the oscillation of~$u$ and the data on~$(x+\cu_n)\cap\cu_m$.

\smallskip

Stability and Proposition~\ref{p.good.scale.mathcal.E} give, on~$\mathcal G_{n+2,z}(1,\nf s8)$,
\begin{equation}
	\mathcal E_{s/3,\infty,1}(y+\cu_{n-2};\a_L,\a_0)
	+\mathcal E_{s/6,\infty,2}(y+\cu_{n-2};\a_L,\a_0)
	\leq C\mathcal E_{s/8,\infty,2}(z+\cu_{n+2};\a_L,\a_0)
	\leq C\,.
	\label{e.harmonic.approximation.local.control}
\end{equation}
By Definition~\ref{d.mathcal.E}, for every~$t\in(0,1]$ and~$q\in[1,\infty]$,
\begin{equation}
	\mathcal E_{t,\infty,q}(y+\cu_{n-2},n-3;\a_L,\a_0)
	\leq3^t\mathcal E_{t,\infty,q}(y+\cu_{n-2};\a_L,\a_0)\,.
	\label{e.harmonic.approximation.local.ellipticity}
\end{equation}
This follows by dropping the nonnegative terms above~$n-3$ and shifting the normalization by one scale; the supremum endpoint has the same comparison.

\smallskip

Let~$v_{\g}$ solve~$-\a_0\Delta v_{\g}=\nabla\cdot\g$ in~$y+\cu_{n-2}$ with trace~$u$. We apply~\citet[Proposition~2.21, (2.119)]{ASD} directly, retaining its unrenormalized fractional test norm during the following duality calculation. Its parameters are~$p=2$ and~$(s_1,\vartheta,s_2)=(s/3,s/2,s)$; the sums stop at~$n-3$. Thus~\eqref{e.harmonic.approximation.local.control} and~\eqref{e.harmonic.approximation.local.ellipticity} also bound the truncated errors in~(2.119). For~$w=u-v_{\g}\in H^1_0(y+\cu_{n-2})$, solve~$-\Delta\psi=w$ in~$y+\cu_{n-2}$, with zero trace. Cube~$H^2$ regularity and rescaling give
\begin{equation*}
	\|\nabla\psi\|_{\underline L^2(y+\cu_{n-2})}
	\leq C3^{n-2}\|w\|_{\underline L^2(y+\cu_{n-2})},\qquad
	\|\nabla^2\psi\|_{\underline L^2(y+\cu_{n-2})}
	\leq C\|w\|_{\underline L^2(y+\cu_{n-2})}\,.
\end{equation*}
For~$0<\vartheta\leq1/8$, the unrenormalized test-function bound is
\begin{multline*}
	\left(\fint_{y+\cu_{n-2}}\int_{y+\cu_{n-2}}
	\frac{|\nabla\psi(p)-\nabla\psi(q)|_2^2}
	{|p-q|_2^{d+2\vartheta}}\,dp\,dq\right)^{1/2}
	+3^{-(n-2)\vartheta}\|\nabla\psi\|_{\underline L^2(y+\cu_{n-2})}
	\\
	\leq C3^{(n-2)(1-\vartheta)}\|w\|_{\underline L^2(y+\cu_{n-2})}\,.
\end{multline*}
To see uniformity, extend~$\nabla\psi$ in~$H^1$ from the unit cube. For displacements of size at most a fixed dimensional constant, its squared translation difference is bounded by the squared displacement times the~$H^1$ norm. The radial integral is~$C/(2-2\vartheta)$, bounded in the stated range. Smooth approximation permits this vector field as a dual test. The identity~$\fint_{y+\cu_{n-2}}w^2=\fint_{y+\cu_{n-2}}\nabla w\cdot\nabla\psi$ therefore bounds~$\|w\|_{\underline L^2(y+\cu_{n-2})}$ by~$C3^{(n-2)(1-\vartheta)}$ times the unrenormalized negative norm of~$\nabla w$, with no inverse power of~$\vartheta$.

\smallskip

In~(2.119), the~$p=2$ local energy sum is at most
\begin{equation*}
	\left(\sum_{k=-\infty}^{n-3}
	3^{-2(s/2-s/3)(n-3-k)}\right)^{1/2}
	\|\a_L^{1/2}\nabla u\|_{\underline L^2(y+\cu_{n-2})}
	\leq Cs^{-1/2}\|\a_L^{1/2}\nabla u\|_{\underline L^2(y+\cu_{n-2})}\,.
\end{equation*}
At each scale, the volume-weighted average of the squared local energies is the squared energy norm on~$y+\cu_{n-2}$, since the subcubes partition~$y+\cu_{n-2}$. Their number therefore contributes no scale loss, and the energy coefficient~$s^{-1}$ becomes~$s^{-3/2}$. The forcing coefficient is~$s^{-9/2}(s-s/2)^{-1}$, and converting the unrenormalized~$H^s$ seminorm to our seminorm contributes~$s^{-1/2}$, giving~$Cs^{-6}$. The factor~$1+\mathcal E_{s/6,\infty,2}(y+\cu_{n-2})^2$ is bounded on the event. The zero-trace difference~$v_{\g}-v$ satisfies the ordinary energy and Poincar\'e estimates, giving~$\|v_{\g}-v\|_{\underline L^2(y+\cu_{n-2})} \leq Cs^{-1/2}\a_0^{-1}3^{(1+s)n}[\g]_{\underline H^s(y+\cu_{n-2})}$. Combining these bounds and passing to~$(x+\cu_{n-4})\cap\cu_m\subset y+\cu_{n-2}$ by the bounded volume ratio yields
\begin{align}
	\|u-v\|_{\underline L^2((x+\cu_{n-4})\cap\cu_m)}
	\indc_{\mathcal G_{n+2,z}(1,s/8)}
	&\leq Cs^{-3/2}3^n
	\mathcal E_{s/8,\infty,2}(z+\cu_{n+2};\a_L,\a_0)
	\a_0^{-1/2}\|\a_L^{1/2}\nabla u\|_{\underline L^2(y+\cu_{n-2})}
	\notag \\ & \qquad
	+Cs^{-6}\a_0^{-1}3^{(1+s)n}
	[\g]_{\underline H^s((x+\cu_n)\cap\cu_m)}
	\,.
	\label{e.harmonic.approximation.coarse.graining}
\end{align}

\smallskip

Partition~$y+\cu_{n-2}$ into its~$9^d$ grandchildren and apply Lemma~\ref{l.centered.energy.good.scale} at their centers. Taking square roots gives
\begin{align}
	\a_0^{-1/2}\|\a_L^{1/2}\nabla u\|_{\underline L^2(y+\cu_{n-2})}
	\indc_{\mathcal G_{n+2,z}(1,s/8)}
	&\leq C3^{-n}\|u-(u)_{(x+\cu_n)\cap\cu_m}\|_{\underline L^2((x+\cu_n)\cap\cu_m)}
	\notag \\ & \qquad
	+C|(\nabla h)_{(x+\cu_n)\cap\cu_m}|_2
	\indc_{\{\overline{(x+\cu_n)\cap\cu_m}\cap\partial\cu_m\neq\varnothing\}}
	\notag \\ & \qquad
	+Cs^{-6}\a_0^{-1}3^{sn}[\g]_{\underline H^s((x+\cu_n)\cap\cu_m)}
	\notag \\ & \qquad
	+Cs^{-2}3^{sn}[\nabla h]_{\underline H^s((x+\cu_n)\cap\cu_m)}
	\indc_{\{\overline{(x+\cu_n)\cap\cu_m}\cap\partial\cu_m\neq\varnothing\}}
	\,.
	\label{e.harmonic.approximation.energy}
\end{align}
Insert~\eqref{e.harmonic.approximation.energy} in~\eqref{e.harmonic.approximation.coarse.graining}. The forcing and trace powers are respectively~$s^{-3/2}s^{-6}=s^{-15/2}$ and~$s^{-3/2}s^{-2}=s^{-7/2}$. Their remaining error factors are bounded by~$C$ on the event; the direct forcing term~$s^{-6}$ is bounded by~$s^{-15/2}$. The mean slope term~$Cs^{-3/2}3^n\mathcal E\,|(\nabla h)_{(x+\cu_n)\cap\cu_m}|_2$ is bounded by the mean slope term in~\eqref{e.good.scale.homog.GMC}, which allows an additional~$s^{-3/2}$ inside the parentheses. This proves the unchanged conclusion for every admissible~$y$.
\end{proof}

The next lemma adapts~\citet[Lemma~4.11]{ASD}.

\begin{lemma}[{Excess decay on good scales}]
\label{l.excess.decay.good.scales.GMC}
There exists a constant~$C(d)<\infty$ such that, for every~$s\in [512\gmcdelta^2,\nf14]$ and~$\ep \in [8s^{-1}\gmcdelta^2, 1]$, every~$k\in\N$ and~$L,m,n\in\N_0$ with~$k\leq n$ and~$L\geq m \geq n+5$, every~$x,z \in \cu_m$ such that~$x \in (z+\cu_{n-3})\cap\cu_m$, and every~$u\in H^1(\cu_m)$,~$\g \in H^{s}(\cu_m;\Rd)$, and~$h \in C^{1,\nf12}(\cu_m)$ satisfying
\begin{equation*}
	\left\{
	\begin{aligned}
		& -\nabla \cdot \a_L \nabla u = \nabla \cdot \g  & \text{in} & \ \cu_m\,, \\
		& u = h & \text{on} & \ \partial \cu_m\,,
	\end{aligned}
	\right.
\end{equation*}
we have the estimate
\begin{align}
	\lefteqn{
	E(u,(x+\cu_{n-k})\cap\cu_m)
	\indc_{\mathcal{G}_{n+2,z}(\ep,\nf s8)}
	} \qquad &
	\notag \\ &
	\leq
	C \bigl( 3^{-\frac{1}{2}k} +
	3^{(1+\frac{d}{2}) k}
	s^{-\nf32}\ep \bigr)
	E(u,(x+\cu_{n})\cap\cu_m)
	\notag \\ & \qquad
	+
	C3^{(1+\frac{d}{2}) k}
	s^{-\nf32}
	\mathcal{E}_{\nf s8,\infty,2}(z+\cu_{n+2} ; \a_L , (\b_{L,n+2})_{z+\cu_{n+2}})
	\Bigl( \bigl| \nabla \ell (u,(x+\cu_{n})\cap\cu_m) \bigr|
	\notag \\ & \qquad
	+ s^{-\nf32} \bigl| (\nabla h) _{(x+\cu_{n})\cap\cu_m} \bigr| \indc_{\{\overline{(x+\cu_{n})\cap\cu_m} \cap \partial \cu_m \neq \emptyset \}}
	\Bigr)
	\notag \\ & \qquad
	+
	Cs^{-\nf{15}2} 3^{(1+\frac{d}{2}) k} (\b_{L,n+2})_{z+\cu_{n+2}}^{-1} 3^{sn} [ \g ]_{\underline{H}^s((x+\cu_{n})\cap\cu_m)}
	\notag \\ & \qquad
	+
	Cs^{-3} 3^{(1+\frac{d}{2}) k} 3^{\frac{1}{2} n} [ \nabla h ]_{W^{\nf12,\infty}((x+\cu_{n})\cap\cu_m)}
	\indc_{\{\overline{(x+\cu_{n})\cap\cu_m} \cap \partial \cu_m \neq \emptyset\}}
	\,.
	\label{e.excess.decay.one.step}
\end{align}
\end{lemma}
\begin{proof}
We compare~$u$ with the harmonic function that has the same boundary values on a cube of side~$3^{n-2}$ near~$x$. Step~1 bounds the difference by Lemma~\ref{l.harmonic.approximation.good.scales.GMC}, Step~2 records the excess decay of the harmonic function, and Step~3 transfers this decay to~$u$. The heart of the matter is Step~3: on the good event the coarse-graining error is at most~$C\ep$, so that the error of the comparison contributes~$C3^{(1+\frac{d}{2})k}s^{-\nf32}\ep\,E(u,(x+\cu_{n})\cap\cu_m)$, plus terms involving the slope of~$u$ and the data.

\smallskip

Fix~$x,z\in\cu_m$ as in the statement. We may assume~$k\geq6$, since the case~$k<6$ follows by increasing~$C(d)$. Choose~$y\in\cu_m$ so that~$(x+\cu_{n-4})\cap\cu_m\subseteq y+\cu_{n-2}\subseteq (x+\cu_{n-1})\cap\cu_m$.

\smallskip

\emph{Step 1.} We approximate~$u$ by a harmonic function. Let~$v$ be the harmonic function in~$y+\cu_{n-2}$ with boundary data~$u$:
\begin{equation*}
	\left\{
	\begin{aligned}
		& -\Delta v = 0  & \text{in} & \ y+\cu_{n-2}\,,
		\\
		&
		v = u & \text{on} & \ \partial (y+\cu_{n-2})\,.
	\end{aligned}
	\right.
\end{equation*}

\smallskip

In particular, the centered energy estimate used to prove the comparison is Lemma~\ref{l.centered.energy.good.scale}. It is valid at tolerance~$1$ for the same arbitrary comparison cube. The factors from the comparison are exactly~$s^{-3/2}$,~$s^{-15/2}$ and~$s^{-7/2}$, as in~\eqref{e.good.scale.homog.GMC}; the smaller event here introduces no additional boundary restriction.

\smallskip

Since~$\ep \leq 1$, the event~$\mathcal{G}_{n+2,z}(\ep,\nf s8)$ implies~$\mathcal{G}_{n+2,z}(1,\nf s8)$. Lemma~\ref{l.harmonic.approximation.good.scales.GMC}, applied in~$y + \cu_{n-2}$ on the good event~$\mathcal{G}_{n+2,z}(1,\nf s8)$, therefore yields~\eqref{e.good.scale.homog.GMC} with~$\indc_{\mathcal{G}_{n+2,z}(\ep,\nf s8)}$ in place of~$\indc_{\mathcal{G}_{n+2,z}(1,\nf s8)}$ on the left side.

\smallskip

\emph{Step 2.} We use the excess decay of harmonic functions. The function~$v$ is harmonic in~$(x+\cu_{n-4})\cap\cu_m$. If~$\overline{(x+\cu_{n-4})\cap\cu_m} \cap \partial \cu_m \neq \emptyset$, then moreover~$v = u = h$ on~$\partial (y+\cu_{n-2}) \cap \partial \cu_m$, with~$h \in C^{1,\nf12}$. The Schauder estimates for harmonic functions, combined with odd reflection across the flat parts of the boundary, therefore give~$v \in C^{1,\nf12}((x+\cu_{n-5})\cap\cu_m)$ with the estimate, where~$[\,\cdot\,]_{C^{0,\nf12}(U)}$ denotes the H\"older seminorm of order~$\nf12$ on~$U$,
\begin{align*}
	[\nabla v]_{C^{0,\nf12}((x+\cu_{n-5})\cap\cu_m)}
	&\leq
	C 3^{-\frac{1}{2} n} E(v,(x+\cu_{n-4})\cap\cu_m)
	\notag \\ & \qquad
	+ C[\nabla h]_{C^{0,\nf12}((x+\cu_{n-4})\cap\cu_m)} \indc_{\{ \overline{(x+\cu_{n})\cap\cu_m} \cap \partial \cu_m \neq \emptyset\}}
	\,.
\end{align*}
It follows that, for every~$k \geq 6$,
\begin{align*}
	E(v,(x+\cu_{n-k})\cap\cu_m)
	&\leq
	C 3^{-\frac{1}{2} k} E(v,(x+\cu_{n-4})\cap\cu_m)
	\notag \\ & \qquad
	+ C 3^{\frac{1}{2}(n-k)} [\nabla h]_{C^{0,\nf12}((x+\cu_{n-2})\cap\cu_m)} \indc_{\{ \overline{(x+\cu_{n})\cap\cu_m} \cap \partial \cu_m \neq \emptyset\}}
	\,.
\end{align*}

\smallskip

\emph{Step 3.} We conclude. By the definition of the excess~\eqref{e.excess.def} and the triangle inequality, for~$k \geq 6$,
\begin{align*}
	E(u,(x+\cu_{n-k})\cap\cu_m)
	&
	\leq
	E(v,(x+\cu_{n-k})\cap\cu_m)
	+
	3^{-(n-k)}\| u - v \|_{\underline{L}^2((x+\cu_{n-k})\cap\cu_m)}
	\notag \\ &
	\leq
	C 3^{-\frac{1}{2} k} E(v,(x+\cu_{n-4})\cap\cu_m)
	+
	3^{-(n-k)}\| u - v \|_{\underline{L}^2((x+\cu_{n-k})\cap\cu_m)}
	\notag \\ & \qquad
	+
	C 3^{\frac{1}{2}(n-k)} [\nabla h]_{C^{0,\nf12}((x+\cu_{n})\cap\cu_m)} \indc_{\{ \overline{(x+\cu_{n})\cap\cu_m} \cap \partial \cu_m \neq \emptyset\}}
	\,.
\end{align*}
For the~$\| u - v \|$ term, since~$(x+\cu_{n-k})\cap\cu_m \subset (x+\cu_{n-4})\cap\cu_m$ and~$|(x+\cu_{n})\cap\cu_m|/|(x+\cu_{n-k})\cap\cu_m| \leq C 3^{dk}$, we have
\begin{equation*}
	3^{-(n-k)}\| u - v \|_{\underline{L}^2((x+\cu_{n-k})\cap\cu_m)} \leq C 3^{(\frac{d}{2} + 1)k} 3^{-n}\| u - v \|_{\underline{L}^2((x+\cu_{n-4})\cap\cu_m)}
	\,.
\end{equation*}
The excess of~$v$ also satisfies~$E(v,(x+\cu_{n-4})\cap\cu_m)\leq CE(u,(x+\cu_{n})\cap\cu_m) +C3^{-n}\|u-v\|_{\underline L^2((x+\cu_{n-4})\cap\cu_m)}$. Thus
\begin{align*}
	E(u,(x+\cu_{n-k})\cap\cu_m)
	&
	\leq
	C 3^{-\frac{1}{2} k} E(u,(x+\cu_{n})\cap\cu_m) + C 3^{(\frac{d}{2} + 1)k} 3^{-n}\| u - v \|_{\underline{L}^2((x+\cu_{n-4})\cap\cu_m)}
	\notag \\ & \qquad
	+
	C 3^{\frac{1}{2}(n-k)} [\nabla h]_{C^{0,\nf12}((x+\cu_{n})\cap\cu_m)}
	\indc_{\{ \overline{(x+\cu_{n})\cap\cu_m} \cap \partial \cu_m \neq \emptyset\}}
	\,.
\end{align*}
To estimate the middle term, we use the bound from Step~1 together with the inequality
\begin{equation*}
	3^{-n} \| u - (u)_{(x+\cu_{n})\cap\cu_m} \|_{\underline{L}^2((x+\cu_{n})\cap\cu_m)}
	\leq
	E(u,(x+\cu_{n})\cap\cu_m) + C \bigl| \nabla \ell (u,(x+\cu_{n})\cap\cu_m) \bigr|
	\,.
\end{equation*}
We multiply the resulting inequality by~$\indc_{\mathcal{G}_{n+2,z}(\ep,\nf s8)}$. By~\eqref{e.mathcalE.bound.for.good} and the definition of~$\mathcal{G}_{n+2,z}(\ep,\nf s8)$, we have~$\mathcal{E}_{\nf s8,\infty,2}(z+\cu_{n+2} ; \a_L, (\b_{L,n+2})_{z+\cu_{n+2}}) \indc_{\mathcal{G}_{n+2,z}(\ep,\nf s8)} \leq C\ep$, so the part of the bound on~$\| u - v \|$ involving~$E(u,(x+\cu_{n})\cap\cu_m)$ contributes~$C 3^{(1+\frac{d}{2})k}s^{-\nf32}\ep \, E(u,(x+\cu_{n})\cap\cu_m)$, the second term in the coefficient of~$E(u,(x+\cu_{n})\cap\cu_m)$ in~\eqref{e.excess.decay.one.step}. In the part involving~$|\nabla \ell(u,(x+\cu_{n})\cap\cu_m)|$, we keep the factor~$\mathcal{E}_{\nf s8,\infty,2}(z+\cu_{n+2} ; \a_L, (\b_{L,n+2})_{z+\cu_{n+2}})$ instead of replacing it by~$\ep$; this produces the second and third lines of~\eqref{e.excess.decay.one.step}. Finally, the averaged seminorm satisfies
\begin{equation*}
	3^{sn}[\nabla h]_{\underline H^s((x+\cu_{n})\cap\cu_m)}
	\leq Cs^{\nf12}3^{\frac{1}{2}n}
	[\nabla h]_{W^{\nf12,\infty}((x+\cu_{n})\cap\cu_m)}\,.
\end{equation*}
The forcing and boundary terms from Step~1 therefore give the last two lines of~\eqref{e.excess.decay.one.step}.
\end{proof}

The following lemma iterates the excess decay over a range of scales. For the sets~$U_j$ used in it, the scale-normalized excess and a best affine approximation are
\begin{equation*}
	E(u,U_j)\coloneqq3^{-j}\min_{\ell\ \mathrm{affine}}\|u-\ell\|_{\underline L^2(U_j)}\,,
	\qquad \ell(u,U_j)\in\mathop{\rm argmin}_{\ell\ \mathrm{affine}}\|u-\ell\|_{\underline L^2(U_j)}\,.
\end{equation*}
The cube contained in~$U_j$ makes the quadratic form on affine functions nondegenerate, so the minimizer exists and is unique. The minimizing constant satisfies~$(\ell(u,U_j))_{U_j}=(u)_{U_j}$.

\begin{lemma}
\label{l.iteration.lemma.GMC}
Let~$k \in \N$ and~$\theta \in (0,1)$ with~$\theta^k \in (0, \nf 35)$. Let~$n,m \in \Z$ with~$n < m$. Let~$\{ U_j \}_{j \in \Z}$ be bounded measurable sets such that, for every~$j\leq m$, we have~$U_{j-1} \subset U_j$ and~$U_j$ contains a translate of~$\cu_{j-2}$ and is contained in a translate of~$\cu_j$. Let~$u \in L^2(U_m)$. Assume that~$\mathcal{B} \subseteq \Z \cap [n,m]$ and that~$\{\ep_j\}_{j \in \Z \cap [n,m]}$ and~$\{\delta_j\}_{j \in \Z \cap [n,m]}$ are sequences of nonnegative numbers such that
\begin{equation}
	E(u,U_{j-k})
	\leq \theta^k E(u,U_{j}) + \ep_j |\nabla \ell(u,U_j)| + \delta_j\,, \qquad \forall j \in \Z \cap [n,m] \setminus \mathcal{B}\,.
	\label{e.Ej.decay.assumption}
\end{equation}
Then there exists a constant~$C(d)<\infty$ such that, with~$A\coloneqq C(k+1)(|\mathcal B|+1)+C\sum_{j=n}^m\ep_j$, we have
\begin{equation}
	3^{-n} \| u - (u)_{U_n} \|_{\underline{L}^2(U_n)}
	\leq
	e^A \biggl(3^{-m} \| u - (u)_{U_m} \|_{\underline{L}^2(U_{m})}
	+ \sum_{j=n}^m \delta_j\biggr)\,.
	\label{e.iteration.slope.lemma}
\end{equation}
Moreover, for every~$M\geq0$ such that~$\ep_j\leq M$ for every~$j\in[n,m]\cap\Z$,
\begin{equation}
	E(u,U_n)
	\leq
	\theta^{-C(k+1)(|\mathcal{B}|+1)}e^A
	\biggl(
	\theta^{m-n} E(u,U_m)
	+M3^{-m} \| u - (u)_{U_m} \|_{\underline{L}^2(U_{m})}
	+(1+M)\sum_{j=n}^m \delta_j
	\biggr)\,.
	\label{e.excess.decay.lemma}
\end{equation}
\end{lemma}

\begin{proof}
The inclusions assumed for~$U_{j-1}$ and~$U_j$, together with the definition of~$\ell(u,U_j)$, give, for every~$j\leq m$,
\begin{equation}
	|\nabla\ell(u,U_j)-\nabla\ell(u,U_{j-1})|
	\leq C(E(u,U_j)+E(u,U_{j-1}))\,.
	\label{e.iteration.gradient.stability.GMC}
\end{equation}
Let~$[a,b]\cap\Z$ be a maximal integer interval contained in~$([n,m]\cap\Z)\setminus\mathcal B$, and partition it into residue classes modulo~$k$. Summing~\eqref{e.Ej.decay.assumption} in each class and using~$\theta^k<3/5$ gives, for~$i\in[a,b]\cap\Z$,
\begin{equation}
	\sum_{j=i}^{b}E(u,U_j)
	\leq C3^{(3+\frac{3d}{2})k}E(u,U_b)
	+C\sum_{j=i+k}^{b}\ep_j|\nabla\ell(u,U_j)|
	+C\sum_{j=i+k}^{b}\delta_j\,.
	\label{e.iteration.excess.sum.GMC}
\end{equation}
To obtain~\eqref{e.iteration.excess.sum.GMC}, apply~\eqref{e.Ej.decay.assumption} at~$j+k$ for~$i\leq j\leq b-k$. The remaining at most~$k$ terms satisfy~$E(u,U_{b-r})\leq C3^{(1+d/2)r}E(u,U_b)$ by nesting and the cube sandwich. Since~$1-\theta^k>2/5$, the excess sum on the right can be absorbed. The sums are empty when~$i+k>b$.

\smallskip

For a fixed~$i\in[a,b]$, telescoping~\eqref{e.iteration.gradient.stability.GMC} and using the shifted sum at every~$r\in[i,b]$ give
\begin{equation*}
	E(u,U_r)+|\nabla\ell(u,U_r)|\leq B_i+C\sum_{j=r+1}^b\ep_j\bigl(E(u,U_j)+|\nabla\ell(u,U_j)|\bigr)\,,
\end{equation*}
where
\begin{equation*}
	B_i=C3^{Ck}\left(E(u,U_b)+|\nabla\ell(u,U_b)|+\sum_{j=i}^b\delta_j\right)\,.
\end{equation*}
There is no~$\ep_r\bigl(E(u,U_r)+|\nabla\ell(u,U_r)|\bigr)$ term, because~$k\geq1$. Backward induction therefore gives the product form of the discrete Gronwall inequality,
\begin{equation*}
	E(u,U_i)+|\nabla\ell(u,U_i)|\leq B_i\prod_{j=i+1}^b(1+C\ep_j)
	\leq B_i\exp\left(C\sum_{j=i+1}^b\ep_j\right)\,.
\end{equation*}
Only nonnegativity of~$\ep_j$ is used. In particular,
\begin{equation}
	E(u,U_i)+|\nabla\ell(u,U_i)|
	\leq
	\exp\biggl(C(k+1)+C\sum_{j=i}^{b}\ep_j\biggr)
	\biggl(E(u,U_b)+|\nabla\ell(u,U_b)|+\sum_{j=i}^{b}\delta_j\biggr)\,.
	\label{e.iteration.interval.GMC}
\end{equation}

\smallskip

For every~$j\leq m$, inequality~\eqref{e.iteration.gradient.stability.GMC} and the definition of the excess imply~$E(u,U_{j-1})+|\nabla\ell(u,U_{j-1})|\leq C(E(u,U_j)+|\nabla\ell(u,U_j)|)$. There are at most~$|\mathcal B|+1$ maximal integer intervals contained in~$([n,m]\cap\Z)\setminus\mathcal B$, and each interval has at most~$k$ incomplete residue classes. Iterating~\eqref{e.iteration.interval.GMC} over these intervals and across the indices in~$\mathcal B$ yields
\begin{equation}
	E(u,U_n)+|\nabla\ell(u,U_n)|
	\leq e^A
	\biggl(E(u,U_m)+|\nabla\ell(u,U_m)|+\sum_{j=n}^{m}\delta_j\biggr)\,.
	\label{e.iteration.combined.GMC}
\end{equation}
The definition of the best affine approximation gives
\begin{equation*}
	E(u,U_m)+|\nabla\ell(u,U_m)|
	\leq C3^{-m}\|u-(u)_{U_m}\|_{\underline L^2(U_m)}
	\,,
\end{equation*}
and, conversely,
\begin{equation*}
	3^{-n}\|u-(u)_{U_n}\|_{\underline L^2(U_n)}
	\leq C(E(u,U_n)+|\nabla\ell(u,U_n)|)\,.
\end{equation*}
Together with~\eqref{e.iteration.combined.GMC}, these inequalities prove~\eqref{e.iteration.slope.lemma}.

\smallskip

To prove~\eqref{e.excess.decay.lemma}, we again sum~\eqref{e.Ej.decay.assumption} separately over the residue classes modulo~$k$, but now keep track of the powers of~$\theta$. For every maximal integer interval~$[a,b]\cap\Z$ contained in~$([n,m]\cap\Z)\setminus\mathcal B$ and every~$i\in[a,b]\cap\Z$,
\begin{equation*}
	E(u,U_i)
	\leq
	e^{C(k+1)}\theta^{-C(k+1)}
	\biggl(
	\theta^{b-i}E(u,U_b)
	+M\max_{i\leq j\leq b}|\nabla\ell(u,U_j)|
	+\sum_{j=i}^{b}\delta_j
	\biggr)\,.
\end{equation*}
The gradient term contributes at most
\begin{equation*}
	M\max_{i\leq j\leq b}|\nabla\ell(u,U_j)|\sum_{r\geq0}\theta^{kr}
	\leq\frac52M\max_{i\leq j\leq b}|\nabla\ell(u,U_j)|\,.
\end{equation*}
At a bad index, quasi-monotonicity gives~$E(u,U_{j-1})\leq C\theta^{-1}(\theta E(u,U_j))$; each incomplete residue class loses at most~$\theta^{-Ck}e^{Ck}$. There are at most~$|\mathcal B|+1$ good intervals. Composing these inequalities therefore introduces a factor of at most~$\theta^{-C(k+1)(|\mathcal B|+1)}e^{C(k+1)(|\mathcal B|+1)}$. The first conclusion, applied on every subinterval, bounds all the gradients by~$e^A(3^{-m}\|u-(u)_{U_m}\|_{\underline L^2(U_m)}+\sum_{j=n}^m\delta_j)$, with its exponent no larger than~$A$. Increasing the dimensional constant in~$A$ absorbs the preceding interval losses and gives~\eqref{e.excess.decay.lemma}.
\end{proof}

\subsection{The large-scale H\"older estimate}
\label{ss.regularity.iteration}

In this subsection we prove Proposition~\ref{p.Holder.regularity}. Lemma~\ref{l.iteration.lemma.GMC} iterates the excess decay estimate~\eqref{e.excess.decay.one.step} of Lemma~\ref{l.excess.decay.good.scales.GMC} from scale~$m$ down to scale~$n$, up to a factor controlled by the number of bad scales between~$n$ and~$m$ and by the sum of the coarse-graining errors on the good scales. The minimal scale~$\mathcal L_L(\gamma,m)$ is built from the random variables~$M$ and~$N$ of Lemmas~\ref{l.sum.the.errors} and~\ref{l.min.scale.good.scale}. For every~$n\leq m-\mathcal L_L(\gamma,m)$, both the number of bad scales and the sum of the coarse-graining errors are at most~$C$ plus a small multiple of~$(1-\gamma)(m-n)$, and hence the factor is at most~$C3^{\frac14(1-\gamma)(m-n)}$.

\subsubsection{The minimal scale}
\label{ss.minimal.scale.sep}

The iteration requires bounds for
\begin{equation}
	\sum_{k=n}^m
	\mathcal{E}_{\nf1{16},\infty,2}(\cu_k ; \a_L,(\b_{L,k})_{\cu_k}) \indc_{\mathcal{G}_k(\ep,\nf1{16}) }
	\label{e.sum.of.eps.k}
\end{equation}
and for the number of bad scales between~$n$ and~$m$,~$\sum_{k=n}^m\bigl(1-\indc_{\mathcal{G}_k(\ep,\nf1{16})}\bigr)$. Propositions~\ref{p.good.scale.mathcal.E} and~\ref{p.density.of.good.scales} control~\eqref{e.sum.of.eps.k} and this number, respectively.

\smallskip

For~$k\in\N_0$,~$z\in\Rd$, and~$s\in(0,\nf12]$, define
\begin{align}
	\mathbf{G}_{k,z}(s)
	&
	\coloneqq
	\sup_{0\leq j \leq k} \sup_{0\leq l \leq j-2} 3^{-\frac{1}{2}s (k-l)}
	\max_{z' \in z+\bigl(3^l \Zd \cap (\cu_j \setminus \cu_{j-1})\bigr)}
	\min\!\left\{
	\max_{|e| = 1}J(z'+\cu_l,e,e\,; \ahom_l^{-1}\a_l),1
	\right\}^{\nf12}
	\notag \\ & \qquad
	+
	\sup_{0\leq j \leq k}
	3^{-\frac{1}{8}s(k-j)} \| \g_{k,j} \|_{L^\infty(z+\cu_{k})}
	+
	3^{-\frac{1}{8}sk}\|\g_0\|_{L^\infty(z+\cu_k)}
	+
	\sum_{j=k}^\infty
	3^{k}
	\| \nabla \g_j \|_{L^\infty(z+\cu_k)}
	\,.
	\label{e.def.bfG.mhq}
\end{align}
Here and below, a supremum over an empty set of integer indices is zero.

\begin{lemma}
\label{l.sum.the.errors}
Let~$s \in (0,\nf 12]$. There exists~$C(d)<\infty$ such that, if~$s \geq C \gmcdelta^2 \left| \log \gmcdelta \right|$, then, for every~$m,n \in\N_0$ with~$n\leq m$ and~$z \in \Rd$,
\begin{equation}
	\avsum_{k=n}^{m}
	\mathbf{G}_{k,z}(s)
	\leq
	Cs^{-\nf 72}   \gmcdelta + \O_{\Gamma_2}\bigl( C s^{-\nf 72}  \gmcdelta \left| \log \gmcdelta \right| ^{\nf 12}  (m-n+1)^{-\nf12} \bigr)
	\,.
	\label{e.after.good.scales.uniform}
\end{equation}
Moreover, for every~$\lambda > 0$ satisfying~$\lambda \geq Cs^{-\nf 72}  \gmcdelta \left| \log \gmcdelta \right|^{\nf 12} $ and every~$m\in\N_0$, there exists a random variable~$M(\lambda,s,m)\in\{-1,0,\ldots,m\}$ such that
\begin{equation}
	0\leq m-M(\lambda,s,m)
	\leq
	1+\O_{\Gamma_1}\bigl(
	Cs^{-7}\gmcdelta^2|\log\gmcdelta|\lambda^{-2}
	\bigr)
	\label{e.M.lambda.m.bound}
\end{equation}
and, for every~$n \in \N_0$ satisfying~$n \leq M(\lambda,s,m)$, we have
\begin{equation}
	\max_{z \in 3^n \Zd \cap \cu_m}
	\sum_{j=n}^{m}
	\mathbf{G}_{j,z}(s)
	\leq
	\lambda (m-n) \,.
	\label{e.after.good.scales}
\end{equation}
\end{lemma}
\begin{proof}
The proof has four steps. Step~1 bounds the sum over the scales of the terms of~$\mathbf G_{k,z}(s)$ involving the fields~$\g_j$. Step~2 bounds the sum of the term involving~$J$ by a sum over the annuli~$\cu_j\setminus\cu_{j-1}$ and estimates the moments of each summand. Step~3 sums these estimates and proves~\eqref{e.after.good.scales.uniform}. In Steps~1 and~3, the main point is the independence of the summands at different scales (in Step~3, within each residue class of scales): it makes the random part of each sum of order~$(m-n+1)^{\nf12}$ rather than~$m-n+1$. Step~4 defines~$M(\lambda,s,m)$ and deduces~\eqref{e.M.lambda.m.bound} from~\eqref{e.after.good.scales.uniform} by a union bound over centers and scales.

\smallskip

\emph{Step 1.} We estimate the terms of~$\mathbf G_{k,z}(s)$ involving the fields~$\g_j$. By stationarity it suffices to prove~\eqref{e.after.good.scales.uniform} for~$z=0$. The triangle inequality gives
\begin{align*}
	&
	\sum_{k=n}^m
	\left(
	\sup_{0\leq j \leq k}
	3^{-\frac{1}{8}s(k-j)} \| \g_{k,j} \|_{L^\infty(\cu_{k})}
	+3^{-\frac{1}{8}sk}\|\g_0\|_{L^\infty(\cu_k)}
	\right)
	\notag \\ & \qquad
	\leq
	2\sum_{k=n}^m
	\sup_{0\leq j \leq k}
	3^{-\frac{1}{8} s(k-j)}
	\sum_{i=j}^k
	\|  \mathbf{g}_i  \|_{L^{\infty}(\cu_k)}
	\notag \\ & \qquad
	\leq
	2\sum_{k=n}^m
	\sum_{i=0}^k
	\sup_{0\leq j \leq i}
	3^{-\frac{1}{8} s(k-j)}
	\|  \mathbf{g}_i  \|_{L^{\infty}(\cu_k)}
	\notag \\ & \qquad
	=
	2\sum_{k=n}^m
	\sum_{i=0}^k
	3^{-\frac{1}{8} s(k-i)}
	\|  \mathbf{g}_i  \|_{L^{\infty}(\cu_k)}
	\notag \\ & \qquad
	=
	2\biggl(
	\sum_{i=n}^m
	\sum_{k=i}^m
	+
	\sum_{k=n}^m \sum_{i=0}^{n-1}  \biggr)
	3^{-\frac{1}{8} s(k-i)}
	\|  \mathbf{g}_i  \|_{L^{\infty}(\cu_k)}
	\,.
\end{align*}
By~\eqref{e.same.law.shift} and~\eqref{e.coefficient.field.regularity}, together with~\eqref{e.max.scaling.O.Gamma.Two}, we obtain~$\|  \mathbf{g}_i  \|_{L^{\infty}(\cu_k)} \leq \O_{\Gamma_2}\bigl( C \gmcdelta (k-i+1)^{\nf 12} \bigr)$ and therefore, by independence, also
\begin{multline*}
	\sum_{k=n}^m
	\left(
	\sup_{0\leq j \leq k}
	3^{-\frac{1}{8}s(k-j)} \| \g_{k,j} \|_{L^\infty(\cu_{k})}
	+3^{-\frac{1}{8}sk}\|\g_0\|_{L^\infty(\cu_k)}
	\right)
	\\
	\leq
	C s^{-\nf 32} \gmcdelta(m-n+1)
	+
	\O_{\Gamma_2}\bigl( C s^{-\nf 32}  \gmcdelta (m-n+1)^{\nf 12} \bigr)
	\,.
\end{multline*}
For the gradient term,
\begin{align*}
	\sum_{k=n}^m
	\sum_{l=k}^\infty
	3^{k}
	\|  \nabla \g_l  \|_{L^{\infty}(\cu_k)}
	&
	=
	\sum_{l=n}^\infty
	\sum_{k=n}^{m\wedge l}
	3^{k}
	\|  \nabla \g_l  \|_{L^{\infty}(\cu_k)}
	\notag \\ &
	\leq
	\frac 32
	\sum_{l=n}^\infty
	3^{m\wedge l }
	\|  \nabla \g_l  \|_{L^{\infty}(\cu_l)}
	\notag \\ &
	=
	\frac 32
	\sum_{l=n}^{m-1}
	3^{l}
	\|  \nabla \g_l  \|_{L^{\infty}(\cu_l)}
	+
	\frac 32
	\sum_{l=m}^\infty
	3^{m}
	\|  \nabla \g_l  \|_{L^{\infty}(\cu_l)}
	\,.
\end{align*}
Again by~\eqref{e.same.law.shift},~\eqref{e.coefficient.field.regularity}, and independence,
\begin{equation*}
	\sum_{k=n}^m
	\sum_{l=k}^\infty
	3^{k}
	\|  \nabla \g_l  \|_{L^{\infty}(\cu_k)}
	\leq
	C \gmcdelta(m-n+1) + \O_{\Gamma_2}\bigl(C \gmcdelta (m-n+1)^{\nf12}\bigr)
	\,.
\end{equation*}
Thus
\begin{multline}
	\sum_{k=n}^m
	\sum_{l=k}^\infty
	3^{k}
	\|  \nabla \g_l  \|_{L^{\infty}(\cu_k)}
	+
	\sum_{k=n}^m
	\left(
	\sup_{0\leq j \leq k}
	3^{-\frac{1}{8}s(k-j)} \| \g_{k,j} \|_{L^\infty(\cu_{k})}
	+3^{-\frac{1}{8}sk}\|\g_0\|_{L^\infty(\cu_k)}
	\right)
	\\
	\leq
	C s^{-\nf 32} \gmcdelta(m-n+1) +
	\O_{\Gamma_2}\bigl( C s^{-\nf 32}  \gmcdelta (m-n+1)^{\nf 12} \bigr)
	\,.
	\label{e.sum.the.errors.fields}
\end{multline}

\smallskip

\emph{Step 2.} We estimate the term of~$\mathbf G_{k,z}(s)$ involving~$J$. Bounding the supremum over~$j$ by a sum, interchanging the sums, and summing the geometric series in~$k$, we obtain
\begin{align}
	\lefteqn{
	\sum_{k=n}^m
	\sup_{0\leq j \leq k} \sup_{0\leq l \leq j-2} 3^{-\frac{1}{2}s (k-l)}
	\max_{z \in 3^l \Zd \cap (\cu_j \setminus \cu_{j-1})}
	\min\!\left\{\max_{|e|=1}J(z+\cu_l,e,e\,;\ahom_l^{-1}\a_l),1\right\}^{\nf12}
	} \quad &
	\notag \\ &
	\leq
	\sum_{j=0}^{m}
	\sum_{k=n \vee j}^m
	\sup_{0\leq l \leq j-2} 3^{-\frac{1}{4}s (k-l)}
	\max_{z \in 3^l \Zd \cap (\cu_j \setminus \cu_{j-1})}
	\min\Bigl\{  \max_{|e| = 1} J(z+\cu_l,e,e\,; \ahom_{l}^{-1} \a_{l})\, ,  1 \Bigr\}^{\nf 12}
	\notag \\ &
	\leq
	4s^{-1}
	\sum_{j=n}^{m}
	\sup_{0\leq l \leq j-2} 3^{-\frac{1}{4}s (j-l)}
	\max_{z \in 3^l \Zd \cap (\cu_j \setminus \cu_{j-1})}
	\min\Bigl\{  \max_{|e| = 1} J(z+\cu_l,e,e\,; \ahom_{l}^{-1} \a_{l})\, ,  1 \Bigr\}^{\nf 12}
	\notag \\ & \qquad
	+
	4s^{-1}
	\sum_{j=0}^{n-1}
	\sup_{0\leq l \leq j-2} 3^{-\frac{1}{4}s (n-l)}
	\max_{z \in 3^l \Zd \cap (\cu_j \setminus \cu_{j-1})}
	\min\Bigl\{  \max_{|e| = 1} J(z+\cu_l,e,e\,; \ahom_{l}^{-1} \a_{l})\, ,  1 \Bigr\}^{\nf 12}
	\,.
	\label{e.split.J.for.minimal.scale}
\end{align}
By Proposition~\ref{p.coarse.grained.bound}, for~$\xi  \in [1, C^{-1} \gmcdelta^{-2} \left| \log \gmcdelta\right|^{-1}]$,
\begin{equation*}
	\E \Bigl[  \max_{|e| = 1} \bigl( J(\cu_l,e,e\,; \ahom_{l}^{-1} \a_{l})
	\bigr) ^{\nf \xi2}
	\Bigr]^{\nf1\xi}
	\leq
	C\gmcdelta  \log(2+\xi)^{\nf 12} \xi^{\nf 12}
	\,.
\end{equation*}
Hence, for every~$\xi \in [1,\infty)$,
\begin{equation*}
	\E\biggl[
	\min\Bigl\{  \max_{|e| = 1} J(z+\cu_l,e,e\,; \ahom_{l}^{-1} \a_{l})\, ,  1 \Bigr\}^{\nf \xi2}  \biggr]^{\nf 1\xi}
	\leq
	C\gmcdelta  \left| \log \gmcdelta \right|^{\nf 12} \xi^{\nf 12} \,,
\end{equation*}
which yields
\begin{equation*}
	\min\Bigl\{  \max_{|e| = 1} J(z+\cu_l,e,e\,; \ahom_{l}^{-1} \a_{l})\, ,  1 \Bigr\}^{\nf 12}
	\leq \O_{\Gamma_2} \bigl( C\gmcdelta  \left| \log \gmcdelta \right|^{\nf 12}  \bigr)
	\,.
\end{equation*}
Moreover, by H\"older's inequality, the condition~$s \geq C \gmcdelta^{2} \left| \log \gmcdelta\right|$, a union bound, and stationarity,
\begin{align*}
	\lefteqn{
	\E\biggl[ \sup_{0\leq l \leq j-2} 3^{-\frac{1}{4}s (j-l)}  \max_{z \in 3^l \Zd \cap (\cu_j \setminus \cu_{j-1})}
	\min\Bigl\{  \max_{|e| = 1} J(z+\cu_l,e,e\,; \ahom_{l}^{-1} \a_{l})\, ,  1 \Bigr\}^{\nf 12} \biggr]
	} \qquad &
	\notag \\ &
	\leq
	\E\biggl[ \sup_{0\leq l \leq j-2} 3^{-2d (j-l)}  \max_{z \in 3^l \Zd \cap (\cu_j \setminus \cu_{j-1})}
	\min\Bigl\{  \max_{|e| = 1} J(z+\cu_l,e,e\,; \ahom_{l}^{-1} \a_{l})\, ,  1 \Bigr\}^{\nf {4d}{s} } \biggr]^{\nf {s}{8d}}
	\notag \\ &
	\leq
	\biggl( \sum_{l=0}^{j-2} 3^{-d (j-l)}
	\E\Bigl[  \max_{|e| = 1} J(\cu_l,e,e\,; \ahom_{l}^{-1} \a_{l})^{\nf {4d}{s} } \Bigr]\biggr)^{\! \nf {s}{8d}}
	\notag \\ &
	\leq
	C s^{-\nf 12} \left| \log s \right|^{\nf 12}
	\gmcdelta
	\,,
\end{align*}
and, using~\eqref{e.max.scaling.O.Gamma.Two},
\begin{equation*}
	\sup_{0\leq l \leq j-2} 3^{-\frac{1}{4}s (j-l)}
	\max_{z \in 3^l \Zd \cap (\cu_j \setminus \cu_{j-1})}
	\min\Bigl\{  \max_{|e| = 1} J(z+\cu_l,e,e\,; \ahom_{l}^{-1} \a_{l})\, ,  1 \Bigr\}^{\nf 12}
	\leq
	\O_{\Gamma_2} \bigl( Cs^{-\nf 32} \gmcdelta  \left| \log \gmcdelta \right|^{\nf 12} \bigr) \,.
\end{equation*}

\smallskip

\emph{Step 3.} We prove~\eqref{e.after.good.scales.uniform}. We split the indices~$j$ into the finitely many residue classes modulo the integer~$r(d)$ from the proof of Lemma~\ref{l.discounted.J.local}; as shown there, the summands with indices in the same class are independent. In each class, we apply~\eqref{e.CLT.scaling.O.Gamma.Two} to the centered summands, which are~$\O_{\Gamma_2}(Cs^{-\nf{3}{2}}\gmcdelta|\log\gmcdelta|^{\nf{1}{2}})$ by the last display, and combining the result with the preceding estimates we obtain
\begin{align*}
	\lefteqn{
	\sum_{j=n}^{m}
	\sup_{0\leq l \leq j-2} 3^{-\frac{1}{4}s (j-l)}
	\max_{z \in 3^l \Zd \cap (\cu_j \setminus \cu_{j-1})}
	\min\Bigl\{  \max_{|e| = 1} J(z+\cu_l,e,e\,; \ahom_{l}^{-1} \a_{l})\, ,  1 \Bigr\}^{\nf 12}
	}  \qquad &
	\notag \\ &
	\leq
	C s^{-\nf 12} \left| \log s \right|^{\nf 12}  \gmcdelta (m-n+1)+ \O_{\Gamma_2}\bigl( C s^{-\nf 32} \gmcdelta \left| \log \gmcdelta \right|^{\nf 12} (m-n+1)^{\nf 12} \bigr)
	\,.
\end{align*}
Applying the last display before Step~3 to the annuli with~$j<n$ and summing the resulting geometric series gives
\begin{multline*}
	\sum_{j=0}^{n-1}
	\sup_{0\leq l \leq j-2} 3^{-\frac{1}{4}s (n-l)}
	\max_{z \in 3^l \Zd \cap (\cu_j \setminus \cu_{j-1})}
	\min\Bigl\{  \max_{|e| = 1} J(z+\cu_l,e,e\,; \ahom_{l}^{-1} \a_{l})\, ,  1 \Bigr\}^{\nf 12}
	\\
	\leq
	\O_{\Gamma_2}\bigl( C s^{-\nf 52}  \gmcdelta  \left| \log \gmcdelta \right|^{\nf 12} \bigr)
	\,.
\end{multline*}
The last two displays, together with~\eqref{e.split.J.for.minimal.scale} and~\eqref{e.sum.the.errors.fields}, give~\eqref{e.after.good.scales.uniform}.

\smallskip

\emph{Step 4.} We define~$M(\lambda,s,m)$ and prove~\eqref{e.M.lambda.m.bound}. Define
\begin{align*}
	M(\lambda,s,m)
	\coloneqq
	\min\!\left(
	\biggl\{ n \in \N_0 \; : \;  n \leq m-1\,, \ \max_{z\in3^{n} \Zd \cap \cu_m}
	\sum_{j=n}^{m}
	\mathbf{G}_{j,z}(s)
	>
	\lambda (m-n)
	\biggr\}\cup\{m\}\right)-1
	\,.
\end{align*}
The averaged estimate~\eqref{e.after.good.scales.uniform}, stationarity, and the sub-Gaussian tail in that estimate show that, for every~$0\leq k\leq m-1$,
\begin{equation*}
	\P\!\left[
	\max_{z\in3^k\Zd\cap\cu_m}
	\sum_{j=k}^m\mathbf G_{j,z}(s)>\lambda(m-k)
	\right]
	\leq
	\exp\!\left(
	C(m-k)-\frac{s^7\lambda^2(m-k)}{C\gmcdelta^2|\log\gmcdelta|}
	\right)\,.
\end{equation*}
Indeed, there are at most~$C3^{d(m-k)}$ centers, and we apply~\eqref{e.after.good.scales.uniform}, translated, at each of them. The same bound is valid for~$k=m$ because its right side is one. Consequently, a union bound gives, for every~$\lambda \geq 2 Cs^{-\nf 72}   \gmcdelta \left| \log \gmcdelta \right|^{\nf 12} $ and~$0\leq n\leq m-1$,
\begin{align}
	\P \bigl[ M(\lambda,s,m) \leq n \bigr]
	&
	\leq
	\sum_{k=0}^{n+1}
	\P \biggl[
	\max_{z\in3^k \Zd \cap \cu_m}
	\sum_{j=k}^{m}
	\mathbf{G}_{j,z}(s)
	>
	\lambda (m-k)
	\biggr]
	\notag \\ &
	\leq
	\sum_{k=0}^{n+1}
	\exp \biggl( C (m-k)
	-\frac{s^7 \lambda^2  (m-k) }{C \gmcdelta^{2} \left| \log \gmcdelta \right| }
	\biggr)
	\leq
	C\exp \biggl(
	-\frac{s^7 \lambda^2  (m-n-1) }{C\gmcdelta^{2} \left| \log \gmcdelta \right|}
	\biggr)\,.
	\label{e.M.lambda.m.tail}
\end{align}
This implies~\eqref{e.M.lambda.m.bound} and completes the proof.
\end{proof}

\begin{lemma}
\label{l.min.scale.good.scale}
There exists~$C(d)<\infty$ such that, for every~$\lambda,\ep ,s \in (0,1]$ satisfying
\begin{equation}
	\lambda \geq C s^{-8}   \ep^{-2} \gmcdelta^{2} \left| \log \gmcdelta\right|
	\label{e.holder.induction.smallness}
\end{equation}
and for every~$m \in \N_0$, there exists a random variable~$N(\lambda,\ep,s,m)\in\{-1,0,\ldots,m\}$ such that
\begin{equation}
	0\leq m-N(\lambda,\ep,s,m)
	\leq
	1+\O_{\Gamma_1}\bigl(
	Cs^{-8}\ep^{-2}\lambda^{-1}\gmcdelta^2|\log\gmcdelta|
	\bigr)
	\label{e.N.theta.m.bound.goodscales}
\end{equation}
and, for every~$n \in \N_0$ with~$n \leq N(\lambda,\ep,s,m)$,
\begin{equation}
	\min_{z\in 3^n\Zd \cap \cu_m}
	\sum_{j=n}^{m}
	\indc{\bigl\{ \mathcal{G}_{j,z}(\ep,s) \bigr\}}
	>
	(1-\lambda) (m-n)  \,.
	\label{e.goodscales.minscale}
\end{equation}
Equivalently, every such~$n$ satisfies
\begin{equation*}
	\max_{z\in3^n\Zd\cap\cu_m}
	\sum_{j=n}^m
	\bigl(1-\indc{\{\mathcal G_{j,z}(\ep,s)\}}\bigr)
	<1+\lambda(m-n)\,.
\end{equation*}
\end{lemma}
\begin{proof}
Define the random variable~$N(\lambda,\ep,s,m)$ by
\begin{align*}
	N(\lambda,\ep,s,m)
	&\coloneqq
	\min\!\Biggl(
	\Biggl\{n\in\N_0:n\leq m,\quad
	\min_{z\in3^n\Zd\cap\cu_m}
	\sum_{j=n}^{m}\indc{\{\mathcal G_{j,z}(\ep,s)\}}
	\\
	&\hspace{47mm}\leq(1-\lambda)(m-n)
	\Biggr\}\cup\{m+1\}\Biggr)-1
	\,.
\end{align*}
The definition of~$N(\lambda,\ep,s,m)$ gives~\eqref{e.goodscales.minscale}. On the event in the next display, the proportion of bad scales among the~$m-n+1$ scales from~$n$ to~$m$ is at least~$\lambda/2$. Proposition~\ref{p.density.of.good.scales}, applied with~$\lambda/2$ in place of~$\theta$, therefore gives
\begin{equation*}
	\P \Biggl[
	\sum_{j=n}^{m}
	(1-\indc_{\mathcal{G}_{j,z}(\ep,s)})
	>
	\lambda (m-n)  \Biggr]
	\leq
	\exp \biggl(
	-\frac{s^8\ep^2\lambda(m-n+1)}
	{C\gmcdelta^2|\log\gmcdelta|}
	\biggr)
	\,.
\end{equation*}
By a union bound, stationarity, and the condition~\eqref{e.holder.induction.smallness} with sufficiently large~$C$, we obtain
\begin{align*}
	\P \Biggl[
	\max_{z\in 3^n\Zd \cap \cu_m}
	\sum_{j=n}^{m}
	(1-\indc_{\mathcal{G}_{j,z}(\ep,s)})
	>
	\lambda (m-n)  \Biggr]
	&
	\leq
	3^{d(m-n)}
	\P \Biggl[
	\sum_{j=n}^{m}
	(1-\indc_{\mathcal{G}_{j}(\ep,s)})
	>
	\lambda (m-n)  \Biggr]
	\notag \\ &
	\leq
	\exp \biggl(
	2d(m-n)
	-\frac{s^8\ep^2\lambda(m-n+1)}
	{C\gmcdelta^2|\log\gmcdelta|}
	\biggr)
	\notag \\ &
	\leq
	\exp \biggl(
	-\frac{s^8\ep^2\lambda(m-n+1)}
	{C\gmcdelta^2|\log\gmcdelta|}
	\biggr)
	\,.
\end{align*}
By another union bound, for~$0\leq n\leq m-1$,
\begin{align*}
	\P \bigl[ N(\lambda,\ep,s,m) \leq n \bigr]
	&
	\leq
	\sum_{k=0}^{n+1}
	\P \Biggl[
	\min_{z\in 3^k\Zd \cap \cu_m}
	\sum_{j=k}^{m}
	\indc_{\mathcal{G}_{j,z}(\ep,s) }
	\leq
	(1-\lambda) (m-k)  \Biggr]
	\notag \\ &
	\leq
	C\exp \biggl(
	-\frac{s^8\ep^2\lambda(m-n)}
	{C\gmcdelta^2|\log\gmcdelta|}
	\biggr)
	\,.
\end{align*}
This implies~\eqref{e.N.theta.m.bound.goodscales}, completing the proof.
\end{proof}

\begin{proof}[{Proof of Proposition~\ref{p.Holder.regularity}}]
The proof has seven steps. Step~1 fixes the parameters, and Step~2 defines the minimal scale~$\mathcal L_L(\gamma,m)$, proves~\eqref{e.mathcal.L.bound}, and records the bounds~\eqref{e.after.good.scales.applied} for~$n\leq m-\mathcal L_L(\gamma,m)$. Step~3 turns Lemma~\ref{l.excess.decay.good.scales.GMC} into the one-step inequality~\eqref{e.ready.for.iteration.lemma} on the good scales, and Step~4 iterates it by Lemma~\ref{l.iteration.lemma.GMC}. Steps~5,~6 and~7 deduce~\eqref{e.Holder.estimate.boxes.local},~\eqref{e.energy.density.estimate} and~\eqref{e.excess.decay.estimate}, respectively. The heart of the matter is Step~4: above the minimal scale, the number of bad scales and the sum of the errors~$\ep_{j,z}$ are at most~$C+C\lambda(m-n)$, so the factor produced by the iteration is at most~$C\exp(C\lambda(m-n))$; the choice of~$\lambda$ in Step~1 makes this at most~$C3^{\frac{1}{4}(1-\gamma)(m-n)}$.

\smallskip

\emph{Step 1.} We fix the parameters. Fix~$\gamma$ in the range allowed by the statement. Choose the dimensional constants in the order indicated below and set
\begin{equation}
	s_0\coloneqq \nf{1}{32}\,,
	\qquad
	\lambda\coloneqq C_1^{-1}(1-\gamma)\,,
	\qquad
	\ep\coloneqq C_2^{-1}(1-\gamma)^{\nf12}\,,
	\qquad
	k\coloneqq\left\lceil4\log_3\!\left(4s_0^{-\nf32}C_3\right)\right\rceil\,,
	\label{e.choices.hepk.Holder.proof}
\end{equation}
where~$C_3(d)$ is the constant of Lemma~\ref{l.excess.decay.good.scales.GMC}. Choose~$C_1(d)\geq2$ sufficiently large for both~\eqref{e.exp.lambda.bound} below and the absorption in the proof of~\eqref{e.energy.density.estimate} in Step~6, and then choose~$C_2(d)$ so that~$C_1C_2^{-8}\leq1$ and
\begin{equation*}
	C_3
	\left(3^{-\frac{1}{2}k}+3^{(1+\frac{d}{2})k}s_0^{-\nf32}C_2^{-1}\right)
	\leq\frac12\,3^{-\frac{1}{4}k}\,.
\end{equation*}
Finally, choose the constant~$C$ of the statement large enough, depending on~$C_1$ and~$C_2$. The lower bound on~$1-\gamma$ then ensures, for all sufficiently small~$\gmcdelta$, that the hypotheses of Lemma~\ref{l.excess.decay.good.scales.GMC}, Proposition~\ref{p.good.scale.mathcal.E}, Lemma~\ref{l.sum.the.errors}, and Lemma~\ref{l.min.scale.good.scale} hold.

\smallskip

\emph{Step 2.} We define the minimal scale. For~$m\in\N_0$, define
\begin{equation}
	\mathcal{L}_L(\gamma,m) \coloneqq
	\max\bigl\{ m - M(\lambda,s_0,m) \,,\, m - N(\lambda,\ep,s_0,m) \bigr\} + k + 5\,,
	\label{e.L.def.reg}
\end{equation}
where~$M$ and~$N$ are supplied by Lemmas~\ref{l.sum.the.errors} and~\ref{l.min.scale.good.scale}. They do not involve~$L$, so~$\mathcal L_L(\gamma,m)$ is the same random variable for every~$L\geq m$.  The summand~$k+5$ ensures that every scale used in Lemma~\ref{l.excess.decay.good.scales.GMC} is at most~$m-5$. Estimates~\eqref{e.M.lambda.m.bound} and~\eqref{e.N.theta.m.bound.goodscales} yield
\begin{align*}
	\mathcal L_L(\gamma,m)
	&\leq C+\O_{\Gamma_1}\bigl( Cs_0^{-7}\gmcdelta^2\left|\log\gmcdelta\right|\lambda^{-2}
	+ Cs_0^{-8}\ep^{-2}\lambda^{-1}\gmcdelta^2\left|\log\gmcdelta\right| \bigr)
	\notag \\ &
	\leq C+\O_{\Gamma_1}\bigl( CC_1^2(1-\gamma)^{-2}\gmcdelta^2\left|\log\gmcdelta\right| \bigr)\,.
\end{align*}
This is~\eqref{e.mathcal.L.bound}. Finally,~\eqref{e.after.good.scales} and~\eqref{e.goodscales.minscale} imply, whenever~$n\leq m-\mathcal L_L(\gamma,m)$,
\begin{equation}
	\max_{z\in 3^n\Zd\cap\cu_m} \sum_{j=n}^{m} \mathbf{G}_{j,z}(s_0) \leq \lambda(m-n)
	\qand
	\max_{z\in 3^n\Zd\cap\cu_m} \sum_{j=n}^{m}
	\bigl(1 - \indc{\bigl\{\mathcal{G}_{j,z}(\ep,s_0)\bigr\}}\bigr)
	<1+\lambda(m-n)\,.
	\label{e.after.good.scales.applied}
\end{equation}

\smallskip

\emph{Step 3.} We prove the excess decay estimate~\eqref{e.excess.decay.one.step} on good scales. Fix~$m\in\N_0$ with~$L\geq m$, and let~$u\in H^1(\cu_m)$ be a solution of~\eqref{e.eq.for.u.holder.estimate}. Fix~$n\in\N_0$ with~$n \leq m-\mathcal{L}_L(\gamma,m)$ and~$z\in 3^n\Zd\cap\cu_m$. We apply Lemma~\ref{l.excess.decay.good.scales.GMC} with~$s = 8s_0 = \nf14$ at the scales~$j\in\N_0\cap[n+k,m-5]$, for which its hypotheses hold. For each such~$j$, the conclusion~\eqref{e.excess.decay.one.step}, with~$j$ in place of~$n$ and with~$x=z$, gives
\begin{align}
	\lefteqn{
	E(u,(z+\cu_{j-k})\cap\cu_m)\indc_{\mathcal{G}_{j+2,z}(\ep,s_0)}
	} \qquad &
	\notag \\ &
	\leq
	C_3\bigl(3^{-\frac{1}{2} k}
	+ 3^{(1+\frac{d}{2})k}s_0^{-\nf32}\ep\bigr) E(u,(z+\cu_{j})\cap\cu_m)
	\notag \\ & \qquad
	+ C_3 3^{(1+\frac{d}{2})k} s_0^{-\nf32}
	\mathcal{E}_{s_0,\infty,2}(z+\cu_{j+2};\a_L,(\b_{L,j+2})_{z+\cu_{j+2}})
	\Bigl(\bigl|\nabla\ell(u,(z+\cu_{j})\cap\cu_m)\bigr|
	\notag \\ & \qquad\qquad
	+ s_0^{-\nf32}\bigl|(\nabla h)_{(z+\cu_{j})\cap\cu_m}\bigr|
	\indc_{\{\overline{(z+\cu_{j})\cap\cu_m}\cap\partial\cu_m\neq\emptyset\}}\Bigr)
	\indc_{\mathcal{G}_{j+2,z}(\ep,s_0)}
	\notag \\ & \qquad
	+ Cs_0^{-\nf{15}2} 3^{(1+\frac{d}{2})k}
	(\b_{L,j+2})_{z+\cu_{j+2}}^{-1} 3^{8s_0 j}[\g]_{\underline{H}^{8s_0}((z+\cu_{j})\cap\cu_m)}
	\notag \\ & \qquad
	+ Cs_0^{-3} 3^{(1+\frac{d}{2})k}
	3^{\frac{1}{2} j}[\nabla h]_{W^{\nf12,\infty}((z+\cu_{j})\cap\cu_m)}
	\indc_{\{\overline{(z+\cu_{j})\cap\cu_m}\cap\partial\cu_m\neq\emptyset\}}\,.
	\label{e.one.step.expanded}
\end{align}
By our choice of~$k$ and~$C_2$, we have~$C_3(3^{-\frac{1}{2}k} + 3^{(1+\frac{d}{2})k}s_0^{-\nf32}\ep) \leq \frac{1}{2}3^{-\frac{1}{4}k}$. Definition~\ref{d.bLm} and the fact that an average lies between the infimum and the supremum give
\begin{align*}
	&\max\!\left\{
	\frac{(\b_{L,j+2})_{z+\cu_{j+2}}}{(\b_{L,m})_{\cu_m}},\frac{(\b_{L,m})_{\cu_m}}{(\b_{L,j+2})_{z+\cu_{j+2}}}
	\right\}
	\\
	&\quad\leq
	\max\!\left\{
	\frac{\ahom_{j+2}}{\ahom_m},\frac{\ahom_m}{\ahom_{j+2}}
	\right\}
	\exp\!\left(
	\sum_{i=j+3}^m\|\g_i-\tau^2\|_{L^\infty(z+\cu_{j+2})}
	+\osc_{\cu_m\cup(z+\cu_{j+2})}\sum_{i=m+1}^L\g_i
	\right)\,.
\end{align*}
Equations~\eqref{e.g.m.j} and~\eqref{e.def.bfG.mhq} give
\begin{equation*}
	\sum_{i=j+3}^m\|\g_i-\tau^2\|_{L^\infty(z+\cu_{j+2})}
	\leq \tau^2(m-j-2)+C\sum_{i=j+2}^m\mathbf G_{i,z}(s_0)\,.
\end{equation*}
By the mean-value theorem and~\eqref{e.def.bfG.mhq},
\begin{align*}
	\osc_{\cu_m\cup(z+\cu_{j+2})}\sum_{i=m+1}^L\g_i
	&\leq C\sum_{i=m+1}^\infty 3^m
	\left(
	\|\nabla\g_i\|_{L^\infty(z+\cu_m)}
	+\|\nabla\g_i\|_{L^\infty(\cu_m)}
	\right)
	\leq C\bigl(\mathbf G_{m,z}(s_0)+\mathbf G_{m,0}(s_0)\bigr)\,.
\end{align*}
Applying~\eqref{e.after.good.scales.applied} at the centers~$z$ and~$0$, we obtain~$\sum_{i=j+2}^m\mathbf G_{i,z}(s_0)+\mathbf G_{m,0}(s_0)\leq C\lambda(m-n)$. The parameter choice in Step~1 gives~$\tau^2\leq C\gmcdelta^2\leq C\lambda$. Together with~\eqref{e.annealed.ordering}, the preceding estimates show that, for every~$j\in\N_0\cap[n,m-5]$,
\begin{equation}
	C^{-1}\exp\bigl(-C\lambda(m-n)\bigr)
	\leq (\b_{L,j+2})_{z+\cu_{j+2}}(\b_{L,m})_{\cu_m}^{-1}
	\leq C\exp\bigl(C\lambda(m-n)\bigr)\,.
	\label{e.ratio.of.bs}
\end{equation}
Moreover,
\begin{equation*}
	3^{8s_0 j}[\g]_{\underline{H}^{8s_0}((z+\cu_{j})\cap\cu_m)}
	\leq C3^{\frac{1}{2} j}[\g]_{W^{\nf12,\infty}((z+\cu_{j})\cap\cu_m)}
	\leq C3^{-\frac{1}{2}(m-j)}\bigl(3^{\frac{1}{2} m}[\g]_{W^{\nf12,\infty}(\cu_m)}\bigr)\,,
\end{equation*}
so, using~\eqref{e.ratio.of.bs},
\begin{multline}
	(\b_{L,j+2})_{z+\cu_{j+2}}^{-1} 3^{8s_0 j}[\g]_{\underline{H}^{8s_0}((z+\cu_{j})\cap\cu_m)}
	\\
	\leq C\exp\bigl(C\lambda(m-n)\bigr) 3^{-\frac{1}{2}(m-j)}
	\bigl((\b_{L,m})_{\cu_m}^{-1} 3^{\frac{1}{2} m}[\g]_{W^{\nf12,\infty}(\cu_m)}\bigr)\,.
	\label{e.forcing.j.bound}
\end{multline}
The boundary datum satisfies
\begin{equation}
	3^{\frac{1}{2} j}[\nabla h]_{W^{\nf12,\infty}((z+\cu_{j})\cap\cu_m)}
	\leq C3^{-\frac{1}{2}(m-j)}\bigl(3^{\frac{1}{2} m}\|\nabla h\|_{\underline{W}^{\nf12,\infty}(\cu_m)}\bigr)\,.
	\label{e.boundary.j.bound}
\end{equation}
Furthermore, we use the following elementary inequality: if~$D\geq0$ and~$0\leq x\leq\ep 3^{s_0D/16}$, then~$3^{-s_0D/2}x\leq 3^{-s_0D/2}\min\{x,1\}+\ep^8$. Indeed, the second term is needed only when~$x>1$, in which case~$3^{-s_0D/16}<\ep$. By~\eqref{e.good.sum.J}, on~$\mathcal G_{j+2,z}(\ep,s_0)$ it applies to the square root~$x$ of every quantity~$J$ on the right side of~\eqref{e.mathcalE.bound.for.good}, with~$D$ the difference between the scale~$j+2$ and the scale of that quantity.  Proposition~\ref{p.good.scale.mathcal.E} and~\eqref{e.def.bfG.mhq} therefore give
\begin{equation}
	\mathcal{E}_{s_0,\infty,2}(z+\cu_{j+2};\a_L,(\b_{L,j+2})_{z+\cu_{j+2}})
	\indc_{\mathcal{G}_{j+2,z}(\ep,s_0)}
	\leq C\min\bigl\{\ep\,,\,
	s_0^{-1}\gmcdelta^2+\ep^8+\mathbf{G}_{j+2,z}(s_0)\bigr\}\,.
	\label{e.mathcalE.bound.applied}
\end{equation}
Inserting these bounds, including~\eqref{e.forcing.j.bound},~\eqref{e.boundary.j.bound} and~\eqref{e.mathcalE.bound.applied}, into~\eqref{e.one.step.expanded} yields, for every~$z\in 3^n\Zd\cap\cu_m$ and~$j\in\N_0\cap[n+k,m-5]$,
\begin{align}
	\lefteqn{
	E(u,(z+\cu_{j-k})\cap\cu_m)\indc_{\mathcal{G}_{j+2,z}(\ep,s_0)}
	} \qquad & \notag \\ &
	\leq 3^{-\frac{1}{4}k} E(u,(z+\cu_{j})\cap\cu_m)
	+ \ep_{j,z}\bigl|\nabla\ell(u,(z+\cu_{j})\cap\cu_m)\bigr|
	+ \delta_{j,z}\,,
	\label{e.ready.for.iteration.lemma}
\end{align}
where we define
\begin{equation}
	\ep_{j,z} \coloneqq
	Cs_0^{-\nf32}\min\bigl\{\ep\,,\,\gmcdelta^2+\ep^8+\mathbf{G}_{j+2,z}(s_0)\bigr\}
	\indc_{\mathcal{G}_{j+2,z}(\ep,s_0)}
	\label{e.ep.j.z.def}
\end{equation}
and
\begin{align*}
	\delta_{j,z} &\coloneqq
	C3^{-\frac{1}{2}(m-j)}\exp\bigl(C\lambda(m-n)\bigr)
	(\b_{L,m})_{\cu_m}^{-1} 3^{\frac{1}{2} m}[\g]_{\underline{W}^{\nf12,\infty}(\cu_m)}
	\notag \\ & \qquad
	+ C\Bigl(\ep_{j,z}\|\nabla h\|_{L^\infty(\cu_m)}
	+ 3^{-\frac{1}{2}(m-j)} 3^{\frac{1}{2} m}\|\nabla h\|_{\underline{W}^{\nf12,\infty}(\cu_m)}
	\Bigr)\indc_{\{\overline{(z+\cu_{j})\cap\cu_m}\cap\partial\cu_m\neq\emptyset\}} \,.
\end{align*}
The geometric factor~$3^{-\frac{1}{2}(m-j)}$ in the last term comes from~\eqref{e.boundary.j.bound}. For~$j\leq m'$, since the sets~$(z+\cu_j)\cap\cu_m$ are nested, if~$\overline{(z+\cu_j)\cap\cu_m}$ meets~$\partial\cu_m$, then~$\overline{(z+\cu_{m'})\cap\cu_m}$ does as well. By~\eqref{e.after.good.scales.applied}, we have
\begin{equation}
	\sum_{j=n}^{m-2}\ep_{j,z}
	\leq Cs_0^{-\nf32}\biggl((\gmcdelta^2+\ep^8)(m-n)
	+ \sum_{j=n}^{m-2}\mathbf{G}_{j+2,z}(s_0)\biggr)
	\leq C\lambda(m-n)\,.
	\label{e.ep.j.z.sum}
\end{equation}
Therefore, we also obtain, for every~$m',n'\in\N_0$ with~$n\leq n'\leq m'\leq m-3$,
\begin{align}
	\sum_{j=n'}^{m'}\delta_{j,z}
	&\leq C\exp\bigl(C\lambda(m-n)\bigr)(\b_{L,m})_{\cu_m}^{-1}
	3^{\frac{1}{2} m'}[\g]_{\underline{W}^{\nf12,\infty}(\cu_m)}
	\notag \\ & \qquad
	+ C\Bigl(\lambda(m-n)\|\nabla h\|_{L^\infty(\cu_m)}
	+ 3^{\frac{1}{2} m'}\|\nabla h\|_{\underline{W}^{\nf12,\infty}(\cu_m)}
	\Bigr)\indc_{\{\overline{(z+\cu_{m'})\cap\cu_m}\cap\partial\cu_m\neq\emptyset\}}
	\,.
	\label{e.delta.j.z.sum}
\end{align}

\smallskip

\emph{Step 4.} We set up the iteration. The inequality~\eqref{e.ready.for.iteration.lemma} is the key hypothesis~\eqref{e.Ej.decay.assumption} of Lemma~\ref{l.iteration.lemma.GMC} (with step~$k$ and~$\theta=3^{-\nf14}$). For~$z\in 3^n\Zd\cap\cu_m$, we define the set~$\mathcal B$ of indices by
\begin{equation}
	\mathcal{B} \coloneqq
	\bigl\{j\in\N_0\cap[n,m-2] : j<n+k
	\ \text{or}\ \mathcal{G}_{j+2,z}(\ep,s_0) \text{ fails}\bigr\}\,.
	\label{e.bad.bad.z}
\end{equation}
By~\eqref{e.after.good.scales.applied} and the definition~\eqref{e.bad.bad.z},
\begin{equation}
	|\mathcal B|
	\leq k+\sum_{j=n+k}^{m-2}\bigl(1-\indc{\bigl\{\mathcal{G}_{j+2,z}(\ep,s_0)\bigr\}}\bigr)
	\leq k+1+\lambda(m-n)\,.
	\label{e.Bz.bound}
\end{equation}
This and~\eqref{e.ep.j.z.sum} imply, for the constant~$A$ in Lemma~\ref{l.iteration.lemma.GMC}, that~$e^A\leq C\exp\bigl(C\lambda(m-n)\bigr)$. Since~$k$ is a dimensional constant,~\eqref{e.Bz.bound} shows that~$3^{C(k+1)(|\mathcal B|+1)/4}e^A$ satisfies the same bound, with a larger~$C$. We apply~\eqref{e.excess.decay.lemma} with~$M=C\ep$, which is admissible because~$\ep_{j,z}\leq C\ep$ by~\eqref{e.ep.j.z.def}. For every~$n',m'\in\N_0$ with~$n\leq n'\leq m'\leq m-5$, Lemma~\ref{l.iteration.lemma.GMC} gives the following two inequalities (for~$n'=m'$ they are immediate). The first inequality, from~\eqref{e.iteration.slope.lemma}, is
\begin{multline}
	3^{-n'}\|u-(u)_{(z+\cu_{n'})\cap\cu_m}\|_{\underline{L}^2((z+\cu_{n'})\cap\cu_m)}
	\\
	\leq C\exp\bigl(C\lambda(m-n)\bigr)
	\biggl(3^{-m'}\|u-(u)_{(z+\cu_{m'})\cap\cu_m}\|_{\underline{L}^2((z+\cu_{m'})\cap\cu_m)}
	+ \sum_{j=n'}^{m'}\delta_{j,z}\biggr)
	\label{e.iteration.slope.applied}
\end{multline}
and the second inequality, the excess decay estimate, is
\begin{align}
	E(u,(z+\cu_{n'})\cap\cu_m)
	&
	\leq C\exp\bigl(C\lambda(m-n)\bigr)
	\biggl(3^{-\frac{1}{4}(m'-n')} E(u,(z+\cu_{m'})\cap\cu_m)
	\notag \\ & \qquad
	+ \ep\, 3^{-m'}\|u-(u)_{(z+\cu_{m'})\cap\cu_m}\|_{\underline{L}^2((z+\cu_{m'})\cap\cu_m)}
	+ \sum_{j=n'}^{m'}\delta_{j,z}\biggr)\,.
	\label{e.excess.decay.applied}
\end{align}

\smallskip

\emph{Step 5.} We prove the large-scale H\"older estimate~\eqref{e.Campanato.full.domain}. We show that, for every~$n,n',m'\in\N_0$ with~$n\leq m-\mathcal L_L(\gamma,m)$ and~$n\leq n'\leq m'\leq m-5$ and every~$z\in 3^n\Zd\cap\cu_m$, we have the estimate
\begin{align}
	\lefteqn{3^{-n'}\|u-(u)_{(z+\cu_{n'})\cap\cu_m}\|_{\underline{L}^2((z+\cu_{n'})\cap\cu_m)}} \qquad &
	\notag \\ &
	\leq C3^{\frac{1}{2}(1-\gamma)(m-n)}
	\Bigl(3^{-m'}\|u-(u)_{(z+\cu_{m'})\cap\cu_m}\|_{\underline{L}^2((z+\cu_{m'})\cap\cu_m)}
	+ (\b_{L,m})_{\cu_m}^{-1} 3^{\frac{1}{2} m}[\g]_{\underline{W}^{\nf12,\infty}(\cu_m)}\Bigr)
	\notag \\ & \qquad
	+ C3^{\frac{1}{2}(1-\gamma)(m-n)} 3^{\frac{1}{2} m}\|\nabla h\|_{\underline{W}^{\nf12,\infty}(\cu_m)}
	\indc_{\{\overline{(z+\cu_{m'})\cap\cu_m}\cap\partial\cu_m\neq\emptyset\}}
	\,.
	\label{e.Campanato.full.domain}
\end{align}
By the choice of~$C_1$ and the identity~$\lambda=C_1^{-1}(1-\gamma)$,
\begin{equation}
	(1+\lambda(m-n))\exp\bigl(C\lambda(m-n)\bigr)
	\leq C\exp\bigl(2CC_1^{-1}(1-\gamma)(m-n)\bigr)
	\leq C3^{\frac{1}{4}(1-\gamma)(m-n)}\,.
	\label{e.exp.lambda.bound}
\end{equation}
Combining~\eqref{e.iteration.slope.applied},~\eqref{e.delta.j.z.sum}, and~\eqref{e.exp.lambda.bound} gives~\eqref{e.Campanato.full.domain}. For each~$j\in\N_0$ with~$j\leq n$ and each~$y\in3^j\Zd\cap((x+\cu_n)\cap\cu_m)$, we apply this estimate with~$j$ in place of~$n$, with~$n'=j$,~$m'=m-5$, and~$z=y$. If~$x\in\cu_{m-1}$, then~$y\in x+\cu_n$ and~$n\leq m-\mathcal L_L(\gamma,m)\leq m-k-5$, so~$\overline{(y+\cu_{m-5})\cap\cu_m}\subset\cu_m$ and the boundary indicator in~\eqref{e.Campanato.full.domain} vanishes; hence this indicator is always bounded by~$\indc_{\{x\notin\cu_{m-1}\}}$. Together with~$3^{-\frac{1}{2}(1+\gamma)(m-j)}\leq 3^{-\gamma(m-j)}$, this proves~\eqref{e.Holder.estimate.boxes.local}.

\smallskip

\emph{Step 6.} We prove the energy estimate~\eqref{e.energy.density.estimate}. We use Lemma~\ref{l.centered.energy.good.scale} at the nearest good anchor to convert the large-scale H\"older estimate~\eqref{e.Campanato.full.domain} into a gradient bound. The coarse-graining estimates of Lemma~\ref{l.coarse.graining.RHS} bound the oscillation of~$u$ at the larger good scale by the energy and the forcing.

\smallskip

Fix~$x\in\cu_m$. If~$m-n\leq2k+14$, then~$|(x+\cu_n)\cap\cu_m|\geq c3^{nd}\geq c|\cu_m|$, and~\eqref{e.energy.density.estimate} follows directly from the energy in~$\cu_m$. We may therefore assume that~$m-n>2k+14$. Then, since~$\lambda\leq\nf12$, the interval~$[n+2,m-5]$ contains more than~$k+1+\lambda(m-n)$ integers, so by~\eqref{e.Bz.bound} it contains a scale outside the set~$\mathcal B$ of Step~4, for the center~$z$ chosen below. Let~$y$ be the coordinatewise clamp of~$x$ to the centers of the scale-$(n+1)$ cubes contained in~$\cu_m$; if~$x\in\cu_{m-1}$, then no coordinate is clamped and~$y=x$. Let~$z\in 3^n\Zd\cap\cu_m$ be a nearest neighbor of~$y$ in~$3^n\Zd\cap\cu_m$, chosen coordinatewise. We define the elements of~$[n+2,m-5]\setminus\mathcal B$ nearest to~$n$ and to~$m$ by
\begin{equation*}
	n' \coloneqq \min\bigl\{j\in\N_0 : j\in[n+2,m-5]\setminus\mathcal B\bigr\}
	\quad\text{and}\quad
	m' \coloneqq \max\bigl\{j\in\N_0 : j\in[n',m-5]\setminus\mathcal B\bigr\}\,,
\end{equation*}
respectively. Since~$|\mathcal B|\leq k+1+\lambda(m-n)$ by~\eqref{e.Bz.bound}, we have~$\max\{n'-n,m-m'\}\leq C+\lambda(m-n)$. As in Step~5, the boundary indicator in~\eqref{e.Campanato.full.domain} vanishes if~$x\in\cu_{m-1}$, so it is bounded by~$\indc_{\{x\notin\cu_{m-1}\}}$. We also fix~$y'\in\cu_m\cap(z+\cu_{m'+1})$ such that~$(z+\cu_{m'})\cap\cu_m\subseteq y'+\cu_{m'}\subseteq\cu_m$. By~\eqref{e.bound.Lambdas.by.Es} and the last assertion of Proposition~\ref{p.good.scale.mathcal.E}, for every~$j\in[n,m-5]\setminus \mathcal B$ we have
\begin{equation}
	C^{-1}(\b_{L,j+2})_{z+\cu_{j+2}}
	\leq \lambda_{s_0,2}(z+\cu_{j+2};\a_L)
	\leq \Lambda_{s_0,2}(z+\cu_{j+2};\a_L)
	\leq C(\b_{L,j+2})_{z+\cu_{j+2}}\,.
	\label{e.ellipticity.good.scale.bounds}
\end{equation}
By Lemma~\ref{l.lambdas.stability},~\eqref{e.big.Lambda.stability},~\eqref{e.ellipticity.good.scale.bounds} at~$j=m'$ and~\eqref{e.ratio.of.bs}, the ellipticity quantities at the larger good scale satisfy, for~$q\in\{1,2\}$,
\begin{equation}
	\Lambda_{\nf 14,q}(y'+\cu_{m'};\a_L)(\b_{L,m})_{\cu_m}^{-1}
	+ \lambda_{\nf 18,q}^{-1}(y'+\cu_{m'};\a_L)(\b_{L,m})_{\cu_m}
	\leq
	C\exp\bigl(C\lambda(m-n)\bigr)\,.
	\label{e.Lambdas.in.reg.bounds}
\end{equation}
Moreover,
\begin{equation*}
	\left(\frac{|(z+\cu_{m})\cap\cu_m|}{|(z+\cu_{m'})\cap\cu_m|}\right)^{\nf12}
	\leq C3^{d(m-m')/2}
	\leq C\exp\bigl(C\lambda(m-n)\bigr)\,.
\end{equation*}
Let~$\a_0'=(\b_{L,n'+2})_{z+\cu_{n'+2}}$. The estimate near the boundary is applied at the available good anchor~$n'+2$. We have~$n'\geq n+2$ and~$\mathcal G_{n'+2,z}(\ep,s_0) \subseteq\mathcal G_{n'+2,z}(1,1/32)$. Moreover~$(x+\cu_n)\cap\cu_m\subseteq (z+\cu_{n+1})\cap\cu_m$. In a coordinate where clamping to~$\cu_m$ occurs, the coordinate of~$y$ is~$\pm(3^m-3^{n+1})/2$, which belongs to~$3^n\Z$, so that~$z_i=y_i$. In an unclamped coordinate,~$y_i=x_i$ and~$|z_i-x_i|\leq3^n/2$. These two observations give the inclusion.

\smallskip

Cover~$(x+\cu_n)\cap\cu_m$ by scale-$(n'-4)$ cubes with centers in~$(z+\cu_{n'-1})\cap\cu_m$. Only~$C(d)$ such cubes are needed, since~$n'-4\geq n-2$. Apply Lemma~\ref{l.centered.energy.good.scale} with~$(n,x,z)=(n',z,z)$ and~$s=1/4$ on each cube. The normalized volumes introduce a factor of at most~$C3^{d(n'-n)}$. If~$x\in\cu_{m-1}$, every boundary indicator for these windows vanishes, since~$n'\leq m-5$ and~$|z-x|_\infty\leq3^n/2$. Consequently
\begin{align*}
	\|\a_L^{1/2}\nabla u\|_{\underline L^2((x+\cu_n)\cap\cu_m)}
	&\leq C3^{d(n'-n)/2}
	(\a_0')^{1/2}3^{-n'}
	\|u-(u)_{(z+\cu_{n'})\cap\cu_m}\|_{\underline L^2((z+\cu_{n'})\cap\cu_m)}
	\notag \\ & \qquad
	+C3^{d(n'-n)/2}(\a_0')^{-1/2}3^{n'/4}
	[\g]_{\underline H^{1/4}((z+\cu_{n'})\cap\cu_m)}
	\notag \\ & \qquad
	+C3^{d(n'-n)/2}(\a_0')^{1/2}
	\bigl(|(\nabla h)_{(z+\cu_{n'})\cap\cu_m}|_2
	\notag \\ & \qquad\qquad
	+3^{n'/4}[\nabla h]_{\underline H^{1/4}((z+\cu_{n'})\cap\cu_m)}\bigr)
	\indc_{\{x\notin\cu_{m-1}\}}
	\,.
\end{align*}
The~$s$ powers are dimensional constants at~$s=1/4$. By~\eqref{e.ratio.of.bs} and~$n'-n\leq C+\lambda(m-n)$, the volume factor and both ratios~$(\a_0')^{\pm1/2}(\b_{L,m})_{\cu_m}^{\mp1/2}$ are bounded by~$C\exp(C\lambda(m-n))$. The scaled fractional embedding gives~$3^{n'/4}[F]_{\underline H^{1/4}((z+\cu_{n'})\cap\cu_m)} \leq C3^{m/2}[F]_{\underline W^{1/2,\infty}(\cu_m)}$. Thus
\begin{align*}
	e^{-C\lambda(m-n)}
	\|\a_L^{1/2}\nabla u\|_{\underline L^2((x+\cu_n)\cap\cu_m)}
	&\leq C(\b_{L,m})_{\cu_m}^{1/2}3^{-n'}
	\|u-(u)_{(z+\cu_{n'})\cap\cu_m}\|_{\underline L^2((z+\cu_{n'})\cap\cu_m)}
	\notag \\ & \qquad
	+C(\b_{L,m})_{\cu_m}^{-1/2}3^{m/2}
	[\g]_{\underline W^{1/2,\infty}(\cu_m)}
	\notag \\ & \qquad
	+C(\b_{L,m})_{\cu_m}^{1/2}3^{m/2}
	\|\nabla h\|_{\underline W^{1/2,\infty}(\cu_m)}
	\indc_{\{x\notin\cu_{m-1}\}}
	\,.
\end{align*}
We apply~\eqref{e.Campanato.full.domain} at~$n'$ and~$m'$ to the first term. Its large-scale input is~$3^{m-m'}$ times the quantity controlled by the following detachment estimate, and~$3^{m-m'}\leq C\exp(C\lambda(m-n))$ by the choice of~$m'$. Furthermore,~\eqref{e.nabla.u.detach} and the fractional Poincar\'e inequality give
\begin{align*}
	3^{-m}\|u-(u)_{(z+\cu_{m'})\cap\cu_m}\|_{\underline{L}^2((z+\cu_{m'})\cap\cu_m)}
	&\leq C3^{-m}\|u-(u)_{y'+\cu_{m'}}\|_{\underline{L}^2(y'+\cu_{m'})}
	\\ &\leq C3^{-\frac14 m}\|\nabla u\|_{\Besov{-\nf14}{2}{2}(y'+\cu_{m'})}\,.
\end{align*}
We apply Lemma~\ref{l.coarse.graining.RHS} and~\eqref{e.Lambdas.in.reg.bounds} to obtain
\begin{align*}
	\lefteqn{
	3^{-\frac14 m}\|\nabla u\|_{\Besov{-\nf14}{2}{2}(y'+\cu_{m'})}
	} \qquad &
	\notag \\ &
	\leq
	C\lambda_{\nf 18,2}^{-\nf12}(y'+\cu_{m'};\a_L)
	\|\a_L^{\nf12}\nabla u\|_{\underline{L}^2(y'+\cu_{m'})}
	+ C\lambda_{\nf 18,2}^{-1}(y'+\cu_{m'};\a_L)
	3^{\frac14 m}[\g]_{\underline{H}^{\nf14}(y'+\cu_{m'})}
	\notag \\
	&\leq C\exp\bigl(C\lambda(m-n)\bigr)
	\Bigl((\b_{L,m})_{\cu_m}^{-\nf12}\|\a_L^{\nf12}\nabla u\|_{\underline{L}^2(\cu_m)}
	+ (\b_{L,m})_{\cu_m}^{-1} 3^{\frac{1}{2} m}[\g]_{\underline{W}^{\nf12,\infty}(\cu_m)}\Bigr)\,.
\end{align*}
Combining these estimates with~\eqref{e.Campanato.full.domain}, applied at~$n'$, and~\eqref{e.exp.lambda.bound}, and enlarging~$C_1$ to include the volume and coefficient losses above, proves~\eqref{e.energy.density.estimate}.

\smallskip

\emph{Step 7.} We prove the excess decay estimate~\eqref{e.excess.decay.estimate}. Fix~$m,n,l$ satisfying~$\mathcal L_L(\gamma,m)\leq m-n\leq(1-\gamma)^{-1}$ and~$n\leq l\leq m-5$, and fix~$x\in\cu_m$. We may assume~$l>n+2$, since otherwise~\eqref{e.excess.decay.estimate} follows, after enlarging~$C$, by comparing the volumes of the sets involved. Choose~$z\in 3^n\Zd\cap\cu_m$ such that
\begin{equation*}
	(x+\cu_n)\cap\cu_m \subseteq (z+\cu_{n+1})\cap\cu_m \subseteq (z+\cu_{l-1})\cap\cu_m \subseteq (x+\cu_l)\cap\cu_m\,.
\end{equation*}
As in Step~5, the boundary indicator in~\eqref{e.delta.j.z.sum} vanishes when~$x\in\cu_{m-1}$, so it is bounded by~$\indc_{\{x\notin\cu_{m-1}\}}$. Using~$E(u,(x+\cu_n)\cap\cu_m)\leq CE(u,(z+\cu_{n+1})\cap\cu_m)$,~$E(u,(z+\cu_{l-1})\cap\cu_m)\leq CE(u,(x+\cu_l)\cap\cu_m)$, and
\begin{equation*}
	\|u-(u)_{(z+\cu_{l-1})\cap\cu_m}\|_{\underline{L}^2((z+\cu_{l-1})\cap\cu_m)}
	\leq C\|u-(u)_{(x+\cu_l)\cap\cu_m}\|_{\underline{L}^2((x+\cu_l)\cap\cu_m)}\,,
\end{equation*}
we obtain~\eqref{e.excess.decay.estimate} from~\eqref{e.excess.decay.applied},~\eqref{e.delta.j.z.sum}, and the assumption~$\mathcal L_L(\gamma,m)\leq(1-\gamma)^{-1}$, using that~$\exp\bigl(C\lambda(m-n)\bigr)\leq C$, which follows from the condition~$m-n\leq(1-\gamma)^{-1} = C_1\lambda^{-1}$.
\end{proof}

\subsection{Regularity uniformly in the cutoff}
\label{ss.cutoff.uniform.regularity}

We remove the restriction~$L\geq m$ from Proposition~\ref{p.Holder.regularity}. Recall from Definition~\ref{d.bLm} that~$\b_{L,n}=\ahom_L$ for~$n\geq L$.  For~$n\leq m$, we shall use the identity
\begin{equation}
	\frac{\a_{m\wedge L}}{\a_{n\wedge L}}
	=
	\exp\left(
	\sum_{i=n+1}^{m}
	(\g_i-\tau^2)\indc_{\{i\leq L\}}
	\right)
	\label{e.cutoff.regularity.centered.ratio}
\end{equation}
and, by~\eqref{e.annealed.ordering},
\begin{equation}
	1
	\leq
	\frac{\ahom_{n\wedge L}}{\ahom_{m\wedge L}}
	\leq
	\exp\bigl(2\tau^2((m\wedge L)-(n\wedge L))\bigr)
	\leq
	\exp\bigl(2\tau^2(m-n)\bigr)\,.
	\label{e.cutoff.regularity.D.ordering}
\end{equation}

\smallskip

In this subsection, for the cutoff~$L$, we modify the good event~$\mathcal G_m(\ep,s)$ of Definition~\ref{d.good.events} by replacing only the quantity~$J$ in~\eqref{e.good.sum.J} by~$J(z+\cu_n,e,e\,;\ahom_{n\wedge L}^{-1}\a_{n\wedge L})$. The conditions~\eqref{e.good.gj} and~\eqref{e.good.discounted.am.vs.aj} are left unchanged, and translated events are defined as before.  Likewise, we modify the quantity~$\mathbf G_{k,z}(s)$ in~\eqref{e.def.bfG.mhq}, for the cutoff~$L$, by replacing its term involving~$J$ with
\begin{equation*}
	\min\!\left\{
	\max_{|e|=1}J\left(
	z'+\cu_l,e,e\,;
	\ahom_{l\wedge L}^{-1}\a_{l\wedge L}
	\right),1
	\right\}^{\nf12}\,.
\end{equation*}
The terms involving the fields~$\g_j$ are unchanged. We retain the original notation. The modified event agrees with the original event when~$L\geq m$, and the modified quantity agrees with the original quantity when~$L\geq k$. Thus the arguments in the previous subsections still apply.

\begin{lemma}[Moments of~$J$ uniformly in the cutoff]
\label{l.cutoff.regularity.response}
There exists~$C(d)<\infty$ such that, for every~$1\leq p\leq C^{-1}\gmcdelta^{-2}|\log\gmcdelta|^{-1}$, and every~$L,n\in\N_0$,
\begin{equation}
	\left\|
	\max_{|e|=1}
	J\left(
	\cu_n,e,e\,;
	\ahom_{n\wedge L}^{-1}\a_{n\wedge L}
	\right)
	\right\|_{L^p(\P)}
	\leq
	Cp\log(2+p)\gmcdelta^2\,.
	\label{e.cutoff.regularity.response.moment}
\end{equation}
\end{lemma}

\begin{proof}
For~$n\leq L$, the estimate follows from Proposition~\ref{p.coarse.grained.bound} and the homogeneity identity
\begin{equation*}
	J(U,e,e\,;\ahom_n^{-1}\a_n)
	=
	J(U,\ahom_n^{-\nf12}e,\ahom_n^{\nf12}e\,;\a_n)\,.
\end{equation*}
If~$n>L$, partition~$\cu_n$ into the cubes~$z+\cu_L$,~$z\in3^L\Zd\cap\cu_n$.  By~\eqref{e.subadditivity.J},
\begin{equation*}
	J(\cu_n,e,e\,;\ahom_L^{-1}\a_L)
	\leq
	\avsum_{z\in3^L\Zd\cap\cu_n}J(z+\cu_L,e,e\,;\ahom_L^{-1}\a_L)\,.
\end{equation*}
By Minkowski's inequality and stationarity, the~$L^p$ norm of the right side is at most that of~$\max_{|e|=1}J(\cu_L,e,e\,;\ahom_L^{-1}\a_L)$, which is covered by the case~$n\leq L$.  The maximum over~$e$ is handled as in Proposition~\ref{p.coarse.grained.bound}.
\end{proof}

\begin{proposition}[Good-scale estimates uniformly in the cutoff]
\label{p.cutoff.regularity.good.scales}
There exists~$C(d)<\infty$ such that, for every~$L\in\N_0$, the following assertions hold, with~$\mathcal G_m$ and~$\mathbf G_{k,z}$ modified as above for the cutoff~$L$.

\smallskip

First, if~$s,\ep,\theta\in(0,1]$ and
\begin{equation}
	\theta
	\geq
	Cs^{-6}\ep^{-2}\gmcdelta^2|\log\gmcdelta|\,,
	\label{e.cutoff.regularity.density.parameters}
\end{equation}
then, for every~$m_0,M\in\N_0$,
\begin{equation}
	\P\left[
	\avsum_{m=m_0}^{m_0+M}
	\indc_{\mathcal G_m(\ep,s)}
	\leq
	1-\theta
	\right]
	\leq
	\exp\left(
	-\frac{s^6\ep^2\theta}
	{C\gmcdelta^2|\log\gmcdelta|}(M+1)
	\right)\,.
	\label{e.cutoff.regularity.density}
\end{equation}

\smallskip

Second, if~$s\in(0,\nf12]$ and~$s\geq C\gmcdelta^2|\log\gmcdelta|$, then, for every~$0\leq n\leq m$ and~$z\in\Rd$,
\begin{equation}
	\avsum_{k=n}^{m}\mathbf G_{k,z}(s)
	\leq
	Cs^{-\nf72}\gmcdelta
	+
	\O_{\Gamma_2}\left(
	Cs^{-\nf72}\gmcdelta|\log\gmcdelta|^{\nf12}
	(m-n+1)^{-\nf12}
	\right)\,.
	\label{e.cutoff.regularity.error.average}
\end{equation}
Consequently, if~$\lambda\in(0,1]$ and~$\lambda\geq Cs^{-\nf72}\gmcdelta|\log\gmcdelta|^{\nf12}$, there exists~$M(\lambda,s,m)\in\{-1,0,\ldots,m\}$ such that
\begin{equation}
	0\leq m-M(\lambda,s,m)
	\leq
	1+\O_{\Gamma_1}\left(
	Cs^{-7}\gmcdelta^2|\log\gmcdelta|\lambda^{-2}
	\right)
	\label{e.cutoff.regularity.error.stopping.tail}
\end{equation}
and, for every~$n\in\N_0$ with~$n\leq M(\lambda,s,m)$,
\begin{equation}
	\max_{z\in3^n\Zd\cap\cu_m}
	\sum_{j=n}^{m}\mathbf G_{j,z}(s)
	\leq
	\lambda(m-n)\,.
	\label{e.cutoff.regularity.error.stopping}
\end{equation}
If, in addition,~$\ep\in(0,1]$ and~$\lambda\geq Cs^{-8}\ep^{-2}\gmcdelta^2|\log\gmcdelta|$, there is~$N(\lambda,\ep,s,m)\in\{-1,0,\ldots,m\}$ satisfying
\begin{equation}
	0\leq m-N(\lambda,\ep,s,m)
	\leq
	1+\O_{\Gamma_1}\left(
	Cs^{-8}\ep^{-2}\lambda^{-1}
	\gmcdelta^2|\log\gmcdelta|
	\right)
	\label{e.cutoff.regularity.good.stopping.tail}
\end{equation}
and, for every~$n\in\N_0$ with~$n\leq N(\lambda,\ep,s,m)$,
\begin{equation}
	\max_{z\in3^n\Zd\cap\cu_m}
	\sum_{j=n}^{m}
	\left(1-\indc_{\mathcal G_{j,z}(\ep,s)}\right)
	<
	1+\lambda(m-n)\,.
	\label{e.cutoff.regularity.good.stopping}
\end{equation}

\smallskip

Finally, for~$s\in[64\gmcdelta^2,\nf12]$,~$\ep\in[s^{-1}\gmcdelta^2,1]$, every~$m\in\N_0$, and every~$z\in\Rd$,
\begin{equation}
	\mathcal E_{s,\infty,2}\left(
	z+\cu_m\,;\a_L,(\b_{L,m})_{z+\cu_m}
	\right)
	\indc_{\mathcal G_{m,z}(\ep,s)}
	\leq
	C\min\left\{\ep,s^{-1}\gmcdelta^2+\ep^8+\mathbf G_{m,z}(s)
	\right\}\,.
	\label{e.cutoff.regularity.good.E}
\end{equation}
In particular,
\begin{equation}
	\mathcal E_{s,\infty,2}\left(
	z+\cu_m\,;\a_L,(\b_{L,m})_{z+\cu_m}
	\right)
	\indc_{\mathcal G_{m,z}(\ep,s)}
	\leq
	C\ep\,.
	\label{e.cutoff.regularity.good.E.epsilon}
\end{equation}
\end{proposition}

\begin{proof}
The conditions~\eqref{e.good.gj} and~\eqref{e.good.discounted.am.vs.aj} in~$\mathcal G_m$ are those of Definition~\ref{d.good.events}. For the modified condition~\eqref{e.good.sum.J}, the proof of Lemma~\ref{l.discounted.J.local} applies, with the moments of~$J$ bounded by~\eqref{e.cutoff.regularity.response.moment}. Indeed, the coefficient~$\ahom_{n\wedge L}^{-1}\a_{n\wedge L}$ has range of dependence at most~$C3^{n\wedge L}\leq C3^n$, so~$J(z+\cu_n,e,e\,;\ahom_{n\wedge L}^{-1}\a_{n\wedge L})$ depends only on the fields in the~$C3^n$-neighborhood of~$z+\cu_n$.  Since the quantities~$J$ enter only with~$n\leq j-2$, splitting the indices~$j$ into finitely many residue classes, as in that proof, gives the finite-range hypothesis of Proposition~\ref{p.concentration.for.scales}.  The moment exponent in that proof is at least~$Cd/s$, and~\eqref{e.cutoff.regularity.density.parameters} ensures both the smallness condition required there and~$\theta\geq p^{-1}$.  This proves~\eqref{e.cutoff.regularity.density}.

\smallskip

For~\eqref{e.cutoff.regularity.error.average}, the terms involving the fields~$\g_j$ are the same as those in the proof of Lemma~\ref{l.sum.the.errors}.  For the term involving~$J$, we estimate the maximum over~$z$ as in Step~2 of that proof, using the moment of order~$8d/s$ of~$\min\{\max_{|e|=1}J(z+\cu_n,e,e\,;\ahom_{n\wedge L}^{-1}\a_{n\wedge L}),1\}^{\nf12}$, that is, the moment of order~$4d/s$ of the same quantity without the square root.  If the dimensional constant in the lower bound on~$s$ is large enough, this order lies in the range allowed by~\eqref{e.cutoff.regularity.response.moment}. With the same residue classes as above, the sum of these terms is then bounded exactly as in Step~3 of that proof, and the rearrangement~\eqref{e.split.J.for.minimal.scale} and the bound for the annuli below scale~$n$ are unchanged.  This proves~\eqref{e.cutoff.regularity.error.average}.  The union bound over translated centers used to prove~\eqref{e.M.lambda.m.tail} gives~\eqref{e.cutoff.regularity.error.stopping.tail} and~\eqref{e.cutoff.regularity.error.stopping}.  Applying the proof of Lemma~\ref{l.min.scale.good.scale} to~\eqref{e.cutoff.regularity.density} gives~\eqref{e.cutoff.regularity.good.stopping.tail} and~\eqref{e.cutoff.regularity.good.stopping}.

\smallskip

It remains to prove~\eqref{e.cutoff.regularity.good.E}.  By~\eqref{e.good.sum.J}, on~$\mathcal G_{m,z}(\ep,s)$ the elementary inequality of Step~3 of the proof of Proposition~\ref{p.Holder.regularity}, with~$s$ in place of~$s_0$, applies to the square root of every quantity~$J$ in the annular decomposition.  For~$m\leq L$, that inequality and~\eqref{e.mathcalE.bound.for.good} give~\eqref{e.cutoff.regularity.good.E}.  Suppose that~$m>L$.  We use the annular decomposition~\eqref{e.mathcalE.annular.decomp.pre} again and, at each nonnegative scale~$n$, compare~$\ahom_L^{-1}\a_L$ with~$\ahom_{n\wedge L}^{-1}\a_{n\wedge L}$.  The two coefficients agree if~$n\geq L$. If~$n<L$, we apply Lemma~\ref{l.J.sensitivity}, with the ratios of the two coefficients given by~\eqref{e.cutoff.regularity.centered.ratio}; the conditions~\eqref{e.good.gj} and~\eqref{e.good.discounted.am.vs.aj} in~$\mathcal G_m$ bound these ratios exactly as in~\eqref{e.J.sensitivity.applied.reg}--\eqref{e.J.localization.reg}, and the ratios of effective diffusivities are bounded by~\eqref{e.cutoff.regularity.D.ordering}.  The contribution of the negative scales is bounded by~\eqref{e.good.scale.subunit.tail}, with the same proof. Combining these estimates with the preceding inequality for these quantities, as in Step~3 of the proof of Proposition~\ref{p.good.scale.mathcal.E}, gives~\eqref{e.cutoff.regularity.good.E}.  The definition of~$\mathcal G_m(\ep,s)$ and~\eqref{e.what.e.good.gj.gives} give~\eqref{e.cutoff.regularity.good.E.epsilon}.
\end{proof}

\begin{lemma}[Regularity at a good scale, uniformly in the cutoff]
\label{l.cutoff.regularity.good.scale.estimates}
Let~$L\in\N_0$.  Then the following hold.
\begin{itemize}
\item \underline{\emph{Harmonic approximation}}: in the setting of Lemma~\ref{l.harmonic.approximation.good.scales.GMC}, with its hypothesis~$L\geq m$ removed, its conclusion holds with the good event~$\mathcal G_{n+2,z}$ modified as above for the cutoff~$L$.

\item \underline{\emph{Excess decay}}: in the setting of Lemma~\ref{l.excess.decay.good.scales.GMC}, with its hypothesis~$L\geq m$ removed, its conclusion holds with the same modification.
\end{itemize}
\end{lemma}

\begin{proof}
Use~$\a_0=(\b_{L,n+2})_{z+\cu_{n+2}}$ and Proposition~\ref{p.cutoff.regularity.good.scales} for the bound on the anchor error. The arbitrary-depth off-grid bound~\eqref{e.boundary.energy.offgrid.error} concerns the same coefficient~$\a_L$ and the same scalar~$\a_0$ on both sides. Consequently its ellipticity consequences, the full-datum competitor, both face-tile estimates and the window absorption in Lemma~\ref{l.centered.energy.good.scale} apply at every cutoff, including~$L<n+2$. The~$L^2$ comparison in Lemma~\ref{l.harmonic.approximation.good.scales.GMC} then gives the same three powers of~$s$ and the same powers of~$\a_0$. The harmonic comparison remains on the arbitrary cube allowed by that lemma. Multiplying the Laplacian by~$\a_0>0$ leaves its harmonic functions unchanged, so the interior and boundary regularity estimates for the comparison function and the transfer in Lemma~\ref{l.excess.decay.good.scales.GMC} apply with the cutoff-modified event.
\end{proof}

\begin{proposition}[Cutoff-uniform large-scale regularity]
\label{p.cutoff.Holder.regularity}
There exists~$C(d)<\infty$ such that, if~$\gmcdelta\leq C^{-1}$, then for every~$\gamma\in[\nf12,1-C\gmcdelta|\log\gmcdelta|^{\nf12}]$ and every~$L,m\in\N_0$, there is an~$\N$-valued random variable~$\mathcal L_L(\gamma,m)$ satisfying, for every~$k\in\N$,
\begin{equation}
	\P\left[
	\mathcal L_L(\gamma,m)>k
	\right]
	\leq
	C\exp\left(
	-\frac{(1-\gamma)^2(k-C)_+}
	{C\gmcdelta^2|\log\gmcdelta|}
	\right)\,.
	\label{e.cutoff.Holder.minimal.scale}
\end{equation}
For every~$\g\in C^{0,\nf12}(\cu_m;\Rd)$,~$h\in C^{1,\nf12}(\cu_m)$, and~$u\in H^1(\cu_m)$ satisfying~\eqref{e.eq.for.u.holder.estimate}, all three conclusions~\eqref{e.Holder.estimate.boxes.local},~\eqref{e.energy.density.estimate}, and~\eqref{e.excess.decay.estimate} hold with~$\mathcal L_L(\gamma,m)$ as the minimal scale. If~$L\geq m$, one may take
\begin{equation}
	\mathcal L_L(\gamma,m)=\mathcal L_m(\gamma,m)\,.
	\label{e.cutoff.Holder.replacements}
\end{equation}
The same statement holds on every translate~$y+\cu_m$, with the tail in~\eqref{e.cutoff.Holder.minimal.scale} uniform in~$y$.
\end{proposition}

\begin{proof}
If~$L\geq m$, then this is Proposition~\ref{p.Holder.regularity}.  We may therefore suppose~$L<m$.  Choose~$s_0,\lambda,\ep$, and the step size~$k$ exactly as in Step~1 of the proof of Proposition~\ref{p.Holder.regularity}.  Define~$\mathcal L_L(\gamma,m)$ by~\eqref{e.L.def.reg}. Here~$M$ and~$N$ are the random variables of Proposition~\ref{p.cutoff.regularity.good.scales} for the cutoff~$L$, not those of Lemmas~\ref{l.sum.the.errors} and~\ref{l.min.scale.good.scale}. Equations~\eqref{e.cutoff.regularity.error.stopping.tail} and~\eqref{e.cutoff.regularity.good.stopping.tail}, with the parameter choices in that proof, give~\eqref{e.cutoff.Holder.minimal.scale}. Equations~\eqref{e.cutoff.regularity.error.stopping} and~\eqref{e.cutoff.regularity.good.stopping} give the two conclusions in~\eqref{e.after.good.scales.applied} for the modified~$\mathbf G$ and~$\mathcal G$.

\smallskip

For~$n\leq j\leq m-5$ and~$z\in3^n\Zd\cap\cu_m$, the calculation leading to~\eqref{e.ratio.of.bs}, with all sums over the fields and all deterministic sums stopped at~$L$, gives
\begin{equation}
	C^{-1}\exp\bigl(-C\lambda(m-n)\bigr)
	\leq
	\frac{
	(\b_{L,j+2})_{z+\cu_{j+2}}
	}{
	(\b_{L,m})_{\cu_m}
	}
	\leq
	C\exp\bigl(C\lambda(m-n)\bigr)\,.
	\label{e.cutoff.Holder.b.ratio}
\end{equation}
Indeed, if~$j+2\geq L$, the numerator and the denominator both equal~$\ahom_L$.  If~$j+2<L<m$, only the fields~$\g_{j+3},\ldots,\g_L$ occur; they are bounded by~\eqref{e.cutoff.regularity.error.stopping}, while~\eqref{e.cutoff.regularity.D.ordering} controls the ratio of effective diffusivities.

\smallskip

We use Lemma~\ref{l.cutoff.regularity.good.scale.estimates} in place of Lemma~\ref{l.excess.decay.good.scales.GMC} in Steps~3--5. With~\eqref{e.cutoff.Holder.b.ratio} in place of~\eqref{e.ratio.of.bs}, the forcing bounds and coefficient comparisons there are unchanged. In Step~6, apply the cutoff-uniform version of Lemma~\ref{l.centered.energy.good.scale} at the nearest good anchor~$n'+2$, as in Step~6 of the proof of Proposition~\ref{p.Holder.regularity}. The volume loss is~$C3^{d(n'-n)/2}$ and the cutoff coefficient ratios are bounded by~$C\exp(C\lambda(m-n))$ by~\eqref{e.cutoff.Holder.b.ratio}. Proposition~\ref{p.cutoff.regularity.good.scales} supplies the same error and stopping bounds. The Campanato estimate at~$n'$ and the detachment estimate at~$m'$ therefore prove the energy estimate. Step~7 uses this unchanged iteration and proves the excess estimate. This proves all three conclusions with~$\mathcal L_L(\gamma,m)$ as the minimal scale.

\end{proof}

\begin{proposition}[Stability under multiplication by a nearly constant function]
\label{p.cutoff.Holder.bounded.multiplier}
There are~$\ep_*(d)>0$,~$c(d)>0$, and~$C(d)<\infty$ such that the following holds if~$\gmcdelta\leq c$. Let~$L\in\N_0$,~$m\in\Z$,~$z\in\Rd$, and~$j\in\N$. Let~$z'+\cu_{m-j}$ be a cube of side~$3^{m-j}$ contained in the concentric cube~$z+\tfrac12\cu_m$ of half the side of~$z+\cu_m$. Outside an event of probability at most
\begin{equation}
	C\exp\left(-\frac{c}{\gmcdelta^2|\log\gmcdelta|^2}\right)\,,
	\label{e.cutoff.Holder.bounded.multiplier.probability}
\end{equation}
the following estimates hold for every positive~$\vartheta\in C(z+\cu_m)$ satisfying
\begin{equation}
	\inf_{b>0}\|b^{-1}\vartheta-1\|_{L^\infty(z+\cu_m)}\leq\ep_*
	\label{e.cutoff.Holder.bounded.multiplier.assumption}
\end{equation}
and every weakly~$\a_L\vartheta$-harmonic function~$h\in H^1(z+\cu_m)$:
\begin{equation*}
	\osc_{z'+\cu_{m-j}}h\leq C3^{-cj}\osc_{z+\cu_m}h\,,
\end{equation*}
and
\begin{equation*}
	\sup_{z'+\cu_{m-j}}|h-(h)_{z+\cu_m}|
	\leq C\|h-(h)_{z+\cu_m}\|_{\underline L^2(z+\cu_m)}\,.
\end{equation*}
\end{proposition}

\begin{proof}
We bound both conclusions by the normalized~$L^2$ oscillation~$\|h-(h)_{z+\cu_m}\|_{\underline L^2(z+\cu_m)}$, combining an estimate below the unit scale with one above it. Below the unit scale, comparison with harmonic functions on small balls bounds point values and H\"older quotients of~$h$ by powers of~$1+\|\nabla\log\a_L\|_{L^\infty(z+\cu_m)}$; no derivative of~$\vartheta$ enters. Above the unit scale, we run the large-scale iteration of the proof of Proposition~\ref{p.cutoff.Holder.regularity} for the coefficient~$\a_L\vartheta$. The main point is that the multiplier changes the coarse-graining errors by at most~$C\ep_*^{\nf12}$ per scale, which is absorbed by running the iteration with exponent~$3/4$ and concluding with exponent~$1/2$. Finally, we bound~$\|\nabla\log\a_L\|_{L^\infty(z+\cu_m)}$ and combine the two estimates; the exceptional events depend on neither~$\vartheta$ nor~$h$.

\smallskip

Multiplying~$\vartheta$ by a positive constant does not change the equation. The infimum in~\eqref{e.cutoff.Holder.bounded.multiplier.assumption} allows a normalization with~$\|\vartheta-1\|_\infty\leq2\ep_*$. We choose~$\ep_*$ small enough for all the estimates below with this bound. Write~$H=\|h-(h)_{z+\cu_m}\|_{\underline L^2(z+\cu_m)}$.

\smallskip

We first record a local estimate that works in every dimension~$d\geq2$. Suppose that a scalar coefficient~$a$ on a Euclidean ball~$V$ satisfies~$\|a/a_0-1\|_\infty\leq\eta$, where~$a_0>0$ and~$\eta$ is a sufficiently small dimensional constant. Let~$v$ be the harmonic function with boundary values~$h$ on~$V$. Testing the equation for~$h-v$, and using the Dirichlet principle for~$v$, gives
\begin{equation*}
	\|\nabla(h-v)\|_{L^2(V)}\leq C\eta\|\nabla h\|_{L^2(V)}\,.
\end{equation*}
The mean-square gradient of~$v$ decreases on concentric balls. Consequently,
\begin{equation*}
	\|\nabla h\|_{\underline L^2(V/2)}
	\leq(1+C2^{d/2}\eta)\|\nabla h\|_{\underline L^2(V)}\,.
\end{equation*}
Choose~$\eta$ so that the factor is at most~$2^{1/2}$. Iteration, Poincar\'e and telescoping ball averages give H\"older exponent~$1/2$ on the inner ball. Caccioppoli on a slightly larger ball controls the initial gradient by the normalized~$L^2$ norm of~$h-c$, for every constant~$c$. This also bounds point values relative to~$c$, and the telescoping averages define a continuous representative of~$h$.

\smallskip

Let~$D=\|\nabla\log\a_L\|_{L^\infty(z+\cu_m)}$. Around every point of~$w+\tfrac12\cu_0$, where~$w+2\cu_0\subseteq z+\cu_m$, choose a ball of radius~$c(d)/(1+D)$. The mean-value theorem and~$\|\vartheta-1\|_\infty\leq2\ep_*$ make~$\a_L\vartheta$, divided by its value of~$\a_L$ at the center, within~$\eta$ of one on twice that ball. The preceding comparison and Caccioppoli estimate, followed by the volume comparison with~$w+\cu_0$, give
\begin{equation*}
	\sup_{w+\tfrac12\cu_0}|h-(h)_{w+\cu_0}|
	\leq C(1+D)^{d/2}\|h-(h)_{w+\cu_0}\|_{\underline L^2(w+\cu_0)}\,.
\end{equation*}
For~$x,y\in w+\tfrac12\cu_0$ the same argument gives
\begin{equation*}
	|h(x)-h(y)|
	\leq C(1+D)^{(d+1)/2}|x-y|^{1/2}
	\|h-(h)_{w+\cu_0}\|_{\underline L^2(w+\cu_0)}\,.
\end{equation*}
For distances larger than~$c(d)/(1+D)$ use the point-value bound and~$\bigl(|x-y|(1+D)/c(d)\bigr)^{\nf12}\geq1$. Thus this estimate covers every distance below one. On a cube of side~$\ell$, the rescaled bounds have~$1+\ell D$ in place of~$1+D$ and~$|x-y|/\ell$ in place of~$|x-y|$. No derivative of~$\vartheta$ is needed.

\smallskip

We next run the interior large-scale iteration for~$\a_L\vartheta$. On every comparison cube~$w+\cu_n\subseteq z+\cu_m$, Lemma~\ref{l.J.sensitivity} and Minkowski's inequality in the weighted subcube sum give
\begin{equation}
	\mathcal E_{s_0,\infty,2}(w+\cu_n;\a_L\vartheta,\a_0)
	\leq(1+C\ep_*)\mathcal E_{s_0,\infty,2}(w+\cu_n;\a_L,\a_0)
	+C\ep_*^{\nf12}\,,\qquad s_0=1/32\,.
	\label{e.cutoff.Holder.bounded.multiplier.error}
\end{equation}
Here~$\a_0$ is the scalar in Lemma~\ref{l.cutoff.regularity.good.scale.estimates}. Cube averages of the long-scale coefficient change by factors between~$1-2\ep_*$ and~$1+2\ep_*$, and their ratios by factors~$1+O(\ep_*)$. The proofs of Lemmas~\ref{l.centered.energy.good.scale},~\ref{l.harmonic.approximation.good.scales.GMC} and~\ref{l.excess.decay.good.scales.GMC} use the good event only to bound the anchor error and the resulting ellipticity ratios. Their interior estimates therefore apply to every positive scalar coefficient with these bounds; for~$\a_L\vartheta$ they are supplied by~\eqref{e.cutoff.Holder.bounded.multiplier.error}. Apply the deterministic interior comparison of Proposition~\ref{p.general.coarse.graining} to the product coefficient; its Caccioppoli factors are comparable to the original ones. In that comparison, therefore, the excess contraction and slope-error coefficient change by at most~$C\ep_*^{\nf12}$. Fix the block length and other parameters of the interior proof of Proposition~\ref{p.cutoff.Holder.regularity} with exponent~$3/4$, before choosing~$\ep_*$. The perturbation preserves a strict contraction on good blocks and adds at most~$C\ep_*^{\nf12}$ times the number of scales to the error sum. More explicitly, on an interval~$[a,b]$ the new error sum is at most
\begin{equation*}
	(1+C\ep_*^{\nf12})\sum_{j=a}^b\ep_j
	+C\ep_*^{\nf12}(b-a+1)\,.
\end{equation*}
Choose~$\ep_*$ so small that this additional exponential loss uses less than the margin between exponents~$3/4$ and~$1/2$. Lemma~\ref{l.iteration.lemma.GMC} then gives exponent~$1/2$ on the same good-scale event, uniformly in~$\vartheta$. All comparisons are interior; no boundary estimate is used.

\smallskip

Here is the covering needed for that iteration. In each coordinate use the four translates with shifts~$0,3^m/4,3^m/2,3^{m+1}/4$. At every scale~$3^n\leq3^m$, their relative phases are the four quarter periods, since powers of three permute these phases. A cube of sufficiently smaller side than~$3^n$ is therefore contained in a cube of one of the~$4^d$ grids at that scale, at distance at least a fixed fraction of~$3^n$ from its boundary. Fix a dimensional integer~$a\geq4$. For each grid retain the top cubes of side~$3^{m-a}$ meeting~$z+\tfrac12\cu_m$ whose comparison neighborhoods lie in~$z+\cu_m$; their inner halves cover~$z+\tfrac12\cu_m$, and every point of~$z+\tfrac12\cu_m$ lies in one of them, at distance at least a fixed fraction of~$3^{m-a}$ from its boundary. There are only dimensionally many such top cubes. Every smaller comparison cube used by the iteration stays in~$z+\cu_m$. In particular, unit cubes~$w+\cu_0$ with~$w+2\cu_0\subseteq z+\cu_m$ and their inner halves cover~$z+\tfrac12\cu_m$ when~$m$ is sufficiently large.

\smallskip

Choose a fixed~$j_0$ larger than the additive stopping constant. On each retained top cube require the minimal depth for exponent~$3/4$ to be at most~$j_0$. The tail~\eqref{e.cutoff.Holder.minimal.scale} and the finite union over the top cubes bound the failure probability by~\eqref{e.cutoff.Holder.bounded.multiplier.probability}, after changing dimensional constants. The perturbed iteration, volume comparison at the top cubes, and telescoping cube averages give
\begin{equation*}
	\|h-(h)_{w+\ell\cu_0}\|_{\underline L^2(w+\ell\cu_0)}
	\leq C(3^{-m}\ell)^{\nf12}H
\end{equation*}
on the interior cubes~$w+\ell\cu_0$ of side~$1\leq\ell\leq3^{m-a}$ used in the covering. The finitely many initial scales and the at most~$j_0$ missing scales contribute only a dimensional constant. Differences of their averages obey the same Campanato bound, and their averages differ from~$(h)_{z+\cu_m}$ by at most~$CH$.

\smallskip

The local point-value factors are also controlled on an event independent of~$\vartheta$. Covering~$z+\cu_m$ at scale~$3^k$, Assumption~\ref{a.g2} and scaling give a sub-Gaussian bound for~$\|\nabla\g_k\|_{L^\infty(z+\cu_m)}$ with mean at most~$C\gmcdelta3^{-k}\sqrt{(m-k)_++1}$ and tail scale~$C\gmcdelta3^{-k}$. Summing these bounds gives
\begin{equation*}
	\P[D>C\gmcdelta\sqrt{m_++1}+t]
	\leq C\exp\left(-\frac{ct^2}{\gmcdelta^2}\right)\,.
\end{equation*}
Take~$t=|\log\gmcdelta|^{-1}$. Its failure probability has the form~\eqref{e.cutoff.Holder.bounded.multiplier.probability}; on the complementary event,~$(1+D)^{(d+1)/2}\leq C3^{m/8}$ for~$m\geq0$. Thus the local bound on unit cubes and the Campanato estimate bound~$|h-(h)_{z+\cu_m}|$ on~$z+\tfrac12\cu_m$ by~$CH$. For~$z'+\cu_{m-j}$ of side~$3^{m-j}\geq1$, comparisons of averages and the unit-scale point errors give
\begin{equation*}
	\osc_{z'+\cu_{m-j}}h\leq C\bigl(3^{-j/2}+3^{-3m/8}\bigr)H
	\leq C3^{-cj}H\,,
\end{equation*}
since~$j\leq m$. A target cube of side below a fixed dimensional fraction of one lies in a single inner half of a unit cube from the covering. Larger targets of side below one are covered by dimensionally many overlapping inner halves; their point errors and the differences of the unit-cube averages obey the same bound because the side~$3^{m-j}$ is then bounded below. If~$3^{m-j}<1$, the local H\"older estimate therefore gives
\begin{equation*}
	\osc_{z'+\cu_{m-j}}h\leq C3^{-3m/8}3^{\frac12(m-j)}H
	\leq C3^{-cj}H\,.
\end{equation*}
These estimates hold for every fixed~$c\leq3/8$, and we then decrease~$c$ to match the probability bound. The finitely many small values of~$j$ follow from~$\osc_{z'+\cu_{m-j}}h\leq\osc_{z+\cu_m}h$ after increasing~$C$.

\smallskip

For~$m\leq a+j_0+6$, including all~$m<0$, use the rescaled local estimates directly on interior cubes of side comparable to~$3^m$. On the same gradient event,~$3^mD$ is dimensionally bounded, so both their point-value constant and their H\"older constant are bounded. This gives the same two conclusions. Finally~$H\leq\osc_{z+\cu_m}h$ whenever that oscillation is finite, and otherwise its asserted upper bound is automatic. Taking the union of the two exceptional events completes the proof. Neither event depends on~$\vartheta$ or~$h$.
\end{proof}

\subsection{Dirichlet homogenization}

The estimates above compare a solution with an affine function.  We now compare it instead with the solution of the homogenized Dirichlet problem in the unit cube, after rescaling the cutoff coefficient~$\a_L$ from a cube at least as large as the cutoff scale.

\begin{theorem}[Quantitative homogenization of the Dirichlet problem]
\label{t.cutoff.Dirichlet.homogenization}
For~$L,M\in\N_0$ with~$L\leq M$, recall the rescaled coefficient of the introduction,
\begin{equation}
	A_{L,M}(x)
	\coloneqq
	\ahom_L^{-1}\a_L(3^Mx)\,,
	\qquad x\in \cu_0\,.
	\label{e.rescaled.cutoff.coefficient.Dirichlet}
\end{equation}
For every~$\vartheta\in(0,1)$ and~$q\in[1,\infty)$, there are~$\gmcdelta_0(\vartheta,q,d)>0$ and~$C(\vartheta,q,d)<\infty$ such that, if~$\gmcdelta\leq\gmcdelta_0$, then there are random variables~$\mathcal Z_{L,M}\geq1$,~$0\leq L\leq M$, satisfying
\begin{equation}
	\sup_{M\in\N_0}\sup_{0\leq L\leq M}
	\|\mathcal Z_{L,M}\|_{L^q(\P)}
	\leq C\,.
	\label{e.Dirichlet.random.factor.moment}
\end{equation}
For~$f\in L^2(\cu_0)$ and~$h\in H^2(\cu_0)$, recall that~$u_{L,M}$ and~$u_{\mathrm{hom}}$ are the unique weak solutions in~$H^1(\cu_0)$ of
\begin{equation}
	\left\{
	\begin{aligned}
		-\nabla\cdot(A_{L,M}\nabla u_{L,M})&=f
		&&\text{in }\cu_0\,,\\
		u_{L,M}&=h
		&&\text{on }\partial \cu_0\,,
	\end{aligned}
	\right.
	\qquad
	\left\{
	\begin{aligned}
		-\Delta u_{\mathrm{hom}}&=f
		&&\text{in }\cu_0\,,\\
		u_{\mathrm{hom}}&=h
		&&\text{on }\partial \cu_0\,,
	\end{aligned}
	\right.
	\,.
	\label{e.cutoff.Dirichlet.problems}
\end{equation}
Then, almost surely, for every such~$f$ and~$h$, the homogenization error~$\mathcal D_{L,M}$ of the introduction satisfies
\begin{equation}
	\mathcal D_{L,M}
	\leq
	\mathcal Z_{L,M}\gmcdelta^\vartheta
	\left(
	\|f\|_{L^2(\cu_0)}+\|h\|_{H^2(\cu_0)}
	\right)\,.
	\label{e.cutoff.Dirichlet.homogenization}
\end{equation}
If~$f=0$, then the factor~$\gmcdelta^\vartheta$ on the right side may be replaced by~$C\gmcdelta$, with the same random variable~$\mathcal Z_{L,M}$.
\end{theorem}

\begin{proof}
We rescale the two Dirichlet problems to the cube~$\cu_M$ and apply the general coarse-graining estimate of Proposition~\ref{p.general.coarse.graining} at the scale~$M-k$. It bounds the errors of the gradient and of the flux in~$\Hminusuls{-s}(\cu_0)$ by~$3^{s_1k}$ times coarse-graining errors whose moments are of order~$\gmcdelta$, plus a term of order~$3^{-s_2k}$ coming from the source~$f$. The main step is the choice of~$k$ that balances these two terms; by~\eqref{e.Dirichlet.parameter.balance}, the resulting power of~$\gmcdelta$ is at most~$\gmcdelta^\vartheta$. A duality argument with the Dirichlet Laplacian then gives the~$L^2$ estimate. These fractional-dual estimates use the unrestricted tests of~\eqref{e.hatted.negative.norm.convention}.

\smallskip

Fix~$s_2\in(0,1)$, and then choose~$0<s_1<s<\min\{s_2,1/2\}$ so that
\begin{equation}
	\frac{s_2}{s_1+s_2}>\vartheta\,.
	\label{e.Dirichlet.parameter.balance}
\end{equation}
Fix a finite exponent~$r\geq\max\{4q,8d/s_1\}$ so large that~$(s_1/2)\log(2+r)\geq c_P$, where~$c_P$ is the logarithmic-condition constant in Proposition~\ref{p.coarse.grained.bound}. We use input moments through~$4r$, leaving room for squares and products before taking~$L^r$ norms. These orders depend on~$\vartheta$ through~$s_1$, as well as on~$q$ and~$d$, and are fixed before the disorder threshold.

\smallskip

To check applicability explicitly, let~$\xi=4r$. Let~$C_P$ be the constant in the first conclusion of Proposition~\ref{p.coarse.grained.bound}, and~$c_L$ the moment-range constant in Lemma~\ref{l.multiscale.response.large.cubes}. Fix
\begin{equation*}
	K\geq\max\left\{1,C_P\xi\log(2+\xi),\frac{2\xi}{c_Ls_1}\right\}\,,
	\qquad \delta_1=K\gmcdelta^2\,.
\end{equation*}
Decrease~$\gmcdelta_0(\vartheta,q,d)\leq e^{-1}$ so that~$C_P\xi\gmcdelta_0\leq1$ and~$K\gmcdelta_0^2<\min\{1/2,c_Ps_1/2\}$. Then~$\xi\leq C_P^{-1}\gmcdelta^{-2}|\log\gmcdelta|^{-1}$, using~$\gmcdelta^2|\log\gmcdelta|\leq\gmcdelta$. The proposition supplies~$\mathcal S(\infty,\xi,\delta_1)$, while~$\gmcdelta^2\leq\delta_1<1$. At both orders~$s_1$ and~$s_1/2$, the large-cube lemma has~$4d/(s_1/2)\leq\xi$ and~$\xi\leq c_L(s_1/2)\gmcdelta^{-2}\delta_1$ by the choice of~$K$. Its bound is~$C\sqrt{\delta_1}=C\gmcdelta$. The lower-order estimates follow by monotonicity of probability norms. Thus, uniformly in~$0\leq L\leq M$, these results give
\begin{equation}
	\left\|
	\mathcal E_{s_1,\infty,1}
	(\cu_M;\a_L,\ahom_L)
	\right\|_{L^r(\P)}
	+
	\left\|
	\mathcal E_{\nf{s_1}{2},\infty,2}
	(\cu_M;\a_L,\ahom_L)
	\right\|_{L^r(\P)}
	\leq C\gmcdelta\,.
	\label{e.Dirichlet.full.response.moments}
\end{equation}

\smallskip

Write~$f$ as the divergence of a field in~$H^1(\cu_0;\Rd)$ of norm at most~$C\|f\|_{L^2(\cu_0)}$, for instance minus the gradient of the solution of the Dirichlet problem for the Laplacian in~$\cu_0$ with source~$f$. Rescale to~$\cu_M$ by setting~$U(y)\coloneqq u_{L,M}(3^{-M}y)$ and~$V(y)\coloneqq u_{\mathrm{hom}}(3^{-M}y)$, and let~$G(y)$ be~$\ahom_L3^{-M}$ times the field above evaluated at~$3^{-M}y$. Then
\begin{equation*}
	-\nabla\cdot(\a_L\nabla U)=\nabla\cdot G\,,
	\qquad
	-\ahom_L\Delta V=\nabla\cdot G\,,
	\qquad
	U-V\in H_0^1(\cu_M)\,.
\end{equation*}

\smallskip

Let~$k\in\N$ and apply Proposition~\ref{p.general.coarse.graining}, with~$p=2$ and~$n=M-k$. Definition~\ref{d.mathcal.E} gives
\begin{equation*}
	\mathcal E_{s_1,\infty,1}
	(\cu_M,M-k;\a_L,\ahom_L)
	\leq
	3^{s_1k}
	\mathcal E_{s_1,\infty,1}
	(\cu_M;\a_L,\ahom_L)\,,
\end{equation*}
and
\begin{equation*}
	\mathcal E_{\nf{s_1}{2},\infty,2}
	(\cu_M,M-k;\a_L,\ahom_L)^2
	\leq
	3^{s_1k}
	\mathcal E_{\nf{s_1}{2},\infty,2}
	(\cu_M;\a_L,\ahom_L)^2\,.
\end{equation*}
The geometric weights in the sum over~$k$ in~\eqref{e.general.coarse.graining.estimate.ASD}, together with~\eqref{e.average.normalized.norm.partition}, bound that sum by~$C\|\a_L^{\nf12}\nabla U\|_{\underline L^2(\cu_M)}$. Moreover,~\eqref{e.cg.RHS} and~\eqref{e.bound.Lambdas.by.Es}, using the fixed input orders chosen above, give a random variable~$Y_{L,M}\geq1$ with~$\|Y_{L,M}\|_{L^r(\P)}\leq C(\vartheta,q,d)$ uniformly in~$L\leq M$, such that
\begin{equation}
	\|A_{L,M}^{\nf12}\nabla u_{L,M}\|_{L^2(\cu_0)}
	\leq
	Y_{L,M}
	\left(
	\|f\|_{L^2(\cu_0)}+\|h\|_{H^2(\cu_0)}
	\right)\,.
	\label{e.Dirichlet.energy.factor}
\end{equation}
Rescaling~\eqref{e.general.coarse.graining.estimate.ASD} back to~$\cu_0$ therefore yields
\begin{align}
	\lefteqn{
	\|\nabla u_{L,M}-\nabla u_{\mathrm{hom}}\|_{\Hminusuls{-s}(\cu_0)}
	+
	\|A_{L,M}\nabla u_{L,M}
	-\nabla u_{\mathrm{hom}}\|_{\Hminusuls{-s}(\cu_0)}
	} \qquad & \notag\\
	&\leq
	C\left(
	3^{s_1k}\mathcal E_1Y_{L,M}
	+
	3^{-s_2k}
	\bigl(1+3^{s_1k}\mathcal E_2^2\bigr)
	\right)
	\left(
	\|f\|_{L^2(\cu_0)}+\|h\|_{H^2(\cu_0)}
	\right)\,,
	\label{e.Dirichlet.prebalance}
\end{align}
where~$\mathcal E_1$ and~$\mathcal E_2$ denote the two random variables in~\eqref{e.Dirichlet.full.response.moments}.

\smallskip

Take~$k\coloneqq\lceil\log_3(\gmcdelta^{-1})/(s_1+s_2)\rceil$. With~$\mathcal E_1,\mathcal E_2$, and~$Y_{L,M}$ as above, define
\begin{equation*}
	\mathcal Z_{L,M}
	\coloneqq 1
	+C\gmcdelta^{-\vartheta}
	\left(
	3^{s_1k}\mathcal E_1Y_{L,M}
	+3^{-s_2k}\bigl(1+3^{s_1k}\mathcal E_2^2\bigr)
	\right)
	+C\gmcdelta^{-1}\mathcal E_1Y_{L,M}\,.
\end{equation*}
Equations~\eqref{e.Dirichlet.parameter.balance},~\eqref{e.Dirichlet.full.response.moments}, and~\eqref{e.Dirichlet.energy.factor}, followed by H\"older's inequality, give~\eqref{e.Dirichlet.random.factor.moment}.  Substituting this choice of~$k$ and the definition of~$\mathcal Z_{L,M}$ into~\eqref{e.Dirichlet.prebalance} proves~\eqref{e.cutoff.Dirichlet.homogenization} for the two~$H^{-1}$ terms of~$\mathcal D_{L,M}$: restriction to zero-boundary tests and the embedding~$H^1_0(\cu_0)\hookrightarrow H^s(\cu_0)$ bound each coordinate dual, and the finite sum over coordinates changes only~$C(d)$.  If~$f=0$, we apply~\eqref{e.Dirichlet.prebalance} again, with~$k=1$.  The~$3^{-s_2k}$ term is then absent, and the last term in the definition of~$\mathcal Z_{L,M}$ gives the factor~$C\gmcdelta$.

\smallskip

Finally, the difference~$u_{L,M}-u_{\mathrm{hom}}$ belongs to~$H_0^1(\cu_0)$. For every~$w\in H_0^1(\cu_0)$, duality and the spectral theorem for the Dirichlet Laplacian on~$\cu_0$ give~$\|w\|_{L^2(\cu_0)}\leq C\|\nabla w\|_{\Hminusuls{-s}(\cu_0)}$: test with the gradient of the image of~$w$ under the inverse Dirichlet Laplacian, which is bounded in~$H^s(\cu_0)$ by~$C\|w\|_{L^2(\cu_0)}$. Applying the unrestricted fractional-dual estimate~\eqref{e.Dirichlet.prebalance} to the difference gives the~$L^2$ term and completes the proof.
\end{proof}

\begin{theorem}[Quantitative homogenization at a fixed cutoff]
\label{t.fixed.cutoff.Dirichlet.algebraic}
There exists~$\alpha_{\mathrm{hom}}(d)>0$ with the following property. Fix~$L\in\N_0$. For~$M\in\N_0$ with~$M\geq L$, let~$A_{L,M}$ be as in~\eqref{e.rescaled.cutoff.coefficient.Dirichlet}. For every~$q\in[1,\infty)$, there is a random variable~$\mathcal Y_{L,q}\geq1$ satisfying
\begin{equation}
	\|\mathcal Y_{L,q}\|_{L^q(\P)}
	\leq C(L,q,d,\gmcdelta)
	\label{e.fixed.cutoff.Dirichlet.random.factor}
\end{equation}
such that, almost surely, for every~$M\geq L$,~$f\in L^2(\cu_0)$ and~$h\in H^2(\cu_0)$, the homogenization error~$\mathcal D_{L,M}$ of the unique weak solutions~$u_{L,M}$ and~$u_{\mathrm{hom}}$ of~\eqref{e.cutoff.Dirichlet.problems} satisfies
\begin{equation}
	\mathcal D_{L,M}
	\leq
	\mathcal Y_{L,q}\,
	3^{-\alpha_{\mathrm{hom}}(M-L)}
	\left(
	\|f\|_{L^2(\cu_0)}+\|h\|_{H^2(\cu_0)}
	\right)\,.
	\label{e.fixed.cutoff.Dirichlet.algebraic}
\end{equation}
\end{theorem}

\begin{proof}
Consider the coefficient~$A_{L,L}$. For fixed~$L$, this is a stationary positive scalar coefficient with range of dependence at most~$\sqrt d$, and its homogenized matrix is the identity.  It satisfies the logarithmic coefficient bound
\begin{equation}
	\|\log A_{L,L}\|_{L^\infty(\cu_0)}
	\leq\O_{\Gamma_2}\bigl(C\gmcdelta(L+1)\bigr)\,.
	\label{e.fixed.cutoff.logarithmic.input}
\end{equation}
Indeed, this follows from~\eqref{e.infrared.approx.cutoffs.field}, Assumption~\ref{a.g2}, and the bound~$|\log \ahom_L|\leq C\gmcdelta^2(L+1)$ from Theorem~\ref{t.renormalized.diffusivities}.

\smallskip

We apply the results of~\citet[Theorem~B, pp.~9--10; Theorem~4.1, pp.~73--74; and Proposition~D.3, pp.~145--146]{AK.HC} to~$A_{L,L}$.  Stationarity and finite range give the mixing assumption of that paper with~$\beta=0$ and~$\nu=d/2$, and with Gaussian concentration.  Proposition~D.3 of that paper, applied to~\eqref{e.fixed.cutoff.logarithmic.input} with exponent~$\nf14$, gives the coarse-grained ellipticity assumption of that paper above a random scale. The probability that this scale exceeds~$\sqrt d\,t$ is at most~$\exp(C_L-c_L\log^2t)$, with positive constants~$C_L$ and~$c_L$ that may depend on~$L$; this tail bound implies that all polynomial moments of this scale are finite.

\smallskip

\citet[Corollary~4.3, p.~76]{AK.HC}, with exponent~$5/16$, gives an estimate uniform over all subcubes. Together with Definition~\ref{d.mathcal.E}, it yields a number~$\alpha_{\mathrm{hom}}(d)>0$ and, for every fixed~$q<\infty$, a random variable~$\mathcal W_{L,q}$ with finite~$q$th moment such that, for all~$N\in\N_0$,~$r\in\{1,2\}$, and~$s\in\{\nf{3}{16},\nf38\}$,
\begin{equation}
	\mathcal E_{s,\infty,r}(\cu_N;A_{L,L},1)
	\leq
	\mathcal W_{L,q}\,3^{-2\alpha_{\mathrm{hom}}N}\,.
	\label{e.fixed.cutoff.response.algebraic}
\end{equation}
Indeed, by the estimate of~\citet[Corollary~4.3, p.~76]{AK.HC} for subcubes, the quantity~$\max_{|e|=1}J(\,\cdot\,,e,e\,;A_{L,L})$ on a subcube at scale~$k\leq N$ is bounded by a random multiple of~$3^{-4\alpha_{\mathrm{hom}}N}3^{\frac{5}{16}(N-k)}$. After taking square roots, the weighted sum over scales in Definition~\ref{d.mathcal.E} converges, since~$s>\nf5{32}$ for both values of~$s$ above. Below that random scale, we use~\eqref{e.ord.bounds.truncated} and the extrema of~$A_{L,L}^{\pm1}$ instead. By a union bound over unit cubes and~\eqref{e.fixed.cutoff.logarithmic.input}, the moments of these extrema on~$\cu_N$ grow at most like~$C(L,q)\exp(C(L,q)\sqrt{N+1})$. H\"older's inequality, the Gaussian-type tail bound for that scale, and the analogous tail for the scale in~\citet[Corollary~4.3]{AK.HC} therefore allow us to absorb these scales into~$\mathcal W_{L,q}$ while keeping its moments finite.

\smallskip

We use the estimate leading to~\eqref{e.Dirichlet.prebalance}, after rescaling~$\cu_M$ first by~$3^L$ and then by~$3^{M-L}$.  Apply Proposition~\ref{p.general.coarse.graining} with~$p=2$,~$s_1=\nf38$,~$s=\nf7{16}$, and~$s_2=\nf34$, so that~$sp^{\prime}=7/8<1$ and both~$s_1$ and~$s_1/2$ are covered by~\eqref{e.fixed.cutoff.response.algebraic}, and choose~$k\coloneqq\max\{1,\lfloor2\alpha_{\mathrm{hom}}(M-L)/(s_1+s_2)\rfloor\}$. In~\eqref{e.Dirichlet.prebalance}, we use~\eqref{e.fixed.cutoff.response.algebraic}, applied at the scale~$M-L$, in place of~\eqref{e.Dirichlet.full.response.moments}, and we construct~$Y_{L,M}$ in~\eqref{e.Dirichlet.energy.factor} from~\eqref{e.cg.RHS},~\eqref{e.bound.Lambdas.by.Es}, and now~\eqref{e.fixed.cutoff.response.algebraic}. These give~$\sup_{M\geq L}\|Y_{L,M}\|_{L^r}<\infty$ for every fixed~$r$, with constants allowed to depend on~$L,r,d,\gmcdelta$. Substituting~\eqref{e.fixed.cutoff.response.algebraic} and the choice of~$k$ in~\eqref{e.Dirichlet.prebalance}, and using the bound~$Y_{L,M}\leq3^{\alpha_{\mathrm{hom}}(M-L)/25}\sup_{M'\geq L}3^{-\alpha_{\mathrm{hom}}(M'-L)/25}Y_{L,M'}$, we obtain the negative-Sobolev part of~\eqref{e.fixed.cutoff.Dirichlet.algebraic} with the random factor
\begin{equation*}
	\mathcal Y_{L,q}
	\coloneqq
	C(1+\mathcal W_{L,4q})^2\Bigl(1+\sup_{M\geq L}3^{-\alpha_{\mathrm{hom}}(M-L)/25}\,Y_{L,M}\Bigr)\,.
\end{equation*}
Since~$Y_{L,M}$ is constructed with a moment larger than~$4q$, the supremum has a finite~$4q$th moment, because its~$4q$th power is at most the sum over~$M\geq L$ of~$3^{-4q\alpha_{\mathrm{hom}}(M-L)/25}Y_{L,M}^{4q}$. H\"older's inequality shows that~$\mathcal Y_{L,q}$ has a finite~$q$th moment, which gives~\eqref{e.fixed.cutoff.Dirichlet.random.factor}. As in the proof of Theorem~\ref{t.cutoff.Dirichlet.homogenization}, restriction of the unrestricted fractional-dual estimate to zero-boundary tests gives the two negative-norm terms of~$\mathcal D_{L,M}$, while duality with the Dirichlet Laplacian gives the~$L^2$ term.
\end{proof}

\subsection{Limiting coefficient and regularity}

The limiting coefficient field is the almost-sure local limit
\begin{equation}
	\widetilde\a_L(x)
	\coloneqq
	\frac{\a_L(x)}{\a_L(0)}\,,
	\qquad
	\a(x)
	\coloneqq
	\lim_{L\to\infty}\widetilde\a_L(x)
	=
	\exp\Bigl(\sum_{k=0}^\infty(\g_k(x)-\g_k(0))\Bigr)\,.
	\label{e.fixed.clock.coefficient}
\end{equation}
Dividing by~$\a_L(0)$ is necessary, since the strong law gives~$(L+1)^{-1}\log\a_L(0)\to-\tau^2$.  This division does not affect the generator of~$X^{(L)}$, and the reversible measures~$(\a_L/\a_L(0))\,dx$ converge locally to~$\a\,dx$. By Assumption~\ref{a.g2} and scaling, the coefficient ratios converge locally in~$C^{1,1}$. The following lemma gives bounds uniform in~$L$.

\begin{lemma}[The coefficients divided by their value at the origin]
	\label{l.finite.cutoff.coefficient.convergence}
	There is~$\gmcdelta_0(d)>0$ such that, if~$\gmcdelta\leq\gmcdelta_0(d)$, then there are~$\kappa\in(0,1)$ and an almost surely finite random variable~$C_\omega$ such that, almost surely, for every~$L\in\N_0$ and~$x\in\Rd$,
	\begin{equation}
		\sup_{\b\in\{\a\}\cup\{\widetilde\a_L\}_{L\in\N_0}}
		\left(\b(x)+\b(x)^{-1}+|\nabla\b(x)|\right)
		\leq C_\omega(1+|x|)^\kappa\,.
		\label{e.finite.cutoff.coefficient.polynomial.growth}
	\end{equation}
	Moreover,
	\begin{equation}
		\widetilde\a_L\longrightarrow\a
		\quad\text{locally in }C^{1,1}\,,
		\qquad
		\nabla\log\widetilde\a_L
		\longrightarrow\nabla\log\a
		\quad\text{locally in }C^{0,1}\,,
		\label{e.finite.cutoff.coefficient.C11.convergence}
	\end{equation}
	and
	\begin{equation}
		\sup_{L\in\N_0}
		\left(
		|\nabla\log\widetilde\a_L(x)|
		+|\nabla^2\log\widetilde\a_L(x)|
		\right)
		\leq
		C_\omega\left(1+\sqrt{\log(2+|x|)}\right)\,.
		\label{e.finite.cutoff.log.derivative.growth}
	\end{equation}
\end{lemma}

\begin{proof}
	Fix a sufficiently large dimensional constant~$C$, and cover each ball~$B_{C2^n}$ by cubes of side~$3^k$. Stationarity,~\eqref{e.same.law.shift},~\eqref{e.coefficient.field.regularity}, and the Borel--Cantelli lemma give a finite random variable~$C_\omega$ such that, for all~$n,k\in\N_0$ and~$r\in\{1,2\}$,
	\begin{equation}
		\|\nabla^r\g_k\|_{L^\infty(B_{C2^n})}
		\leq
		C_\omega3^{-rk}\sqrt{n+k+1}\,.
		\label{e.finite.cutoff.derivative.Borel.Cantelli}
	\end{equation}
	Summing over~$k$ proves the local convergence~\eqref{e.finite.cutoff.coefficient.C11.convergence} and the bound~\eqref{e.finite.cutoff.log.derivative.growth}.

\smallskip

    It remains to obtain the bound on the coefficient itself uniformly in the cutoff~$L$. Let~$K_n=\max\{k\in\N_0:3^k\leq C2^n\}$, so~$K_n+1\leq C(d)(n+1)$. For a fixed lattice point~$z$, the independence of the fields and Assumptions~\ref{a.g1}--\ref{a.g2} give centered, independent sub-Gaussian variables~$\g_k(z)$. Hence
\begin{equation*}
	\E\exp\left(t\sum_{k=0}^j\g_k(z)\right)
	\leq\exp(Ct^2\gmcdelta^2(j+1))\,,
	\qquad
	\P\left[\left|\sum_{k=0}^j\g_k(z)\right|>u\right]
	\leq2\exp\left(-\frac{u^2}{C\gmcdelta^2(j+1)}\right)\,.
\end{equation*}
The second bound follows by optimizing the exponential Markov bound for both signs. There are at most~$C2^{dn}(n+1)$ pairs~$(z,j)$ with~$z\in\mathbb Z^d\cap B_{C2^n}$ and~$0\leq j\leq K_n$. With~$u=A\gmcdelta(n+1)$ their union therefore has probability at most~$C(n+1)2^{dn}\exp(-cA^2(n+1))$. Fix a dimensional~$A$ so large that this is summable in~$n$. Borel--Cantelli bounds every such unanchored sum by~$A\gmcdelta(n+1)+C_\omega$. Since~$0$ is one of the lattice points, subtracting its sum gives the anchored bound with twice this constant. Absorbing the finitely many initial rows and the additive deterministic constant into~$C_\omega$ yields
	\begin{equation}
		\max_{z\in\mathbb Z^d\cap B_{C2^n},\ 0\leq j\leq K_n}
		\left|
		\sum_{k=0}^j(\g_k(z)-\g_k(0))
		\right|
		\leq
		C\gmcdelta n+C_\omega\,.
		\label{e.finite.cutoff.value.maximal.bound}
	\end{equation}
	For~$x\in B_{2^n}$, we apply~\eqref{e.finite.cutoff.value.maximal.bound} at a nearest lattice point and control the difference using~\eqref{e.finite.cutoff.derivative.Borel.Cantelli}; for the scales~$k$ with~$3^k\leq C2^n$, this difference is bounded by~$C_\omega\sqrt{n+1}$. For the other scales, the mean-value theorem along the segment from~$0$ to~$x$ gives
	\begin{equation*}
		\sum_{3^k>C2^n}|\g_k(x)-\g_k(0)|
		\leq
		C_\omega2^n
		\sum_{3^k>C2^n}3^{-k}\sqrt{n+k+1}
		\leq
		C_\omega\sqrt{n+1}\,.
	\end{equation*}
	Both estimates hold for every partial sum of the series. Consequently,
	\begin{equation*}
		\sup_{L\in\N_0}\sup_{x\in B_{2^n}}
		|\log\widetilde\a_L(x)|
		\leq
		C\gmcdelta n+C_\omega\sqrt{n+1}+C_\omega\,.
	\end{equation*}
	Choosing~$\gmcdelta_0(d)$ small, bounding the square-root term by an arbitrarily small multiple of~$n$ plus a random constant, and exponentiating, we obtain~\eqref{e.finite.cutoff.coefficient.polynomial.growth} for~$\widetilde\a_L$ and its reciprocal.  Since~$\nabla\widetilde\a_L=\widetilde\a_L\nabla\log\widetilde\a_L$, the bound~\eqref{e.finite.cutoff.log.derivative.growth} gives the gradient bound after a final arbitrarily small increase of~$\kappa$.  Passing to the local limit proves the assertions for~$\a$.
\end{proof}

Recall that, for a large~$C_0(d)<\infty$ and~$\gmcdelta$ small, the regularity exponent~$\gamma_{\mathrm{reg}}\coloneqq1-C_0\gmcdelta|\log\gmcdelta|^{\nf12}$ of the introduction lies in~$(\nf12,1)$.

\begin{theorem}[Large-scale H\"older regularity and a Liouville theorem]
\label{t.large.scale.Holder.multifractal}
There are~$\gmcdelta_0(d)>0$ and~$C(d)<\infty$ such that the following holds if~$\gmcdelta\leq\gmcdelta_0(d)$. For~$L=\infty$, the coefficient~$\a_L$ stands for the limiting coefficient~$\a$ of~\eqref{e.fixed.clock.coefficient}.
\begin{itemize}
\item \underline{\emph{Large-scale H\"older estimate}}: Let~$\gamma\in[\nf12,\gamma_{\mathrm{reg}}]$. For every~$L\in\N_0\cup\{\infty\}$ and~$m\in\N$, there is an~$\N$-valued random variable~$\mathcal L_L(\gamma,m)$ satisfying
\begin{equation}
	\P\left[\mathcal L_L(\gamma,m)>k\right]
	\leq
	C\exp\left(
	-\frac{(1-\gamma)^2(k-C)_+}
	{C\gmcdelta^2|\log\gmcdelta|}
	\right)
	,	\qquad k\in\N\,.
	\label{e.limiting.Holder.minimal.scale}
\end{equation}
Moreover, almost surely, for every~$L$ and~$m$, if~$u\in H^1(\cu_m)$ is~$\a_L$-harmonic in~$\cu_m$, then, for every~$n\in\N_0$ satisfying~$n\leq m-\mathcal L_L(\gamma,m)$ and every~$z\in3^n\Zd$ such that~$z+\cu_n\subseteq\cu_{m-1}$,
\begin{equation}
	\left\|
	u-(u)_{z+\cu_n}
	\right\|_{\underline L^2(z+\cu_n)}
	\leq
	C3^{-\gamma(m-n)}
	\left\|
	u-(u)_{\cu_m}
	\right\|_{\underline L^2(\cu_m)}
	\label{e.large.scale.Holder.multifractal}
\end{equation}
and
\begin{equation}
	\left\|\a_L^{\nf12}\nabla u\right\|_{\underline L^2(z+\cu_n)}
	\leq
	C3^{(1-\gamma)(m-n)}
	\left\|\a_L^{\nf12}\nabla u\right\|_{\underline L^2(\cu_m)}\,.
	\label{e.large.scale.energy.multifractal}
\end{equation}

\item \underline{\emph{Liouville theorem}}: Almost surely, for every~$L\in\N_0\cup\{\infty\}$, if~$u\in H^1_{\mathrm{loc}}(\Rd)$ satisfies
\begin{equation}
	-\nabla\cdot(\a_L\nabla u)=0
	\quad\text{in }\Rd\,,
	\qquad
	\liminf_{R\to\infty}
	R^{-\gamma_{\mathrm{reg}}}
	\inf_{c\in\mathbb R}
	\|u-c\|_{\underline L^2(B_R)}=0\,,
	\label{e.multifractal.Liouville.power.growth}
\end{equation}
then~$u$ is constant.
\end{itemize}
\end{theorem}

\begin{proof}
Fix~$\gamma\in[\nf12,\gamma_{\mathrm{reg}}]$. We take~$C_0$ at least as large as the constant~$C$ of Proposition~\ref{p.cutoff.Holder.regularity}, and~$\gmcdelta_0\leq C^{-1}$; then that proposition applies with exponent~$\gamma$. For finite~$L$, the two estimates are special cases of Proposition~\ref{p.cutoff.Holder.regularity}, for~$\a_L$-harmonic functions and interior cubes: we take~$\g=0$,~$x=z$, and~$j=n$ in~\eqref{e.Holder.estimate.boxes.local}, and use~\eqref{e.energy.density.estimate}. Approximating the trace in~$H^{\nf{1}{2}}(\partial\cu_m)$ by smooth functions extends the conclusion from smooth to arbitrary~$H^1$ boundary data.

\smallskip

For~$L=\infty$, fix~$m$ and, for~$N\geq m$, let~$u_N\in u+H_0^1(\cu_m)$ be~$\widetilde\a_N$-harmonic. By~\eqref{e.cutoff.Holder.replacements}, Proposition~\ref{p.cutoff.Holder.regularity} supplies the same random variable~$\mathcal L_N(\gamma,m)$ for every~$N\geq m$, and multiplication of the coefficient by~$\a_N(0)^{-1}$ changes neither the equation nor the estimate~\eqref{e.large.scale.energy.multifractal}. Lemma~\ref{l.finite.cutoff.coefficient.convergence} and the energy identity for~$u_N-u$ give~$\|\nabla u_N-\nabla u\|_{L^2(\cu_m)}\to0$. The large-scale H\"older and energy estimates therefore pass to the limit. Taking~$\mathcal L_\infty(\gamma,m)=\mathcal L_N(\gamma,m)$, which is the same for every~$N\geq m$, also gives~\eqref{e.limiting.Holder.minimal.scale}.

\smallskip

For the Liouville assertion, fix~$L\in\N_0\cup\{\infty\}$ and suppose that~$u$ satisfies~\eqref{e.multifractal.Liouville.power.growth}. By~\eqref{e.limiting.Holder.minimal.scale} and Borel--Cantelli,~$\mathcal L_L(\gamma_{\mathrm{reg}},m)\leq m$ for all sufficiently large~$m$. Choose a sequence of radii~$R\to\infty$ and constants~$c=c(R)\in\mathbb R$ such that~$R^{-\gamma_{\mathrm{reg}}}\|u-c\|_{\underline L^2(B_{R})}\to0$, and let~$m$ be the largest integer such that~$\cu_{m}\subseteq B_{R}$. Since~$3^{m}$ is comparable to~$R$ and~$(u)_{\cu_{m}}$ minimizes~$c\mapsto\|u-c\|_{\underline L^2(\cu_{m})}$, we have~$3^{-\gamma_{\mathrm{reg}}m}\|u-(u)_{\cu_m}\|_{\underline L^2(\cu_m)}\longrightarrow0$. For every fixed~$z\in\Zd$, we apply~\eqref{e.large.scale.Holder.multifractal} with~$\gamma=\gamma_{\mathrm{reg}}$ and~$n=0$ on~$z+\cu_0$, which is contained in~$\cu_{m-1}$ once~$R$ is large enough, and let~$R\to\infty$. Thus~$\nabla u=0$ almost everywhere in~$\Rd$, so~$u$ is constant.
\end{proof}

\section{The scaling limit}
\label{s.scaling.limit}

In this section we state Theorem~\ref{t.A} in the precise form in which we prove it. Sections~\ref{sec:mfd-uniqueness} and~\ref{sec:mfd-process-convergence} prove part~(i) of Theorem~\ref{t.scaling.limit} and the convergence of the measures in part~(ii), and Section~\ref{lim:sec-process-measure} proves the remaining assertions.

\smallskip

To prove almost-sure convergence, we realize all rescaled processes using one family of independent fields~$(\g_j)_{j\in\Z}$, with~$\g_j$ having the law of~$\g_0(3^{-j}\,\cdot)$. Rescaling space by the factor~$3^{-N}$ turns the field with scale~$3^k$ into a field with scale~$3^{k-N}$. Accordingly, we take the rescaled potential at cutoff~$N$ to be~$H+\sum_{j=0}^N\g_{-j}$, where~$H(x)\coloneqq\sum_{j=1}^\infty(\g_j(x)-\g_j(0))$ collects the fields with scales larger than one; the series converges locally in~$C^{1,1}$. Increasing~$N$ adds the field with scale~$3^{-N}$ and leaves the others unchanged. Let~$K^\omega_{N,x}$ be the law on continuous paths of the diffusion started at~$x$ and constructed from the elliptic resolvents of the generator~$\ahom_N^{-1}\bigl(\Delta+\nabla(H+\sum_{j=0}^N\g_{-j})\cdot\nabla\bigr)$. Its stochastic differential description is
\begin{equation}
	d\widetilde Z^{(N)}_t=\ahom_N^{-1}\nabla\Bigl(H+\sum_{j=0}^N\g_{-j}\Bigr)(\widetilde Z^{(N)}_t)\,dt+\sqrt{2/\ahom_N}\,dW_t\,,
	\qquad \widetilde Z^{(N)}_0=x\,.
	\label{e.intro.rescaled.generator}
\end{equation}
For each~$N$, the fields~$\g_{k-N}$ with~$k\in\N_0$ have the same joint law as the fields~$\g_k(3^N\,\cdot)$. Thus~$K^\omega_{N,x}$, viewed as a random probability measure, has the same distribution as the law of~$(3^{-N}X_{3^{2N}t/\ahom_N})_{t\geq0}$ conditional on the environment, with~$X_0=3^Nx$. The diffusion~\eqref{e.intro.rescaled.generator} is reversible with respect to
\begin{equation*}
	\mu_N(dx)\coloneqq\exp\Bigl(H(x)+\sum_{j=0}^N\bigl(\g_{-j}(x)-\tau^2\bigr)\Bigr)\,dx\,.
\end{equation*}
The measure~$\mu_N^0\coloneqq e^{-H}\mu_N$ is stationary and has the law of~$\a_N(3^Nx)\,dx$. For Gaussian fields, it is the standard approximation of Gaussian multiplicative chaos: its density is obtained by exponentiating~$\sum_{j=0}^N\g_{-j}$ and dividing by the expectation of that exponential. We write~$d_{\mathrm{weak}}$ for a bounded complete metric on the probability measures on~$C([0,\infty),\Rd)$ that induces weak convergence.

\begin{theorem}[Scaling limit]
\label{t.scaling.limit}
Assume~\ref{a.g1}--\ref{a.g4}. There exists~$\gmcdelta_0(d)>0$ such that, if~$0<\gmcdelta\leq\gmcdelta_0$, then there exists a continuous probability kernel~$x\mapsto K^\omega_x$ such that the following statements hold.
\begin{enumerate}[label=(\roman*)]
\item \underline{\emph{Convergence of the path laws}}: Almost surely, for every compact set~$F\subseteq\Rd$,
\begin{equation}
	\lim_{N\to\infty}\sup_{x\in F}
	d_{\mathrm{weak}}(K^\omega_{N,x},K^\omega_x)=0\,.
	\label{eq:mfd-main-as}
\end{equation}
Almost surely, the laws~$(K^\omega_x)_{x\in\Rd}$ define a continuous, conservative, strong Markov process~$Z$, started at~$x$ under~$K_x^\omega$. On the same event, its semigroup maps bounded continuous functions to continuous functions, and, for every bounded open set~$U$, the mean exit time of~$Z$ from~$U$ under~$K_x^\omega$ is bounded uniformly in~$x\in U$. For every~$x\in\Rd$, the annealed law of~$(3^{-N}X_{3^{2N}t/\ahom_N})_{t\geq0}$, with~$X_0=3^Nx$, converges to~$\E[K^\omega_x]$.

\item \underline{\emph{The reversible measure}}: Almost surely, the measures~$\mu_N^0$ and~$\mu_N$ converge locally weakly to measures~$\mu^0$ and~$\mu=e^H\mu^0$, respectively. Both limits are locally finite, nonatomic, of full support, and singular with respect to Lebesgue measure. The measure~$\mu^0$ is stationary and not deterministic, and it satisfies~$\E[\mu^0(dx)]=dx$. Almost surely, the process~$Z$ is reversible with respect to~$\mu$, which is invariant: for every nonnegative Borel function~$f$ and every~$t\geq0$,
\begin{equation*}
	\int\biggl(\int f(Z_t)\,dK_x^\omega\biggr)\mu(dx)=\int f\,d\mu\,.
\end{equation*}

\item \underline{\emph{Anomalous scaling}}: There exist constants~$C\geq1$ and~$\eta>0$, depending only on~$d$ and the law of~$\g_0$, such that, for all~$\ell,m\in\N_0$ with~$\ell\leq m$,
\begin{equation*}
	\frac{\ahom_m}{\ahom_\ell}\leq C3^{-\eta(m-\ell)}\,.
\end{equation*}
For a physical side length~$r>0$, write~$\sigma_r\coloneqq\inf\{t\geq0:Z_t-Z_0\notin r\cu_0\}$, as in Theorem~\ref{t.A}. For every~$x\in\Rd$ and every integer~$k\geq0$,
\begin{equation}
	\E\left[\int\sigma_{3^{-k}}\,dK_x^\omega\right]\leq C3^{-(2+\eta)k}\,.
	\label{e.intro.limit.exit}
\end{equation}
For every~$x\in\Rd$, almost surely in~$\omega$, the law~$K_x^\omega$ is mutually singular with every mixture of Brownian laws with arbitrary covariance matrices, and~$K_x^\omega$-almost every path satisfies
\begin{equation}
	\sup\{\gamma\geq0:|Z_t-Z_0|=O(t^\gamma)\text{ as }t\downarrow0\}
	\leq\nf1{(2+\eta)}\,.
	\label{e.intro.limit.Holder}
\end{equation}
In dimension two, one can take~$\eta=\tau^2/\log3$.

\item \underline{\emph{Non-Gaussian marginals}}: Almost surely, for every~$x\in\Rd$ and every~$t>0$, the law of~$Z_t$ under~$K_x^\omega$ is absolutely continuous with respect to~$\mu$, and it is not a Gaussian law.
\end{enumerate}
\end{theorem}

Theorem~\ref{t.A} follows, with~$\mathbf P_x\coloneqq K^\omega_x$: for its exit-time bound, apply~\eqref{e.intro.limit.exit} with~$3^{-k-1}<r\leq3^{-k}$.

\section{Uniqueness of the limit form}
\label{sec:mfd-uniqueness}

In this section we prove that the energies~$\int A_N|\nabla u|^2\,dx$ of the rescaled processes have a unique limit: all subsequential limits of the Dirichlet forms on cubes coincide (Theorem~\ref{mfd:thm-c1}). Section~\ref{sec:mfd-process-convergence} combines this limit with the limit of the reversible measures~$\mu_N$ to construct the limiting diffusion.

\smallskip

Compactness of the energies comes from their insensitivity to fields of small scale, and gives subsequential limits. To identify them, we adapt the normalization comparison of~\citet[Section~3.4]{DGComparison} and the local replacement arguments of~\citet[Section~6]{GM} and~\citet[Section~4.5]{DGUniqueness}, using harmonic extensions in place of paths. We compare two subsequential limit forms~$\mathscr E$ and~$\mathscr F$. They bound each other, and we call the factors of the largest multiple of~$\mathscr E$ below~$\mathscr F$ and of the smallest multiple of~$\mathscr E$ above~$\mathscr F$ the \textbf{optimal comparison constants}. These constants are deterministic. If they differed, then, at a positive proportion of scales, on the cubes where the estimates of Section~\ref{mfd:sec-branches} hold, the minimal~$\mathscr F$-energy for boundary values close to an affine function would be smaller, by a fixed fraction of the gap between the constants, than the upper constant times the minimal~$\mathscr E$-energy, or the symmetric statement would hold. Replacing an~$\mathscr E$-harmonic function on such cubes by the~$\mathscr F$-harmonic function with the same boundary values would then improve one of the constants. The difficulty is quantitative: both the energy saved by a replacement and the errors it must absorb are proportional to the gap between the optimal comparison constants. To obtain errors of this size, we apply the Efron--Stein inequality to the difference between the minimal energy of one form and a multiple of the minimal energy of the other, rather than to each energy separately.

When the optimal comparison constants agree,~$\mathscr F=c\mathscr E$ for a deterministic~$c>0$, which we call the \textbf{proportionality constant}. With the same reversible measure, the diffusion associated with~$\mathscr F$ runs~$c$ times as fast as the diffusion associated with~$\mathscr E$. Thus proportionality leaves a possible change of time. The effective-diffusivity normalization fixes the speed and forces~$c=1$.

\smallskip

The section is organized as follows.
\begin{itemize}
\item In Section~\ref{mfd:sec-environment} we introduce the rescaled coefficients and energies.

\item In Section~\ref{mfd:sec-local} we bound, uniformly in~$N$, the energy that solutions carry in small balls, including solutions of Neumann problems whose datum is approximated by a volume source. These bounds control the change in the energy on a fixed cube when a field of small scale is resampled.

\item In Section~\ref{mfd:sec-resampling} we use this control and the Efron--Stein inequality to prove compactness of the finite-volume energies.

\item In Section~\ref{mfd:sec-local-form} we show that, along subsequences, the Green operators of the Dirichlet problems converge and that their limits define Dirichlet forms. These limit forms are local and consistent under restriction to subcubes, and multiplying the coefficient by a positive continuous function multiplies their energy measures by the same function.

\item In Section~\ref{mfd:sec-branches} we bound, on a single event that does not depend on the solution, the number of cubes of each chain of nested cubes on which the regularity, extension and harmonic-approximation estimates of Lemma~\ref{mfd:in-deterministic} fail. Integrating this bound against the energy measure of a solution selects cubes that carry a fixed fraction of its energy, on which these estimates hold and the solution is close to an affine function.

\item In Section~\ref{mfd:sec-compare} we prove that the optimal comparison constants of two subsequential limit forms agree, so that the two forms are proportional. The required concentration estimate for differences of minimal energies is Proposition~\ref{mfd:prop-conc}.

\item In Section~\ref{mfd:sec-scalar} we show that the proportionality constant is one. The coarse-grained matrices of~$e^{-H}A_N$ on large cubes approach the identity uniformly in~$N$, and both limit forms inherit this normalization.
\end{itemize}

Throughout, Assumptions~\ref{a.g1}--\ref{a.g4} hold. Before reducing~$\gmcdelta$, we fix all exponents and tolerances and a finite list of moment orders. Bounds in~$L^p(\P)$ for the orders in this list are uniform in~$N$. When a result below is stated for every finite~$p$, the moment orders needed for that application, including the larger orders used in H\"older's inequality, spatial maxima and resampling, are fixed first; the disorder threshold then depends on this finite list. This convention applies in particular to Lemmas~\ref{mfd:lem-neumann-error},~\ref{mfd:lem-prefix-limit} and~\ref{mfd:lem-band}, and Proposition~\ref{mfd:prop-conc}. Such a result does not give a single~$\gmcdelta>0$ that works for all moment orders.

\subsection{Rescaled coefficients and energies}
\label{mfd:sec-environment}

In this subsection we write the Dirichlet form of the rescaled process with index~$N$ in terms of the relabelled fields, and we define the forms, minimal energies and Dirichlet Green operators on cubes. The coefficient~$A_N$ is the product of the positive factor~$e^H$, determined by the fields at scales larger than one, and a stationary coefficient, as in~\eqref{eq:mfd-normalization}; Lemma~\ref{mfd:lem-infrared} bounds the moments of~$H$. The forms and the Dirichlet Green operators are defined in~\eqref{eq:mfd-energy} and~\eqref{eq:mfd-quadratic-response}.

\smallskip

At fixed~$N$, relabel the fields by replacing~$\g_i$ with~$\g_{i-N}(3^{-N}\,\cdot)$, and let~$\a_N$ be the cutoff coefficient of~\eqref{e.intro.cutoff.coefficients} formed from the relabelled fields. The diffusion with law~$K^\omega_{N,x}$ has Dirichlet form~$\int A_N|\nabla u|^2\,dx$ on~$L^2(\mu_N)$, where~$A_N\coloneqq\ahom_N^{-1}\,d\mu_N/dx$ (see~\eqref{e.intro.rescaled.generator}). In terms of the relabelled fields,
\begin{equation}
	A_N(x)=e^{H(x)}\frac{\a_N(3^Nx)}{\ahom_N}\,,
	\qquad A_N^0(x)\coloneqq e^{-H(x)}A_N(x)\,.
	\label{eq:mfd-normalization}
\end{equation}
If~$H$ is replaced by one of its partial sums, then~$A_N$ becomes a spatially constant multiple of a stationary cutoff coefficient, so the results of Sections~\ref{s.self.similar} and~\ref{s.Holder.regularity} apply.

\begin{lemma}[Moments of the field~$H$]
\label{mfd:lem-infrared}
\begin{itemize}
\item For every cube~$w+R\cu_0$ there is~$C_{w,R}<\infty$ such that, for every~$p\geq2$ and~$\lambda\geq0$, the field~$H$ and its partial sums~$\sum_{j=1}^{L'}(\g_j-\g_j(0))$ satisfy
\begin{equation*}
	\bigl\|\|H\|_{C^{1,1}(\overline{w+R\cu_0})}\bigr\|_{L^p}\leq C_{w,R}\gmcdelta\sqrt p\,,
	\qquad
	\E e^{\lambda\|H\|_{C^{1,1}(\overline{w+R\cu_0})}}\leq2e^{C_{w,R}\lambda^2\gmcdelta^2}
	\,,
\end{equation*}
with the same bounds for the partial sums, uniformly in the number~$L'$ of terms.
\item The partial sums converge to~$H$ locally in~$C^{1,1}$ almost surely.
\end{itemize}
\end{lemma}
\begin{proof}
As in the proof of Lemma~\ref{l.finite.cutoff.coefficient.convergence}, scaling, the mean value theorem and~\eqref{e.coefficient.field.regularity} give~$\|\|\g_j-\g_j(0)\|_{C^{1,1}(\overline{w+R\cu_0})}\|_{L^p}\leq C_{w,R}\gmcdelta\sqrt p\,3^{-j}$ for~$j\geq1$. Sum over~$j$ and expand the exponential. Since the sum of these norms has finite mean, the partial sums also converge in~$C^{1,1}$ almost surely.
\end{proof}

On a bounded cube~$w+R\cu_0$, let
\begin{equation}
	\mathscr E_N^{w+R\cu_0}(u,v)\coloneqq\int_{w+R\cu_0}A_N\nabla u\cdot\nabla v\,dx\,,
	\qquad \Gamma_N(u)\coloneqq A_N|\nabla u|^2\,dx\,.
	\label{eq:mfd-energy}
\end{equation}
The form~$\mathscr E_N^{w+R\cu_0}$ has domain~$H^1_0(w+R\cu_0)$, and we use the same integral for functions in~$H^1(w+R\cu_0)$ when prescribing boundary values. Let~$\Lambda_{N,w+R\cu_0}(g)$ be the minimum of~$\mathscr E_N^{w+R\cu_0}(u)$ over~$u\in H^1(w+R\cu_0)$ with~$u-g\in H^1_0(w+R\cu_0)$. Let~$G_N^{w+R\cu_0}$ be the Dirichlet Green operator, the inverse of~$-\nabla\cdot(A_N\nabla)$ on~$L^2(w+R\cu_0,dx)$ with zero Dirichlet boundary values, so that
\begin{equation}
	\langle f,G_N^{w+R\cu_0}f\rangle
	=\sup_{v\in H^1_0(w+R\cu_0)}\left(2\int_{w+R\cu_0}fv\,dx-\mathscr E_N^{w+R\cu_0}(v)\right)\,.
	\label{eq:mfd-quadratic-response}
\end{equation}
For~$p\in\Rd$, the quantity~$R^d\,p\cdot \a_*^{-1}(w+R\cu_0;A_N)p$ is the energy~$\mathscr E_N^{w+R\cu_0}(v)$ of a solution~$v$ of~$-\nabla\cdot(A_N\nabla v)=0$ in~$w+R\cu_0$ with boundary flux~$A_N\nabla v\cdot\mathbf n=p\cdot\mathbf n$, and it has the variational formula~\eqref{e.variational.astar}. Since~$d\mu_N/dx=\ahom_N A_N$, the form~$\mathscr E_N^{w+R\cu_0}$, regarded as a form on~$L^2(\mu_N)$, has generator~$u\mapsto\nabla\cdot(A_N\nabla u)/(\ahom_NA_N)$, the generator of the diffusion~\eqref{e.intro.rescaled.generator} killed on leaving~$w+R\cu_0$; its resolvent is studied in Section~\ref{mfd:sec-speed}.

\smallskip

For~$k\in\Z$ and~$z\in\Rd$, let~$\Theta_{k,z}$ replace the field~$\g_i$ by~$\g_{i-k}(z+3^{-k}\,\cdot)$, for every~$i\in\Z$.

\smallskip

For a cube~$w+3^{-\ell}\cu_0$ with integer~$\ell\leq N$, including negative~$\ell$, which give cubes of side larger than one, let
\begin{equation*}
	s_N(w+3^{-\ell}\cu_0)\coloneqq
	\frac{e^{(N-\ell+1)\tau^2}\ahom_{N-\ell}}{e^{(N+1)\tau^2}\ahom_N}
	\exp\left(H(w)+
	\begin{cases}
		\displaystyle\sum_{j=0}^{\ell-1}\g_{-j}(w),&\ell\geq0\,,\\
		\displaystyle-\sum_{j=\ell}^{-1}\g_{-j}(w),&\ell<0\,.
	\end{cases}
	\right)\,.
\end{equation*}

\subsection{Local estimates and approximation of the Neumann datum}
\label{mfd:sec-local}

In this subsection we bound, uniformly in~$N$, the energy that solutions of~$-\nabla\cdot(A_N\nabla u)=f$ carry in small balls. For the Dirichlet problem this is the Morrey bound~\eqref{eq:mfd-4} of Proposition~\ref{mfd:prop-growth}. For the Neumann problem with a bounded source, we reflect solutions across the faces of a cube (Proposition~\ref{mfd:prop-folded-iteration}) to obtain the same bound (Proposition~\ref{mfd:cor-neumann-source}), and we approximate the Neumann datum of a constant flux by such a source (Lemma~\ref{mfd:lem-neumann-error}). The bounds follow from the coarse-grained inequalities of Section~\ref{s.coarse.graining}, the moments of Proposition~\ref{p.coarse.grained.bound} and the regularity estimates of Proposition~\ref{p.cutoff.Holder.regularity}. The Morrey bound and the coercivity bound~\eqref{eq:mfd-1} are used again in Section~\ref{s.tightness}.

\smallskip

At radii above the smallest field scale~$3^{-N}$, we apply Proposition~\ref{p.cutoff.Holder.regularity} after replacing~$H$ by a partial sum. At smaller radii, elliptic estimates for coefficients with small oscillation complete the bounds. At these small radii we use Lemma~\ref{mfd:lem-extremes} for the coefficient and its logarithmic gradient. The moments of these bounds grow with~$N$, and positive powers of~$3^{-N}$ absorb this growth.

\smallskip

\subsubsection{Coercivity and extension}

The coarse-grained bounds of Section~\ref{s.self.similar} also cover scales below the smallest field scale. We first apply them on a rescaled cube. We call a cube~$z+3^{-k}\cu_0$ of a triadic \textbf{mesh}, a partition of a cube or of~$\Rd$ into congruent cubes, a \textbf{cell}.

\begin{lemma}[Ellipticity moments on a cell]
\label{mfd:lem-cell-ellipticity}
Fix~$s\in(0,1]$,~$q\in\{1,2\}$ and a bounded cube~$w+R\cu_0$.
\begin{itemize}
\item For every~$p<\infty$ and all~$\gmcdelta$ sufficiently small depending only on~$d$,~$s$,~$q$ and~$p$, uniformly for~$0\leq k\leq N$ and cubes~$z+3^{-k}\cu_0\subset w+R\cu_0$,
\begin{equation}
	\bigl\|\Lambda_{s,q}(z+3^{-k}\cu_0;A_N)+\lambda_{s,q}^{-1}(z+3^{-k}\cu_0;A_N)\bigr\|_{L^p}
	\leq C_{p,s,w,R}\exp\bigl(C(p+p^2)\gmcdelta^2k\bigr)\,.
	\label{eq:mfd-cell-moment}
\end{equation}
\item In particular, on a fixed cell~$z+3^{-k}\cu_0$ these moments are bounded uniformly in~$N$.
\end{itemize}
\end{lemma}
\begin{proof}
Choose a finite moment order~$\xi\geq4p$ large enough for Proposition~\ref{p.coarse.grained.bound} at exponent~$s/2$, then decrease~$\gmcdelta$ so that its remaining hypotheses hold. The estimate~\eqref{e.homogenize.at.scale.m.E}, the comparison of the weighted sums at exponents~$s/2$ and~$s$, and~\eqref{e.bound.Lambdas.by.Es} give
\begin{equation*}
	\sup_m\bigl\|\Lambda_{s,q}(\cu_0;A_m^0)
	+\lambda_{s,q}^{-1}(\cu_0;A_m^0)\bigr\|_{L^{2p}}<\infty\,.
\end{equation*}

\smallskip

By definition of~$A_N$, the coefficient~$A_N(z+3^{-k}y)$ is the product of three factors: the scalar~$\frac{\ahom_{N-k}}{\ahom_N}\exp\bigl(H(z)+\sum_{j<k}(\g_{-j}(z)-\tau^2)\bigr)$, the ratio of~$\exp\bigl(H+\sum_{j<k}\g_{-j}\bigr)$ at~$z+3^{-k}y$ to its value at~$z$, and a function of~$y$ with the law of~$A_{N-k}^0$. Summing the bounds on the gradients of the fields as in the proof of Lemma~\ref{mfd:lem-infrared} bounds the oscillation of~$H+\sum_{j<k}\g_{-j}$ over the cube~$z+3^{-k}\cu_0$ by a sub-Gaussian variable with parameter~$C_{w,R}\gmcdelta$. The scalar involves the values of~$k$ independent fields. By~\eqref{e.annealed.ordering.ratio}, for every fixed~$P<\infty$ the sum of the scalar and its inverse, multiplied by the exponential of this oscillation, has~$L^P$ norm at most a constant depending on~$P$ and~$w+R\cu_0$ times~$e^{C(P+P^2)\gmcdelta^2k}$. Monotonicity and homogeneity of the coarse-grained ellipticity constants, followed by H\"older's inequality, prove~\eqref{eq:mfd-cell-moment}.
\end{proof}

\begin{lemma}[Coercivity of the form~$\mathscr E_N^{z+r\cu_0}$]
\label{mfd:lem-coercivity}
Let~$z+r\cu_0$ be a fixed cube with~$r\leq1$.
\begin{itemize}
\item For every~$N\geq0$ there is a random variable~$Q_N$ such that, for every~$v\in H^1_0(z+r\cu_0)$ and every~$u\in H^1(z+r\cu_0)$,
\begin{equation}
	\|v\|_{H^{\nf34}(z+r\cu_0)}^2\leq Q_N\mathscr E_N^{z+r\cu_0}(v)
	\quad\text{and}\quad
	\|u-(u)_{z+r\cu_0}\|_{H^{\nf34}(z+r\cu_0)}^2\leq Q_N\mathscr E_N^{z+r\cu_0}(u)\,.
	\label{eq:mfd-1}
\end{equation}
\item For every~$p<\infty$ there is~$\gmcdelta_0(d,p)>0$ such that, if~$\gmcdelta\leq\gmcdelta_0$, then~$\sup_N\|Q_N\|_{L^p}<\infty$.
\end{itemize}
\end{lemma}

\begin{proof}
On~$\cu_0$, take~$s=\nf18$. By the functional inequalities~\eqref{e.besov.grad.poincare},~\eqref{e.nabla.u.detach} and~\eqref{e.CG.Poincare.trace.zero}, followed by~$B^{\nf{7}{8}}_{2,\infty}\subset H^{\nf34}$, we obtain~\eqref{eq:mfd-1} with~$Q_N\coloneqq C\lambda_{1/8,1}^{-1}(\cu_0;A_N)$. Lemma~\ref{mfd:lem-cell-ellipticity} bounds its~$L^p$ norms uniformly in~$N$. For a fixed smaller cube, rescale the same inequalities and use~\eqref{eq:mfd-cell-moment}. If its side is not triadic, then the remaining dilation lies in~$[1,3]$. Apply the same argument after replacing~$\gmcdelta$ by~$C_d\gmcdelta$, with constants depending on~$z+r\cu_0$.
\end{proof}

Fix~$\beta\in(\nf12,1)$ before choosing~$\gmcdelta$, and let~$\upsilon\coloneqq\frac12(\beta-\nf12)$. For a bounded cube~$w+R\cu_0$, fix finitely many grids, each consisting of the cubes~$z+r\,\cu_0$ with~$r=3^{-k}$ for an integer~$k\geq0$ and~$z$ in~$z_0+r\mathbb Z^d$ for a point~$z_0$ fixed for each grid. We call their cubes contained in~$w+R\cu_0$ \textbf{grid cubes}.

\begin{lemma}[Boundary extension and bounds on grid cubes]
\label{mfd:lem-extension}
\begin{enumerate}[label=(\roman*)]
\item There is~$C(d,\beta)<\infty$ such that, for every cube~$z+r\cu_0$ with~$r\leq1$ and every~$g\in C^\beta(\partial(z+r\cu_0))$,
\begin{equation}
	\Lambda_{N,z+r\cu_0}(g)\leq C\,\Lambda_{\upsilon/2,2}(z+r\cu_0;A_N)\,r^{d-2}\,\|g(z+r\,\cdot)\|_{C^\beta(\partial \cu_0)/\mathbb R}^2\,.
	\label{eq:mfd-2}
\end{equation}
\item For every~$\kappa>0$ and~$p<\infty$ there is~$\gmcdelta_0(d,\beta,\kappa,p)>0$ such that, if~$\gmcdelta\leq\gmcdelta_0$, then there are random variables~$K_N$ with~$\sup_N\|K_N\|_{L^p}<\infty$ such that, for every~$N$ and every grid cube~$z+r\cu_0$ of side~$r\geq3^{-N}$,
\begin{equation}
	\Lambda_{\upsilon/2,2}(z+r\cu_0;A_N)+\lambda_{\upsilon/2,2}^{-1}(z+r\cu_0;A_N)\leq K_Nr^{-\kappa}\,.
	\label{eq:mfd-3}
\end{equation}
\end{enumerate}
\end{lemma}

\begin{proof}
\emph{Step 1.} We prove~\eqref{eq:mfd-2}. The change of variables~$x=z+ry$ multiplies the unnormalized energy by~$r^{d-2}$ and leaves the volume-normalized constant~$\Lambda_{\upsilon/2,2}$ unchanged, so we may take~$z+r\cu_0=\cu_0$. We construct an extension~$g_{\mathrm{ext}}$ of~$g$ directly. First extend~$g$ to a~$\beta$-H\"older function~$\tilde g$ on~$\Rd$, and let~$\tilde g_k$ be its mollification at scale~$3^{-k}$. Use a smooth Whitney partition~$\{\chi_k\}$ with bounded overlap, subordinate to the layers where~$\operatorname{dist}(\cdot,\partial\cu_0)$ is comparable to~$3^{-k}$, satisfying~$|D^i\chi_k|\leq C3^{ik}$ for~$i=1,2$ and~$\sum_k\chi_k=1$. Let~$g_{\mathrm{ext}}\coloneqq\sum_k\chi_k\tilde g_k$. The mollifications~$\tilde g_k$ on overlapping layers differ by at most~$C[g]_{C^\beta}3^{-k\beta}$. Differentiating~$g_{\mathrm{ext}}$ and using~$\sum_kD^i\chi_k=0$ shows that~$g_{\mathrm{ext}}$ is smooth inside the cube, has continuous boundary values~$g$, and satisfies
\begin{equation*}
	|\nabla g_{\mathrm{ext}}(x)|\leq C[g]_{C^\beta}\operatorname{dist}(x,\partial\cu_0)^{\beta-1}\,,
	\qquad |D^2g_{\mathrm{ext}}(x)|\leq C[g]_{C^\beta}\operatorname{dist}(x,\partial\cu_0)^{\beta-2}\,.
\end{equation*}
For the fractional norm of~$\nabla g_{\mathrm{ext}}$, split the pairs~$(x,y)$ according to whether~$|x-y|$ is smaller than~$\operatorname{dist}(x,\partial\cu_0)/2$. On the closer pairs use the Hessian bound; on the remaining pairs use the gradient bound. Both contributions are bounded by
\begin{equation*}
	C[g]_{C^\beta}^2\int_{\cu_0}\operatorname{dist}(x,\partial\cu_0)^{2\beta-2-2\upsilon}\,dx
	\leq C(d,\beta)[g]_{C^\beta}^2\,,
\end{equation*}
since~$2\beta-2-2\upsilon>-1$. The same gradient bound gives its~$L^2$ norm. Thus~$g_{\mathrm{ext}}\in H^{1+\upsilon}(\cu_0)$ and~$\|\nabla g_{\mathrm{ext}}\|_{H^\upsilon}\leq C(d,\beta)[g]_{C^\beta}$. Applying~\eqref{e.cg.RHS} with zero right side,~$s=\upsilon$ and boundary datum~$g_{\mathrm{ext}}$, we find that the~$A_N$-harmonic function~$v$ with~$v=g_{\mathrm{ext}}$ on~$\partial \cu_0$ satisfies
\begin{equation*}
	\Lambda_{N,\cu_0}(g)=\|A_N^{\nf12}\nabla v\|_{L^2}^2
	\leq C\upsilon^{-3}\Lambda_{\upsilon/2,2}(\cu_0;A_N)
	\|\nabla g_{\mathrm{ext}}\|_{H^\upsilon}^2\,.
\end{equation*}
Subtracting a constant from~$g$ does not change~$\Lambda_{N,\cu_0}(g)$ or the H\"older seminorm of~$g$. The extension just constructed therefore gives~$\|\nabla g_{\mathrm{ext}}\|_{H^\upsilon}\leq C\|g\|_{C^\beta(\partial \cu_0)/\mathbb R}$. Since~$\upsilon$ depends only on~$\beta$, this proves~\eqref{eq:mfd-2}.

\smallskip

\emph{Step 2.} We prove~\eqref{eq:mfd-3}. Since~$L^p$ norms increase with~$p$, we may assume that~$p>1$ and~$p\kappa>d$. Apply Lemma~\ref{mfd:lem-cell-ellipticity} with exponents~$\upsilon/2$ and~$2$. There are at most~$C3^{dk}$ grid cubes at depth~$k$, so taking the maximum over them introduces a factor~$C3^{\frac{d}{p}k}$ in~$L^p$. Reduce~$\gmcdelta$ so that~$C(p+p^2)\gmcdelta^2<\frac12(\kappa-d/p)\log3$. For~$K_N$ defined as the supremum, over~$k\leq N$ and grid cubes~$z+3^{-k}\cu_0$ of depth~$k$, of~$3^{-\kappa k}\bigl(\Lambda_{\upsilon/2,2}(z+3^{-k}\cu_0;A_N)+\lambda^{-1}_{\upsilon/2,2}(z+3^{-k}\cu_0;A_N)\bigr)$, Minkowski's inequality gives
\begin{equation*}
	\|K_N\|_{L^p}
	\leq C_p\sum_{k\geq0}3^{-(\kappa-\frac{d}{p})k}
	e^{C(p+p^2)\gmcdelta^2k}\leq C_p\,.
\end{equation*}
The right side does not depend on~$N$, so~$\sup_N\|K_N\|_{L^p}<\infty$. By the definition of~$K_N$, every grid cube~$z+r\cu_0$ of side~$r=3^{-k}\geq3^{-N}$ satisfies~\eqref{eq:mfd-3}.
\end{proof}

\subsubsection{Morrey bounds for the Dirichlet problem}

To apply the regularity estimates with a bounded volume source, we first write that source as a divergence.

For~$f\in L^\infty$ supported in a cube of side~$R$, let~$g$ be the gradient of the Newtonian potential of~$f$, that is, the convolution of~$f$ with the gradient of the fundamental solution of the Laplacian.

\begin{lemma}[The gradient of the Newtonian potential]
\label{mfd:lem-primitive}
\begin{itemize}
\item The field~$g$ satisfies~$\nabla\cdot g=f$.
\item There is~$C(d)<\infty$ such that, for all~$x$ and~$y$ with~$|x-y|\leq R$,
\begin{equation}
	|g(x)-g(y)|\leq C\|f\|_\infty|x-y|\bigl(1+\log(R/|x-y|)\bigr)\leq CR^{\nf12}\|f\|_\infty|x-y|^{\nf12}\,.
	\label{eq:mfd-primitive}
\end{equation}
\item In particular,~$[g]_{W^{\nf12,\infty}}\leq CR^{\nf12}\|f\|_\infty$ on that cube.
\end{itemize}
\end{lemma}

\begin{proof}
The kernel of the convolution and its gradient are bounded by~$C|z|^{1-d}$ and~$C|z|^{-d}$. If~$x$ is farther than~$2\sqrt dR$ from the supporting cube, the derivative bound gives~$|g(x)-g(y)|\leq C\|f\|_\infty|x-y|$ directly. Otherwise the support lies in~$B_{CR}(x)$. We split the convolution according to whether the integration variable~$z$ satisfies~$|z-x|\leq2|x-y|$. The part with~$|z-x|\leq2|x-y|$ is bounded by~$C\|f\|_\infty|x-y|$, both at~$x$ and at~$y$. For~$|z-x|>2|x-y|$, the mean value theorem bounds the difference of the kernel at~$x-z$ and at~$y-z$ by~$C|x-y||x-z|^{-d}$. The integral over~$2|x-y|<|z-x|<CR$ is at most~$C|x-y|\|f\|_\infty\log(CR/|x-y|)$. The second inequality in~\eqref{eq:mfd-primitive} follows from~$(|x-y|/R)^{\nf12}(1+\log(R/|x-y|))\leq C$.
\end{proof}

For a fixed cube~$w+R\cu_0$, set~$D_N\coloneqq3^{-N}\|\nabla\log A_N\|_{L^\infty(w+R\cu_0)}$.

\begin{lemma}[Pointwise extremes of the coefficient]
\label{mfd:lem-extremes}
For every~$p$ there are~$C_p<\infty$ depending on~$p$ and~$w+R\cu_0$,~$C_d<\infty$ and~$c_d>0$ such that, if~$\gmcdelta\leq c_d/p$, then
\begin{equation*}
	\|D_N\|_{L^p}\leq C_p(1+N)^{\nf12}\,,
\end{equation*}
and
\begin{equation}
	\bigl\|\|A_N\|_{L^\infty(w+R\cu_0)}+\|A_N^{-1}\|_{L^\infty(w+R\cu_0)}\bigr\|_{L^p}
	\leq C_pe^{(C_d\gmcdelta+C_p\gmcdelta^2)N}\,.
	\label{eq:mfd-6}
\end{equation}
\end{lemma}

The linear term~$C_d\gmcdelta N$ in the exponent arises from taking the maximum over the cubes of side~$3^{-N}$ in~$w+R\cu_0$.

\begin{proof}
\emph{Step 1.} We bound~$D_N$. We have~$\nabla\log A_N=\nabla H+\sum_{j=0}^N\nabla\g_{-j}$ and, in law,~$\|\nabla\g_{-j}\|_{L^\infty(w+R\cu_0)}$ is at most~$3^j$ times the maximum, over at most~$C3^{dj}$ cells, of the regularity norms on unit cells from~\eqref{e.coefficient.field.regularity}. The maximum of~$n$ sub-Gaussian variables with parameter~$\gmcdelta$ has~$L^p$ norm at most~$C\gmcdelta(\sqrt{\log n}+\sqrt p)$. Therefore
\begin{equation*}
	\|D_N\|_{L^p}\leq3^{-N}\bigl\|\|\nabla H\|_{L^\infty(w+R\cu_0)}\bigr\|_{L^p}+\sum_{j\leq N}3^{j-N}C_p\gmcdelta\sqrt{1+j}\leq C_p\sqrt{1+N}\,.
\end{equation*}
This is the first estimate of the lemma.

\smallskip

\emph{Step 2.} We bound~$\|A_N\|_{L^\infty(w+R\cu_0)}$ and~$\|A_N^{-1}\|_{L^\infty(w+R\cu_0)}$. The deterministic factor~$\ahom_N^{-1}e^{-(N+1)\tau^2}$ of~$A_N$ satisfies~$1\leq e^{(N+1)\tau^2}\ahom_N\leq e^{C\gmcdelta^2(N+1)}$ by~\eqref{e.ord.bounds.truncated}, Proposition~\ref{p.homogenized.coefficient.reciprocal.lower} and~\eqref{e.tau.bound}, and~$e^{\|H\|_\infty}$ has moments of every order by Lemma~\ref{mfd:lem-infrared}. By H\"older's inequality, it therefore suffices to bound the moments of~$\exp(\sup_{w+R\cu_0}|\sum_{j\leq N}\g_{-j}|)$. Cover~$w+R\cu_0$ by at most~$C3^{dN}$ cubes~$x_i+3^{-N}\cu_0$. On each of them,~$\sup_{x_i+3^{-N}\cu_0}|\sum_{j\leq N}\g_{-j}|\leq|\sum_{j\leq N}\g_{-j}(x_i)|+\sum_{j\leq N}\operatorname{osc}_{x_i+3^{-N}\cu_0}\g_{-j}$, where the first term is sub-Gaussian with parameter~$C\gmcdelta\sqrt{N+1}$ and the second is sub-Gaussian with parameter~$C\gmcdelta$. H\"older's inequality gives~$\E\exp(q\sup_{x_i+3^{-N}\cu_0}|\sum_{j\leq N}\g_{-j}|)\leq Ce^{Cq^2\gmcdelta^2(N+1)}$. For~$q\geq p$, the inequality~$\exp(p\sup_{w+R\cu_0}|\sum_{j\leq N}\g_{-j}|)\leq(\sum_i\exp(q\sup_{x_i+3^{-N}\cu_0}|\sum_{j\leq N}\g_{-j}|))^{\nf pq}$ and Jensen's inequality bound~$\E\exp(p\sup_{w+R\cu_0}|\sum_{j\leq N}\g_{-j}|)$ by~$(C3^{dN}e^{Cq^2\gmcdelta^2(N+1)})^{\nf pq}$. We choose~$q=\gmcdelta^{-1}\sqrt{d\log3/C}$, which is admissible when~$\gmcdelta\leq c_d/p$; the exponent is then at most~$C_dp\gmcdelta(N+1)$.
\end{proof}

Fix~$t\in(d-1,d)$,~$\alpha\in(0,1)$, a finite list of moment orders and a cube~$y_0+r_0\cu_0$ with~$r_0\leq1$. For~$f\in L^\infty(y_0+r_0\cu_0)$ and~$\varphi\in C^2(\overline{y_0+r_0\cu_0})$, let~$u_N\in H^1(y_0+r_0\cu_0)$ solve~$-\nabla\cdot(A_N\nabla u_N)=f$ in~$y_0+r_0\cu_0$ with~$u_N=\varphi$ on~$\partial(y_0+r_0\cu_0)$.

\begin{proposition}[Morrey bounds at all radii]
\label{mfd:prop-growth}
There is~$\gmcdelta_0(d,t,\alpha,\text{list})>0$ such that, if~$\gmcdelta\leq\gmcdelta_0$, then there are random variables~$K_N$ with~$\sup_N\|K_N\|_{L^p}<\infty$ for every order~$p$ in the list, such that, almost surely, the following estimates hold simultaneously for all~$f$ and~$\varphi$ as above.
\begin{enumerate}[label=(\roman*)]
\item \underline{\emph{Morrey bound}}: for every~$x\in y_0+r_0\cu_0$ and~$r\in(0,1]$,
\begin{equation}
	\Gamma_N(u_N)(B_r(x)\cap(y_0+r_0\cu_0))\leq K_N\bigl(\|f\|_\infty+\|\varphi\|_{C^2}\bigr)^2r^t\,.
	\label{eq:mfd-4}
\end{equation}
\item \underline{\emph{H\"older bound}}:~$\|u_N\|_{C^\alpha(\overline{y_0+r_0\cu_0})}\leq K_N(\|f\|_\infty+\|\varphi\|_{C^2})$. If~$\varphi=0$, then
\begin{equation}
	\|G_N^{y_0+r_0\cu_0}f\|_{C^\alpha(\overline{y_0+r_0\cu_0})}\leq K_N\|f\|_\infty\,.
	\label{eq:mfd-5}
\end{equation}
\end{enumerate}
\end{proposition}

\begin{proof}
The proof proceeds in seven steps, and the quantity we track is the energy~$\Gamma_N(u_N)(B_r(x)\cap(y_0+r_0\cu_0))$ of small balls. Step~1 replaces~$H$ by a partial sum of its series, so that the coefficient becomes a multiple of a cutoff coefficient. Steps~2 and~3 prove the Morrey bound~\eqref{eq:mfd-4} at radii~$r\geq3^{-N}$ from the regularity estimates for cutoff coefficients, and Steps~4 and~5 prove it at radii~$r\leq3^{-N}$, away from~$\partial(y_0+r_0\cu_0)$ and near~$\partial(y_0+r_0\cu_0)$, from estimates for coefficients with small oscillation. The main difficulty lies at these small radii: there the moments of the coefficient bounds of Lemma~\ref{mfd:lem-extremes} grow with~$N$, and positive powers of~$3^{-N}$ absorb this growth. Step~6 removes the truncation, and Step~7 proves the H\"older bounds.

\smallskip

Rescale so that~$r_0=1$. In the coordinates~$y=3^Nx$ of the relabelled fields, the domain is~$3^Ny_0+\cu_N$, a deterministic translate of~$\cu_N$. By linearity, it suffices to consider data with~$\|f\|_\infty+\|\varphi\|_{C^2}\leq1$. Throughout the proof,~$K_N$ denotes a random variable with~$\sup_N\|K_N\|_{L^p}<\infty$ for every order~$p$ in the list, which may change from line to line. Choose the exponents as follows. Take~$\alpha_M\in(\nf12,1)$ with~$\alpha_M>\max\{\alpha,(t+2-d)/2\}$, and note that~$d-2+2\alpha_M>t$. Choose~$t_1$ with~$\max\{t,d-2+2\alpha\}<t_1<d-2+2\alpha_M$, and require~$\gmcdelta$ to be so small that~$\alpha_M\leq1-C\gmcdelta|\log\gmcdelta|^{\nf12}$. We also take~$\gmcdelta$ so small that the tail condition~\eqref{eq:mfd-tail-rate} of Step~3 holds for every order in the list and for the finitely many larger orders used in H\"older's inequality, and that the moment absorptions in Steps~4,~5 and~7 hold; once the exponents and the orders are fixed, each of these is a condition on~$\gmcdelta$ alone.

\smallskip

\emph{Step 1.} We replace~$H$ by a partial sum~$\sum_{j=1}^{L'}(\g_j-\g_j(0))$ of its series, which makes the coefficient a random scalar multiple of a cutoff coefficient. In Steps~1--5 we keep the notation~$H$,~$A_N$,~$u_N$ and~$\Gamma_N$ for the truncated objects; the truncation is removed in Step~6. By~\eqref{eq:mfd-normalization} and~\eqref{e.intro.cutoff.coefficients},
\begin{equation*}
	A_N(3^{-N}y)=\ahom_N^{-1}e^{L'\tau^2}\exp\Bigl(-\sum_{i=N+1}^{N+L'}\g_i(0)\Bigr)\a_{N+L'}(y)\,,
\end{equation*}
so~$A_N(3^{-N}\cdot)$ is a random scalar multiple of the cutoff coefficient~$\a_{N+L'}$. In the coordinates~$y$, the rescaled function~$u_N(3^{-N}\cdot)$ solves~$-\nabla\cdot(\a_{N+L'}\nabla u_N(3^{-N}\cdot))=\nabla\cdot g$, where~$g$ is the vector field constructed in Lemma~\ref{mfd:lem-primitive} for the source~$f(3^{-N}\cdot)$ extended by zero, divided by~$3^{2N}$ and by this random scalar. Moreover,~$u_N(3^{-N}\cdot)=\varphi(3^{-N}\cdot)\in C^{1,\nf12}$ on the boundary. Since~$N+L'\geq N$, the estimate~\eqref{e.energy.density.estimate} applies with~$m=N$; we write~$\mathcal L\coloneqq\mathcal L_{N+L'}(\alpha_M,N)$ for the minimal scale in that estimate.

\smallskip

\emph{Step 2.} We prove the Morrey bound~\eqref{eq:mfd-4} at the radii~$3^{-N}\leq r\leq3^{-\mathcal L}$. For~$x\in y_0+\cu_0$ and~$0\leq n\leq N-\mathcal L$, we have~$B_{3^{n-N}/2}(x)\cap(y_0+\cu_0)\subset3^{-N}\bigl((3^Nx+\cu_n)\cap(3^Ny_0+\cu_N)\bigr)$. Substituting the definitions of the random scalar and of~$\b_{N+L',N}$, the constants~$L'\tau^2$ cancel, and the random scalar times the average~$(\b_{N+L',N})_{3^Ny_0+\cu_N}$ equals~$\fint_{y_0+\cu_0}e^H\,dx$. We square~\eqref{e.energy.density.estimate}, multiply it by the volume, which is at most~$3^{dn}$, and return to the~$x$-coordinates. Using Lemma~\ref{mfd:lem-primitive} with~$R=3^N$ and bounding~$\fint_{y_0+\cu_0}e^H\,dx$ and its inverse by~$e^{\|H\|_{L^\infty(y_0+\cu_0)}}$, we obtain, for~$r=3^{n-N}/2$,
\begin{multline}
	\Gamma_N(u_N)(B_r(x)\cap(y_0+\cu_0))
	\\
	\leq Cr^{d-2+2\alpha_M}\Bigl(\Gamma_N(u_N)(y_0+\cu_0)+e^{\|H\|_{L^\infty(y_0+\cu_0)}}\bigl(\|f\|_\infty^2+\|\varphi\|_{C^2}^2\bigr)\Bigr)\,.
	\label{eq:mfd-macro}
\end{multline}

\smallskip

For balls~$B_r(x)$ meeting~$\partial(y_0+\cu_0)$, the energy estimate used here is the boundary version obtained in Step~6 of the proof of Proposition~\ref{p.Holder.regularity}, from Lemma~\ref{l.centered.energy.good.scale} at its nearest good anchor, together with its cutoff extension in Proposition~\ref{p.cutoff.Holder.regularity}. The minimal scale~$\mathcal L$ constructed in Propositions~\ref{p.Holder.regularity} and~\ref{p.cutoff.Holder.regularity} is independent of~$x$. Thus the same calculation gives~\eqref{eq:mfd-macro} near the boundary, with the scalar cancellation of Step~2.

\smallskip

The factor~$e^{\|H\|_{L^\infty(y_0+\cu_0)}}$ has moments bounded uniformly in~$N$ and~$L'$ by Lemma~\ref{mfd:lem-infrared}.

\smallskip

\emph{Step 3.} We treat the remaining large radii~$3^{-\mathcal L}<r\leq1$ and bound the total energy. For these radii,~$\Gamma_N(u_N)(B_r\cap(y_0+\cu_0))\leq\Gamma_N(u_N)(y_0+\cu_0)\leq3^{t_1\mathcal L}r^{t_1}\Gamma_N(u_N)(y_0+\cu_0)$. Together with~\eqref{eq:mfd-macro}, this gives the bound~$K_N(\|f\|_\infty+\|\varphi\|_{C^2})^2r^{t_1}$ at every radius~$r\geq3^{-N}$, also when~$\mathcal L>N$. For an order~$p$ in the list and for the larger orders used in H\"older's inequality,~\eqref{e.cutoff.Holder.minimal.scale} gives~$\E3^{pt_1\mathcal L}<\infty$ as soon as
\begin{equation}
	\frac{(1-\alpha_M)^2}{C\gmcdelta^2|\log\gmcdelta|}>pt_1\log3\,,
	\label{eq:mfd-tail-rate}
\end{equation}
which is a condition on~$\gmcdelta$ alone once~$\alpha_M$,~$p$ and~$t_1$ are fixed. The total energy has moments: by~\eqref{eq:mfd-quadratic-response},~\eqref{eq:mfd-1} and~\eqref{eq:mfd-2},
\begin{equation*}
	\Gamma_N(u_N)(y_0+\cu_0)\leq2\Lambda_{N,y_0+\cu_0}(\varphi)+2\langle f,G_N^{y_0+\cu_0}f\rangle\leq2\Lambda_{N,y_0+\cu_0}(\varphi)+2K_N\|f\|_{L^2}^2\,,
\end{equation*}
uniformly in the truncation, since the truncated coefficient satisfies the same estimates.

\smallskip

\emph{Step 4.} We prove the Morrey bound for~$r\leq3^{-N}$ on balls away from~$\partial(y_0+\cu_0)$. Let~$p_1<\infty$ satisfy~$d-2d/p_1>t_1$ and~$1-d/p_1>\alpha$, and let~$r_{N,*}$ be~$3^{-N}(1+D_N)^{-1}$ times a small constant depending on~$d$ and~$p_1$, chosen below. For~$x_0\in y_0+\cu_0$ with~$B_{2r_{N,*}}(x_0)\subset y_0+\cu_0$, the oscillation of~$\log A_N$ on~$B_{2r_{N,*}}(x_0)$ is at most~$2r_{N,*}D_N/3^{-N}$, hence at most twice the small constant. Dividing the equation by~$A_N(x_0)$ gives an equation with coefficient~$A_N/A_N(x_0)$, whose logarithm is bounded by twice the small constant on the ball, and source~$f/A_N(x_0)$. By~\citet[Theorem~2, equation~(48)]{Meyers1963}, the admissible gradient exponent tends to infinity as the ellipticity ratio tends to one. Choose the small constant so that the admissible exponent exceeds~$p_1$, subtract the mean of~$u_N$, and apply that estimate with initial gradient exponent~$2$ and a finite source exponent greater than~$d$. After scaling to the unit ball, this gives
\begin{equation*}
	\Bigl(\fint_{B_{r_{N,*}}(x_0)}|\nabla u_N|^{p_1}\Bigr)^{\nf1{p_1}}\leq C\Bigl(\fint_{B_{2r_{N,*}}(x_0)}|\nabla u_N|^2\Bigr)^{\nf12}+Cr_{N,*}A_N(x_0)^{-1}\|f\|_\infty\,.
\end{equation*}
H\"older's inequality on~$B_r(x_0)$ and the fact that~$A_N$ is within a fixed factor of~$A_N(x_0)$ on~$B_{2r_{N,*}}(x_0)$ yield, for every~$r\leq r_{N,*}$,
\begin{equation}
	\Gamma_N(u_N)(B_r)\leq C\Bigl(\frac r{r_{N,*}}\Bigr)^{d-2d/p_1}\Gamma_N(u_N)(B_{2r_{N,*}})
	+C\|A_N^{-1}\|_{L^\infty(y_0+\cu_0)}\|f\|_\infty^2r_{N,*}^{2+\nf{2d}{p_1}}r^{d-2d/p_1}\,.
	\label{eq:mfd-micro}
\end{equation}
We divide by~$r^t$. Since~$d-2d/p_1>t_1>t$, the right side is largest at~$r=r_{N,*}$, so
\begin{equation*}
	\Gamma_N(u_N)(B_r)r^{-t}\leq C(1+D_N)^t3^{tN}\Gamma_N(u_N)(B_{C3^{-N}})+C\|A_N^{-1}\|_{L^\infty(y_0+\cu_0)}\|f\|_\infty^23^{-(d+2-t)N}\,.
\end{equation*}
For~$r_{N,*}\leq r\leq3^{-N}$, we use~$\Gamma_N(u_N)(B_r)\leq\Gamma_N(u_N)(B_{C3^{-N}})\leq(3^{-N}/r_{N,*})^tr^t3^{tN}\Gamma_N(u_N)(B_{C3^{-N}})$. Together with Step~3 at~$r=C3^{-N}$, this gives, for every~$r\leq3^{-N}$,
\begin{equation*}
	\Gamma_N(u_N)(B_r)\leq C(1+D_N)^t3^{-(t_1-t)N}K_Nr^t
	+C\|A_N^{-1}\|_{L^\infty(y_0+\cu_0)}\|f\|^2_\infty3^{-(d+2-t)N}r^t\,.
\end{equation*}
By~\eqref{eq:mfd-6}, the positive powers~$3^{-(t_1-t)N}$ and~$3^{-(d+2-t)N}$ absorb the polynomially growing moments of~$D_N$ and the moments of~$\|A_N^{-1}\|_{L^\infty(y_0+\cu_0)}$, once~$\gmcdelta$ is small enough for the finitely many larger moment orders used in H\"older's inequality.

\smallskip

\emph{Step 5.} We prove the Morrey bound for~$r\leq3^{-N}$ near the faces of~$y_0+\cu_0$, where the Dirichlet condition holds. Let~$x_0$ be within distance~$2r_{N,*}$ of a face, say~$\{x_d=0\}$. The function~$u_N-\varphi$ vanishes on the face and solves~$-\nabla\cdot(A_N\nabla(u_N-\varphi))=f+\nabla\cdot(A_N\nabla\varphi)$, where~$\|\nabla\cdot(A_N\nabla\varphi)\|_\infty\leq CA_N(x_0)(\|D^2\varphi\|_\infty+r_{N,*}^{-1}\|\nabla\varphi\|_\infty)$ on the ball around~$x_0$ on which~$\log A_N$ has small oscillation, as in Step~4. We reflect~$u_N-\varphi$ oddly,~$A_N$ evenly and the source oddly across the face; the reflected function solves the reflected equation, with no additional source concentrated on the face. Indeed, for a compactly supported test function on a ball centered on the face, the difference of the test function and its reflection, restricted to~$y_0+\cu_0$, lies in~$H^1$, has zero trace on the face, and vanishes near the other faces, so it is an admissible test function for the original equation. The map that folds the reflected domain back onto~$y_0+\cu_0$ is~$1$-Lipschitz, so the oscillation of the logarithm of the reflected coefficient on a ball~$B$ centered on the face is at most that of~$\log A_N$ on~$B\cap\overline{y_0+\cu_0}$. Hence~\eqref{eq:mfd-micro} applies to the reflected function and coefficient, since the proof of~\eqref{eq:mfd-micro} uses only the boundedness and the small oscillation of the coefficient, not its differentiability across the face. Since the energy measure of the reflected function on a symmetric set~$B$ is twice the energy measure of~$u_N-\varphi$ on~$B\cap(y_0+\cu_0)$, and
\begin{equation*}
	\Gamma_N(u_N-\varphi)(B_{C3^{-N}}\cap(y_0+\cu_0))
	\leq2\Gamma_N(u_N)(B_{C3^{-N}}\cap(y_0+\cu_0))
	+C\|A_N\|_{L^\infty(y_0+\cu_0)}\|\nabla\varphi\|_\infty^23^{-dN}\,,
\end{equation*}
we obtain
\begin{align*}
	\Gamma_N(u_N)(B_r\cap(y_0+\cu_0))
	&\leq C(r/r_{N,*})^{d-2d/p_1}\Gamma_N(u_N)(B_{C3^{-N}}\cap(y_0+\cu_0))\notag\\
	&\qquad+C\|A_N^{-1}\|_{L^\infty(y_0+\cu_0)}\|f\|_\infty^2r_{N,*}^{2+\nf{2d}{p_1}}r^{d-2d/p_1}\notag\\
	&\qquad+C\|A_N\|_{L^\infty(y_0+\cu_0)}\|\varphi\|_{C^2}^2(1+D_N)^dr_{N,*}^{\nf{2d}{p_1}}r^{d-2d/p_1}\,.
\end{align*}
Dividing the last term by~$r^t$ and maximizing at~$r=r_{N,*}$ gives~$C\|A_N\|_{L^\infty(y_0+\cu_0)}\|\varphi\|_{C^2}^23^{-(d-t)N}(1+D_N)^t$. The positive power~$3^{-(d-t)N}$ absorbs the moments of~$\|A_N\|_{L^\infty(y_0+\cu_0)}(1+D_N)^t$ as in Step~4. Near an edge or a corner of~$y_0+\cu_0$, we reflect across each face meeting~$B_{2r_{N,*}}(x_0)$. On a set invariant under these reflections, the energy of the reflected solution is two to the power of the number of these faces times the energy of the solution in the part of the set inside~$y_0+\cu_0$. Divide by this factor to obtain the estimates on the original cube.

\smallskip

\emph{Step 6.} We remove the truncation. By Lemma~\ref{mfd:lem-infrared}, the truncated coefficients converge uniformly on~$\overline{y_0+\cu_0}$ to the untruncated~$A_N$ almost surely. For fixed~$N$, they are therefore bounded away from zero once the truncation is sufficiently large. Testing the difference of the two weak equations with the difference of the solutions, which lies in~$H^1_0(y_0+\cu_0)$, bounds its energy by the uniform distance between the coefficients times the product of the~$L^2$ norms of the gradient of the untruncated solution and of the gradient of the difference. The lower bound on the coefficients then shows that the gradients converge in~$L^2(y_0+\cu_0)$. Consequently, for every~$x$ and~$r$, the truncated energies converge to~$\Gamma_N(u_N)(B_r\cap(y_0+\cu_0))$. By Steps~2--5, these energies satisfy~\eqref{eq:mfd-4} with random constants whose~$L^p$ norms are bounded uniformly in~$N$ and in the truncation, for every order~$p$ in the list. The limit inferior of these random constants gives a constant for the untruncated estimate, and Fatou's lemma preserves these moment bounds.

\smallskip

\emph{Step 7.} We prove the H\"older bounds. At radii above~$3^{-N}$,~\eqref{e.Holder.estimate.boxes.local} is a Campanato bound at every~$x\in y_0+\cu_0$, with exponent~$\alpha_M>\alpha$ and the same source norm; the large radii~$3^{-\mathcal L}<r\leq1$ are handled as in Step~3, and the~$L^2$ oscillation on~$y_0+\cu_0$ is bounded by the energy through~\eqref{eq:mfd-1}. At radii below~$3^{-N}$, the~$W^{1,p_1}$ bound of Step~4 and Morrey's inequality bound the oscillation on~$B_{r_{N,*}}$. After division by~$r^\alpha$, the energy on the ball of radius~$C3^{-N}$, the volume source and the boundary datum contribute, respectively,
\begin{align*}
	&CK_N^{\nf12}\|A_N^{-1}\|_{L^\infty(y_0+\cu_0)}^{\nf12}
	3^{-((t_1-d+2)/2-\alpha)N}(1+D_N)^C\,,\\
	&C\|A_N^{-1}\|_{L^\infty(y_0+\cu_0)}\|f\|_\infty
	3^{-(2-\alpha)N}(1+D_N)^C\,,\\
	&C\|\varphi\|_{C^2}
	3^{-(1-\alpha)N}(1+D_N)^C\,.
\end{align*}
All three powers of~$3^{-N}$ are positive by the choice of the parameters. Chaining balls of radius~$r_{N,*}$ inside a cube of side~$3^{-N}$ introduces another fixed power of~$1+D_N$. Above the smallest field scale, telescoping the cube averages produces the geometric sum~$\sum_{r_1\leq r}r_1^{\alpha_M}\leq Cr^{\alpha_M}$ over triadic radii~$r_1$. With the same absorption of moments as in Step~4, this proves the H\"older bounds. Approximate a general bounded source by smooth sources. Remove the truncation as in Step~6: the equicontinuity just proved and the convergence in~$L^2$ give uniform convergence of the truncated solutions.
\end{proof}

\smallskip

\begin{remark}[Morrey and H\"older bounds on larger cubes]
\label{mfd:rem-large-cube-growth}
For the coefficient~$A_N$, both conclusions of Proposition~\ref{mfd:prop-growth} hold on every fixed cube~$w+R\cu_0$, with moment constants depending on~$w+R\cu_0$. The disorder threshold depends only on~$d$,~$t$,~$\alpha$ and the prescribed moment orders. To see this, choose~$j\geq1$ so that~$R\leq3^j$. With~$\Theta_{k,z}$ as in Section~\ref{mfd:sec-environment}, the scale shift~$(\Theta_{-j,0}\omega)_\ell(y)=\g_{\ell+j}^\omega(3^jy)$ preserves the law of the fields. Splitting the series defining~$H$ at~$j$ gives
\begin{equation}
	A_N^\omega(3^jy)=s_N(3^j\cu_0)A_{N+j}^{\Theta_{-j,0}\omega}(y)\,.
	\label{eq:mfd-large-cube-transfer}
\end{equation}
By~\eqref{e.annealed.ordering}, both~$s_N(3^j\cu_0)$ and~$s_N(3^j\cu_0)^{-1}$ are at most~$\exp(j\tau^2+|\sum_{i=1}^j\g_i(0)|)$, uniformly in~$N$. This exponential has every finite moment by~\ref{a.g2}.

\smallskip

If~$u_N$ solves the Dirichlet problem on~$w+R\cu_0$, then~$v(y)=u_N(3^jy)$ solves it on~$3^{-j}(w+R\cu_0)$, with coefficient~$A_{N+j}^{\Theta_{-j,0}\omega}$, source~$3^{2j}s_N(3^j\cu_0)^{-1}f(3^j\,\cdot)$ and boundary datum~$\varphi(3^j\,\cdot)$. The weak equation follows by changing variables in its test integrals, and the energy satisfies
\begin{equation*}
	\Gamma_N^\omega(u_N)(B_r(x)\cap(w+R\cu_0))
	=s_N(3^j\cu_0)3^{j(d-2)}
	\Gamma_{N+j}^{\Theta_{-j,0}\omega}(v)
	(B_{r/3^j}(x/3^j)\cap3^{-j}(w+R\cu_0))\,.
\end{equation*}
Apply Proposition~\ref{mfd:prop-growth} on the cube~$3^{-j}(w+R\cu_0)$ of side at most one. Since~$t>d-2$ and~$3^j\geq1$, the radius factor~$3^{j(d-2-t)}$ is at most one. Also~$\|\varphi(3^j\,\cdot)\|_{C^2}\leq(1+3^j+3^{2j})\|\varphi\|_{C^2}$, and~$s_N(3^j\cu_0)\max\{1,s_N(3^j\cu_0)^{-1}\}^2=\max\{s_N(3^j\cu_0),s_N(3^j\cu_0)^{-1}\}$. The energy and H\"older bounds are consequently transferred by multiplying the random constant~$K_{N+j}$ of this application by~$(1+3^j+3^{2j})^2\exp(j\tau^2+|\sum_{i=1}^j\g_i(0)|)$. For each prescribed order~$p$, use Proposition~\ref{mfd:prop-growth} on~$3^{-j}(w+R\cu_0)$ at order~$2p$ before choosing the disorder threshold. H\"older's inequality and preservation of the law of the fields then bound the product in~$L^p$, uniformly in~$N$. This proves both the Morrey and H\"older estimates on~$w+R\cu_0$ for~$A_N$; their moment constants may depend on the fixed position of~$w+R\cu_0$.
\end{remark}

\subsubsection{Approximating the Neumann datum by a volume source}

We approximate the functional~$v\mapsto\int_{(\tfrac12,\dots,\tfrac12)+\cu_0}p\cdot\nabla v$, which encodes the Neumann datum~$p\cdot\mathbf n$, by a bounded volume source near the faces.

\smallskip

Let~$p\in\Rd$ with~$|p|=1$, and let
\begin{equation*}
	L_p(v)\coloneqq\int_{(\tfrac12,\dots,\tfrac12)+\cu_0}p\cdot\nabla v\,dx=\int_{\partial((\tfrac12,\dots,\tfrac12)+\cu_0)}(p\cdot\mathbf n)\operatorname{Tr}v\,d\mathcal H^{d-1}\,,
\end{equation*}
where~$\mathbf n$ is the outer normal. For~$0<\ep<\nf18$, let~$f_\ep$ be the sum over~$i=1,\dots,d$ of~$p_i$ times the difference of two nonnegative smooth functions, constant in the tangential coordinates and of integral one in the normal coordinate, the first supported in the layer~$\{\ep<1-x_i<2\ep\}$ near the face~$\{x_i=1\}$ and the second in the layer~$\{\ep<x_i<2\ep\}$ near the face~$\{x_i=0\}$, both obtained from one fixed profile by dilation by~$\ep$, and let
\begin{equation}
	L_{p,\ep}(v)\coloneqq\int_{(\tfrac12,\dots,\tfrac12)+\cu_0}f_\ep v\,dx\,.
	\label{eq:mfd-bump}
\end{equation}
Each summand of~$f_\ep$ has integral zero,~$f_\ep$ belongs to~$C^\infty(\overline{(\tfrac12,\dots,\tfrac12)+\cu_0})$ with~$\|f_\ep\|_\infty\leq C\ep^{-1}$, and~$L_{p,\ep}(1)=0$.

Let~$v_N,v_{N,\ep}$ be the mean-zero Riesz representers of~$L_p,L_{p,\ep}$ with respect to~$\mathscr E_N^{(\tfrac12,\dots,\tfrac12)+\cu_0}$, so that~$\mathscr E_N^{(\tfrac12,\dots,\tfrac12)+\cu_0}(v_N)=p\cdot \a_*^{-1}((\tfrac12,\dots,\tfrac12)+\cu_0;A_N)p$.

\begin{lemma}[Error of the volume-source approximation of the Neumann datum]
\label{mfd:lem-neumann-error}
\begin{itemize}
\item For every~$v\in H^1((\tfrac12,\dots,\tfrac12)+\cu_0)$,
\begin{equation}
	|L_p(v)-L_{p,\ep}(v)|
	\leq C\ep^{\nf14}\|v-(v)_{(\tfrac12,\dots,\tfrac12)+\cu_0}\|_{H^{\nf34}((\tfrac12,\dots,\tfrac12)+\cu_0)}\,.
	\label{eq:mfd-7}
\end{equation}
\item With~$Q_N$ from~\eqref{eq:mfd-1}, this gives
\begin{equation}
	\mathscr E_N^{(\tfrac12,\dots,\tfrac12)+\cu_0}(v_N-v_{N,\ep})\leq CQ_N\ep^{\nf12}\,.
	\label{eq:mfd-8}
\end{equation}
\item For every moment order~$q$ in the finite list fixed at the beginning of Section~\ref{sec:mfd-uniqueness}, if~$\gmcdelta$ is small enough, then
\begin{equation}
	\|\mathscr E_N^{(\tfrac12,\dots,\tfrac12)+\cu_0}(v_N)-\mathscr E_N^{(\tfrac12,\dots,\tfrac12)+\cu_0}(v_{N,\ep})\|_{L^q}\leq C_q\ep^{\nf14}
	\quad\text{and}\quad \sup_{N,\ep}\|\mathscr E_N^{(\tfrac12,\dots,\tfrac12)+\cu_0}(v_{N,\ep})\|_{L^q}<\infty\,.
	\label{eq:mfd-9}
\end{equation}
\end{itemize}
\end{lemma}
\begin{proof}
Subtract~$(v)_{(\tfrac12,\dots,\tfrac12)+\cu_0}$ and extend to~$H^{\nf34}(\Rd)$. Fourier transformation in the normal variable and Cauchy--Schwarz give, on each face,
\begin{equation*}
	\|v(s,\cdot)-v(0,\cdot)\|_{L^2(\mathbb R^{d-1})}^2
	\leq\|v\|_{H^{\nf34}(\Rd)}^2
	\int_{\mathbb R}\frac{\min\{4,s^2\vartheta^2\}}
	{(1+\vartheta^2)^{\nf34}}\,d\vartheta
	\leq Cs^{\nf12}\|v\|_{H^{\nf34}(\Rd)}^2\,.
\end{equation*}
Integrate against the bump profile and sum over the faces to get~\eqref{eq:mfd-7}. Test the difference of the two weak equations with~$v_N-v_{N,\ep}$ and use coercivity to obtain~\eqref{eq:mfd-8}. The reverse triangle inequality for the energy norm then gives
\begin{equation*}
	|\mathscr E_N^{(\tfrac12,\dots,\tfrac12)+\cu_0}(v_N)-\mathscr E_N^{(\tfrac12,\dots,\tfrac12)+\cu_0}(v_{N,\ep})|
	\leq CQ_N^{\nf12}\ep^{\nf14}\mathscr E_N^{(\tfrac12,\dots,\tfrac12)+\cu_0}(v_N)^{\nf12}
	+CQ_N\ep^{\nf12}\,.
\end{equation*}
The energy~$\mathscr E_N^{(\tfrac12,\dots,\tfrac12)+\cu_0}(v_N)=p\cdot \a_*^{-1}((\tfrac12,\dots,\tfrac12)+\cu_0;A_N)p$ has the required moments uniformly in~$N$ by Lemma~\ref{mfd:lem-cell-ellipticity}. H\"older's inequality at the doubled orders proves~\eqref{eq:mfd-9}.

\smallskip

We also record the transfer used on a fixed physical cube~$w+R\bigl((\tfrac12,\dots,\tfrac12)+\cu_0\bigr)$. For~$p\neq0$, let~$f_h^{p/|p|}$ be the source~$f_\ep$ of~\eqref{eq:mfd-bump} for the vector~$p/|p|$ and~$\ep=h$. Define
\begin{align*}
	&f_{w,R,h}(w+Ry)\coloneqq\frac{|p|}{R}f_h^{p/|p|}(y)\,,\\
	&L_{w,R,p}(u)\coloneqq\int_{w+R\bigl((\tfrac12,\dots,\tfrac12)+\cu_0\bigr)}p\cdot\nabla u\,,\\
	&L_{w,R,p,h}(u)\coloneqq\int_{w+R\bigl((\tfrac12,\dots,\tfrac12)+\cu_0\bigr)}f_{w,R,h}u\,.
\end{align*}
For~$p=0$ take both functionals and the source to be zero. With~$\widetilde u(y)\coloneqq u(w+Ry)$, change of variables gives
\begin{equation*}
	L_{w,R,p}(u)=|p|R^{d-1}L_{p/|p|}(\widetilde u),\qquad
	L_{w,R,p,h}(u)=|p|R^{d-1}L_{p/|p|,h}(\widetilde u)\,.
\end{equation*}
In particular,~$\int_{w+R\bigl((\tfrac12,\dots,\tfrac12)+\cu_0\bigr)}f_{w,R,h}=0$ and~$\|f_{w,R,h}\|_\infty\leq C_{w,R,p}h^{-1}$, and the fixed dilation of the fractional norm in~\eqref{eq:mfd-7} gives
\begin{equation*}
	|L_{w,R,p}(u)-L_{w,R,p,h}(u)|
	\leq C_{w,R,p}h^{\nf14}\|u-(u)_{w+R\bigl((\tfrac12,\dots,\tfrac12)+\cu_0\bigr)}\|_{H^{\nf34}(w+R\bigl((\tfrac12,\dots,\tfrac12)+\cu_0\bigr))}\,.
\end{equation*}

\smallskip

The coercivity bound needed here also transfers to fixed larger cubes. For~$R\leq1$ use Lemma~\ref{mfd:lem-coercivity}. For~$R>1$, choose~$k\geq1$ with~$r=R/3^k\in(1/3,1]$ and use the change of variables~$x=z+3^ky$ from the cube~$r\cu_0$, where~$z$ is the center of~$w+R\bigl((\tfrac12,\dots,\tfrac12)+\cu_0\bigr)$. The environment~$(\Theta_{-k,z}\omega)_j(y)=\g^\omega_{j+k}(z+3^ky)$ has the original law. The coefficient at cutoff~$N+k$ in this environment is bounded above by~$F_{w,R}(\omega)A_N(z+3^ky)$, and the same bound holds with~$A_N^0$ on the right, where we may take
\begin{multline*}
	F_{w,R}\coloneqq\exp\Bigl(k\tau^2+\sum_{i=1}^k\|\g^\omega_i\|_{C(\overline{w+R\bigl((\tfrac12,\dots,\tfrac12)+\cu_0\bigr)})}\\
	+\|H(\Theta_{-k,z}\omega)\|_{C(\overline{r\cu_0})}
	+\|H(\omega)\|_{C(\overline{w+R\bigl((\tfrac12,\dots,\tfrac12)+\cu_0\bigr)})}\Bigr)\,.
\end{multline*}
The normalization ratio is bounded by~\eqref{e.annealed.ordering}, namely~$\ahom_N\leq e^{2k\tau^2}\ahom_{N+k}$; the other factors come from the fields~$\g_1,\dots,\g_k$ and from the two restrictions of~$H$. Each has every fixed exponential moment on its fixed compact set by~\eqref{e.coefficient.field.regularity} and Lemma~\ref{mfd:lem-infrared}. Pull back the mean-zero function, and likewise a function with zero trace (the~$H^1_0$ clause of Lemma~\ref{mfd:lem-coercivity}); use that lemma on~$r\cu_0$ and scale the energy by~$3^{k(d-2)}$. The fixed scaling of the fractional norm gives a physical coercivity constant
\begin{equation*}
	Q_{w+R((\tfrac12,\dots,\tfrac12)+\cu_0),N}\leq C_{w,R}F_{w,R}\,Q_{r\cu_0,N+k}(\Theta_{-k,z}\omega)\,.
\end{equation*}
The constants in the moment bounds may depend on~$w$ and~$R$; the disorder threshold is the one for the finite doubled moment orders on~$r\cu_0$. H\"older's inequality and preservation of the law of the fields give the required uniform-in-~$N$ moments, for both~$A_N$ and~$A_N^0$. For a triadic~$R$,~$r=1$. Only a first-moment bound for~$Q_{r\cu_0,N+k}$ is needed for a fixed polynomial tail: the union bound at~$\sqrt S$ gives~$\P[Q_{w+R((\tfrac12,\dots,\tfrac12)+\cu_0),N}>S]\leq M_{w,R}S^{-1/2}$ after absorbing~$C_{w,R}$ into~$F_{w,R}$. This uses no independence between~$F_{w,R}$ and~$Q_{r\cu_0,N+k}(\Theta_{-k,z}\omega)$.

\smallskip

Let~$Y_N$ and~$Y_{N,h}$ be the energies of the mean-zero Riesz representers of~$L_{w,R,p}$ and~$L_{w,R,p,h}$. These representers and the preceding estimate give
\begin{equation*}
	|Y_N-Y_{N,h}|
	\leq 2C_{w,R,p}h^{\nf14}Q_{w+R((\tfrac12,\dots,\tfrac12)+\cu_0),N}^{\nf12}Y_N^{\nf12}
	+C_{w,R,p}^2h^{\nf12}Q_{w+R((\tfrac12,\dots,\tfrac12)+\cu_0),N}\,.
\end{equation*}
The trace bound~$|L_{w,R,p}(u)|\leq C_{w,R,p}\|u-(u)_{w+R\bigl((\tfrac12,\dots,\tfrac12)+\cu_0\bigr)}\|_{H^{\nf34}(w+R\bigl((\tfrac12,\dots,\tfrac12)+\cu_0\bigr))}$ gives~$Y_N\leq C_{w,R,p}Q_{w+R((\tfrac12,\dots,\tfrac12)+\cu_0),N}$ by the variational formula. Consequently, for every moment order~$q$ in the finite list,
\begin{equation*}
	\sup_N\|Y_N-Y_{N,h}\|_{L^q}\leq C_{q,w,R,p}h^{\nf14}\quad\text{and}\quad
	\sup_{N,\,0<h<1/8}\|Y_{N,h}\|_{L^q}<\infty\,.
\end{equation*}
These are the fixed-cube versions of~\eqref{eq:mfd-9}. For the identity coefficient in place of~$A_N$, the energy~$Y_N$ on~$w+R\bigl((\tfrac12,\dots,\tfrac12)+\cu_0\bigr)$ is~$R^d\,p\cdot p$.
\end{proof}

\subsubsection{Reflection across the faces}

Let~$A$ be a scalar coefficient, and let~$u$ be a solution with zero conormal derivative and bounded source~$f$. Reflect~$u$ and~$f$ evenly across a set of faces. This gives a solution~$u\circ T$ of the equation with coefficient~$A\circ T$ and source~$f\circ T$, where~$T$ folds the reflected domain onto the original cube. The reflected equation has no additional source concentrated on the folding planes. On a set invariant under these reflections, the energy of~$u\circ T$ is two to the power of the number of reflected faces times the energy of~$u$ in the part of the set inside the original cube.

\begin{lemma}[Effect of folding on~$\mathcal E_{s,\infty,2}$]
\label{mfd:thm-fold}
Let~$T$ fold the cube~$z+\cu_m$ across coordinate planes through~$z$, as above. For every continuous positive coefficient~$A$ on its closure, scalar~$A_0>0$,~$s\in(0,\nf12)$ and~$n\leq m$,
\begin{equation}
	\mathcal E_{s,\infty,2}(z+\cu_m,n;A\circ T,A_0)^2
	\leq\left(1+\frac{3d}{3^{1-2s}-1}\right)
	\mathcal E_{s,\infty,2}(z+\cu_m,n;A,A_0)^2\,.
	\label{eq:mfd-11}
\end{equation}
\end{lemma}
\begin{proof}
On a cube meeting no folding plane,~$T$ reflects some of the coordinates, so changing variables by~$T$ shows that, for each~$n$, the quantity~$\sup_{|e|=1}J(\,\cdot\,,e,e;(A\circ T)/A_0)$ on that cube is at most the maximum of~$\sup_{|e|=1}J(z'+\cu_n,e,e;A/A_0)$ over the cubes~$z'+\cu_n\subset z+\cu_m$ with~$z'\in z+3^n\Zd$. Subdivide each cube that meets a folding plane triadically, and repeat the subdivision on the subcubes that still meet one. After~$j$ subdivisions, the subcubes that still meet a folding plane fill at most a fraction~$d3^{-j}$ of the cube. By subadditivity~\eqref{e.subadditivity.J}, the corresponding maxima satisfy
\begin{multline*}
	\max_{\substack{z'\in z+3^n\Zd\\z'+\cu_n\subset z+\cu_m}}\sup_{|e|=1}J\bigl(z'+\cu_n,e,e;(A\circ T)/A_0\bigr)
	\leq\max_{\substack{z'\in z+3^n\Zd\\z'+\cu_n\subset z+\cu_m}}\sup_{|e|=1}J\bigl(z'+\cu_n,e,e;A/A_0\bigr)\\
	+\sum_{j\geq1}3d\,3^{-j}\max_{\substack{z'\in z+3^{n-j}\Zd\\z'+\cu_{n-j}\subset z+\cu_m}}\sup_{|e|=1}J\bigl(z'+\cu_{n-j},e,e;A/A_0\bigr)\,.
\end{multline*}
Insert this inequality into Definition~\ref{d.mathcal.E} and sum the geometric series~$\sum_{j\geq1}3^{-(1-2s)j}$ to obtain~\eqref{eq:mfd-11}.
\end{proof}

\begin{proposition}[Energy decay for the folded coefficient]
\label{mfd:prop-folded-iteration}
There are dimensional constants~$C<\infty$ and~$M\geq1$ with the following property. Let~$\gmcdelta\leq C^{-1}$ and~$\gamma\in[1/2,1)$ satisfy~$1-\gamma\geq CM\gmcdelta|\log\gmcdelta|^{1/2}$, let~$m,L\in\N_0$ with~$m\leq L$, and let~$T$ fold the cube~$z+\cu_m$ across coordinate planes through~$z$. Let~$\gamma_*\coloneqq1-(1-\gamma)/M$. Define~$\mathcal L_L(\gamma,m)$ by~\eqref{e.L.def.reg} for the cube~$z+\cu_m$, with all parameter choices in~\eqref{e.choices.hepk.Holder.proof} evaluated at~$\gamma_*$; in particular,~$\lambda_*=\lambda/M$ and~$\ep_*=\ep/\sqrt M$, where~$\lambda$ and~$\ep$ denote the choices at~$\gamma$.
\begin{itemize}
\item If~$u\in H^1(z+\cu_m)$ solves~$-\nabla\cdot((\a_L\circ T)\nabla u)=\nabla\cdot\tilde g$, with~$\tilde g\in C^{0,\nf12}(z+\cu_m;\Rd)$, then, for~$0\leq n\leq m-\mathcal L_L(\gamma,m)$,
\begin{multline}
	\|(\a_L\circ T)^{\nf12}\nabla u\|_{\underline L^2(z+\cu_n)}
	\\
	\leq C3^{(1-\gamma)(m-n)}
	\left(\|(\a_L\circ T)^{\nf12}\nabla u\|_{\underline L^2(z+\cu_m)}
	+(\b_{L,m})_{z+\cu_m}^{-\nf12}3^{\frac{1}{2}m}[\tilde g]_{W^{\nf12,\infty}(z+\cu_m)}\right)\,.
	\label{eq:mfd-folded-energy}
\end{multline}
\item For every~$h\geq0$, this minimal scale satisfies
\begin{equation*}
	\P[\mathcal L_L(\gamma,m)>h]
	\leq C\exp\left(-\frac{(1-\gamma)^2(h-C)_+}
	{CM^2\gmcdelta^2|\log\gmcdelta|}\right)\,,
\end{equation*}
uniformly in~$L$,~$m$ and the deterministic translate~$z$.
\end{itemize}
\end{proposition}
\begin{proof}
Apply the argument proving~\eqref{e.energy.density.estimate} at~$z$, using the original averages~$(\b_{L,j+2})_{z+\cu_{j+2}}$. The harmonic approximation and excess decay use the good event through~$\mathcal E_{s,\infty,2}$, as in the proof of Lemma~\ref{l.cutoff.regularity.good.scale.estimates}. By~\eqref{eq:mfd-11}, folding multiplies this quantity by at most a dimensional constant. The bounds with~$q=1$ at the larger exponents follow from the bound with~$q=2$ at~$s_0=\nf1{32}$ by weighted Cauchy--Schwarz. Choose a dimensional~$M\geq1$ large enough to absorb these factors, and evaluate the parameter choices~\eqref{e.choices.hepk.Holder.proof} at~$\gamma_*=1-(1-\gamma)/M$, so that~$\lambda_*=\lambda/M$ and~$\ep_*=\ep/\sqrt M$, as in the statement. This uses the same good events at the stricter parameter~$\ep_*$. The error sum, the bound on bad scales and the ratios of the original averages are the bounds~\eqref{e.ep.j.z.sum},~\eqref{e.Bz.bound} and~\eqref{e.ratio.of.bs} with these parameter choices. They supply the iteration of Lemma~\ref{l.iteration.lemma.GMC};~\eqref{e.bound.Lambdas.by.Es} supplies its coarse-grained ellipticity inputs. All its cubes are centered at~$z$ and invariant under the reflections. The stricter choices give the asserted estimate with exponent~$\gamma$.

\smallskip

We record the tail of the minimal scale with this reduced tolerance. Each of the two scales in~\eqref{e.L.def.reg} has tail length bounded by a constant times
\begin{equation*}
	s_0^{-7}\gmcdelta^2|\log\gmcdelta|\lambda_*^{-2}
	+s_0^{-8}\ep_*^{-2}\lambda_*^{-1}
	\gmcdelta^2|\log\gmcdelta|
	\leq \frac{CM^2\gmcdelta^2|\log\gmcdelta|}{(1-\gamma)^2}\,.
\end{equation*}
Indeed, both~$\lambda_*^{-2}$ and~$\ep_*^{-2}\lambda_*^{-1}$ gain the factor~$M^2$. Equations~\eqref{e.M.lambda.m.bound} and~\eqref{e.N.theta.m.bound.goodscales}, or~\eqref{e.cutoff.regularity.error.stopping.tail} and~\eqref{e.cutoff.regularity.good.stopping.tail} for the cutoff version, therefore give
\begin{equation*}
	\P[\mathcal L_L(\gamma,m)>h]
	\leq C\exp\left(-\frac{(1-\gamma)^2(h-C)_+}
	{CM^2\gmcdelta^2|\log\gmcdelta|}\right)\,.
\end{equation*}
The constants are uniform in~$L$,~$m$ and~$z$. The hypothesis on~$1-\gamma$ places~$\gamma_*$ in the range of Proposition~\ref{p.cutoff.Holder.regularity}, because~$C$ may be taken at least as large as the constant in that proposition. These are parameter choices in the existing estimates; the definitions of the good events are unchanged.
\end{proof}

For a fixed cube~$w+R\cu_0$ and~$f\in L^\infty(w+R\cu_0)$ with~$\int_{w+R\cu_0}f=0$, let~$u_N\in H^1(w+R\cu_0)$ be the mean-zero solution of the Neumann problem
\begin{equation*}
	\int_{w+R\cu_0}A_N\nabla u_N\cdot\nabla v=\int_{w+R\cu_0}fv
	\,,
\end{equation*}
for every~$v\in H^1(w+R\cu_0)$.

\begin{proposition}[Bounded Neumann sources]
\label{mfd:cor-neumann-source}
Let~$\alpha\in(0,1)$ and~$t\in(d-1,d)$, and fix a finite list of moment orders. There is~$\gmcdelta_0(d,\alpha,t,\text{list})>0$ such that, if~$\gmcdelta\leq\gmcdelta_0$, then on every fixed cube~$w+R\cu_0$ there are random variables~$K_N$ with~$\sup_N\|K_N\|_{L^p}<\infty$ for every order~$p$ in the list, such that, almost surely, for every cutoff~$N$ and every mean-zero~$f\in L^\infty(w+R\cu_0)$, the solutions just defined satisfy
\begin{equation*}
	\|u_N\|_{C^\alpha(\overline{w+R\cu_0})}\leq K_N\|f\|_\infty\,,
\end{equation*}
and, for every~$x\in w+R\cu_0$ and~$0<r\leq1$,
\begin{equation*}
	\Gamma_N(u_N)(B_r(x)\cap(w+R\cu_0))\leq K_N\|f\|_\infty^2r^t\,.
\end{equation*}
\end{proposition}
\begin{proof}
We track the energy~$\Gamma_N(u_N)(B_r(x)\cap(w+R\cu_0))$ of small balls. Step~1 rescales the problem to the stationary coefficient. Step~2 joins each target ball to a ball of fixed radius by a chain of at most~$d+1$ applications of Proposition~\ref{mfd:prop-folded-iteration}, reflecting evenly across each face that the chain comes close to. The random factors in these applications depend on the target ball; the main work, in Steps~3 and~4, is to bound them, up to small negative powers of the radii, by random variables that do not depend on the target and whose moments are bounded uniformly in~$N$; this is done by summing over meshes of centers and root depths, and of reference averages. Together with the coercivity estimate~\eqref{eq:mfd-1} for the total energy, this gives the Morrey bound for~$r\geq3^{-N}$. Step~5 treats the radii below~$3^{-N}$ and proves the H\"older bound; the last paragraphs remove the truncation of~$H$ and pass to general fixed cubes.

\smallskip

We first work on~$(\tfrac12,\dots,\tfrac12)+\cu_0$. Choose~$\max\{t,d-2+2\alpha\}<t_1<t_0<d$ and~$0<\theta<d+2-t_0$. Fix all moment orders to be used below before choosing~$\gmcdelta$. These include the given list, the larger orders for products with coercivity constants and the auxiliary orders used for the two meshes. When this proposition is used in Proposition~\ref{mfd:prop-16}, include~$3p$; when it is used at order~$P$ and at order~$2P$ in Lemma~\ref{mfd:lem-band}, include~$3P$ and~$6P$, respectively. This is still a finite list. The auxiliary orders are~$\xi$ in Step~3 and~$q$ in Step~4, and the conditions on~$\gmcdelta$ imposed in those steps depend only on these orders and on quantities fixed before~$\gmcdelta$.

\smallskip

\emph{Step 1.} We specify the physical normalization. Temporarily replace~$H$ by~$H_{L'}\coloneqq\sum_{i=1}^{L'}(\g_i-\g_i(0))$, including~$L'=0$. Reindexing the stationary fields gives
\begin{equation*}
	A_N^{(L')}(y)=c_{N,L'}\a_{N+L'}(3^Ny),\qquad
	c_{N,L'}\coloneqq\ahom_N^{-1}
	\exp\left(L'\tau^2-\sum_{i=1}^{L'}\g_i(0)\right)\,.
\end{equation*}
Here the coefficient~$\a_{N+L'}$ on the right is taken in the reindexed environment. Choose a measurable representative of~$f$ bounded everywhere by its essential supremum, changing it on a null set if necessary; this preserves its integral and the functional~$v\mapsto\int fv$. Under~$X=3^Ny$ and~$U(X)\coloneqq3^Nu_N(3^{-N}X)$, its equation has source~$(c_{N,L'}3^N)^{-1}f(3^{-N}X)$. Even reflection replaces this source by its composition with the coordinate fold. On a reflection-invariant set, the stationary energy of~$U$ is~$2^{|I|}3^{dN}c_{N,L'}^{-1}$ times the physical energy in its intersection with~$(\tfrac12,\dots,\tfrac12)+\cu_0$, where~$I$ is the set of reflected coordinates.

\smallskip

Thus, in physical coordinates, the average in the source term of~\eqref{eq:mfd-folded-energy} is multiplied by~$c_{N,L'}$. At depth~$k$ and center~$z$ we use for this product the ball average
\begin{equation*}
	\bar s_{N,k}(z)\coloneqq\fint_{B_{3^{-k}/2}(z)}
	\frac{\ahom_{N-k}}{\ahom_N}
	\exp\left(H_{L'}(x)+\sum_{j=0}^{k-1}\g_{-j}(x)-k\tau^2\right)\,dx\,,
	\qquad 0\leq k\leq N\,.
\end{equation*}
Use the fixed shift of depth prescribed by the outer cube in the iteration: for a root (Step~2) of radius~$3^{-k_i}/2$ we use~$\bar s_{N,(k_i-\ell)_+}(z_i)$ for this product, where~$\ell$ is the integer fixed first in Step~2. Its averaging radius is comparable to~$3^{-k_i}/2$ with a fixed factor. We write~$\bar s_i\coloneqq \bar s_{N,(k_i-\ell)_+}(z_i)$. The factor~$c_{N,L'}$ has already been included.

\smallskip

\emph{Step 2.} We construct the reflected chains. Fix~$\ell$ first, and then fix~$L_*>10$ with~$L_*\geq2\sqrt d\,3^\ell$. Call the radii~$3^{-k}/2$ with~$k\in\Z$ half-triadic, and fix a half-triadic~$R_*>0$ with~$R_*<(100L_*)^{-1}$ and~$R_*\leq3^{-\ell}/2$. All three are fixed independently of~$N$ and~$L'$. At a center~$y$, let~$a_i\coloneqq\min\{y_i,1-y_i\}$ and let~$\bar y_i\in\{0,1\}$ be such that the face~$\{x_i=\bar y_i\}$ is nearer to~$y$. If~$I$ is the set of activated coordinates, then the current center~$z_I$ has coordinate~$\bar y_i$ for~$i\in I$ and coordinate~$y_i$ otherwise. Start with the least half-triadic radius~$r_1$ whose ball about the current center contains the target ball~$B_r(x)$. Stop if~$r_1\geq R_*/24$.

\smallskip

If an inactive coordinate~$i$ satisfies~$a_i\leq100L_*r_1$, then activate~$i$, project the center to the face~$\{x_i=\bar y_i\}$, and replace~$r_1$ by the least half-triadic radius at least~$r_1+a_i$. The old ball lies in the new one, and the radius increases by at most the factor~$3(1+100L_*)$. Otherwise choose the largest half-triadic~$r_2$ such that~$r_2\leq R_*$ and~$4L_*r_2\leq a_i$ for every~$i\notin I$. Then~$8r_1<r_2$. Reflect across every activated face and apply Proposition~\ref{mfd:prop-folded-iteration} with~$\gamma\coloneqq(t_0-d+2)/2$ between radii comparable to~$r_1$ and~$r_2$. The ball of radius~$r_2$ about the current center is the \textbf{root} of this application; for the~$i$th application we write~$z_i$ for its center and~$r_2=3^{-k_i}/2$, and call~$k_i$ the depth of the root. If~$r_2=R_*$, stop. Otherwise activate a nearest inactive face and take the least half-triadic radius at least~$L_*r_2+a_i$ about the newly projected center. The condition~$4L_*r_2\leq a_i$ for~$i\notin I$ keeps every other face outside the outer cube of this application.

\smallskip

Each transition activates a new coordinate. There are at most~$d$ transitions and~$d+1$ energy-decay applications, including at edges and corners. The intervening containments bound the enlargement factors by fixed constants. Telescoping the ratios in the energy-decay applications bounds their product by~$C(r/R_*)^{t_0}$. No assumption that the target meets only one face is made.

\smallskip

In each application let~$Y_i$ be the minimal scale of Proposition~\ref{mfd:prop-folded-iteration} at its reduced tolerance, for the original, translated stationary field, with its fixed outer-scale buffers. If the requested depth is less than~$Y_i$, use energy monotonicity instead: increasing a fixed~$c\geq t_0\log3$ makes~$e^{cY_i}\bigl(r_1/(3^{-k_i}/2)\bigr)^{t_0}\geq1$ in that case. Thus the same one-step inequality holds at every requested depth. The local Newtonian primitive of the evenly reflected bounded source has half-H\"older seminorm bounded by~$C(3^{-k_i}/2)^{\nf12}\|f\|_\infty$, by~\eqref{eq:mfd-primitive}. If instead a primitive supported on the fixed outer cube is used, its local bound has the extra factor~$1+\log(2\cdot3^{k_i})$. Keeping this larger bound, composition of the chain gives
\begin{multline*}
	\Gamma_N(u_N)(B_r(x)\cap((\tfrac12,\dots,\tfrac12)+\cu_0))
	\\
	\leq C e^{c\sum_iY_i}r^{t_0}
	\Bigl[R_*^{-t_0}\mathscr E_N^{(\tfrac12,\dots,\tfrac12)+\cu_0}(u_N)
	+\|f\|_\infty^2\sum_i
	\bar s_i^{-1}(3^{-k_i}/2)^{d+2-t_0}\bigl(1+\log(2\cdot3^{k_i})\bigr)^2\Bigr]\,.
\end{multline*}
Here and below constants absorb the at most~$2^d$ reflection multiplicity and the bounded number of enlargements. It remains to bound the random factors~$e^{c\sum_iY_i}$ and~$\bar s_i^{-1}$, up to small negative powers of the radii, by random variables that do not depend on the chain and whose moments are bounded uniformly in~$N$; Steps~3 and~4 do this.

\smallskip

\emph{Step 3.} We bound the factors~$e^{c\sum_iY_i}$ uniformly over a mesh of chains. At target depth~$n$, use spacing~$3^{-(n+\ell)}$. Enlarge the target by this spacing to contain every target centered in the corresponding mesh cell, and choose an interior mesh point even for targets touching the boundary. Catalogue all permutations of the coordinates and all choices of the depths of the~$d+1$ roots in~$\{0,\dots,n+\ell\}$. Activated sets are initial segments of the permutation, so every chain in Step~2 is represented. This catalogue has cardinality at most~$C_\ell(n+\ell+1)^{d+1}3^{dn}$. For a root of depth~$k\leq N$, the allowance is the minimal scale on the outer stationary cube at scale~$N-k+\ell$, plus~$\ell$ and the fixed iteration buffer; let it be zero when~$k>N$.

\smallskip

By the tail bound of Proposition~\ref{mfd:prop-folded-iteration}, there are fixed~$C_1$ and~$C$ such that each allowance~$Y$ satisfies
\begin{equation*}
	\P[Y>m]\leq C_1 e^{-\varkappa(m-C_1)_+},\qquad
	\varkappa\coloneqq\frac{\bigl(1-\frac{t_0-d+2}{2}\bigr)^2}{CM^2\gmcdelta^2|\log\gmcdelta|}\,.
\end{equation*}
This is uniform in every deterministic translated center, cutoff and root depth. For~$\zeta>0$ and~$\varkappa\geq2\zeta$, summing the tail gives
\begin{equation*}
	\E e^{\zeta Y}\leq e^{\zeta(C_1+1)}
	\left(1+C_1\sum_{m\geq0}(e^\zeta-1)e^{-(\varkappa-\zeta)m}\right)
	\leq C(C_1,\zeta)\,.
\end{equation*}
Choose~$\xi$ larger than all the moment orders needed for the products below and such that~$\xi(t_0-t_1)>d$, and then decrease~$\gmcdelta$ so that~$\varkappa\geq2(d+1)\xi c$. H\"older's inequality over the~$d+1$ applications gives a~$\xi$th moment of~$Ce^{c\sum_iY_i}$ that is uniform in~$N$ and in the chain. It requires no independence between the applications. Let
\begin{equation*}
	\mathsf m_N\coloneqq\sup_{n\geq0}\,\max_{\text{catalogue at depth }n}\,3^{-(t_0-t_1)n}Ce^{c\sum_iY_i}\,.
\end{equation*}
Then
\begin{equation*}
	\E\mathsf m_N^\xi\leq C\sum_{n\geq0}
	(n+\ell+1)^{d+1}3^{(d-\xi(t_0-t_1))n}<\infty\,.
\end{equation*}
The bound is uniform in~$N$. In particular, the factor~$Ce^{c\sum_iY_i}$ for a target ball of radius~$r$ is at most~$C\mathsf m_Nr^{-(t_0-t_1)}$.

\smallskip

\emph{Step 4.} We bound the averages~$\bar s_i^{-1}$ uniformly over a mesh of centers. Write~$s_{N,k}(x)$ for the integrand defining~$\bar s_{N,k}$. At spacing comparable to~$3^{-k}/2$, use the statistic
\begin{equation*}
	V_{N,k}(y)\coloneqq e^{\operatorname{osc}_{B_{3\cdot3^{-k}/2}(y)}
	(H_{L'}+\sum_{j<k}\g_{-j})}
	\bigl(s_{N,k}(y)+s_{N,k}(y)^{-1}\bigr)\,.
\end{equation*}
Every average~$\bar s_{N,k}(z)$ with~$|z-y|\leq3^{-k}/2$ and its inverse are bounded by~$V_{N,k}(y)$. The pointwise positive and negative moment bounds and the bounded-scale oscillation estimate used in Lemma~\ref{mfd:lem-cell-ellipticity} give, for every fixed auxiliary order~$q$ and~$k\leq N$, the bound~$\|V_{N,k}(y)\|_{L^q}\leq C_q e^{C(q+q^2)\gmcdelta^2 k}$. These bounds apply to the averages~$\bar s_{N,k}$; in particular, they include the ratio~$\ahom_{N-k}/\ahom_N$. Choose~$q$ at least as large as~$\xi$ and with~$q\theta>d$, and then decrease~$\gmcdelta$ until~$C(q+q^2)\gmcdelta^2<\tfrac12(\theta-d/q)\log3$. There are at most~$C3^{dk}$ points at depth~$k$. Consequently, for
\begin{equation*}
	V_N\coloneqq\sup_{k\leq N,\,y\text{ in the depth-}k\text{ mesh}}
	(3^{-k}/2)^\theta V_{N,k}(y)\,,
\end{equation*}
Minkowski's inequality bounds~$\|V_N\|_{L^q}$ by
\begin{equation*}
	C_q\sum_{k\geq0}3^{-(\theta-d/q)k}
	e^{C(q+q^2)\gmcdelta^2 k}<\infty\,,
\end{equation*}
uniformly in~$N$ and~$L'$. Thus~$\bar s_i^{-1}\leq CV_N(3^{-k_i}/2)^{-\theta}$, including the fixed depth shift in Step~1. The choice~$\theta<d+2-t_0$ gives~$\sup_{0<r_2\leq R_*}r_2^{d+2-t_0-\theta}(1+\log(1/r_2))^2<\infty$. With at most~$d+1$ source terms, Steps~2--4 therefore imply, for~$r\geq3^{-N}$,
\begin{equation*}
	\Gamma_N(u_N)(B_r(x)\cap((\tfrac12,\dots,\tfrac12)+\cu_0))
	\leq C\mathsf m_N r^{t_1}
	\bigl(\mathscr E_N^{(\tfrac12,\dots,\tfrac12)+\cu_0}(u_N)+V_N\|f\|_\infty^2\bigr)\,.
\end{equation*}
The mean-zero coercivity estimate~\eqref{eq:mfd-1} bounds~$\mathscr E_N^{(\tfrac12,\dots,\tfrac12)+\cu_0}(u_N)$ by~$CQ_N\|f\|_\infty^2$. Taking~$\mathsf m_N,V_N$ and~$Q_N$ at twice the requested orders bounds the resulting product in every requested~$L^p$. Here~$Q_N$ denotes the coercivity constant; the final constant in the statement is enlarged to include the resulting products.

\smallskip

The constants can be chosen uniformly in~$L'$. On a fixed compact set containing all cubes used above, the proof of Lemma~\ref{mfd:lem-infrared} gives~$\bigl\|\sum_{i\geq1}\|\g_i-\g_i(0)\|_{C^{1,1}}\bigr\|_{L^p}\leq C\gmcdelta\sqrt p$ and finite exponential moments at every fixed order. This sum bounds~$H$,~$0$ and all~$H_{L'}$ simultaneously. In particular, the averages~$\bar s_{N,k}$ and the inverse ellipticity constants for~$H$ and for every~$H_{L'}$ compare with those for~$L'=0$ by the exponential of this sum. Use the doubled orders in these products; the contributions of~$\nabla H$ and~$\nabla H_{L'}$ to~$D_N$ are bounded by the same sum as well. The minimal scales~$Y_i$ are those of the reindexed stationary coefficient. Every root actually used in a chain has~$k\geq \ell$, hence its outer scale~$N-k+\ell$ is at most~$N\leq N+L'$. Equation~\eqref{e.cutoff.Holder.replacements} therefore uses the same stopping variable for all~$L'$. Catalogue entries that occur in no chain use the same minimal scale. Thus no union over infinitely many truncations~$H_{L'}$ is needed. The radius~$R_*$ and all geometric constants were fixed before~$L'$. For the stationary coefficient~$A_N^0$, the translates~$w'+\cu_0$ use the law-preserving translation of the fields, and their constants~$\mathsf m_N$,~$V_N$ and~$Q_N$ are computed in the translated fields. Every spatial maximum needed inside a fixed cube is included in the displayed mesh sums. A fixed finite collection of such translates can be included in their catalogues without changing the exponent conditions.

\smallskip

\emph{Step 5.} We treat the radii below~$3^{-N}$ and pass to~$H$. Reflect evenly across every face meeting~$B_{2r_{N,*}}(x)$, with~$r_{N,*}$ as below. The coordinate fold is~$1$-Lipschitz, so it preserves the required small oscillation of~$\log A_N$; the reflected weak equation has the evenly reflected bounded source. Take~$p_1$ so large that~$d-2d/p_1>t_1$ and~$1-d/p_1>\alpha$, and apply~\eqref{eq:mfd-micro} at the radius~$r_{N,*}$ of Step~4 of the proof of Proposition~\ref{mfd:prop-growth}. The bound at the end of Step~4 of this proof, at~$r=C3^{-N}$, bounds the energy on~$B_{2r_{N,*}}$ in~\eqref{eq:mfd-micro}. Dividing by~$r^t$ gives the two factors~$3^{-(t_1-t)N}$ and~$3^{-(d+2-t)N}$ from Step~4 of Proposition~\ref{mfd:prop-growth}; their strictly positive exponents absorb the moments of Lemma~\ref{mfd:lem-extremes} at the finitely many larger orders used in the products. Divide by the reflection multiplicity to return to~$(\tfrac12,\dots,\tfrac12)+\cu_0$.

\smallskip

For the H\"older estimate, the local coarse-grained Poincar\'e inequality~\eqref{e.CG.Poincare.mean.zero} applied to the same energy bound gives the Campanato exponent~$t_1+2-d$ minus an arbitrarily small loss. Choose that loss so that the exponent is larger than~$2\alpha$. The maximum over a mesh of the inverse coarse-grained ellipticity constants in this inequality has the same convergent moment summation as~\eqref{eq:mfd-3}. Below~$3^{-N}$ use the preceding~$W^{1,p_1}$ estimate and Morrey's inequality. The first two factors in Step~7 of Proposition~\ref{mfd:prop-growth} have positive exponents because~$t_1>d-2+2\alpha$; there is no boundary-datum term. Telescope the cube averages and use the mean-zero condition to control the supremum norm. Enlarging the random constant gives both asserted estimates. Every constant constructed above is a statistic of the coefficient, independent of~$f$ and of the chosen solution.

\smallskip

On the event of Lemma~\ref{mfd:lem-infrared} on which~$H_{L'}\to H$, the truncated coefficients converge uniformly on~$\overline{(\tfrac12,\dots,\tfrac12)+\cu_0}$. Compare the weak equations on the mean-zero space, testing with the difference of the solutions. Coercivity yields gradient convergence; the H\"older bound yields uniform convergence. Fatou's lemma preserves the moment bounds. One can use the same moment-bounded coefficient statistics for every source before the passage~$L'\to\infty$, so the final constant is source independent. Countable intersections make the event common to all cutoffs and to all cubes of the meshes of Steps~3 and~4; the deterministic weak-equation estimates then apply to every bounded mean-zero source.

\smallskip

Finally consider a fixed cube~$w+R\cu_0$. For~$A_N$, choose the integer~$j$ with~$3^{j-1}<R\leq3^j$, transport the stationary environment by translation and that scale shift, and run the preceding argument on the remaining cube of side in a fixed bounded interval. For~$A_N^0$ and the truncations~$H_{L'}$, a shift of a cube with~$R>1$ to bounded size deletes the fields of the largest scales, which only~$A_N$ contains. We therefore run the preceding argument directly on~$w+R\cu_0$ in these cases. Its chains have roots of radius at most~$R_*$, so their outer scales do not exceed the cutoff; the catalogue of Step~3 is larger by a factor depending only on~$R$; and the total energy is bounded by the coercivity bound for larger cubes in the proof of Lemma~\ref{mfd:lem-neumann-error}. For~$A_N$, we have~$A_N(w+Ry)=s_N(w+\cu_j)\widetilde A(y)$, where~$\widetilde A(y)\coloneqq A_{N+j}^{\Theta_{-j,w}\omega}(3^{-j}Ry)$, and~$\widetilde u(y)\coloneqq u_N(w+Ry)$. The source of~$\widetilde u$ is~$R^2s_N(w+\cu_j)^{-1}f(w+Ry)$, and the energy of~$u_N$ is~$R^{d-2}s_N(w+\cu_j)$ times that of~$\widetilde u$. Thus the transported Morrey constant is multiplied by~$R^{d+2-t}s_N(w+\cu_j)^{-1}$, and the H\"older constant by~$R^2\max\{1,R^{-\alpha}\}s_N(w+\cu_j)^{-1}$. The factor~$s_N(w+\cu_j)^{-1}$ is bounded by an exponential of~$|H(w)|$ and of the finitely many fields~$\g_i(w)$ entering~$s_N(w+\cu_j)$; it has every required moment uniformly in~$N$ and~$L'$ by Lemma~\ref{mfd:lem-infrared} and Assumption~\ref{a.g2}. Use twice the requested orders in this final product. If~$N+j<0$, which holds only for the finitely many~$N<-j$, the microscopic argument of Step~5 covers these cutoffs. The exponent conditions and the disorder threshold are the same; the final moment bounds may depend on the fixed cube.

\smallskip

The proof also gives both estimates for~$A_N^0$ and for the coefficients obtained with~$H_{L'}$, with the same threshold and moment bounds, uniformly in~$L'$. The solution and energy are computed for each coefficient. On each fixed cube one~$K_N$ suffices for all these coefficients and all bounded mean-zero sources. For~$A_N^0$, stationarity gives the bounds on each fixed translate of the cube, and taking the maximum of the constants handles every fixed finite collection of translates.
\end{proof}

Applying Proposition~\ref{mfd:cor-neumann-source} to the mean-zero source~$f_\ep$ of~\eqref{eq:mfd-bump}, with~$\|f_\ep\|_\infty\leq C\ep^{-1}$, gives, with one random constant for all~$\ep$,
\begin{equation}
	\Gamma_N(v_{N,\ep})(B_r(x)\cap((\tfrac12,\dots,\tfrac12)+\cu_0))
	\leq K_N\ep^{-2}r^t\,.
	\label{eq:mfd-13}
\end{equation}

\subsection{Compactness of the finite-volume energies}
\label{mfd:sec-resampling}

In this subsection we prove that the finite-volume energies form a compact family in~$L^p(\P)$ (Proposition~\ref{mfd:prop-response-compact}). The effect of a field on a finite-volume energy decays as a power of its scale (Lemma~\ref{mfd:lem-15}). We delete the field on thin strips along the faces of a grid of boxes. By the Morrey bound~\eqref{eq:mfd-4}, the strips and each box carry energy at most a positive power of the scale of the field (Lemma~\ref{mfd:lem-strips}); the restrictions of the field to the parts of different boxes outside the strips are independent, and the Efron--Stein inequality bounds the fluctuations. Summing over the fields that we do not retain, we approximate each finite-volume energy, uniformly in the cutoff, by its conditional expectation given finitely many fields (Proposition~\ref{mfd:prop-16}).

\smallskip

On a fixed cube~$w+R\cu_0$, by a \textbf{finite-volume energy}~$\mathcal R_N$ we mean either the minimal energy~$\Lambda_{N,w+R\cu_0}(\varphi)$ with a nonconstant smooth boundary datum~$\varphi$, or the quadratic form~$\langle f,G_N^{w+R\cu_0}f\rangle$ of the Dirichlet Green operator for a nonzero smooth source~$f$. Write~$\mathcal R(g)$ for the same energy computed with the coefficient whose logarithm is~$g$, so that~$\mathcal R_N=\mathcal R(\log A_N)$.

\subsubsection{Local perturbations of the coefficient}

Let~$A$ be bounded above and below by positive constants, and let~$g$ be bounded and supported in~$U\subset w+R\cu_0$. Write~$\mathscr E^g(v)\coloneqq\int_{w+R\cu_0}e^gA|\nabla v|^2$, with energy measure~$\Gamma_{\mathscr E^g}(v)\coloneqq e^gA|\nabla v|^2\,dx$, and~$\mathscr E\coloneqq \mathscr E^0$. Let~$u$ and~$u^g$ be the minimizers for the coefficients~$A$ and~$e^gA$ with the same Dirichlet boundary values, or the solutions for these coefficients with the same linear right-hand side and zero Dirichlet or zero conormal boundary condition; in the conormal case, the right-hand side vanishes on constants and the solutions have mean zero. Denote their energies by~$\mathcal R$ and~$\mathcal R^g$.

\begin{lemma}[Multiplicative perturbations]
\label{mfd:cor-14}
For the coefficients and solutions just defined,
\begin{equation}
	\begin{aligned}
		&\mathscr E^g(u^g-u)\leq\|g\|_\infty^2e^{3\|g\|_\infty}\Gamma_{\mathscr E}(u)(U),\qquad
		\Gamma_{\mathscr E^g}(u^g)(U)\leq4e^{\|g\|_\infty}\Gamma_{\mathscr E}(u)(U)\,,\\
		&|\mathcal R^g-\mathcal R|\leq2\|g\|_\infty e^{4\|g\|_\infty}\Gamma_{\mathscr E}(u)(U)\,.
	\end{aligned}
	\label{eq:mfd-14}
\end{equation}
For two log coefficients~$g_1$ and~$g_2$ and each of these energies,
\begin{equation}
	e^{-\|g_1-g_2\|_\infty}\mathcal R(g_2)\leq\mathcal R(g_1)
	\leq e^{\|g_1-g_2\|_\infty}\mathcal R(g_2)\,.
	\label{eq:mfd-mult}
\end{equation}
\end{lemma}
\begin{proof}
Subtracting the Euler equations gives~$\mathscr E^g(u^g-u,\phi)=-\int_U(e^g-1)\,d\Gamma_{\mathscr E}(u,\phi)$. Test with~$u^g-u$ and apply Cauchy--Schwarz to obtain
\begin{equation}
	\mathscr E^g(u^g-u)^{\nf12}\leq
	(e^{\|g\|_\infty/2}-e^{-\|g\|_\infty/2})\Gamma_{\mathscr E}(u)(U)^{\nf12}\,.
	\label{eq:mfd-perturbation-sharp}
\end{equation}
This proves the first bound in~\eqref{eq:mfd-14}; the triangle inequality in~$L^2(U,e^gA\,dx)$ gives the second estimate. The Euler equations also give
\begin{equation*}
	\mathcal R^g-\mathcal R
	=\pm\left[\int_U(e^g-1)\,d\Gamma_{\mathscr E}(u)-\mathscr E^g(u^g-u)\right]\,,
\end{equation*}
with plus sign when~$\mathcal R$ is a minimal energy with given boundary values and minus sign when it is the energy of a solution with given right-hand side. Using~$|e^g-1|\leq\|g\|_\infty e^{\|g\|_\infty}$ proves the third bound. Monotonicity in the coefficient and scaling by constants give~\eqref{eq:mfd-mult}. These arguments use only the Euler equations and the Cauchy--Schwarz inequality for energy measures, so they also apply to the regular strongly local forms constructed in Section~\ref{mfd:sec-local-form}.
\end{proof}

\subsubsection{Resampling one random field}

We bound the change in a finite-volume energy when a field~$\g_{-j}$ with~$j\geq0$ is replaced by an independent copy. The field has scale~$r\coloneqq3^{-j}$ and range of dependence~$\sqrt d\,r$, so its restrictions to boxes at mutual distance larger than~$\sqrt d\,r$ are independent. We partition~$w+R\cu_0$ into boxes of side~$r_1$ with~$r\ll r_1\ll1$. Delete the field on strips of width of order~$r$ around the faces of the boxes, and apply the Efron--Stein inequality to the independent restrictions of the field to the remaining interiors of the boxes. The energy in the strips controls the error caused by deleting the field there. The largest energy in one box controls the differences~$X-X^{(i)}$ in the Efron--Stein inequality~\eqref{eq:mfd-ES} below. By the Morrey bound~\eqref{eq:mfd-4}, both are bounded by a positive power of the scale~$r$. Exponential moments of the sup norm of the field enter through Lemma~\ref{mfd:cor-14}; they grow more slowly than any negative power of the scale~$r$, so this decay absorbs them.

For~$r>0$ and~$4\sqrt d\,r<r_1\leq1$, partition a fixed cube~$w+R\cu_0$ into boxes~$z+r_1\cu_0$, and call \textbf{strips} the points within distance~$4\sqrt d\,r$ of a box face.

\begin{lemma}[Mass near the box faces and in a single box]
\label{mfd:lem-strips}
Let~$t>d-1$ and~$K<\infty$, and let~$\nu$ be a finite positive measure on~$w+R\cu_0$ such that~$\nu(B_s(x)\cap(w+R\cu_0))\leq Ks^t$ for every~$x\in w+R\cu_0$ and~$0<s\leq1$. There is~$C(d,t,w,R)<\infty$ such that
\begin{equation*}
	\nu(\text{strips})\leq CKr_1^{-1}r^{t-d+1}\,,
\end{equation*}
and
\begin{equation*}
	\max_{z+r_1\cu_0\ \mathrm{box}}\nu(z+r_1\cu_0)\leq CKr_1^t\,.
\end{equation*}
With~$r_1=r^{(t-d+1)/(t+1)}$ and sufficiently small~$r$, both bounds are~$CKr^b$, where~$b=t(t-d+1)/(t+1)>0$.
\end{lemma}

\begin{proof}
There are at most~$C_{w,R}r_1^{-1}$ face hyperplanes. Each strip is covered by~$C_{w,R}r^{1-d}$ balls of radius~$4\sqrt d\,r$, and each box by a dimensional number of balls of radius~$r_1$. Apply the bound~$\nu(B_s(x)\cap(w+R\cu_0))\leq Ks^t$ and use~$r_1^{-1}r^{t-d+1}=r_1^t=r^b$ for the stated choice of~$r_1$.
\end{proof}

The sup norm on~$w+R\cu_0$ of a fine-scale field grows only logarithmically in the inverse of its scale, and its exponential moments grow more slowly than any power of the inverse scale. By~\eqref{e.coefficient.field.regularity}, the variable~$\|\g_0\|_{L^\infty(\cu_0)}$ is sub-Gaussian with parameter~$C\gmcdelta$. Let~$j\geq0$ and~$r\coloneqq3^{-j}$, and let~$S\coloneqq\|\g_{-j}\|_{L^\infty(w+R\cu_0)}$. The field~$\g_{-j}$ has the law of~$\g_0(3^j\,\cdot)$. Cover~$3^j(w+R\cu_0)$ by at most~$C_{w,R}r^{-d}$ unit cubes to obtain~$\P[S>s]\leq2C_{w,R}r^{-d}\exp(-s^2/C\gmcdelta^2)$. The union bound and~$(s+c)^2\geq s^2+c^2$ show that the excess of~$S$ over~$\gmcdelta\sqrt{C_{w,R}(1+|\log r|)}$ has the tail~$\exp(-s^2/C\gmcdelta^2)$. Thus the excess is sub-Gaussian with parameter~$C\gmcdelta$, uniformly in~$r$. Consequently, with constants also allowed to depend on the fixed cube~$w+R\cu_0$, for all~$k\geq0$,~$\lambda\geq0$,~$p<\infty$ and~$\gamma>0$, and every~$0<\gmcdelta\leq1$,
\begin{equation}
	\|S^ke^{\lambda S}\|_{L^p}
	\leq C(p,k,\lambda)\,\gmcdelta^{k}(1+|\log r|)^{\nf{k}{2}}
	\exp\bigl(C\lambda\gmcdelta\sqrt{1+|\log r|}+C_p\lambda^2\gmcdelta^2\bigr)
	\leq C(p,k,\lambda,\gamma)\,\gmcdelta^{k}r^{-\gamma}\,.
	\label{eq:mfd-layernorm}
\end{equation}
The splitting of~$S$ into a deterministic logarithmic part and a sub-Gaussian excess is needed: a union bound applied directly to the exponential moment would produce the factor~$C_{w,R}r^{-d}$.

\smallskip

The fluctuations are estimated with the following moment form of the Efron--Stein inequality. Let~$\xi_1,\dots,\xi_m$ be independent random elements, let~$p\geq2$, and let~$X=X(\xi_1,\dots,\xi_m)\in L^p$ be measurable. Write~$X^{(i)}$ for the variable obtained by replacing~$\xi_i$ with an independent copy, and~$\E_i$ for the integral over this copy alone. Then
\begin{equation}
	\|X-\E X\|_{L^p}\leq C_p\Bigl\|\Bigl(\sum_i\E_i\bigl[(X-X^{(i)})^2\bigr]\Bigr)^{\nf12}\Bigr\|_{L^p}\,.
	\label{eq:mfd-ES}
\end{equation}
This is a consequence of~\citet*[Theorem~2]{BBLM2005}, who proved two one-sided bounds in terms of variance proxies; both proxies are dominated by the sum on the right side of~\eqref{eq:mfd-ES}. The inequality also holds conditionally on every~$\sigma$-field independent of the random elements and their copies.

\begin{lemma}[Resampling one random field]
\label{mfd:lem-15}
Let~$0<\gmcdelta\leq1$,~$t\in(d-1,d)$ and~$p\geq2$, and let~$\mathcal R_N$ be a finite-volume energy. Assume that, for every~$N$, the function~$u_N$ defining~$\mathcal R_N$ (the minimizer with the given boundary values, or the solution with zero Dirichlet boundary values) satisfies~\eqref{eq:mfd-4} on~$w+R\cu_0$ with the constant~$K_N$, and that
\begin{equation*}
	\sup_N\|K_N\|_{L^{3p}}+\sup_N\|\mathcal R_N\|_{L^{3p}}<\infty\,.
\end{equation*}
For~$0\leq j\leq N$, let~$\mathcal R_N^{(j)}$ be the finite-volume energy obtained after replacing the field~$\g_{-j}$ by an independent copy. Let~$b$ be the exponent of Lemma~\ref{mfd:lem-strips}. Then there exists~$C_p<\infty$, depending on~$p$,~$d$,~$t$,~$w+R\cu_0$, on the boundary datum or source defining~$\mathcal R_N$ and on the two suprema in the previous display, such that, for all integers~$0\leq j\leq N$,
\begin{equation}
	\|\mathcal R_N-\mathcal R_N^{(j)}\|_{L^p}\leq C_p\,\gmcdelta\,3^{-bj/8}\,.
	\label{eq:mfd-15}
\end{equation}
\end{lemma}

\begin{proof}
Let~$r\coloneqq3^{-j}$ and let~$S\coloneqq\|\g_{-j}\|_\infty$ on~$w+R\cu_0$, as in~\eqref{eq:mfd-layernorm}. A superscript~$(j)$ denotes the corresponding quantity after this field is resampled. Throughout, expectations over this field are conditional on the other fields. For the finitely many scales at which the strip construction of Lemma~\ref{mfd:lem-strips} does not apply, the bound~$|\mathcal R_N-\mathcal R_N^{(j)}|\leq(S+S^{(j)})e^{S+S^{(j)}}\mathcal R_N$ and~\eqref{eq:mfd-layernorm} suffice.

\smallskip

At the remaining scales, partition~$w+R\cu_0$ into boxes of side~$r_1=r^{(t-d+1)/(t+1)}$. Choose smooth cutoff functions~$\chi_i$, one for each box, supported on interior subboxes that are separated by more than~$\sqrt d\,r$, such that~$\chi=\sum_i\chi_i$ equals one outside the strips of width~$4\sqrt d\,r$. Let~$\widehat{\mathcal R}_N$ be the energy after deleting~$(1-\chi)\g_{-j}$. Let~$\widehat\nu$ be the energy measure of the corresponding solution. The restrictions of the field to these subboxes are independent. Lemmas~\ref{mfd:cor-14} and~\ref{mfd:lem-strips} give~$\widehat\nu(w+R\cu_0)\leq e^S\mathcal R_N$,
\begin{equation*}
	|\mathcal R_N-\widehat{\mathcal R}_N|\leq CSe^{CS}K_Nr^b\,,
\end{equation*}
and
\begin{equation*}
	\max_U\widehat\nu(U)\leq Ce^{CS}K_Nr^b\,,
\end{equation*}
where the maximum is over the interior subboxes~$U$. For the last bound, apply~\eqref{eq:mfd-perturbation-sharp} on the strips and the inequality that~$\widehat\nu(U)$ is at most~$2e^S\Gamma_N(u_N)(U)$ plus twice the energy, for the form with the deleted field, of the difference between the minimizers before and after the deletion. The second term measures the change in the minimizer when replacing~$\g_{-j}$ by~$\chi\g_{-j}$.

\smallskip

Resampling the field on an interior subbox~$U$ changes~$\widehat{\mathcal R}_N$ by at most~$C(S+S^{(j)})e^{C(S+S^{(j)})}\widehat\nu(U)$. Since~$\sum_U\widehat\nu(U)^2\leq\max_U\widehat\nu(U)\widehat\nu(w+R\cu_0)$, the conditional Efron--Stein inequality~\eqref{eq:mfd-ES} and conditional Jensen give
\begin{equation*}
	\|\widehat{\mathcal R}_N-\E_j\widehat{\mathcal R}_N\|_p
	\leq C_pr^{\nf{b}{2}}
	\bigl\|(S+S^{(j)})e^{C(S+S^{(j)})}K_N^{\nf{1}{2}}\mathcal R_N^{\nf{1}{2}}\bigr\|_p
	\leq C_p\gmcdelta r^{\nf{3b}{8}}\,.
\end{equation*}
The last step uses H\"older's inequality with three factors and~\eqref{eq:mfd-layernorm}, at a total loss of~$r^{-\nf{b}{8}}$. The errors~$\mathcal R_N-\widehat{\mathcal R}_N$ and~$\widehat{\mathcal R}_N^{(j)}-\mathcal R_N^{(j)}$ from deleting the field on the strips are at most~$C_p\gmcdelta r^{\nf{7b}{8}}$ in~$L^p$. The energy~$\widehat{\mathcal R}_N^{(j)}$, computed with the cutoff function applied to the resampled field, has the same conditional mean as~$\widehat{\mathcal R}_N$, so
\begin{equation}
	\mathcal R_N-\mathcal R_N^{(j)}=(\mathcal R_N-\widehat{\mathcal R}_N)+(\widehat{\mathcal R}_N-\E_j\widehat{\mathcal R}_N)
	-(\widehat{\mathcal R}_N^{(j)}-\E_j\widehat{\mathcal R}_N^{(j)})+(\widehat{\mathcal R}_N^{(j)}-\mathcal R_N^{(j)})\,.
	\label{eq:mfd-fourterm}
\end{equation}
Since~$r^{\nf{3b}{8}}\leq3^{-bj/8}$, this proves~\eqref{eq:mfd-15}.
\end{proof}

The following martingale bound adds up the effects of the fields on which we do not condition.

\begin{lemma}[Summing resampling bounds]
\label{mfd:lem-resampling-sum}
Let~$(\xi_j)_{j\in I}$ be a countable independent family and let~$X\in L^p$ with~$1\leq p<\infty$ be measurable with respect to this family. Write~$X^{(j)}$ for the result of replacing~$\xi_j$ by an independent copy. For every~$I_0\subseteq I$,
\begin{equation*}
	\bigl\|X-\E[X\mid\sigma(\xi_j:j\in I_0)]\bigr\|_p
	\leq\sum_{j\notin I_0}\|X-X^{(j)}\|_p\,.
\end{equation*}
\end{lemma}

\begin{proof}
Condition first on the retained variables, and then reveal the omitted ones one at a time, in some enumeration. At the step that reveals~$\xi_j$, the martingale difference is the conditional expectation of~$X-X^{(j)}$ given the variables revealed so far. Conditional Jensen bounds its~$L^p$ norm by~$\|X-X^{(j)}\|_p$. Sum and use~$L^p$ martingale convergence.
\end{proof}

We apply this bound to approximate an energy by conditioning on finitely many fields.

\begin{proposition}[Approximation by conditioning on finitely many fields]
\label{mfd:prop-16}
Let~$\ell\geq0$ be an integer, and let~$b$ be the exponent of Lemma~\ref{mfd:lem-strips}.
\begin{itemize}
\item Let~$\mathcal R_N$ be a finite-volume energy satisfying the hypotheses of Lemma~\ref{mfd:lem-15} for every~$j\geq0$. Then there exists~$C_p<\infty$ such that, for every~$N\geq \ell$,
\begin{equation}
	\bigl\|\mathcal R_N-\E[\mathcal R_N\mid\g_{-j},\ |j|\leq\ell]\bigr\|_{L^p}
	\leq C_p\,\gmcdelta\,3^{-b\ell/8}\,.
	\label{eq:mfd-16}
\end{equation}
\item On~$(\tfrac12,\dots,\tfrac12)+\cu_0$, suppose instead that~$\mathcal R_N=|(\tfrac12,\dots,\tfrac12)+\cu_0|\,e\cdot \a_*^{-1}((\tfrac12,\dots,\tfrac12)+\cu_0;A_N)e$ is the Neumann energy for a Euclidean unit vector~$e$, let~$p\geq2$ and~$\gmcdelta\leq\gmcdelta_0(d,t,p)$, and assume
\begin{itemize}
\item the bound~\eqref{eq:mfd-13}, with~$K_N$ bounded in~$L^{3p}$,
\item the bound on~$\sup_{N,\ep}\|\mathscr E_N^{(\tfrac12,\dots,\tfrac12)+\cu_0}(v_{N,\ep})\|$ in~\eqref{eq:mfd-9} at order~$3p$, and
\item the approximation-error bound in~\eqref{eq:mfd-9} at order~$p$.
\end{itemize}
The threshold is chosen so that the moment range of Proposition~\ref{p.coarse.grained.bound} contains~$4p$ and so that Lemma~\ref{mfd:lem-cell-ellipticity} applies at order~$4p$; that lemma supplies the uniform~$L^{4p}$ bound for~$\mathcal R_N$. Under these conditions,~\eqref{eq:mfd-16} holds with right side~$C_p\gmcdelta^{\nf12}3^{-b\ell/256}$.
\end{itemize}
\end{proposition}

\begin{proof}
For~$j\geq0$,~\eqref{eq:mfd-15} bounds the resampling error by~$C_p\gmcdelta3^{-bj/8}$, and the error is zero for~$j>N$. For~$j<0$, the field~$H$ changes by the difference of the expressions~$\g_{-j}-\g_{-j}(0)$ for the field and for its independent copy. This difference has norm of size~$C\gmcdelta3^j$ and sub-Gaussian tails by Lemma~\ref{mfd:lem-infrared}. Equation~\eqref{eq:mfd-mult} bounds the resampling error by~$C_p\gmcdelta3^j$. Sum both geometric series using Lemma~\ref{mfd:lem-resampling-sum}.

\smallskip

For the Neumann energy on~$(\tfrac12,\dots,\tfrac12)+\cu_0$, first use the volume source~$f_\ep$. The proof of Lemma~\ref{mfd:lem-15} applies because Lemma~\ref{mfd:cor-14} also covers zero-conormal problems with a given right-hand side. Here~$C_p$ depends on the uniform~$L^{3p}$ bound for the source-independent~$K_N$ of~\eqref{eq:mfd-13}, the uniform~$L^{3p}$ bounds for the smoothed energies, and the uniform~$L^p$ bound for the approximation error in~\eqref{eq:mfd-9} after division by~$\ep^{\nf14}$. Proposition~\ref{p.coarse.grained.bound} and Lemma~\ref{mfd:lem-cell-ellipticity} are used at the order~$4p$ as well. All these orders are fixed before decreasing~$\gmcdelta$; the moment orders are those specified in the statement.

\smallskip

Choose a numerical bound~$M$ for the~$L^{3p}$ norms of~$K_N$ and the smoothed energies, the~$L^{4p}$ norms of~$\mathcal R_N$, and the~$L^p$ approximation norm after division by~$\ep^{\nf14}$. The last two inputs are supplied by Lemma~\ref{mfd:lem-cell-ellipticity} and~\eqref{eq:mfd-9}. Let~$C_{\rm resp}\coloneqq C(d)$ be the constant in Proposition~\ref{p.coarse.grained.bound}, and decrease~$\gmcdelta$ so that the admissible range of moment orders in that proposition contains~$4p$, that is,~$4p\leq C_{\rm resp}^{-1}\gmcdelta^{-2}|\log\gmcdelta|^{-1}$. The constants and disorder threshold are chosen after these bounds and the fixed bump profile, and before~$N$ and~$\ell$. Here~$p\geq2$ as in Lemma~\ref{mfd:lem-15};~$C_p$ depends on~$d,t,p,M,C_{\rm resp}$ and the fixed data and bump profile.

\smallskip

By~\eqref{eq:mfd-13}, the Morrey constant for the solution with source~$f_\ep$ is~$K_N\ep^{-2}$. This factor enters deletion of the field on strips linearly and the Efron--Stein estimate through its square root. Together with~\eqref{eq:mfd-9} and H\"older's inequality at the stated orders, this gives
\begin{equation*}
	\|\mathcal R_N-\E[\mathcal R_N\mid\g_{-j},\ |j|\leq\ell]\|_p
	\leq C_p\bigl(\ep^{\nf14}
	+\gmcdelta\ep^{-2}3^{-b\ell/8}
	+\gmcdelta3^{-\ell}\bigr)\,.
\end{equation*}
Choose~$\ep=3^{-b\ell/32}/16$ to obtain~$C_p3^{-b\ell/128}$.

\smallskip

For the other bound compare with~$|(\tfrac12,\dots,\tfrac12)+\cu_0|\,e\cdot e=1$. Let~$\mathcal R_N^0$ be the Neumann energy for~$A_N^0$. The physical dilation~$X=3^Nx$ maps~$(\tfrac12,\dots,\tfrac12)+\cu_0$ to the translated cube of side~$3^N$ at its matching stationary cutoff~$N$. Reindex the stationary environment by this dilation and translation. Write~$J_N^{\rm ref}$ for the maximum over~$|e|=1$ on the right of~\eqref{e.coarse.grained.matrix.deviation.by.J} for this cube, the coefficient~$\ahom_N^{-1}\a_N$ and~$\b=\Id$. Proposition~\ref{p.coarse.grained.bound}, at the finite orders up to~$4p$, gives~$\|J_N^{\rm ref}\|_{L^{4p}}\leq C_p\gmcdelta^2$. The matrix-deviation estimate~\eqref{e.coarse.grained.matrix.deviation.by.J} gives
\begin{equation*}
	|\mathcal R_N^0-1|^2\leq C J_N^{\rm ref}(1+J_N^{\rm ref}),\qquad
	\|\mathcal R_N^0-1\|_p\leq C_p\gmcdelta\,.
\end{equation*}
Restoring~$H$ changes the energy by at most~$(e^{\|H\|_{L^\infty((\tfrac12,\dots,\tfrac12)+\cu_0)}}-1)\mathcal R_N^0$; H\"older's inequality and Lemma~\ref{mfd:lem-infrared} bound this change by~$C_p\gmcdelta$. Conditional expectation is an~$L^p$ contraction, so
\begin{equation*}
	\|\mathcal R_N-\E[\mathcal R_N\mid\g_{-j},\ |j|\leq\ell]\|_p
	\leq2\|\mathcal R_N-1\|_p\leq C_p\gmcdelta\,.
\end{equation*}
Taking the geometric mean with the preceding decay bound proves the asserted~$C_p\gmcdelta^{\nf12}3^{-b\ell/256}$ estimate.

\smallskip

For the cell terms in Lemma~\ref{mfd:lem-band}, apply this estimate to the normalized Neumann energy on the unit cube in the translated environment at the cell's own cutoff. Translation between~$\cu_0$ and~$(\tfrac12,\dots,\tfrac12)+\cu_0$ preserves the Neumann energy. More generally, if~$A_N(w+Ry)=s_N(w+R\cu_0)\widetilde A_m(y)$ for a coefficient~$\widetilde A_m$ at cutoff~$m$, then the Neumann energy on~$w+R\cu_0$ is~$R^d s_N(w+R\cu_0)^{-1}$ times the Neumann energy for~$\widetilde A_m$ on~$\cu_0$. The conditioning estimate applies after removing this scalar and at cutoff~$m$. Its finite moment inputs are chosen at~$3P$ and~$6P$ when the Neumann energy is used at~$P$ and~$2P$. For general fixed physical cubes we use the transported volume loads in the proof of Lemma~\ref{mfd:lem-neumann-error}, followed by Proposition~\ref{mfd:prop-as-response-bank}.
\end{proof}

\subsubsection{Compactness of conditional expectations}

On a cube~$z+r\cu_0$, let~$\mathcal R_N$ denote the minimal energy~$\Lambda_{N,z+r\cu_0}(\varphi)$ with boundary values~$\varphi$, the quadratic form~$\langle f,G_N^{z+r\cu_0}f\rangle$ of the Dirichlet Green operator with source~$f$, or the Neumann energy~$r^d\,e\cdot \a_*^{-1}(z+r\cu_0;A_N)e$ for a constant flux~$e\in\Rd$.

\begin{proposition}[Strong compactness of the finite-volume energies]
\label{mfd:prop-response-compact}
Assume
\begin{itemize}
\item the moment orders satisfy~$1\leq p<q<\infty$ and~$\sup_N\|\mathcal R_N\|_{L^q}<\infty$,
\item for some~$a>0$ and~$C_{z,r}<\infty$, every pair of integers~$0\leq\ell\leq N$ satisfies
\begin{equation*}
	\bigl\|\mathcal R_N-\E[\mathcal R_N\mid\g_{-j},\ |j|\leq\ell]\bigr\|_{L^p}
	\leq C_{z,r}\gmcdelta3^{-a\ell}\,.
\end{equation*}
\end{itemize}
Here~$C_{z,r}$ may depend on the fixed cube, data and disorder, but not on~$N$ or~$\ell$. For a countable collection we assume that the exponent~$a$ is common, while~$C_{z,r}$ and the moment bound may depend on its member. Then the following statements hold.
\begin{enumerate}[label=(\roman*)]
\item The sequence~$(\mathcal R_N)_N$ is relatively compact in~$L^p(\P)$.
\item The same holds for~$\langle f,G_N^{z+r\cu_0}f'\rangle$ with two sources~$f$ and~$f'$.
\item For every countable collection of cubes, sources, boundary data and fixed integer shifts~$N\mapsto N+j$ of the cutoff~$N$, a common subsequence can be chosen along which all the corresponding energies converge in~$L^p(\P)$.
\end{enumerate}
\end{proposition}

For the minimal energy and the Dirichlet Green quadratic form the approximation is supplied by Proposition~\ref{mfd:prop-16}; for the Neumann energy on a fixed cube it follows from the approximation of~$L_{w,R,e}$ in the proof of Lemma~\ref{mfd:lem-neumann-error}, as shown below.

\begin{proof}
Fix~$\ell$ and~$N\geq\ell$. The log coefficient on~$\overline{z+r\cu_0}$ is the sum of the retained potential~$\sum_{-\ell\leq j<0}(\g_{-j}-\g_{-j}(0))+\sum_{0\leq j\leq\ell}\g_{-j}$ and a remaining potential independent of the retained fields. As a function of~$v\in C(\overline{z+r\cu_0})$, take the expectation of the energy for the log coefficient~$v$ plus that remaining potential. Independence identifies its value at the retained potential with the conditional expectation of~$\mathcal R_N$. By~\eqref{eq:mfd-mult}, values at~$v,\psi$ compare by~$e^{\|v-\psi\|_\infty}$. Choose a fixed radius which bounds the retained potential with positive probability. Integration on this event bounds the value at~$0$ by~$\E\mathcal R_N$ times the exponential of the radius divided by that probability. The functions are therefore uniformly bounded and equicontinuous on every bounded set. A diagonal extraction on a countable dense subset of~$C(\overline{z+r\cu_0})$ gives pointwise convergence everywhere. Conditional Jensen gives a uniform~$L^q$ bound, so convergence at the retained random potential is in~$L^p$.

\smallskip

For every~$\ep>0$, choose~$\ell$ so that~$C_{z,r}\gmcdelta3^{-a\ell}<\ep/2$. A finite~$\ep/2$-net of the conditional expectations is an~$\ep$-net of the energies with~$N\geq\ell$. The finitely many energies with~$N<\ell$ form a compact set. This proves~(i).

\smallskip

For the Neumann energy on a general fixed cube~$z+r\cu_0$, use the energies~$Y_{N,h}$ defined in the proof of Lemma~\ref{mfd:lem-neumann-error}. The source-independent estimate of Proposition~\ref{mfd:cor-neumann-source} has Morrey constant at most~$C_{z,r,p}K_Nh^{-2}$ for this source, while the energies~$Y_{N,h}$ have uniformly bounded prescribed moments by the estimates for~$Y_{N,h}$ in that proof. Lemma~\ref{mfd:cor-14}, the proof of Lemma~\ref{mfd:lem-15} and Lemma~\ref{mfd:lem-resampling-sum} therefore give
\begin{equation*}
	\|Y_{N,h}-\E[Y_{N,h}\mid\g_{-j},\ |j|\leq\ell]\|_{L^p}
	\leq C(p,z,r)\gmcdelta\bigl(h^{-2}3^{-b\ell/8}+3^{-\ell}\bigr)\,.
\end{equation*}
The required finite moment orders are fixed before~$\gmcdelta$, independently of~$h$. The load approximation and conditional-expectation contraction yield
\begin{equation*}
	\|Y_N-\E[Y_N\mid\g_{-j},\ |j|\leq\ell]\|_{L^p}
	\leq C(p,z,r)\bigl(h^{\nf14}
	+\gmcdelta h^{-2}3^{-b\ell/8}+\gmcdelta3^{-\ell}\bigr)\,.
\end{equation*}
Take~$h=3^{-b\ell/32}/16$. Since~$0<b<1$ and~$0<\gmcdelta\leq1$, this is at most~$C(p,z,r)3^{-b\ell/128}$. Thus the second hypothesis of the proposition holds with~$a=b/128$ and~$C_{z,r}=C(p,z,r)/\gmcdelta$, or with any smaller positive~$a$. No bound on~$C_{z,r}$ uniform as~$\gmcdelta\downarrow0$ is asserted here. The retained-potential argument applies to this fixed linear functional vanishing on constants. This proves the general fixed-cube Neumann case; for~$p<2$ apply the resampling estimates at order two and downgrade the norm.

\smallskip

For a countable collection, choose a common exponent~$a=b/256$. Apply the assumptions to the quadratic forms~$\langle f\pm f',G_N^{z+r\cu_0}(f\pm f')\rangle$; the chosen countable linear source family contains~$f+f'$ and~$f-f'$. Polarization gives~(ii), and a diagonal extraction over the countable collection gives~(iii). Translation and every fixed cutoff shift are included in that collection.
\end{proof}

\subsection{Subsequential limit forms}
\label{mfd:sec-local-form}

We construct limit forms on~$L^2(z+\cu_n,dx)$ for~$z\in\mathbb Q^d$ and~$n\in\mathbb Z$, along a subsequence. They are regular, strongly local Dirichlet forms (Propositions~\ref{mfd:prop-killed-inverse},~\ref{mfd:prop-regularity} and~\ref{mfd:prop-locality}), consistent under restriction to subcubes (Proposition~\ref{mfd:prop-killed-consistency}), and multiplying the coefficient by a positive continuous function multiplies their energy measures by the same function (Proposition~\ref{mfd:prop-21}). The minimal energies with prescribed boundary values converge to the minimal energy of the limit form on the cube over functions with the same boundary values (Proposition~\ref{mfd:prop-boundary}), so that harmonic replacement is available for the limit forms (Proposition~\ref{mfd:prop-gluing}). Since~$A_N$ has no upper bound uniform in~$N$, we localize with harmonic interpolants of cutoff functions (Lemma~\ref{mfd:lem-cutoffs}).

\smallskip

The cubes~$z+\cu_n$ have rational centers and triadic side lengths. They include the cubes~$\cu_j$ with~$j\geq0$ and all their triadic subcubes. Let~$\beta$ be the exponent of Lemma~\ref{mfd:lem-extension}. Fix~$\beta<\alpha<1$,~$d-1<t<d$ and~$0<\kappa<2\alpha-1$, with~$\alpha$ so close to~$1$ that~\eqref{eq:mfd-33} of Lemma~\ref{mfd:lem-affine} holds for some~$\gamma\in(0,1)$ and~$\zeta>0$, and with~$2(1-\alpha)+\kappa<\nf1{24}$. Choose a countable~$\mathbb Q$-linear family~$\mathcal D\subset C_c^\infty(\Rd)$ such that, for each cube, the functions of~$\mathcal D$ supported in the cube are dense in~$L^2$ and in~$C_c^\infty$ on that cube. This also gives density in~$C^\alpha_c$ for the~$C^\beta$ norm. For each finite union of closed rational cubes contained in a finite union of open rational cubes, include in~$\mathcal D$ a smooth cutoff function, that is, a function with values in~$[0,1]$ equal to one near the first set and zero outside the second. Also include the smooth functions used to build the cutoff functions~$\chi_{N,r}$ of Lemma~\ref{mfd:lem-cutoffs}(iii).

\smallskip

The \textbf{harmonic interpolant} of a function on a mesh is the function equal, on each cube of the mesh, to the~$A_N$-harmonic extension of the restriction of that function to the boundary of the cube.

\smallskip

We extract a subsequence jointly for the finite-volume energies for these cubes and data, the solutions, the environment and the random constants in~\eqref{eq:mfd-1}--\eqref{eq:mfd-5}, together with the coercivity constants~$Q_{z+\cu_n,N}$ in the proof of Lemma~\ref{mfd:lem-neumann-error} for the cubes~$z+\cu_n$. Include as single random variables each constant of~\eqref{eq:mfd-3}, which is valid on all grid cubes at once, the constants of Proposition~\ref{mfd:prop-growth} and Remark~\ref{mfd:rem-large-cube-growth}, which do not depend on the source, the corresponding constants for~$A_N^0$ on the cells~$z+\cu_n$ (see the discussion after Proposition~\ref{mfd:prop-boundary}), and~$3^{-\kappa N}\|A_N\|_{L^\infty(z+\cu_n)}$, whose moments are bounded by Lemma~\ref{mfd:lem-extremes}. Proposition~\ref{mfd:prop-response-compact} gives~$L^p$ convergence of the energies along a subsequence. The moment and H\"older bounds give tightness of all the remaining coordinates.

\smallskip

Apply Prokhorov's theorem and Skorokhod's representation theorem~\citep[Theorems~16.3 and~4.30]{Kallenberg} in the countable product of the Polish coordinate spaces. We obtain a further subsequence and a representation on a probability space where all coordinates converge almost surely. The energy limits retain their joint law with the environment and their measurability with respect to the random fields.

\smallskip

We work on the event of probability one on which all coordinates converge. The constants just included are bounded along the sequence in the Skorokhod representation. The bounds~\eqref{eq:mfd-3} hold for every rational mesh, with a constant depending on the mesh, and the bounds~\eqref{eq:mfd-4} hold for all sources and boundary data. Fix~$z_0\in\mathbb Q^d$ and~$n_0\in\mathbb Z$ and abbreviate~$\mathscr E_N\coloneqq \mathscr E_N^{z_0+\cu_{n_0}}$ and~$G_N\coloneqq G_N^{z_0+\cu_{n_0}}$.

\subsubsection{The Dirichlet Green operator and the limit form}

\begin{proposition}[Convergence of the Dirichlet Green operators and the limit form]
\label{mfd:prop-killed-inverse}
Along the sequence in the Skorokhod representation, the following statements hold.
\begin{enumerate}[label=(\roman*)]
\item \underline{\emph{Norm convergence}}: the operators~$G_N$ converge in operator norm on~$L^2(z_0+\cu_{n_0},dx)$ to a compact, positive, self-adjoint, injective operator~$G$.
\item \underline{\emph{Limit form}}: the form
\begin{equation}
	\mathscr E(u)=\sup_{f\in L^2(z_0+\cu_{n_0})}\bigl\{2\langle f,u\rangle-\langle f,Gf\rangle\bigr\}\,,
	\qquad D(\mathscr E)=\{u:\mathscr E(u)<\infty\}=\operatorname{Ran}G^{\nf12}\,,
	\label{eq:mfd-18}
\end{equation}
is a densely defined closed Dirichlet form on~$L^2(z_0+\cu_{n_0},dx)$.
\item \underline{\emph{Mosco convergence}}: if~$u_N\rightharpoonup u$ weakly in~$L^2(z_0+\cu_{n_0})$, then~$\mathscr E(u)\leq\liminf \mathscr E_N(u_N)$ (the liminf inequality). Every~$u\in D(\mathscr E)$ has a \textbf{recovery sequence}, that is, a sequence~$u_N\in H^1_0(z_0+\cu_{n_0})$ with~$u_N\to u$ in~$L^2(z_0+\cu_{n_0})$ and~$\mathscr E_N(u_N)\to \mathscr E(u)$.
\item \underline{\emph{Density}}: the functions~$Gf$ with~$f\in C_c^\infty(z_0+\cu_{n_0})$ are dense in~$D(\mathscr E)$ for the form norm.
\end{enumerate}
\end{proposition}

\begin{proof}
\emph{Step 1.} We show that~$G_N$ converges in operator norm. By~\eqref{eq:mfd-1} and the energy identity,
\begin{equation*}
	\|G_Nf\|_{H^{\nf34}}^2\leq Q_N\mathscr E_N(G_Nf)=Q_N\langle f,G_Nf\rangle\leq Q_N\|f\|_2\|G_Nf\|_2\,,
\end{equation*}
so~$\|G_N\|_{L^2\to H^{\nf34}}\leq\sup_NQ_N<\infty$ along the sequence in the Skorokhod representation. By Rellich's theorem the images of the unit ball under all~$G_N$ lie in a common precompact subset of~$L^2$. Given~$\ep>0$, choose a finite~$\ep$-net for these images consisting of elements of~$\mathcal D$, and let~$\Pi$ be the orthogonal projection onto its span. Then~$\|(I-\Pi)G_N\|\leq\ep$. By self-adjointness,~$\|G_N(I-\Pi)\|=\|(I-\Pi)G_N\|\leq\ep$. Hence~$\|G_N-G_M\|\leq\|\Pi(G_N-G_M)\Pi\|+4\ep$. The operator~$\Pi G_N\Pi$ is determined by the finite matrix of the bilinear form~$\langle\cdot,G_N\cdot\rangle$ on the elements of the net. Polarization expresses each entry as a real linear combination of quadratic forms~$\langle h,G_Nh\rangle$ with~$h\in\mathcal D$, which converge by the choice of the sequence in the Skorokhod representation. Hence~$(G_N)$ is Cauchy in operator norm.

\smallskip

\emph{Step 2.} We show that~$G$ is injective. For~$v\in H^1_0(z_0+\cu_{n_0})$ and~$s\in\mathbb R$, the variational formula~\eqref{eq:mfd-quadratic-response} gives~$\langle f,G_Nf\rangle\geq2s\langle f,v\rangle-s^2\mathscr E_N(v)$, and optimizing in~$s$ gives~$|\langle f,v\rangle|^2\leq\langle f,G_Nf\rangle \mathscr E_N(v)$. Let~$f\in L^2(z_0+\cu_{n_0})$ satisfy~$Gf=0$, let~$\phi\in\mathcal D\cap C_c^\infty(z_0+\cu_{n_0})$, and let~$v_{N,r}$ be the harmonic interpolant of~$\phi$ on the triadic subdivision of~$z_0+\cu_{n_0}$ into cubes of side~$r<\operatorname{dist}(\operatorname{supp}\phi,\partial(z_0+\cu_{n_0}))$. By~\eqref{eq:mfd-2} with~$\|\phi(z+r\,\cdot)\|_{C^\beta(\partial \cu_0)/\mathbb R}\leq Cr\|\nabla\phi\|_\infty$ for the centers~$z$ of the mesh cells,
\begin{equation*}
	\mathscr E_N(v_{N,r})=\sum_z\Lambda_{N,z+r\cu_0}(\phi|_{\partial(z+r\cu_0)})
	\leq Cr^{-d}\max_z\Lambda_{\upsilon/2,2}(z+r\cu_0;A_N)\,r^{d-2}\,r^2\|\nabla\phi\|_\infty^2\,,
\end{equation*}
which is bounded along the sequence in the Skorokhod representation at fixed~$r$, since only finitely many cells~$z+r\cu_0$ occur. The maximum principle on each cell gives~$\|v_{N,r}-\phi\|_\infty\leq Cr\|\nabla\phi\|_\infty$. Hence
\begin{equation*}
	|\langle f,\phi\rangle|\leq|\langle f,v_{N,r}\rangle|+Cr\|f\|_1\|\nabla\phi\|_\infty\,,
\end{equation*}
and the first term tends to zero as~$N\to\infty$, since~$\langle f,G_Nf\rangle\to\langle f,Gf\rangle=0$. Letting~$r\downarrow0$ gives~$f\perp C_c^\infty(z_0+\cu_{n_0})$, so~$f=0$.

\smallskip

\emph{Step 3.} The spectral theorem gives~$\mathscr E(u)=\|G^{-\nf12}u\|_2^2$ with domain~$\operatorname{Ran}G^{\nf12}$, a dense closed form domain. For every~$\lambda>0$, Step~1 gives
\begin{equation*}
	(G_N^{-1}+\lambda I)^{-1}=G_N(I+\lambda G_N)^{-1}
	\longrightarrow G(I+\lambda G)^{-1}=(G^{-1}+\lambda I)^{-1}
\end{equation*}
in operator norm. The resolvent criterion~\citep[Theorem~2.4]{KuwaeShioya2003} therefore gives Mosco convergence on~$L^2(z_0+\cu_{n_0},dx)$. Applying the liminf inequality to normal contractions of recovery sequences proves the Markov property of~$\mathscr E$. Finally,~$\mathscr E(u,Gf)=\langle u,f\rangle$; if~$u$ is orthogonal in the form norm to~$G\mathcal D$, then~$\langle(I+G)u,f\rangle=0$ for every~$f\in\mathcal D$, so~$u=0$. This proves~(ii)--(iv).
\end{proof}

Since~$\mathbb Q^d\times\mathbb Z$ is countable, the proposition gives, along the sequence in the Skorokhod representation, a limit form~$\mathscr E^{z+\cu_n}$ on every cube~$z+\cu_n$ with~$z\in\mathbb Q^d$ and~$n\in\mathbb Z$; we write~$\mathscr E\coloneqq \mathscr E^{z_0+\cu_{n_0}}$ when the cube is fixed.

\smallskip

The comparison arguments of Section~\ref{mfd:sec-compare} need the quadratic forms~$\langle f,G_N^{z_0+\cu_{n_0}}f\rangle$ to be bounded above and below with positive probability, uniformly in~$N$ and in the extracted limit. Both bounds follow from moment bounds.

\begin{lemma}[Positive and negative moments of the quadratic forms]
\label{mfd:lem-response-reciprocal}
Fix a cube~$z+R\cu_0$ of side at most one and a nonzero function~$f\in C_c^\infty(z+R\cu_0)$. For every~$p<\infty$, if~$\gmcdelta$ is sufficiently small, then
\begin{equation*}
	\sup_N\E\bigl[\langle f,G_N^{z+R\cu_0}f\rangle^p+\langle f,G_N^{z+R\cu_0}f\rangle^{-p}\bigr]<\infty\,.
\end{equation*}
The same bounds hold for every extracted limit, uniformly over the subsequence and over all values of~$\gmcdelta$ in a fixed, sufficiently small interval.
\end{lemma}

\begin{proof}
The coercivity bound~\eqref{eq:mfd-1} gives~$\langle f,G_N^{z+R\cu_0}f\rangle\leq Q_N\|f\|_2^2$. For the negative moments, fix a mesh of side~$r$ and let~$v_{N,r}$ be the harmonic interpolant of~$f$ on it, as in Step~2 of the proof of Proposition~\ref{mfd:prop-killed-inverse}. The maximum principle gives~$\|v_{N,r}-f\|_\infty\leq Cr\|\nabla f\|_\infty$. Fix~$r$ so small that~$\langle f,v_{N,r}\rangle\geq\frac{1}{2}\|f\|_2^2$ for every~$N$. The energy of~$v_{N,r}$ has a uniformly bounded~$p$th moment by~\eqref{eq:mfd-2}--\eqref{eq:mfd-3} on this fixed finite mesh. Optimize the variational formula~\eqref{eq:mfd-quadratic-response} over the multiples of~$v_{N,r}$ to obtain
\begin{equation*}
	\langle f,G_N^{z+R\cu_0}f\rangle^{-1}\leq\frac{\mathscr E_N^{z+R\cu_0}(v_{N,r})}{\langle f,v_{N,r}\rangle^2}
	\leq4\|f\|_2^{-4}\mathscr E_N^{z+R\cu_0}(v_{N,r})\,.
\end{equation*}
Fatou's lemma passes both moment bounds to each extracted limit.
\end{proof}

\subsubsection{Traces and localization}

We will apply the following trace estimate to the energy measures of the harmonic interpolants of cutoff functions (Lemma~\ref{mfd:lem-cutoffs}).

\begin{lemma}[Measure trace estimate]
\label{mfd:lem-19}
Let~$t>d-1$ and~$K<\infty$, and let~$\nu$ be a finite positive Borel measure on~$\Rd$ supported in~$\overline{z_0+\cu_{n_0}}$ such that~$\nu(B_r(x))\leq Kr^t$ for every~$x\in\overline{z_0+\cu_{n_0}}$ and~$0<r\leq1$. Then
\begin{itemize}
\item The map~$v\mapsto v$, from smooth functions on~$\overline{z_0+\cu_{n_0}}$ to~$L^2(\nu)$, extends by completion to~$H^{\nf12}(z_0+\cu_{n_0})$.
\item There exists~$C(d,t,z_0+\cu_{n_0})<\infty$ such that, for every~$v\in H^{\nf12}(z_0+\cu_{n_0})$,
\begin{equation}
	\int|v|^2\,d\nu\leq C(d,t,z_0+\cu_{n_0})\bigl(K+\nu(z_0+\cu_{n_0})\bigr)\|v\|_{H^{\nf12}(z_0+\cu_{n_0})}^2\,.
	\label{eq:mfd-19}
\end{equation}
\item Moreover, if~$v_N\to v$ in~$L^2(z_0+\cu_{n_0})$ and~$\sup_N\|v_N\|_{H^{\nf34}(z_0+\cu_{n_0})}<\infty$, then~$v_N\to v$ in~$H^{\nf12}(z_0+\cu_{n_0})$.
\end{itemize}
\end{lemma}

\begin{proof}
The Bessel kernel of order one is bounded by~$C|z|^{1-d}$ for~$|z|\leq1$ and is bounded outside the unit ball. The assumption~$\nu(B_r(x))\leq Kr^t$ and~$t>d-1$ therefore give
\begin{equation*}
	\sup_{x\in\operatorname{supp}\nu}\int\bigl(|x-y|^{1-d}\vee1\bigr)\,d\nu(y)
	\leq C\left(\nu(\overline{z_0+\cu_{n_0}})+K\sum_{n\geq0}2^{-n(t-d+1)}\right)
	\leq C\bigl(K+\nu(\overline{z_0+\cu_{n_0}})\bigr)\,.
\end{equation*}
Regard the Bessel potential of order~$\nf12$ as a map from~$L^2(dx)$ to~$L^2(\nu)$. The composition of this map with its adjoint has the Bessel kernel of order one as its kernel. This kernel is bounded by the integrand above. By Schur's test, the squared norm of the map is at most~$C(K+\nu(\overline{z_0+\cu_{n_0}}))$. The Bessel-potential description of~$H^{\nf12}(\Rd)$ and a bounded extension operator from the cube give~\eqref{eq:mfd-19} for smooth~$v$.

\smallskip

The boundary of~$z_0+\cu_{n_0}$ is~$\nu$-null: cover it with~$O(r^{1-d})$ balls. Hence~$\nu(\overline{z_0+\cu_{n_0}})=\nu(z_0+\cu_{n_0})$. Density of smooth functions in~$H^{\nf12}(z_0+\cu_{n_0})$ gives the asserted completion, and convergence in~$L^2(dx)$ identifies it with the ordinary integral when~$\nu$ has a bounded density with respect to Lebesgue measure. Finally, the same extension operator and Sobolev interpolation give
\begin{equation*}
	\|w\|_{H^{\nf12}(z_0+\cu_{n_0})}
	\leq C\|w\|_{L^2(z_0+\cu_{n_0})}^{\nf13}\|w\|_{H^{\nf34}(z_0+\cu_{n_0})}^{\nf23}\,.
\end{equation*}
Weak compactness first puts~$v$ in~$H^{\nf34}(z_0+\cu_{n_0})$; applying this inequality to~$w=v_N-v$ proves the last assertion.
\end{proof}

Along the sequence in the Skorokhod representation, recovery sequences and solutions with zero Dirichlet boundary values are bounded in~$H^{\nf34}$ by~\eqref{eq:mfd-1}, and the energy measures of the cutoff functions of Lemma~\ref{mfd:lem-cutoffs} below satisfy~\eqref{eq:mfd-4} with a common constant. Thus~\eqref{eq:mfd-19} controls~$\int(v_N-v)^2\,d\Gamma_N(\chi_N)$ uniformly in~$N$. We use this estimate when localizing recovery sequences in Proposition~\ref{mfd:prop-locality}.

\subsubsection{Regularity, locality, and restriction to subcubes}

To obtain energy measures for~$\mathscr E$, we show that~$\mathscr E$ is regular and strongly local. Both properties are proved by multiplying recovery sequences by cutoff functions. Smooth cutoff functions are not suitable, since~$A_N$ has no upper bound uniform in~$N$; we use instead the harmonic interpolants of cutoff functions, whose energies, constants in the Morrey bounds and H\"older norms are controlled cell by cell by~\eqref{eq:mfd-2}--\eqref{eq:mfd-5}.

Let~$F\subset O$ be finite unions of closed and open rational cubes, respectively, with~$F$ compact and~$\overline O\subset z_0+\cu_{n_0}$. Let~$\chi_N$ be the harmonic interpolant, on a mesh of cubes~$z+\cu_n$ with~$z\in\mathbb Q^d$ and~$n\in\mathbb Z$, of a smooth cutoff function from~$\mathcal D$, with values in~$[0,1]$, equal to one near~$F$ and compactly supported in~$O$. Choose the mesh so fine that every cell meeting the set where this function lies strictly between~$0$ and~$1$ is contained in~$O\setminus F$.

For~$r=3^{n_0-j}$ with~$j\in\N_0$, let~$\chi_{N,r}$ be the harmonic interpolant on the triadic mesh of~$z_0+\cu_{n_0}$ of side~$r$ of a smooth function~$\Rd\to[0,1]$ that vanishes within distance~$r$ of~$\partial(z_0+\cu_{n_0})$, equals one at distance at least~$3r$ from~$\partial(z_0+\cu_{n_0})$ and has gradient at most~$C/r$ on~$\Rd$.

\begin{lemma}[Harmonic cutoff functions]
\label{mfd:lem-cutoffs}
\begin{enumerate}[label=(\roman*)]
\item The function~$\chi_N$
\begin{itemize}
\item belongs to~$H^1_0(z_0+\cu_{n_0})\cap C(\overline{z_0+\cu_{n_0}})$,
\item satisfies~$0\leq\chi_N\leq1$,
\item equals one on a neighborhood of~$F$ and zero outside~$O$,
\item has~$\nabla\chi_N=0$ on every cell on which the smooth cutoff function is constant.
\end{itemize}
\item Along the sequence in the Skorokhod representation, the energy~$\mathscr E_N(\chi_N)$, the constant in the bound~\eqref{eq:mfd-4} for~$\Gamma_N(\chi_N)$ at all radii, and the~$C^\alpha$ norm of~$\chi_N$ are bounded. In particular, along a subsequence,~$\chi_N$ converges uniformly to a continuous function~$\chi$ that agrees with the original cutoff function on every cell where that function is constant.
\item With~$K_N$ from~\eqref{eq:mfd-3}, for every~$r=3^{n_0-j}$ with~$j\in\N_0$, we have~$0\leq\chi_{N,r}\leq1$ and
\begin{equation*}
	\mathscr E_N(\chi_{N,r})\leq C\bigl(K_N+3^{-\kappa N}\|A_N\|_{L^\infty(z_0+\cu_{n_0})}\bigr)r^{-1-\kappa}\,.
\end{equation*}
\end{enumerate}
\end{lemma}

\begin{proof}
On a cell on which the cutoff function is constant the harmonic extension is that constant, by uniqueness, and the maximum principle gives~$0\leq\chi_N\leq1$. Continuity across faces follows from the matching traces, and continuity up to the boundary from the H\"older bound for the cell solutions. At a fixed mesh there are finitely many cells, each with bounded~$\Lambda_{\upsilon/2,2}(z+\cu_n;A_N)$ and bounded constants in the Morrey and H\"older bounds along the sequence in the Skorokhod representation. Sum over the cells to obtain~(ii). The Arzel\`a--Ascoli theorem gives uniform convergence.

\smallskip

For~(iii), the maximum principle gives~$0\leq\chi_{N,r}\leq1$ as in~(i). There are~$O(r^{-(d-1)})$ cells on which the function defining~$\chi_{N,r}$ is not constant. On each such cell the rescaled function has~$C^\beta$ norm at most~$C$. If~$r\geq3^{-N}$, then~\eqref{eq:mfd-2} and~\eqref{eq:mfd-3} bound the energy on each such cell by~$C\Lambda_{\upsilon/2,2}(z+r\cu_0;A_N)r^{d-2}\leq CK_Nr^{d-2-\kappa}$. If~$r<3^{-N}$, the smooth function defining~$\chi_{N,r}$ is itself a competitor in the minimization on each cell, which gives
\begin{equation*}
	\mathscr E_N(\chi_{N,r})\leq C\|A_N\|_{L^\infty(z_0+\cu_{n_0})}r^{-1}
	\leq C\bigl(3^{-\kappa N}\|A_N\|_{L^\infty(z_0+\cu_{n_0})}\bigr)r^{-1-\kappa}\,.
\end{equation*}
Both the factor~$3^{-\kappa N}\|A_N\|_{L^\infty(z_0+\cu_{n_0})}$ and the constant of~\eqref{eq:mfd-3} were included in the joint extraction, and their sum bounds the energy for all widths~$r$.
\end{proof}

\begin{proposition}[Regularity]
\label{mfd:prop-regularity}
The Dirichlet form~$\mathscr E$ of Proposition~\ref{mfd:prop-killed-inverse} is regular.
\end{proposition}

\begin{proof}
For~$f\in\mathcal D$, the H\"older bound~\eqref{eq:mfd-5} gives~$u_N=G_Nf\to Gf$ uniformly on~$\overline{z_0+\cu_{n_0}}$, with zero boundary values. Let~$\chi_{N,r}$ be the cutoff function of Lemma~\ref{mfd:lem-cutoffs}(iii). Covering the boundary layer~$\{x\in z_0+\cu_{n_0}:\operatorname{dist}(x,\partial(z_0+\cu_{n_0}))<3r\}$ by~$O(r^{1-d})$ balls and using~\eqref{eq:mfd-4}--\eqref{eq:mfd-5} gives
\begin{align}
	\mathscr E_N(u_N(1-\chi_{N,r}))
	&\leq 2\Gamma_N(u_N)(\{\operatorname{dist}(\cdot,\partial(z_0+\cu_{n_0}))<3r\})
	+2\sup_{\operatorname{dist}(\cdot,\partial(z_0+\cu_{n_0}))<3r}|u_N|^2\mathscr E_N(\chi_{N,r})
	\notag \\ &
	\leq C(\omega,f)\bigl(r^{t-d+1}+r^{2\alpha-1-\kappa}\bigr)\,.
	\label{eq:mfd-20}
\end{align}
At fixed~$r$, pass to a subsequence along which~$\chi_{N,r}\to\chi_r$ uniformly. By~\eqref{eq:mfd-20} and the liminf inequality,~$Gf\chi_r\in D(\mathscr E)\cap C_c(z_0+\cu_{n_0})$ converges to~$Gf$ in the form norm as~$r\downarrow0$. Since the functions~$Gf$ are dense in~$D(\mathscr E)$ by Proposition~\ref{mfd:prop-killed-inverse}(iv),~$D(\mathscr E)\cap C_c(z_0+\cu_{n_0})$ is dense in~$D(\mathscr E)$ for the form norm.

\smallskip

For density in the uniform norm, take~$\phi\in\mathcal D\cap C_c^\infty(z_0+\cu_{n_0})$ and its harmonic interpolants on a mesh of side~$r$. Their bounded energies and H\"older norms give a uniform subsequential limit~$v_r\in D(\mathscr E)\cap C_c(z_0+\cu_{n_0})$. The maximum principle on each cell gives~$\|v_r-\phi\|_\infty\leq Cr\|\nabla\phi\|_\infty$. Since these~$\phi$ are dense in~$C_0(z_0+\cu_{n_0})$, the space~$D(\mathscr E)\cap C_c(z_0+\cu_{n_0})$ is uniformly dense as well.
\end{proof}

\begin{proposition}[Strong locality on~$D(\mathscr E)\cap C_c(z_0+\cu_{n_0})$]
\label{mfd:prop-locality}
Let~$u,v\in D(\mathscr E)\cap C_c(z_0+\cu_{n_0})$, and suppose that~$u$ equals a constant~$c$ on a neighborhood of~$\operatorname{supp}v$. Then~$\mathscr E(u,v)=0$.
\end{proposition}

\begin{proof}
Choose two cutoff functions from~$\mathcal D$ and let~$\psi_N,\chi_N$ be their harmonic interpolants from Lemma~\ref{mfd:lem-cutoffs}, with~$\psi_N=1$ near~$\operatorname{supp}v$,~$\chi_N=1$ near~$\operatorname{supp}\psi_N$, and~$\operatorname{supp}\chi_N\subset\{u=c\}^\circ$. For recovery sequences~$u_N,v_N$ of~$u$ and~$v$, consider the sequences~$c\chi_N+(1-\chi_N)u_N$ and~$\psi_Nv_N$. These are still recovery sequences. Indeed, they converge in~$L^2$ to~$u$ and~$v$, and the only extra gradient terms are controlled by
\begin{equation*}
	\int(u_N-c)^2\,d\Gamma_N(\chi_N)
	+\int v_N^2\,d\Gamma_N(\psi_N)\longrightarrow0\,.
\end{equation*}
On the support of~$\nabla\chi_N$ we have~$u=c$, and on that of~$\nabla\psi_N$ we have~$v=0$. Apply the trace estimate~\eqref{eq:mfd-19} together with the~$H^{\nf12}$ convergence of the recovery sequences. Young's inequality gives~$\limsup \mathscr E_N(c\chi_N+(1-\chi_N)u_N)\leq(1+\ep)\mathscr E(u)$ for every~$\ep>0$, and likewise for~$\psi_Nv_N$. The liminf inequality of Proposition~\ref{mfd:prop-killed-inverse}(iii) gives equality. Since~$c\chi_N+(1-\chi_N)u_N=c$ near~$\operatorname{supp}\psi_N$, we have~$\mathscr E_N(c\chi_N+(1-\chi_N)u_N,\psi_Nv_N)=0$. Applying the liminf inequality to their sum and difference yields~$\mathscr E(u\pm v)\leq \mathscr E(u)+\mathscr E(v)$, hence~$\mathscr E(u,v)=0$.
\end{proof}

By regularity and Proposition~\ref{mfd:prop-locality}, uniqueness in the Beurling--Deny representation~\citep[Section~3.2]{FOT} makes the jump and killing measures vanish. Thus~$\mathscr E$ is strongly local, and its energy measures~$\Gamma_{\mathscr E}(u,v)$ satisfy the chain rule, Leibniz rule, Cauchy--Schwarz inequality and locality on open sets~\citep[Theorem~3.2.2 and Lemma~3.2.5]{FOT}.

\smallskip

We normalize the energy measure by~$\int\phi\,d\Gamma_{\mathscr E}(u)=\mathscr E(u,u\phi)-\frac{1}{2}\mathscr E(u^2,\phi)$ for bounded~$u\in D(\mathscr E)$ and~$\phi\in D(\mathscr E)\cap C_c(z_0+\cu_{n_0})$, so that the corresponding measure at finite~$N$ is~$\Gamma_N(u)=A_N|\nabla u|^2\,dx$, and~$\Gamma_{\mathscr E}(u)(z_0+\cu_{n_0})=\mathscr E(u)$ and~$\Gamma_{\mathscr E}(u,v)(z_0+\cu_{n_0})=\mathscr E(u,v)$. The measure~$\Gamma_{\mathscr E}(u)$ is supported in~$\operatorname{supp}u$.

\smallskip

\begin{proposition}[Consistency under restriction to subcubes]
\label{mfd:prop-killed-consistency}
Let~$z+\cu_n\subset z_0+\cu_{n_0}$ be an open cube with~$z\in\mathbb Q^d$ and~$n\in\mathbb Z$, and let~$\mathscr E^{z+\cu_n}$ and~$\mathscr E^{z_0+\cu_{n_0}}$ be the limit forms on~$z+\cu_n$ and~$z_0+\cu_{n_0}$ along the sequence in the Skorokhod representation. Then
\begin{itemize}
\item Extension by zero to~$z_0+\cu_{n_0}$ sends~$D(\mathscr E^{z+\cu_n})$ onto the~$\mathscr E^{z_0+\cu_{n_0}}$-closure of~$\{w\in D(\mathscr E^{z_0+\cu_{n_0}})\cap C_c(z_0+\cu_{n_0}):\operatorname{supp}w\subset z+\cu_n\}$.
\item For every~$u\in D(\mathscr E^{z+\cu_n})$, the energy~$\mathscr E^{z+\cu_n}(u)$ equals the~$\mathscr E^{z_0+\cu_{n_0}}$-energy of the extension of~$u$ by zero.
\end{itemize}
\end{proposition}

\begin{proof}
Zero extension of a recovery sequence on~$z+\cu_n$ shows that the~$\mathscr E^{z_0+\cu_{n_0}}$-energy of the extension of~$u$ by zero is at most~$\mathscr E^{z+\cu_n}(u)$. Conversely, fix~$w\in D(\mathscr E^{z_0+\cu_{n_0}})\cap C_c(z+\cu_n)$. Localize a recovery sequence on~$z_0+\cu_{n_0}$ by a harmonic cutoff function equal to one near~$\operatorname{supp}w$ and supported in~$z+\cu_n$. The trace argument in Proposition~\ref{mfd:prop-locality} shows that this is still a recovery sequence. The liminf inequality on~$z+\cu_n$ (Proposition~\ref{mfd:prop-killed-inverse}(iii)) gives~$\mathscr E^{z+\cu_n}(w|_{z+\cu_n})\leq \mathscr E^{z_0+\cu_{n_0}}(w)$. Hence equality holds on~$D(\mathscr E^{z_0+\cu_{n_0}})\cap C_c(z+\cu_n)$. Regularity of~$\mathscr E^{z+\cu_n}$ (Proposition~\ref{mfd:prop-regularity}) and the first inequality extend equality to all of~$D(\mathscr E^{z+\cu_n})$. The resulting isometry has closed image, and this image is exactly the stated closure.
\end{proof}

\subsubsection{Continuous weights and minimal energies with given boundary values}

For~$\varrho\in C(\overline{z_0+\cu_{n_0}})$ with~$\varrho>0$, let~$\mathscr E_N^{\varrho}(u)\coloneqq\int_{z_0+\cu_{n_0}}\varrho\,d\Gamma_N(u)$ for~$u\in H^1_0(z_0+\cu_{n_0})$, and let~$G_N^{\varrho}$ be its Dirichlet Green operator.

\begin{lemma}[Form order implies measure order]
\label{mfd:lem-sincos}
Let~$\mathscr E,\mathscr F$ be regular strongly local Dirichlet forms on~$L^2(z_0+\cu_{n_0},dx)$ with~$D(\mathscr E)=D(\mathscr F)$, let~$O\subset z_0+\cu_{n_0}$ be open, and let~$C<\infty$ be such that~$\mathscr F(v)\leq C\mathscr E(v)$ for every~$v\in D(\mathscr E)\cap C_c(z_0+\cu_{n_0})$ with~$\operatorname{supp}v\subset O$. Then~$\Gamma_{\mathscr F}(u)\leq C\Gamma_{\mathscr E}(u)$ on~$O$ for every~$u\in D(\mathscr E)\cap C_c(z_0+\cu_{n_0})$.
\end{lemma}

\begin{proof}
Let~$v\in D(\mathscr E)\cap C_c(z_0+\cu_{n_0})$ with~$\operatorname{supp}v\subset O$, and let~$k\geq1$. The functions~$k^{-1}v\sin(ku)$ and~$k^{-1}v\cos(ku)=k^{-1}v(\cos(ku)-1)+k^{-1}v$ belong to~$D(\mathscr E)\cap C_c(z_0+\cu_{n_0})$ and they are supported in~$\operatorname{supp}v$. For either form, the chain and Leibniz rules show that the cross terms in the energy measures of these two functions cancel in the sum, so that
\begin{equation*}
	\mathscr E(k^{-1}v\sin ku)+\mathscr E(k^{-1}v\cos ku)=\int v^2\,d\Gamma_{\mathscr E}(u)+k^{-2}\mathscr E(v)\,,
\end{equation*}
and likewise for~$\mathscr F$. Applying~$\mathscr F\leq C\mathscr E$ to both functions and letting~$k\to\infty$ gives~$\int v^2\,d\Gamma_{\mathscr F}(u)\leq C\int v^2\,d\Gamma_{\mathscr E}(u)$. For nonnegative~$\phi\in C_c(O)$, choose by regularity a function in~$D(\mathscr E)\cap C_c(z_0+\cu_{n_0})$ with values in~$[0,1]$ that equals one on~$\operatorname{supp}\phi$ and is supported in~$O$. Also choose functions in~$D(\mathscr E)\cap C_c(z_0+\cu_{n_0})$ converging uniformly to~$\phi^{\nf12}$. Then the products of the latter with the former are admissible and their squares converge uniformly to~$\phi$. Hence~$\int\phi\,d\Gamma_{\mathscr F}(u)\leq C\int\phi\,d\Gamma_{\mathscr E}(u)$.
\end{proof}

\begin{proposition}[Continuous weights]
\label{mfd:prop-21}
Let~$\varrho\in C(\overline{z_0+\cu_{n_0}})$ be positive. Along the sequence in the Skorokhod representation, the forms~$\mathscr E_N^\varrho$ Mosco-converge to a form~$\mathscr E^\varrho$ such that~$D(\mathscr E^\varrho)=D(\mathscr E)$ and, for all~$u,v\in D(\mathscr E)$,
\begin{equation}
	\mathscr E^{\varrho}(u,v)=\int_{z_0+\cu_{n_0}}\varrho\,d\Gamma_{\mathscr E}(u,v)\,.
	\label{eq:mfd-21}
\end{equation}
\end{proposition}

\begin{proof}
The bounds~$(\inf\varrho)\mathscr E_N\leq \mathscr E_N^\varrho\leq(\sup\varrho)\mathscr E_N$ show that the forms~$\mathscr E_N^\varrho$ satisfy the~$H^{\nf34}$ coercivity bound~\eqref{eq:mfd-1} up to a constant factor, and that the order reverses for their Green operators:~$(\sup\varrho)^{-1}G_N\leq G_N^\varrho\leq(\inf\varrho)^{-1}G_N$. By the finite-rank compactness argument of Proposition~\ref{mfd:prop-killed-inverse}, every subsequence therefore has a further subsequence along which~$G_N^\varrho$ converges in norm on all cubes~$z+\cu_n$ with~$z\in\mathbb Q^d$ and~$n\in\mathbb Z$. The reversed order makes each limit injective. The same proposition constructs their Mosco limit~$\mathscr F$, with~$D(\mathscr F)=D(\mathscr E)$ and equivalent form norms. Therefore~$\mathscr F$ is regular, and~$D(\mathscr F)\cap C_c(z_0+\cu_{n_0})=D(\mathscr E)\cap C_c(z_0+\cu_{n_0})$. The proofs of locality and consistency under restriction to smaller cubes also apply: recovery sequences retain their~$H^{\nf{3}{4}}$ bound, and the energy measures of the cutoff functions are at most~$(\sup\varrho)\Gamma_N(\chi_N)$. Thus~$\mathscr F$ is strongly local and its restrictions to smaller cubes satisfy~$(\inf_{z+\cu_n}\varrho)\mathscr E^{z+\cu_n}\leq \mathscr F^{z+\cu_n}\leq(\sup_{z+\cu_n}\varrho)\mathscr E^{z+\cu_n}$ for~$z\in\mathbb Q^d$ and~$n\in\mathbb Z$. By Lemma~\ref{mfd:lem-sincos}, the same comparison holds for energy measures on~$z+\cu_n$. For~$u\in D(\mathscr E)\cap C_c(z_0+\cu_{n_0})$, take a finite partition of unity near~$\operatorname{supp}u$ subordinate to cubes on which~$\operatorname{osc}\varrho\leq\ep$. Integrate against this partition of unity to obtain~$|\mathscr F(u)-\int\varrho\,d\Gamma_{\mathscr E}(u)|\leq\ep \mathscr E(u)$. Let~$\ep\downarrow0$, extend by density in the form norm and polarize. Every subsequential limit form is therefore the form in~\eqref{eq:mfd-21}, which proves convergence of the whole sequence in the Skorokhod representation.
\end{proof}

We also obtain convergence in probability of the weighted Dirichlet Green operators along the subsequence extracted at the beginning of this subsection. Indeed, given a subsequence, we represent the weight, the limit form~$\mathscr E$ and the preceding coordinates jointly on one probability space. The comparison~\eqref{eq:mfd-mult} absorbs the error in the uniform approximation of the weight. The proof of Proposition~\ref{mfd:prop-21} identifies every limit with~\eqref{eq:mfd-21}. The subsequence criterion for convergence in probability then gives the assertion, also for random positive continuous weights and for countable families of them.

\begin{lemma}[Bounded weights and deletion of a field]
\label{mfd:lem-borel-weights}
Let~$g$ be a bounded Borel function on~$z_0+\cu_{n_0}$. Then
\begin{itemize}
\item The formula~$\mathscr E^g(u,v)\coloneqq\int e^g\,d\Gamma_{\mathscr E}(u,v)$ defines a regular strongly local Dirichlet form with domain~$D(\mathscr E)$.
\item Its energy measures are~$\Gamma_{\mathscr E^g}(u,v)=e^g\,\Gamma_{\mathscr E}(u,v)$.
\item For every field~$\g_{-j}$ with~$j\geq0$, the form~$e^{-\g_{-j}}\mathscr E$, given by~$(u,v)\mapsto\int e^{-\g_{-j}}\,d\Gamma_{\mathscr E}(u,v)$, is measurable with respect to the remaining fields.
\end{itemize}
\end{lemma}

\begin{proof}
The bounds~$e^{-\|g\|_\infty}\mathscr E\leq \mathscr E^g\leq e^{\|g\|_\infty}\mathscr E$ show that~$\mathscr E^g$ is closed on~$D(\mathscr E)$, with a form norm equivalent to that of~$\mathscr E$, so that~$\mathscr E^g$ is regular, like~$\mathscr E$. The chain rule for energy measures gives the Markov property, and the locality of energy measures gives strong locality~\citep[Theorem~3.2.2 and Lemma~3.2.5]{FOT}. The defining identity of the energy measure and the Leibniz rule give the formula for~$\Gamma_{\mathscr E^g}$.

\smallskip

For the deletion, the weight~$e^{-\g_{-j}}$ is continuous and bounded above and below on each fixed cube. The coefficient~$e^{-\g_{-j}}A_N$ is independent of the deleted field for all sufficiently large~$N$. Proposition~\ref{mfd:prop-21}, applied to random weights as explained after its proof, identifies the limit of the Dirichlet Green operators of~$e^{-\g_{-j}}A_N$ with the Dirichlet Green operator of~$e^{-\g_{-j}}\mathscr E$. Hence each quadratic form~$\langle f,G^{z+\cu_n}f\rangle$ of this Dirichlet Green operator, with~$f\in\mathcal D$,~$z\in\mathbb Q^d$ and~$n\in\mathbb Z$, is a limit in probability of functions of the remaining fields. It is therefore measurable with respect to the remaining fields. These quadratic forms determine the form~$e^{-\g_{-j}}\mathscr E$.
\end{proof}

Therefore, conditionally on the remaining fields, the forms corresponding to coefficients in which the deleted field is multiplied by the cutoff functions~$\chi_i$ of the proof of Lemma~\ref{mfd:lem-15} can be defined by weighting the fixed form~$e^{-\g_{-j}}\mathscr E$ with the independent restrictions of the deleted field. The first part of the lemma allows these weights to be discontinuous.

\smallskip

\begin{corollary}[Energy-measure convergence]
\label{mfd:cor-energy-measures}
Along the sequence in the Skorokhod representation, the following statements hold, with weak convergence of finite measures on~$\overline{z_0+\cu_{n_0}}$.
\begin{enumerate}[label=(\alph*)]
\item If~$u_N\in H^1_0(z_0+\cu_{n_0})$ is a recovery sequence for~$u\in D(\mathscr E)$, then~$\Gamma_N(u_N)\rightharpoonup\Gamma_{\mathscr E}(u)$.
\item If~$u_N\to u$ in~$L^2$ and~$\sup_N\mathscr E_N(u_N)<\infty$, then every weak subsequential limit of~$\Gamma_N(u_N)$ is at least~$\Gamma_{\mathscr E}(u)$.
\item If~$(u_N)$ is as in~(b) and~$(v_N)$ is a recovery sequence for~$v\in D(\mathscr E)$, then~$\Gamma_N(u_N,v_N)\rightharpoonup\Gamma_{\mathscr E}(u,v)$.
\end{enumerate}
\end{corollary}

\begin{proof}
Let~$\phi\in C(\overline{z_0+\cu_{n_0}})$ be nonnegative. For a recovery sequence, the liminf inequalities of Proposition~\ref{mfd:prop-21} for the weighted forms with the positive weights~$1\pm\ep\phi$, where~$\ep>0$ and~$\ep\|\phi\|_\infty<1$, together with~$\mathscr E_N(u_N)\to \mathscr E(u)$, give both inequalities needed for~(a). For~(b), apply the liminf inequality for the weight~$\phi+\ep$, and let~$\ep\downarrow0$ along a subsequence on which the measures converge weakly.

\smallskip

For~(c), Cauchy--Schwarz bounds the total variations of the mutual energy measures. Pass to a joint subsequential limit of~$\Gamma_N(u_N)$ and~$\Gamma_N(u_N,v_N)$. Apply~(b) to~$u_N+sv_N$ and~(a) to~$v_N$. For every rational~$s$, this shows that the measure
\begin{equation*}
	\lim_N\Gamma_N(u_N)-\Gamma_{\mathscr E}(u)
	+2s\bigl(\lim_N\Gamma_N(u_N,v_N)-\Gamma_{\mathscr E}(u,v)\bigr)
\end{equation*}
is nonnegative. Testing against~$0\leq\phi\in C(\overline{z_0+\cu_{n_0}})$ and letting~$s\to\pm\infty$ identifies the limit of~$\Gamma_N(u_N,v_N)$ with~$\Gamma_{\mathscr E}(u,v)$. Every subsequential limit is therefore the same.
\end{proof}

For the regular strongly local form~$\mathscr E$, the one-dimensional energy-image-density theorem and the absolutely continuous chain rule~\citep[Proposition~1.16 and its proof]{Fitzsimmons2004} give~$\Gamma_{\mathscr E}(v)(\{v=c\})=0$ and~$\Gamma_{\mathscr E}(\Phi(v))=\Phi'(v)^2\Gamma_{\mathscr E}(v)$ for Lipschitz~$\Phi$ with~$\Phi(0)=0$, using quasi-continuous representatives. The reference normalizes the energy measure as~$2\Gamma_{\mathscr E}$, which does not affect either assertion.

\begin{lemma}[Zero-trace differences lie in~$D(\mathscr E^{z+\cu_n})$]
\label{mfd:lem-truncation}
Let~$z+\cu_n\subset z_0+\cu_{n_0}$ be a cube with~$z\in\mathbb Q^d$ and~$n\in\mathbb Z$, and let~$v\in D(\mathscr E^{z_0+\cu_{n_0}})\cap C(\overline{z_0+\cu_{n_0}})$ vanish on~$\partial(z+\cu_n)$. Then
\begin{itemize}
\item The function~$v_{z+\cu_n}\coloneqq v\indc_{\overline{z+\cu_n}}$ belongs to the space~$D(\mathscr E^{z+\cu_n})\subset D(\mathscr E^{z_0+\cu_{n_0}})$ of functions extended by zero.
\item Its energy is~$\mathscr E^{z_0+\cu_{n_0}}(v_{z+\cu_n})=\Gamma_{\mathscr E}(v)(z+\cu_n)$.
\end{itemize}
\end{lemma}

\begin{proof}
Since~$v=0$ on~$\partial(z+\cu_n)$, for every~$\ep>0$ the set~$\{|v|\geq\ep\}\cap\overline{z+\cu_n}$ is compact in~$z+\cu_n$. Truncate~$v_{z+\cu_n}$ at height~$\ep$, that is, replace it by~$\operatorname{sgn}(v_{z+\cu_n})(|v_{z+\cu_n}|-\ep)_+$, and multiply by a cutoff function in~$D(\mathscr E^{z_0+\cu_{n_0}})\cap C_c(z_0+\cu_{n_0})$ equal to one near this set and supported in~$z+\cu_n$. The result belongs to~$D(\mathscr E^{z+\cu_n})$ by consistency under restriction to smaller cubes (Proposition~\ref{mfd:prop-killed-consistency}). The chain rule and locality show that its energy measure is the restriction of~$\Gamma_{\mathscr E}(v)$ to~$(z+\cu_n)\cap\{|v|>\ep\}$. The energy of the difference of the truncations at two heights is the mass of~$\Gamma_{\mathscr E}(v)$ on the part of~$z+\cu_n$ where~$|v|$ lies between these heights. Since the level sets of~$v$ carry no energy, the truncations are Cauchy in the form norm and converge uniformly to~$v_{z+\cu_n}$. Since~$\mathscr E^{z+\cu_n}$ is closed,~$v_{z+\cu_n}\in D(\mathscr E^{z+\cu_n})$. Letting~$\ep\downarrow0$ in the identity for the energy measure of the truncation gives~$\mathscr E(v_{z+\cu_n})=\Gamma_{\mathscr E}(v)(z+\cu_n)$, with~$\Gamma_{\mathscr E}(v_{z+\cu_n})(\partial(z+\cu_n))=0$.
\end{proof}

We now identify the limits of the minimal energies~$\Lambda_{N,z+\cu_n}$ with given boundary values on cubes~$z+\cu_n$ with~$z\in\mathbb Q^d$ and~$n\in\mathbb Z$. The candidate is the energy in~$z+\cu_n$ of the limit~$U$ of the harmonic interpolants of a fixed smooth datum on the~$3^d$ cells around~$z+\cu_n$.

Let~$z+r\cu_0=z+\cu_n$ with~$z\in\mathbb Q^d$ and~$n\in\mathbb Z$, and let~$z+3r\cu_0$ be the cube formed by the~$3^d$ congruent cells around~$z+r\cu_0$. For a boundary datum~$g$ on~$\partial(z+r\cu_0)$ from~$\mathcal D$, fix a smooth extension of~$g$ with compact support in~$z+3r\cu_0$. Let~$U_N$ equal, on each of the~$3^d$ cells, the~$A_N$-harmonic extension of the restriction of this extension to the cell boundary.

\begin{proposition}[Limits of the minimal energies with given boundary values]
\label{mfd:prop-boundary}
Along the sequence in the Skorokhod representation,~$U_N$ converges uniformly on~$\overline{z+3r\cu_0}$ to a function~$U$ with the following properties.
\begin{enumerate}[label=(\roman*)]
\item \underline{\emph{Boundary values}}: the function~$U$ belongs to~$D(\mathscr E^{z+3r\cu_0})\cap C(\overline{z+3r\cu_0})$, and~$U=g$ on~$\partial(z+r\cu_0)$.
\item \underline{\emph{No energy on the faces}}: the energy measure of~$U$ satisfies~$\Gamma_{\mathscr E}(U)(\partial(z+r\cu_0))=0$.
\item \underline{\emph{Harmonicity}}:~$\mathscr E^{z+3r\cu_0}(U,\phi)=0$ for every~$\phi\in D(\mathscr E^{z+r\cu_0})$, extended by zero to~$z+3r\cu_0$.
\item \underline{\emph{Minimality}}: the minimal energies converge, and
\begin{multline*}
	\lim_N\Lambda_{N,z+r\cu_0}(g)=\Gamma_{\mathscr E}(U)(z+r\cu_0)\\
	=\min\bigl\{\Gamma_{\mathscr E}(V)(z+r\cu_0):V\in D(\mathscr E^{z+3r\cu_0})\cap C(\overline{z+3r\cu_0})\,,
	\; V=g\text{ on }\partial(z+r\cu_0)\bigr\}\,.
\end{multline*}
\end{enumerate}
\end{proposition}

\begin{proof}
The energy, H\"older norm and constant in the Morrey bound for~$U_N$ are bounded by summing the corresponding cell bounds. Hence every subsequence has a further subsequence with~$U_N\to U$ uniformly, where~$U\in D(\mathscr E^{z+3r\cu_0})$ and~$U$ equals the fixed extension of~$g$ on every cell boundary. We show that~$(U_N)$ is a recovery sequence for~$U$.

\smallskip

Let~$W_N$ be a recovery sequence for~$U$, let~$\Sigma$ be the union of the cell boundaries, and take on each cell the cutoff function of Lemma~\ref{mfd:lem-cutoffs}(iii) with~$s=3^{-j}r$,~$j\in\N_0$, in place of~$r$. Their sum~$\chi_{N,s}$ vanishes near~$\Sigma$ and equals one away from its~$3s$-neighborhood. Use~$\chi_{N,s}W_N+(1-\chi_{N,s})U_N$, which has the same cell-boundary values as~$U_N$, as a competitor. Cellwise minimality and convexity give, for~$\ep>0$,
\begin{equation*}
	\mathscr E_N(U_N)
	\leq(1+\ep)\bigl[\mathscr E_N(W_N)
	+\Gamma_N(U_N)(\{\operatorname{dist}(\cdot,\Sigma)<3s\})\bigr]
	+(1+\ep^{-1})
	\int(W_N-U_N)^2\,d\Gamma_N(\chi_{N,s})\,.
\end{equation*}
At fixed~$s$, the last term tends to zero by~\eqref{eq:mfd-19} and convergence in~$H^{\nf12}$. The~$\Gamma_N(U_N)$-measure of the set of points at distance less than~$3s$ from~$\Sigma$ is at most~$Cs^{t-d+1}$ by the Morrey bound at every radius and a covering of the finitely many faces. Let first~$N\to\infty$, then~$s\downarrow0$ and~$\ep\downarrow0$. The liminf inequality gives~$\mathscr E_N(U_N)\to \mathscr E(U)$.

\smallskip

Corollary~\ref{mfd:cor-energy-measures}(a) now gives~$\Gamma_N(U_N)\rightharpoonup\Gamma_{\mathscr E}(U)$. The common Morrey bound passes to the limit. Hence~$\Gamma_{\mathscr E}(U)(\Sigma)=0$ and~$\Lambda_{N,z+r\cu_0}(g)=\Gamma_N(U_N)(z+r\cu_0)\to\Gamma_{\mathscr E}(U)(z+r\cu_0)$. For~$\phi\in D(\mathscr E^{z+r\cu_0})$, a recovery sequence~$\phi_N$ on~$z+r\cu_0$, extended by zero, is also a recovery sequence on~$z+3r\cu_0$ by consistency under restriction to smaller cubes (Proposition~\ref{mfd:prop-killed-consistency}). The identity~$\mathscr E_N(U_N,\phi_N)=0$ at finite~$N$ passes to~$\mathscr E(U,\phi)=0$ by Corollary~\ref{mfd:cor-energy-measures}(c).

\smallskip

If~$V=g$ on~$\partial(z+r\cu_0)$, then by Lemma~\ref{mfd:lem-truncation} the function~$v=(V-U)\indc_{\overline{z+r\cu_0}}$ lies in~$D(\mathscr E^{z+r\cu_0})$. Its energy measure does not charge~$\partial(z+r\cu_0)$. Harmonicity therefore yields~$\Gamma_{\mathscr E}(V)(z+r\cu_0)=\Gamma_{\mathscr E}(U)(z+r\cu_0)+\mathscr E(v)\geq\Gamma_{\mathscr E}(U)(z+r\cu_0)$. The same argument on each cell shows that the difference of two possible limits has zero energy on every cell. Since the energy measures do not charge the faces, coercivity forces this difference to vanish on~$z+3r\cu_0$. Thus the uniform convergence and all the conclusions hold along the full sequence in the Skorokhod representation. The argument applies to every fixed finite partition into cubes with a common smooth datum vanishing on the outer boundary.
\end{proof}

We will also use this boundary identification after removing the factor~$e^H$. On the same realization on~$z+3r\cu_0$, let~$\varrho\coloneqq e^{-H}$ and let~$\mathscr E^{\varrho,z+3r\cu_0}$ be the form obtained from the Dirichlet Green operators on~$z+3r\cu_0$ by Proposition~\ref{mfd:prop-21}. Its energy measure is~$\Gamma_{\mathscr E^\varrho}=\varrho\Gamma_{\mathscr E}$, its domain is~$D(\mathscr E^{z+3r\cu_0})$, and its form norm is equivalent to that of~$\mathscr E^{z+3r\cu_0}$. Thus regularity and consistency under restriction to subcubes are preserved.

\smallskip

We give the cell bounds for~$A_N^0$ also on cells larger than one. The unit-scale bounds follow from the proof of Proposition~\ref{mfd:prop-growth} with~$H$ replaced by~$0$. On a fixed larger cell, join the harmonic extensions of the smooth datum on a finite mesh of cubes of side at most one. Their energy bounds control the energy of this competitor and hence the total energy of the cell minimizer. For a bounded volume source, the same bound follows from~\eqref{eq:mfd-quadratic-response} and the coercivity bound for larger cubes in the proof of Lemma~\ref{mfd:lem-neumann-error}.

\smallskip

For a cell~$z+\cu_j$ with~$j\geq0$, let~$f\in L^\infty(z+\cu_j)$ and~$\varphi\in C^2(\overline{z+\cu_j})$, and let~$u_N\in H^1(z+\cu_j)$ solve~$-\nabla\cdot(A_N^0\nabla u_N)=f$ in~$z+\cu_j$ with~$u_N=\varphi$ on~$\partial(z+\cu_j)$. Use the coordinates~$y=3^Nx$ and the relabelled fields of Step~1 of the proof of Proposition~\ref{mfd:prop-growth}. The whole cell becomes~$3^Nz+\cu_{N+j}$, and we have~$A_N^0(3^{-N}y)=\ahom_N^{-1}\a_N(y)$ and~$\b_{N,N+j}=\ahom_N$. Apply Proposition~\ref{p.cutoff.Holder.regularity} with cutoff~$L=N$, outer scale~$m=N+j$ and the exponent~$\alpha_M$ chosen in the growth proof. Its energy and oscillation estimates apply on the whole cell, including balls meeting its faces. The coefficient~$\ahom_N$ cancels with the scalar~$\ahom_N^{-1}$ when returning to the original energy. Lemma~\ref{mfd:lem-primitive} supplies the vector source when a volume source is present. Steps~2--3 of that proof therefore give, for every~$x\in z+\cu_j$ and every~$3^{-N}\leq r\leq1$, with~$\max\{t,d-2+2\alpha\}<t_1<d-2+2\alpha_M$ and~$\Gamma_N^0(u)\coloneqq A_N^0|\nabla u|^2\,dx$,
\begin{equation}
	\Gamma_N^0(u_N)(B_r(x)\cap (z+\cu_j))
	\leq K_N\bigl(\|f\|_\infty+\|\varphi\|_{C^2}\bigr)^2r^{t_1}
	\,.
	\label{eq:mfd-zero-large-cell-growth}
\end{equation}
The total energy bound treats radii beyond the regularity range, as in Step~3. The tail~\eqref{e.cutoff.Holder.minimal.scale} is uniform in~$N$,~$j$ and the deterministic translation; choose the moment orders for the total energy and~$3^{t_1\mathcal L_N(\alpha_M,N+j)}$ before decreasing~$\gmcdelta$, so that~\eqref{eq:mfd-tail-rate} holds at these orders. This bounds the listed moments of~$K_N$ uniformly in~$N$.

\smallskip

At radii below~$3^{-N}$, retain Steps~4--5 of the proof of Proposition~\ref{mfd:prop-growth}, with~$H=0$. Combine the small-oscillation estimate of those steps with~\eqref{eq:mfd-zero-large-cell-growth}; the positive powers~$3^{-(t_1-t)N}$,~$3^{-(d+2-t)N}$ and~$3^{-(d-t)N}$ absorb the coefficient and log-gradient moments at the finitely many enlarged orders. The odd reflections give the same bound at faces, edges and corners. Step~7, using the whole-cell oscillation estimate above, gives the~$C^\alpha$ bound. Together, these steps give the Morrey bound~\eqref{eq:mfd-4} for~$\Gamma_N^0$ at every radius~$r\leq1$ and the H\"older bounds of Proposition~\ref{mfd:prop-growth}(ii) for~$A_N^0$ on~$z+\cu_j$, with one random constant~$K_N$ that does not depend on~$f$ and~$\varphi$ and whose moments of the listed orders are bounded uniformly in~$N$. This is the whole-cell argument used in Step~3 of Lemma~\ref{tight:lem-static}, with~$\alpha_M$ in place of~$7/8$. Stationarity makes the moment constants uniform over deterministic translations of a fixed cell.

\smallskip

These constants, for the~$3^d$ cells around~$z+r\cu_0$, are among the constants for~$A_N^0$ included in the joint extraction of the subsequence. The cutoff energy measures are bounded by~$(\sup_{z+3r\cu_0}\varrho)\Gamma_N(\chi_{N,s})$. The same comparison with a recovery sequence for~$\mathscr E^{\varrho,z+3r\cu_0}$ makes the weighted harmonic interpolants a recovery sequence: the face-layer error is~$Cs^{t-d+1}$, and the cutoff term vanishes by~\eqref{eq:mfd-19}. The weighted energy-measure convergence then gives zero energy on the faces, harmonicity and minimality. In particular, with~$\Lambda^0_{N,z+r\cu_0}(g)$ the minimal energy for~$A_N^0$,
\begin{multline}
	\Lambda_{N,z+r\cu_0}^0(g) \\
	\longrightarrow\min\bigl\{\Gamma_{\mathscr E^\varrho}(V)(z+r\cu_0):
	V\in D(\mathscr E^{\varrho,z+3r\cu_0})\cap C(\overline{z+3r\cu_0}),\ V=g\text{ on }\partial(z+r\cu_0)\bigr\}\,.
	\label{eq:mfd-weighted-boundary-minimum}
\end{multline}
For a random~$H$, include it in the joint representation and use~\eqref{eq:mfd-mult} to absorb its uniform approximation. The subsequence criterion then gives convergence in probability on the original field space, after uniqueness. This applies in particular to the uses of~$e^{-H}$ below.

Let~$\beta$ and~$\alpha$ be as fixed in Section~\ref{mfd:sec-local-form}, so that~$1/2<\beta<\alpha<1$, and let~$\upsilon$ be the exponent of Lemma~\ref{mfd:lem-extension}.

For the gluing assertions~(ii) and~(iii) below, let~$w+\cu_n$ be a cube with~$w\in\mathbb Q^d$ and~$n\in\mathbb Z$, and let~$z_i+\cu_{n_i}$ for~$i=1,\dots,m$ be disjoint cells with~$z_i\in\mathbb Q^d$ and~$n_i\in\mathbb Z$ such that the closure of each~$z_i+\cu_{n_i}$ and the concentric cube~$z_i+\cu_{n_i+1}$ of three times its side length lie in~$w+\cu_n$. Let~$v\in D(\mathscr E^{w+\cu_n})\cap C(\overline{w+\cu_n})$ satisfy~$v|_{\partial(z_i+\cu_{n_i})}\in C^\alpha$ for every~$i$. Define~$v'$ to equal~$v$ outside~$\bigcup_i(z_i+\cu_{n_i})$ and, on each~$z_i+\cu_{n_i}$, to equal the harmonic extension~$U_i$ of~$v|_{\partial(z_i+\cu_{n_i})}$.

\begin{proposition}[H\"older boundary data, gluing, harmonic replacement]
\label{mfd:prop-gluing}
The following statements hold along the sequence in the Skorokhod representation.
\begin{enumerate}[label=(\roman*)]
\item Let~$z+r\cu_0$ and~$z+3r\cu_0$ be as in Proposition~\ref{mfd:prop-boundary}. For every~$g\in C^\alpha(\partial(z+r\cu_0))$:
\begin{itemize}
\item the limit~$\Lambda_{\mathscr E,z+r\cu_0}(g)=\lim_N\Lambda_{N,z+r\cu_0}(g)$ exists and satisfies, with~$C$ the constant of~\eqref{eq:mfd-2},
\begin{equation*}
	\Lambda_{\mathscr E,z+r\cu_0}(g)\leq C\bigl(\sup_N\Lambda_{\upsilon/2,2}(z+r\cu_0;A_N)\bigr)r^{d-2}\|g(z+r\,\cdot)\|_{C^\beta/\mathbb R}^2\,,
\end{equation*}
\item all conclusions of Proposition~\ref{mfd:prop-boundary} hold for~$g$, with a harmonic extension~$U\in D(\mathscr E^{z+3r\cu_0})\cap C(\overline{z+3r\cu_0})$.
\end{itemize}
\item For every choice of cells and~$v$ as above,~$v'\in D(\mathscr E^{w+\cu_n})\cap C(\overline{w+\cu_n})$ and
\begin{equation*}
	\mathscr E^{w+\cu_n}(v')=\mathscr E^{w+\cu_n}(v)-\sum_i\bigl[\Gamma_{\mathscr E}(v)(z_i+\cu_{n_i})-\Lambda_{\mathscr E,z_i+\cu_{n_i}}(v|_{\partial(z_i+\cu_{n_i})})\bigr]\leq \mathscr E^{w+\cu_n}(v)\,.
\end{equation*}
\item If moreover~$\Gamma_{\mathscr E}(v)$ charges no face~$\partial(z_i+\cu_{n_i})$, and the cells~$z_i+\cu_{n_i}$ partition a cube~$w'+R'\cu_0\subset w+\cu_n$, then~$\Gamma_{\mathscr E}(v')(w'+R'\cu_0)=\sum_i\Lambda_{\mathscr E,z_i+\cu_{n_i}}(v|_{\partial(z_i+\cu_{n_i})})$.
\end{enumerate}
\end{proposition}

\begin{proof}
(i) Extend~$g$ to a function in~$C_c^\alpha(z+3r\cu_0)$ and approximate it in~$C^\beta$ by functions in~$\mathcal D$. On the fixed finite partition, apply the bound~\eqref{eq:mfd-2} for H\"older boundary data and the maximum principle to two approximations. The bounds show that the energy~$\mathscr E_N$ of the difference of the corresponding interpolants is at most~$C$ times the squared~$C^\beta$ norm of the difference of the approximations. The sup norm of the difference of the interpolants is at most the sup norm of the difference of the approximations. The interpolants converge by Proposition~\ref{mfd:prop-boundary}, and the same bounds hold for their limits by the liminf inequality. Hence the corresponding harmonic extensions converge in energy and uniformly to a limit~$U$ with boundary values~$g$. Harmonicity passes to the limit by energy convergence. The energy on the faces also remains zero, since, by the Cauchy--Schwarz inequality for energy measures, energy measures depend continuously in total variation on the function in the form norm. Minimality follows again from harmonicity and Lemma~\ref{mfd:lem-truncation}. Since~$\sqrt{\Lambda_{N,z+r\cu_0}}$ is a seminorm of the datum, its reverse triangle inequality and~\eqref{eq:mfd-2} transfer the convergence of~$\Lambda_{N,z+r\cu_0}$ and the asserted bound from the approximations to~$g$. The argument holds for each realization and covers random boundary data.

\smallskip

(ii) Consistency under restriction to subcubes extends each~$U_i$ to~$D(\mathscr E^{w+\cu_n})$, and Lemma~\ref{mfd:lem-truncation} puts~$(U_i-v)\indc_{\overline{z_i+\cu_{n_i}}}$ in~$D(\mathscr E^{z_i+\cu_{n_i}})$ with zero energy on~$\partial(z_i+\cu_{n_i})$. Adding these differences to~$v$ shows that~$v'$ belongs to the form domain. Locality identifies its energy measure with~$\Gamma_{\mathscr E}(U_i)$ inside~$z_i+\cu_{n_i}$ and with~$\Gamma_{\mathscr E}(v)$ elsewhere, including the faces by Cauchy--Schwarz. Integrating gives the identity. Cell minimality gives the inequality.

\smallskip

(iii) The same identification and the assumption that~$\Gamma_{\mathscr E}(v)$ charges no face give the sum over the partition. For a common smooth datum, the energy measure of the limit does not charge the faces by Proposition~\ref{mfd:prop-boundary}; for a common~$C^\alpha$ datum, this follows from the approximation in energy in~(i).
\end{proof}

\begin{lemma}[Finite partitions with prescribed boundary values]
\label{mfd:lem-skeleton}
Let~$\alpha$ be as fixed in Section~\ref{mfd:sec-local-form}, let the cubes~$z_i+\cu_{n_i}$ with~$z_i\in\mathbb Q^d$ and~$n_i\in\mathbb Z$ form a finite partition of a cube~$w+\cu_n$ with~$w\in\mathbb Q^d$ and~$n\in\mathbb Z$, and let~$g\in C^\alpha(\overline{w+\cu_n})$ vanish on~$\partial(w+\cu_n)$. Let~$v$ be the function equal, on each cell~$z_i+\cu_{n_i}$, to the limiting harmonic extension of~$g|_{\partial(z_i+\cu_{n_i})}$. Then~$v\in D(\mathscr E^{w+\cu_n})\cap C(\overline{w+\cu_n})$,
\begin{equation*}
	\mathscr E^{w+\cu_n}(v)\leq\sum_i\Lambda_{\mathscr E,z_i+\cu_{n_i}}(g|_{\partial(z_i+\cu_{n_i})})\,,
\end{equation*}
and
\begin{equation*}
	\|v-g\|_\infty\leq\sup_i\operatorname{osc}_{z_i+\cu_{n_i}}g\,.
\end{equation*}
\end{lemma}

\begin{proof}
For smooth data compactly supported in~$w+\cu_n$, apply the last paragraph of the proof of Proposition~\ref{mfd:prop-boundary} to this finite partition. The harmonic interpolants converge uniformly and in energy, and their cell energies add. For general~$g$, use~$\beta$ as fixed in Section~\ref{mfd:sec-local-form}, with~$\beta<\alpha$. The zero extension of~$g$ can be approximated in~$C^\beta$ by~$C_c^\infty(w+\cu_n)$ functions. The bound~\eqref{eq:mfd-2} for H\"older boundary data and the maximum principle, as in Proposition~\ref{mfd:prop-gluing}(i), give convergence in energy and uniformly. This proves the energy bound, also for cells touching~$\partial(w+\cu_n)$. The maximum principle on each cell gives~$\|v-g\|_\infty\leq\sup_i\operatorname{osc}_{z_i+\cu_{n_i}}g$.
\end{proof}

\subsection{Estimates along chains of nested cubes}
\label{mfd:sec-branches}

In this subsection we bound the number of cubes of each chain of nested cubes on which the regularity, extension and harmonic-approximation estimates of Lemma~\ref{mfd:in-deterministic} fail, on a single event that does not depend on the solution (Proposition~\ref{mfd:prop-allchain}). The energy measure~$\Gamma_{\mathscr E}(u)$ and the events on which these estimates hold depend on the same fields. We therefore cover the failure of the estimates by events that depend only on the fields in a window of scales around that of the cube, with probabilities decaying exponentially in its width (Lemma~\ref{mfd:lem-witness}). At weak disorder, the resulting bound shows that, averaged over generations, the cubes on which the estimates fail carry an arbitrarily small proportion of the energy (inequality~\eqref{eq:mfd-24}).

\smallskip

Let~$\mathscr E$ be a subsequential limit form along the subsequence of cutoffs of Section~\ref{mfd:sec-local-form}, realized on the Skorokhod probability space.

\subsubsection{Continuous versions of the conditions}

Fix~$s\in(0,\nf1{32}]$ before choosing~$\gmcdelta$; Lemma~\ref{l.harmonic.approximation.good.scales.GMC} is applied with parameter~$8s\in(0,\nf14]$. Let~$\mathsf F_{m,z}$,~$\mathsf M_{m,z}$ and~$\mathsf R_{m,z}$ be the left sides of~\eqref{e.good.gj},~\eqref{e.good.discounted.am.vs.aj} and~\eqref{e.good.sum.J}, translated to~$z$. Here~$\mathsf M_{m,z}$ is a scalar statistic. For~$\gmcdelta=0$, they equal~$0$,~$2$ and~$0$, respectively. We use the quantity~$\mathbf G_{m,z}(s)$ of~\eqref{e.def.bfG.mhq}, suppressing its fixed argument~$s$, and the quantity~$\mathcal J_n(P)$ of~\eqref{e.normalized.J.maximum.definition}.

\smallskip

To measure failure of the bounds at center~$z$ and scale~$m$, apply continuous piecewise linear functions to~$\mathsf F_{m,z}$,~$\mathsf M_{m,z}$ and~$\mathsf R_{m,z}$, and let~$\mathsf Z_{m,z}$ be their sum. Each function increases linearly from zero to one between its lower and upper thresholds, is zero below the lower threshold and one above the upper threshold, and lies between the indicators of exceeding these thresholds. The threshold pairs are~$(\ep/2,\ep)$,~$(6,12)$ and~$(\ep^2/4,\ep^2)$, respectively, with the primitive tolerance~$\ep\in(0,1)$ fixed in the next subsubsection. Thus~$\mathsf Z_{m,z}$ takes values in~$[0,3]$; it is a continuous version of the number of the three conditions~\eqref{e.good.gj},~\eqref{e.good.discounted.am.vs.aj} and~\eqref{e.good.sum.J} that fail at center~$z$ and scale~$m$. The probability bounds of Lemmas~\ref{l.bad.scales.for.nabla.gj},~\ref{l.exp.goodscales} and~\ref{l.discounted.J.local} are applied at the lower thresholds. We call~$\mathsf F_{m,z}$,~$\mathsf M_{m,z}$,~$\mathsf R_{m,z}$,~$\mathsf Z_{m,z}$ and~$\mathbf G_{m,z}$ the primitive quantities. Every quantity~$\mathcal J$ entering~$\mathsf R$ or~$\mathbf G$ at stationary scale~$n$ is computed for~$\ahom_n^{-1}\a_n$, hence for~$A_n^0$ after translation and rescaling to the unit cube. The cutoff~$n$ of this coefficient is part of their definition. The normalized quantities on a physical cell, defined next, are different quantities. Where~$\mathsf Z_{m,z}<1$, each quantity lies below its upper threshold. Apply the deterministic good-scale estimates at these upper thresholds; in particular, they use the threshold~$12$ for~$\mathsf M_{m,z}$ in place of~$6$. Only the constants change, provided~$\tau^2\leq s\log3/16$: the product bound~$12$ gives~$e^{\|\g_{m,n}\|_\infty}\leq12\cdot3^{s(m-n)/8}$, while~\eqref{e.rho.m.n.bound} gives~$e^{|(\log(\ahom_n/\ahom_m)-(m-n)\tau^2)|}\leq3^{s(m-n)/16}$. The squared coefficient ratio in~\eqref{e.combined.coefficient.ratio.reg}, multiplied by the bound~\eqref{e.good.sum.J} and its weight, therefore still decays like~$3^{-s(m-n)}$.

\subsubsection{The good event of a cell}

We specify the quantities on a cell after dividing out its scalar~$s_N$. Use~$\beta$ fixed in Section~\ref{mfd:sec-local-form} with~$\upsilon$ as in Lemma~\ref{mfd:lem-extension}. The exponent pairs used below are~$(\upsilon/2,2)$ for boundary extension and~$(1,1)$ for the mean-zero Poincar\'e bound; the homogenization error uses~$(s,\infty,2)$. These exponents are fixed before the moment orders and the disorder threshold.

\smallskip

For a cube~$w+3^{-\ell}\cu_0$ with integer~$\ell\leq N$, including the negative depths of fixed enlargements, we use the cell scalar~$s_N(w+3^{-\ell}\cu_0)$ of Section~\ref{mfd:sec-environment}. For~$\ell\geq0$ this is exactly the scalar of Lemma~\ref{mfd:lem-local-normalizations}. For~$\ell<0$, its definition removes the values at~$w$ of the fields~$\g_1,\dots,\g_{-\ell}$. In the coordinates~$y$ with~$x=w+3^{-\ell}y$ the coefficient is~$s_N(w+3^{-\ell}\cu_0)^{-1}A_N(w+3^{-\ell}y)$. Retain the three ellipticity coordinates
\begin{equation*}
	\frac{\lambda_{\upsilon/2,2}(w+3^{-\ell}\cu_0;A_N)}{s_N(w+3^{-\ell}\cu_0)}\,,\qquad
	\frac{\Lambda_{\upsilon/2,2}(w+3^{-\ell}\cu_0;A_N)}{s_N(w+3^{-\ell}\cu_0)}\,, \qquad
	\frac{s_N(w+3^{-\ell}\cu_0)}{\lambda_{\upsilon/2,2}(w+3^{-\ell}\cu_0;A_N)}\,.
\end{equation*}
Also retain the matrices~$s_N(w+3^{-\ell}\cu_0)^{-1}\a(w+3^{-\ell}\cu_0;A_N)$ and~$s_N(w+3^{-\ell}\cu_0)\a_*^{-1}(w+3^{-\ell}\cu_0;A_N)$, entry by entry. For scalar comparisons let~$w'+3^{-\ell'}\cu_0$ range over the cubes of the catalogue fixed below, in particular~$z+3^{-k}\cu_0$, its concentric enlargement~$z+3^{1-k}\cu_0$ and the \textbf{comparison cube}, the cube concentric with~$z+3^{-k}\cu_0$ of side~$3^{\lfloor\gamma H_1\rfloor+4-k}$, with~$\gamma$ and~$H_1$ of Lemma~\ref{mfd:lem-affine}. For each such cube retain~$\mathcal E_{s,\infty,2}(\cu_0,0;s_N(w'+3^{-\ell'}\cu_0)^{-1}A_N(w'+3^{-\ell'}\,\cdot),\Id)$ and~$s_N(z+3^{-k}\cu_0)/s_N(w'+3^{-\ell'}\cu_0)$. The ellipticity coordinates have value~$1$, the matrices have value~$\Id$, the error has value~$0$ and the ratio has value~$1$ when~$\gmcdelta=0$. The finite tests use these normalized coordinates and their deviations from these values. By Lemma~\ref{mfd:lem-local-normalizations}, for~$\ell\geq0$ the coefficient~$s_N(w+3^{-\ell}\cu_0)^{-1}A_N(w+3^{-\ell}\,\cdot)$ in these coordinates is the coefficient~$A_{N-\ell}$ in the translated and rescaled environment; it includes the centered increments of the fields of scale larger than the side of~$w+3^{-\ell}\cu_0$. For~$\ell<0$, the definition of~$s_N$ gives the same identification. In particular, its matrices are not those of~$A_{N-\ell}^0$. The estimates for the primitive quantities and for these normalized arrays are proved separately in Lemma~\ref{mfd:lem-band}.

\smallskip

The pair~$(1,1)$ requires no additional random test. Weighted Cauchy--Schwarz in~\eqref{e.coarse.grained.ellipticity} gives, for every cube and coefficient,
\begin{equation*}
	\Lambda_{1,1}\leq C_\upsilon\Lambda_{\upsilon/2,2},\qquad
	\lambda_{1,1}^{-1}\leq C_\upsilon\lambda_{\upsilon/2,2}^{-1}\,,
	\qquad
	C_\upsilon=\frac{(1-3^{-1})^2}
	{(1-3^{-2+\upsilon})(1-3^{-\upsilon})}\,.
\end{equation*}
Indeed, apply Cauchy--Schwarz to~$\sum_{j\geq0}3^{-j}M_j^{\nf12} =\sum_{j\geq0}3^{-(1-\upsilon/2)j}(3^{-\upsilon j}M_j)^{\nf12}$, where~$M_j$ is the maximum matrix norm at relative depth~$j$ in either definition. This supplies both required exponent pairs on every enlargement from the same normalized ellipticity arrays.

\smallskip

Fix a rational cube and its triadic subdivision. For a cell~$z+3^{-k}\cu_0$, use a fixed finite catalogue of concentric enlargements and the cubes of finitely many translated triadic grids needed for the harmonic comparison. A catalogue cube has depth~$\ell=k-i$ for a fixed integer offset~$i$ and center~$z+3^{-\ell}v$ for a fixed grid shift~$v$. Include~$z+3^{-k}\cu_0$, its concentric enlargement of side~$3\cdot3^{-k}$ and the comparison cube. Choose the finite grids, with shifts~$v\in\{-\nf12,0,\nf12\}^d$ for the catalogue cubes of side at most~$3\cdot3^{-k}$, so that every cube~$x+\rho\cu_0$ with~$x\in\overline{z+3^{-k}\cu_0}$ and~$0<\rho\leq3^{-k}$ is contained in a \textbf{descendant}, a cube of the triadic subdivision of a catalogue cube, of side at most~$9\rho$. The catalogue cubes of side at most~$3\cdot3^{-k}$, that is, all except the comparison cube, lie in the concentric cube of side~$6\cdot3^{-k}$. Items~(i) and~(ii) of Lemma~\ref{mfd:in-deterministic} use only these cubes, so for them it suffices that the cube of the lemma contain their closures; the comparison cube is used only in item~(iii). At descendant depth~$D$, use all descendant centers and the finitely many catalogue centers needed for the scalar comparisons. Their number is at most~$C_d3^{dD}$.

\smallskip

Fix integers~$k_1\geq0$ and~$k_0\geq\max\{k_1,1\}$. At each observation center~$w$, the two prefixes are the sums over the physical depths~$j=\ell-k_1,\ldots,\ell+D$, that is, over the depths of a catalogue cube and its descendants down to relative depth~$D$, together with the~$k_1$ coarser depths~$\ell-k_1,\ldots,\ell-1$. The prefix conditions are imposed only for~$D\geq k_0$, so each prefix has at most~$D+k_1+1\leq3D$ terms. Their summands are, respectively, the primitive quantities~$\mathsf Z_{N-j,\,3^Nw}$ and~$\mathbf G_{N-j,\,3^Nw}$ in the stationary coordinates of~\eqref{eq:mfd-normalization}. Use these prefixes when~$\ell+D\leq N$, so that every index~$N-j$ is nonnegative; otherwise let them be zero and use the estimates below the smallest field scale.

\smallskip

Choose fixed numbers
\begin{equation*}
	0<\lambda_{\rm cut}<\lambda_{\rm lim}/4<\lambda_{\rm lim}<\lambda_{\rm det}<1\,,
	\qquad 0<\ep_{\rm ell}<1,\qquad
	\ep_{\rm hom}>0,\qquad c_{\rm det}>0\,.
\end{equation*}
Fix the primitive tolerance~$\ep\in(0,1)$ sufficiently small for the deterministic estimates at the chosen regularity exponents, and then choose~$\lambda_{\rm det}$ sufficiently small for those estimates, and~$\lambda_{\rm lim}$ and~$\lambda_{\rm cut}$ satisfying the displayed inequalities. In particular, after reducing~$\gmcdelta$, the quantity~$C(\lambda_{\rm det}+\gmcdelta^2+\ep^8)$ is smaller than the gap between the iteration exponent and~$\alpha$. Choose~$c_{\rm det}$ small enough that the hypothesis of the harmonic comparison of Lemma~\ref{l.cutoff.regularity.good.scale.estimates} holds even at error~$2\ep_{\rm hom}c_{\rm det}$; this leaves room for passage from the limiting error to a finite cutoff. The finite-cutoff prefix threshold is~$\lambda_{\rm cut}D$; the limiting prefix threshold is~$\lambda_{\rm lim}D$, strictly below the threshold~$\lambda_{\rm det}D$ used by the iteration of Lemma~\ref{l.iteration.lemma.GMC}. The \textbf{good event}~$\mathcal G_{\mathscr E}(z+3^{-k}\cu_0)$ is defined using the limiting arrays of these exact finite-cutoff quantities along the subsequence defining~$\mathscr E$. It requires, for every catalogue cube, every~$D\geq k_0$ and every observation center, both limiting prefixes to be smaller than~$\lambda_{\rm lim}D$, and the following finite tests:
\begin{enumerate}[label=(\alph*),leftmargin=*]
\item On every catalogue enlargement~$w+3^{-\ell}\cu_0$,~$\ep_{\rm ell}\leq\lim_N\lambda_{\upsilon/2,2}(w+3^{-\ell}\cu_0;A_N)/s_N(w+3^{-\ell}\cu_0)$ and~$\lim_N\Lambda_{\upsilon/2,2}(w+3^{-\ell}\cu_0;A_N)/s_N(w+3^{-\ell}\cu_0)\leq\ep_{\rm ell}^{-1}$. Consequently, the normalized~$(1,1)$ lower and upper bounds are~$\ep_{\rm ell}/C_\upsilon$ and~$C_\upsilon\ep_{\rm ell}^{-1}$. For the limit~$\lim_Ns_N(w+3^{-\ell}\cu_0)^{-1}\a(w+3^{-\ell}\cu_0;A_N)$ these bounds also give, for every~$p\in\Rd$,
\begin{multline*}
	\frac{\ep_{\rm ell}^2}{d}\operatorname{tr}\bigl(\lim_Ns_N(w+3^{-\ell}\cu_0)^{-1}\a(w+3^{-\ell}\cu_0;A_N)\bigr)|p|^2\\
	\leq p\cdot\bigl(\lim_Ns_N(w+3^{-\ell}\cu_0)^{-1}\a(w+3^{-\ell}\cu_0;A_N)\bigr)p\,.
\end{multline*}
\item On the comparison cube~$w'+3^{-\ell'}\cu_0$, the limiting error obeys
\begin{equation*}
	\lim_N\mathcal E_{s,\infty,2}\bigl(\cu_0,0;s_N(w'+3^{-\ell'}\cu_0)^{-1}A_N(w'+3^{-\ell'}\,\cdot),\Id\bigr)
	\leq\ep_{\rm hom}c_{\rm det}\,.
\end{equation*}
\item For every catalogue cube~$w'+3^{-\ell'}\cu_0$ the limiting ratio~$s_N(z+3^{-k}\cu_0)/s_N(w'+3^{-\ell'}\cu_0)$ belongs to~$(\nf12,2)$, strictly inside the closed interval used by the deterministic comparison.
\end{enumerate}
At a finite cutoff, write~$\mathcal G_N^{\rm det}(z+3^{-k}\cu_0)$ for the conditions used by the deterministic estimates: prefix threshold~$\lambda_{\rm det}D$, ellipticity bounds~$\ep_{\rm ell}/2$ and~$2\ep_{\rm ell}^{-1}$, error bound~$2\ep_{\rm hom}c_{\rm det}$ and closed ratio interval~$[\nf12,2]$. Each fixed finite set of these conditions follows eventually from~$\mathcal G_{\mathscr E}(z+3^{-k}\cu_0)$ under almost-sure convergence of the arrays. The smaller prefix threshold~$\lambda_{\rm cut}D$ is used to prove the probability bound passed to the limit at the threshold~$\lambda_{\rm lim}D/4$; eventual satisfaction of the threshold~$\lambda_{\rm cut}D$ is not needed. The catalogue and the finite collection of tests are chosen before~$N$,~$k$, the source and the environment. The ellipticity and matrix conditions of these tests are used for extension and comparison on a good cell. The regularity iteration itself uses only the primitive prefixes, together with its reference-ratio estimates; it does not require transfer of the normalized tests~(a)--(c).

\smallskip

For continuous tests in~(a)--(c), use piecewise linear ramps with the following inner and outer thresholds. For the lower ellipticity coordinate they are~$(1+\ep_{\rm ell})/2$ and~$\ep_{\rm ell}$, with the ramp increasing as the coordinate decreases; for the upper coordinate they are~$(1+\ep_{\rm ell}^{-1})/2$ and~$\ep_{\rm ell}^{-1}$; for the error they are~$\ep_{\rm hom}c_{\rm det}/2$ and~$\ep_{\rm hom}c_{\rm det}$. For the lower ratio test they are~$3/4$ and~$1/2$, and for the upper ratio test~$3/2$ and~$2$. The outer thresholds are the bounds defining~$\mathcal G_{\mathscr E}(z+3^{-k}\cu_0)$; the inner thresholds place the zero-disorder coordinates in the regions where their respective ramps vanish. The deterministic finite-cutoff bounds are weaker, as specified above. Each ramp is zero at the zero-disorder value of its coordinate and one when the corresponding outer condition fails. Finite sums of these ramps, or of absolute values of the centered normalized coordinates, are fixed nonnegative Lipschitz functions of the finite array. Their constants depend only on the displayed thresholds and catalogue.

Fix~$N$ and a cell~$z+3^{-k}\cu_0$ on which~$\mathcal G_N^{\rm det}(z+3^{-k}\cu_0)$ holds, and fix a cube containing the closures of its catalogue cubes. Let~$u\in H^1$ be a weak solution of~$-\nabla\cdot(A_N\nabla u)=f$ in this cube, with~$f\in L^{d/(1-\alpha)}$.

\begin{lemma}[Deterministic consequences of the good event]
\label{mfd:in-deterministic}
Every such solution satisfies the following statements.
\begin{enumerate}[label=(\roman*),leftmargin=*]
\item The local H\"older estimate with source, of the form~\eqref{e.Holder.estimate.boxes.local} on~$z+3^{-k}\cu_0$ with~$z+3^{1-k}\cu_0$ in place of~$\cu_m$, holds at exponent~$\alpha$, with a fixed deterministic constant. More precisely, if the iteration of Lemma~\ref{l.iteration.lemma.GMC} is run with a fixed exponent in~$(\alpha,1)$, then the constant in its conclusion after~$D$ scales is at most~$3^{Ck_0+C(\lambda_{\rm det}+\gmcdelta^2+\ep^8)D}$, where~$C(\lambda_{\rm det}+\gmcdelta^2+\ep^8)$ is smaller than the difference between this exponent and~$\alpha$.
\item The bound~\eqref{eq:mfd-2} for H\"older boundary values holds on~$z+3^{-k}\cu_0$ with the factor~$\Lambda_{\upsilon/2,2}(z+3^{-k}\cu_0;A_N)$ of Lemma~\ref{mfd:lem-extension}, which lies between~$\ep_{\rm ell}/2$ and~$2\ep_{\rm ell}^{-1}$ times the finite-cutoff scalar~$s_N(z+3^{-k}\cu_0)$ of Section~\ref{mfd:sec-environment}; this scalar is the factor by which the coefficient on~$z+3^{-k}\cu_0$, rescaled to unit size, is a multiple of a copy of~$A^0_{N-k}$, up to the oscillation factor of Lemma~\ref{mfd:lem-cell-ellipticity}.
\item The interior harmonic comparison of Lemma~\ref{l.cutoff.regularity.good.scale.estimates} holds with tolerance~$\ep_{\rm hom}$ and the source term normalized by~$s_N(w'+3^{-\ell'}\cu_0)$ on the comparison cube~$w'+3^{-\ell'}\cu_0$.
\end{enumerate}
\end{lemma}

\begin{proof}
On~$\mathcal G_N^{\rm det}(z+3^{-k}\cu_0)$, we apply the iteration of Lemma~\ref{l.iteration.lemma.GMC} with the bad set consisting of the~$k_1$ additional indices and the indices of the scales at which one of the conditions fails. This set has at most~$2k_0+\lambda_{\rm det}D$ elements. By~\eqref{e.cutoff.regularity.good.E} the exponent in the constant of the conclusion is at most~$C(\gmcdelta^2+\ep^8+\lambda_{\rm det})D+Ck_0$. This gives~(i). Condition~(a) gives both the~$(\upsilon/2,2)$ bound for extension and, by its weighted Cauchy--Schwarz comparison, the~$(1,1)$ bound for mean-zero Poincar\'e. In particular, on a catalogue enlargement~$w+3^{-\ell}\cu_0$, the bounds~$\operatorname{tr}\bigl(s_N(w+3^{-\ell}\cu_0)^{-1}\a(w+3^{-\ell}\cu_0;A_N)\bigr)\leq d\Lambda_{\upsilon/2,2}(w+3^{-\ell}\cu_0;A_N)/s_N(w+3^{-\ell}\cu_0)$ and~$p\cdot s_N(w+3^{-\ell}\cu_0)^{-1}\a(w+3^{-\ell}\cu_0;A_N)p\geq\lambda_{\upsilon/2,2}(w+3^{-\ell}\cu_0;A_N)|p|^2/s_N(w+3^{-\ell}\cu_0)$ give trace coercivity with constant~$\ep_{\rm ell}^2/(4d)$ at a finite cutoff, and~$\ep_{\rm ell}^2/d$ in the limit. Condition~(a), the estimate~\eqref{e.cg.RHS} and the trace extension constructed in the proof of Lemma~\ref{mfd:lem-extension}, with~$0<\upsilon<\beta-\nf12$, give~(ii). To obtain~(iii), use the deterministic comparison argument in the proof of Lemma~\ref{l.cutoff.regularity.good.scale.estimates}, with comparison coefficient~$s_N(w'+3^{-\ell'}\cu_0)\Id$. In the coordinates~$y$ with~$x=w'+3^{-\ell'}y$, divide both coefficient and source by~$s_N(w'+3^{-\ell'}\cu_0)$; the comparison coefficient is then~$\Id$. Use regularity parameter~$8s\leq\nf14$, so that the error index is~$s$. The proof after the anchor error estimate uses only that error, positive scalar comparison and the ellipticity bounds, and applies to this scalar as well. Condition~(b) and the choice of~$c_{\rm det}$ give the tolerance~$\ep_{\rm hom}$; scaling back gives the factor~$3^{-2\ell'}s_N(w'+3^{-\ell'}\cu_0)^{-1}$ on the source. Condition~(c) replaces~$s_N(w'+3^{-\ell'}\cu_0)$ by~$s_N(z+3^{-k}\cu_0)$ up to a fixed factor when required. These bounds depend only on the coefficient and the deterministic geometry, so the same bounds apply to every weak solution.
\end{proof}

\begin{lemma}[Passage to the limit in sums over finitely many scales]
\label{mfd:lem-prefix-limit}
Along the subsequence that defines~$\mathscr E$, the following statements hold for each fixed finite interval of scales.
\begin{enumerate}[label=(\roman*),leftmargin=*]
\item The quantities~$\mathsf Z_{N-n,3^Nz}$ and~$\mathbf G_{N-n,3^Nz}$, for fixed~$n$ and~$z$, and each of the finitely many coarse-grained quantities in (a)--(c) converge in~$L^p$, for each of the finitely many orders~$p$ fixed at the start of Section~\ref{sec:mfd-uniqueness}. The moment bounds are uniform in the deterministic position of the interval of scales.
\item Let~$\varkappa,D>0$. If, at every cutoff, the sum of~$\mathsf Z_{N-n,3^Nz}$ or of~$\mathbf G_{N-n,3^Nz}$ over the interval of scales exceeds a given threshold with probability at most~$e^{-\varkappa D}$, then the limit of this sum exceeds every larger threshold with probability at most~$e^{-\varkappa D}$.
\end{enumerate}
\end{lemma}

\begin{proof}
We first prove~(i). For each requested~$p$, fix an auxiliary order~$P\geq2p$ with~$\nf dP$ below the margin of the geometric weights, and then decrease~$\gmcdelta$ so that~$C_P\gmcdelta^2$ is below the remaining margin, where~$C_P$ is the constant in the bound~$e^{C_P\gmcdelta^2j}$ on the~$L^P$ norm of~$\mathcal J$ obtained below. Split the terms that involve only the fields into finite sums at finitely many depths and geometrically weighted tails. By~\eqref{e.coefficient.field.regularity} the tails are bounded in~$L^p$ by their geometric weights times the moments of the fields. The other parts are geometrically weighted maxima of~$J$. For the~$\mathcal J$-terms of~$\mathsf R$ and~$\mathbf G$, the cube scale is the stationary cutoff of its coefficient~$\ahom_m^{-1}\a_m$. Translation and rescaling to~$\cu_0$ give the coefficient~$A_m^0$ in its own reindexed environment, with no remaining cell scalar or normalization ratio. On~$\cu_0$ the exact decomposition is
\begin{equation*}
	J(\cu_0,e,e;A)
	=\nf12\Lambda_{\cu_0;A}(e\cdot x)
	+\nf12\,e\cdot\a_*^{-1}(\cu_0;A)e-|e|^2\,.
\end{equation*}
Both the minimal energy with affine boundary values and the Neumann energy occur. Include these energies, for~$A_m^0$ in the translated environment, in the countable extraction of Section~\ref{mfd:sec-local-form}. For the finite coordinates in~(a)--(c), include their actual normalized matrices, scalar ratios and geometrically weighted descendant arrays in that extraction as well.

\smallskip

We truncate the geometrically weighted sums at a large relative depth. This extraction gives convergence in measure of the finitely many retained quantities; the uniform bounds at the larger fixed moment orders upgrade this to convergence in~$L^p$. At a relative depth~$j$ beyond the truncation there are at most~$C3^{dj}$ cells. By Proposition~\ref{p.coarse.grained.bound} and the bound~\eqref{e.annealed.ordering.ratio} on the ratios~$\ahom_m\ahom_n^{-1}$, the quantity~$\mathcal J$ of each cell has~$L^P$ norm at most of order~$e^{C_P\gmcdelta^2j}$, and taking a spatial maximum at the auxiliary order~$P$ introduces a factor~$3^{\frac{d}{P}j}$. The tail is then exponentially small in the truncation depth, uniformly in the cutoff.

\smallskip

Below the smallest field scale~$3^{-N}$, the maximum of the coefficient on a cell whose side equals that scale controls every smaller descendant by the direct bounds~\eqref{e.direct.J.bound.subunit}, so the finer depths require no further spatial maximum. The continuous piecewise linear functions and minima are Lipschitz, and square roots satisfy~$|\sqrt x-\sqrt y|\leq\sqrt{|x-y|}$. Hence the quantities built from these are Cauchy in~$L^p$.

\smallskip

For~(ii), on the event that the limiting sum exceeds the larger threshold, a subsequence of the approximating sums that converges almost surely eventually exceeds some threshold strictly between the two thresholds, and each approximating sum exceeds that threshold with probability at most~$e^{-\varkappa D}$. Fatou's lemma applied to the indicators gives the bound for the limit.
\end{proof}

A cell~$z+3^{-k}\cu_0$ is \textbf{good} for~$\mathscr E$ when the event~$\mathcal G_{\mathscr E}(z+3^{-k}\cu_0)$ specified above holds, and \textbf{bad} otherwise. Refine the subsequence defining~$\mathscr E$ jointly for the deterministic reference ratios and the countable family of normalized arrays. In particular,~$s_N(z+3^{-k}\cu_0)\to s_{\mathscr E}(z+3^{-k}\cu_0)>0$, with the scalar of Section~\ref{mfd:sec-compare}; for a negative-depth enlargement use the definition of~$s_N$ for~$\ell<0$. Lemma~\ref{mfd:lem-prefix-limit} gives convergence in~$L^p$ for each of the fixed orders. A further diagonal refinement gives almost-sure convergence at all the countably many depths and observation centers. Each fixed finite set of conditions defining~$\mathcal G_N^{\rm det}(z+3^{-k}\cu_0)$ then holds eventually on~$\mathcal G_{\mathscr E}(z+3^{-k}\cu_0)$, using~$\lambda_{\rm lim}<\lambda_{\rm det}$, the weaker finite-cutoff ellipticity and error bounds, and the open limiting ratio interval.

\smallskip

Apply Lemma~\ref{mfd:in-deterministic} at each fixed depth to the finite-cutoff solutions converging to~$u=G_{\mathscr E}f$ with smooth~$f$, on a cube that contains the closures of the catalogue cubes that the estimate uses. Pass their inequalities to the limit along recovery sequences, use Corollary~\ref{mfd:cor-energy-measures}(a) for their local energy measures, and intersect over the depths. This gives, on~$\mathcal G_{\mathscr E}(z+3^{-k}\cu_0)$, the local H\"older estimate of Lemma~\ref{mfd:in-deterministic}(i), the bound~\eqref{eq:mfd-2} for H\"older boundary values, and harmonic approximation with tolerance~$\ep_{\rm hom}$. The error constant was chosen above to absorb the factor two in the finite-cutoff bound. The comparison argument of Section~\ref{mfd:sec-compare} applies regularity only to the solutions~$G_{\mathscr E}f$ with smooth~$f$. For the second limit form~$\mathscr F$ introduced there, the argument uses only the maximum principle for harmonic extensions.

\subsubsection{Events depending on finitely many fields}

The conditions depend on all the fields. We approximate the quantities in the conditions by functions of the fields with indices in a finite interval around the scale of the cube, and then cover the failure of~$\mathcal G_{\mathscr E}(z+3^{-k}\cu_0)$ by events that depend on such intervals and have exponentially small probabilities.

Let~$X_{n,z}$ denote~$\mathbf G_{n,z}$, one of the three summands of~$\mathsf Z_{n,z}$, or a fixed finite nonnegative combination of the absolute values of the centered normalized coordinates specified in~(a)--(c), or their fixed continuous ramps. For the quantities in~(a)--(c),~$X_{n,z}$ refers to the cell of physical depth~$k=N-n$ and physical center~$3^{-N}z$ in the relabelled coordinates; all catalogue depths must be at most~$N$. The window about stationary index~$n$ is therefore the window about the original index~$n-N=-k$ of the fields~$\g_i$.

\begin{lemma}[Approximation by functions of finitely many fields]
\label{mfd:lem-band}
There exist~$a,c>0$, depending only on the fixed exponents and catalogue and not on the moment order, with the following property.
\begin{itemize}
\item Let~$1\leq p<\infty$. Assume that the auxiliary orders needed for~$p$ are included in the finite list fixed at the start of Section~\ref{sec:mfd-uniqueness}, and that~$\gmcdelta$ is sufficiently small, depending on that list and the fixed thresholds.
\item Then there is~$C_p<\infty$, depending on~$d$, the regularity and iteration exponents, the moment orders, the catalogue and the ramp thresholds, and uniform in the valid scale, deterministic center and cutoff, such that, for every~$h\geq1$ and every deterministic~$n$ and~$z$ for which these quantities are defined at the chosen cutoff, there is a random variable~$\widetilde X_{n,z}^{(h)}$, measurable with respect to the fields with indices in an interval of length at most~$C(h+1)$ about~$n$, such that
\begin{equation}
	\bigl\|X_{n,z}-\widetilde X_{n,z}^{(h)}\bigr\|_{L^p}\leq C_p\gmcdelta^c3^{-ah}\,.
	\label{eq:mfd-22}
\end{equation}
\end{itemize}
\end{lemma}

\begin{proof}
The proof proceeds in four steps, and the quantity we track is the error~$\|X_{n,z}-\widetilde X_{n,z}^{(h)}\|_{L^p}$. In Steps~1--3 we truncate the geometrically weighted sums defining~$X_{n,z}$ at a depth proportional to~$h$, so that the retained terms are, or are approximated by, functions of the fields in the window: Step~1 treats the terms that involve only the fields, Step~2 the product term~$\mathsf M$, and Step~3 the energy terms of~$\mathsf R$ and~$\mathbf G$ and the normalized coordinates in~(a)--(c). Step~4 obtains the factor~$\gmcdelta^c$ in~\eqref{eq:mfd-22}. The heart of the matter is Step~3, where Proposition~\ref{mfd:prop-16} approximates each of these energies by its conditional expectation given finitely many fields, and the geometric depth weights absorb the spatial maxima over the~$C3^{dr}$ cells at relative depth~$r$.

\smallskip

\emph{Step 1.} We treat the terms that involve only the fields. Before the geometric weights are applied, the depth-$j$ sum in~$\mathsf F_{n,z}$ has~$L^p$ norm at most~$C_p\gmcdelta(j+1)^{\nf32}$: we cover~$z+\cu_{n+1+j}$ by cubes of the scale of each field, use the sub-Gaussian bound for the maximum from~\eqref{e.coefficient.field.regularity}, and sum the~$2j+1$ terms. Fix a constant~$c_4>0$, depending only on the fixed exponents, and discard the depths~$j>c_4h$, whose total contribution is at most~$C_p\gmcdelta3^{-\frac{sc_4}{16}h}$. The tail terms of~$\mathbf G$ (see~\eqref{e.def.bfG.mhq}) carry the explicit weight~$3^{k-i}$.

\smallskip

\emph{Step 2.} We treat the product term. We retain the depths~$j\leq\lfloor h/2\rfloor$ and stop the second product at~$n+h$. Both~$\mathsf M$ and its approximant are at least~$2$, the approximant is not larger, and the removed depths contribute only their excess over~$2$. Fix~$j$. The first product factor at depth~$j$ is the exponential of~$\sup_{x\in z+\cu_{n+1+j}}\sum_{i=(n-j)\vee0}^{n+j}|\g_i(x)|$. We cover the cube by at most~$C_d3^{d(2j+1)}$ cubes of side~$3^{(n-j)\vee0}$. On each such cube, this sum is at most the sum over~$i$ of the value of~$|\g_i|$ at the center plus~$C_d3^{(n-j)\vee0}$ times the gradient norm of~$\g_i$ on the cube; these terms are independent in~$i$ and sub-Gaussian with parameter~$C_d\gmcdelta$. We take the exponential moment before the spatial maximum, which bounds the moment of order~$\xi$ on each cube by~$\exp(C_d\xi\gmcdelta(j+1)+C_d\xi^2\gmcdelta^2(j+1))$, and a maximum over the covering cubes, taken in~$L^q$ with~$q\geq p$, introduces a factor equal to their number to the power~$\nf1q$. With~$e^x-1\leq xe^x$, choosing~$q$ large depending on~$d,s,p$ and then~$\gmcdelta$ small shows that~$3^{-\frac{s}{8}j}$ times the~$L^p$ norm of the excess of the first product factor over one is at most~$C_p\gmcdelta(j+1)3^{-\frac{s}{16}j}$.

\smallskip

The second product factor is the exponential of~$4\sup_x\sum_{i\geq n+j}|\g_i(x)-\g_i(z)|$; this supremum is at most a sum, over~$i\geq n+j$, of sub-Gaussian variables multiplied by~$C_d3^{n+j-i}$. Hence its excess over one has~$L^p$ norm at most~$C_p\gmcdelta$, and the part with~$i>n+h$ contributes~$C_p\gmcdelta3^{j-h}$.

\smallskip

Summing over~$j\leq h/2$ gives~$C_p\gmcdelta3^{-\frac12h}$, and the discarded depths~$j>h/2$ contribute~$\sum_{j>h/2}3^{-\frac{s}{8}j}$ times the~$L^p$ norm of the sum of the two excesses at depth~$j$, which is at most~$C_p\gmcdelta3^{-ah}$.

\smallskip

\emph{Step 3.} We treat the energy terms, fixing their moment orders first. Choose an auxiliary~$P\geq\max\{2,2p\}$ so large that~$d/P\leq\min\{s,\upsilon/2\}/32$. Increase it, if necessary, to dominate the fixed spatial and comparison exponents. Fix moment orders~$3P$ and~$6P$ for the Morrey and smoothed Neumann bounds, the orders~$4P$ and~$8P$ for the Neumann energies in Proposition~\ref{mfd:prop-16}, and the finitely many higher product and mesh orders required by those estimates, before reducing~$\gmcdelta$. Reduce~$\gmcdelta$ until twelve times each selected order is at most~$C^{-1}\gmcdelta^{-2}|\log\gmcdelta|^{-1}$ and the geometric moment conditions below hold. Lemma~\ref{mfd:lem-15}, both clauses of Proposition~\ref{mfd:prop-16}, and their moment inputs hold with~$H$ replaced by~$0$ or by its partial sums~$H_{L'}$, uniformly in~$L'$. Their proofs use only the increments~$\g_i-\g_i(0)$ on a fixed compact set. Lemma~\ref{mfd:lem-infrared} bounds these increments uniformly for~$0$, every~$H_{L'}$ and~$H$. The proof of Proposition~\ref{mfd:cor-neumann-source} gives its Neumann bounds for these coefficients with one random constant, uniformly in~$L'$, and gives the translated estimates for~$A_N^0$. Thus the subsequent perturbation estimates have the same moment inputs.

\smallskip

For~$\mathsf R$ and the energy part of~$\mathbf G$, every~$\mathcal J$ is computed for~$A_n^0$ at its own stationary scale~$n$, as in the preceding proof. Apply both clauses of Proposition~\ref{mfd:prop-16} to the affine Dirichlet and constant-flux Neumann terms, using the finite polarization family of the vectors~$e_i$ and~$e_i+e_j$, normalizing these vectors and restoring their quadratic factors by homogeneity. The slower rate gives a nonnegative approximation~$\widehat{\mathcal J}$ with
\begin{equation*}
	\|\mathcal J-\widehat{\mathcal J}\|_{L^P}
	\leq C_P\gmcdelta^{\nf12}3^{-b\ell/256}\,.
\end{equation*}
The approximant can be obtained by conditioning~$\mathcal J$ itself; nonnegativity is preserved. When~$\ell$ exceeds the cutoff~$n$ of~$\mathcal J$, the quantity~$\mathcal J$ already depends only on the retained fields and the error is zero. No scalar~$s_N$ enters this estimate.

\smallskip

Retain relative depths~$r\leq\lfloor h/4\rfloor$ and use~$\ell=\lfloor h/2\rfloor$ for each retained energy. At depth~$r$ there are at most~$C3^{dr}$ cells, so the maximum over them introduces a factor~$C3^{dr/P}$. The retained errors and the discarded tail of~$\mathsf R$ are bounded by
\begin{equation*}
	C_P\gmcdelta^{\nf12}3^{-bh/512}
	\sum_{r\geq2}3^{-(s/8-d/P)r}
	+C_P\gmcdelta^2
	\sum_{r>h/4}3^{-(s/8-d/P)r}\,.
\end{equation*}
For~$\mathbf G$ use~$|\sqrt{\min\{x,1\}}-\sqrt{\min\{y,1\}}| \leq\sqrt{|x-y|}$ and its stronger depth weight. The square root halves the error exponents. Hence, for these two energy families,
\begin{equation*}
	\|X-\widetilde X\|_{L^p}
	\leq C_P\gmcdelta^{\nf14}3^{-\min\{b/1024,s/64\}h}\,.
\end{equation*}
Small~$h$ is absorbed into the constant. The retained original field indices lie in an interval of length~$C(h+1)$ about the parent scale. In the original indexing of the fields~$\g_i$, a parent of depth~$k$ has this interval about~$-k$. Fixed catalogue offsets only enlarge it by a fixed amount.

\smallskip

For the normalized coordinates in~(a)--(c), first write the coarse ellipticity sums, both coarse-grained matrices and the error using the coefficient of the rescaled cell. Split the sums at relative depth~$\lfloor h/4\rfloor$. On a retained descendant, separate the reference ratio between that descendant and the catalogue cube. This ratio is part of the normalized coefficient decomposition; its positive and inverse moments satisfy~$\|\rho+\rho^{-1}\|_{L^\xi}\leq C_\xi e^{C_\xi\gmcdelta^2r}$. Retain the intervening coarse fields and truncate the remaining centered infrared increments. Apply the moment bounds for the energies on the unit cube at~$2P$ and H\"older at~$2P$ to these products, using the already selected~$6P$ Morrey inputs. After the spatial maximum, the bound is multiplied by~$3^{-\vartheta r}3^{dr/P}e^{C_P\gmcdelta^2r}$, where~$\vartheta>0$ is the fixed depth weight of the quantity being treated. Choose~$P$ with~$d/P\leq \vartheta/4$ and then~$\gmcdelta$ with~$C_P\gmcdelta^2/\log3\leq \vartheta/4$. This leaves~$3^{-\vartheta r/2}$. The omitted infrared increments contribute at most~$C_P\gmcdelta3^{-c_2h}$ for a fixed~$c_2>0$. Below the smallest field scale, the direct bound~\eqref{e.direct.J.bound.subunit} on each terminal cell controls all its finer descendants, so there is one spatial maximum at that terminal depth.

\smallskip

The upper ellipticity sum is approximated first, then the inverse lower ellipticity sum. For the lower coefficient, truncate the reciprocal at~$\frac12\,3^{-ah/2}$ and use the moments at order~$4P$ to control the small-denominator part. This loses at most half the decay exponent~$a$; it does not change its independence from~$p$. Polarization supplies the two normalized matrices. Apply the same descendant argument to the error with comparison coefficient~$s_N(w'+3^{-\ell'}\cu_0)\Id$. The finite scalar ratios use the explicit formula defining~$s_N$. These are exactly the three ellipticity coordinates, the two matrices, the error and the ratio specified above. The finite family of exponents and catalogue cubes permits one common field interval, the minimum of their positive decay rates and the minimum of their disorder thresholds. A fixed Lipschitz function of the approximated coordinates then gives the asserted band approximation for every finite test. In particular the~$(1,1)$ consequences require no extra approximation: their bounds follow deterministically from the~$(\upsilon/2,2)$ tests on each enlargement.

\smallskip

\emph{Step 4.} We obtain the factor~$\gmcdelta^c$ in~\eqref{eq:mfd-22}. The deviations of the quantities defining the conditions from their values at~$\gmcdelta=0$ satisfy~$\|\mathsf F\|_p+\|\mathsf M-2\|_p\leq C_p\gmcdelta$, while~$\|\mathsf R\|_p\leq C_p\gmcdelta^2$ and~$\|\mathbf G\|_p\leq C_p\gmcdelta^{\nf12}$. The first bound follows from Steps~1 and~2. For~$\mathsf R$, group the cubes of~\eqref{e.good.sum.J} over which the supremum is taken by~$r=m-n$. There are at most~$C_d3^{dr}$ such cubes, the quantity~$\mathcal J$ of each has~$L^q$ norm~$C_q\gmcdelta^2$ by Proposition~\ref{p.coarse.grained.bound}, and~$\sum_{r\geq2}3^{-(\frac s8-\frac dq)r}<\infty$ for a fixed~$q\geq p$ large enough. For~$\mathbf G$, we take~$m=n$ in Lemma~\ref{l.sum.the.errors} and obtain~$C_p\gmcdelta(1+|\log\gmcdelta|^{\nf12})\leq C_p\gmcdelta^{\nf12}$. If a truncation gives only~$C_p3^{-ah}$, without a positive power of~$\gmcdelta$, we combine it with a bound~$C_p\gmcdelta^{c_3}$, for some~$c_3>0$, on the distance of~$X$ from its deterministic value at~$\gmcdelta=0$. Conditional expectation onto the fields with indices in the interval of length~$C(h+1)$ about~$n$ is a contraction in~$L^p$. Thus the distance of~$X$ from its conditional expectation is at most twice its distance from the approximation. The same bound holds with the value at~$\gmcdelta=0$ in place of the approximation. The minimum of the two resulting bounds is at most their geometric mean~$2C_p\gmcdelta^{\nf {c_3}2}3^{-\frac{a}{2}h}$. The centered moment bound for the normalized coordinates is~$\max_{Y}\|Y-Y^{(0)}\|_{L^p}\leq C_p\gmcdelta^{c_1}$ with~$c_1>0$, where~$Y$ runs over the fixed finite array,~$Y^{(0)}$ is its unit-coefficient value and~$c_1$ is fixed before~$p$. To obtain it, apply~\eqref{e.bound.Lambdas.by.Es} with comparison~$\Id$ after translation, scaling and removal of the scalar of the catalogue cube. The small-error moments of~\eqref{e.homogenize.at.scale.m.E}, together with the centered infrared and coarse-field oscillation bounds, control the right-hand side uniformly in the cutoff. The order~$\lambda\leq\Lambda$ turns the upper bounds for~$\Lambda$ and~$\lambda^{-1}$ into two-sided bounds, so their deviations from~$1$ are bounded by~$C(\mathcal E_{\upsilon/2,\infty,2}+\mathcal E_{\upsilon/2,\infty,2}^2)$. The spectral estimate in the proof of~\eqref{e.bound.Lambdas.by.Es} gives the corresponding deviations of both matrices from~$\Id$. The error uses its own small-error bound; the formula defining~$s_N$ and the ordering of~$\ahom$ give the ratio deviation. All products use the larger orders already selected. Below the smallest field scale, use the direct terminal-cell bound as above. The fixed Lipschitz ramps inherit these bounds. Take the minimum of the finitely many positive decay and disorder-power exponents. These exponents are fixed by the geometric weights and the regularity parameters before~$p$ is chosen; only the auxiliary moment orders, disorder threshold and prefactor depend on~$p$.
\end{proof}

\begin{lemma}[Sums of~$\mathsf Z$ and~$\mathbf G$ over scales at a fixed center]
\label{mfd:in-prefix}
Let~$\lambda>0$ and~$\varkappa<\infty$. If~$\gmcdelta$ is sufficiently small, then there is~$C<\infty$ such that, for every deterministic center~$z$, every interval~$n,\ldots,n+k-1$ of scales, in the coordinates of~\eqref{eq:mfd-normalization}, every cutoff~$N$, with all quantities evaluated at~$N$, and each~$X\in\{\mathsf Z,\mathbf G\}$,
\begin{equation*}
	\P\Bigl[\sum_{i=0}^{k-1}X_{n+i,z}>\frac{\lambda k}{4}\Bigr]\leq Ce^{-\varkappa k}\,.
\end{equation*}
\end{lemma}

\begin{proof}
For~$X=\mathsf Z$ we apply Lemmas~\ref{l.bad.scales.for.nabla.gj},~\ref{l.exp.goodscales} and~\ref{l.discounted.J.local} at the lower thresholds of the continuous piecewise linear functions. The exponential rates in these lemmas tend to infinity as~$\gmcdelta\downarrow0$. For~$X=\mathbf G$ we apply Lemma~\ref{l.sum.the.errors}. Its estimate~\eqref{e.after.good.scales.uniform} concerns the average of the approximation errors over the interval, with an~$\mathcal O_{\Gamma_2}$ fluctuation of size~$(m-n+1)^{-\nf12}$. After multiplication by the length~$m-n+1$, the sum exceeds~$\lambda(m-n+1)$ with probability at most~$\exp(-c\lambda^2(m-n+1)\gmcdelta^{-2}|\log\gmcdelta|^{-1})$. The bounds are uniform in the position of the interval by the stationarity in Assumption~\ref{a.g1} and by Proposition~\ref{p.cutoff.Holder.regularity}.
\end{proof}

For use with the limiting good event, fix a physical depth~$j$, an observation center, an integer~$h\geq0$ and a window of the fields~$\g_i$ with indices in an interval of length at most~$C(h+1)$ about~$-j$. Write~$X_N$ for a primitive summand or finite ramp above, pulled back from the stationary coordinates at cutoff~$N$, and~$X_{\mathscr E}$ for its limit. The window is the same for all sufficiently large~$N$. Lemma~\ref{mfd:lem-prefix-limit}(i) and contraction give
\begin{align*}
	\bigl\|\E[X_N\mid\text{fields in the window}]-\E[X_{\mathscr E}\mid\text{fields in the window}]\bigr\|_{L^p} \leq\|X_N-X_{\mathscr E}\|_{L^p}\longrightarrow0\,.
\end{align*}
By Lemma~\ref{mfd:lem-band}, the distance of~$X_N$ from its conditional expectation is at most~$2C_p\gmcdelta^c3^{-ah}$. Passing to the limit gives the same bound for~$X_{\mathscr E}$; the centered moment bounds pass to the limit as well. Thus the limiting quantities have band approximants~$\E[X_{\mathscr E}\mid\text{fields in the window}]$, with the same decay exponents and at most twice the constant. For the primitive prefixes, apply Lemma~\ref{mfd:in-prefix} with parameter~$\lambda_{\rm cut}$ to the interval of depths~$\ell-k_1,\ldots,\ell+D$ of the prefix. Its length is at most~$3D$ since~$D\geq k_0\geq\max\{k_1,1\}$, so the resulting bound applies at threshold~$\lambda_{\rm cut}D<\lambda_{\rm lim}D/4$. Lemma~\ref{mfd:lem-prefix-limit}(ii) then gives the exponential bound at~$\lambda_{\rm lim}D/4$ for the limiting sum. The rate can be chosen arbitrarily large before the moment orders and disorder threshold are fixed.

\begin{lemma}[Covering the limiting bad event]
\label{mfd:lem-witness}
Use the limiting arrays defining~$\mathcal G_{\mathscr E}(z+3^{-n}\cu_0)$. Let~$\varkappa<\infty$. If~$\gmcdelta$ is sufficiently small, then, for every cell~$z+3^{-n}\cu_0$ at triadic depth~$n$, there exist events~$W_{z+3^{-n}\cu_0,h}$ for~$h\geq1$ such that
\begin{itemize}
\item the event~$W_{z+3^{-n}\cu_0,h}$ is measurable with respect to the fields~$\g_{-j}$ whose indices~$j$ lie between~$n-h$ and~$n+2h$,
\item its probability satisfies~$\P(W_{z+3^{-n}\cu_0,h})\leq e^{-\varkappa h}$,
\item up to a null set,
\begin{equation*}
	\{\mathcal G_{\mathscr E}(z+3^{-n}\cu_0)\text{ fails}\}\subseteq\bigcup_{h\geq1}W_{z+3^{-n}\cu_0,h}\,.
\end{equation*}
\end{itemize}
The same conclusion holds for the failure of one additional scalar or matrix condition on a quantity whose centered~$L^p$ norm is at most~$C_p\gmcdelta$ and whose approximation by a function of the fields with indices in an interval of length at most~$C(h+1)$ about~$n$ has error at most~$C_p\gmcdelta3^{-ah}$ in~$L^p$, for every~$h\geq1$, at one order~$p$, large enough depending on~$a$,~$\varkappa$ and the catalogue.
\end{lemma}

\begin{proof}
Outside a common null set, all dyadic conditional-expectation approximants converge to their limiting quantities. This follows from the summable~$L^p$ errors and a countable intersection over the cells, centers and finite tests. If~$\mathcal G_{\mathscr E}(z+3^{-n}\cu_0)$ fails, then either one of the finite tests~(a)--(c) fails, or there are a prefix length~$D\geq k_0$, one of at most~$Ce^{CD}$ choices of a cube and a center, and a limiting primitive prefix sum~$S$ with~$S\geq\lambda_{\rm lim}D$. The~$k_1$ coarser depths only change the constant in the bounds below. Approximate~$S$ by~$S^{(\ell)}$, the sum of the conditional expectations of its summands on windows of radius~$2^\ell D$. The limiting band bound above gives~$\|S-S^{(\ell)}\|_{L^p}\leq C_pD\gmcdelta^c3^{-a2^\ell D}$. If~$S\geq\lambda_{\rm lim}D$, then either~$S^{(0)}>\lambda_{\rm lim}D/2$ or, for some~$\ell\geq1$,~$|S^{(\ell)}-S^{(\ell-1)}|>\lambda_{\rm lim}D2^{-\ell-2}$. Indeed, the increment allowances sum to~$\lambda_{\rm lim}D/4$, and the approximants converge almost surely to~$S$.

\smallskip

The first event implies either~$S>\lambda_{\rm lim}D/4$ or~$|S-S^{(0)}|>\lambda_{\rm lim}D/4$. Use the limiting prefix bound proved before the lemma, with a rate larger than~$3\varkappa+C$, and Markov's inequality. Its probability is at most~$Ce^{-(3\varkappa+C)D}+(4C_p\gmcdelta^c/\lambda_{\rm lim})^p3^{-apD}$. The~$\ell$th increment event has probability at most~$(C_p\gmcdelta^c)^p2^{p(\ell+C)}3^{-ap2^{\ell-1}D}$.

\smallskip

Take the union over the~$Ce^{CD}$ spatial choices. For each~$h$, combine the events for pairs~$(D,\ell)$ with~$2^\ell D\asymp h$; there are polynomially many such pairs. Their field windows lie within an interval of radius~$Ch$ about~$n$. After relabeling~$h$, enlarge this interval to~$[n-h,n+2h]$. Choose~$p$ so that the decay rate~$ap\log3$ absorbs the prescribed rate, the spatial count and this fixed enlargement, and then decrease~$\gmcdelta$ to absorb the prefactors.

\smallskip

For each finite test~(a)--(c), let~$T_{\mathscr E}\in[0,1]$ be its specified ramp. Failure of that test implies~$T_{\mathscr E}=1$, including either endpoint of the open ratio interval. The limiting centered bound gives~$\|T_{\mathscr E}\|_{L^p}\leq C_p\gmcdelta^{c_1}$, and its conditional-expectation error is at most~$C_p\gmcdelta^c3^{-ah}$. Start at a fixed window large enough to include the catalogue offsets. The probability that this first approximant exceeds~$1/2$ is small by its~$L^p$ norm; decrease~$\gmcdelta$ so that it has the required exponential bound at that fixed radius. Telescope the approximants with increment allowances~$2^{-\ell-2}$ and apply Markov's inequality as for the prefixes. The finite union of these covers gives the required events for~(a)--(c), after starting with a slightly larger rate. The same argument applies to the additional centered condition in the statement, after dividing by its fixed positive failure threshold.
\end{proof}

Fix a subdivision factor~$3^{H_1}$ and a bounded cube~$w+R\cu_0$. A descendant at~$3^{H_1}$-adic generation~$n$ has side length~$3^{-H_1n}$ times the side length of~$w+R\cu_0$. When~$w+R\cu_0$ has triadic side length~$R=3^{-k_*}$, its descendants at generation~$n$ have physical triadic depth~$k_*+H_1n$. A \textbf{chain} of length~$m$ is a sequence of~$m$ cubes of the~$3^{H_1}$-adic subdivision of~$w+R\cu_0$, each a child of the preceding one and the first a child of~$w+R\cu_0$.

\begin{proposition}[Bad cells along every chain]
\label{mfd:prop-allchain}
Fix a subdivision factor~$3^{H_1}$ and a bounded cube~$w+R\cu_0$, and let~$\theta\in(0,1)$ and~$b_*>0$. There exists~$\gmcdelta_0>0$, depending on~$d$,~$H_1$,~$\theta$,~$b_*$ and the parameters of the good event, such that, if~$\gmcdelta\leq\gmcdelta_0$, then there is a random variable~$B_\omega<\infty$ such that
\begin{itemize}
\item for every~$\lambda>0$,
\begin{equation*}
	\P(B_\omega>\lambda)\leq(1-e^{-b_*})^{-1}\exp\Bigl(-\frac{b_*\lambda}{1-\theta}\Bigr)\,,
\end{equation*}
\item almost surely, for every~$m\geq1$ and every chain~$\pi$ of length~$m$,
\begin{equation}
	\#\{\text{bad cells on }\pi\}\leq\theta m+B_\omega\,.
	\label{eq:mfd-23}
\end{equation}
\end{itemize}
\end{proposition}

\begin{proof}
Apply Lemma~\ref{mfd:lem-witness} with rate~$\varkappa$. On a fixed chain, successive generations have field-index centers spaced by~$H_1$; for a triadic root the center is~$-(k_*+H_1n)$. The rare-interval estimate holds on every integer interval, so a fixed root offset does not change it. Enlarge the interval of field indices determining~$W_{z+3^{-k}\cu_0,h}$ to the symmetric interval of radius~$2h$ about its depth, replace all events at other depths or odd radii by the empty event, and apply Lemma~\ref{l.concentration.rare.intervals}. In that lemma, the rate is~$\varkappa/3$ and the density is at least~$\theta/H_1$. Thus, once~$\varkappa\theta\geq CH_1$, the probability that this chain has at least~$\theta m$ bad cells is at most~$e^{-c\varkappa\theta m}$. There are~$3^{H_1dm}$ chains. Choose~$\varkappa$ with~$c\varkappa\theta\geq d\log3^{H_1}+b_*$. The probability that some chain of length~$m$ has at least~$\theta m$ bad cells is then at most~$e^{-b_*m}$. These choices precede the moment order and the threshold for~$\gmcdelta$ in Lemma~\ref{mfd:lem-witness}.

\smallskip

Let~$B_\omega\coloneqq\sup_m\max_{|\pi|=m}(\#\{\text{bad cells on }\pi\}-\theta m)_+$. If~$B_\omega>\lambda$, then for some~$m>\lambda/(1-\theta)$ there is a chain of length~$m$ with more than~$\theta m$ bad cells. Summing~$e^{-b_*m}$ over these~$m$ proves the stated tail bound and~\eqref{eq:mfd-23}.
\end{proof}

Every point of~$w+R\cu_0$ lies in exactly one chain of half-open cubes. Integrating~\eqref{eq:mfd-23} over these chains gives, on the same event, for every finite Borel measure~$\nu$ on~$w+R\cu_0$,
\begin{equation}
	\sum_{j=1}^m\nu(\text{union of the bad cells at~$3^{H_1}$-adic generation }j)
	\leq(\theta m+B_\omega)\nu(w+R\cu_0)\,.
	\label{eq:mfd-24}
\end{equation}
In particular this holds for the random energy measure~$\Gamma_{\mathscr E}(u)$.

\subsection{Comparison of two subsequential limits}
\label{mfd:sec-compare}

In this subsection we show that two subsequential limit forms differ only by a deterministic scalar factor (Theorem~\ref{mfd:thm-prop}). Let~$\mathscr E$ and~$\mathscr F$ be two such forms, obtained along two subsequences of cutoffs, realized in the same Skorokhod representation (Section~\ref{mfd:sec-local-form}), so that both are functions of the same fields. Their domains are not assumed to be equal. We compare them on the cubes~$z+\cu_n$ with~$z\in\mathbb Q^d$ and~$n\in\mathbb Z$.

\smallskip

The preliminary two-sided comparison (Theorem~\ref{mfd:thm-C0}) adapts the argument of~\citet[Section~3.4]{DGComparison} on the ratios of the normalizations of the two forms. In that argument, comparisons built from paths around annuli, together with upper and lower bounds on crossing distances, exclude ratios of the two normalizations that are too large or too small. Here harmonic extensions play the role of the paths. Proposition~\ref{mfd:prop-density} glues harmonic extensions on cubes at the scales where the ratio of the two normalizations is at most a given threshold, and shows that if this happens on a positive density of scales, then~$\mathscr F$ is bounded by~$\mathscr E$ times a constant multiple of the threshold. For a small threshold, the variational formula~\eqref{eq:mfd-18} shows that this contradicts the moment bounds of Lemma~\ref{mfd:lem-response-reciprocal}. This excludes ratios that are too small. Exchanging~$\mathscr E$ and~$\mathscr F$ excludes ratios that are too large. The proof of Theorem~\ref{mfd:thm-C0} then applies the same construction on the remaining scales, where the ratio is bounded above and below, to compare the domains and the energies in both directions.

\smallskip

As in the metric uniqueness proofs of~\citet[Section~6]{GM} and~\citet[Section~4.5]{DGUniqueness}, we consider the optimal upper and lower comparison constants. Those proofs shortened geodesic segments to obtain probability bounds incompatible with a strict difference between the constants. They had to control which annuli a field-dependent geodesic visits. The probability comparison and counting argument in~\citet[Section~4.2]{DGUniqueness} provided this control after suitable smooth changes of the Gaussian free field.

\smallskip

Here Proposition~\ref{mfd:prop-allchain} bounds, along every chain of nested cubes and on a single event that does not depend on the solution, the number of bad cells of the chain. Its integrated form~\eqref{eq:mfd-24} therefore applies to the energy measure of the solution being replaced. We use harmonic extensions on disjoint cubes in place of shorter path segments. Since the solution is close to affine on good cells, the comparison of the minimal energies of~$\mathscr E$ and~$\mathscr F$ with affine boundary values transfers to its boundary values. Summing, over these cubes, the amount by which the minimal~$\mathscr F$-energy falls below the upper constant times the minimal~$\mathscr E$-energy improves one of the comparison constants for the forms themselves.

\smallskip

In this subsection a cell~$z+3^{-k}\cu_0$ is called good exactly when both~$\mathcal G_{\mathscr E}(z+3^{-k}\cu_0)$ and~$\mathcal G_{\mathscr F}(z+3^{-k}\cu_0)$ hold. Applying Proposition~\ref{mfd:prop-allchain} to each form with half the allowed density of bad cells and adding the two counts gives~\eqref{eq:mfd-23} for this notion of bad cell. On a good cell we thus have the regularity and coercivity estimates for~$\mathscr E$ and the bound~\eqref{eq:mfd-2} for H\"older boundary values for~$\mathscr F$, and the same with the roles of~$\mathscr E$ and~$\mathscr F$ exchanged.

\smallskip

The comparison is expressed through the scalars~$s_{\mathscr E}(z+3^{-k}\cu_0)$ and~$s_{\mathscr F}(z+3^{-k}\cu_0)$ of the two forms. At the cell~$z+3^{-k}\cu_0$ of side~$3^{-k}$ and center~$z$, the scalar~$s_{\mathscr E}(z+3^{-k}\cu_0)$ of~$\mathscr E$ is the limit of the cell scalar~$s_N(z+3^{-k}\cu_0)$ of Section~\ref{mfd:sec-environment} as~$N\to\infty$ along the subsequence that defines~$\mathscr E$. We define~$s_{\mathscr F}(z+3^{-k}\cu_0)$ in the same way, along the subsequence that defines~$\mathscr F$. In~$s_N(z+3^{-k}\cu_0)$, only the ratio~$\ahom_{N-k}/\ahom_N$ depends on~$N$, and the other factors are the same for the two forms. After a further diagonal extraction, these limits exist for all~$k$ at once: since~$|\ahom_{N}\ahom_{N-k}^{-1}-1|\leq C\gmcdelta^2ke^{C\gmcdelta^2k}$ by~\eqref{e.annealed.ordering.ratio}, this ratio lies in a compact subinterval of~$(0,\infty)$ that depends only on~$k$. The ratio~$\rho_k\coloneqq\frac{s_{\mathscr F}(z+3^{-k}\cu_0)}{s_{\mathscr E}(z+3^{-k}\cu_0)}$ of these two scalars is deterministic and does not depend on~$z$.

\smallskip

Choose the parameters in the following order. First fix the geometric and regularity parameters of this subsection, then~$H_1$, the moment orders, and a first threshold~$\gmcdelta_*$ for~$\gmcdelta$, which determines the comparison constant~$C_0$ of Theorem~\ref{mfd:thm-C0}. With~$C_0$ fixed, we tighten the tolerance~$\ep_{\rm hom}$ of the harmonic comparison and add to the good event the comparison of the minimal energies of~$\mathscr E$ and~$\mathscr F$ with affine boundary values, of Corollary~\ref{mfd:cor-32}, whose error is proportional to the difference between the optimal comparison constants. This requires only a further decrease of~$\gmcdelta$. The thresholds for sums over intervals of scales and~$H_1$ stay fixed. None of these choices depends on the difference of the optimal constants.

\subsubsection{Comparison from a positive density of small ratios}

We first show that a positive density of scales at which~$\rho_k\leq a$ forces~$\mathscr F\leq C_*a\mathscr E$. Fix~$\kappa_0$ with~$8(2(1-\alpha)+\kappa)<\kappa_0<\nf13$, which is possible by the choice of~$\alpha$ and~$\kappa$ in Section~\ref{mfd:sec-local-form}, and choose an integer~$D\geq\nf{8d}{\kappa_0}$, used in the stopping rule of the proof of Proposition~\ref{mfd:prop-density}, and~$\theta\leq\nf{\kappa_0}{8}$. By the choice of~$\kappa_0$ we have~$2(1-\alpha)+\kappa<\kappa_0/8$, which is equivalent to~$d-2+2\alpha-\kappa>d-\nf{\kappa_0}{8}$. Then choose~$H_1$ so large that~$3^{H_1(\frac{\kappa_0}{3}-\frac{\kappa_0}{8})}\geq2(6d)^{\nf{\kappa_0}{3}}$. Here we use the~$3^{H_1}$-adic grid, translated as in the proof of Proposition~\ref{mfd:prop-density}, with cell sides~$3^{-H_1n}$, so~$n$ is measured from side length one and the physical triadic depth is~$H_1n$. Thus~$\rho_{H_1n}$ is the ratio at that physical depth. In a~$3^{H_1}$-adic subdivision starting from a root of side~$R=3^{-k_*}$, generation~$n$ instead has physical depth~$k_*+H_1n$. For a root of side~$3^{-H_1n_0}$ this shifts every index by the fixed integer~$n_0$, which does not affect upper densities. For a child~$z+r\cu_0$ of side~$r$ of a cell~$z'+3^{H_1}r\cu_0$ of side~$3^{H_1}r$, call it \textbf{padded} if its distance from~$\partial(z'+3^{H_1}r\cu_0)$ is at least~$3r$. Then the concentric cube of side~$6r$, which contains the concentric enlargement of side~$3r$ and the catalogue cubes of~$z+r\cu_0$ of side at most~$3r$, has its closure in~$z'+3^{H_1}r\cu_0$. At most~$6d\,3^{H_1(d-1)}$ children of each cell are not padded. For every real number~$c$, a set~$S\subseteq\mathbb N$ has \textbf{upper density} at least~$c$ if~$\limsup_mm^{-1}\#(S\cap\{1,\ldots,m\})\geq c$.

\smallskip

On a good padded cell, the regularity of~$\mathscr E$ and the bound~\eqref{eq:mfd-2} for H\"older boundary values of~$\mathscr F$ combine as follows.

Let~$z_0+\cu_{j_0}$ be a cube with~$z_0\in\mathbb Q^d$ and~$j_0\in\mathbb Z$ and~$f$ a smooth function, and let~$u\coloneqq G_{\mathscr E}^{z_0+\cu_{j_0}}f$, where~$G_{\mathscr E}^{z_0+\cu_{j_0}}$ is the Dirichlet Green operator of~$\mathscr E$ on~$z_0+\cu_{j_0}$.

In the next lemma,~$C$ depends only on~$d$, the fixed regularity and iteration exponents,~$k_0$, the primitive and trace constants, the displayed good-event thresholds and ratio bounds, and the finite enlargement and grid geometry. It is independent of the cell, depth, source, extracted subsequence and environment. It uses only the primitive iteration thresholds, the extension exponent, the ellipticity bounds and the scalar comparisons on~$z+r\cu_0$ and~$z+3r\cu_0$, and is independent of~$\ep_{\rm hom}$ and~$c_{\rm det}$. The random constant~$K$ may also depend on~$z_0+\cu_{j_0}$,~$f$ and the environment.

\begin{lemma}[Extension on good cells]
\label{mfd:lem-goodext}
\begin{enumerate}[label=(\roman*),leftmargin=*]
\item If the cell~$z+r\cu_0$ of side~$r=3^{-k}$ is a good child of~$z'+3^{H_1}r\cu_0\subseteq z_0+\cu_{j_0}$ whose concentric enlargement of side~$3r$ has its closure in~$z'+3^{H_1}r\cu_0$, and the closures of the catalogue cubes of~$z+r\cu_0$ of side at most~$3r$ lie in~$z_0+\cu_{j_0}$, then
\begin{multline*}
	\inf_{c\in\R}\|u(z+r\,\cdot)-c\|_{C^\alpha(\cu_0)}
	\\
	\leq Cr^{\nf{(2-d)}{2}}s_{\mathscr E}(z+r\cu_0)^{-\nf12}\sqrt{\Gamma_{\mathscr E}(u)(z'+3^{H_1}r\cu_0)}
	+Cr^2s_{\mathscr E}(z+r\cu_0)^{-1}\|f\|_\infty\,,
\end{multline*}
and
\begin{equation*}
	\Lambda_{\mathscr F,z+r\cu_0}(u|_{\partial(z+r\cu_0)})
	\leq C\rho_k\bigl[\Gamma_{\mathscr E}(u)(z'+3^{H_1}r\cu_0)+s_{\mathscr E}(z+r\cu_0)^{-1}r^{d+2}\|f\|_\infty^2\bigr]\,.
\end{equation*}
\item There is a random constant~$K<\infty$ such that~$\Lambda_{\mathscr F,z+r\cu_0}(u|_{\partial(z+r\cu_0)})\leq Kr^{d-2+2\alpha-\kappa}$ for every cell~$z+r\cu_0\subseteq z_0+\cu_{j_0}$ of side~$r$.
\end{enumerate}
\end{lemma}

\begin{proof}
For~(i), the local H\"older estimate of Lemma~\ref{mfd:in-deterministic}(i) on the good cell, in the cube~$z_0+\cu_{j_0}$ that contains its catalogue cubes of side at most~$3r$, bounds the oscillation on~$z+r\cu_0$ by the normalized~$L^2$ oscillation on the concentric cube of side~$3r$ and the source (Section~\ref{mfd:sec-branches}). The coarse-grained Poincar\'e inequality~\eqref{e.CG.Poincare.mean.zero} on that cube bounds its oscillation by its energy. Since its closure lies in~$z'+3^{H_1}r\cu_0$, this energy is at most~$\Gamma_{\mathscr E}(u)(z'+3^{H_1}r\cu_0)$, giving the factor~$r^{\nf{(2-d)}{2}}$ in unnormalized units. The source is controlled by~\eqref{eq:mfd-5}: by Lemma~\ref{mfd:lem-primitive}, convolving~$f$ with the gradient of the Newtonian kernel gives a vector field with~$C^{\nf12}$ seminorm at most~$Cr^{\nf12}\|f\|_\infty$ on a cube of side~$r$. Scaling gives the term~$r^2s_{\mathscr E}(z+r\cu_0)^{-1}\|f\|_\infty$. For the form~$\mathscr F$ on~$z+r\cu_0$, the coefficient factor in~\eqref{eq:mfd-2} is comparable to~$s_{\mathscr F}(z+r\cu_0)=\rho_ks_{\mathscr E}(z+r\cu_0)$. Inserting the first estimate into that bound gives the second estimate: we square the first bound and multiply by~$s_{\mathscr F}(z+r\cu_0)r^{d-2}$.

\smallskip

For~(ii), the global bound~\eqref{eq:mfd-5} for~$u$ gives~$\|u(z+r\,\cdot)-c\|_{C^\beta(\partial \cu_0)}\leq r^\alpha[u]_\alpha$, and the bound~\eqref{eq:mfd-3} on grid cubes bounds the factor in~\eqref{eq:mfd-2} for~$\mathscr F$ on~$z+r\cu_0$ by~$Kr^{-\kappa}$. Then~\eqref{eq:mfd-2} gives the bound~$Kr^{-\kappa}r^{d-2}r^{2\alpha}[u]_\alpha^2$.
\end{proof}

\begin{proposition}[Comparison from a positive density of small ratios]
\label{mfd:prop-density}
Assume that the parameters are chosen as above. There exists a deterministic~$C_*<\infty$, depending on~$d$,~$\kappa_0$,~$H_1$ and the constants of the good event, such that the following holds for every~$a>0$. Let~$S$ be a deterministic set of~$3^{H_1}$-adic depths of upper density at least~$\kappa_0$ with~$\rho_{H_1n}\leq a$ for every~$n\in S$. Then, on every cube~$z+\cu_j$ with~$z\in\mathbb Q^d$ and~$j\in\mathbb Z$:
\begin{itemize}
\item the domains satisfy~$D(\mathscr E)\subseteq D(\mathscr F)$,
\item for every~$v\in D(\mathscr E)$,
\begin{equation}
	\mathscr F(v)\leq C_*a\mathscr E(v)\,.
	\label{eq:mfd-25}
\end{equation}
\end{itemize}
\end{proposition}

\begin{proof}
Fix~$z_0\in\mathbb Q^d$ and~$j_0\in\mathbb Z$. Translate the~$3^{H_1}$-adic grid so that its origin is a corner of~$z_0+\cu_{j_0}$. Since the side of~$z_0+\cu_{j_0}$ is triadic, every sufficiently fine~$3^{H_1}$-adic mesh partitions~$z_0+\cu_{j_0}$ exactly. Fix~$u=G_{\mathscr E}f$ with smooth~$f$, fix~$\ep_{\rm Leb}>0$, and let~$\nu\coloneqq\Gamma_{\mathscr E}(u)+\ep_{\rm Leb}\,dx$. By~\eqref{eq:mfd-3} with~$\kappa<2$, the source term in Lemma~\ref{mfd:lem-goodext} is at most~$Kr^{d+2-\kappa}\|f\|_\infty^2\leq\ep_{\rm Leb} r^d$ for cells of side below a number~$s_0$ depending on~$u$,~$\ep_{\rm Leb}$ and the environment.

\smallskip

Fix a~$3^{H_1}$-adic mesh of side at most~$s_0$ that partitions~$z_0+\cu_{j_0}$; all cells below are subcubes of the finitely many cubes of this mesh, and we apply Proposition~\ref{mfd:prop-allchain} to each of them.

\smallskip

\emph{Step 1.} We stop the subdivision at good cells that carry a fixed fraction of the mass of their parents, and bound the sum, over these cells~$z+r\cu_0$, of the minimal~$\mathscr F$-energy~$\Lambda_{\mathscr F,z+r\cu_0}(u|_{\partial(z+r\cu_0)})$ with the boundary values of~$u$. Within each cube of the mesh of side at most~$s_0$ we follow the subdivision tree and stop at the first child~$z+r\cu_0$ that lies at a~$3^{H_1}$-adic depth in~$S$, is good and padded, and satisfies~$\nu(z'+3^{H_1}r\cu_0)\leq3^{H_1D}\nu(z+r\cu_0)$. By Lemma~\ref{mfd:lem-goodext}, each cell at which the subdivision stops satisfies~$\Lambda_{\mathscr F,z+r\cu_0}(u|_{\partial(z+r\cu_0)})\leq Ca[\nu(z'+3^{H_1}r\cu_0)+\ep_{\rm Leb} r^d]\leq Ca3^{H_1D}\nu(z+r\cu_0)$. The interiors of these cells are disjoint. Hence the sum of~$\Lambda_{\mathscr F,z+r\cu_0}(u|_{\partial(z+r\cu_0)})$ over them is at most~$C_*a\nu(z_0+\cu_{j_0})$ with~$C_*=C3^{H_1D}$.

\smallskip

\emph{Step 2.} We show that the sum of~$\Lambda_{\mathscr F,z+r\cu_0}(u|_{\partial(z+r\cu_0)})$ over the remaining cells, those on chains along which the subdivision never stops, tends to zero. There are arbitrarily large depths~$m$ at which at least~$\frac{3\kappa_0}{4}m$ of the~$3^{H_1}$-adic depths below the mesh of side at most~$s_0$ belong to~$S$. Along a chain on which the subdivision has not stopped by such a depth~$m$, every step to a good padded child, at a~$3^{H_1}$-adic depth in~$S$, fails the condition~$\nu(z'+3^{H_1}r\cu_0)\leq3^{H_1D}\nu(z+r\cu_0)$. Since~$\nu$ is monotone along the chain and drops by more than the factor~$3^{H_1D}$ at each of these steps, while every leaf satisfies~$\nu(\text{leaf})\geq\ep_{\rm Leb}|\text{leaf}|=\ep_{\rm Leb}3^{-H_1dm}$ on the~$3^{H_1}$-adic mesh of depth~$m$, the number of these steps is at most
\begin{equation*}
	\frac{d}{D}m+\frac{\log(\nu(z_0+\cu_{j_0})/\ep_{\rm Leb})}{D\log3^{H_1}}
	\leq\frac{\kappa_0}{8}m+O_{z_0+\cu_{j_0},\ep_{\rm Leb},u,\omega}(1)\,.
\end{equation*}
By~\eqref{eq:mfd-23}, at most~$\theta m+B_\omega\leq\frac{\kappa_0}{8}m+B_\omega$ of the steps at~$3^{H_1}$-adic depths in~$S$ are bad. For all sufficiently large such~$m$, the chain therefore has at least~$\frac{3\kappa_0}{4}m-\frac{\kappa_0}{4}m-O(1)\geq\frac{\kappa_0}{3}m$ steps at~$3^{H_1}$-adic depths in~$S$ to a child that is not padded. At a step to such a child there are at most~$6d\,3^{H_1(d-1)}$ choices of the child. At every other step there are at most~$3^{H_1d}$. Summing over the at most~$2^m$ positions of the steps to such children, we find that the number of leaves on chains along which the subdivision does not stop is at most
\begin{equation*}
	C_{z_0+\cu_{j_0}}2^m3^{H_1dm}(6d/3^{H_1})^{\nf{\kappa_0m}{3}}\leq C(z_0+\cu_{j_0},\ep_{\rm Leb},u,\omega)3^{H_1(d-\frac{\kappa_0}{8})m}\,,
\end{equation*}
where we used~$3^{H_1(\frac{\kappa_0}{3}-\frac{\kappa_0}{8})}\geq2(6d)^{\nf{\kappa_0}{3}}$. By Lemma~\ref{mfd:lem-goodext}(ii), the sum of~$\Lambda_{\mathscr F,z+r\cu_0}(u|_{\partial(z+r\cu_0)})$ over them is at most~$K3^{H_1(d-\frac{\kappa_0}{8})m}$ times~$3^{-H_1(d-2+2\alpha-\kappa)m}$, which tends to zero because~$d-2+2\alpha-\kappa>d-\nf{\kappa_0}{8}$.

\smallskip

\emph{Step 3.} We show that~$u\in D(\mathscr F)$ and~$\mathscr F(u)\leq C_*a\mathscr E(u)$. At a finite depth~$m$, the cells at which the subdivision stops and the remaining leaves form a finite partition of~$z_0+\cu_{j_0}$ into cells, on whose boundaries we take the values of~$u$. By Lemma~\ref{mfd:lem-skeleton}, the~$\mathscr F$-harmonic extensions of these boundary values glue to a function~$v_m\in D(\mathscr F)$ with~$\mathscr F(v_m)\leq C_*a\nu(z_0+\cu_{j_0})+o(1)$. By the maximum principle, each extension lies between the extreme boundary values of~$u$ on its cell, so~$\|v_m-u\|_\infty$ is at most the modulus of continuity of~$u$ at scale~$s_0$. For fixed~$s_0$ we take~$m$ large and then let~$s_0\to0$. Then~$v_m\to u$ uniformly, and the closedness of~$\mathscr F$ gives~$\mathscr F(u)\leq C_*a(\mathscr E(u)+\ep_{\rm Leb}3^{dj_0})$. Letting~$\ep_{\rm Leb}\downarrow0$ gives~$u\in D(\mathscr F)$ and~$\mathscr F(u)\leq C_*a\mathscr E(u)$.

\smallskip

\emph{Step 4.} We extend the inequality to~$D(\mathscr E)$. Let~$v\in D(\mathscr E)$. By Proposition~\ref{mfd:prop-killed-inverse}, there are smooth~$f_n$ with~$G_{\mathscr E}f_n\to v$ in the form norm of~$\mathscr E$. By Step~3, the sequence~$(G_{\mathscr E}f_n)$ is Cauchy in the form norm of~$\mathscr F$, and the closedness of~$\mathscr F$ identifies its limit with~$v$. This gives~$v\in D(\mathscr F)$ and~\eqref{eq:mfd-25}.
\end{proof}

\begin{theorem}[Two-sided comparison]
\label{mfd:thm-C0}
There exists a deterministic~$C_0\geq1$ such that, for every pair of subsequential limit forms and every~$0<\gmcdelta\leq\gmcdelta_*$, where~$\gmcdelta_*$ is the first threshold for~$\gmcdelta$ chosen at the beginning of this subsection, on every cube~$z+\cu_n$ with~$z\in\mathbb Q^d$ and~$n\in\mathbb Z$,
\begin{equation}
	D(\mathscr E)=D(\mathscr F)\quad\text{and}\quad C_0^{-1}\mathscr E\leq \mathscr F\leq C_0\mathscr E\,.
	\label{eq:mfd-26}
\end{equation}
\end{theorem}

\begin{proof}
Fix a nonzero smooth source~$f_0$ on the unit cube. By the positive and negative moments of~$\langle f_0,G_{\mathscr E}f_0\rangle$ and~$\langle f_0,G_{\mathscr F}f_0\rangle$ in Lemma~\ref{mfd:lem-response-reciprocal}, choose a deterministic~$M\geq1$, uniform over the subsequential limit forms and over~$0<\gmcdelta\leq\gmcdelta_*$, so that, with positive probability, both~$\langle f_0,G_{\mathscr E}f_0\rangle$ and~$\langle f_0,G_{\mathscr F}f_0\rangle$ lie in~$[M^{-1},M]$. It suffices to choose~$M$ so that each of the four tail probabilities is less than~$\nf18$. Let~$a_0\coloneqq(2C_*M^2)^{-1}$ and let~$S_-$ be the set of~$3^{H_1}$-adic depths~$n$ with~$\rho_{H_1n}\leq a_0$. Since the ratios~$\rho_k$ are deterministic, so is~$S_-$. If~$S_-$ has upper density at least~$\kappa_0$, then Proposition~\ref{mfd:prop-density} gives~$D(\mathscr E)\subseteq D(\mathscr F)$ and~$\mathscr F\leq C_*a_0\mathscr E$ on~$D(\mathscr E)$. The variational formula~\eqref{eq:mfd-18} then gives, almost surely,
\begin{align*}
	\langle f_0,G_{\mathscr F}f_0\rangle
	&=\sup_{v\in D(\mathscr F)}\bigl\{2\langle f_0,v\rangle-\mathscr F(v)\bigr\}
	\geq\sup_{v\in D(\mathscr E)}\bigl\{2\langle f_0,v\rangle-C_*a_0\mathscr E(v)\bigr\}
	\notag \\ &
	=(C_*a_0)^{-1}\langle f_0,G_{\mathscr E}f_0\rangle=2M^2\langle f_0,G_{\mathscr E}f_0\rangle\,,
\end{align*}
which contradicts the choice of~$M$. Hence~$S_-$ has upper density less than~$\kappa_0$. Exchanging~$\mathscr E$ and~$\mathscr F$ gives the same bound for the set~$S_+$ of~$3^{H_1}$-adic depths with~$\rho_{H_1n}\geq a_0^{-1}$. The complement of~$S_-\cup S_+$ therefore has lower density greater than~$1-2\kappa_0>\kappa_0$, since~$\kappa_0<\nf13$, and on this complement~$a_0<\rho_{H_1n}<a_0^{-1}$. Proposition~\ref{mfd:prop-density}, applied in both directions with~$a=a_0^{-1}$, gives~$\mathscr F\leq C_*a_0^{-1}\mathscr E$ and~$\mathscr E\leq C_*a_0^{-1}\mathscr F$ with equal domains. Thus~\eqref{eq:mfd-26} holds with~$C_0=C_*/a_0=2C_*^2M^2$.
\end{proof}

\subsubsection{The optimal constants are deterministic}

For~$a>0$ we write~$a\mathscr E\leq \mathscr F$ if~$a\mathscr E^{z+\cu_n}\leq \mathscr F^{z+\cu_n}$ on~$D(\mathscr E^{z+\cu_n})=D(\mathscr F^{z+\cu_n})$ for every~$z\in\mathbb Q^d$ and~$n\in\mathbb Z$, and we optimize over all these cubes:
\begin{equation*}
	c_-\coloneqq\sup\{a:\ a\mathscr E\leq \mathscr F\}\qquad\text{and}\qquad c_+\coloneqq\inf\{a:\ \mathscr F\leq a\mathscr E\}\,.
\end{equation*}
The numbers~$c_-$ and~$c_+$ are the optimal lower and upper comparison constants. By~\eqref{eq:mfd-26},~$C_0^{-1}\leq c_-\leq c_+\leq C_0$, and~$c_-\mathscr E\leq \mathscr F\leq c_+\mathscr E$ by monotone limits. The comparison constants are defined for each realization; the next lemma shows that they are deterministic, which Theorem~\ref{mfd:thm-prop} uses to conclude that the proportionality constant is deterministic.

\begin{lemma}[The optimal comparison constants are deterministic]
\label{mfd:lem-endpoints}
The optimal comparison constants~$c_-$ and~$c_+$ are almost surely constant.
\end{lemma}

\begin{proof}
The comparison constants are measurable: on each cube~$z+\cu_n$ with~$z\in\mathbb Q^d$ and~$n\in\mathbb Z$,~$a\mathscr E\leq \mathscr F$ is equivalent, by the variational formula~\eqref{eq:mfd-18}, to~$\langle f,G_{\mathscr F}f\rangle\leq a^{-1}\langle f,G_{\mathscr E}f\rangle$ for all sources~$f$ in the countable dense family.

\smallskip

Resampling one field multiplies both energy measures by the same positive continuous weight, by~\eqref{eq:mfd-21}. By Lemma~\ref{mfd:lem-sincos}, an inequality between the two forms is equivalent to the same inequality between their energy measures. Multiplying both measures by a common weight preserves these inequalities in both directions. Thus resampling one field leaves~$c_-$ and~$c_+$ unchanged. Since these constants are bounded, Lemma~\ref{mfd:lem-resampling-sum}, with no fields retained, gives~$\|c_--\E c_-\|_1=\|c_+-\E c_+\|_1=0$.
\end{proof}

\subsubsection{Concentration of the difference between the energies}

Assume~$c_+>c_-$. We prove a concentration bound for the difference of the two energies, with an error proportional to~$c_+-c_-$.

\smallskip

We first check that the minimal energies of~$\mathscr E$ and~$\mathscr F$ with given boundary values on a cell are taken over the same class of admissible functions. Let~$z+3^{-k}\cu_0$ be a cell and let~$z_0+\cu_{j_0}\supset z+3^{-k}\cu_0$ be a larger cube with~$z_0\in\mathbb Q^d$ and~$j_0\in\mathbb Z$. A boundary value~$\ell$ on~$\partial(z+3^{-k}\cu_0)$, affine or of class~$C^\beta$ with~$\beta>\nf12$, is realized by an~$\mathscr E$-harmonic extension through a neighborhood of~$\partial(z+3^{-k}\cu_0)$, which belongs to~$D(\mathscr E^{z_0+\cu_{j_0}})$. By~\eqref{eq:mfd-26} this extension lies in~$D(\mathscr F^{z_0+\cu_{j_0}})$. Its restriction to~$z+3^{-k}\cu_0$, together with the identification of the minimal energies in Propositions~\ref{mfd:prop-boundary} and~\ref{mfd:prop-gluing}, places it in the admissible class of both~$\Lambda_{\mathscr E,z+3^{-k}\cu_0}(\ell)$ and~$\Lambda_{\mathscr F,z+3^{-k}\cu_0}(\ell)$. Thus both minimal energies below are taken over the same class of functions with the prescribed boundary values, on which~$c_-\Lambda_{\mathscr E,z+3^{-k}\cu_0}\leq\Lambda_{\mathscr F,z+3^{-k}\cu_0}\leq c_+\Lambda_{\mathscr E,z+3^{-k}\cu_0}$ as quadratic forms.

\begin{lemma}[Difference of two minimizers]
\label{mfd:lem-diff}
Let~$u_{\mathscr E}$ and~$u_{\mathscr F}$ be the minimizers of~$\mathscr E$ and~$\mathscr F$ in~$z+3^{-k}\cu_0$ with the same affine boundary values, let~$c\in[c_-,c_+]$, and let~$w\coloneqq u_{\mathscr F}-u_{\mathscr E}$. Then~$\mathscr F(w,v)=-(\mathscr F-c\mathscr E)(u_{\mathscr E},v)$ for every~$v\in D(\mathscr E^{z+3^{-k}\cu_0})$, and
\begin{equation}
	\mathscr E(w)\leq\Bigl(\frac{c_+-c_-}{c_-}\Bigr)^2\mathscr E(u_{\mathscr E})\,.
	\label{eq:mfd-27}
\end{equation}
Moreover, for all~$u,v$ in the common domain and every Borel set~$U$,
\begin{equation*}
	|\Gamma_{\mathscr F-c\mathscr E}(u,v)|(U)\leq(c_+-c_-)\sqrt{\Gamma_{\mathscr E}(u)(U)\Gamma_{\mathscr E}(v)(U)}\,.
\end{equation*}
\end{lemma}

\begin{proof}
Subtracting the Euler equations~$\mathscr F(u_{\mathscr F},v)=0=\mathscr E(u_{\mathscr E},v)$ gives~$\mathscr F(w,v)=\mathscr F(u_{\mathscr F},v)-\mathscr F(u_{\mathscr E},v)=-\mathscr F(u_{\mathscr E},v)+c\mathscr E(u_{\mathscr E},v)=-(\mathscr F-c\mathscr E)(u_{\mathscr E},v)$. Since~$(c_--c)\mathscr E\leq \mathscr F-c\mathscr E\leq(c_+-c)\mathscr E$ and~$\max\{|c_--c|,|c_+-c|\}\leq c_+-c_-$, the symmetric form~$\mathscr F-c\mathscr E$ satisfies~$|(\mathscr F-c\mathscr E)(u,v)|\leq(c_+-c_-)\sqrt{\mathscr E(u)\mathscr E(v)}$. Taking~$v=w$ gives~$c_-\mathscr E(w)\leq \mathscr F(w)=-(\mathscr F-c\mathscr E)(u_{\mathscr E},w)\leq(c_+-c_-)\sqrt{\mathscr E(u_{\mathscr E})\mathscr E(w)}$, which is~\eqref{eq:mfd-27}.

\smallskip

Let~$z_0+\cu_{j_0}$ be the cube on which the forms are considered and let~$V\subseteq z_0+\cu_{j_0}$ be a Borel set. Lemma~\ref{mfd:lem-sincos}, applied to~$(\mathscr E,\mathscr F)$ with constant~$c_+$ and to~$(\mathscr F,\mathscr E)$ with constant~$c_-^{-1}$, gives~$c_-\Gamma_{\mathscr E}(\psi)\leq\Gamma_{\mathscr F}(\psi)\leq c_+\Gamma_{\mathscr E}(\psi)$ for~$\psi\in D(\mathscr E)\cap C_c(z_0+\cu_{j_0})$. Since~$\Gamma_{\mathscr F-c\mathscr E}(\psi)=\Gamma_{\mathscr F}(\psi)-c\Gamma_{\mathscr E}(\psi)$ and~$\max\{|c_--c|,|c_+-c|\}\leq c_+-c_-$, it follows that~$|\Gamma_{\mathscr F-c\mathscr E}(\psi)(V)|\leq(c_+-c_-)\Gamma_{\mathscr E}(\psi)(V)$. This extends to every~$\psi$ in the common domain, because~$D(\mathscr E)\cap C_c(z_0+\cu_{j_0})$ is dense in the form norm, which is equivalent to that of~$\mathscr F$ by~\eqref{eq:mfd-26}, and the energy measures~$\Gamma_{\mathscr E}(\psi)$ and~$\Gamma_{\mathscr F}(\psi)$ depend continuously in total variation on~$\psi$ in the form norm. Polarization therefore gives, for every~$s>0$,
\begin{equation*}
	4s|\Gamma_{\mathscr F-c\mathscr E}(u,v)(V)|
	\leq2(c_+-c_-)\bigl[\Gamma_{\mathscr E}(u)(V)+s^2\Gamma_{\mathscr E}(v)(V)\bigr]\,.
\end{equation*}
Optimizing in~$s$ gives~$|\Gamma_{\mathscr F-c\mathscr E}(u,v)(V)|\leq(c_+-c_-)\sqrt{\Gamma_{\mathscr E}(u)(V)\Gamma_{\mathscr E}(v)(V)}$. Sum this bound over a finite Borel partition of~$U$, apply Cauchy--Schwarz, and take the supremum over the partitions to obtain the total-variation bound.
\end{proof}

We normalize the minimal energies with affine boundary values on a cell. Fix~$c\in[c_-,c_+]$. Let~$\ell$ range over the affine boundary values in a fixed finite set that contains the coordinate functions and determines quadratic forms by polarization. Let~$T_{z+3^{-k}\cu_0}\coloneqq\sum_{i=1}^d\Lambda_{\mathscr E,z+3^{-k}\cu_0}(x_i)$ and, for one of these boundary values,~$\mathcal N\coloneqq(\Lambda_{\mathscr F,z+3^{-k}\cu_0}(\ell)-c\Lambda_{\mathscr E,z+3^{-k}\cu_0}(\ell))/T_{z+3^{-k}\cu_0}$. By polarization, these quantities give the entries of the matrix~$B_k$ in~\eqref{eq:mfd-31}. Let
\begin{equation*}
	\Psi\coloneqq\frac{1}{T_{z+3^{-k}\cu_0}}\sum_{\ell}\bigl[\Gamma_{\mathscr E}(u_{\mathscr E}^{\ell})+\Gamma_{\mathscr E}(u_{\mathscr F}^{\ell})\bigr]\,,
	\qquad
	\Psi_w\coloneqq\frac{1}{T_{z+3^{-k}\cu_0}}\sum_{\ell}\Gamma_{\mathscr E}(w^{\ell})\,,
\end{equation*}
where the sums run over the fixed affine boundary values and~$w^\ell\coloneqq u_{\mathscr F}^\ell-u_{\mathscr E}^\ell$. Then~$\Psi(z+3^{-k}\cu_0)\leq C(C_0)$ by~\eqref{eq:mfd-26}, and~$\Psi_w(z+3^{-k}\cu_0)\leq C(C_0)(c_+-c_-)^2$ by~\eqref{eq:mfd-27}. The measure~$\Psi$ inherits the Morrey bound~\eqref{eq:mfd-4} of the energy measures of the minimizers. For a Borel set~$U\subseteq z+3^{-k}\cu_0$ and a bounded~$g$ supported in~$U$, let~$\mathcal N^g$ denote the quantity~$\mathcal N$ computed for the forms~$\mathscr E^g$ and~$\mathscr F^g$.

\begin{lemma}[Perturbation of the normalized difference~$\mathcal N$]
\label{mfd:lem-relvar}
There exists~$C(C_0)<\infty$ such that, for every Borel set~$U\subseteq z+3^{-k}\cu_0$ and every bounded~$g$ supported in~$U$,
\begin{equation}
	|\mathcal N^g-\mathcal N|\leq C\|g\|_\infty e^{C\|g\|_\infty}
	\bigl[(c_+-c_-)\Psi(U)+\sqrt{\Psi(U)\Psi_w(U)}\bigr]\,.
	\label{eq:mfd-30}
\end{equation}
\end{lemma}

\begin{proof}
Let~$u_{\mathscr E}^g$ and~$u_{\mathscr F}^g$ minimize the weighted forms~$\mathscr E^g$ and~$\mathscr F^g$ with the same boundary values as~$u_{\mathscr E}$ and~$u_{\mathscr F}$, respectively. Let~$\widetilde{\mathcal N}_{z+3^{-k}\cu_0}\coloneqq\Lambda_{\mathscr F,z+3^{-k}\cu_0}(\ell)-c\Lambda_{\mathscr E,z+3^{-k}\cu_0}(\ell)$,~$v_{\mathscr E}\coloneqq u_{\mathscr E}^g-u_{\mathscr E}$ and~$v_{\mathscr F}\coloneqq u_{\mathscr F}^g-u_{\mathscr F}$. Subtract the weighted Euler equations to obtain, for every test function~$\phi$ with zero boundary values on~$\partial(z+3^{-k}\cu_0)$,
\begin{equation}
	\mathscr F^g(v_{\mathscr F}-v_{\mathscr E},\phi)=-\int_U(e^g-1)
	[d\Gamma_{\mathscr F}(w,\phi)+d\Gamma_{\mathscr F-c\mathscr E}(u_{\mathscr E},\phi)]-(\mathscr F^g-c\mathscr E^g)(v_{\mathscr E},\phi)\,.
	\label{eq:mfd-28}
\end{equation}
Test~\eqref{eq:mfd-28} with~$v_{\mathscr F}-v_{\mathscr E}$ and use Lemma~\ref{mfd:lem-diff} and the perturbation bound~\eqref{eq:mfd-14} for~$v_{\mathscr E}$. Cauchy--Schwarz gives
\begin{equation}
	\mathscr E^g(v_{\mathscr F}-v_{\mathscr E})\leq C\|g\|_\infty^2e^{C\|g\|_\infty}
	[\Gamma_{\mathscr E}(w)(U)+(c_+-c_-)^2\Gamma_{\mathscr E}(u_{\mathscr E})(U)]\,.
	\label{eq:mfd-29}
\end{equation}
As in the proof of Lemma~\ref{mfd:cor-14},
\begin{equation*}
	\widetilde{\mathcal N}_{z+3^{-k}\cu_0}^g-\widetilde{\mathcal N}_{z+3^{-k}\cu_0}
	=\int_U(e^g-1)\,d\bigl[\Gamma_{\mathscr F}(u_{\mathscr F})-c\Gamma_{\mathscr E}(u_{\mathscr E})\bigr]
	-\bigl[\mathscr F^g(v_{\mathscr F})-c\mathscr E^g(v_{\mathscr E})\bigr]\,.
\end{equation*}
Use the bilinear identities
\begin{align*}
	\Gamma_{\mathscr F}(u_{\mathscr F})-c\Gamma_{\mathscr E}(u_{\mathscr E})
	&=\Gamma_{\mathscr F-c\mathscr E}(u_{\mathscr E})+\Gamma_{\mathscr F}(w,u_{\mathscr F}+u_{\mathscr E})\,,\\
	\mathscr F^g(v_{\mathscr F})-c\mathscr E^g(v_{\mathscr E})
	&=(\mathscr F^g-c\mathscr E^g)(v_{\mathscr E})+\mathscr F^g(v_{\mathscr F}-v_{\mathscr E},v_{\mathscr F}+v_{\mathscr E})\,.
\end{align*}
Lemma~\ref{mfd:lem-diff} bounds the first line, in total variation on~$U$, by~$CT_{z+3^{-k}\cu_0}[(c_+-c_-)\Psi(U)+\sqrt{\Psi(U)\Psi_w(U)}]$. For the second,~\eqref{eq:mfd-14} for each form and~\eqref{eq:mfd-29} give
\begin{equation*}
	\mathscr E^g(v_{\mathscr E})+\mathscr E^g(v_{\mathscr F}+v_{\mathscr E})\leq C\|g\|_\infty^2e^{C\|g\|_\infty}T_{z+3^{-k}\cu_0}\Psi(U)\,,
\end{equation*}
and
\begin{equation*}
	\mathscr E^g(v_{\mathscr F}-v_{\mathscr E})\leq C\|g\|_\infty^2e^{C\|g\|_\infty}T_{z+3^{-k}\cu_0}[\Psi_w(U)+(c_+-c_-)^2\Psi(U)]\,.
\end{equation*}
{\emergencystretch=8em Cauchy--Schwarz and~$|(\mathscr F^g-c\mathscr E^g)(v_{\mathscr E})|\leq(c_+-c_-)\mathscr E^g(v_{\mathscr E})$ thus bound the correction by~$C\|g\|_\infty^2e^{C\|g\|_\infty}T_{z+3^{-k}\cu_0}[(c_+-c_-)\Psi(U)+\sqrt{\Psi(U)\Psi_w(U)}]$. Absorbing~$\|g\|_\infty^2$ into~$\|g\|_\infty e^{C\|g\|_\infty}$ proves the required bound for~$|\widetilde{\mathcal N}_{z+3^{-k}\cu_0}^g-\allowbreak \widetilde{\mathcal N}_{z+3^{-k}\cu_0}|/T_{z+3^{-k}\cu_0}$.\par}

\smallskip

Finally,~$|\widetilde{\mathcal N}_{z+3^{-k}\cu_0}|\leq C(c_+-c_-)T_{z+3^{-k}\cu_0}$,~$e^{-\|g\|_\infty}T_{z+3^{-k}\cu_0}\leq T_{z+3^{-k}\cu_0}^g\leq e^{\|g\|_\infty}T_{z+3^{-k}\cu_0}$, and~\eqref{eq:mfd-14} summed over the coordinate data gives
\begin{equation*}
	|T_{z+3^{-k}\cu_0}^g-T_{z+3^{-k}\cu_0}|\leq C\|g\|_\infty e^{C\|g\|_\infty}T_{z+3^{-k}\cu_0}\Psi(U)\,.
\end{equation*}
Apply these bounds to~$\mathcal N^g-\mathcal N=(\widetilde{\mathcal N}_{z+3^{-k}\cu_0}^g-\widetilde{\mathcal N}_{z+3^{-k}\cu_0})/T_{z+3^{-k}\cu_0}^g +\widetilde{\mathcal N}_{z+3^{-k}\cu_0}(T_{z+3^{-k}\cu_0}-T_{z+3^{-k}\cu_0}^g)/(T_{z+3^{-k}\cu_0}T_{z+3^{-k}\cu_0}^g)$.
\end{proof}

On a cell~$z+3^{-k}\cu_0$ at depth~$k$, define the coarse-grained matrices~$A_{\mathscr E}$ and~$A_{\mathscr F}$ by~$\Lambda_{\mathscr E,z+3^{-k}\cu_0}(y\cdot x)=3^{-kd}\,y\cdot A_{\mathscr E}y$ for~$y\in\Rd$, and similarly for~$\mathscr F$. Set~$c_k\coloneqq\E[\operatorname{tr}A_{\mathscr F}/\operatorname{tr}A_{\mathscr E}]$ and~$B_k\coloneqq(A_{\mathscr F}-c_kA_{\mathscr E})/\operatorname{tr}A_{\mathscr E}$.

\begin{proposition}[Concentration of the normalized difference~$B_k$]
\label{mfd:prop-conc}
\begin{itemize}
\item the constant~$c_k$ lies in~$[c_-,c_+]$,
\item the matrix~$B_k$ has mean~$\E B_k=0$,
\item Moreover, there exists~$a>0$, depending only on the dimension and on the regularity exponents fixed earlier, and not on~$C_0$ or~$p$, such that the following holds for every~$1\leq p<\infty$. If~$\gmcdelta$ is sufficiently small, depending on~$p$ and~$C_0$, then there is~$C_p(C_0)<\infty$, independent of~$k$, such that, for every~$k$ and every integer~$h\geq0$,
\begin{equation}
	\begin{aligned}
		&\|B_k\|_{L^p}\leq C_p(C_0)\gmcdelta(c_+-c_-)\,,\\
		&\bigl\|B_k-\E[B_k\mid\g_{-j}:|j-k|\leq h]\bigr\|_{L^p}\leq C_p(C_0)\gmcdelta(c_+-c_-)3^{-a h}\,,
	\end{aligned}
	\label{eq:mfd-31}
\end{equation}
\end{itemize}
\end{proposition}

\begin{proof}

\emph{Step 1.} The form order gives~$c_k\in[c_-,c_+]$. By Assumptions~\ref{a.g1} and~\ref{a.g3}, the joint law of~$A_{\mathscr E}$ and~$A_{\mathscr F}$, after their common scalar factor is removed, is invariant under signed permutations of the coordinates. Thus~$\E B_k$ is a scalar matrix. Its trace is zero by the definition of~$c_k$, so~$\E B_k=0$.

\smallskip

\emph{Step 2.} We bound the moments of the constant in the Morrey bound~\eqref{eq:mfd-4} for~$\Psi$, uniformly in~$k$. Rescale a cell to the unit cube and divide~$\mathscr E$ by~$s_{\mathscr E}(z+3^{-k}\cu_0)$ and~$\mathscr F$ by~$s_{\mathscr F}(z+3^{-k}\cu_0)$. The fields of scale smaller than the cell then have their original law, and the sum of the fields of larger scale, with its value at the cell center subtracted, has the moments of Lemma~\ref{mfd:lem-infrared}, uniformly in~$k$. At finite cutoff, let~$u$ minimize the normalized energy with affine boundary values of slope~$e_1$. Test the Neumann variational formula for boundary flux~$e_1\cdot\mathbf n$ with scalar multiples of~$u$. The product of the normalized Dirichlet and Neumann energies is then at least~$(\int_{\cu_0}e_1\cdot\nabla u)^2=1$. This inequality passes to the extracted limits. The proof of Lemma~\ref{mfd:lem-cell-ellipticity}, applied to the normalized cell, bounds the positive moments of the normalized Dirichlet and Neumann energies, uniformly in the cutoff and in~$k$, and Fatou's lemma passes these bounds to the extracted limits. By the product inequality, each of the two energies is at least the reciprocal of the other, so the same bounds control the reciprocals. For some deterministic~$M<\infty$, both normalized minimal energies with slope~$e_1$ lie in~$[M^{-1},M]$ with positive probability. On that event, the order of the original minimal energies and their scaling give that~$\rho_k$ times the ratio of the two normalized energies lies in~$[c_-,c_+]$, so that~$c_-M^{-2}\leq \rho_k\leq c_+M^2$. The last bounds hold deterministically because~$\rho_k$ is deterministic. Also the reciprocal of the sum of the normalized minimal energies over the coordinate directions is at most the normalized energy of the Neumann problem for~$\mathscr E$ with boundary flux~$e_1\cdot\mathbf n$. H\"older's inequality, the Morrey bounds for the two minimizers and these reciprocal moments now give the moments, of the finitely many orders fixed at the start of Section~\ref{sec:mfd-uniqueness}, of the constant in the Morrey bound for~$\Psi$, uniformly in~$k$.

\smallskip

\emph{Step 3.} Let~$p\geq2$, and consider one of the fields of scale smaller than the cell, of scale~$r$. Use the cutoff functions and strips of the proof of Lemma~\ref{mfd:lem-15}, with boxes of side~$r^{(t-d+1)/(t+1)}$ and strips of width~$4\sqrt d\,r$. The restrictions of the field to the separated interior subboxes are independent by~\ref{a.g1}, and by Lemma~\ref{mfd:lem-borel-weights} the whole field can be removed from the limit forms. With~$b$ from Lemma~\ref{mfd:lem-strips}, the original measure satisfies~$\Psi(\mathrm{strips})\leq Kr^b$. By~\eqref{eq:mfd-30}, the change of~$\mathcal N$ caused by deleting the field on the strips is at most~$C(c_+-c_-)Ye^{CY}(Kr^b+K^{\nf{1}{2}}r^{\nf{b}{2}})$, where~$Y$ is the sup norm of the field. As in the proof of Lemma~\ref{mfd:lem-15},~\eqref{eq:mfd-14} controls the changed minimizers, so the energy measures after the cutoff function is applied to the field satisfy~$\widehat\Psi(z+3^{-k}\cu_0)\leq C$,~$\widehat\Psi_w(z+3^{-k}\cu_0)\leq C(c_+-c_-)^2$ and~$\max_U\widehat\Psi(U)\leq Ce^{CY}Kr^b$. The bounds on~$\widehat\Psi(z+3^{-k}\cu_0)$ and~$\widehat\Psi_w(z+3^{-k}\cu_0)$ follow from the order between the two forms, which is preserved when the same cutoff function is applied to the field in both forms. The new ingredient in the Efron--Stein estimate is the bound on the sum of squares
\begin{align*}
	\sum_U[(c_+-c_-)\widehat\Psi(U)
	+\sqrt{\widehat\Psi(U)\widehat\Psi_w(U)}]^2
	& \leq2\max_U\widehat\Psi(U)
	[(c_+-c_-)^2\widehat\Psi(z+3^{-k}\cu_0)+\widehat\Psi_w(z+3^{-k}\cu_0)]
	\notag \\ & \leq C(c_+-c_-)^2\max_U\widehat\Psi(U)\,.
\end{align*}
Thus the first and last terms of~\eqref{eq:mfd-fourterm} in Lemma~\ref{mfd:lem-15}, and the conditional Efron--Stein estimate, all keep the factor~$c_+-c_-$. By~\eqref{eq:mfd-layernorm}, the norms of the fields are absorbed after reducing the exponent~$b$ by a fixed fraction. For an entry~$X$ of~$B_k$, resampling the field at relative index~$j$ consequently gives
\begin{equation*}
	\|X-X^{(j)}\|_p\leq C_p(C_0)\gmcdelta(c_+-c_-)
	\begin{cases}3^{-a j},&j\geq0,\\3^j,&j<0.\end{cases}
\end{equation*}
For~$j<0$ this follows directly from~\eqref{eq:mfd-30} and the bound of Step~2 on the sum of the fields of larger scale, with their value at the cell center subtracted. Choose~$a$ below a fixed fraction of~$b$. Apply Lemma~\ref{mfd:lem-resampling-sum} with no fields retained and then with the fields within distance~$h$ of~$k$ retained. Centering, polarization and the two geometric sums give~\eqref{eq:mfd-31}.
\end{proof}

By the concentration bound, we can add to the good event the requirement that the minimal energies of~$\mathscr F$ and~$c_k\mathscr E$ with affine boundary values differ by at most a small multiple of~$c_+-c_-$ times the minimal energy of~$\mathscr E$ with the same boundary values.

\begin{corollary}[Minimal energies with affine boundary values on good cells]
\label{mfd:cor-32}
Let~$\varpi>0$. The good event of Section~\ref{mfd:sec-branches} can be strengthened by one additional condition, with a threshold for~$\gmcdelta$ and a bad-cell density~$\theta$ that do not depend on~$c_+-c_-$, such that, on every good cell~$z+3^{-k}\cu_0$ at depth~$k$ and for every affine boundary value~$\ell$,
\begin{equation}
	|\Lambda_{\mathscr F,z+3^{-k}\cu_0}(\ell)-c_k\Lambda_{\mathscr E,z+3^{-k}\cu_0}(\ell)|\leq\varpi(c_+-c_-)\Lambda_{\mathscr E,z+3^{-k}\cu_0}(\ell)\,.
	\label{eq:mfd-32}
\end{equation}
\end{corollary}

\begin{proof}
Apply the last assertion of Lemma~\ref{mfd:lem-witness} to~$B_k/(c_+-c_-)$, with the additional condition~$|B_k/(c_+-c_-)|\leq\varpi c_0$, where~$c_0>0$ is the deterministic constant of the next paragraph. Its hypotheses are precisely~\eqref{eq:mfd-31}. Reserve half the allowed density of bad cells for this condition. Choose the rate~$\varkappa$ of the events~$W_{z+3^{-k}\cu_0,h}$ of Lemma~\ref{mfd:lem-witness} larger than the exponential growth rate~$d\log3^{H_1}$ of the number of chains. Then choose the required moment order and decrease~$\gmcdelta$. Proposition~\ref{mfd:prop-allchain} gives the strengthened good event, with choices independent of~$c_+-c_-$.

\smallskip

We pass from~$B_k$ to~\eqref{eq:mfd-32}. Condition~(a) of the good event bounds the condition number of the normalized matrix~$A_{\mathscr E}$: the coarse-grained ellipticity constants~\eqref{e.coarse.grained.ellipticity} bound~$|A_{\mathscr E}|$ and~$|A_{\mathscr E}^{-1}|$. Hence~$\lambda_{\min}(A_{\mathscr E})\geq c_0\operatorname{tr}A_{\mathscr E}$ with a deterministic~$c_0>0$. For an affine function with slope~$y$, we obtain~$|y\cdot(A_{\mathscr F}-c_kA_{\mathscr E})y|\leq|B_k|\operatorname{tr}A_{\mathscr E}|y|^2\leq\varpi c_0(c_+-c_-)\operatorname{tr}A_{\mathscr E}|y|^2\leq\varpi(c_+-c_-)\,y\cdot A_{\mathscr E}y$, which is~\eqref{eq:mfd-32}.
\end{proof}

\subsubsection{Affine approximation of the boundary values}

For the next lemma, fix~$\nf12<\beta<\alpha<1$,~$0<\gamma<1$,~$\zeta>0$ and~$\ep_{\rm aff}>0$ such that
\begin{equation}
	\frac{d+\zeta}{2}
	-\gamma\left(\frac{d-2}{2}+\alpha
	+\frac{\alpha-\beta}{\alpha+d/2}(2-\alpha)\right)<0\,.
	\label{eq:mfd-33}
\end{equation}

Let~$\mathscr E$ be a subsequential limit form along the subsequence of Section~\ref{mfd:sec-local-form}, realized in the Skorokhod representation, let~$z_0\in\mathbb Q^d$ and~$j_0\in\mathbb Z$, let~$f\in C_c^\infty(z_0+\cu_{j_0})$ and~$u\coloneqq G_{\mathscr E}^{z_0+\cu_{j_0}}f$, and let~$\ep_{\rm Leb}>0$ and~$\nu\coloneqq\Gamma_{\mathscr E}(u)+\ep_{\rm Leb}\,dx$.

Choose~$H_1$ large depending on~$d$,~$\alpha$,~$\beta$,~$\gamma$,~$\zeta$,~$\ep_{\rm aff}$ and the constants of the good event fixed before~$H_1$, with~$\lfloor\gamma H_1\rfloor+6<H_1$. Then choose~$\ep_{\rm hom}$ small depending on these and~$H_1$, and~$\gmcdelta$ small depending on these,~$H_1$ and~$\ep_{\rm hom}$. These choices are independent of~$f$ and~$\ep_{\rm Leb}$.

\begin{lemma}[Affine approximation of the boundary values]
\label{mfd:lem-affine}
Almost surely, there is~$s_0>0$, which may depend on~$f$,~$\ep_{\rm Leb}$ and the environment, such that the following holds for every good child~$z+r\cu_0$, of side~$r<s_0$, of a cell~$z'+3^{H_1}r\cu_0\subset z_0+\cu_{j_0}$. If the comparison cube of~$z+r\cu_0$, of side~$3^{\lfloor\gamma H_1\rfloor+4}r$, has its closure in~$z'+3^{H_1}r\cu_0$, and~$\nu(z'+3^{H_1}r\cu_0)\leq3^{H_1(d+\zeta)}\nu(z+r\cu_0)$, then there is an affine function~$\ell_{z+r\cu_0}$ such that
\begin{equation}
	\Lambda_{\mathscr E,z+r\cu_0}(u|_{\partial(z+r\cu_0)}-\ell_{z+r\cu_0})\leq\ep_{\rm aff}\nu(z+r\cu_0)\,.
	\label{eq:mfd-34}
\end{equation}
\end{lemma}

\begin{proof}
The constants~$C$ below depend only on~$d$,~$\alpha$,~$\beta$,~$\gamma$,~$\zeta$ and the constants of the good event fixed before~$H_1$ (as in Lemma~\ref{mfd:lem-goodext}); they do not depend on~$H_1$,~$\ep_{\rm hom}$,~$\gmcdelta$,~$f$,~$\ep_{\rm Leb}$ or the environment. Take~$H_1$ with~$\lfloor\gamma H_1\rfloor+6<H_1$ and let~$r_1\coloneqq3^{\lfloor\gamma H_1\rfloor+4}r/729$. Condition~(c) of the good event compares~$s_{\mathscr E}$ at the sides~$729r_1$ and~$r$ within a factor two. The mass bound on the parent and the estimate of the oscillation by the energy bound the oscillation of~$u$ on the cube of side~$r_1$ concentric with~$z+r\cu_0$ by~$Cr_1^{\nf{(2-d)}{2}}s_{\mathscr E}(z+r\cu_0)^{-\nf12}3^{H_1\frac{d+\zeta}{2}}\sqrt{\nu(z+r\cu_0)}$, together with the source term~$Cr_1^2s_{\mathscr E}(z+r\cu_0)^{-1}\|f\|_\infty$. Condition~(b) on the comparison cube, of side~$729r_1$, and Lemma~\ref{l.cutoff.regularity.good.scale.estimates}, with the Laplace-harmonic comparison function on the concentric cube of side~$9r_1$ and its error measured on the concentric cube of side~$r_1$, show that the~$\underline L^2$ distance on the cube of side~$r_1$ between~$u$ and the Laplace-harmonic function on the cube of side~$9r_1$ with the boundary values of~$u$ is at most
\begin{equation*}
	C\ep_{\rm hom}r_1^{(2-d)/2}s_{\mathscr E}(z+r\cu_0)^{-\nf12}
	3^{H_1(d+\zeta)/2}\sqrt{\nu(z+r\cu_0)}
	+Cr_1^2s_{\mathscr E}(z+r\cu_0)^{-1}\|f\|_\infty\,.
\end{equation*}
The estimates at finite cutoff pass to the extracted solutions by uniform convergence and Corollary~\ref{mfd:cor-energy-measures}, since the energy measure~$\Gamma_{\mathscr E}(u)$ does not charge the faces of~$z'+3^{H_1}r\cu_0$. Let~$\ell_{z+r\cu_0}$ be the affine function that best approximates~$u$ in~$L^2(z+r\cu_0)$; we use it for both norms below. Approximating the Laplace-harmonic function by its first-order Taylor polynomial gains a factor~$(r/r_1)^2$ in~$L^2(z+r\cu_0)$. Restricting the error of the harmonic approximation from the concentric cube of side~$r_1$ to~$z+r\cu_0$ introduces a factor~$(r_1/r)^{\nf d2}$ in volume-normalized~$L^2$. The local H\"older bound on~$z+r\cu_0$ gains a factor~$(r/r_1)^\alpha$, since the affine projection is bounded on~$C^\alpha$. Since~$3^{H_1\gamma} r/27<r_1\leq 3^{H_1\gamma} r/9$, after division by~$r^{\nf{(2-d)}{2}}\sqrt{\nu(z+r\cu_0)/s_{\mathscr E}(z+r\cu_0)}$ and reducing~$\ep_{\rm hom}$ after fixing~$H_1$, the rescaled~$L^2$ norm of~$u-\ell_{z+r\cu_0}$ on~$z+r\cu_0$ is therefore at most
\begin{equation*}
	C3^{H_1(\frac{d+\zeta}{2}-\gamma\frac{d+2}{2})}
	+C\ep_{\rm hom}3^{H_1(\frac{d+\zeta}{2}+\gamma)}\,,
\end{equation*}
and its rescaled~$C^\alpha$ norm on~$z+r\cu_0$ is at most~$C3^{H_1(\frac{d+\zeta}{2}-\gamma(\frac{d-2}{2}+\alpha))}$. Interpolate between these two bounds with the exponent~$(\alpha-\beta)/(\alpha+d/2)$ and apply the bound~\eqref{eq:mfd-2} for H\"older boundary values on the good cell, where the factor in it is comparable to~$s_{\mathscr E}(z+r\cu_0)$. The spatial and coefficient factors cancel, giving
\begin{align*}
	\Lambda_{\mathscr E,z+r\cu_0}(u|_{\partial(z+r\cu_0)}-\ell_{z+r\cu_0})
	&\leq C\bigl(3^{H_1(\frac{d+\zeta}{2}-\gamma\frac{d+2}{2})}
	+\ep_{\rm hom}3^{H_1(\frac{d+\zeta}{2}+\gamma)}\bigr)^{\nf{2(\alpha-\beta)}{(\alpha+d/2)}}\\
	&\qquad\cdot3^{2H_1(1-\frac{\alpha-\beta}{\alpha+d/2})(\frac{d+\zeta}{2}-\gamma(\frac{d-2}{2}+\alpha))}\nu(z+r\cu_0)\,.
\end{align*}
At~$\ep_{\rm hom}=0$, half the exponent of~$3^{H_1}$ is the left side of~\eqref{eq:mfd-33}, which is negative. First choose~$H_1$ to make this contribution smaller than~$\frac{\ep_{\rm aff}}{4}$. Next choose~$\ep_{\rm hom}$ so small that the contribution of the error of the harmonic approximation is smaller than~$\frac{\ep_{\rm aff}}{4}$. Then decrease~$\gmcdelta$ so that the enlarged good event has the required density of bad cells along every chain. After the same normalization, the source terms are at most~$C_{3^{H_1}}\ep_{\rm Leb}^{-\nf12}s_{\mathscr E}(z+r\cu_0)^{-\nf12}r\|f\|_{C^1}\leq K_{3^{H_1},\omega}\ep_{\rm Leb}^{-\nf12}r^{1-\nf\kappa2}\|f\|_{C^1}$, by~$\nu(z+r\cu_0)\geq\ep_{\rm Leb} r^d$ and~\eqref{eq:mfd-3}. They tend to zero because~$\kappa<2$. A small~$s_0$ makes them smaller than~$\frac{\ep_{\rm aff}}{4}$. The left side of~\eqref{eq:mfd-33} tends to~$-(1-\beta)/(1+d/2)<0$ as~$\alpha,\gamma\uparrow1$ and~$\zeta\downarrow0$, so the requirement on~$\alpha$ in Section~\ref{mfd:sec-local-form} and the choices of~$\gamma$ and~$\zeta$ can be met in every dimension.
\end{proof}

\subsubsection{Strict improvement of a comparison constant}

We now combine concentration for affine boundary values with the affine approximation of the solution's boundary values to rule out~$c_+>c_-$. Let~$\mathcal U\coloneqq\{n:c_{H_1n}\leq(c_-+c_+)/2\}$. Either this set or its complement has upper density at least~$\nf12$. Fix~$\varpi=\nf18$ in~\eqref{eq:mfd-32}. On a good cell~$z+3^{-k}\cu_0$ at a~$3^{H_1}$-adic depth in~$\mathcal U$, for every affine~$\ell$,
\begin{align*}
	P_{z+3^{-k}\cu_0}(\ell)&\coloneqq c_+\Lambda_{\mathscr E,z+3^{-k}\cu_0}(\ell)-\Lambda_{\mathscr F,z+3^{-k}\cu_0}(\ell)
	\notag \\
	&\geq(c_+-c_k-\varpi(c_+-c_-))\Lambda_{\mathscr E,z+3^{-k}\cu_0}(\ell)\geq\frac38(c_+-c_-)\Lambda_{\mathscr E,z+3^{-k}\cu_0}(\ell)\,.
\end{align*}
Since both minimal energies are taken over the same class of admissible functions,~$0\leq P_{z+3^{-k}\cu_0}\leq(c_+-c_-)\Lambda_{\mathscr E,z+3^{-k}\cu_0}$ on all admissible boundary values. If~$u$ is replaced on~$z+3^{-k}\cu_0$ by the~$\mathscr F$-harmonic function with the same boundary values~$g$, then the new function has~$\mathscr F$-energy~$\Lambda_{\mathscr F,z+3^{-k}\cu_0}(g)=c_+\Lambda_{\mathscr E,z+3^{-k}\cu_0}(g)-P_{z+3^{-k}\cu_0}(g)$ on~$z+3^{-k}\cu_0$, which is at most~$c_+\Gamma_{\mathscr E}(u)(z+3^{-k}\cu_0)-P_{z+3^{-k}\cu_0}(g)$.

\smallskip

Harmonic replacements will give~$\mathscr F(u)\leq(c_+-3(c_+-c_-)/256)\mathscr E(u)$ for~$u=G_{\mathscr E}f$, once the errors they introduce are controlled. Density extends this inequality to the common domain, contradicting the optimality of~$c_+$. If the complementary set of~$3^{H_1}$-adic depths has upper density at least~$\nf12$, then we interchange the forms and improve the optimal upper constant~$1/c_-$ instead. Both the energy saved by the replacements and the approximation error are proportional to~$c_+-c_-$, so the tolerances can be chosen independently of its size.

\smallskip

With~$\beta$ and~$\alpha$ as fixed in Section~\ref{mfd:sec-local-form}, choose~$\gamma$ and~$\zeta$ satisfying~\eqref{eq:mfd-33}. For~$H_1$ as in Lemma~\ref{mfd:lem-affine}, let~$b_1\coloneqq H_1-\lfloor\gamma H_1\rfloor-4$. Then~$b_1\geq1$,~$H_1\geq b_1+2$ and~$3^{-b_1}3^{H_1}r$ is the side of the comparison cube of~$z+r\cu_0$, used in conditions~(b)--(c) in Lemma~\ref{mfd:lem-affine}. With the dimensional constant~$C_d\geq1$ of the next lemma, the closure of this cube lies in~$z'+3^{H_1}r\cu_0$ for every child at distance at least~$C_d3^{-b_1}$ times the side of~$z'+3^{H_1}r\cu_0$ from~$\partial(z'+3^{H_1}r\cu_0)$, while~$C_d3^{-b_1}\leq81C_d3^{H_1(\gamma-1)}$ tends to zero as~$H_1$ grows. Averaging over finitely many translated grids makes the mass of the children close to the boundary of their parent small, for every measure.

\begin{lemma}[Shifted rational grids]
\label{mfd:lem-shifts}
Let~$H_1\geq b_1+2$. Translate the nested~$3^{H_1}$-adic grid by the shifts~$\mathbf k/(3^{b_1}+1)$ with~$\mathbf k\in\{0,\ldots,3^{b_1}\}^d$. Consider the children at depth~$n$ whose closure is within distance~$C_d3^{-b_1}3^{H_1(1-n)}$ of the boundary of their parent. Let~$C$ be the dimensional constant in the counting bounds below; it may be larger than~$C_d$. Then
\begin{itemize}
\item For every~$x$ and every~$n\geq1$,
\begin{multline*}
	(3^{b_1}+1)^{-d}\# \Bigl\{\mathbf k:\ \text{the depth-$n$ child of }x\text{ in the grid shifted by }\mathbf k/(3^{b_1}+1)
	\\
	\text{ \quad is within distance }C_d3^{-b_1}3^{H_1(1-n)}\text{ of the boundary of its parent} \Bigr\}
	\\
	\leq2d\Bigl(C_d3^{-b_1}+\frac1{3^{H_1}}+\frac1{3^{b_1}+1}\Bigr)
	\leq C3^{-b_1}\,.
\end{multline*}
\item Consequently, for every finite Borel measure~$\nu$ on a cube~$z_0+\cu_{j_0}$, the average over~$\mathbf k$ of the mass of the cells at depth~$n$ within distance~$C_d3^{-b_1}3^{H_1(1-n)}$ of the boundary of their parent is at most~$C3^{-b_1}\nu(z_0+\cu_{j_0})$.
\end{itemize}
\end{lemma}

\begin{proof}
Fix a coordinate~$i$. Relative to its parent at depth~$n-1$, the point~$x$ has fractional coordinate~$\{3^{H_1(n-1)}(x_i-\mathbf k_i/(3^{b_1}+1))\}$. Since~$\gcd(3^{H_1},3^{b_1}+1)=1$, multiplication by~$3^{H_1(n-1)}$ permutes the residues modulo~$3^{b_1}+1$. As~$\mathbf k_i/(3^{b_1}+1)$ ranges over the~$3^{b_1}+1$ shifts, these fractional coordinates form an equally spaced set of~$3^{b_1}+1$ points with a common rotation. If the child is within this distance of the boundary of its parent in coordinate~$i$, then~$x_i$ is within~$(C_d3^{-b_1}3^{H_1}+1)3^{-H_1n}$ of a face of the parent. Its fractional coordinate then lies in two intervals of total length at most~$2C_d3^{-b_1}+2/3^{H_1}$. A rotated equally spaced set has at most~$2(3^{b_1}+1)(C_d3^{-b_1}+1/3^{H_1})+2$ points in these two intervals. We divide by~$3^{b_1}+1$, sum over the~$d$ coordinates, and integrate the pointwise bound against~$\nu$.
\end{proof}

\begin{lemma}[Good cells carry a fixed fraction of the mass]
\label{mfd:lem-mass}
Suppose that~$\mathcal U$ has upper density at least~$\nf12$. Fix~$z_0\in\mathbb Q^d$,~$j_0\in\mathbb Z$,~$f\in C_c^\infty(z_0+\cu_{j_0})$ and~$\ep_{\rm Leb}>0$, let~$u\coloneqq G_{\mathscr E}f$ and~$\nu\coloneqq\Gamma_{\mathscr E}(u)+\ep_{\rm Leb}\,dx$, let~$\ep_{\rm coll}>0$, and assume that~$\theta+3^{-H_1\zeta}+C3^{-b_1}+\ep_{\rm coll}\leq\nf18$, with~$C$ the counting constant of Lemma~\ref{mfd:lem-shifts}. Then
\begin{itemize}
\item At arbitrarily large~$3^{H_1}$-adic depths~$n\in\mathcal U$, there are a shift~$\mathbf k\in\{0,\ldots,3^{b_1}\}^d$ and a finite family of cells~$z_i+r\cu_0$, indexed by~$i$, at depth~$n$ in the grid shifted by~$\mathbf k/(3^{b_1}+1)$ such that~$\nu(\bigcup_i(z_i+r\cu_0))\geq\nu(z_0+\cu_{j_0})/8$.
\item Every cell~$z_i+r\cu_0$ of such a family, of side~$r$ and with parent~$z_i'+3^{H_1}r\cu_0$,
\begin{itemize}
\item is good,
\item satisfies~$\nu(z_i'+3^{H_1}r\cu_0)\leq3^{H_1(d+\zeta)}\nu(z_i+r\cu_0)$,
\item has its concentric enlargement of side~$3^{\lfloor\gamma H_1\rfloor+4}r$ with closure in~$z_i'+3^{H_1}r\cu_0$.
\end{itemize}
\end{itemize}
\end{lemma}

\begin{proof}
In each shifted grid and at each depth, the children~$z+r\cu_0$ that violate~$\nu(z'+3^{H_1}r\cu_0)\leq3^{H_1(d+\zeta)}\nu(z+r\cu_0)$ have total mass at most~$3^{-H_1\zeta}\nu(z_0+\cu_{j_0})$, since each has mass less than~$3^{-H_1(d+\zeta)}\nu(z'+3^{H_1}r\cu_0)$ and each parent has at most~$3^{H_1d}$ children. By~\eqref{eq:mfd-24} for each of the~$(3^{b_1}+1)^d$ grids, the mass of the bad cells over the first~$m$ of the~$3^{H_1}$-adic depths is at most~$(\theta m+B_\omega^{\mathbf k})\nu(z_0+\cu_{j_0})$. The maximum over the finitely many grids of the random constants~$B^{\mathbf k}_\omega$ is finite on the intersection of the events of Proposition~\ref{mfd:prop-allchain} for these grids. After division by~$m$, this maximum is absorbed in the~$o(1)$ below. By Lemma~\ref{mfd:lem-shifts}, the mass of the cells within distance~$C_d3^{-b_1}3^{H_1(1-n)}$ of the boundary of their parent, averaged over~$\mathbf k$, is at most~$C3^{-b_1}\nu(z_0+\cu_{j_0})$ at every depth.

\smallskip

We omit a neighborhood of~$\partial(z_0+\cu_{j_0})$ of mass at most~$\ep_{\rm coll}\nu(z_0+\cu_{j_0})$. At sufficiently large~$3^{H_1}$-adic depths, the parents of the remaining points lie inside~$z_0+\cu_{j_0}$. Along depths~$m$ with~$\#(\mathcal U\cap\{1,\ldots,m\})\geq m/2-o(m)$, we sum the mass of the good cells with enlargement in the parent that satisfy the bound on their parents over~$n\in\mathcal U\cap\{1,\ldots,m\}$ and~$\mathbf k\in\{0,\ldots,3^{b_1}\}^d$, and divide by~$m(3^{b_1}+1)^d$. This normalized sum is at least
\begin{equation*}
	\Bigl[\frac12-\ep_{\rm coll}-\theta-3^{-H_1\zeta}-C3^{-b_1}
	-o(1)\Bigr]\nu(z_0+\cu_{j_0})>\frac{\nu(z_0+\cu_{j_0})}{4}
\end{equation*}
for large~$m$. A finite average that exceeds~$\nu(z_0+\cu_{j_0})/4$ has a summand that exceeds it, and since every fixed finite set of depths contributes~$o(m)$, such depths can be taken arbitrarily large.
\end{proof}

Let~$z_0\in\mathbb Q^d$ and~$j_0\in\mathbb Z$, let~$u\in D(\mathscr E^{z_0+\cu_{j_0}})\cap C(\overline{z_0+\cu_{j_0}})$, and let~$z_i+r_i\cu_0$, indexed by~$i$, be a finite disjoint family of cells of side at most~$r_{\max}$, whose threefold enlargements lie in~$z_0+\cu_{j_0}$. Assume that~$u|_{\partial(z_i+r_i\cu_0)}\in C^\alpha$ for every~$i$, with~$\alpha$ as fixed in Section~\ref{mfd:sec-local-form}. Define~$u^{(z_i,r_i)_i}$ to equal the~$\mathscr E$-harmonic extension of~$u|_{\partial(z_i+r_i\cu_0)}$ on each~$z_i+r_i\cu_0$ and to equal~$u$ elsewhere.

\begin{lemma}[Harmonic replacement]
\label{mfd:lem-replace}
\begin{itemize}
\item The difference~$u-u^{(z_i,r_i)_i}$ is the~$\mathscr E$-orthogonal projection of~$u$ onto~$\bigoplus_iD(\mathscr E^{z_i+r_i\cu_0})$.
\item The energy splits as~$\mathscr E(u)=\mathscr E(u^{(z_i,r_i)_i})+\mathscr E(u-u^{(z_i,r_i)_i})$.
\item If moreover~$u=G_{\mathscr E}f$, then~$\mathscr E(u-u^{(z_i,r_i)_i})\leq Kr_{\max}^{\nf32}\|f\|_2^2$, where~$K<\infty$ depends on~$z_0+\cu_{j_0}$ and on the environment through~\eqref{eq:mfd-1}, but not on~$f$,~$r_{\max}$ or the family.
\item In particular, if~$u=G_{\mathscr E}f$, then~$\sum_i[\Gamma_{\mathscr E}(u)(z_i+r_i\cu_0)-\Lambda_{\mathscr E,z_i+r_i\cu_0}(u|_{\partial(z_i+r_i\cu_0)})] =\mathscr E(u-u^{(z_i,r_i)_i})\to0$ as~$r_{\max}\to0$, uniformly over the choice of the family.
\end{itemize}
\end{lemma}

\begin{proof}
The Euler equations and strong locality identify~$w=u-u^{(z_i,r_i)_i}$ as the stated orthogonal projection and give both energy identities. Recovery sequences and~\eqref{eq:mfd-1} give~$\|v\|_{H^{\nf34}(z_0+\cu_{j_0})}^2\leq K\mathscr E(v)$; they also show that~$D(\mathscr E^{z_i+r_i\cu_0})\subset H^{\nf34}(z_i+r_i\cu_0)$ on each cell and that its functions have zero trace on~$\partial(z_i+r_i\cu_0)$. The fractional Poincar\'e inequality for zero-trace functions, rescaled to~$z_i+r_i\cu_0$, gives
\begin{equation*}
	\|w\|_2^2
	\leq Cr_{\max}^{\nf32}\sum_i[w]_{H^{\nf34}(z_i+r_i\cu_0)}^2
	\leq Cr_{\max}^{\nf32}[w]_{H^{\nf34}(z_0+\cu_{j_0})}^2
	\leq Kr_{\max}^{\nf32}\mathscr E(w)\,.
\end{equation*}
Here the seminorm is the Gagliardo double integral. Finally, for~$u=G_{\mathscr E}f$,~$\mathscr E(w)=\mathscr E(u,w)=\langle f,w\rangle\leq\|f\|_2\|w\|_2$, which proves the estimate.
\end{proof}

\begin{theorem}[Proportionality]
\label{mfd:thm-prop}
The optimal comparison constants coincide:~$c_+=c_-$. Hence there exists a deterministic~$c\in[C_0^{-1},C_0]$ such that
\begin{equation}
	\mathscr F=c\mathscr E\,.
	\label{eq:mfd-36}
\end{equation}
\end{theorem}

\begin{proof}
Suppose for a contradiction that~$c_+>c_-$.

\smallskip

\emph{Step 1.} We show that~$\mathcal U$ having upper density at least~$\nf12$ contradicts the definition of~$c_+$. Fix~$z_0\in\mathbb Q^d$ and~$j_0\in\mathbb Z$, let~$u\coloneqq G_{\mathscr E}f$ with smooth~$f$, fix~$\ep_{\rm Leb}>0$, let~$\nu\coloneqq\Gamma_{\mathscr E}(u)+\ep_{\rm Leb}\,dx$ and choose~$\ep_{\rm aff}\leq3/1024$ in Lemma~\ref{mfd:lem-affine}. Choose the parameters in the order above, enlarging~$H_1$ so that~$\theta+3^{-H_1\zeta}+C3^{-b_1}+\ep_{\rm coll}\leq\nf18$ with~$b_1=H_1-\lfloor\gamma H_1\rfloor-4$, before fixing~$\ep_{\rm hom}$ and~$\gmcdelta$. Let the cells~$z_i+r_i\cu_0$ form the family of Lemma~\ref{mfd:lem-mass} at a sufficiently large~$3^{H_1}$-adic depth, with maximal side~$r_{\max}$. For each cell~$z_i+r_i\cu_0$ of this family, the triangle inequalities for~$\sqrt{P_{z_i+r_i\cu_0}}$ and~$\sqrt{\Lambda_{\mathscr E,z_i+r_i\cu_0}}$, applied to~$u|_{\partial(z_i+r_i\cu_0)}=\ell_{z_i+r_i\cu_0}+(u|_{\partial(z_i+r_i\cu_0)}-\ell_{z_i+r_i\cu_0})$, give~$P_{z_i+r_i\cu_0}(u|_{\partial(z_i+r_i\cu_0)})\geq(c_+-c_-)[\frac{3}{16}\Lambda_{\mathscr E,z_i+r_i\cu_0}(u|_{\partial(z_i+r_i\cu_0)})-4\ep_{\rm aff}\nu(z_i+r_i\cu_0)]$, by~\eqref{eq:mfd-34}. Summing over~$i$ and using Lemma~\ref{mfd:lem-replace}, we obtain
\begin{align*}
	\sum_iP_{z_i+r_i\cu_0}(u|_{\partial(z_i+r_i\cu_0)})
	&\geq(c_+-c_-)\Bigl[\frac{3}{16}\sum_i\Lambda_{\mathscr E,z_i+r_i\cu_0}(u|_{\partial(z_i+r_i\cu_0)})-4\ep_{\rm aff}\nu(z_0+\cu_{j_0})\Bigr]\\
	&\geq(c_+-c_-)\Bigl[\frac{3}{16}\Gamma_{\mathscr E}(u)\bigl(\textstyle\bigcup_i(z_i+r_i\cu_0)\bigr)
	-\frac{3}{16}Kr_{\max}^{\nf32}\|f\|_2^2\\
	&\qquad-4\ep_{\rm aff}(\mathscr E(u)+\ep_{\rm Leb}3^{dj_0})\Bigr]\,.
\end{align*}
Since~$\nu(\bigcup_i(z_i+r_i\cu_0))\geq\nu(z_0+\cu_{j_0})/8$, we have~$\Gamma_{\mathscr E}(u)(\bigcup_i(z_i+r_i\cu_0))\geq \mathscr E(u)/8-7\ep_{\rm Leb}3^{dj_0}/8$, and with~$\ep_{\rm aff}\leq3/1024$ this gives
\begin{equation*}
	\sum_iP_{z_i+r_i\cu_0}(u|_{\partial(z_i+r_i\cu_0)})\geq(c_+-c_-)\bigl[\frac{3}{256}\mathscr E(u)-C\ep_{\rm Leb}3^{dj_0}-o_{r_{\max}}(1)\bigr]\,.
\end{equation*}
Now~$u\in D(\mathscr F)$ by~\eqref{eq:mfd-26}. Let~$v$ be the function equal to the~$\mathscr F$-harmonic extension of~$u|_{\partial(z_i+r_i\cu_0)}$ on each~$z_i+r_i\cu_0$ and to~$u$ elsewhere. By strong locality and the bound~$\Gamma_{\mathscr F}\leq c_+\Gamma_{\mathscr E}$ off the cells,
\begin{align*}
	\mathscr F(v)
	&\leq c_+\Gamma_{\mathscr E}(u)\bigl((z_0+\cu_{j_0})\setminus\textstyle\bigcup_i(z_i+r_i\cu_0)\bigr)+\sum_i\Lambda_{\mathscr F,z_i+r_i\cu_0}(u|_{\partial(z_i+r_i\cu_0)})
	\notag \\ &
	=c_+\mathscr E(u)-\sum_i\bigl[c_+\Gamma_{\mathscr E}(u)(z_i+r_i\cu_0)-\Lambda_{\mathscr F,z_i+r_i\cu_0}(u|_{\partial(z_i+r_i\cu_0)})\bigr]\\ &
	\leq c_+\mathscr E(u)-\sum_iP_{z_i+r_i\cu_0}(u|_{\partial(z_i+r_i\cu_0)})\,,
\end{align*}
where the last inequality uses~$\Gamma_{\mathscr E}(u)(z_i+r_i\cu_0)\geq \Lambda_{\mathscr E,z_i+r_i\cu_0}(u|_{\partial(z_i+r_i\cu_0)})$. As~$r_{\max}\to0$ along the~$3^{H_1}$-adic depths of Lemma~\ref{mfd:lem-mass}, the maximum principle bounds~$\|v-u\|_\infty$ by the modulus of continuity of~$u$ at scale~$r_{\max}$, which tends to zero. Thus~$v\to u$ in~$L^2$. Lower semicontinuity of~$\mathscr F$ gives~$\mathscr F(u)\leq(c_+-3(c_+-c_-)/256)\mathscr E(u)+C(c_+-c_-)\ep_{\rm Leb}3^{dj_0}$.

\smallskip

The tolerances and~$\ep_{\rm coll}$ do not depend on~$\ep_{\rm Leb}$. Letting~$\ep_{\rm Leb}\downarrow0$ gives~$\mathscr F(u)\leq(c_+-3(c_+-c_-)/256)\mathscr E(u)$ for every~$u=G_{\mathscr E}^{z_0+\cu_{j_0}}f$ with smooth~$f$, on every cube~$z_0+\cu_{j_0}$ with~$z_0\in\mathbb Q^d$ and~$j_0\in\mathbb Z$. Step~4 of the proof of Proposition~\ref{mfd:prop-density} extends the bound to~$D(\mathscr E^{z_0+\cu_{j_0}})$ for every such cube. This contradicts the definition of~$c_+$ as an infimum over all these cubes.

\smallskip

\emph{Step 2.} We show that the complement of~$\mathcal U$ having upper density at least~$\nf12$ contradicts the definition of~$c_-$. On good cells at~$3^{H_1}$-adic depths outside~$\mathcal U$, we have~$\Lambda_{\mathscr F,z+3^{-k}\cu_0}(\ell)\geq(c_k-\varpi(c_+-c_-))\Lambda_{\mathscr E,z+3^{-k}\cu_0}(\ell)\geq(c_-+3(c_+-c_-)/8)\Lambda_{\mathscr E,z+3^{-k}\cu_0}(\ell)$. We regard~$\mathscr F$ as the first form and~$\mathscr E$ as the second. The comparison constants of~$\mathscr E$ relative to~$\mathscr F$ are~$1/c_+\leq1/c_-$, with difference~$1/c_--1/c_+$, and the gain at affine boundary values is at least~$(3/8)(1/c_--1/c_+)$, because~$c_-+3(c_+-c_-)/8\leq c_+$. The argument of Step~1, with the solutions~$G_{\mathscr F}f$ for smooth~$f$, improves the optimal upper comparison constant~$1/c_-$ of this pair, which contradicts the optimality of~$c_-$.

\smallskip

Hence~$c_+=c_-$, so~$\mathscr F=c_-\mathscr E$, and the constant in~\eqref{eq:mfd-36} is~$c_-$, which is deterministic by Lemma~\ref{mfd:lem-endpoints}.
\end{proof}

We now show that the proportionality constant~$c$ in Theorem~\ref{mfd:thm-prop} is one, using the fact that the homogenized matrix of the stationary coefficients~$A_N^0$ is the identity.

\subsection{Determination of the proportionality constant}
\label{mfd:sec-scalar}

In this subsection we show that the proportionality constant~$c$ of Theorem~\ref{mfd:thm-prop} is one, so that all subsequential limit forms coincide (Theorem~\ref{mfd:thm-c1}). Homogenization fixes the constant: the coarse-grained matrices of the stationary coefficients~$A_N^0$ on large cubes approach the identity in expectation, uniformly in~$N$ (Theorem~\ref{mfd:thm-eta}), and~\eqref{eq:mfd-37} states this for the minimal energies with affine boundary values. After removing the factor~$e^H$ from the two limit forms, both inherit this normalization, so the proportionality constant is one.

\smallskip

The coefficients~$A_N^0$ of~\eqref{eq:mfd-normalization} are stationary, have range of dependence at most~$\sqrt d$, and have homogenized matrix~$\Id$.

\begin{theorem}[Uniform qualitative homogenization]
\label{mfd:thm-eta}
There exists~$\gmcdelta_0>0$ such that, if~$0<\gmcdelta\leq\gmcdelta_0$, then, for every unit vector~$e$,
\begin{equation*}
	\lim_{k\to\infty}\sup_N\E J(\cu_k,e,e;A_N^0)=0\,.
\end{equation*}
\end{theorem}
\begin{proof}
The proof proceeds in three steps, and the quantity we track is the sum~$S_N(k)$ of the expectations~$\E J(\cu_k,e_j,e_j;A_N^0)$ over the coordinate directions, defined in Step~1. Step~1 proves uniform moment bounds and chooses cutoffs~$N_k$ for which~$S_{N_k}$ changes little between the scales~$k$ and~$4k$. Step~2 extracts a subsequence of these cutoffs along which the expectations of the coarse-grained matrices~$A_{N_k}^0(\cu_n)$ converge to a deterministic matrix~$M_1$, and uses the finite range of dependence to show that, along it, the matrices~$(A_{N_k}^0)_*^{-1}(\cu_n)$ concentrate around a deterministic matrix~$M_2$. Step~3 uses the estimates of Section~\ref{s.self.similar} to show that~$M_1=M_2^{-1}$, so that~$M_1=M_2=\Id$. The heart of the matter is Step~3, where the small gaps in the subadditivity inequality obtained in Step~1 force the two limits to agree.

\smallskip

\emph{Step 1.} We obtain uniform moment bounds and choose cutoffs for which the expected energy changes little across a long interval of scales. Fix a sufficiently large even integer~$\xi$. Proposition~\ref{p.coarse.grained.bound}, subadditivity~\eqref{e.subadditivity.J} and stationarity give~$\sup_{N,k\geq0}\|\sup_{|e|=1}J(\cu_k,e,e;A_N^0)\|_{L^\xi}<\infty$. Moreover, the family
\begin{equation*}
	\left\{\Lambda_{1/4}(\cu_k;A_N^0)+\lambda_{1/4}(\cu_k;A_N^0)^{-1}:N,k\geq0\right\}
\end{equation*}
is tight. This follows from Lemma~\ref{l.multiscale.response.large.cubes}, with~$s=\nf18$,~$r=2$ and the initial error~$C\xi\log(2+\xi)\gmcdelta^2$, followed by~\eqref{e.bound.Lambdas.by.Es} and weighted Cauchy--Schwarz. Choose~$\xi$ large enough for that lemma and then reduce~$\gmcdelta$.

\smallskip

Write
\begin{equation*}
	S_N(k)\coloneqq\sum_{j=1}^d\E J(\cu_k,e_j,e_j;A_N^0)
	=\frac12\E\operatorname{tr}\bigl(A_N^0(\cu_k)+(A_N^0)_*^{-1}(\cu_k)-2\Id\bigr)\,.
\end{equation*}
By~\eqref{e.subadditivity.a} and~\eqref{e.annealed.matrices.def},~$\E A_N^0(\cu_k)$ and~$\E(A_N^0)_*^{-1}(\cu_k)$ decrease to~$\Id$. Thus~$S_N(k)\geq0$ decreases in~$k$, and the limit of~$\sup_N S_N(k)$ as~$k\to\infty$ exists. Choose~$N_k$ with~$S_{N_k}(4k)\geq\sup_N S_N(4k)-k^{-1}$. Then
\begin{equation*}
	\sup_{k\leq n\leq4k}\bigl(S_{N_k}(n)-S_{N_k}(4k)\bigr)
	\leq\sup_N S_N(k)-\sup_N S_N(4k)+k^{-1}\longrightarrow0\,.
\end{equation*}

\smallskip

\emph{Step 2.} We show that the coarse-grained matrices~$(A_{N_k}^0)_*^{-1}(\cu_n)$ concentrate and that the expectations of the coarse-grained matrices~$A_{N_k}^0(\cu_n)$ converge along these scales. By subadditivity, the average of~$(A_{N_k}^0)_*^{-1}(z+\cu_k)$ over the scale-$k$ subcubes of~$\cu_n$, minus~$(A_{N_k}^0)_*^{-1}(\cu_n)$, is positive semidefinite and has expected trace at most~$2(S_{N_k}(k)-S_{N_k}(n))$. This difference tends to zero in~$L^1$, uniformly for~$k\leq n\leq4k$. The moment bound of Step~1 and~\eqref{e.coarse.grained.matrix.deviation.by.J} upgrade this to~$L^2$. For~$n\geq2k$, the variance of the average is at most~$C3^{-d(n-k)}$, because the finite range of dependence splits the subcubes into a dimensional number of classes of independent cubes. Consequently, after extracting a subsequence, there are deterministic matrices~$M_1,M_2\geq \Id$ such that
\begin{equation*}
	\sup_{k\leq n\leq4k}\bigl\|\E A_{N_k}^0(\cu_n)-M_1\bigr\|\longrightarrow0\,,
\end{equation*}
and
\begin{equation*}
	\sup_{2k\leq n\leq4k}\bigl\|(A_{N_k}^0)_*^{-1}(\cu_n)-M_2\bigr\|_{L^2}\longrightarrow0\,.
\end{equation*}
For the first convergence, extract a limit of~$\E A_{N_k}^0(\cu_{4k})$ ; the difference~$\E A_{N_k}^0(\cu_n)-\E A_{N_k}^0(\cu_{4k})$ is positive semidefinite with trace at most~$2(S_{N_k}(n)-S_{N_k}(4k))$, which tends to zero by Step~1.

\smallskip

\emph{Step 3.} We use the estimates of Section~\ref{s.self.similar} to show that both limiting matrices equal the identity. Fix~$p\in\Rd$ and let~$v\coloneqq v(\cdot,\cu_{4k},p,M_2^{-1}p;A_{N_k}^0)$. Apply~\eqref{e.Besov.norms.gradient} with~$s=\nf14$ and~$L_0=2k$. The expected gaps in the subadditivity inequality~\eqref{e.subadditivity.a} tend to zero uniformly over~$2k\leq n\leq4k$ by Step~1. The mean-square differences of the coarse-grained matrices~$(A_{N_k}^0)_*^{-1}(\cu_n)$ tend to zero by Step~2 and stationarity. The weights of the two sums in~\eqref{e.Besov.norms.gradient} are summable, and~$(A_{N_k}^0)_*^{-1}(\cu_{4k})M_2^{-1}p-p\to0$ in~$L^2$. Since the coarse-grained ellipticity constants are tight, we obtain~$3^{-4k}\|\nabla v\|_{\Hminusul(\cu_{4k})}\to0$ in probability. Now~\eqref{e.J.L.bound}, with the same tight constants and the vanishing gaps in the subadditivity inequality, gives~$J(\cu_{4k-1},p,M_2^{-1}p;A_{N_k}^0)\to0$ in probability. The moment bounds make these energies uniformly integrable, so their expectations tend to zero. By~\eqref{e.coarse.grained.matrices}, the limit of these expectations is also~$\frac12p\cdot(M_1-M_2^{-1})p$. Thus~$M_1=M_2^{-1}$. Since~$M_1,M_2\geq \Id$, both equal~$\Id$, and~$S_{N_k}(4k)\to0$. The choice of~$N_k$ then shows that the limit of~$\sup_NS_N(k)$ is zero.
\end{proof}

Recall from the proof of Proposition~\ref{mfd:prop-boundary} the minimal energy~$\Lambda^0_{N,w+R\cu_0}$ with given boundary values for~$A_N^0$. The bound~\eqref{e.coarse.grained.matrix.deviation.by.J} and Cauchy--Schwarz turn Theorem~\ref{mfd:thm-eta} into the statement that, for every~$p\in\Rd$,
\begin{equation}
	\lim_{k\to\infty}\sup_N\E\left|
	\frac{\Lambda^0_{N,3^k\cu_0}(p\cdot x)}{3^{dk}}-|p|^2\right|=0
	\,.
	\label{eq:mfd-37}
\end{equation}

\begin{theorem}[Identification of the scalar]
\label{mfd:thm-c1}
The scalar~$c$ in~\eqref{eq:mfd-36} equals one. Consequently, all subsequential limit forms coincide.
\end{theorem}

\begin{proof}
Let~$\mathscr E^\varrho\coloneqq e^{-H}\mathscr E$ and~$\mathscr F^\varrho\coloneqq e^{-H}\mathscr F$ be the forms obtained by removing the factor~$e^H$ with the rule~\eqref{eq:mfd-21} for continuous weights, applied to the random continuous weight~$e^{-H}$ on a compact cube. They are subsequential limits of the same normalized sequence, whose coefficient is~$A_N^0$. Fix~$k\in\Z$ and~$p\neq0$. By~\eqref{eq:mfd-weighted-boundary-minimum}, with~$\varrho=e^{-H}$ and~$\mathscr E^{\varrho,3^{k+1}\cu_0}=\mathscr E^\varrho$ on the cube~$3^{k+1}\cu_0$,~$\Lambda_{N,3^k\cu_0}^0(p\cdot x)$ converges along the subsequence that defines~$\mathscr E$ to~$\Lambda_{\mathscr E^\varrho,3^k\cu_0}(p\cdot x)$. The uniform moments of these energies give uniform integrability. Hence
\begin{equation*}
	\E\Bigl|\frac{\Lambda_{\mathscr E^\varrho,3^k\cu_0}(p\cdot x)}{3^{dk}}-|p|^2\Bigr|
	\leq\sup_N\E\Bigl|\frac{\Lambda^0_{N,3^k\cu_0}(p\cdot x)}{3^{dk}}-|p|^2\Bigr|\,.
\end{equation*}
The same holds for~$\mathscr F^\varrho$ along its subsequence. Letting~$k\to\infty$,~\eqref{eq:mfd-37} shows that both~$\E|\Lambda_{\mathscr E^\varrho,3^k\cu_0}(p\cdot x)/3^{dk}-|p|^2|$ and its analogue for~$\mathscr F^\varrho$ tend to zero. But~$\mathscr F^\varrho=c\mathscr E^\varrho$ gives~$\Lambda_{\mathscr F^\varrho,3^k\cu_0}=c\Lambda_{\mathscr E^\varrho,3^k\cu_0}$ on every cube, so~$c|p|^2=|p|^2$ and~$c=1$.

\smallskip

Since all subsequential limit forms are therefore equal, and the minimal energies with prescribed boundary values and the Dirichlet Green quadratic forms at finite~$N$ are strongly relatively compact by Proposition~\ref{mfd:prop-response-compact}, every subsequence of these quantities has a further subsequence converging to the same limit, identified by~\eqref{eq:mfd-18} and Proposition~\ref{mfd:prop-boundary}. Hence the minimal energies with prescribed boundary values and the Dirichlet Green quadratic forms, for the countable family of cubes and data fixed in Section~\ref{mfd:sec-local-form}, converge along the full sequence of cutoffs. Convergence of the Neumann energies on general triadic cubes is established separately by the comparison with the smoothed volume sources in Proposition~\ref{mfd:prop-as-response-bank}(iv). This comparison uses no Neumann limit form; the uniform higher moments of the Neumann energies also give convergence in~$L^p$.
\end{proof}
\section{Convergence of the process}
\label{sec:mfd-process-convergence}

In this section we prove Theorem~\ref{t.scaling.limit}(i): almost surely, the quenched laws of the rescaled processes converge, uniformly for starting points in compact sets. The limit form of Section~\ref{sec:mfd-uniqueness} and the limit~$\mu$ of the reversible measures~$\mu_N$ determine the limiting resolvents. We first prove convergence in probability, and we then upgrade it to almost-sure convergence.

\smallskip

The difficulty in the first stage is that the densities of~$\mu_N$ are not bounded uniformly in~$N$, so the sources in the resolvent equations of the rescaled processes are random measures rather than functions. We control them through bounds, uniform in~$N$, on the masses of~$\mu_N$ in small balls. The difficulty in the second stage is that convergence in probability gives no almost-sure statement along the full sequence of cutoffs. We obtain one by adapting the comparison argument of~\citet[Section~3.4]{DGComparison}, strengthened by~\citet[Lemma~3.11]{Devlin2024}, to the finite-volume energies at two cutoffs. If the ratios of the deterministic normalizations at the two cutoffs were at most a fixed number less than one on a positive proportion of scales, harmonic replacement would bound the minimal energy for a fixed nonconstant affine boundary value at a finite cutoff by a fixed factor smaller than one times its positive limit, contradicting convergence in probability.

\smallskip

The section is organized as follows.
\begin{itemize}
\item In Section~\ref{s.fixed.coefficient} we construct the diffusions~$X$ and~$X^{(L)}$ from their elliptic resolvents, from every starting point.

\item In Section~\ref{mfd:sec-speed} we prove that the measures~$\mu_N$ converge and that their masses in small balls are bounded uniformly in~$N$. These bounds let us pass to the limit in the resolvent equations of the rescaled processes killed on leaving a cube. Mollifying the sources and applying the H\"older estimates of Section~\ref{sec:mfd-uniqueness} gives convergence in probability of the killed resolvents, uniformly in the starting point.

\item In Section~\ref{s.tightness} we prove tightness of the path laws.

\item In Section~\ref{mfd:sec-paths} we combine tightness with the convergence of the killed resolvents to identify the limiting path laws, which gives convergence in probability.

\item In Section~\ref{mfd:sec-finite-comparison} we compare the finite-volume energies at two cutoffs and obtain, by the Borel--Cantelli lemma, almost-sure convergence of the finite-volume energies.

\item In Sections~\ref{mfd:sec-pathwise-bounds} and~\ref{mfd:sec-as-paths} we prove bounds that hold for all cutoffs on one event of probability one, and we carry the almost-sure convergence over to the Dirichlet forms, the killed resolvents and the path laws.
\end{itemize}
The conventions of Section~\ref{sec:mfd-uniqueness} remain in force, except that Sections~\ref{s.fixed.coefficient} and~\ref{s.tightness} work with the diffusions~$X$ and~$X^{(L)}$ and their coefficients~$\a$ and~$\a_L$.

\subsection{Construction of the diffusions}
\label{s.fixed.coefficient}

In this subsection we construct the diffusions~$X$ and~$X^{(L)}$ with coefficients~$\a$ and~$\a_L$ from every starting point. We identify their resolvents killed on leaving a cube with those of the Dirichlet form~$\int\a|\nabla f|^2\,dx$ with zero boundary values; Section~\ref{mfd:sec-speed} uses this to write the killed resolvents of the rescaled processes as minimizers. The construction starts from the elliptic resolvents: weak solutions on increasing cubes converge to a resolvent, which defines a Feller semigroup, and moment bounds for its transition kernels give continuous paths. By Lemma~\ref{l.finite.cutoff.coefficient.convergence}, the coefficient~$\a$ is positive and locally~$C^{1,1}$. The same holds for~$\a_L$. The logarithmic derivative bound~\eqref{e.finite.cutoff.log.derivative.growth} gives at most linear growth of the drift, uniformly in~$L$. For either coefficient,~$|\nabla\log\a(x)|_2\leq K(1+|x|_2)$, with a random~$K<\infty$ independent of~$L$. The weak solutions of~$\lambda\a u-\nabla\cdot(\a\nabla u)=\a f$ on increasing bounded cubes, with zero boundary values, converge to the whole-space resolvent. To see the decay and uniqueness, let~$w(x)\coloneqq1+|x|_2^2$. Direct differentiation gives~$|\a^{-1}\nabla\cdot(\a\nabla w^{-\gamma})|\leq C\gamma w^{-\gamma}$ and~$\a^{-1}\nabla\cdot(\a\nabla w^\gamma)\leq C\gamma w^\gamma$ for~$0<\gamma\leq1$, where~$C\coloneqq2d+8+3K$. Fix~$\lambda>0$ and choose~$C\gamma\leq\lambda/2$. If~$f$ is supported in~$\{|x|_2\leq\mathsf r\}$, then~$2\lambda^{-1}\|f\|_\infty(1+\mathsf r^2)^\gamma w^{-\gamma}$ is a supersolution of~$\lambda-\a^{-1}\nabla\cdot(\a\nabla\,\cdot)$ with source~$|f|$. Comparison on each exhausting cube bounds the absolute value of its zero-boundary solution by this function, which tends to zero at infinity. For uniqueness among bounded solutions, compare the difference with~$\ep w^\gamma$ on large cubes and then let~$\ep\downarrow0$. Contraction extends the resolvent to~$C_0(\Rd)$. For~$f\in C_0(\Rd)$, uniform approximation by compactly supported data and the local energy estimates pass the weak equation to~$R_\lambda f$, so it is the unique bounded weak solution. Consequently, for every smooth~$v\in C_0(\Rd)$ with~$\a^{-1}\nabla\cdot(\a\nabla v)\in C_0(\Rd)$, uniqueness gives~$R_\lambda(\lambda v-\a^{-1}\nabla\cdot(\a\nabla v))=v$. The maximum principle gives positivity and~$\|R_\lambda\|\leq\lambda^{-1}$; uniqueness gives the resolvent identity. In particular,~$\lambda R_\lambda\phi\to\phi$ for~$\phi\in C_c^\infty(\Rd)$. Contraction extends this convergence to~$C_0(\Rd)$. Hille--Yosida~\citep[Theorem~19.11]{Kallenberg} therefore gives a Feller semigroup~$P_t$, with generator~$\a^{-1}\nabla\cdot(\a\nabla)$ on~$C_c^\infty(\Rd)$. Its positive operators have transition kernels by the Riesz representation theorem. These kernels have total mass one: the same resolvent identity puts~$w^{-\gamma}$ in the generator domain, and the bound~$\|\a^{-1}\nabla\cdot(\a\nabla w^{-\gamma})\|_\infty\leq C\gamma$ gives~$|P_tw^{-\gamma}-w^{-\gamma}|\leq Ct\gamma$. Letting~$\gamma\downarrow0$ proves~$P_t1=1$.

\smallskip

The weak resolvent equation also supplies the bounds needed for continuous paths, before any path realization is used. Use the regularization~$\Phi_{\mathsf r}(y)\coloneqq\frac{\mathsf r^2|y-x|_2^2}{\mathsf r^2+|y-x|_2^2}$, a concave function of the squared distance, where~$\mathsf r\geq1$. Both~$\Phi_{\mathsf r}-\mathsf r^2$ and~$\Phi_{\mathsf r}^2-\mathsf r^4$ belong to~$C_0(\Rd)$, as do their images under~$\a^{-1}\nabla\cdot(\a\nabla\,\cdot)$. The preceding weak-resolvent identity places them in the generator domain; adding back the constants is justified by~$P_t1=1$. Concavity in the squared distance and direct differentiation give~$\a^{-1}\nabla\cdot(\a\nabla\Phi_{\mathsf r})\leq C(w(x)+\Phi_{\mathsf r})$ and~$\a^{-1}\nabla\cdot(\a\nabla(\Phi_{\mathsf r}^2))\leq C(w(x)\Phi_{\mathsf r}+\Phi_{\mathsf r}^2)$, with~$C$ independent of~$\mathsf r$. The semigroup identity and Gronwall's inequality first bound~$P_t\Phi_{\mathsf r}(x)$ by~$Cw(x)te^{Ct}$; integrating this bound in the equation for~$P_t(\Phi_{\mathsf r}^2)(x)$ gives~$Cw(x)^2t^2e^{Ct}$. Letting~$\mathsf r\to\infty$ and using Fatou's lemma yields the second-moment bound and the fourth-moment bound~$\int_{\Rd}|y-x|_2^4P_t(x,dy)\leq Cw(x)^2t^2e^{Ct}$. In particular, for every compact~$F\subset\Rd$ and~$r>0$,~$\sup_{x\in F}P_t(x,\{|y-x|_2\geq r\})=O(t^2)=o(t)$. Kolmogorov extension constructs the finite-dimensional laws on a countable dense set of times. The semigroup identity and the fourth-moment bound give, on every time interval~$[0,T]$,~$\mathbf E_x|Z_t-Z_s|_2^4\leq C(x,T)|t-s|^2$ for~$s,t\in[0,T]$ and the coordinate process~$Z$, so the Kolmogorov continuity criterion extends these laws to continuous paths. We denote the resulting conservative diffusions for~$\a$ and~$\a_L$ by~$X$ and~$X^{(L)}$, respectively. To obtain the strong Markov property, approximate a finite stopping time~$S$ from above by~$S_n\coloneqq2^{-n}\lceil2^nS\rceil$. The Markov identity holds at each~$S_n$. Path continuity and continuity of~$P_tf$ for~$f\in C_0(\Rd)$ allow passage to~$S$; a monotone-class argument extends the identity to bounded Borel tests. Localization gives the identity on~$\{S<\infty\}$.

\smallskip

The semigroup is symmetric on~$L^2(\Rd,\a\,dx)$, and its resolvents agree with those of the closure of the form~$f\mapsto\int_{\Rd}\a|\nabla f|^2\,dx$ on~$C_c^\infty(\Rd)$. This follows by testing the weak equation and exhausting~$\Rd$ by cubes. Its generator is~$\a^{-1}\nabla\cdot(\a\nabla)$. Restricting the form to a cube~$z+r\cu_0$ with zero boundary values identifies its resolvents with those of the diffusion killed on leaving the cube, initially as~$\a\,dx$-almost-everywhere versions~\citep[Theorems~4.4.2--4.4.3]{FOT}. The identification holds at every interior starting point as follows. For~$f\in C_c^\infty(z+r\cu_0)$, let~$u$ be the zero-boundary solution of~$\lambda u-\a^{-1}\nabla\cdot(\a\nabla u)=f$. It is~$C^2$ in~$z+r\cu_0$ and continuous on~$\overline{z+r\cu_0}$, with value zero on~$\partial(z+r\cu_0)$. Take cubes~$z+r_j\cu_0$ with~$r_j<r$ increasing to~$r$, and compactly supported~$C^2$ functions agreeing with~$u$ near~$\overline{z+r_j\cu_0}$. Smooth approximation in the generator identity gives Dynkin's formula for these functions. Stopping at~$t\wedge\tau_{z+r_j\cu_0}$, where~$\tau_{z+r_j\cu_0}$ is the exit time of the diffusion from~$z+r_j\cu_0$, yields, for every~$x\in z+r_j\cu_0$,
\begin{equation*}
	u(x)=\mathbf E_x\left[\int_0^{t\wedge\tau_{z+r_j\cu_0}}e^{-\lambda s}f(X_s)\,ds
	+e^{-\lambda(t\wedge\tau_{z+r_j\cu_0})}u(X_{t\wedge\tau_{z+r_j\cu_0}})\right]\,.
\end{equation*}
Path continuity and the zero boundary value allow~$j\to\infty$; boundedness and the exponential factor then allow~$t\to\infty$. Thus~$u(x)=\mathbf E_x\int_0^{\tau_{z+r\cu_0}}e^{-\lambda s}f(X_s)\,ds$ for every~$x\in z+r\cu_0$. Approximation by uniformly bounded smooth data, using local elliptic regularity and dominated convergence, extends this identity to bounded continuous data. The same argument applies to~$\a_L$.

\smallskip

These laws also describe the usual SDE construction. Since~$\nabla\log\a$ is locally Lipschitz, there is a unique maximal strong solution from every~$x\in\Rd$ of
\begin{equation}
	dX_t=\sqrt2\,dW_t+\nabla\log\a(X_t)\,dt\,,
	\label{e.fixed.clock.sde}
\end{equation}
where~$W$ is standard Brownian motion on~$\Rd$. Fix a bounded cube~$z+r\cu_0$, a number~$\lambda>0$ and a function~$f\in C_c^\infty(z+r\cu_0)$, and let~$u$ solve~$\lambda\a u-\nabla\cdot(\a\nabla u)=\a f$ in~$z+r\cu_0$, with zero boundary values. Interior elliptic regularity and It\^o's formula, stopped first strictly inside~$z+r\cu_0$, give, for every~$x\in z+r\cu_0$, after exhaustion and letting the deterministic time horizon tend to infinity,
\begin{equation*}
	u(x)=\mathbf E_x^{\mathrm{SDE}}
	\left[\int_0^{\tau_{z+r\cu_0}}e^{-\lambda t}f(X_t)\,dt\right]\,.
\end{equation*}
Here~$\mathbf E_x^{\mathrm{SDE}}$ denotes expectation for the maximal solution of~\eqref{e.fixed.clock.sde}. Thus its killed resolvents agree with those of the analytic construction. Uniqueness of Laplace transforms identifies the killed semigroups, and the Markov property identifies their path laws. Exhausting~$\Rd$ identifies the laws up to their lifetimes, for every starting point. The same argument applies to~$\a_L$ and to the space and time rescalings in~\eqref{e.intro.rescaled.generator}; Proposition~\ref{tight:prop-tightness} makes these lifetimes infinite. For an open set~$U$, write~$\tau_U\coloneqq\inf\{t\geq0:X_t\notin U\}$ and let~$\tau_U^{(L)}$ be the exit time of~$X^{(L)}$ from~$U$. The diffusion~$X^{(L)}$ has generator~\eqref{e.intro.cutoff.generator}.

\subsection{The reversible measure and resolvent convergence in probability}
\label{mfd:sec-speed}

In this subsection we prove that the reversible measures~$\mu_N$ converge locally weakly, almost surely, to a locally finite measure~$\mu$, with mass at most a random constant times~$r^{d-\ep}$ in balls of radius~$r\leq1$, for every~$\ep\in(0,1)$ when~$\gmcdelta$ is small enough, depending on~$\ep$ (Proposition~\ref{mfd:prop-chaos-growth}). We also prove that the killed resolvents of the rescaled processes converge in probability, uniformly on each cube, to the minimizer of a functional built from the limit form and~$\mu$ (Propositions~\ref{mfd:prop-speed-resolvent} and~\ref{mfd:prop-uniform-resolvent}). The sources in the resolvent equations are the measures~$\mu_N$, whose densities are not bounded uniformly in~$N$. The mass bound lets us mollify them and apply the H\"older estimate for bounded sources.

\smallskip

Write~$\mathscr E$ for the unique limit form; it is measurable with respect to the random fields. Along the sequence in each Skorokhod representation of Section~\ref{mfd:sec-local-form}, the forms~$\mathscr E_N^{z+\cu_n}$ Mosco-converge to~$\mathscr E^{z+\cu_n}$ on~$L^2(z+\cu_n,dx)$, for every~$z\in\mathbb Q^d$ and~$n\in\mathbb Z$, and their coercivity constants are bounded.

\subsubsection{The chaos martingale and the mass of small balls}

The measures~$\mu_N^0$ and the reversible measures~$\mu_N=e^H\mu_N^0$ are defined in the paragraph following~\eqref{e.intro.rescaled.generator}.

\begin{lemma}[Positive moments of the chaos]
\label{mfd:lem-chaos-moments}
Let~$p\geq1$ be an integer. There exist~$C_p<\infty$ and~$c_p>0$ such that, for every bounded region, there is~$C<\infty$ such that, if~$\gmcdelta\leq c_p$, then, for every~$N\geq0$ and every cube~$z+r\cu_0$ of side~$r\leq1$ contained in this region,
\begin{equation}
	\E\bigl[\mu_N(z+r\cu_0)^p\bigr]\leq C\,r^{dp-C_p\gmcdelta^2}\,.
	\label{eq:mfd-38}
\end{equation}
Moreover,~$N\mapsto \mu_N^0(z+r\cu_0)$ is a nonnegative martingale, and~$\mu_N(z+r\cu_0)$ is a nonnegative martingale conditionally on~$H$.
\end{lemma}
\begin{proof}
Since~$\E e^{\g_{-j}(x)}=e^{\tau^2}$ for every~$x$ and the fields are independent, integration against fixed nonnegative test functions shows that~$\mu_N^0$ is a measure-valued martingale. The sum~$H$ of the fields with scales larger than one is independent of the fields at scales at most one, so~$\mu_N=e^H\mu_N^0$ is a martingale given~$H$. On a bounded region the factor~$e^{pH}$ is controlled by~$\E e^{p\sup H}<\infty$, which follows from the exponential moments of the~$C^{1,1}$ norm of~$H$ on compact cubes in Lemma~\ref{mfd:lem-infrared}. It therefore suffices to prove~\eqref{eq:mfd-38} for~$\mu_N^0$.

\smallskip

For~$p\geq2$ we expand
\begin{equation*}
	\E \mu_N^0(z+r\cu_0)^p=\int_{(z+r\cu_0)^p}\E\prod_{i=1}^p\prod_{j=0}^N
	e^{\g_{-j}(x_i)-\tau^2}\,dx\,.
\end{equation*}
The field~$\g_{-j}$ has range of dependence~$C3^{-j}$ by Assumption~\ref{a.g1}. If all pairwise distances~$|x_i-x_k|$ exceed that range, then the~$p$ factors are independent and the normalized expectation equals one. For each of the remaining fields, H\"older's inequality gives~$\E\prod_ie^{\g_{-j}(x_i)-\tau^2}\leq\E e^{p\g_0(0)}e^{-p\tau^2}\leq e^{Cp^2\gmcdelta^2}$ by~\eqref{e.coefficient.field.regularity}. The number of fields for which this separation condition fails is at most~$C+\sum_{i<k}\log_3^+(C/|x_i-x_k|)$, so the integrand is bounded by~$C\prod_{i<k}(|x_i-x_k|\wedge1)^{-a}$ with~$a\coloneqq Cp^2\gmcdelta^2$. There are~$p(p-1)/2$ pair factors. H\"older's inequality with exponent~$p(p-1)/2$ and the bound~$\int_{(z+r\cu_0)^2}|x-y|^{-ap(p-1)/2}\,dx\,dy\leq Cr^{2d-ap(p-1)/2}$, valid when~$ap(p-1)/2<d$, give~$\E \mu_N^0(z+r\cu_0)^p\leq C_pr^{dp-ap(p-1)/2}$. Since~$ap(p-1)/2\leq Cp^4\gmcdelta^2$, the condition~$ap(p-1)/2<d$ holds once~$\gmcdelta\leq c_p$, and~\eqref{eq:mfd-38} follows.
\end{proof}

\begin{proposition}[Mass of small balls and the limiting measure]
\label{mfd:prop-chaos-growth}
Let~$\ep\in(0,1)$, let~$p$ be an integer with~$p\ep>d$, and assume that~$\gmcdelta\leq c_p$ and~$C_p\gmcdelta^2<p\ep-d$, with~$c_p$ and~$C_p$ as in Lemma~\ref{mfd:lem-chaos-moments}. Then, almost surely,~$\mu_N$ converges locally weakly to a locally finite measure~$\mu$, and, for every bounded region, there is a random constant~$K_\mu$ with~$\E K_\mu^p<\infty$ such that, for every~$x$ in this region, every~$0<r\leq1$ and every~$N\geq0$,
\begin{equation}
	\mu_N(B_r(x))\leq K_\mu r^{d-\ep}
	\quad\text{and}\quad
	\mu(B_r(x))\leq K_\mu r^{d-\ep}\,.
	\label{eq:mfd-39}
\end{equation}
Almost surely, the measure~$\mu$ has full support and no atoms, and~$\mu(\partial(w+R\cu_0))=0$ for every cube~$w+R\cu_0$.
\end{proposition}

\begin{proof}
For~$z\in2^{-k}\Zd$ such that~$z+2^{-k}\cu_0$ lies in a fixed bounded region, Doob's inequality~\citep[Proposition~7.16]{Kallenberg} and~\eqref{eq:mfd-38} give~$\E\sup_N\mu_N(z+2^{-k}\cu_0)^p\leq(p/(p-1))^p\sup_N\E\mu_N(z+2^{-k}\cu_0)^p \leq C2^{-k(dp-C_p\gmcdelta^2)}$. Let~$C$ be the constant in the covering below and let
\begin{equation*}
	K_\mu^p\coloneqq C^p\sum_{k\geq0}\ \sum_{z\in2^{-k}\Zd,\ z+2^{-k}\cu_0\text{ in the region}}
	2^{kp(d-\ep)}\sup_N\mu_N(z+2^{-k}\cu_0)^p\,.
\end{equation*}
There are~$C2^{dk}$ cubes of side~$2^{-k}$ in the region, so~$\E K_\mu^p\leq C\sum_k2^{k(d+C_p\gmcdelta^2-p\ep)}<\infty$ by the choice of~$p$ and then of~$\gmcdelta$. A ball~$B_r(x)$ with~$2^{-k-1}<r\leq2^{-k}$ is covered by a number of cubes of side~$2^{-k}$ that depends only on~$d$. This gives the first bound in~\eqref{eq:mfd-39} with one random constant for all~$N$,~$x$ and~$r\leq1$.

\smallskip

Since~$p>1$, the bound~\eqref{eq:mfd-38} at this order gives uniform integrability of~$\mu_N(\phi)$ for every~$\phi\in C_c(\Rd)$. The martingale convergence theorem~\citep[Theorem~7.18 and Corollary~7.22]{Kallenberg} gives almost-sure and~$L^p$ convergence against every function in a countable dense subset of~$C_c(\Rd)$. The local mass bound above gives vague relative compactness~\citep[Theorem~A2.3]{Kallenberg}, so these limits determine a single locally finite measure~$\mu$, and~$\mu_N\to\mu$ locally weakly almost surely. The bound for~$\mu$ follows from the portmanteau theorem applied to open balls. The measure~$\mu$ has no atoms by~\eqref{eq:mfd-39}. To prove that~$\mu(\partial(w+R\cu_0))=0$, cover~$\partial(w+R\cu_0)$ by~$Cr^{1-d}$ balls of radius~$r$. This gives~$\mu(\partial(w+R\cu_0))\leq CK_\mu r^{1-\ep}\to0$.

\smallskip

It remains to prove that~$\mu$ has full support. Let~$\mu^0\coloneqq e^{-H}\mu$. Fix an open cube~$w+R\cu_0$ and~$k\geq0$, and write~$\mu^0=e^{\sum_{j=0}^k\g_{-j}-(k+1)\tau^2}\mu^{0,(>k)}$, where~$\mu^{0,(>k)}$ is the limiting measure built in the same way from the fields~$\g_{-j}$ with~$j>k$ only. Since~$\E \mu^{0,(>k)}(\partial(w+R\cu_0))=|\partial(w+R\cu_0)|=0$ and the factor in front of~$\mu^{0,(>k)}$ is a fixed continuous positive function, weak convergence shows that the two sides of this identity give the same mass to~$w+R\cu_0$. Hence the events~$\{\mu^0(w+R\cu_0)=0\}$ and~$\{\mu^{0,(>k)}(w+R\cu_0)=0\}$ coincide up to null sets for every~$k$, so~$\{\mu^0(w+R\cu_0)=0\}$ is a tail event of the independent fields. Its probability is therefore~$0$ or~$1$, and it is~$0$ because~$\E \mu^0(w+R\cu_0)=R^d>0$ by uniform integrability. The positive continuous factor~$e^H$ gives the full support of~$\mu$.
\end{proof}

\subsubsection{Traces on the reversible measures}

Fix~$z\in\mathbb Q^d$ and~$n\in\mathbb Z$ for the rest of this subsection. We choose~$\ep<1$ in~\eqref{eq:mfd-39}, so that the trace inequality~\eqref{eq:mfd-19} applies with~$t=d-\ep>d-1$ to every~$\mu_N$ and to~$\mu$, with the constant~$C\cdot(K_\mu+\mu_N(z+\cu_n))$. We extend functions in~$H^1_0(z+\cu_n)$ by zero to~$\Rd$. For such a function~$v$, the zero extension satisfies~$\|v\|_{H^{\nf34}(\Rd)}\leq C\|v\|_{H^{\nf34}(z+\cu_n)}$, so~\eqref{eq:mfd-1} gives
\begin{equation}
	\|v\|_{H^{\nf34}(\Rd)}^2\leq CQ_N\mathscr E_N^{z+\cu_n}(v)\,,
	\label{eq:mfd-D-zeroext}
\end{equation}
for every~$v\in H^1_0(z+\cu_n)$, with a constant that is bounded along the sequence in the Skorokhod representation. Traces on~$\mu$ are defined by density: by~\eqref{eq:mfd-19}, the map sending a smooth function on~$\overline{z+\cu_n}$ to its class in~$L^2(\mu)$ is bounded from~$H^{\nf12}$ to~$L^2(\mu)$, and it extends by continuity.

\begin{lemma}[Varying-measure convergence]
\label{mfd:lem-varying-trace}
Let~$u_N\in H^1_0(z+\cu_n)$ satisfy~$\sup_N\mathscr E_N^{z+\cu_n}(u_N)<\infty$ along the sequence in the Skorokhod representation and~$u_N\to u$ in~$L^2(z+\cu_n,dx)$. Then, for every bounded continuous function~$f$ on~$\overline{z+\cu_n}$,
\begin{equation}
	\int u_N^2\,d\mu_N\longrightarrow\int u^2\,d\mu
	\quad\text{and}\quad
	\int fu_N\,d\mu_N\longrightarrow\int fu\,d\mu\,.
	\label{eq:mfd-40}
\end{equation}
\end{lemma}

\begin{proof}
By~\eqref{eq:mfd-D-zeroext} the zero extensions are bounded in~$H^{\nf34}(\Rd)$, and the interpolation inequality~$\|w\|_{H^{\nf12}}\leq\|w\|_{L^2}^{\nf{1}{3}}\|w\|_{H^{\nf34}}^{\nf{2}{3}}$ upgrades the convergence in~$L^2$ to strong convergence in~$H^{\nf12}(\Rd)$. We write
\begin{equation*}
	\int u_N^2\,d\mu_N-\int u^2\,d\mu
	=\int(u_N-u)(u_N+u)\,d\mu_N+\Bigl(\int u^2\,d\mu_N-\int u^2\,d\mu\Bigr)\,.
\end{equation*}
The first term tends to zero by Cauchy--Schwarz and~\eqref{eq:mfd-19}. For the second, we approximate~$u$ in~$H^{\nf12}$ by smooth functions~$\phi$. By~\eqref{eq:mfd-19}, the resulting error is small uniformly over all~$\mu_N$ and~$\mu$, and~$\int\phi^2\,d\mu_N\to\int\phi^2\,d\mu$ by weak convergence, since~$\mu(\partial(z+\cu_n))=0$. The proof of the second limit is the same.
\end{proof}

\subsubsection{The killed resolvent as a minimizer}

At finite cutoff the generator is symmetric in~$L^2(\mu_N)$ with form~$\mathscr E_N^{z+\cu_n}$. Since the density of~$\mu_N$ is bounded on~$\overline{z+\cu_n}$ at every finite~$N$, it follows that, for~$\lambda>0$ and bounded~$f$, the killed resolvent~$u_N=R^{z+\cu_n}_{N,\lambda}f$ is the unique minimizer over~$H^1_0(z+\cu_n)$ of~$\mathscr E_N^{z+\cu_n}(u)+\lambda\int u^2\,d\mu_N-2\int fu\,d\mu_N$. It satisfies the sub-Markov bound~$\|u_N\|_\infty\leq\|f\|_\infty/\lambda$, the weak equation~$\mathscr E_N^{z+\cu_n}(u_N,\phi)=\int(f-\lambda u_N)\phi\,d\mu_N$, and, by testing with~$u_N$, the energy bound~$\mathscr E_N^{z+\cu_n}(u_N)\leq\|f\|_\infty^2\mu_N(z+\cu_n)/\lambda$.

\begin{proposition}[Passage of the minimization to the limit]
\label{mfd:prop-speed-resolvent}
Let~$f\in C(\overline{z+\cu_n})$ and~$\lambda>0$, and assume that, along the sequence in the Skorokhod representation,~$\mathscr E_N^{z+\cu_n}$ Mosco-converges to~$\mathscr E^{z+\cu_n}$ and that~$\mathscr E^{z+\cu_n}$ satisfies the coercivity bound~\eqref{eq:mfd-1} with a finite constant. Then the functional
\begin{equation}
	u\longmapsto \mathscr E^{z+\cu_n}(u)+\lambda\int u^2\,d\mu-2\int fu\,d\mu
	\label{eq:mfd-41}
\end{equation}
has a unique minimizer, denoted~$R^{z+\cu_n}_\lambda f$, in the energy Hilbert space~$(D(\mathscr E^{z+\cu_n}),\mathscr E^{z+\cu_n})$, and every~$L^2(dx)$ cluster point of the killed resolvents~$u_N=R^{z+\cu_n}_{N,\lambda}f$ equals~$R^{z+\cu_n}_\lambda f$.
\end{proposition}

\begin{proof}
The energy space is complete because~$\mathscr E^{z+\cu_n}(u)\geq\|G^{z+\cu_n}\|^{-1}\|u\|_{L^2}^2$ by the variational formula~\eqref{eq:mfd-18}. The trace from~$D(\mathscr E^{z+\cu_n})$ to~$L^2(\mu)$ is bounded, since~$\int u^2\,d\mu\leq C(K_\mu+\mu(z+\cu_n))\|u\|_{H^{\nf12}}^2\leq CM\mathscr E^{z+\cu_n}(u)$ with a random constant~$M$, by~\eqref{eq:mfd-19} and the coercivity of~$\mathscr E^{z+\cu_n}$. The functional~\eqref{eq:mfd-41} is therefore continuous, strictly convex and coercive on the energy space, and the Riesz representation theorem gives a unique minimizer~$R^{z+\cu_n}_\lambda f$, which satisfies~$\mathscr E^{z+\cu_n}(R^{z+\cu_n}_\lambda f,v)+\lambda\int R^{z+\cu_n}_\lambda f\,v\,d\mu=\int fv\,d\mu$ for every~$v\in D(\mathscr E^{z+\cu_n})$.

\smallskip

The energies~$\mathscr E_N^{z+\cu_n}(u_N)$ are bounded, so, by~\eqref{eq:mfd-1}, a subsequence converges strongly in~$L^2(dx)$. The liminf inequality of Mosco convergence and Lemma~\ref{mfd:lem-varying-trace} show that the~$\liminf$ of the finite-cutoff functionals at~$u_N$ is at least the value of~\eqref{eq:mfd-41} at the limit. Since~$u_N$ minimizes the finite-cutoff functional, recovery sequences for every~$v\in D(\mathscr E^{z+\cu_n})$, together with~\eqref{eq:mfd-40}, show that this~$\liminf$ is at most the value of~\eqref{eq:mfd-41} at~$v$. Hence the limit minimizes~\eqref{eq:mfd-41}, so it equals~$R^{z+\cu_n}_\lambda f$.
\end{proof}

\subsubsection{Uniform convergence of the killed resolvents}

To obtain uniform convergence we need H\"older bounds on the killed resolvents that are uniform in~$N$. The source in the resolvent equation is the measure~$(f-\lambda u_N)\mu_N$, whose density is not bounded uniformly in~$N$, so we mollify the source before applying the H\"older estimate for bounded sources.

For the next proposition, assume the hypotheses of Proposition~\ref{mfd:prop-speed-resolvent}. Assume also that solutions with zero Dirichlet boundary values and bounded sources satisfy the H\"older estimate~\eqref{eq:mfd-5} with~$\alpha=1/2$ up to the boundary, with the supremum norm of the source on the right side and with a constant whose moments of the finitely many orders fixed at the start of Section~\ref{sec:mfd-uniqueness} are finite. Take~$\ep<1/[8(d+2)]$ in~\eqref{eq:mfd-39}.

\begin{proposition}[Uniform convergence of killed resolvents]
\label{mfd:prop-uniform-resolvent}
Under these assumptions, for every bounded continuous~$f$ and every~$\lambda>0$,
\begin{equation}
	\sup_{\overline{z+\cu_n}}\bigl|R^{z+\cu_n}_{N,\lambda}f-R^{z+\cu_n}_\lambda f\bigr|\longrightarrow0
	\quad\text{in probability}\,,
	\label{eq:mfd-42}
\end{equation}
and the limit vanishes on~$\partial(z+\cu_n)$. There is a random constant~$M_{z+\cu_n}$ with finite moments of these orders such that, for all~$\lambda>0$ and~$f\in C(\overline{z+\cu_n})$,
\begin{equation}
	\|R^{z+\cu_n}_\lambda f\|_{C^{\nf14}(\overline{z+\cu_n})}\leq M_{z+\cu_n}\|f\|_\infty\,.
	\label{eq:mfd-uniform-lambda}
\end{equation}
\end{proposition}

\begin{proof}
Let~$u_N\coloneqq R^{z+\cu_n}_{N,\lambda}f$ and~$\nu_N\coloneqq(f-\lambda u_N)\mu_N$, so that~$|\nu_N|\leq2\|f\|_\infty\mu_N$ and~$\mathscr E_N^{z+\cu_n}(u_N,\phi)=\int\phi\,d\nu_N$. Fix~$0<r\leq1$, let~$s\coloneqq r^{d+2}\leq1$, and let~$\psi_s$ be a mollifier at scale~$s$.

\smallskip

\emph{Step 1.} We bound the mollified source and the error of the mollification. By~\eqref{eq:mfd-39},
\begin{equation*}
	|\nu_N*\psi_s(x)|\leq Cs^{-d}|\nu_N|(B_s(x))\leq CK_\mu s^{-\ep}\|f\|_\infty\,.
\end{equation*}
For the error, we first use that, by Schur's test,~$\|\nu_N\|_{H^{-\nf14}}^2$ is at most~$\|\nu_N\|_{TV}$ times the supremum over~$x$ of the Bessel potential of order~$1/2$ of~$|\nu_N|$ at~$x$. Splitting the Bessel potential into dyadic annuli and using~\eqref{eq:mfd-39} bounds this supremum, since~$\ep<1/2$. Hence~$\|\nu_N\|_{H^{-\nf14}}\leq M\|f\|_\infty$. Then the bound~$\sup_\vartheta(1+|\vartheta|^2)^{-\nf12}\min(1,s|\vartheta|)^2\leq Cs$ gives
\begin{equation*}
	\|\nu_N-\nu_N*\psi_s\|_{H^{-\nf34}}\leq Cs^{\nf12}\|\nu_N\|_{H^{-\nf14}}
	\leq Ms^{\nf12}\|f\|_\infty\,.
\end{equation*}

\smallskip

\emph{Step 2.} We compare~$u_N$ with the solution for the mollified source. Let~$v_{N,s}\in H^1_0(z+\cu_n)$ solve~$\mathscr E_N^{z+\cu_n}(v_{N,s},\phi)=\int\phi\,(\nu_N*\psi_s)\,dx$. Its source is a bounded function of size~$Ms^{-\ep}\|f\|_\infty$, so~\eqref{eq:mfd-5} with~$\alpha=1/2$ gives~$\|v_{N,s}\|_{C^{\nf12}(\overline{z+\cu_n})}\leq Ms^{-\ep}\|f\|_\infty$. The difference of the two equations and~\eqref{eq:mfd-D-zeroext} give
\begin{align*}
	\mathscr E_N^{z+\cu_n}(u_N-v_{N,s})=\langle\nu_N-\nu_N*\psi_s,u_N-v_{N,s}\rangle
	&
	\leq Ms^{\nf12}\|f\|_\infty\|u_N-v_{N,s}\|_{H^{\nf34}}
	\notag \\ &
	\leq Ms^{\nf12}\|f\|_\infty \mathscr E_N^{z+\cu_n}(u_N-v_{N,s})^{\nf12}\,,
\end{align*}
hence~$\|u_N-v_{N,s}\|_{L^2}\leq Ms^{\nf12}\|f\|_\infty$.

\smallskip

\emph{Step 3.} We prove a H\"older bound for~$u_N$. By Steps~1 and~2, the mean-square oscillation of~$u_N$ on~$B_r(x)\cap(z+\cu_n)$ is at most
\begin{equation*}
	M\bigl(r^{\nf12}s^{-\ep}+r^{-\nf{d}{2}}s^{\nf12}\bigr)\|f\|_\infty
	=M\bigl(r^{\nf{1}{2}-(d+2)\ep}+r\bigr)\|f\|_\infty\,,
\end{equation*}
and~$1/2-(d+2)\ep>3/8>1/4$. Using zero extensions across~$\partial(z+\cu_n)$, Campanato's criterion gives~$\|u_N\|_{C^{\nf14}(\overline{z+\cu_n})}\leq M\|f\|_\infty$ and~$u_N=0$ on~$\partial(z+\cu_n)$. Here~$M$ is a product of finitely many random constants with finite moments of these orders, and it is therefore tight.

\smallskip

\emph{Step 4.} We prove convergence in probability. Fix two cutoff subsequences, and apply the Skorokhod representation of Section~\ref{mfd:sec-local-form} to both of them jointly; the variables represented jointly include the fields~$(\g_j)_{j\in\Z}$, the reversible measures, the quadratic forms~$\langle h,G_N^{y+\cu_\ell}h\rangle$ with~$h\in\mathcal D$,~$y\in\mathbb Q^d$ and~$\ell\in\mathbb Z$, and the constants in Steps~1--3. The constants have tight laws. The limit forms of the two sequences coincide by Theorem~\ref{mfd:thm-c1}, and their limiting reversible measures are the same. The Arzel\`a--Ascoli theorem and Proposition~\ref{mfd:prop-speed-resolvent} therefore identify both uniform limits with the unique minimizer of the same functional.

\smallskip

Consequently, every pair of cutoff subsequences has a further pair of subsequences along which the joint laws of the two killed resolvents converge to a law supported on the diagonal of~$C(\overline{z+\cu_n})^2$. Since~$C(\overline{z+\cu_n})$ is Polish, the Gy\"ongy--Krylov criterion~\citep[Lemma~1.1]{GyongyKrylov1996} gives convergence in probability on the original probability space. The joint representation with these quadratic forms and the reversible measure identifies the limit with~$R^{z+\cu_n}_\lambda f$. The constant in Step~3 is independent of~$\lambda$ and~$f$, which gives~\eqref{eq:mfd-uniform-lambda} and the zero boundary values.
\end{proof}

\subsection{Tightness}
\label{s.tightness}

In this subsection we prove Proposition~\ref{tight:prop-tightness} (tightness in expectation and in probability). Section~\ref{mfd:sec-paths} uses both forms of tightness, together with the compact containment estimate~\eqref{tight:eq-containment}, to identify the limiting path law from its resolvents. The early-exit bound~\eqref{tight:eq-fast-exit} is used again in Section~\ref{lim:sec-process-measure} to bound mean exit times from below.

\smallskip

We first bound the probability of an early exit from a cube. After the rescaling introduced below, a cube of side~$3^m$ becomes~$\cu_0$ and time is measured in units of~$3^{2m}/\ahom_{m\wedge L}$. The process started in~$\cu_{-2}$ leaves~$\cu_0$ before time~$t$ with probability at most~$M\sqrt t$, where~$M$ has moments bounded uniformly in the cutoff, the cube size, and the position of the cube.

\smallskip

The proof bounds the Laplace transform of the exit time by one minus the solution of a resolvent problem with a cutoff function as source. A local maximum principle for subsolutions, with the Laplace parameter equal to~$1/t$, then gives the exit-probability bound. Since pointwise ellipticity is not uniform in the cutoff and the reversible measure becomes singular in the limit, we prove this maximum principle by Moser iteration with respect to the reversible measure. The iteration uses the mass bounds, the~$H^{\nf34}$ coercivity bound, and cutoff functions of Lemma~\ref{tight:lem-static}. The energy measures of these cutoff functions satisfy an upper bound with the same power of the ball radius as the upper mass bound for the reversible measure. The exit-probability bound shows that, with high probability, the paths stay in a large cube up to a fixed time; applied at stopping times, it also controls their modulus of continuity.

\smallskip

Throughout this subsection we assume~\ref{a.g1}--\ref{a.g4} and that~$\gmcdelta$ is small enough, depending on the moment orders fixed below. For the diffusion~$X^{(L)}$ of~\eqref{e.intro.cutoff.coefficients}--\eqref{e.intro.cutoff.generator}, let~$Z_t^{(L)}\coloneqq3^{-L}X_{3^{2L}t/\ahom_L}^{(L)}$ for~$t\geq0$. We use rescaled coordinates for~$Z^{(L)}$ and physical coordinates for~$X^{(L)}$. All statements of this subsection also hold for the rescaled diffusion~$3^{-L}X_{3^{2L}t/\ahom_L}$ without cutoff.

\smallskip

For~$L,m\in\N_0$ and~$z\in\Rd$, consider the cube~$z+\cu_m$ in physical coordinates, and let
\begin{equation*}
	j\coloneqq m\wedge L,\qquad
	c\coloneqq\left(\frac{\a_L}{\a_j}\right)(z),\qquad
	d\mu\coloneqq c^{-1}\a_L(z+3^m x)\,dx,\qquad
	A\coloneqq\ahom_j^{-1}c^{-1}\a_L(z+3^m x)\,.
\end{equation*}
If~$m\geq L$, then~$c=1$. In the local coordinates~$x\mapsto z+3^mx$, which rescale~$z+\cu_m$ to~$\cu_0$, and in units of time
\begin{equation}
	\frac{3^{2m}}{\ahom_{m\wedge L}}\,,
	\label{tight:eq-local-time}
\end{equation}
the generator of~$X^{(L)}$ becomes~$(d\mu/dx)^{-1}\nabla\cdot(A\nabla)$. The generator is symmetric on~$L^2(\mu)$, with Dirichlet form and local energy
\begin{equation*}
	\mathscr E_A(v,w)\coloneqq\int A\nabla v\cdot\nabla w
	\qquad\text{and}\qquad
	\mathscr E_A(v;V)\coloneqq\int_V A|\nabla v|^2\,.
\end{equation*}
The indexed quantity~$\mathcal E_{s,\infty,q}$ of Definition~\ref{d.mathcal.E} denotes the coarse-graining error. The factor~$c$ normalizes the product of the factors contributed by~$\g_i$, with~$j<i\leq L$, to have value one at~$z$. By~\eqref{e.sensitivity.field}, the oscillation of the sum of these fields on the cube has uniformly bounded exponential moments, and therefore the moment bounds for the random constants below are uniform in~$L$,~$m$ and~$z$. For the coefficient~$\a$ without cutoff we use instead~$j\coloneqq m$ and~$c\coloneqq\a(z)/\a_m(z)$, so that the density of~$\mu$ is~$c^{-1}\a(z+3^mx)$, and the unit of time is~$3^{2m}/\ahom_m$. This density is then~$\a_m(z+3^mx)$ multiplied by~$\exp(\sum_{i>m}(\g_i(z+3^mx)-\g_i(z)))$, and this exponential factor satisfies the same oscillation bound.

\smallskip

Fix a bounded reference cube, a moment order~$q\in[1,\infty)$ and finitely many concentric pairs~$x_0+R_0\cu_0\Subset x_0+R_1\cu_0$ inside the reference cube. For each pair and each~$s\in[0,1]$, let~$x_0+R\cu_0$ be the concentric cube with side length~$R=(1-s)R_0+sR_1$. The nested families in the next lemma consist of these cubes with~$s=k2^{-n}$, for~$n\geq0$ and~$0\leq k\leq2^n$.

\begin{lemma}[Local mass, coercivity and cutoff-function estimates]
\label{tight:lem-static}
There exists a threshold~$\gmcdelta_0>0$, which may depend on~$q$, such that, if~$0<\gmcdelta\leq\gmcdelta_0$, then there is a deterministic~$b_2<\infty$, independent of~$q$, such that, for every~$L,m\in\N_0$ and~$z\in\Rd$, there is a random constant~$M=M_{L,m,z}\geq1$ satisfying
\begin{equation*}
	\sup_{L,m,z}\E\bigl[M_{L,m,z}^q\bigr]<\infty
\end{equation*}
such that, after rescaling~$z+\cu_m$ to~$\cu_0$ as in~\eqref{tight:eq-local-time}, the following statements hold.
\begin{enumerate}[label=(\roman*)]
\item \underline{\emph{Mass bounds}}: for every ball~$B_r(x)$ contained in the reference cube with~$0<r\leq1$,
\begin{equation}
	M^{-1}r^{d+\nf12}\leq\mu(B_r(x))\leq Mr^{d-\nf12}\,.
	\label{tight:eq-mass}
\end{equation}

\item \underline{\emph{Coercivity in~$H^{\nf34}$}}: for each fixed interior cube~$x_0+R\cu_0$ of the family with~$n=0$ and every~$v\in H^1(x_0+R\cu_0)$,
\begin{equation}
	\|v\|_{H^{\nf34}(x_0+R\cu_0)}^2
	\leq M\bigl(\mathscr E_A(v;x_0+R\cu_0)+\|v\|_{L^2(x_0+R\cu_0,dx)}^2\bigr)\,,
	\label{tight:eq-coercivity}
\end{equation}
and the last term can be omitted if~$v\in H_0^1(x_0+R\cu_0)$. For~$s=k2^{-n}$, the same assertion holds with~$M2^{b_2n}$ in place of~$M$.

\item \underline{\emph{Harmonic cutoff functions}}: between consecutive cubes~$x_0+R\cu_0\Subset x_0+R'\cu_0$ of the family, with~$s=k2^{-n}$ and~$s=(k+1)2^{-n}$ respectively, with gap~$a\asymp2^{-n}$ there is a cutoff function~$\chi$, equal to one on the smaller cube, vanishing near the boundary of the larger cube, and satisfying~$0\leq\chi\leq1$, such that, for every ball~$B_r(x)$ with~$0<r\leq1$,
\begin{equation}
	\Gamma_\chi(B_r(x))
	\coloneqq\int_{B_r(x)}A|\nabla\chi|^2
	\leq Ma^{-b_2}r^{d-\nf12}\,.
	\label{tight:eq-cutoff-measure}
\end{equation}
\end{enumerate}
\end{lemma}

When needed, we replace~$M$ by a fixed power of~$M$; this is possible if the lemma is applied from the outset with a correspondingly larger moment order~$q$. For a finite collection of cubes we take the maximum of their constants; the resulting moment bounds remain uniform in~$L$,~$m$ and~$z$.

\begin{proof}
The three estimates are local versions of estimates proved in Sections~\ref{sec:mfd-uniqueness} and~\ref{mfd:sec-speed} for the rescaled coefficient~$A_N$ of~\eqref{eq:mfd-normalization}: Lemma~\ref{mfd:lem-chaos-moments} bounds the positive moments of the mass of small cubes under the multiplicative chaos, Lemma~\ref{mfd:lem-coercivity} bounds the~$H^{\nf34}$ norm by the energy, Proposition~\ref{mfd:prop-growth} bounds the energy measure of balls at every radius for solutions of Dirichlet problems, and Lemma~\ref{mfd:lem-cutoffs} constructs cutoff functions as harmonic interpolants. We give the additional arguments needed to pass to the rescaling of~$z+\cu_m$ to~$\cu_0$ and to obtain one random constant for all the nested cubes.

\smallskip

\emph{Step 1.} We prove the mass bounds~\eqref{tight:eq-mass}. For~$m\leq L$, reindexing the fields in reverse order identifies~$\int_{\cu_0}\a_m(z+3^mx)\,dx$ in law with~$\mu_m^0(\cu_0)$ from the paragraph after~\eqref{e.intro.rescaled.generator}; the fields above scale~$m$ contribute the oscillation factor described before the lemma, whose positive and negative moments are uniformly bounded. Lemma~\ref{mfd:lem-chaos-moments}, applied to the martingale~$N\mapsto\mu_N^0(\cu_0)$, bounds the positive moments of~$\mu_N^0(\cu_0)$ uniformly in~$N$, including the field~$\g_0$ at scale one. The bound on the second moment and~$\E \mu_N^0(\cu_0)=1$ give~$\inf_N\P(\mu_N^0(\cu_0)\geq\nf12)>0$, by the Cauchy--Schwarz inequality. Fix~$i\in\N$ and choose two sufficiently separated triadic subcubes of~$\cu_0$ of side~$3^{-i}$. Since the fields are independent and have a finite range of dependence, for every~$N\geq i$ we have
\begin{equation*}
	\mu_N^0(\cu_0)\geq 3^{-di}Y\bigl(\mu_{N-i}^{0,(1)}(\cu_0)+\mu_{N-i}^{0,(2)}(\cu_0)\bigr)\,,
\end{equation*}
where~$\mu_{N-i}^{0,(1)}$ and~$\mu_{N-i}^{0,(2)}$ are independent copies of~$\mu_{N-i}^0$, independent of the random variable~$Y>0$. Here~$Y$ is the infimum over~$\cu_0$ of the part of the density of~$\mu_N^0$ contributed by the~$i$ coarsest fields. By Assumption~\ref{a.g2}, the variable~$-\log Y$ has a Gaussian upper tail. Consequently, the function~$\mathsf h(s)\coloneqq\sup_N\P(\mu_N^0(\cu_0)<e^{-s})$ satisfies~$\mathsf h(s)\leq Ce^{-cs^2}+\mathsf h(3s/4)^2$ for all sufficiently large~$s$, where the finitely many~$N<i$ are absorbed into the first term. The lower bound on~$\P(\mu_N^0(\cu_0)\geq\nf12)$ gives~$\mathsf h(s)\leq\theta<1$ for large~$s$. Letting~$s\to\infty$ in the recursion gives~$\mathsf h(s)\to0$. Iterating the recursion on the intervals~$[(\nf43)^ns_0,(\nf43)^{n+1}s_0]$ then gives~$\mathsf h(s)\leq C\exp(-cs^{\nf32})$: since~$2(\nf34)^{\nf32}>1$, squaring the bound at~$\nf{3s}4$ gives a stronger bound at~$s$, and the Gaussian term in the recursion is of smaller order. Hence, for every~$p<\infty$,
\begin{equation}
	\sup_N\E\bigl[\mu_N^0(\cu_0)^{-p}\bigr]<\infty\,.
	\label{tight:eq-negative-mass}
\end{equation}

\smallskip

For~$m>L$, we have~$c=1$. Scaling by~$3^L$ and using stationarity shows that~$\mu(\cu_0)$ has the law of~$3^{-d(m-L)}\mu_L^0(\cu_{m-L})$. Partition~$\cu_{m-L}$ into its~$3^{d(m-L)}$ unit cubes. This normalized mass is the average of their masses, each with the law of~$\mu_L^0(\cu_0)$. For every~$p\geq1$, convexity of~$v\mapsto v^p$ and~$v\mapsto v^{-p}$ on~$(0,\infty)$ bounds both moments of the average by those of one unit mass. Thus the positive and negative moment bounds for~$\mu(\cu_0)$ also hold above the cutoff. Independence of the unit masses is not required.

\smallskip

We now pass to the triadic subcubes of the reference cube. Let~$x_1+\cu_{-n}$ be such a subcube, of side~$3^{-n}$, and split the density of~$\mu$ into the factors from fields of scale below and above the physical side of~$x_1+\cu_{-n}$. The product of the factors from fields of scale between the physical side of~$x_1+\cu_{-n}$ and~$3^j$, and its reciprocal, have moments of order~$p$ bounded by~$C_p3^{C_p\gmcdelta^2n}$, and the oscillation of the fields above scale~$j$ has a uniformly bounded exponential moment. Combining the positive moments of~$\mu_N^0(\cu_0)$ with~\eqref{tight:eq-negative-mass} therefore gives, for every~$p<\infty$,
\begin{equation*}
	\E\left[
	\left(\frac{\mu(x_1+\cu_{-n})}{3^{-nd}}\right)^p
	+\left(\frac{\mu(x_1+\cu_{-n})}{3^{-nd}}\right)^{-p}\right]
	\leq C_p3^{C_p\gmcdelta^2n}\,.
\end{equation*}
The reference cube contains~$O(3^{dn})$ triadic subcubes of side~$3^{-n}$. Choose~$p>q+2d$, and then decrease~$\gmcdelta$ so that~$p/2>d+C_p\gmcdelta^2$. Markov's inequality and a summation over~$n$ prove~\eqref{tight:eq-mass} for triadic cubes, with the asserted moments, and a finite collection of translated triadic grids gives the bounds for balls.

\smallskip

\emph{Step 2.} We prove the~$H^{\nf34}$ coercivity bound~\eqref{tight:eq-coercivity}. If~$m\leq L$, it is the estimate of Lemma~\ref{mfd:lem-coercivity}, applied with~$N=m$ and with the value at~$z$ subtracted from each field of scale larger than~$3^m$; the proof is uniform in the number of these fields, by~\eqref{e.sensitivity.field}.

\smallskip

For~$m\geq L$, fix a moment order~$p\geq\max\{q,64d\}$, let~$\xi\coloneqq2p$, and choose the disorder threshold afterwards. Write~$C_{\rm resp}$ for the constant of Proposition~\ref{p.coarse.grained.bound} and~$c_1$ for the admissibility constant of Lemma~\ref{l.multiscale.response.large.cubes}. Let~$\delta_1\coloneqq(1+C_{\rm resp}\xi\log(2+\xi)+\nf{32\xi}{c_1})\gmcdelta^2$. For small enough~$\gmcdelta$, we have~$\gmcdelta^2\leq\delta_1<1$ and~$\xi$ is in the range of Proposition~\ref{p.coarse.grained.bound}. That proposition bounds the coarse-grained defect at every scale by~$\delta_1$, and hence supplies the induction hypothesis of Lemma~\ref{l.multiscale.response.large.cubes} with~$m_0=k$ at each application on a physical cube of scale~$k\geq L$. Its lower moment requirement holds because~$\xi\geq128d=32\cdot4d$, and its upper requirement holds because~$\xi\leq c_1\gmcdelta^{-2}\delta_1/32$. Applying it with~$(\nf1{16},1)$ and~$(\nf1{32},2)$ gives, uniformly in~$k\geq L$ and every deterministic center~$y$,
\begin{equation}
	\|\mathcal E_{\nf1{16},\infty,1}(y+\cu_k;\a_L,\ahom_L)\|_{L^{2p}}
	+\|\mathcal E_{\nf1{32},\infty,2}(y+\cu_k;\a_L,\ahom_L)\|_{L^{2p}}
	\leq C_p\gmcdelta\,.
	\label{eq:mfd-static-above-cutoff}
\end{equation}
Stationarity gives the same constants after translation. By~\eqref{e.bound.Lambdas.by.Es}, this bounds~$\ahom_L\lambda_{\nf1{16},1}^{-1}(y+\cu_k;\a_L)$ in~$L^p$ uniformly in~$L,k,y$. The squared coarse Poincar\'e inequality~\eqref{e.besov.grad.poincare} therefore uses this random factor times the energy divided by~$\ahom_L$, which is precisely the normalization of~$A$ in~\eqref{tight:eq-local-time}. Sobolev embedding gives~\eqref{tight:eq-coercivity}; the smaller physical cubes use the case~$m\leq L$ treated above. For~$v\in H^1_0(x_0+R\cu_0)$, use~\eqref{e.CG.Poincare.trace.zero}.

\smallskip

For triadic subcubes of side~$3^{-n}$ in the local coordinates, the fields of scale between the physical side of the subcube and~$3^j$ contribute scalar factors. Together with~\eqref{e.annealed.ordering}, they multiply the constant by a random factor with~$L^p$ norm at most~$C_p3^{C_p\gmcdelta^2n}$. For subcubes of physical side less than one, the ellipticity bounds on unit cubes give a fixed power of~$3^n$ instead, with a random factor whose moments are bounded by Lemma~\ref{mfd:lem-extremes}. Taking the maximum over the~$O(3^{dn})$ triadic cubes of side~$3^{-n}$ and summing over~$n$ with a sufficiently large geometric weight gives one constant~$M$ for all these meshes. The same summation over the~$O(2^n)$ cubes of the nested families gives~$M2^{b_2n}$ after enlarging~$b_2$. The comparison~\eqref{e.lambda.stability} of translated cubes allows fixed enlargements and changes of grid.

\smallskip

\emph{Step 3.} We construct the cutoff functions of~\eqref{tight:eq-cutoff-measure}. We use harmonic interpolants as in Lemma~\ref{mfd:lem-cutoffs}. Choose a smooth cutoff function~$f$ that equals one on a neighborhood of the smaller cube and vanishes on a neighborhood of the boundary of the larger cube, and a mesh of side~$r_0\asymp a$ such that the cells meeting the transition region of~$f$ lie between the two prescribed cubes. Let~$\chi$ be the function which, on each cell, is~$A$-harmonic with boundary values~$f$. The maximum principle gives~$0\leq\chi\leq1$, and~$\chi=f$ on cells where~$f$ is constant.

\smallskip

On a physical cell of side~$3^k$ with~$k\geq0$, apply Proposition~\ref{p.cutoff.Holder.regularity} at exponent~$7/8$ and outer scale~$k$. Above the cutoff, the whole-cell Dirichlet energy is bounded by~\eqref{e.cg.RHS} and~\eqref{eq:mfd-static-above-cutoff}, with exponent~$\nf1{16}$. The tail coefficient is~$\ahom_{k\wedge L}$, since~$k\wedge L\leq k$. The gradient estimate then gives energy growth with exponent~$d-1/4$ down to physical radius one; its random factor consists of the squared whole-cell energy constant and~$3^{d\mathcal L}$, where~$\mathcal L$ is the minimal scale of Proposition~\ref{p.cutoff.Holder.regularity}. To bound this factor at a requested order~$p$, obtain the energy constant at order~$4p$ and~$3^{d\mathcal L}$ at order~$2p$ before choosing~$\gmcdelta$.

\smallskip

Below physical radius one, use the small-oscillation estimate and the odd reflections of Steps~4--5 of Proposition~\ref{mfd:prop-growth}. Join it to the preceding bound at a radius proportional to the inverse log-gradient bound. The resulting constant has the form of the preceding random factor plus~$C3^{-k/4}$ times a polynomial of degree at most~$d+6$ in that factor, the bounds on~$\|A\|_{L^\infty}$ and~$\|A^{-1}\|_{L^\infty}$ and the log-gradient bound. For the requested order~$p$, take those bounds at order~$p(d+6)$ with geometric factor~$3^{k/(8(d+6))}$, and take the energy and stopping moments at the enlarged orders just described. The factor~$3^{-k/4}$ absorbs~$3^{k/8}$. H\"older's inequality thus gives a cell constant bounded in~$L^p$ uniformly in~$L$ and~$k$, including~$k>L$. This proves the cell estimate at all radii, including faces and corners. For~$k\leq L$, subtract the values of the remaining fields at the cell center as in Proposition~\ref{mfd:prop-growth}; their oscillation factors and reciprocals have the same uniform moments.

\smallskip

It remains to express the cell estimates after rescaling~$z+\cu_m$ to~$\cu_0$ as in~\eqref{tight:eq-local-time}. Let~$x_1+r_0\cu_0$ be a cell of side~$r_0\coloneqq3^{-n}$ with~$n\leq m$ and physical center~$z+3^mx_1$, and let~$k\coloneqq m-n$. In the coordinates of the cell, the coefficient~$A$ is the normalized coefficient of the cell multiplied by two factors: the exponential~$\exp(\sum_{j<i\leq L}(\g_i(z+3^mx_1)-\g_i(z)))$, which is controlled by the oscillation bound for these fields, and the scalar~$(\ahom_{k\wedge L}/\ahom_j)(\a_j(z+3^mx_1)/\a_{k\wedge L}(z+3^mx_1))$. By~\eqref{e.annealed.ordering} and~\eqref{e.sensitivity.field}, this scalar and its reciprocal have~$L^p$ norms at most~$C_p3^{C_p\gmcdelta^2n}$. After subtracting a constant from~$f$, its rescaled~$C^2$ norm on the cell is at most~$Cr_0\|\nabla f\|_{C^1}$, where~$\|\nabla f\|_{C^1}\coloneqq\|\nabla f\|_\infty+\|D^2f\|_\infty$. The estimate on the cell therefore gives a random constant~$M_{x_1+r_0\cu_0}$ with~$\|M_{x_1+r_0\cu_0}\|_{L^p}\leq C_pr_0^{-C_p\gmcdelta^2}$ such that, for every~$0<r\leq r_0$,
\begin{equation*}
	\Gamma_{\chi}(B_r(x)\cap(x_1+r_0\cu_0))
	\leq M_{x_1+r_0\cu_0}\|\nabla f\|_{C^1}^2r_0^{\nf12}r^{d-\nf12}\,.
\end{equation*}
At each mesh size there are~$O(r_0^{-d})$ cells. We take the maximum of the constants~$M_{x_1+r_0\cu_0}$ over these cells and sum over the meshes with a sufficiently large geometric weight. This proves~\eqref{tight:eq-cutoff-measure} when the cells have physical side at least one, or equivalently~$r_0\geq3^{-m}$.

\smallskip

We give the microscopic estimate when~$r_0<3^{-m}$. On the fixed reference cube, which contains the cells, let
\begin{equation*}
	D_m\coloneqq3^{-m}\|\nabla\log A\|_{L^\infty}\,,
	\qquad \Upsilon_m\coloneqq(1+D_m)^d\|A\|_{L^\infty}\,.
\end{equation*}
The Dirichlet principle on each cell~$x_1+r_0\cu_0$ gives~$\int_{x_1+r_0\cu_0}A|\nabla\chi|^2\leq\|A\|_\infty\|\nabla f\|_\infty^2r_0^d$. The difference~$\chi-f$ belongs to~$H^1_0(x_1+r_0\cu_0)$; its energy is at most four times this bound and~$-\nabla\cdot(A\nabla(\chi-f))=\nabla\cdot(A\nabla f)$. Choose~$r_1\coloneqq c_0r_0/(1+D_m)$ with a sufficiently small dimensional~$c_0$. Since~$r_0\|\nabla\log A\|_\infty\leq D_m$, the oscillation of~$\log A$ on~$B_{2r_1}$ is at most~$4c_0$. Near the boundary of the cell, reflect~$\chi-f$ and its source oddly and~$A$ evenly across every face met by the ball. The reflected weak equation has no additional face source. This includes edges and corners, with at most~$2^d$ reflected copies.

\smallskip

After dividing the coefficient by its value at the center, choose~$c_0$ so that the small-oscillation Meyers estimate~\eqref{eq:mfd-micro} applies with exponent~$4d$. The source term is bounded by~$C\|\nabla f\|_{C^1}$ after multiplication by~$r_1$. H\"older's inequality and the whole-cell energy bound give, for~$r\leq r_1$,
\begin{equation*}
	\int_{B_r\cap(x_1+r_0\cu_0)}A|\nabla(\chi-f)|^2
	\leq C\|A\|_\infty\|\nabla f\|_{C^1}^2
	\bigl(r_0^dr_1^{-d+\nf12}+r_1^{\nf12}\bigr)r^{d-\nf12}
	\leq C\Upsilon_m\|\nabla f\|_{C^1}^2r^{d-\nf12}\,.
\end{equation*}
The coefficient value used in the normalization cancels against the energy; no exponential of~$D_m$ occurs. Add the energy of~$f$ to obtain the same estimate for~$\chi$. For~$r_1<r\leq r_0$, use the whole-cell energy and~$r_0/r\leq C(1+D_m)$. A ball of radius~$r\leq r_0$ meets only a dimensional number of cells; for~$r>r_0$, summing their energies gives~$C\|A\|_\infty\|\nabla f\|_\infty^2r^d$. Since~$r\leq1$, we have proved, for all centers and radii,
\begin{equation}
	\Gamma_\chi(B_r(x))
	\leq C\Upsilon_m\|\nabla f\|_{C^1}^2r^{d-\nf12}\,.
	\label{eq:mfd-static-microscopic}
\end{equation}

\smallskip

The proof of Lemma~\ref{mfd:lem-extremes} also applies to this local normalization. The physical log-gradient is a sum of fields at scales~$3^i$, weighted by~$3^{-i}$; the maximum over the~$O(3^{d(m-i)})$ unit cells at that scale gives~$\|D_m\|_{L^p}\leq C_p\sqrt{1+m}$. For~$j=m\wedge L$, the density is the product of the fields through~$j$ times the coarse oscillation factor. The deterministic factor~$\ahom_j^{-1}e^{-(j+1)\tau^2}$ is at most one. The maximum of the fine sum over the physical cube, estimated as in Step~2 of that lemma, gives
\begin{equation*}
	\bigl\|\|A\|_{L^\infty}\bigr\|_{L^p}
	\leq C_p\exp\bigl((C_d\gmcdelta+C_p\gmcdelta^2)m\bigr)\,.
\end{equation*}
These bounds are uniform in~$L,m,z$: when~$m>L$, there are only~$L+1\leq m+1$ fields, and when~$m\leq L$, the remaining oscillation factor has uniformly bounded exponential moments. For the uncut coefficient take~$j=m$ and use the same bound on the infinite coarse tail. Fix~$q\geq1$ first, use H\"older's inequality at the required larger orders, and then decrease~$\gmcdelta$ so that~$C_d\gmcdelta+C_q\gmcdelta^2<\tfrac12\log3$. It follows that
\begin{equation*}
	\|\Upsilon_m\|_{L^q}\leq C_q(1+m)^{d/2}
	e^{(C_d\gmcdelta+C_q\gmcdelta^2)m}\quad\text{and}\quad
	\sup_{L,m,z}\E\bigl[(3^{-m}\Upsilon_m)^q\bigr]<\infty\,.
\end{equation*}
The factor~$\Upsilon_m$ is independent of the mesh. We have~$\Upsilon_m\leq(3^{-m}\Upsilon_m)r_0^{-1}$, and~$\|\nabla f\|_{C^1}\leq Ca^{-2}$ with~$r_0\asymp a$. Thus~\eqref{eq:mfd-static-microscopic} implies~\eqref{tight:eq-cutoff-measure} in the case~$r_0<3^{-m}$ with~$b_2\geq5$, independently of~$q$. Finally, the estimates for~$\a$ follow also by locally smooth convergence of the normalized coarse tails and Fatou's lemma, as in Proposition~\ref{mfd:prop-growth}.
\end{proof}

The regularity estimates of Section~\ref{s.Holder.regularity} concern solutions, and the coefficient~$A$ has no ellipticity bounds uniform in~$L$, so we derive the local maximum principle for subsolutions from Lemma~\ref{tight:lem-static} alone.

\begin{lemma}[Local bound for subsolutions]
\label{tight:lem-subharmonic}
After rescaling a cube~$z+\cu_m$ to~$\cu_0$ as in~\eqref{tight:eq-local-time}, let~$x_0+R_0\cu_0\Subset x_1+R_1\cu_0$ be fixed cubes. The estimates of Lemma~\ref{tight:lem-static} imply that there is a random constant~$M$ such that, for every bounded nonnegative weak subsolution~$w$ of~$\nabla\cdot(A\nabla w)\geq0$ in~$x_1+R_1\cu_0$,
\begin{equation}
	\operatorname*{ess\,sup}_{x_0+R_0\cu_0}w\leq M\left(\int_{x_1+R_1\cu_0} w^2\,d\mu\right)^{\nf12}\,.
	\label{tight:eq-mean-value}
\end{equation}
For every fixed~$q<\infty$, there exists~$\gmcdelta_0>0$ such that, if~$0<\gmcdelta\leq\gmcdelta_0$, then~$\E[M^q]$ is bounded uniformly in~$L$,~$m$ and~$z$.
\end{lemma}

Since~$M$ does not depend on~$w$, the lemma applies with the same constant to every bounded nonnegative solution of~$\nabla\cdot(A\nabla w)=\lambda w\,d\mu/dx$, for every~$\lambda\geq0$. We use this uniformity in~$\lambda$ in Step~1 of the proof of Proposition~\ref{tight:prop-tightness}.

\begin{proof}
By a finite covering it suffices to treat a fixed concentric pair~$x_0+R_0\cu_0\Subset x_0+R_1\cu_0$. In both iterations of Step~4 we use the nested family of Lemma~\ref{tight:lem-static} associated with this pair. We apply Moser iteration with respect to~$\mu$, in which trace inequalities for measures satisfying the upper bound in~\eqref{tight:eq-mass} replace the Sobolev inequality. Step~1 proves the trace inequalities, Steps~2 and~3 use them to bound the~$L^2(dx)$ norm and the integral against the energy measure of a cutoff function by the energy and the~$L^2(\mu)$ norm, and Step~4 iterates the resulting reverse H\"older inequalities for~$w$ on nested cubes.

\smallskip

\emph{Step 1.} We prove two trace inequalities. Let~$\nu$ be a measure satisfying~$\nu(B_r(x))\leq Mr^{d-\nf12}$ for every ball~$B_r(x)$ with~$0<r\leq1$, and let~$p_0\geq2$. We show that, if~$s>\nf14$ and~$t>\nf d2-\nf{(d-\nf12)}{p_0}$, then
\begin{equation}
	\|v\|_{L^2(\nu)}^2\leq CM\|v\|_{H^s}^2
	\quad\text{and}\quad
	\|v\|_{L^{p_0}(\nu)}^2\leq CM\|v\|_{H^t}^2\,.
	\label{tight:eq-trace}
\end{equation}
For a Littlewood--Paley decomposition of~$v$, each block~$v_i$ at frequency~$2^i$ satisfies~$\|v_i\|_{L^2(\nu)}^2\leq CM2^{\nf{i}{2}}\|v_i\|_2^2$. To see this, write~$v_i$ as its convolution with a rapidly decaying reproducing kernel, apply the Cauchy--Schwarz inequality, and integrate against~$\nu$ using the bound on the mass of balls. Interpolation with Bernstein's inequality for that block gives
\begin{equation*}
	\|v_i\|_{L^{p_0}(\nu)}
	\leq CM^{\nf1{p_0}}2^{i(\frac d2-\frac{d-\frac12}{p_0})}\|v_i\|_2\,.
\end{equation*}
Summing over~$i$ with the Cauchy--Schwarz inequality, after extending~$v$ from a cube to the whole space, gives~\eqref{tight:eq-trace} on the fixed cubes. Below we take~$s=\nf12$,~$t=\nf34$ and~$p_0=2+\nf1d$, for which~$\nf d2-\nf{(d-\nf12)}{p_0}=\nf d{(2d+1)}<\nf34$; Section~\ref{lim:sec-process-measure} uses other exponents.

\smallskip

\emph{Step 2.} We bound the~$L^2(dx)$ norm of~$v$ by the energy of~$v$ and its~$L^2(\mu)$ norm. Partition a cube into cubes~$y+r_0\cu_0$ of side~$r_0$, and take the average~$(v)_{y+r_0\cu_0}$ of~$v$ over~$y+r_0\cu_0$ with respect to Lebesgue measure. Apply the first trace inequality, with~$\nf34$ in place of~$\nf12$, to~$v-(v)_{y+r_0\cu_0}$. The fractional Poincar\'e inequality and the lower mass bound give~$\int_{y+r_0\cu_0}|v-(v)_{y+r_0\cu_0}|^2\,d\mu\leq CMr_0[v]_{H^{\nf34}(y+r_0\cu_0)}^2$ and
\begin{equation*}
	\|v\|_{L^2(dx)}^2
	\leq CMr_0^{-\nf12}\|v\|_{L^2(\mu)}^2
	+CM^2r_0^{\nf12}\|v\|_{H^{\nf34}}^2\,.
\end{equation*}
We apply~\eqref{tight:eq-coercivity} to the last term and choose~$r_0$ small enough to absorb the resulting multiple of~$\|v\|_{L^2(dx)}^2$ into the left side. For every~$0<\ep\leq1$ this yields
\begin{equation}
	\|v\|_{L^2(dx)}^2
	\leq\ep\mathscr E_A(v)
	+CM^C\ep^{-C}\|v\|_{L^2(\mu)}^2\,.
	\label{tight:eq-two-norms}
\end{equation}

\smallskip

\emph{Step 3.} We bound the integral of~$v^2$ against the energy measure of a cutoff function by the energy of~$v$ and its~$L^2(\mu)$ norm. Let~$\chi$ be a cutoff function across a gap~$a$ between nested cubes~$x_0+R\cu_0\Subset x_0+R'\cu_0$. We apply the first inequality in~\eqref{tight:eq-trace} to~$\nu=\Gamma_\chi$, which satisfies its hypothesis with~$Ma^{-b_2}$ in place of~$M$ by~\eqref{tight:eq-cutoff-measure}, interpolate~$H^{\nf12}$ between~$H^{\nf34}$ and~$L^2(dx)$, and then use~\eqref{tight:eq-coercivity} and~\eqref{tight:eq-two-norms}. For every~$0<\ep\leq1$ this gives
\begin{equation}
	\int v^2\,d\Gamma_\chi
	\leq\ep\mathscr E_A(v;x_0+R'\cu_0)
	+CM^C\ep^{-b_2}a^{-b_2}\int_{x_0+R'\cu_0} v^2\,d\mu\,.
	\label{tight:eq-cutoff-form}
\end{equation}

\smallskip

\emph{Step 4.} We iterate. For~$p\geq2$, test the subsolution inequality with~$\chi^2w^{p-1}$. The Cauchy--Schwarz inequality and~\eqref{tight:eq-cutoff-form} applied to~$w^{\nf p2}$, with~$\ep$ a small multiple of~$p^{-1}$, give, for nested cubes~$x_0+R\cu_0\Subset x_0+R'\cu_0$,
\begin{equation*}
	\mathscr E_A(w^{\nf p2};x_0+R\cu_0)
	\leq2^{-b_2-2}\mathscr E_A(w^{\nf p2};x_0+R'\cu_0)
	+CM^Cp^{b_2}(R'-R)^{-b_2}\int_{x_0+R'\cu_0}w^p\,d\mu\,.
\end{equation*}
Iterate this inequality over sides increasing from~$R$ to~$R'$ with geometrically decreasing gaps. Hence
\begin{equation*}
	\mathscr E_A(w^{\nf p2};x_0+R\cu_0)
	\leq CM^Cp^{b_2}(R'-R)^{-b_2}\int_{x_0+R'\cu_0}w^p\,d\mu\,.
\end{equation*}
The second trace inequality, the~$H^{\nf34}$ coercivity bound~\eqref{tight:eq-coercivity} and~\eqref{tight:eq-two-norms} give, since~$p_0/2>1$,
\begin{equation*}
	\|w\|_{L^{p_0p/2}(x_0+R\cu_0,\mu)}
	\leq\bigl(CM^Cp^{b_2}(R'-R)^{-b_2}\bigr)^{\nf1p}
	\|w\|_{L^p(x_0+R'\cu_0,\mu)}\,.
\end{equation*}
We iterate with~$p=2(p_0/2)^i$ for~$i=0,1,2,\ldots$ and geometrically decreasing gaps. The product of the resulting constants converges, which proves~\eqref{tight:eq-mean-value}. For continuous functions such as~$1-u$ in the proof of Proposition~\ref{tight:prop-tightness}, the full support of~$\mu$ identifies the essential supremum with the pointwise supremum.
\end{proof}

Let~$K^\omega_{L,x}$ be the quenched law of~$Z^{(L)}$, or of~$3^{-L}X_{3^{2L}t/\ahom_L}$, started at~$x$. Equip~$C([0,\infty),\Rd)$ with the topology of locally uniform convergence.

\begin{proposition}[Tightness]
\label{tight:prop-tightness}
Almost surely these processes are conservative. For every compact set~$F\subset\Rd$, the following statements hold.
\begin{enumerate}[label=(\roman*)]
\item \underline{\emph{Tightness in expectation}}: for every~$\ep>0$, there is a compact set~$\mathcal K\subset C([0,\infty),\Rd)$ such that
\begin{equation}
	\sup_L\E\left[\sup_{x\in F}K^\omega_{L,x}(\mathcal K^c)\right]
	\leq\ep\,.
	\label{tight:eq-uniform-laws}
\end{equation}
In particular, the annealed laws~$\E[K^\omega_{L,x}]$ with~$L\in\N_0$ and~$x\in F$ are tight.

\item \underline{\emph{Tightness in probability}}: for every~$\ep>0$, there is a compact set~$\mathcal M_\ep$ of probability measures on~$C([0,\infty),\Rd)$ such that
\begin{equation}
	\sup_L\P\bigl[K^\omega_{L,x}\notin\mathcal M_\ep\ \text{for some}\ x\in F\bigr]
	\leq\ep\,.
	\label{tight:eq-random-laws}
\end{equation}
\end{enumerate}
\end{proposition}

\begin{proof}
The proof proceeds in four steps. Step~1 bounds the probability that the process, started in the concentric cube of one ninth the side, leaves the cube before a short time~$t$, measured after the rescaling of space and time in~\eqref{tight:eq-local-time}, by a random constant times~$\sqrt t$, with moments of the constant bounded uniformly in the cutoff and in the size and position of the cube; this exit estimate is the heart of the matter. Step~2 applies it to large cubes and gives compact containment and conservativeness, Step~3 applies it at stopping times to bound the modulus of continuity and proves~(i), and Step~4 deduces~(ii) from~(i).

\smallskip

\emph{Step 1.} We show that, after rescaling every cube~$z+\cu_m$ to~$\cu_0$ as in~\eqref{tight:eq-local-time}, the process started in~$\cu_{-2}$ leaves~$\cu_0$ before time~$t$ with probability at most~$M\sqrt t$. Let~$\chi$ be a cutoff function of Lemma~\ref{tight:lem-static} between~$\cu_{-1}$ and~$\cu_0$, so that~$\chi$ takes values in~$[0,1]$, equals one on~$\cu_{-1}$, vanishes near~$\partial\cu_0$, and satisfies~$\mathscr E_A(\chi)\leq M$. For~$\lambda>0$, let~$u$ solve
\begin{equation*}
	\lambda u-(d\mu/dx)^{-1}\nabla\cdot(A\nabla u)=\lambda\chi\quad\text{in }\cu_0\,,
	\qquad u=0\quad\text{on }\partial\cu_0\,.
\end{equation*}
The maximum principle gives~$0\leq u\leq1$. Since~$u$ minimizes~$\mathscr E_A(v)+\lambda\|v-\chi\|_{L^2(\mu)}^2$ over~$v\in H^1_0(\cu_0)$, comparison with~$v=\chi$ gives~$\int_{\cu_0}|u-\chi|^2\,d\mu\leq\mathscr E_A(\chi)/\lambda\leq M/\lambda$. On~$\cu_{-1}$, the nonnegative function~$1-u$ satisfies~$\nabla\cdot(A\nabla(1-u))=\lambda(1-u)\,d\mu/dx\geq0$, and, since~$\chi=1$ on~$\cu_{-1}$, this estimate bounds the norm of~$1-u$ in~$L^2(\cu_{-1},\mu)$. Therefore Lemma~\ref{tight:lem-subharmonic}, whose constant does not depend on~$\lambda$, gives~$\sup_{\cu_{-2}}(1-u)\leq M\lambda^{-\nf12}$. Let~$\widehat X$ be the process obtained by rescaling~$z+\cu_m$ to~$\cu_0$, let~$\mathbf P_x^\omega$ and~$\mathbf E_x^\omega$ be its law and expectation started at~$x$, and let~$\tau_{\cu_0}$ be its exit time from~$\cu_0$. The stopped It\^o formula gives, for every~$x\in\cu_0$,
\begin{equation*}
	u(x)=\mathbf E_x^\omega\int_0^{\tau_{\cu_0}}
	\lambda e^{-\lambda s}\chi(\widehat X_s)\,ds
	\leq1-\mathbf E_x^\omega e^{-\lambda\tau_{\cu_0}}\,.
\end{equation*}
Taking~$\lambda=\nf1t$, we obtain
\begin{equation}
	\sup_{x\in\cu_{-2}}\mathbf P_x^\omega(\tau_{\cu_0}\leq t)
	\leq e\sup_{x\in\cu_{-2}}\mathbf E_x^\omega
	\exp\Bigl(-\frac{\tau_{\cu_0}}{t}\Bigr)
	\leq M\sqrt t\,.
	\label{tight:eq-fast-exit}
\end{equation}
The constant~$M$ does not depend on~$t$, and~$\E M\leq C$. In physical coordinates,~\eqref{tight:eq-fast-exit} says that the probability that~$X^{(L)}$, started at a point of~$z+\cu_{m-2}$, leaves~$z+\cu_m$ before time~$t$ is at most
\begin{equation}
	M\sqrt{t\ahom_{m\wedge L}}\,3^{-m}\,,
	\label{tight:eq-physical-exit}
\end{equation}
and the same bound holds for~$X$ with~$m$ in place of~$m\wedge L$. These bounds concern the processes killed on leaving a bounded cube and do not presume that the processes are conservative.

\smallskip

\emph{Step 2.} We prove compact containment and conservativeness. For triadic~$R\geq1$, a cube of side~$R$ in rescaled coordinates has time unit~$R^23^{2L}/\ahom_L$ at finite cutoff, and, by~\eqref{e.annealed.ordering}, at least~$R^23^{2L}/\ahom_L$ for the diffusion~$X$ without cutoff. After the cube sizes are adjusted by fixed factors,~\eqref{tight:eq-fast-exit} therefore gives, for every~$T>0$ and every~$R\geq1$ larger than a fixed multiple of~$\max_{x\in F}|x|$,
\begin{equation}
	\sup_L\E\left[\sup_{x\in F}K^\omega_{L,x}\left(\sup_{s\leq T}|Z_s|>R\right)\right]
	\leq C\frac{\sqrt T}{R}\,.
	\label{tight:eq-containment}
\end{equation}
Letting~$R\to\infty$ shows that, for almost every environment and every starting point in~$F$, the process does not explode. The same estimate then gives compact containment for the laws on continuous path space.

\smallskip

\emph{Step 3.} We bound the modulus of continuity and prove~\eqref{tight:eq-uniform-laws}. Fix~$R<\infty$,~$T>0$ and~$r>0$, fix~$k\in\N_0$ such that cubes of side~$3^{-k}$ have diameter at most~$r/4$, and cover the ball~$B_R$ by finitely many cubes~$y+\cu_{-k-2}$. For~$L\geq k$, the time unit of each of the cubes~$y+\cu_{-k}$, measured in units of~$3^{2L}/\ahom_L$, satisfies, by~\eqref{e.annealed.ordering},~$3^{-2k}\ahom_L/\ahom_{L-k}\geq3^{-(2+\frac{2\tau^2}{\log3})k}>0$. We apply~\eqref{tight:eq-fast-exit} to each cube of this finite cover and take the maximum of the finitely many random constants, whose expectation is at most the sum of their expectations. The strong Markov property then gives a constant~$C$, depending on~$R$ and~$r$, such that, for every stopping time~$S\leq T$ and every~$h>0$,
\begin{equation*}
	\sup_{L\geq k}\E\Bigl[\sup_{x\in B_R}K^\omega_{L,x}\Bigl(
	\sup_{0\leq s\leq h}|Z_{S+s}-Z_S|>r/4,\ S<\tau_{B_R}\Bigr)\Bigr]
	\leq C\sqrt h\,.
\end{equation*}

\smallskip

To obtain~\eqref{tight:eq-uniform-laws}, in which the supremum over starting points lies inside the expectation, we first restrict to the event on which the maximum of the random constants of this finite cover is at most~$M_0$. Its complement has probability at most~$C/M_0$, uniformly in~$L\geq k$. On this event, for every starting point in~$B_R$, the probability that the path moves by more than~$r/4$ within time~$h$ is at most~$CM_0\sqrt h$, and, by the strong Markov property, the same bound holds after every stopping time before exit from~$B_R$, conditionally on the path up to that stopping time.

\smallskip

Consider the successive exit times of the path from balls of radius~$r/4$, the first centered at the starting point and each later one at the position of the path at the preceding exit time. Choose~$h_0$ so small that this probability bound, with~$h=h_0$, is at most~$\nf1{16}$. If~$n>2T/h_0$ such exits occur before time~$T$, then at least~$\nf n2$ of the gaps between consecutive exit times are shorter than~$h_0$. The strong Markov property in the fixed environment and a union bound over the choices of these gaps bound the probability of this event by~$2^n16^{-\nf n2}$. Hence the number of exits before time~$T$ is tight, uniformly in the starting point.

\smallskip

For a fixed number of exits, the probability that some gap is shorter than~$h$ tends to zero as~$h\to0$. Since the path remains in a ball of radius~$r/4$ between consecutive exit times, this gives a uniform bound on the probability that the path oscillates by more than~$r$ on an interval of length~$h$. We choose first~$M_0$, then the number of exits, then~$h$.

\smallskip

For each of the finitely many~$L<k$, the coefficient is locally~$C^{1,1}$. Nonexplosion and continuous dependence on the starting point show that, almost surely, the laws of these diffusions started in~$F$ are uniformly tight, and dominated convergence gives the corresponding bound for the expected supremum. We apply the preceding construction along a sequence of integer time horizons and of spatial tolerances tending to zero, with summable error probabilities. Compact containment and the Arzel\`a--Ascoli theorem then prove~\eqref{tight:eq-uniform-laws}.

\smallskip

\emph{Step 4.} We prove~\eqref{tight:eq-random-laws}. By~(i), for every~$j\geq1$ there is a compact set~$\mathcal K_j$ of paths for which the left side of~\eqref{tight:eq-uniform-laws} is at most~$\ep2^{-2j}$. Let
\begin{equation*}
	\mathcal M_\ep\coloneqq\bigl\{P:P(\mathcal K_j^c)\leq2^{-j}
	\text{ for every }j\geq1\bigr\}\,.
\end{equation*}
The set~$\mathcal M_\ep$ is closed and uniformly tight, hence compact. Markov's inequality and a union bound over~$j$ prove~\eqref{tight:eq-random-laws}.
\end{proof}

\subsection{Convergence in probability of the path laws}
\label{mfd:sec-paths}

In this subsection we prove that the quenched path laws~$K^\omega_{N,x}$ converge in probability, locally uniformly in~$x$, to a continuous random kernel (Proposition~\ref{mfd:prop-quenched-convergence}). The limit is a conservative strong Markov process, symmetric with respect to~$\mu$, and its resolvents killed on leaving a cube are the minimizers of Section~\ref{mfd:sec-speed} (Proposition~\ref{mfd:prop-limit-properties} and Lemma~\ref{mfd:lem-killing}). In particular its mean exit times from bounded sets are finite (Corollary~\ref{mfd:cor-finite-exit}).

\smallskip

We use the killed resolvents to identify the limiting laws on~$C([0,\infty),\Rd)$. The ingredients are compact containment, the Markov identity at finite cutoff, and the fact that discounted integrals of finite-dimensional observations determine a continuous path law.

\smallskip

We record the short-time approximation for the killed resolvents. For~$z\in\mathbb Q^d$,~$n\in\mathbb Z$ and~$f\in C_0(z+\cu_n)$, extend~$f$ by zero. Splitting the occupation integral at a time~$h>0$ gives, for every~$r>0$,
\begin{align}
	\sup_N\E\|\lambda R^{z+\cu_n}_{N,\lambda}f-f\|_\infty
	& \leq\sup_{|x-y|\leq r}|f(x)-f(y)|
	\notag \\ & \qquad
	+2\|f\|_\infty\Bigl(e^{-\lambda h} +\sup_N\E\sup_{x\in\overline{z+\cu_n}}K^\omega_{N,x}
	\bigl(\sup_{s\leq h}|Z_s-x|>r\bigr)\Bigr)\,.
	\label{eq:mfd-resolvent-identity-approximation}
\end{align}
The last probability tends to zero as~$h\downarrow0$ by~\eqref{tight:eq-uniform-laws}. If a path is killed before~$h$ without moving farther than~$r$, its starting point is within~$r$ of the boundary, where~$f$ is small. Proposition~\ref{mfd:prop-uniform-resolvent} passes~\eqref{eq:mfd-resolvent-identity-approximation} to the limiting resolvents. In particular,~$\E\|\lambda R^{z+\cu_n}_\lambda f-f\|_\infty\to0$ as~$\lambda\to\infty$.

\begin{proposition}[Convergence in probability]
\label{mfd:prop-quenched-convergence}
Let~$F\subset\Rd$ be compact and~$\ep_0>0$. Then
\begin{equation*}
	\lim_{N,N'\to\infty}\P\!\left[
	\sup_{x\in F}d_{\mathrm{weak}}(K^\omega_{N,x},K^\omega_{N',x})
	>\ep_0\right]=0\,.
\end{equation*}
Consequently, the kernels~$x\mapsto K^\omega_{N,x}$ converge in probability, locally uniformly in~$x$, to a continuous random probability kernel~$x\mapsto K^\omega_x$ that is measurable with respect to the fields.
\end{proposition}
\begin{proof}
Use the Prokhorov metric~$d_{\mathrm P}$ associated with the complete metric on path space
\begin{equation*}
	\sum_{m\geq1}2^{-m}
	\bigl(1\wedge\sup_{t\leq m}|w(t)-w'(t)|\bigr)\,.
\end{equation*}
Write~$R_{N,\lambda}f(x)\coloneqq\int\!\int_0^\infty e^{-\lambda t}f(Z_t)\,dt\,dK^\omega_{N,x}$ for the whole-space resolvent. If~$y+\cu_\ell$, with~$y\in\mathbb Q^d$ and~$\ell\in\mathbb Z$, contains~$F$, splitting at the exit time gives
\begin{equation}
	\sup_{x\in F}|R_{N,\lambda}f(x)-R^{y+\cu_\ell}_{N,\lambda}f(x)|
	\leq\frac{\|f\|_\infty}{\lambda}
	\left(e^{-\lambda T}+\sup_{x\in F}K^\omega_{N,x}(\tau_{y+\cu_\ell}\leq T)\right)\,.
	\label{eq:mfd-path-localization}
\end{equation}
Choose first~$T$ large and then~$\ell$ large. Compact containment~\eqref{tight:eq-containment} makes the expectation of the right side arbitrarily small, uniformly in~$N$. Proposition~\ref{mfd:prop-uniform-resolvent} therefore makes the whole-space resolvents Cauchy in probability, locally uniformly in~$x$. Contraction and localization give the same conclusion for measurable random continuous sources~$g_N$ converging in probability locally uniformly to~$g$, provided~$\|g_N\|_\infty$ and~$\|g\|_\infty$ are bounded by one deterministic constant. The intermediate sources in the iterated resolvents below have this bound by contraction.

\smallskip

For bounded continuous~$f_1,\ldots,f_k$ and positive integers~$n_1,\ldots,n_k$, consider the bounded continuous path functional
\begin{equation*}
	\int_{(0,\infty)^k}e^{-\sum_{i=1}^k n_i s_i}
	\prod_{i=1}^k f_i\bigl(w(s_1+\cdots+s_i)\bigr)\,ds\,.
\end{equation*}
The Markov property and Fubini identify its expectation with~$R_{N,n_1}(f_1R_{N,n_2}(f_2\cdots R_{N,n_k}f_k))(x)$. Iterating the preceding resolvent convergence makes these expectations Cauchy in probability, uniformly for~$x\in F$. These functionals determine continuous path laws. After changing variables from~$s_i$ to~$e^{-s_i}$, varying the positive integers~$n_i$ gives the integrals of the finite-dimensional observations against every monomial on~$[0,1]^k$. Polynomial approximation and continuity in time then determine these observations. A separating algebra of spatial test functions determines the finite-dimensional laws.

\smallskip

Tightness~\eqref{tight:eq-uniform-laws} now transfers this Cauchy property to~$d_{\mathrm P}$. Indeed, on a compact set of paths carrying all but a small mass, finitely many continuous tests suffice to distinguish laws at a prescribed positive~$d_{\mathrm P}$ distance; the determining functionals approximate these tests there. Markov's inequality controls the mass outside that set, uniformly over~$x\in F$. A finite union bound proves the asserted Cauchy property. The finite-cutoff kernels are continuous in~$x$, by their Feller property and local path tightness. Completeness of the Prokhorov metric and a subsequence with summable error probabilities on a compact exhaustion of~$\Rd$ give a measurable continuous limiting kernel.

\smallskip

Finally, on the compact set of laws~$\mathcal M_\ep$ in~\eqref{tight:eq-random-laws}, every metric inducing the weak topology is uniformly equivalent to~$d_{\mathrm P}$. With probability at least~$1-2\ep$, both finite-cutoff kernels take values there for every~$x\in F$. Letting~$\ep\downarrow0$ transfers the assertion to~$d_{\mathrm{weak}}$.
\end{proof}

\subsubsection{The Markov property of the limit}

\begin{proposition}[Properties of the limit]
\label{mfd:prop-limit-properties}
Almost surely, the limiting kernel~$x\mapsto K^\omega_x$ has the following properties.
\begin{enumerate}[label=(\roman*),leftmargin=*]
\item The map~$x\mapsto K^\omega_x$ is continuous in the weak topology.
\item The laws~$(K^\omega_x)_x$ are the laws of a continuous, conservative, strong Markov process~$Z$.
\item The process is symmetric with respect to~$\mu$: for every~$t\geq0$ and all compactly supported continuous functions~$f$ and~$g$,
\begin{equation*}
	\int g(x)\int f(Z_t)\,dK^\omega_x\,\mu(dx)=\int f(x)\int g(Z_t)\,dK^\omega_x\,\mu(dx)\,.
\end{equation*}
\end{enumerate}
\end{proposition}

\begin{proof}
Continuity of the kernel and infinite lifetime hold by Proposition~\ref{mfd:prop-quenched-convergence}: the limiting laws are probability measures on~$C([0,\infty),\Rd)$. The finite-cutoff Markov identity passes to the limit for bounded continuous path functionals, using locally uniform convergence of the kernels and compact containment. A monotone-class argument and approximation of stopping times from above by dyadic times, which uses the continuity of the paths and of~$x\mapsto K^\omega_x$, then give the strong Markov property.

\smallskip

At finite cutoff the same identity holds for~$K^\omega_{N,x}$ and~$\mu_N$, by the symmetry of the generator in~$L^2(\mu_N)$. For compactly supported continuous~$f$ and~$g$, locally uniform convergence of the kernels and local weak convergence~$\mu_N\to\mu$ (Proposition~\ref{mfd:prop-chaos-growth}) pass this identity to the limit.
\end{proof}

\begin{lemma}[Killed resolvents of the limiting process]
\label{mfd:lem-killing}
Almost surely, for every~$z\in\mathbb Q^d$, every~$n\in\mathbb Z$ and every~$\lambda>0$, the resolvent of the process of Proposition~\ref{mfd:prop-limit-properties} killed on leaving~$z+\cu_n$ is~$R^{z+\cu_n}_\lambda$.
\end{lemma}
\begin{proof}
By Proposition~\ref{mfd:prop-quenched-convergence} and~\eqref{eq:mfd-42}, we choose a subsequence along which the whole-space kernels converge almost surely, locally uniformly in the starting point, and the killed resolvents converge uniformly on~$\overline{z+\cu_n}$, for a countable dense family of data and positive parameters. Fix an environment in this event. For nonnegative data, the discounted occupation functional stopped at~$\tau_{z+\cu_n}$ is lower semicontinuous on continuous path space, so the portmanteau theorem bounds its expectation under~$K^\omega_x$ from above by~$R^{z+\cu_n}_\lambda f(x)$.

\smallskip

For the opposite inequality, first take~$f=1$. Fix~$\ep>0$ and a smaller closed concentric cube~$\overline{z+r\cu_0}\Subset z+\cu_n$ so that, for all sufficiently large~$N$,~$\lambda R^{z+\cu_n}_{N,\lambda}1\leq\ep$ on~$(z+\cu_n)\setminus\overline{z+r\cu_0}$. Uniform convergence and the continuous zero-boundary representative of~$R^{z+\cu_n}_\lambda1$ allow this choice. Let~$\sigma_{h,\overline{z+r\cu_0}}$ be the first time on the grid~$h\N_0$ at which the path is outside~$\overline{z+r\cu_0}$, with~$e^{-\lambda\infty}=0$. The strong Markov property at~$\sigma_{h,\overline{z+r\cu_0}}$ at finite cutoff gives
\begin{equation*}
	(1-\ep)\int e^{-\lambda\sigma_{h,\overline{z+r\cu_0}}}\,dK^\omega_{N,x}
	\leq\int e^{-\lambda\tau_{z+\cu_n}}\,dK^\omega_{N,x}
	=1-\lambda R^{z+\cu_n}_{N,\lambda}1(x)\,.
\end{equation*}
The grid functional on the left is bounded and lower semicontinuous, since~$\overline{z+r\cu_0}$ is closed. Portmanteau passes the inequality to the limiting law. Continuity of paths gives~$e^{-\lambda\tau_{z+\cu_n}}\leq\liminf_{h\downarrow0}e^{-\lambda\sigma_{h,\overline{z+r\cu_0}}}$; Fatou's lemma and then~$\ep\downarrow0$ give~$R^{z+\cu_n}_\lambda1(x)\leq\int\int_0^{\tau_{z+\cu_n}}e^{-\lambda t}\,dt\,dK^\omega_x$. Together with the first inequality this proves equality for~$f=1$. Applying the first inequality to~$\|f\|_\infty+f$ and~$\|f\|_\infty-f$ gives equality for every bounded continuous~$f$. Contraction, continuity in~$\lambda$, and a countable intersection give the identity for all~$\lambda>0$ and all cubes on one event of probability one.

\smallskip

To allow~$f\in C(\overline{z+\cu_n})$, we multiply it by bounded compactly supported cutoff functions increasing to one in~$z+\cu_n$. The resulting data converge in~$L^2(\mu)$, and the variational identity and the trace bound give the convergence of their resolvents in energy. The common~$C^{\nf14}$ bound then gives uniform convergence, and dominated convergence for the path integral weighted by~$e^{-\lambda t}$ proves the identity for~$f$.
\end{proof}

The identification of the killed resolvents also gives the finite mean exit times asserted in Theorem~\ref{t.scaling.limit}.

\begin{corollary}[Finite mean exit times]
\label{mfd:cor-finite-exit}
Almost surely,~$\sup_{x\in U}\int\tau_U\,dK^\omega_x<\infty$ for every bounded open set~$U$.
\end{corollary}
\begin{proof}
Choose~$z\in\mathbb Q^d$ and~$n\in\mathbb Z$ such that~$z+\cu_n$ contains~$\overline U$. Lemma~\ref{mfd:lem-killing} and~\eqref{eq:mfd-uniform-lambda}, with~$f=1$, give, for every~$\lambda>0$,
\begin{equation*}
	\sup_{x\in z+\cu_n}\int\!\int_0^{\tau_{z+\cu_n}}e^{-\lambda t}\,dt\,dK^\omega_x
	=\|R^{z+\cu_n}_\lambda1\|_\infty\leq M_{z+\cu_n}
	\,.
\end{equation*}
Monotone convergence as~$\lambda\downarrow0$ and~$\tau_U\leq\tau_{z+\cu_n}$ prove the assertion.
\end{proof}

\subsection{Almost-sure convergence by comparison of two cutoffs}
\label{mfd:sec-finite-comparison}

In this subsection we prove that the minimal energies with fixed boundary values, and the Dirichlet Green quadratic forms and Neumann energies for fixed data, converge almost surely on every cube along the full sequence of cutoffs (Propositions~\ref{mfd:prop-as-dirichlet} and~\ref{mfd:prop-as-response-bank}). By the Borel--Cantelli lemma, it suffices to compare two cutoffs with error probabilities that are summable in~$N$.

\smallskip

We compare the finite-volume energies at two cutoffs. On a small cube, rescaling to the unit cube writes the ratio of the two energies as a deterministic ratio of normalizations times a ratio of two unit-cube energies. To combine these local comparisons, we join harmonic extensions on disjoint small cubes, using the finite-cutoff counterpart (Lemma~\ref{mfd:lem-finite-good-levels}) of the bound of Proposition~\ref{mfd:prop-allchain} on the number of bad cells along chains of nested cubes.

\smallskip

Our control of the ratios of the two normalizations follows the argument of~\citet[Section~3.4]{DGComparison}. There, two random metrics were compared by joining paths around annuli at a positive proportion of scales into a single path. Upper and lower bounds on crossing distances then excluded ratios of the two normalizations that are too large or too small. \citet[Lemma~3.11]{Devlin2024} strengthened this argument by using convergence in probability, which rules out ratios bounded away from one on a positive proportion of scales. We apply the argument of~\citet[Lemma~3.11]{Devlin2024} to the minimal energy with a fixed nonconstant affine boundary value. If the ratios were at most a fixed number less than one on a positive proportion of scales, then joining harmonic extensions on cubes of those scales would bound this energy at a finite cutoff by a fixed factor less than one times its limit, which is positive. This would contradict the convergence in probability. Interchanging the two cutoffs excludes, in the same way, ratios at least a fixed number greater than one. On the remaining scales the ratios are close to one, and Proposition~\ref{mfd:prop-as-dirichlet} turns the comparison on these scales into bounds with summable error probabilities.

\smallskip

We first record how a minimal energy with given boundary values on a small cube rescales. The maps~$\Theta_{k,z}$ translate and rescale the fields; their increments about the origin are taken only when forming~$H$. We write~$\Lambda^{\Theta_{k,z}\omega}_{m,\cu_0}$ for the minimal energy on~$\cu_0$ with given boundary values at cutoff~$m$ in this environment.

\begin{lemma}[Rescaling of minimal energies on small cubes]
\label{mfd:lem-local-normalizations}
For every integer~$k\geq0$, the limit
\begin{equation*}
	\bar\rho(k)=\lim_{N\to\infty}\frac{\ahom_{N-k}}{\ahom_N}
\end{equation*}
exists and is positive. Let~$z\in\Rd$. For~$N\geq k$, the cell scalar~$s_N(z+3^{-k}\cu_0)$ of Section~\ref{mfd:sec-environment} equals
\begin{equation*}
	\frac{\ahom_{N-k}}{\ahom_N}
	\exp\Bigl(H(z)+\sum_{j=0}^{k-1}\bigl(\g_{-j}(z)-\tau^2\bigr)\Bigr)\,,
\end{equation*}
the value at the center of~$z+3^{-k}\cu_0$ of the factor in~$A_N$ coming from the normalization and the fields whose scales exceed the side length of~$z+3^{-k}\cu_0$. Then, for every~$N\geq k$ and~$g\in H^1(z+3^{-k}\cu_0)$,
\begin{equation*}
	\Lambda_{N,z+3^{-k}\cu_0}(g)
	=3^{-(d-2)k}s_N(z+3^{-k}\cu_0)\,
	\Lambda_{N-k,\cu_0}^{\Theta_{k,z}\omega}(g(z+3^{-k}\,\cdot))\,.
\end{equation*}
For each fixed~$z$, the same identity holds almost surely for the minimal energies of the limit form, with~$\ahom_{N-k}/\ahom_N$ in~$s_N(z+3^{-k}\cu_0)$ replaced by~$\bar\rho(k)$, for every boundary trace on which these are defined.
\end{lemma}
\begin{proof}
Reindexing the fields gives, for every~$N\geq k$,
\begin{equation*}
	A_N(z+3^{-k}y)=s_N(z+3^{-k}\cu_0)A_{N-k}^{\Theta_{k,z}\omega}(y)\,.
\end{equation*}
The increments of~$H$ and of the fields~$\g_{-j}$ with~$j<k$ form the field~$H$ of the new environment; the fields with~$j\geq k$ retain their values. The energy identity follows by changing variables~$x=z+3^{-k}y$. Stationarity, scaling and independence show that~$\Theta_{k,z}$ preserves the joint field law. The same transformed environment serves every cutoff. In particular, dividing the two energies at cutoffs~$N$ and~$M$ by their respective scalars~$s_N(z+3^{-k}\cu_0)$ and~$s_M(z+3^{-k}\cu_0)$ cancels those scalars exactly and gives the two unit-cube energies at cutoffs~$N-k$ and~$M-k$ in this one environment.

\smallskip

To prove that the normalizations converge, fix a nonconstant affine datum and choose~$z\in\mathbb Q^d$, which does not affect~$\bar\rho(k)$. Proposition~\ref{mfd:prop-response-compact} and Theorem~\ref{mfd:thm-c1} give convergence in probability, along the full sequence of cutoffs, of the minimal energies on the cubes~$w+\cu_n$ with~$w\in\mathbb Q^d$ and~$n\in\mathbb Z$ for the data in~$\mathcal D$, and the seminorm inequality for~$\Lambda^{\nf12}$ together with~\eqref{eq:mfd-2} extends this to the affine datum. Hence both minimal energies in the identity converge in probability; for the one on~$\cu_0$ in the environment~$\Theta_{k,z}\omega$ this holds because~$\Theta_{k,z}$ preserves the law of the fields. Their limits are strictly positive: by Lemma~\ref{mfd:lem-coercivity}, each of these minimal energies is at least a positive constant times~$Q_N^{-1}$, and~$\sup_N\|Q_N\|_{L^p}<\infty$. By~\eqref{e.annealed.ordering}, for fixed~$k$ the deterministic ratios~$\ahom_{N-k}/\ahom_N$ lie in~$[1,e^{2\tau^2k}]$, a compact subset of~$(0,\infty)$. Two subsequential limits of these ratios must coincide, since otherwise the same positive minimal energy of the limit form would equal one random variable multiplied by two different constants. Finally, we pass to the limit in the identity for the countable family of smooth traces of Section~\ref{mfd:sec-local-form}, and the bound~\eqref{eq:mfd-2} for H\"older boundary values of Lemma~\ref{mfd:lem-extension} gives the identity for every trace on which the minimal energy of the limit form is defined.
\end{proof}

We compare the two cutoffs on a cell~$z+3^{-k}\cu_0$ through the trace on~$\partial(z+3^{-k}\cu_0)$ of a single function, which is H\"older continuous but otherwise arbitrary. The next lemma reduces this comparison to finitely many smooth traces.

\begin{lemma}[Reduction to finitely many smooth boundary values]
\label{mfd:lem-finite-trace-tests}
Let~$1/2<\beta<\alpha<1$,~$C_1<\infty$ and~$\varpi>0$. There exists a finite family of smooth functions on~$\partial \cu_0$ with the following property. Let~$\Lambda_1$ and~$\Lambda_2$ be minimal energies on~$\cu_0$ with given boundary values such that~$\Lambda_i(g)\leq C_1\|g\|_{C^\beta/\mathbb R}^2$ for~$i\in\{1,2\}$ and every~$g\in C^\beta(\partial \cu_0)$, and such that~$|\Lambda_1(\varphi)-\Lambda_2(\varphi)|\leq\varpi/2$ for every function~$\varphi$ of the family. Then, for every~$g\in C^\alpha(\partial \cu_0)$,
\begin{equation*}
	|\Lambda_1(g)-\Lambda_2(g)|
	\leq\varpi\|g\|_{C^\alpha/\mathbb R}^2\,.
\end{equation*}
\end{lemma}
\begin{proof}
Since~$\Lambda_i^{\nf12}$ is a seminorm, the triangle inequality gives, for all~$g,\varphi\in C^\beta(\partial \cu_0)$,
\begin{equation*}
	|\Lambda_i(g)-\Lambda_i(\varphi)|
	\leq C_1\|g-\varphi\|_{C^\beta/\mathbb R}
	\bigl(\|g\|_{C^\beta/\mathbb R}+\|\varphi\|_{C^\beta/\mathbb R}\bigr)\,.
\end{equation*}
The unit ball of~$C^\alpha/\mathbb R$ is precompact in~$C^\beta/\mathbb R$, so for every tolerance it has a finite net consisting of smooth functions. Take the family to be a sufficiently fine such net, and apply the last display to~$\Lambda_1$ and to~$\Lambda_2$. This gives the assertion for~$g$ in the unit ball, and homogeneity gives it for every~$g$.
\end{proof}

The number of functions in the family and their smooth norms may depend on~$\varpi$; the exponents and the moment orders do not. The constant~$C_1$ is the one in Lemma~\ref{mfd:lem-finite-good-cell}(iii) below, and it is fixed before the family is chosen.

\smallskip

On a good cell we also need, at a finite cutoff, the local H\"older estimate and the bound~\eqref{eq:mfd-2} for H\"older boundary values that Lemma~\ref{mfd:lem-goodext} provides for the limit form. The scalar~$s_N(z+3^{-k}\cu_0)$ of Lemma~\ref{mfd:lem-local-normalizations} is the finite-cutoff analogue of the scalar~$s_{\mathscr E}(z+3^{-k}\cu_0)$ of Section~\ref{mfd:sec-compare}.

Fix~$1/2<\beta<\alpha<1$, a subdivision factor~$3^{H_1}$, a distance factor~$b_0\geq3$ and a rate~$\varkappa>0$ for the exponential bounds below. For~$N\geq k\geq0$,~$z\in\Rd$ and~$r=3^{-k}$, rescaling~$z+r\cu_0$ to~$\cu_0$ gives cutoff~$m=N-k$. In the next lemma,~$A$ denotes either~$A_N$ or~$A_N^0=e^{-H}A_N$,~$s(z+r\cu_0)$ its respective scalar~$s_N(z+r\cu_0)$ of Lemma~\ref{mfd:lem-local-normalizations} or~$s_N^0(z+r\cu_0)=e^{-H(z)}s_N(z+r\cu_0)$, and~$\Gamma(u)$ its energy measure.

\begin{lemma}[Good cells at finite cutoff]
\label{mfd:lem-finite-good-cell}
The thresholds for the sums over intervals of scales, the finite moment orders and~$\gmcdelta_0>0$ can be chosen in this order so that the following holds for~$\gmcdelta\leq\gmcdelta_0$. There are deterministic constants~$C,C_1<\infty$, independent of the cutoff, the cell, the source and the tolerance~$\varpi$ of Lemma~\ref{mfd:lem-finite-trace-tests}, and good events~$\mathcal G_N(z+r\cu_0)$ on which the conclusions below hold for both coefficients simultaneously.
\begin{enumerate}[label=(\roman*),leftmargin=*]
\item Up to a null set, the complement of~$\mathcal G_N(z+r\cu_0)$ is covered by events indexed by integers~$h\geq1$, the event of index~$h$ being measurable with respect to~$\sigma(\g_{-j}:k-h\leq j\leq k+2h)$ and having probability at most~$e^{-\varkappa h}$.
\item Assume~$H_1\geq4$. Suppose~$\mathcal G_N(z+r\cu_0)$ holds and~$z+r\cu_0$ lies in its parent~$z'+3^{H_1}r\cu_0$ of side~$3^{H_1}r$ with~$\operatorname{dist}(z+r\cu_0,\partial(z'+3^{H_1}r\cu_0))\geq b_0r$. For every~$f\in L^\infty(z'+3^{H_1}r\cu_0)$ and every weak solution~$u\in H^1(z'+3^{H_1}r\cu_0)$ of~$-\nabla\cdot(A\nabla u)=f$ on~$z'+3^{H_1}r\cu_0$, there is a constant~$c$ such that
\begin{equation}
	\|u(z+r\,\cdot)-c\|_{C^\alpha(\cu_0)}
	\leq C r^{(2-d)/2}s(z+r\cu_0)^{-\nf12}\Gamma(u)(z'+3^{H_1}r\cu_0)^{\nf12}
	+Cr^2s(z+r\cu_0)^{-1}\|f\|_\infty\,.
	\label{eq:mfd-finite-good-holder}
\end{equation}
\item On~$\mathcal G_N(z+r\cu_0)$, for every~$g\in C^\beta(\partial \cu_0)$, the normalized minimal energy satisfies
\begin{equation*}
	r^{2-d}s(z+r\cu_0)^{-1}\Lambda_{N,z+r\cu_0}\bigl(g((\,\cdot-z)/r)\bigr)
	\leq C_1\|g\|_{C^\beta/\mathbb R}^2\,.
\end{equation*}
\end{enumerate}
\end{lemma}
\begin{proof}
The proof proceeds in four steps. Step~1 defines~$\mathcal G_N(z+r\cu_0)$ and proves~(i) for it. Step~2 bounds the normalized coefficient and its gradient on~$z+r\cu_0$ and its fixed enlargement. Step~3 proves~(ii), by the iteration of Lemma~\ref{l.iteration.lemma.GMC} above the smallest field scale and an interior Sobolev estimate below it. Step~4 adds condition~(a), proves~(iii) and intersects the good events of the two coefficients. The main point is Step~3: the growth of the coefficient bounds of Step~2 below the smallest field scale must be absorbed by the gain~$\alpha_1-\alpha$ of the H\"older exponent, which fixes how small~$s_*$ and~$\tau^2$ must be.

\smallskip

\emph{Step 1.} We define the event~$\mathcal G_N(z+r\cu_0)$. It imposes the conditions of~$\mathcal G_N^{\rm det}(z+r\cu_0)$ on sums over intervals of scales, down to the smallest field scale, with an exponent~$\alpha_1\in(\alpha,1)$ in place of~$\alpha$ in the iteration of Lemma~\ref{l.iteration.lemma.GMC}. We choose the integrability exponent used in Step~3 so that the Sobolev embedding gives a H\"older exponent larger than~$\alpha_1$, then a small test exponent~$s_*\in(0,\nf1{32}]$, and add the root and enlargement tests
\begin{equation*}
	\mathsf F_{m+\iota,\,3^Nz}\leq\ep,\qquad
	\mathsf M_{m+\iota,\,3^Nz}\leq12
	\qquad(\iota=0,1),\qquad m=N-k\,.
\end{equation*}
These are the quantities of Section~\ref{mfd:sec-branches}, evaluated in the relabelled fields~$\g_i^{(N)}(y)\coloneqq\g_{i-N}(3^{-N}y)$ with~$s=s_*$. The index~$m$ belongs to~$z+r\cu_0$, and~$m+1$ to its fixed concentric enlargement~$z+3r\cu_0$; their stationary centers are both~$3^Nz$. In the original indexing of the fields~$\g_i$ their field windows are centered at~$-k$ and~$1-k$, respectively.

\smallskip

For these extra tests use the continuous ramps of Section~\ref{mfd:sec-branches}, with threshold pairs~$(\ep/2,\ep)$ and~$(6,12)$; they are two of the three summands of~$\mathsf Z_{m+\iota,\,3^Nz}$. Their zero-disorder values are zero. Lemma~\ref{mfd:lem-band} approximates them by quantities depending on windows about the same root indices, with exponentially decreasing errors, uniformly in~$N$ and~$k$. Its proof applies verbatim with~$s=s_*$, because Section~\ref{mfd:sec-branches} allows every exponent~$s\in(0,\nf1{32}]$ fixed before~$\gmcdelta$. The proof also gives an~$L^p$ bound of order~$C_p\gmcdelta^c$ for these ramps. The telescoping argument of Lemma~\ref{mfd:lem-witness}, applied to each ramp at this one cutoff, therefore covers its failure by events with a prescribed exponential rate after the moment orders are fixed and~$\gmcdelta$ is reduced. Enlarge and reindex the windows to absorb their fixed widths and the offset~$1$, and increase the rate before taking the finite union. All these events then have the windows and the probability bound in~(i), centered at the physical root depth~$k$, uniformly in~$m$. For the conditions on sums over intervals of~$D\leq m$ scales, we apply, at the single fixed cutoff~$N$, the telescoping argument from the proof of Lemma~\ref{mfd:lem-witness}, which uses dyadically increasing intervals of field indices. Lemmas~\ref{mfd:in-prefix} and~\ref{mfd:lem-band} are uniform at that cutoff, so the union over the at most~$Ce^{CD}$ spatial choices has the same bound. Choose the moment order used in the approximation so large that the resulting exponential rate exceeds both~$\varkappa$ and the exponential growth rate of the number of spatial choices, and then decrease~$\gmcdelta$. This gives~(i) directly at finite cutoff.

\smallskip

\emph{Step 2.} We bound the coefficient on the cell and its fixed enlargement. The term~$j=m$ in the root product test gives
\begin{equation*}
	\sup_{y\in 3^Nz+\cu_{m+1}}\sum_{i=0}^{m}|\g_i^{(N)}(y)|
	\leq\log12+\frac{s_* m}{8}\log3\,.
\end{equation*}
Indeed, the first product in that term is bounded by~$12\cdot3^{s_* m/8}$ on the larger cube~$3^Nz+\cu_{2m+1}$, and includes all indices~$0,\ldots,m$. The term~$j=0$ in the same test bounds
\begin{equation*}
	\sup_{y\in 3^Nz+\cu_{m+1}}
	\sum_{i=m+1}^{\infty}
	|\g_i^{(N)}(y)-\g_i^{(N)}(3^Nz)|
	\leq\frac14\log12\,.
\end{equation*}
For~$A_N^0$ this sum stops at~$i=N$, so the same bound applies. The identity of Lemma~\ref{mfd:lem-local-normalizations}, with the respective scalar~$s(z+r\cu_0)$, expresses the logarithm of the normalized coefficient as the sum of the fine fields~$i=0,\ldots,m$ and these coarse-field increments, minus~$\log(e^{(m+1)\tau^2}\ahom_m)$. Since~$e^{-2\tau^2}\leq\ahom_0\leq1$ by~\eqref{e.ord.bounds.truncated} and Jensen's inequality,~\eqref{e.annealed.ordering} gives~$e^{-2\tau^2(m+1)}\leq\ahom_m\leq1$, hence~$|\log(e^{(m+1)\tau^2}\ahom_m)|\leq\tau^2(m+1)$. Thus the logarithm of the normalized coefficient is bounded by~$C+C(s_*+\tau^2)m$, for both coefficients.

\smallskip

The root field test controls the gradients separately. Its term~$j=m$ bounds the sum of~$3^i|\nabla\g_i^{(N)}|$ for~$0\leq i\leq m$ by~$3^{s_* m/8}$. For~$i=m+j>m$, its term of depth~$j$ gives~$3^i\|\nabla\g_i^{(N)}\|_\infty\leq3^{s_* j/8}$ on~$3^Nz+\cu_{m+1}$. Hence
\begin{equation*}
	\sum_{i=m+1}^{\infty}
	\|\nabla\g_i^{(N)}\|_{L^\infty(3^Nz+\cu_{m+1})}
	\leq3^{-m}\sum_{j\geq1}3^{-(1-s_*/8)j}\,.
\end{equation*}
The series on the right converges, and the same bound holds with every finite upper limit in the sum on the left. In unit coordinates~$x$, the stationary coordinate is~$3^Nz+3^mx$, and the smallest field scale is~$3^{-m}$. For~$A=A_N$, the coefficient in these unit coordinates divided by~$s(z+r\cu_0)$ is~$A_{N-k}^{\Theta_{k,z}\omega}$, and the value and gradient bounds give
\begin{equation*}
	\|(A_{N-k}^{\Theta_{k,z}\omega})^{-1}\|_\infty+\|A_{N-k}^{\Theta_{k,z}\omega}\|_\infty
	+3^{-m}\|\nabla\log A_{N-k}^{\Theta_{k,z}\omega}\|_\infty
	\leq C3^{C(s_*+\tau^2)m}\,.
\end{equation*}
The norms here are on~$3\cu_0$. The normalized rescaling of~$A_N^0$ is~$A_{N-k}^{\Theta_{k,z}\omega}$ with~$H$ replaced by its partial sum~$H_k$ in that environment, so the same bounds hold. They also hold with~$H$ omitted or replaced by each of its partial sums, since all omitted terms were bounded in absolute value.

\smallskip

\emph{Step 3.} We prove~(ii). Above the scale~$3^{-m}$, the iteration of Lemma~\ref{l.iteration.lemma.GMC} on the event controlling sums over intervals of scales gives~\eqref{eq:mfd-finite-good-holder} with exponent~$\alpha_1$. Below~$3^{-m}$, we use the interior Sobolev estimate, with this exponent, of Step~4 of the proof of Proposition~\ref{mfd:prop-growth} on balls of radius~$a3^{-m}/[C(1+3^{-m}\|\nabla\log A_{N-k}^{\Theta_{k,z}\omega}\|_\infty)]$ (for~$A_N^0$, read~$A_{N-k}^{\Theta_{k,z}\omega}$ with~$H_k$ in place of~$H$), with a small constant~$a>0$ depending only on~$d$ and the integrability exponent, on which the oscillation of~$\log A_{N-k}^{\Theta_{k,z}\omega}$ is small. After matching with the estimate at scale~$3^{-m}$, the energy term carries the factor~$3^{-m(\alpha_1-\alpha)}$ and the source term the factor~$3^{-m(2-\alpha)}$, each multiplied only by fixed powers of the quantities in the last display. We choose~$s_*$, and then~$\gmcdelta$ (which bounds~$\tau^2$ by~\eqref{e.tau.bound}), so that~$C(s_*+\tau^2)<\alpha_1-\alpha$, where~$C$ depends only on~$d$ and the integrability exponent. The net exponents of~$3^{-m}$ then remain positive, and~\eqref{eq:mfd-finite-good-holder} holds at all radii with a deterministic constant.

\smallskip

\emph{Step 4.} We prove~(iii) and obtain one good event for both coefficients. First intersect the event of Step~1 with condition~(a) on~$z+r\cu_0$ and its fixed enlargements, using the finite-cutoff ellipticity bounds of~$\mathcal G_N^{\rm det}(z+r\cu_0)$ at the trace exponent of Lemma~\ref{mfd:lem-extension}. These normalized ellipticity conditions have the centered moment and band bounds of Lemma~\ref{mfd:lem-band}, also with~$H$ replaced by zero. Their failure is covered at the fixed cutoff by the finite-test telescoping argument of Lemma~\ref{mfd:lem-witness}. Add these covers to those of Step~1, using a larger rate before the finite union. Thus~(i) covers condition~(a) as well. On that condition, Lemma~\ref{mfd:lem-extension} gives~(iii), with deterministic~$C_1$ independent of~$\varpi$. Finally, we apply the construction separately to the coefficient with and without its factor~$e^H$, and take the intersection of the two good events. When~$H$ is replaced by one of its partial sums, the same bounds hold uniformly in the number of terms, by the summable bounds on the increments of the fields of scale larger than one.
\end{proof}

The comparison of Lemma~\ref{mfd:lem-finite-trace-tests} on a cell holds when the normalized minimal energies at the two cutoffs are close on the finite family of smooth boundary values of that lemma. By the convergence in probability, this comparison fails only with small probability, but the counting argument of Section~\ref{mfd:sec-branches} requires more: the failure must be covered by events of exponentially small probability, each measurable with respect to the fields with indices in a short interval.

\begin{lemma}[Rare events and field indices]
\label{mfd:lem-rare-tests}
There is an absolute constant~$C\geq4$ with the following property. Let~$p\in[1,\infty)$,~$a,C_2>0$ and~$\varkappa>0$ satisfy~$C\varkappa<ap\log3$. Then there is an integer~$h_0\geq1$, depending only on these four parameters, such that the following holds. Let~$n_0\in\Z$ and let~$(X_m)_{m\in\N_0}$ be random variables with values in~$[0,1]$ that converge to zero in probability. For every integer~$h\geq0$, let~$\mathcal F_h$ be the~$\sigma$-field generated by the fields~$\g_j$ with~$n_0-h\leq j\leq n_0+h$. Assume that, for every~$m,h\in\N_0$,
\begin{equation*}
	\|X_m-\E[X_m\mid\mathcal F_h]\|_{L^p}
	\leq C_23^{-ah}\,.
\end{equation*}
Then there is an integer~$m_0\geq0$ such that, for every~$m\geq m_0$, up to a null set,~$\{X_m>1/2\}$ is contained in a union of events indexed by integers~$h\geq h_0$. The event with index~$h$ is measurable with respect to~$\sigma(\g_j:n_0-2h\leq j\leq n_0+h)$ and has probability at most~$e^{-\varkappa h}$.
\end{lemma}
\begin{proof}
Let~$X_m^{(\ell)}\coloneqq\E[X_m\mid\mathcal F_{2^\ell h_0}]$ for~$\ell\geq0$, with~$h_0$ to be chosen. The summable~$L^p$ errors imply~$X_m^{(\ell)}\to X_m$ almost surely. Since~$\sum_{\ell\geq1}2^{-\ell-3}=1/8$, the event~$\{X_m>1/2\}$ is covered, up to a null set, by~$\{X_m^{(0)}\geq3/8\}$ and the events~$\{X_m^{(\ell)}-X_m^{(\ell-1)}\geq2^{-\ell-3}\}$ for~$\ell\geq1$. Markov's inequality bounds the probability of the~$\ell$th increment event by~$(16C_2)^p2^{p\ell}3^{-ap2^{\ell-1}h_0}$. Choose~$h_0$ large enough in terms of~$a,p,C_2,\varkappa$ that this is at most~$e^{-\varkappa2^\ell h_0}$ for every~$\ell\geq1$. Since~$0\leq X_m\leq1$ and~$X_m\to0$ in probability,~$\E X_m\to0$. Choose~$m_0$ so that~$\E X_m\leq e^{-\varkappa h_0}/8$ for~$m\geq m_0$; Markov's inequality then bounds the first event by~$e^{-\varkappa h_0}$. Index these events by~$h=2^\ell h_0$ and use the empty event at the remaining indices. Each event is~$\mathcal F_h$-measurable, hence measurable in the stated interval.
\end{proof}

We verify the conditioning and convergence inputs for the finite trace tests. Let~$\varphi_1,\ldots,\varphi_{n_\varpi}$ be the smooth family of Lemma~\ref{mfd:lem-finite-trace-tests}. Let~$\psi_\varpi$ vanish on~$[0,\varpi/4]$, equal one on~$[\varpi/2,\infty)$ and be linear in between. For~$z\in\Rd$,~$r\coloneqq3^{-k}$ and~$M\geq N\geq k$, let
\begin{multline*}
	X_{N,M}(z+r\cu_0)\coloneqq\psi_\varpi\Bigl(\max_{i\leq n_\varpi}
	\Bigl|r^{2-d}s_N(z+r\cu_0)^{-1}\Lambda_{N,z+r\cu_0}\bigl(\varphi_i((\cdot-z)/r)\bigr)
	\\ -r^{2-d}s_M(z+r\cu_0)^{-1}\Lambda_{M,z+r\cu_0}\bigl(\varphi_i((\cdot-z)/r)\bigr)\Bigr|\Bigr)\,.
\end{multline*}
By Lemma~\ref{mfd:lem-local-normalizations}, these are the minimal energies on~$\cu_0$ at the cutoffs~$N-k$ and~$M-k$, in the same environment~$\Theta_{k,z}\omega$. No moment of either cell scalar enters this test.

\smallskip

Write~$\mathcal F_\ell\coloneqq\sigma(\g_j:|j|\leq\ell)$. Extend~$\varphi_i$ smoothly to~$\overline{\cu_0}$. Proposition~\ref{mfd:prop-growth} supplies its Morrey bound, with a factor~$\|\varphi_i\|_{C^2}^2$, and~\eqref{eq:mfd-2} and Lemma~\ref{mfd:lem-cell-ellipticity} supply uniform moments of its minimal energy. Thus Lemma~\ref{mfd:lem-15} bounds the change caused by resampling~$\g_{-j}$ by~$C_{p,i}\gmcdelta3^{-bj/8}$ in~$L^p$ for~$0\leq j\leq m$, where~$b$ is the exponent of Lemma~\ref{mfd:lem-strips}; for~$j>m$ that field is absent. For~$j\geq1$, resampling~$\g_j$ changes~$H$ by its increment about the origin, whose supremum has sub-Gaussian size~$C\gmcdelta3^{-j}$. The coefficient comparison~\eqref{eq:mfd-mult} and H\"older's inequality with the moments of order~$3p$ of this minimal energy bound this change by~$C_{p,i}\gmcdelta3^{-j}$. For every~$m,\ell\geq0$, Lemma~\ref{mfd:lem-resampling-sum} gives
\begin{equation*}
	\|\Lambda_{m,\cu_0}(\varphi_i)-\E[\Lambda_{m,\cu_0}(\varphi_i)\mid\mathcal F_\ell]\|_{L^p}
	\leq C_{p,i}\gmcdelta\left(\sum_{j=\ell+1}^{m}3^{-bj/8}+\sum_{j>\ell}3^{-j}\right)
	\leq C_{p,i}\gmcdelta3^{-b\ell/8}\,.
\end{equation*}
The first sum is empty when~$\ell\geq m$. This is the proof of the first item of Proposition~\ref{mfd:prop-16}, extended to windows beyond the cutoff; its stated restriction~$m\geq\ell$ is unnecessary here.

\smallskip

Proposition~\ref{mfd:prop-response-compact} and Theorem~\ref{mfd:thm-c1} give full-sequence convergence of the unit-cube energies on the original field space. Their uniform higher moments give~$L^p$ convergence. The seminorm inequality for~$\Lambda_{m,\cu_0}^{\nf12}$ and~\eqref{eq:mfd-2} extend it from the dense trace family~$\mathcal D$ to each~$\varphi_i$. Denote the limits by~$\Lambda_{\mathscr E,\cu_0}(\varphi_i)$ and, for~$m\in\N_0$, let
\begin{equation*}
	X_m\coloneqq\psi_{\varpi/2}\Bigl(\max_{i\leq n_\varpi}|\Lambda_{m,\cu_0}(\varphi_i)-\Lambda_{\mathscr E,\cu_0}(\varphi_i)|\Bigr)\,.
\end{equation*}
Since~$\psi_{\varpi/2}(x)>\ep$ forces~$x>\varpi/8$, we have~$\P[X_m>\ep]\to0$ as~$m\to\infty$, for every~$\ep>0$.

\smallskip

Contraction of conditional expectation, the Lipschitz constant~$8/\varpi$ of~$\psi_{\varpi/2}$, the finite-maximum inequality in~$L^p$ and the~$L^p$ convergence of the energies turn the bound for~$\Lambda_{m,\cu_0}(\varphi_i)$ above into
\begin{equation*}
	\|X_m-\E[X_m\mid\mathcal F_\ell]\|_{L^p}
	\leq \frac{32\gmcdelta}{\varpi}
	\left(\sum_{i=1}^{n_\varpi} C_{p,i}^p\right)^{1/p}3^{-b\ell/8}
\end{equation*}
for every~$m,\ell\geq0$. This is the conditioning hypothesis of Lemma~\ref{mfd:lem-rare-tests} at~$n_0=0$, with~$a=b/8$ and~$C_2$ the constant on the right. Choose~$p\geq2$ so that~$(b/8)p\log3>C(\varkappa+\log2)$ for the counting rate~$\varkappa$, choose the finite higher orders needed for the Morrey bounds, the moments of the minimal energies and the product bounds, and then decrease~$\gmcdelta$. These orders and this disorder threshold are independent of~$\varpi$ and of the smooth net. Only the norm factors~$C_{p,i}$, the size~$n_\varpi$ of the net, and hence~$C_2$ and the initial window depend on~$\varpi$.

\smallskip

Apply Lemma~\ref{mfd:lem-rare-tests} with~$n_0=0$ and at rate~$\varkappa+\log2$ to the sequence~$(X_m)$, which satisfies its hypotheses by the two preceding paragraphs. The resulting~$h_0$ and~$m_0$ are fixed before the cell and the two cutoffs. Outside the two covers at~$m=N-k$ and at~$m=M-k$, the triangle inequality bounds the difference of the two energy vectors by~$3\varpi/8$ in maximum norm. Pulling these covers back by~$\Theta_{k,z}$ therefore covers~$\{X_{N,M}(z+r\cu_0)>1/2\}$; it preserves their probabilities and gives the field interval~$[-k-2h,-k+h]$. Taking the union of the two events at each~$h\geq1$ gives probability at most~$e^{-\varkappa h}$. Thus the same~$h_0,m_0$ work for every cell and pair with~$M\geq N\geq k+m_0$. The countably many cells and cutoffs used below have a common null exception. The comparison for~$A_N^0$ is transferred in Lemma~\ref{mfd:lem-finite-source-comparison}, on the cells at which the subdivision stops.

\smallskip

For an initial cube~$w+\cu_{-k_*}$ of a~$3^{H_1}$-adic subdivision, define the eligible generations by~$I_{k_*}^{a_1,a_2}(N)\coloneqq\{n\in\mathbb N_0:a_1N\leq k_*+H_1n\leq a_2N\}$. Thus the interval~$[a_1N,a_2N]$ always refers to physical triadic depth. Good and bad cells are counted by their generation indices in this set.

Fix~$0<a_1<a_2<1$,~$\theta\in(0,1)$, a subdivision factor~$3^{H_1}$ with~$H_1\geq4$ and a bounded cube~$w+\cu_{-k_*}$. Choose the finite moment orders and the threshold~$\gmcdelta_0$ as in Lemmas~\ref{mfd:lem-finite-good-cell} and~\ref{mfd:lem-rare-tests}, with the rate~$\varkappa$ of Lemma~\ref{mfd:lem-finite-good-cell}(i) large enough for the counting estimate in the proof. These choices precede the tolerance~$\varpi>0$.

For cutoffs~$M\geq N$, a generation~$n\in I_{k_*}^{a_1,a_2}(N)$ of a chain is \textbf{good} if, on its cell of physical depth~$k_*+H_1n$, the good events of Lemma~\ref{mfd:lem-finite-good-cell} hold at both cutoffs and the two separately normalized minimal energies for~$A_N$ satisfy the conclusion of Lemma~\ref{mfd:lem-finite-trace-tests} with tolerance~$\varpi$ for every~$C^\alpha$ trace; otherwise it is \textbf{bad}.

\begin{lemma}[Simultaneous good generations for two cutoffs]
\label{mfd:lem-finite-good-levels}
For every~$\varpi>0$ there are~$N_\varpi,C_\varpi<\infty$ and~$c_\varpi>0$ such that, for each pair~$M\geq N\geq N_\varpi$, outside an event of probability at most~$C_\varpi3^{-c_\varpi N}$, every chain in the~$3^{H_1}$-adic subdivision of~$w+\cu_{-k_*}$ has at most~$\theta N/H_1$ bad generations with index in~$I_{k_*}^{a_1,a_2}(N)$.
\end{lemma}

The comparison for~$A_N^0$ is obtained by the common-weight transfer in Lemma~\ref{mfd:lem-finite-source-comparison}.

\begin{proof}
The two covers obtained from Lemma~\ref{mfd:lem-rare-tests} in the preceding argument and the covers of Lemma~\ref{mfd:lem-finite-good-cell}(i) at the two cutoffs cover each bad generation by events of the form in Lemma~\ref{mfd:lem-finite-good-cell}(i), with a prescribed rate~$\varkappa$, after increasing that rate by a fixed factor. For a cell of physical depth~$k$, the extra tests use stationary indices~$N-k$ and~$N-k+1$ at cutoff~$N$, and~$M-k$ and~$M-k+1$ at cutoff~$M$. Both pairs correspond to the same original field indices~$-k$ and~$1-k$, so their covers have the same field windows at that cell. Their field windows are centered at~$-(k_*+H_1n)$, for~$n\in I_{k_*}^{a_1,a_2}(N)$, and these centers are spaced by~$H_1$. As in the proof of Proposition~\ref{mfd:prop-allchain}, apply Lemma~\ref{l.concentration.rare.intervals} on this interval of field indices, with the events at unselected indices replaced by the empty event. For sufficiently large~$N$, the eligible set has~$(a_2-a_1)N/H_1+O(1)$ elements and its last generation is~$(a_2N-k_*)/H_1+O(1)$. Thus the number of chains ending there is at most~$C_{w+\cu_{-k_*}}3^{da_2N}$, where the fixed root offset is absorbed in~$C_{w+\cu_{-k_*}}$. Choose~$\varkappa$ so that the exponential rate exceeds this growth and take a union bound. If~$\theta\geq a_2-a_1$, use the estimate with a smaller positive value of~$\theta$; it implies the stated bound.
\end{proof}

At physical triadic depth~$k\coloneqq k_*+H_1n$ the minimal energies at the cutoffs~$N$ and~$M$ carry the normalizations~$\ahom_{N-k}/\ahom_N$ and~$\ahom_{M-k}/\ahom_M$ of Lemma~\ref{mfd:lem-local-normalizations}, while the random factor in the scalars~$s_N(z+3^{-k}\cu_0)$ and~$s_M(z+3^{-k}\cu_0)$ is common to both. The ratio~$\rho_{N,M}(k)\coloneqq(\ahom_{N-k}/\ahom_N)/(\ahom_{M-k}/\ahom_M)$ of these two normalizations is therefore, up to the tolerance~$\varpi$ of Lemma~\ref{mfd:lem-finite-trace-tests}, the factor relating the two minimal energies on a good cell of depth~$k$. The next lemma adapts the stopping argument of Proposition~\ref{mfd:prop-density} to two finite cutoffs.

\smallskip

For comparison on an arbitrary cube~$w+\cu_{-k_*}$, first partition it into triadic cubes~$y+\cu_{-\bar k}$ of common depth~$\bar k\coloneqq H_1\lceil k_*/H_1\rceil$. This uses fewer than~$H_1$ preliminary triadic subdivisions and gives finitely many roots, independently of~$N$ and~$M$. Subdivide each of these roots~$3^{H_1}$-adically. All these roots have the same eligible set~$I_{\bar k}^{1/4,3/4}(N)$, and generation~$n$ in each of them has physical depth~$\bar k+H_1n$. For large~$N$, writing~$j\coloneqq\bar k/H_1+n$ identifies this set exactly with~$\{j\in\mathbb N_0:N\leq4H_1j\leq3N\}$. The same preliminary partition and the same eligible set are used in both comparisons below. At these depths~$\rho_{N,M}(\bar k+H_1n)=\rho_{N,M}(H_1j)$. For each fixed cube and all sufficiently large~$N$, its eligible generations therefore correspond to the same set of absolute indices~$j$ with~$N\leq4H_1j\leq3N$.

\begin{lemma}[Comparison of two cutoffs by stopping]
\label{mfd:lem-finite-stopping}
Let~$\vartheta>0$. The geometry, the exponents and the finitely many moment orders can be chosen so that there exists~$\gmcdelta_0>0$ such that, if~$0<\gmcdelta\leq\gmcdelta_0$, then the following holds with constants~$D$,~$H_1$ and~$C$ independent of~$\varpi$,~$c$,~$N$ and~$M$. Let~$w+\cu_{-k_*}$ be a fixed cube and~$g$ a fixed smooth function on~$\partial(w+\cu_{-k_*})$. For every~$\varpi>0$ there exist~$C_\varpi<\infty$ and~$\gamma_\varpi>0$ with the following property. Let~$N\leq M$ be cutoffs and~$c\leq2$. Use the preliminary roots just described, and assume that a deterministic set~$S\subseteq I_{\bar k}^{1/4,3/4}(N)$ satisfies~$\#S\geq\vartheta\#I_{\bar k}^{1/4,3/4}(N)$ and~$\rho_{N,M}(\bar k+H_1n)\leq c$ for every~$n\in S$. Then, outside an event of probability at most~$C_\varpi3^{-\gamma_\varpi N}$,
\begin{equation*}
	\Lambda_{N,w+\cu_{-k_*}}(g)
	\leq c(1+C\varpi3^{H_1D})\Lambda_{M,w+\cu_{-k_*}}(g)
	+C_\varpi3^{-\gamma_\varpi N}\|g\|_{C^2}^2\,.
\end{equation*}
Under these assumptions with~$\rho_{N,M}^{-1}$ in place of~$\rho_{N,M}$, the reverse inequality holds with~$\Lambda_{N,w+\cu_{-k_*}}$ and~$\Lambda_{M,w+\cu_{-k_*}}$ interchanged.
\end{lemma}
\begin{proof}
\emph{Step 1.} We define where the subdivision stops and estimate the energy on the cells at which it stops. Let~$u$ be the~$\mathscr E_M$-harmonic extension of~$g$ to~$w+\cu_{-k_*}$, fix a small~$\zeta>0$, and let~$\nu\coloneqq\Gamma_M(u)+3^{-\zeta N}\,dx$, so that~$\nu(z+r\cu_0)\geq3^{-\zeta N}r^d$ for every cell~$z+r\cu_0$. Use the preliminary partition into roots and their~$3^{H_1}$-adic subdivisions. Apply Lemma~\ref{mfd:lem-finite-good-levels} with~$a_1=1/4$,~$a_2=3/4$ to each of these finitely many roots, and take the union of their exceptional events. The fraction of bad generations is chosen in Step~2. Let~$m\coloneqq\lfloor N/H_1\rfloor$. For large~$N$ the first generation of~$I_{\bar k}^{1/4,3/4}(N)$ is~$\lceil(N/4-\bar k)/H_1\rceil$, and its last generation is~$n_+\coloneqq\lfloor(3N/4-\bar k)/H_1\rfloor$. We follow each root subdivision from generation~$\lceil(N/4-\bar k)/H_1\rceil$ through generation~$n_+$. Stop at the first child~$z+r\cu_0$, along each chain, with parent~$z'+3^{H_1}r\cu_0$ such that its generation~$n$ belongs to~$S$, the cell~$z+r\cu_0$ is good in the sense of Lemma~\ref{mfd:lem-finite-good-levels} and at distance at least~$b_0$ times its side from~$\partial(z'+3^{H_1}r\cu_0)$, and~$\nu(z'+3^{H_1}r\cu_0)\leq3^{H_1D}\nu(z+r\cu_0)$. At such a cell, let~$k\coloneqq\bar k+H_1n$. The local H\"older estimate~\eqref{eq:mfd-finite-good-holder} and the comparison of Lemma~\ref{mfd:lem-finite-trace-tests} give
\begin{equation*}
	\Lambda_{N,z+r\cu_0}(u|_{\partial(z+r\cu_0)})
	\leq \rho_{N,M}(k)\Lambda_{M,z+r\cu_0}(u|_{\partial(z+r\cu_0)})
	+C\varpi \rho_{N,M}(k)\Gamma_M(u)(z'+3^{H_1}r\cu_0)\,.
\end{equation*}
Since~$\Lambda_{M,z+r\cu_0}(u|_{\partial(z+r\cu_0)})$ is bounded by~$\Gamma_M(u)(z+r\cu_0)$, the stopping rule and the disjointness of the cells at which the subdivision stops show that the sum of~$\Lambda_{N,z+r\cu_0}(u|_{\partial(z+r\cu_0)})$ over these cells is at most
\begin{equation*}
	c \mathscr E_M(u)+C\varpi c 3^{H_1D}
	\bigl(\mathscr E_M(u)+3^{-\zeta N}3^{-k_*d}\bigr)\,.
\end{equation*}
This bounds the energy of the competitor of Step~4 on these cells.

\smallskip

\emph{Step 2.} We count the cells at generation~$n_+$ on chains along which the subdivision never stops. The moment estimates of Section~\ref{mfd:sec-local} are uniform in both cutoffs, so, outside an event of probability at most~$C3^{-p\ep N}$, each of the finitely many global constants used below is at most~$3^{\ep N}$. Consider a chain of nested cells along which the subdivision does not stop. At each of its steps at which the generation belongs to~$S$ and the child is good and at distance at least~$b_0$ times its side from the boundary of its parent,~$\nu$ drops by a factor larger than~$3^{H_1D}$. Since~$\nu(w+\cu_{-k_*})$ is controlled by the global energy and~$\nu(z+r\cu_0)\geq3^{-\zeta N}r^d$ for every cell~$z+r\cu_0$ at generation~$n_+$, the number of such steps is at most
\begin{equation*}
	\frac1{DH_1\log3}\Bigl(\log\bigl(3^{\zeta N}\nu(w+\cu_{-k_*})\bigr)
	+\tfrac34 dN\log3+O_{w+\cu_{-k_*}}(1)\Bigr)\,.
\end{equation*}
We choose~$D$ large in terms of~$d$ and~$\vartheta$, and~$\theta$,~$\ep$ and~$\zeta$ small in terms of~$d$,~$\vartheta$ and~$D$. Since~$\#I_{\bar k}^{1/4,3/4}(N)=m/2+O(1)$, each chain has at least~$\vartheta m/2-O(1)$ steps with generation in~$S$. Discarding the steps counted above and the bad steps with generation in~$S$ leaves at least~$\vartheta m/3-O(1)$ steps with generation in~$S$ at which the child is within distance~$b_0$ times its side of the boundary of its parent. Each cell has at most~$C_d3^{H_1(d-1)}$ children within that distance of its boundary, so, for~$H_1$ large, the number of cells at generation~$n_+$ on chains where the subdivision has not stopped is at most~$C_{w+\cu_{-k_*}}3^{H_1(\frac{3d}{4}-\frac{\vartheta}{4})m}$. Thus these remaining cells form an exponentially small fraction of the cells at generation~$n_+$.

\smallskip

\emph{Step 3.} We bound the energies on the remaining cells at generation~$n_+$. Their physical depth is~$\bar k+H_1n_+=3N/4+O(H_1)$, so their side length is comparable to~$3^{-3N/4}$. On such a cell of side~$r$, the bound~\eqref{eq:mfd-2} for H\"older boundary values, together with the bounds~\eqref{eq:mfd-3} on grid cubes, bounds the energy by~$K_N[u]_\alpha^2r^{d-2+2\alpha-\kappa}$. We choose the fixed exponents so that~$2-2\alpha+\kappa<\vartheta/8$. The sum of the energies on these cells is then at most~$C3^{3\ep N}3^{-H_1\vartheta \frac{m}{8}}$, which decays exponentially in~$N$ if~$\ep$ is small.

\smallskip

\emph{Step 4.} We compare the minimal energies with boundary values~$g$. The cells at which the subdivision stops and the remaining cells at generation~$n_+$ form a finite partition of~$w+\cu_{-k_*}$. Joining the~$\mathscr E_N$-harmonic extensions of the traces of~$u$ on the cells of this partition gives a competitor with boundary values~$g$. The Dirichlet principle then gives the asserted upper bound for~$\Lambda_{N,w+\cu_{-k_*}}(g)$, and all excluded events have exponentially small probability. The reverse inequality follows from the same proof with the cutoffs interchanged. The finitely many smaller cutoffs are covered by taking their exceptional event to be the whole probability space and enlarging~$C_\varpi$.
\end{proof}

The same construction compares the energies of solutions with a volume source.

\begin{lemma}[Comparison for volume sources]
\label{mfd:lem-finite-source-comparison}
Assume the hypotheses of Lemma~\ref{mfd:lem-finite-stopping}. There is an event of probability at least~$1-C_\varpi3^{-\gamma_\varpi N}$ on which the following holds. Let~$f\in L^\infty(w+\cu_{-k_*})$ and let~$u_M$ solve~$-\nabla\cdot(A_M\nabla u_M)=f$ on~$w+\cu_{-k_*}$, either with a smooth Dirichlet datum~$g$, in which case we let~$B\coloneqq\|f\|_\infty+\|g\|_{C^2}$, or with zero conormal derivative and mean zero, in which case we assume that~$\int_{w+\cu_{-k_*}}f=0$ and let~$B\coloneqq\|f\|_\infty$. Then there is a trial function~$v_N$, with the same Dirichlet datum in the first case and with mean zero in the second, such that
\begin{equation*}
	\|v_N-u_M\|_\infty\leq C_\varpi3^{-\gamma_\varpi N}B
\end{equation*}
and
\begin{equation*}
	\mathscr E_N^{w+\cu_{-k_*}}(v_N)\leq c(1+C\varpi3^{H_1D})\mathscr E_M^{w+\cu_{-k_*}}(u_M)
	+C_\varpi3^{-\gamma_\varpi N}B^2\,.
\end{equation*}
The same holds with~$A_N$ and~$A_M$ replaced by~$A_N^0=e^{-H}A_N$ and~$A_M^0=e^{-H}A_M$.
\end{lemma}
\begin{proof}
We use the partition of the proof of Lemma~\ref{mfd:lem-finite-stopping}, now with~$\nu\coloneqq\Gamma_M(u_M)+3^{-\zeta N}B^2\,dx$ in place of the measure of that proof. On a good cube~$z+r\cu_0$ whose parent~$z'+3^{H_1}r\cu_0$ satisfies~$\operatorname{dist}(z+r\cu_0,\partial(z'+3^{H_1}r\cu_0))\geq b_0r$, the source term in~\eqref{eq:mfd-finite-good-holder} of Lemma~\ref{mfd:lem-finite-good-cell}(ii) adds to the local comparison the term~$C\varpi \rho_{N,M}(k)s_M(z+r\cu_0)^{-1}r^{d+2}\|f\|_\infty^2$. On the event of the preceding proof, the bounds~\eqref{eq:mfd-3} on grid cubes give~$s_M(z+r\cu_0)^{-1}\leq3^{\ep N}r^{-\kappa}$. Since the side of every cube used in the comparison is at most~$C3^{-\frac14N}$, this source term is bounded by~$C\varpi c3^{-\zeta N}B^2r^d$ if~$\ep+\zeta<(2-\kappa)/4$, so it is absorbed in the lower bound~$3^{-\zeta N}B^2r^d$ for~$\nu(z+r\cu_0)$. The global estimates of Proposition~\ref{mfd:prop-growth} (for cubes of side larger than one, Remark~\ref{mfd:rem-large-cube-growth}) and Proposition~\ref{mfd:cor-neumann-source} bound the energy by~$3^{\ep N}B^2$ and the~$C^\alpha$ norm by~$3^{\ep N}B$, after these constants are added to the finite list of constants required to be at most~$3^{\ep N}$.

\smallskip

On each element of the partition we take the~$A_N$-harmonic extension of the trace of~$u_M$. By the maximum principle, its distance from~$u_M$ is at most twice the oscillation of~$u_M$ on that element. All elements have side at most~$C3^{-\frac14N}$, so this distance is at most~$C3^{-(\frac{\alpha}{4}-\ep)N}B$. In the Neumann case we first subtract the mean. This proves both inequalities.

\smallskip

For~$A_N^0$ and~$A_M^0$, the Dirichlet estimates with source on larger cells are supplied by the whole-cell argument leading to~\eqref{eq:mfd-zero-large-cell-growth}. For the mean-zero Neumann solutions, the proof of Proposition~\ref{mfd:cor-neumann-source} gives the estimates for~$A_N^0$ and~$A_M^0$ on every fixed cube. Its chains use roots of radius at most~$R_*$, so their outer scales do not exceed the cutoff, and the coercivity bound for larger cubes in the proof of Lemma~\ref{mfd:lem-neumann-error} bounds the total energy at the same finitely many moment orders. The local trace comparison transfers separately. On each cube~$z+r\cu_0$ used in the partition,~$r\leq C3^{-N/4}$. Let~$s_j^0(z+r\cu_0)\coloneqq e^{-H(z)}s_j(z+r\cu_0)$ for~$j=N,M$. The variational formula gives, for every trace~$g$,
\begin{align*}
	e^{-\operatorname{osc}_{z+r\cu_0}H}s_j(z+r\cu_0)^{-1}\Lambda_{j,z+r\cu_0}(g)
	&\leq s_j^0(z+r\cu_0)^{-1}\Lambda_{j,z+r\cu_0}^0(g)
	\\ &\leq e^{\operatorname{osc}_{z+r\cu_0}H}s_j(z+r\cu_0)^{-1}\Lambda_{j,z+r\cu_0}(g)\,.
\end{align*}
Both cutoffs use this same weight, and~$s_N^0(z+r\cu_0)/s_M^0(z+r\cu_0)=\rho_{N,M}(k)$. Add~$e^{\|H\|_{C^1(w+\cu_{-k_*})}}$, whose fixed moments are finite, to the constants required to be at most~$3^{\ep N}$. For large~$N$, the oscillation of~$H$ on~$z+r\cu_0$ is at most~$C3^{-N/5}$ on this event. The uniform good-cell trace bound then transfers the finite trace test with an additional tolerance~$C3^{-N/5}$; choose the original tolerance smaller by a fixed factor. Equivalently, a local energy comparison is multiplied by at most~$e^{2\operatorname{osc}_{z+r\cu_0}H}=1+O(3^{-N/5})$. After summing over the cells at which the subdivision stops, its additional error is at most~$C3^{-N/5}\mathscr E_M^0(u_M)$. It is absorbed by the multiplicative tolerance after increasing the deterministic initial cutoff. The additive source and terminal-cell errors retain a positive exponential rate. This proves the two-cutoff comparison for~$A_N^0$ without a convergence assertion uniform over the partial sums of~$H$.
\end{proof}

\begin{proposition}[Almost-sure convergence of minimal energies]
\label{mfd:prop-as-dirichlet}
Let~$w+\cu_{-k_*}$ be a cube and~$g$ a smooth function on~$\partial(w+\cu_{-k_*})$. Then~$\Lambda_{N,w+\cu_{-k_*}}(g)\to\Lambda_{\mathscr E,w+\cu_{-k_*}}(g)$ almost surely as~$N\to\infty$ along the full sequence of integers.
\end{proposition}
\begin{proof}
\emph{Step 1.} We show that the normalizations are close to their limits on most of the eligible generations of the preliminary roots. Let~$I\coloneqq I_{\bar k}^{1/4,3/4}(N)$ and~$\rho_N(k)\coloneqq(\ahom_{N-k}/\ahom_N)/\bar\rho(k)$, and fix~$0<\ep<1$ and a density~$0<\vartheta<1/3$. Suppose that, for infinitely many~$N$, the set~$S_-\coloneqq\{n\in I:\rho_N(\bar k+H_1n)\leq1-\ep\}$ has~$\#S_-\geq\vartheta\#I$. For each such fixed~$N$, the set~$S_-$ is finite. The convergence~$\ahom_{M-k}/\ahom_M\to \bar\rho(k)>0$ at each physical depth~$\bar k+H_1n$ with~$n\in S_-$ therefore gives a cutoff beyond which~$\rho_{N,M}(\bar k+H_1n)\leq1-\ep/2$ for every~$n\in S_-$. Apply Lemma~\ref{mfd:lem-finite-stopping} on~$\cu_0$, using the corresponding set of indices~$\{\bar k/H_1+n:n\in S_-\}$ and a fixed nonconstant affine datum, with~$\varpi$ so small that~$C\varpi3^{H_1D}<\ep/4$, and let~$M\to\infty$ first. Along these~$N$, the minimal energy at cutoff~$N$ is bounded by a fixed factor strictly smaller than one times its limit, up to an error that tends to zero in probability. This contradicts its convergence in probability to a strictly positive limit. The reverse comparison, using~$\rho_{N,M}^{-1}$ and the finite set~$\{n\in I:\rho_N(\bar k+H_1n)\geq1+\ep\}$, excludes in the same way a fraction~$\vartheta$ of large ratios. Hence, for all sufficiently large~$N$, the set~$S\coloneqq\{n\in I:1-\ep\leq \rho_N(\bar k+H_1n)\leq1+\ep\}$ has~$\#S>(1-2\vartheta)\#I$.

\smallskip

\emph{Step 2.} We prove an exponential bound and conclude. For each large~$N$, we apply the two comparisons of Lemma~\ref{mfd:lem-finite-stopping}, with~$S\subseteq I_{\bar k}^{1/4,3/4}(N)$ the set of generation indices found in Step~1, and let~$M\to\infty$. Taking~$\varpi$ small relative to~$\ep$ gives, with constants depending on~$\ep$ and on the datum,
\begin{equation*}
	\P\left[
	|\Lambda_{N,w+\cu_{-k_*}}(g)-\Lambda_{\mathscr E,w+\cu_{-k_*}}(g)|
	>C\ep\Lambda_{\mathscr E,w+\cu_{-k_*}}(g)+C_\ep3^{-\gamma_\ep N}
	\right]\leq C_\ep3^{-\gamma_\ep N}\,.
\end{equation*}
The right side is summable in~$N$, so the Borel--Cantelli lemma for~$\ep=1/j$ with~$j\in\N$ gives the almost-sure convergence.
\end{proof}

The same argument applies to the following quantities on a fixed cube~$z+\cu_n$, computed for the coefficient~$A_N$ or for~$A_N^0$. We write~$\Xi_N$ for any one of them:
\begin{enumerate}[label=(\roman*),leftmargin=*]
\item the minimal energy~$\Lambda_{N,z+\cu_n}(g)$ with a fixed smooth boundary value~$g$;
\item the quadratic form~$\langle f,G_N^{z+\cu_n}f\rangle$ of~\eqref{eq:mfd-quadratic-response} for a fixed bounded source~$f$;
\item the quadratic form~$\sup_{v\in H^1(z+\cu_n)}\{2\int_{z+\cu_n}fv-\mathscr E_N^{z+\cu_n}(v)\}$ for a fixed bounded source~$f$ with~$\int_{z+\cu_n}f=0$;
\item the Neumann energy~$3^{nd}\,p\cdot\a_*^{-1}(z+\cu_n;A_N)p$ of Section~\ref{mfd:sec-environment}, for the constant flux given by a fixed vector~$p$.
\end{enumerate}

\begin{proposition}[Summable error probabilities]
\label{mfd:prop-as-response-bank}
There exists~$\gmcdelta_0>0$ such that, if~$0<\gmcdelta\leq\gmcdelta_0$, then~$\Xi_N$ converges in probability to a limit~$\Xi$, and there is~$C<\infty$ such that, for every~$\ep>0$, there exist~$C_\ep<\infty$,~$c_\ep>0$ and a deterministic~$N_\ep$ such that, for every~$N\geq N_\ep$,
\begin{equation}
	\P\bigl[|\Xi_N-\Xi|>
	C\ep \Xi+C_\ep3^{-c_\ep N}\bigr]
	\leq C_\ep3^{-c_\ep N}\,.
	\label{eq:mfd-summable-response}
\end{equation}
The finitely many moment orders and the threshold~$\gmcdelta_0$ do not depend on~$\ep$. For the coefficients~$A_N^0$, the constants are uniform over deterministic translations of~$z+\cu_n$. In particular, for every countable family of such quantities, all of them converge almost surely along the full sequence of cutoffs.
\end{proposition}
\begin{proof}
For the minimal energies of~$A_N$, the bound~\eqref{eq:mfd-summable-response} was proved in Proposition~\ref{mfd:prop-as-dirichlet}. For~$A_N^0$, Proposition~\ref{mfd:prop-21} and~\eqref{eq:mfd-weighted-boundary-minimum}, with the weight~$e^{-H}$ and the same limit form on the enlarged cube, give convergence in probability. The common-weight transfer in Lemma~\ref{mfd:lem-finite-source-comparison} supplies the two-cutoff comparison. Fix before~$\gmcdelta$ all moment orders needed for the rare tests, the global estimates and the products below; their number and values are independent of the tolerance.

\smallskip

For a volume source~$f$, let~$u_M$ be the solution and write~$\Xi_M\coloneqq \mathscr E_M(u_M)=\int_{z+\cu_n}fu_M$. Let~$c$ and~$\varpi$ be as in Lemma~\ref{mfd:lem-finite-source-comparison}, and insert~$v_N/(c(1+C\varpi3^{H_1D}))$ in the variational formula for~$\Xi_N$. For~$1/2\leq c(1+C\varpi3^{H_1D})\leq3$, the trial-energy and sup-distance bounds give
\begin{equation*}
	\Xi_M\leq c(1+C\varpi3^{H_1D})\,\Xi_N+C_\ep3^{-\gamma N}\|f\|_\infty^2\,.
\end{equation*}
Interchanging the cutoffs gives the reverse comparison. Step~1 of Proposition~\ref{mfd:prop-as-dirichlet} supplies a set of eligible generations, of density larger than~$1/3$, on which the normalization ratio is close to one. For every fixed tolerance and all large~$N$, the resulting comparisons therefore hold with~$c(1+C\varpi3^{H_1D})\leq1+C\ep$, uniformly in all sufficiently large~$M$. The exceptional events are independent of the bounded source, and~$\gamma>0$ depends on the tolerance.

\smallskip

We first record the moment bound used in case~(iii). The mean-zero coercivity constant~$Q_{z+\cu_n,N}$ in the proof of Lemma~\ref{mfd:lem-neumann-error} gives~$0\leq \Xi_N\leq Q_{z+\cu_n,N}\|f\|_{L^2(z+\cu_n)}^2$. Only upper moments of this constant are needed. On a fixed larger triadic cube~$z+\cu_n$ with~$n>0$, this constant can be computed exactly. The fields of~$\Theta_{-n,z}\omega$ have the original joint law, and expansion of the logarithms gives
\begin{equation*}
	A_N(z+3^ny)=s_N(z+\cu_n)A_{N+n}^{\Theta_{-n,z}\omega}(y)\,.
\end{equation*}
The energy scales by~$3^{n(d-2)}s_N(z+\cu_n)$ and the squared~$L^2$ norm by~$3^{nd}$. Thus the physical~$L^2$ coercivity constant is~$3^{2n}s_N(z+\cu_n)^{-1}Q_{N+n}^{\Theta_{-n,z}\omega}$. For the inhomogeneous~$H^{\nf34}$ norm it is at most
\begin{equation*}
	\max\{3^{2n},3^{\frac12n}\}\,s_N(z+\cu_n)^{-1}Q_{N+n}^{\Theta_{-n,z}\omega}\,.
\end{equation*}
By~\eqref{e.annealed.ordering.ratio} and Lemma~\ref{mfd:lem-infrared}, both~$s_N(z+\cu_n)$ and its inverse have every fixed finite moment uniformly in~$N$, with constants depending on~$z+\cu_n$. Take~$Q_{N+n}^{\Theta_{-n,z}\omega}$ at twice the requested order before choosing~$\gmcdelta$, and apply H\"older's inequality; no independence is needed. For~$A_N^0$, the bound~$\mathscr E_N^{0,z+\cu_n}\geq e^{-\sup_{z+\cu_n}H}\mathscr E_N^{z+\cu_n}$ introduces one more finite-moment factor. For sides at most one use Lemma~\ref{mfd:lem-coercivity}; the finitely many cutoffs below the depth of the cube are also bounded by~$\|A_N^{-1}\|_{L^\infty(z+\cu_n)}$ and the ordinary cube Poincar\'e inequality.

\smallskip

Now compare the cutoffs~$N$ and~$M\geq N$, with~$N$ fixed and large. These upper moments make the multiplicative error small in probability as~$\ep\downarrow0$, and the two exceptional probabilities vanish as~$N\to\infty$. This proves the Cauchy property in probability in case~(iii), without a Neumann limit form. In case~(ii), the limit is the quadratic form of the unique Dirichlet limit form. Letting~$M\to\infty$ with~$N$ fixed gives~\eqref{eq:mfd-summable-response} in both cases.

\smallskip

For case~(iv), take the mean-zero source~$f_{z+\cu_n,h}$ constructed in the proof of Lemma~\ref{mfd:lem-neumann-error}. On~$z+\cu_n$ it contains the physical factor~$3^{-n}$ and satisfies~$\|f_{z+\cu_n,h}\|_\infty\leq C_{z+\cu_n,p}h^{-1}$. The comparison of~$L_{w,R,p}$ and~$L_{w,R,p,h}$ there gives an~$L^q$ approximation error at most~$C_{q,z+\cu_n,p}h^{\nf14}$ and uniformly bounded~$L^q$ moments of the energies~$Y_{N,h}$ defined there, at the selected orders. After the preceding change of variables, the Neumann energy is~$3^{nd}s_N(z+\cu_n)^{-1}$ times the Neumann energy on the unit cube; the factor~$3^{-n}$ in the source gives the same scaling for~$Y_{N,h}$. That proof also supplies the coercivity constant~$Q_{z+\cu_n,N}$ for both~$A_N$ and~$A_N^0$ on the same fixed cube, with a uniform polynomial tail~$\P[Q_{z+\cu_n,N}>T]\leq B_{z+\cu_n}T^{-p_0}$, with~$p_0>0$. Here~$B_{z+\cu_n}$ is the constant~$M_{w,R}$ of that proof for the cube~$z+\cu_n$. For the required larger triadic cubes, its construction uses the law-preserving shift to cutoff~$N+n$, the fields~$\g_1,\dots,\g_n$ and the two restrictions of~$H$. The unit-cube first-moment bound already suffices for~$p_0=1/2$. This transfer fixes its disorder threshold before the tolerance.

\smallskip

Let~$Y_{N,h}$ be the quadratic form~(iii) for the smoothed volume source~$f_{z+\cu_n,h}$. The comparison of~$L_{w,R,p}$ and~$L_{w,R,p,h}$ gives
\begin{equation*}
	|\Xi_N-Y_{N,h}|
	\leq C_{z+\cu_n,p}h^{\nf14}Q_{z+\cu_n,N}^{\nf12}\Xi_N^{\nf12}
	+C_{z+\cu_n,p}h^{\nf12}Q_{z+\cu_n,N}\,.
\end{equation*}
For fixed~$h$, case~(iii) and the uniform moment approximation prove that~$(\Xi_N)_N$ is Cauchy in probability. For the rate, use the two-cutoff comparison directly, with~$h=3^{-\gamma N/4}$. Outside the source-independent comparison events and the two events~$\{Q_{z+\cu_n,N}>3^{\gamma N/32}\}$ and~$\{Q_{z+\cu_n,M}>3^{\gamma N/32}\}$, absorb the mixed terms by~$x_1x_2\leq\ep x_1^2+C_\ep x_2^2$. Applying this absorption in both directions, and taking~$\ep\leq1/2$, yields
\begin{equation*}
	|\Xi_N-\Xi_M|
	\leq C\ep \Xi_N+
	C_\ep
	\left(3^{-\gamma N}h^{-2}+3^{\gamma N/32}h^{\nf12}\right)\,.
\end{equation*}
The two deterministic remainders decay as~$3^{-\gamma N/2}$ and~$3^{-3\gamma N/32}$. The probability of either coercivity event is at most~$B_{z+\cu_n}3^{-p_0\gamma N/32}$, uniformly for every sufficiently large~$M$. Hence a common positive rate is~$c_\ep\coloneqq\min\{3\gamma/32,p_0\gamma/32\}$. Let~$M\to\infty$ with~$N$ fixed, using convergence in probability and then a further subsequence converging almost surely. The resulting bound first has~$C\ep \Xi_N$; choosing the tolerance smaller absorbs this term against the limit~$\Xi$ and gives~\eqref{eq:mfd-summable-response}. No centered estimate on the physical cube is used in this comparison.

\smallskip

The finite list of moment orders is chosen before~$\gmcdelta$. Only the comparison geometry, initial field windows, constants and the positive rate~$\gamma$ depend on the tolerance; the coercivity order does not. For~$A_N^0$, translate the fixed cube and its boundary datum or volume source together, and leave the constant flux vector unchanged. Every finite-cutoff quantity~$\Xi_N$ on the translated cube equals the corresponding quantity on the original cube in the translated fields. Choose its limit by the same translation. Stationarity therefore preserves the joint law of~$\Xi_N$ and its limit, so the same constants and deterministic initial cutoff work for every deterministic translate of this fixed cube. They may depend on its side and the fixed data. The final assertion follows by Borel--Cantelli at~$\ep=1/j$ with~$j\in\N$, and a countable intersection.
\end{proof}

\subsection{Almost-sure bounds uniform in the cutoff}
\label{mfd:sec-pathwise-bounds}

In this subsection we prove that the coarse-grained ellipticity bounds (Lemma~\ref{mfd:lem-as-coarse}) and the regularity estimates (Lemma~\ref{mfd:lem-as-regularity}) hold with one finite random constant for all cutoffs, on one event of probability one. On this event, the Green operators converge in norm, the Dirichlet forms Mosco-converge and the killed resolvents converge uniformly (Proposition~\ref{mfd:prop-as-forms}). Let~$\theta>0$ be the proportion chosen in the proof of Lemma~\ref{mfd:lem-as-coarse}. On cells of depth at most~$\theta N$, Proposition~\ref{mfd:prop-as-response-bank} compares the finite-cutoff energies with their limits, with error probabilities summable over the cells and over~$N$. For the deeper cells, the moment bounds imply that, almost surely, the random constants in the estimates are smaller than every fixed exponential in~$N$ once~$N$ is large. The margins in the spatial exponents absorb this growth.

\begin{lemma}[Coarse-grained ellipticity bounds uniform in the cutoff]
\label{mfd:lem-as-coarse}
Let~$w+\cu_n$ be a cube and let~$s\in(0,1]$. If~$\gmcdelta$ is sufficiently small, depending on~$s$ and on the H\"older exponent in~\eqref{eq:mfd-2}, then, almost surely, for~$q\in\{1,2\}$,
\begin{equation*}
	\sup_N\bigl(\Lambda_{s,q}(w+\cu_n;A_N)
	+\lambda_{s,q}(w+\cu_n;A_N)^{-1}\bigr)<\infty\,.
\end{equation*}
In particular, applied with~$s=\nf18$ and~$q=1$, and with~$s=\upsilon/2$ and~$q=2$, the coercivity bound~\eqref{eq:mfd-1} and the bound~\eqref{eq:mfd-2} for H\"older boundary values on~$w+\cu_n$ hold with one finite random constant for all~$N$. The same conclusions hold with~$A_N^0$ in place of~$A_N$.
\end{lemma}
\begin{proof}
It suffices to consider a unit cube.

\smallskip

\emph{Step 1.} We bound the finite-cutoff energies of the cells of depth at most~$\theta N$ by their limits. Denote by~$\mathcal R_m(z)$ the finite family consisting of the minimal energies with the affine boundary values~$x\mapsto e\cdot x$ and the Neumann energies for the fluxes~$e$, for~$e\in\{e_i,e_i+e_j\}$, computed for~$A_m^0$ in the environment shifted to~$z$, at cutoff~$m$. We apply~\eqref{eq:mfd-summable-response} with one fixed small tolerance to this family on every grid cube of depth~$k\leq\theta N$, whose cutoff after rescaling to the unit cube is~$m=N-k$. There are at most~$C3^{dk}$ such cubes. We choose~$\theta>0$ so small that~$d\theta<c(1-\theta)/2$, where~$c$ is the exponent in~\eqref{eq:mfd-summable-response} at this tolerance; then the sum over these cubes of the exceptional probabilities is summable in~$N$. Hence, almost surely, for all sufficiently large~$N$, each of these finite-cutoff energies is bounded by~$C\cdot(1+\mathcal R(\Theta_{k,z}\omega))$, where~$\mathcal R$ is its limit. All the entries are nonnegative, and polarization recovers the coarse-grained matrices from this finite family.

\smallskip

\emph{Step 2.} We bound the coarse-grained matrices of the cells of depth at most~$\theta N$. For a cube~$z+3^{-k}\cu_0$, separating the scalar~$s_N(z+3^{-k}\cu_0)$ of Lemma~\ref{mfd:lem-local-normalizations} gives
\begin{equation*}
	A_N(z+3^{-k}y)=
	s_N(z+3^{-k}\cu_0)\exp\Bigl(\Bigl[H+\sum_{j=0}^{k-1}\g_{-j}\Bigr](z+3^{-k}y)
	-\Bigl[H+\sum_{j=0}^{k-1}\g_{-j}\Bigr](z)\Bigr)
	A_{N-k}^{0,\Theta_{k,z}\omega}(y)\,.
\end{equation*}
The deterministic ratio in~$s_N(z+3^{-k}\cu_0)$ and its inverse are at most~$e^{C\tau^2k}$. Hence the factors in front of~$A_{N-k}^{0,\Theta_{k,z}\omega}$ and their inverses are bounded by
\begin{equation*}
	e^{C\tau^2k}\exp\Bigl(\Bigl|H(z)+\sum_{j=0}^{k-1}\g_{-j}(z)\Bigr|
	+\operatorname{osc}_{z+3^{-k}\cu_0}\Bigl[H+\sum_{j=0}^{k-1}\g_{-j}\Bigr]\Bigr)\,,
\end{equation*}
whose~$L^p$ norm is at most~$C_pe^{C_p\gmcdelta^2 k}$, by the independence of the fields and the bound on their gradients weighted by their scales. The limits~$\mathcal R$ have the same moment bounds as their approximants. We choose~$\zeta<s/4$, a moment order~$p$ with~$p\zeta>2d$, and then~$\gmcdelta$ small enough. A union bound over all~$k$ and~$z$ gives a random~$K<\infty$ such that, for every grid cube of depth~$k$, the product of the expression in the last display and~$1+\mathcal R(\Theta_{k,z}\omega)$ is at most~$K3^{\zeta k}$. Since both coarse-grained matrices are monotone in the coefficient,~$\a(z+3^{-k}\cu_0;A_N)+\a_*^{-1}(z+3^{-k}\cu_0;A_N)\leq CK3^{\zeta k}$ on the cells of depth at most~$\theta N$, with the same~$K$ for all sufficiently large~$N$.

\smallskip

\emph{Step 3.} We bound the matrices of the cells of depth greater than~$\theta N$ and sum over all depths. For~$\theta N<k\leq N$, the proof of Lemma~\ref{mfd:lem-extension}, applied only to the matrices of the cells, gives the bound~$K_N3^{\zeta k}$, where the moments of~$K_N$ of the finitely many orders fixed at the start of Section~\ref{sec:mfd-uniqueness} are bounded uniformly in~$N$. By the Borel--Cantelli lemma,~$K_N\leq3^{\zeta_1 N}$ for all large~$N$, for every fixed~$\zeta_1>0$. The contribution of these cells to either sum defining~$\Lambda_{s,q}$ or~$\lambda_{s,q}$ is at most~$CK_N3^{-(s-\zeta)\theta N}$, which is bounded if~$\zeta_1<(s-\zeta)\theta$. The cells below the smallest field scale contribute at most~$C3^{-sN}(\sup_{w+\cu_n}A_N+\sup_{w+\cu_n}A_N^{-1})$, which is bounded almost surely once~$C_d\gmcdelta+C_1\gmcdelta^2<s\log3$, with the constants of~\eqref{eq:mfd-6} at~$p=1$, by Lemma~\ref{mfd:lem-extremes} and Markov's inequality. The cells of depth at most~$\theta N$ contribute at most~$CK\sum_{k\geq0}3^{-(s-\zeta)k}$. For~$q=1$, weighted Cauchy--Schwarz bounds the square of the sum of the square roots of the matrices by this sum of matrices; for~$q=2$, the definition is already a sum of matrices, with the smaller weight~$3^{-2sk}$. Finally, the coercivity bound and the bound~\eqref{eq:mfd-2} for H\"older boundary values follow from these bounds as in the proofs of Lemmas~\ref{mfd:lem-coercivity} and~\ref{mfd:lem-extension}.
\end{proof}

\begin{lemma}[Regularity uniform in the cutoff]
\label{mfd:lem-as-regularity}
Let~$0<\alpha_1<1$ and~$d-1<t_1<d$ be strictly smaller than the exponents fixed in the moment estimates of Proposition~\ref{mfd:prop-growth} and Proposition~\ref{mfd:cor-neumann-source}. Then, for every fixed cube~$w+R\cu_0$, almost surely the estimates of Proposition~\ref{mfd:prop-growth}, extended to larger cubes by Remark~\ref{mfd:rem-large-cube-growth}, and Proposition~\ref{mfd:cor-neumann-source} with the exponents~$\alpha_1$ and~$t_1$ hold on~$w+R\cu_0$ with one finite random constant for all~$N$.
\end{lemma}
\begin{proof}
We run the iteration of Lemma~\ref{l.iteration.lemma.GMC} using only the primitive conditions on~$\mathsf F$,~$\mathsf M$ and~$\mathsf R$, and the accumulated error~$\mathbf G$, defined in Section~\ref{mfd:sec-branches}. The normalized tests~(a)--(c) of~$\mathcal G_N^{\rm det}(z+3^{-k}\cu_0)$ are not inputs to this iteration. Its reference ratios follow from~\eqref{e.ratio.of.bs}; the actual coefficient's initial energy and ellipticity bounds are supplied separately by Lemma~\ref{mfd:lem-as-coarse}. In this proof,~$\theta$ denotes a smaller proportion, chosen in Steps~2 and~5. We transfer the primitive prefixes above the scale~$3^{-\theta N}$ to their limits (Steps~1--5), and use the uniform moment bounds below that scale (Step~6).

\smallskip

\emph{Step 1.} We pass the primitive quantities to the limit. In~$\mathsf R$ and~$\mathbf G$, the minimal energies with affine boundary values and the Neumann energies at stationary scale~$n$ are computed for~$\ahom_n^{-1}\a_n$. After translation and rescaling to the unit cube they are energies of~$A_n^0$, with no cell scalar and no factor~$e^H$. Define the limits of these energies by Proposition~\ref{mfd:prop-as-response-bank}, and use them to define the limiting primitive quantities. The pure-field parts are taken in the same reindexed environment. These energies converge along the full sequence of cutoffs by Proposition~\ref{mfd:prop-as-response-bank}, so the proof of Lemma~\ref{mfd:lem-prefix-limit}(i) gives their finite-interval~$L^p$ limits along that full sequence. Geometric tail bounds pass by Fatou, and the conditional-expectation argument preceding Lemma~\ref{mfd:lem-witness} gives their limiting band bounds. The matrices computed for~$A_N/s_N(z+3^{-k}\cu_0)$ include the increments of~$H$ about the center and are not identified with these energies of~$A_n^0$. Their almost-sure uniform ellipticity bounds are those of Lemma~\ref{mfd:lem-as-coarse}, and no convergence of their normalized tests is used here.

\smallskip

\emph{Step 2.} We compare the finite-cutoff primitive quantities with their limits. Fix~$\varpi>0$ below the margins between the lower and upper ramp thresholds. Truncate only the primitive quantities at relative depth~$\lfloor\theta_1N\rfloor$, and consider intervals of length at most~$\lfloor\theta N\rfloor$ at all the centers needed in the proof of Proposition~\ref{mfd:cor-neumann-source}, including the face meshes and the centers of the initial balls. For one of these energies at physical depth~$j$, the remaining own cutoff is~$N-j$; its limit is taken in that same translated and reindexed environment. Thus all the retained energies are unit-cube energies of~$A_{N-j}^0$, not energies of~$A_N/s_N(z+3^{-k}\cu_0)$. There are at most~$CN^C3^{C_d(\theta+\theta_1)N}$ of them, with remaining cutoffs at least~$[1-C(\theta+\theta_1)]N-C$. By~\eqref{eq:mfd-summable-response} and a union bound, the probability that some retained energy differs from its limit by more than~$C\varpi$ times the limit plus the exponentially small additive error is at most
\begin{equation*}
	C_\varpi N^C3^{C_d(\theta+\theta_1)N-c_\varpi[1-C(\theta+\theta_1)]N}\,.
\end{equation*}
The discarded primitive tails, at finite cutoff and in the limit, have~$L^p$ norm at most~$C_p3^{-a\lfloor\theta_1N\rfloor}$. The probability that one of them exceeds~$\varpi$ is at most~$C_\varpi N^C3^{C_d\theta N-ap\theta_1N}$. Choose~$\theta_1$ small enough for the first exponent to be negative, and then~$\theta$ small enough for both exponents to be negative. Both bounds are summable in~$N$. Borel--Cantelli gives these comparisons simultaneously for all sufficiently large cutoffs.

\smallskip

\emph{Step 3.} We compare the finite and limiting primitive prefixes. The relative error of the unit-cube energies controls the continuous piecewise linear functions: polarization gives~$|\mathcal J_N-\mathcal J_\infty|\leq C\varpi(1+\mathcal J_\infty)$ on each retained cube at its own cutoff, and weighting gives the same inequality for the condition on~$\mathsf R$. A bounded continuous piecewise linear function of this quantity changes by at most~$C\varpi$, and the square-root terms of the sums over intervals of scales change by at most~$C\sqrt\varpi$. The terms of~$\mathbf G$ involving the fields~$\g_j$ at a retained depth coincide once their finite interval of field indices lies below~$N$, and their omitted tails satisfy the geometric bound of Step~2. Consequently, almost surely for all large~$N$, every retained sum over an interval of scales at finite cutoff is bounded by the corresponding limiting sum plus~$C\sqrt\varpi$ per scale.

\smallskip

\emph{Step 4.} We bound, uniformly over the choices of center and initial radius, the number of scales after which the limiting sums over intervals of scales satisfy the thresholds. For a fixed choice~$\mathsf c$ of center and initial radius, let~$m_{\mathsf c}$ be the least integer~$b$ such that the limiting primitive prefixes at~$\mathsf c$ of every length~$D>b$ satisfy the stricter thresholds of the iteration. Apply Lemma~\ref{mfd:in-prefix} to~$\mathsf Z$ and~$\mathbf G$ at a smaller finite-cutoff threshold and pass the bounds to the limit by Lemma~\ref{mfd:lem-prefix-limit}(ii). Summing over the interval lengths gives~$\P(m_{\mathsf c}>b)\leq Ce^{-\varkappa b}$, uniformly in~$\mathsf c$. At target depth~$n$, the meshes of centers and initial radii used in the proof of Proposition~\ref{mfd:cor-neumann-source} involve at most~$C(n+1)^C3^{C_dn}$ choices of centers, initial radii and faces, edges and corners. We fix a small~$\zeta_2>0$ before decreasing~$\gmcdelta$, and choose the rate~$\varkappa$ so that~$\varkappa\zeta_2>C_d\log3+2$. A union bound over the choices~$\mathsf c$ at depth~$n$ gives
\begin{equation*}
	\P\Bigl[\max_{\mathsf c} m_{\mathsf c}>\zeta_2 n+b\Bigr]
	\leq C(n+1)^C3^{C_dn}e^{-\varkappa(\zeta_2 n+b)}\,.
\end{equation*}
Summing over~$n$ shows that, almost surely, there is one finite~$m_\omega$ such that~$m_{\mathsf c}\leq\zeta_2n+m_\omega$ for every choice~$\mathsf c$ at every depth~$n$.

\smallskip

\emph{Step 5.} We transfer these bounds to finite cutoffs and prove the estimates at radii~$r\geq 3^{-\theta N}$. We take~$\varpi$ so small that the loss per scale of Step~3 is below the strict margin of the iteration. Then, for finite cutoffs and depths~$n\leq\theta N$, the same numbers~$m_{\mathsf c}$ apply with slightly weaker thresholds. The~$\mathsf Z$ prefix bounds the number of failures of the primitive good-scale conditions, since a failed upper threshold contributes one to~$\mathsf Z$; the~$\mathbf G$ prefix bounds the accumulated comparison error. These are the two inputs of~\eqref{e.after.good.scales.applied}. The reference-ratio input is~\eqref{e.ratio.of.bs}, derived from the same primitive conditions and annealed ordering. Thus the iteration applies without imposing~(a)--(c) at these cells. Its initial averages and energy factors are bounded separately as follows. The estimate of Lemma~\ref{mfd:thm-fold} transfers the original field bounds to the reflected coefficients near faces, edges and corners. The reflection across faces is applied at most~$d+1$ times, and the numbers~$m_{\mathsf c}$ of these applications add up, so their total contribution is bounded by a common random constant times~$3^{C\zeta_2 n}$.

\smallskip

We next show that the averages~$\bar s_i$ at the initial radius of each application are bounded by a random constant times a small negative power of the initial radius. Let~$\theta_0$ be the value of~$\theta$ in the proof of Lemma~\ref{mfd:lem-as-coarse}, so that the cells of depth at most~$\theta_0N$ were compared with their limits there, and decrease~$\theta$ to at most~$\theta_0/2$. For an initial radius at depth~$k\leq\theta N$, we split its coarse-grained sum at descendant depth~$\theta_0N$, with~$s=\nf18$ for the coercivity constant and~$s=\upsilon/2$ for the constant in~\eqref{eq:mfd-2}, the two orders of Lemma~\ref{mfd:lem-as-coarse}. For~$0<\zeta<s$, the part with descendant depth at most~$\theta_0N$ is at most~$CK3^{\zeta k}$, and the part with descendant depth greater than~$\theta_0N$ is at most
\begin{equation*}
	CK_N3^{\zeta k}3^{-(s-\zeta)(\theta_0N-k)}
	\leq C3^{\zeta k}3^{(\zeta_1-(s-\zeta)\theta_0/2)N}\,,
\end{equation*}
where~$K_N\leq3^{\zeta_1 N}$ for large~$N$, as in the proof of Lemma~\ref{mfd:lem-as-coarse}. We choose~$\zeta_1<(s-\zeta)\theta_0/2$; the contribution below the smallest field scale satisfies the same bound by Lemma~\ref{mfd:lem-extremes}. The average at its own scale involves only the value of~$H+\sum_{j<k}\g_{-j}$ at the center and the deterministic ratio~$\ahom_{N-k}/\ahom_N$. Hence the energy in the initial ball, the coercivity constant and the constant in~\eqref{eq:mfd-2} are bounded by a random constant times a small negative power of the initial radius.

\smallskip

We apply the iteration of Lemma~\ref{l.iteration.lemma.GMC} with exponents strictly larger than~$\alpha_1$ and~$t_1$, and choose~$\zeta_2$ and~$\zeta$ so small that the margins in the exponents absorb both losses. The iteration gives one random constant for the oscillation and energy bounds at all radii~$r\geq 3^{-\theta N}$.

\smallskip

\emph{Step 6.} We treat the radii~$r\leq 3^{-\theta N}$. We choose exponents~$\alpha_0>\alpha_1$ and~$t_0>t_1$ for which the Morrey bounds at all radii of Proposition~\ref{mfd:prop-growth} (on cubes of side larger than one, Remark~\ref{mfd:rem-large-cube-growth}) and Proposition~\ref{mfd:cor-neumann-source} hold. Their constants satisfy~$K_N\leq3^{\zeta_1 N}$ almost surely for all large~$N$, for every fixed~$\zeta_1>0$. For~$r\leq 3^{-\theta N}$, we have~$K_Nr^{\alpha_0-\alpha_1}\leq3^{(\zeta_1-\theta(\alpha_0-\alpha_1))N}$ and~$K_Nr^{t_0-t_1}\leq3^{(\zeta_1-\theta(t_0-t_1))N}$. We choose~$\zeta_1$ smaller than~$\theta(\alpha_0-\alpha_1)$ and~$\theta(t_0-t_1)$; then these quantities are bounded. In the oscillation estimates at radii larger than~$3^{-\theta N}$, point values are compared with averages over balls of radius~$3^{-\theta N}$, using the estimate at small radii just proved.
\end{proof}

\begin{proposition}[Almost-sure convergence of energies and resolvents]
\label{mfd:prop-as-forms}
Almost surely, for every~$z\in\mathbb Q^d$ and~$n\in\mathbb Z$, the operators~$G_N^{z+\cu_n}$ converge to~$G^{z+\cu_n}$ in operator norm on~$L^2(z+\cu_n,dx)$, the forms~$\mathscr E_N^{z+\cu_n}$ Mosco-converge to~$\mathscr E^{z+\cu_n}$, and, for every~$f\in C(\overline{z+\cu_n})$ and~$\lambda>0$,
\begin{equation*}
	\|R^{z+\cu_n}_{N,\lambda}f-R^{z+\cu_n}_\lambda f\|_{C(\overline{z+\cu_n})}\longrightarrow0\,.
\end{equation*}
On each fixed cube, the~$C^{\nf14}(\overline{z+\cu_n})$ norms of~$R^{z+\cu_n}_{N,\lambda}f$ are bounded by one finite random constant, uniformly in~$N$, in~$\lambda$ in a compact subset of~$(0,\infty)$, and in~$f$ with~$\|f\|_\infty\leq1$.
\end{proposition}
\begin{proof}
Intersect the events of Lemmas~\ref{mfd:lem-as-coarse} and~\ref{mfd:lem-as-regularity}, Proposition~\ref{mfd:prop-chaos-growth}, and Proposition~\ref{mfd:prop-as-response-bank} for a countable dense family of bounded sources and their pairwise sums. Coercivity makes the images of the~$L^2$ unit ball under~$G_N^{z+\cu_n}$ bounded in~$H^{\nf34}$ uniformly in~$N$, and polarization gives convergence of the matrix elements of~$G_N^{z+\cu_n}$ on the dense family. The finite-rank argument of Proposition~\ref{mfd:prop-killed-inverse} and the variational formula~\eqref{eq:mfd-quadratic-response} give norm and Mosco convergence to the limits already identified in probability.

\smallskip

All the constants in Proposition~\ref{mfd:prop-uniform-resolvent} are now bounded along the full sequence. Its variational and Arzel\`a--Ascoli argument therefore proves uniform convergence and the H\"older bound.
\end{proof}

\subsection{Almost-sure convergence of the path laws}
\label{mfd:sec-as-paths}

In this subsection we prove that the quenched path laws converge almost surely, uniformly in the starting point on compact sets (Proposition~\ref{mfd:prop-as-quenched}), and we complete the proof of Theorem~\ref{t.scaling.limit}(i). In an environment on which the killed resolvents converge uniformly (Proposition~\ref{mfd:prop-as-forms}) and the limiting kernel is continuous and conservative (Proposition~\ref{mfd:prop-limit-properties}), the uniform convergence transfers the exit estimates of the limiting kernel to all large cutoffs, which gives tightness of the path laws. Tightness and the convergence of the discounted occupation integrals then identify the limit.

\smallskip

\begin{proposition}[Almost-sure quenched convergence]
\label{mfd:prop-as-quenched}
Almost surely, for every compact~$F\subset\Rd$,
\begin{equation*}
	\lim_{N\to\infty}\sup_{x\in F}
	d_{\mathrm{weak}}(K^\omega_{N,x},K^\omega_x)=0\,.
\end{equation*}
\end{proposition}
\begin{proof}
Intersect the full-probability events of Proposition~\ref{mfd:prop-as-forms}, Proposition~\ref{mfd:prop-limit-properties} and Lemma~\ref{mfd:lem-killing} over the countable family of cubes, and fix an environment in this event. The killed resolvents of~$K^\omega_{N,x}$ converge to those of~$K^\omega_x$, uniformly on each closed cube.

\smallskip

We first obtain path tightness in this fixed environment. The image of a compact starting set under the continuous kernel~$x\mapsto K^\omega_x$ is a compact family of laws on continuous path space. Its paths have infinite lifetime, so their probability of leaving a large cube before a fixed time, or of moving a fixed distance in a sufficiently short time, tends to zero uniformly in the starting point. Uniform convergence of the killed resolvents transfers these bounds to all sufficiently large cutoffs. For compact containment, compare~$R^{w+\cu_n}_{N,1}1$ with its limiting counterpart on a sufficiently large cube~$w+\cu_n$ with~$w\in\mathbb Q^d$ and~$n\in\mathbb Z$: the identity~$1-R^{w+\cu_n}_{N,1}1(x)=\int e^{-\tau_{w+\cu_n}}\,dK^\omega_{N,x}$ bounds exits before a fixed time. For short exits, fix cubes~$z+r\cu_0\Subset y+\cu_\ell$, with~$y\in\mathbb Q^d$ and~$\ell\in\mathbb Z$, and a compact starting set~$F_0\subset z+r\cu_0$. Choose a continuous~$f$ with values in~$[0,1]$, supported in~$y+\cu_\ell$, zero near~$F_0$ and one near~$\partial(z+r\cu_0)$. Continuity of the limiting path kernel, whose image of~$F_0\cup\partial(z+r\cu_0)$ is compact, gives~$\lambda R^{y+\cu_\ell}_\lambda f\to f$ uniformly there. Given~$\ep>0$, first choose~$\lambda$ so that~$\lambda R^{y+\cu_\ell}_\lambda f\leq\ep$ on~$F_0$ and~$\lambda R^{y+\cu_\ell}_\lambda f\geq3/4$ on~$\partial(z+r\cu_0)$. Proposition~\ref{mfd:prop-as-forms} and Lemma~\ref{mfd:lem-killing} then give, for all sufficiently large~$N$,~$u_N\coloneqq\lambda R^{y+\cu_\ell}_{N,\lambda}f\leq2\ep$ on~$F_0$ and~$u_N\geq1/2$ on~$\partial(z+r\cu_0)$. The strong Markov property at~$\tau_{z+r\cu_0}$, applied to the nonnegative occupation integral killed at~$\tau_{y+\cu_\ell}$, gives, for every~$x\in F_0$,
\begin{equation*}
	u_N(x)\geq\int
	e^{-\lambda\tau_{z+r\cu_0}}u_N(Z_{\tau_{z+r\cu_0}})\indc_{\{\tau_{z+r\cu_0}\leq h\}}\,dK^\omega_{N,x}
	\geq\frac12e^{-\lambda h}K^\omega_{N,x}(\tau_{z+r\cu_0}\leq h)
	\,.
\end{equation*}
Consequently,~$\sup_{x\in F_0}K^\omega_{N,x}(\tau_{z+r\cu_0}\leq h)\leq8\ep$ for~$h\leq(\log2)/\lambda$. The choices are in the order~$\ep$, then~$\lambda$, then~$N$, then~$h$. Using small cubes~$z+r\cu_0$ and a finite cover gives the bound for a prescribed displacement, uniformly over a compact starting set.

\smallskip

Apply this short-exit bound at successive exits from small balls, as in Step~3 of the proof of Proposition~\ref{tight:prop-tightness}. The strong Markov property bounds the number of exits on a finite time interval; it also makes the probability of a very short gap tend to zero. Thus the moduli of continuity are tight. Together with compact containment, Arzel\`a--Ascoli gives tightness on continuous path space.

\smallskip

Now let~$N_j\to\infty$ and~$x_j\to x$. Localization as in~\eqref{eq:mfd-path-localization} and killed-resolvent convergence give convergence of the whole-space discounted occupation integrals. Tightness and uniqueness of the Laplace transform give convergence of the one-time laws at every time. The Markov identity then gives convergence of the finite-dimensional laws, and approximation by cylinder functions on compact sets of paths gives~$K^\omega_{N_j,x_j}\to K^\omega_x$. If convergence were not uniform on a compact~$F$, then a sequence of offending starting points would have a convergent subsequence, contradicting this conclusion and continuity of~$x\mapsto K^\omega_x$. This proves the result in~$d_{\mathrm P}$; tightness transfers it to every metric inducing the weak topology.
\end{proof}

\begin{proof}[Proof of Theorem~\ref{t.scaling.limit}(i)]
The relabelling~\eqref{eq:mfd-normalization} and the changes of space and time transform the weak elliptic resolvent equations of the physical process into those defining~$K^\omega_{N,x}$. Uniqueness of the weak solutions identifies the resolvents and hence the semigroups; the resulting finite-dimensional distributions identify the continuous path laws.

\smallskip

Proposition~\ref{mfd:prop-as-quenched} gives the almost-sure convergence~\eqref{eq:mfd-main-as} of the quenched laws, and Proposition~\ref{mfd:prop-quenched-convergence} shows that the limiting kernel is measurable with respect to the fields. Proposition~\ref{mfd:prop-limit-properties} gives a continuous kernel, which is the family of laws of a continuous, conservative, strong Markov process with a~$\mu$-symmetric semigroup. Since the kernel is continuous, the semigroup maps bounded continuous functions to continuous functions. Proposition~\ref{mfd:prop-chaos-growth} shows that~$\mu=\lim_N\mu_N$ is locally finite, has no atoms and has full support, and Corollary~\ref{mfd:cor-finite-exit} gives the finite mean exit times. Finally, integrating the quenched convergence against~$\P$ and using that the path functionals are bounded gives the convergence of the annealed laws.
\end{proof}

\section{The limiting process and measure}
\label{lim:sec-process-measure}

In this section we prove the remaining assertions of Theorem~\ref{t.scaling.limit}: the anomalous scaling of the limiting paths in part~(iii), the properties of the reversible measure in part~(ii), and the non-Gaussian marginals in part~(iv).

\smallskip

The paths of the scaling limit leave small cubes faster than Brownian paths. The reason is the power-law decay of the effective diffusivities between arbitrary scales (Lemma~\ref{lim:lem-strict-decay}): it bounds the mean time to leave a cube of side~$3^{-k}$ by~$C3^{-(2+\eta)k}$ for every cutoff~$L\geq k$, with~$C$ independent of~$L$ and~$k$ (Lemma~\ref{lim:lem-mean-exit}), and these bounds pass to the limit. Consequently, the pointwise H\"older exponent of the limiting paths is at most~$(2+\eta)^{-1}$, and the annealed path laws are mutually singular with every mixture of Brownian laws with arbitrary covariance matrices; for each deterministic starting point, so is the quenched path law in almost every environment. The reversible measure is singular with respect to Lebesgue measure, by the change of measure of multiplicative chaos, and it gives no mass to affine hyperplanes. Since the transition laws are absolutely continuous with respect to it, none of them is Gaussian.

\smallskip

The section is organized as follows.
\begin{itemize}
\item In Section~\ref{lim:sec-exit} we extend the power-law decay of the effective diffusivities to the ratio between two arbitrary scales and bound the mean exit times.

\item In Section~\ref{lim:sec-path-laws} we prove the anomalous scaling and the singularity of the annealed and quenched path laws (Theorem~\ref{lim:thm-nonbrownian} and Corollary~\ref{lim:cor-paths}).

\item In Section~\ref{lim:sec-singularity} we prove that the reversible measure is invariant and singular with respect to Lebesgue measure (Theorem~\ref{lim:thm-measure}) and that no transition law is Gaussian (Theorem~\ref{lim:thm-nongaussian}).
\end{itemize}

Throughout, we assume~\ref{a.g1}--\ref{a.g4} at weak disorder and rescale cubes to~$\cu_0$ as in Section~\ref{s.tightness}, with unit of time~\eqref{tight:eq-local-time}. We use the notation~$\mathbf P_x^\omega$,~$\mathbf E_x^\omega$ of Step~1 of the proof of Proposition~\ref{tight:prop-tightness} for the diffusion under consideration ($X^{(L)}$,~$X$, or their rescalings to~$\cu_0$) started at~$x$.

\subsection{Power-law decay and mean exit times}
\label{lim:sec-exit}

In this subsection we bound the mean time that the diffusions~$X^{(L)}$ and~$X$ need to leave a cube (Lemma~\ref{lim:lem-mean-exit}). After the cube is rescaled to~$\cu_0$ with the unit of time~\eqref{tight:eq-local-time}, this time is at most a random constant whose moments are bounded uniformly in the cutoff and in the size and position of the cube. The mean exit time solves a Dirichlet problem whose source is the reversible measure. The mass and coercivity bounds of Section~\ref{s.tightness} give a Sobolev inequality for this measure, which bounds the solution by an iteration on level sets. To compare cubes of different sizes, we first extend the power-law decay of the effective diffusivities in Proposition~\ref{p.homogenized.coefficient.strict.decay} to the ratio between two arbitrary scales (Lemma~\ref{lim:lem-strict-decay}).

\begin{lemma}[Power-law decay between arbitrary scales]
\label{lim:lem-strict-decay}
There exist constants~$C<\infty$ and~$\eta>0$, depending only on~$d$ and the law of~$\g_0$, such that, for every~$\ell,m\in\N_0$ with~$\ell\leq m$,
\begin{equation}
	\frac{\ahom_m}{\ahom_\ell}\leq C3^{-\eta(m-\ell)}\,.
	\label{lim:eq-strict-decay}
\end{equation}
In dimension two, one can take~$\eta=\tau^2/\log3$.
\end{lemma}

\begin{proof}
In dimension two, use~\eqref{e.renormalized.diffusivities.two.dimensions}. For~$d\geq3$, we start the gluing argument of Proposition~\ref{p.homogenized.coefficient.strict.decay} at scale~$\ell$.

\smallskip

We first prove for~$\spx_\ell^\pi$ the energy bound that Proposition~\ref{p.coarse.grained.bound} gives for cubes. The Whitney cubes of~$\spx_\ell^\pi$ of side~$3^{\ell-i}$ occupy at most a proportion~$C3^{-i}$ of its volume. For~$i\leq\ell$, use on each cube the minimizer of the energy with coefficient~$\a_{\ell-i}$ among functions equal to~$p\cdot x$ on the boundary of the cube. The factors contributed by the fields in~$\a_\ell$ that are absent from~$\a_{\ell-i}$ are independent of this minimizer and have mean one. By~\eqref{e.annealed.ordering}, the expected energy density is bounded by~$C\ahom_\ell e^{2\tau^2i}|p|^2$. For~$i>\ell$, the affine competitor gives the same bound since~$\ahom_\ell\geq e^{-2\tau^2(\ell+1)}$. Join these competitors, which all equal~$p\cdot x$ on the boundaries of the Whitney cubes. Uniformly in~$\ell$,~$\pi$ and~$p$, we obtain
\begin{equation*}
	\E[p\cdot\a_\ell(\spx_\ell^\pi)p]
	\leq C\ahom_\ell|p|^2\sum_{i\geq0}3^{-i}e^{2\tau^2i}
	\leq C\ahom_\ell|p|^2\,.
\end{equation*}

\smallskip

Take~$R$ and~$q_R<1$ from~\eqref{e.strict.decay.qR.gap}. Keep the fields with indices at most~$\ell$ and with indices~$\ell+R,\ell+2R,\ldots$, and omit the others. The induction~\eqref{e.strict.decay.sparse.induction}, started at scale~$\ell$ from the bound for~$\spx_\ell^\pi$ above, and reinsertion of the omitted fields by~\eqref{e.strict.decay.full.sparse.comparison} give~$\ahom_m\leq C\ahom_\ell q_R^{\lfloor(m-\ell)/R\rfloor}$. This proves~\eqref{lim:eq-strict-decay} with~$\eta=-\log q_R/(R\log3)$.
\end{proof}

\begin{lemma}[Mean exit times]
\label{lim:lem-mean-exit}
There exists~$\gmcdelta_0>0$ such that, if~$0<\gmcdelta\leq\gmcdelta_0$, then there exists~$C<\infty$ such that, for every~$L,m\in\N_0$ and every~$z\in\Rd$, the cube~$z+\cu_m$ satisfies
\begin{equation}
	\E\Bigl[\sup_{x\in z+\cu_m}\mathbf E_x^\omega\tau_{z+\cu_m}^{(L)}\Bigr]
	\leq C3^{2m}/\ahom_{m\wedge L}\,.
	\label{lim:eq-mean-exit}
\end{equation}
More generally, there is~$c>0$ such that the following holds. For every~$q<\infty$, there exists~$\gmcdelta_0>0$ such that, if~$0<\gmcdelta\leq\gmcdelta_0$, then there are random constants~$M\geq1$, one for each~$L,m\in\N_0$ and~$z\in\Rd$, with~$\E[M^q]$ bounded uniformly in~$L$,~$m$ and~$z$, such that the process obtained by rescaling~$z+\cu_m$ to~$\cu_0$ satisfies
\begin{equation*}
	\sup_{x\in\cu_0}\mathbf E_x^\omega\tau_{\cu_0}\leq M
\end{equation*}
and
\begin{equation*}
	\inf_{x\in\cu_{-2}}\mathbf E_x^\omega\tau_{\cu_0}\geq cM^{-2}\,.
\end{equation*}
The same assertions hold for the diffusion~$X$ without cutoff, with~$\ahom_m$ in~\eqref{lim:eq-mean-exit}.
\end{lemma}

\begin{proof}
We rescale~$z+\cu_m$ to~$\cu_0$ as in Section~\ref{s.tightness}, for either the diffusion~$X^{(L)}$ or the diffusion~$X$; here~$\mu$ has a positive, locally~$C^{1,1}$ density. The mass bound~\eqref{tight:eq-mass}, the trace inequality~\eqref{tight:eq-trace} and the coercivity bound~\eqref{tight:eq-coercivity} for functions with zero boundary values give, on every cube~$w+r\cu_0$ that is fixed in the rescaled coordinates, and for every~$v\in H_0^1(w+r\cu_0)$,
\begin{equation}
	\|v\|_{L^{\nf{2d}{(d-1)}}(w+r\cu_0,\mu)}^2\leq M\mathscr E_A(v;w+r\cu_0)
	\,.
	\label{lim:eq-killed-sobolev}
\end{equation}
Here~$M$ may depend on~$w$ and~$r$ and has the asserted uniform moment bounds.

\smallskip

Let~$u$ solve~$-\nabla\cdot(A\nabla u)=d\mu/dx$ in~$\cu_0$ with zero boundary values. Testing with~$(u-a)_+$ gives, for~$b>a\geq0$,
\begin{equation*}
	\mu(\{u>b\})\leq\left(\frac M{b-a}\right)^{\nf{2d}{(d-1)}}\mu(\{u>a\})^{\nf{(d+1)}{(d-1)}}\,.
\end{equation*}
Iterate this inequality along the thresholds equal to~$1-2^{-i}$ times a dimensional multiple of~$M\mu(\cu_0)^{\nf1d}$ to obtain~$\|u\|_\infty\leq CM\mu(\cu_0)^{\nf1d}$, which by~\eqref{tight:eq-mass} is bounded by another random constant~$M$ with the asserted moments.

\smallskip

We identify~$u$ with the mean exit time. The semigroup of the diffusion killed on leaving~$\cu_0$ is associated with the part of the Dirichlet form on~$\cu_0$~\citep[Theorems~4.4.2--4.4.3]{FOT}. Interpolate~\eqref{lim:eq-killed-sobolev} with~$L^1(\mu)$ to obtain a Nash inequality with respect to~$\mu$. On each fixed cube, it bounds the killed semigroup from~$L^1(\mu)$ to~$L^\infty(\mu)$ by a constant times~$t^{-d}$. We adapt the eigenfunction construction for symmetric operators in~\citet[Chapter~VII, Section~6, Theorem~6.3 and Proposition~6.4]{Bass} to~$L^2(\cu_0,\mu)$. Interior elliptic regularity of the eigenfunctions, and uniform convergence of the eigenfunction expansion on compact sets of positive times and interior points, give a kernel~$\mathsf p_{\cu_0}(t,x,y)$ that is jointly continuous at positive times and interior points.

\smallskip

Since the drift of the rescaled stochastic differential equation is locally bounded, It\^o's formula applied to the squared distance from the starting point, stopped on exit from a small ball, bounds the probability of leaving that ball before time~$t$ by a quantity tending to zero as~$t\downarrow0$, uniformly over starting points in a compact subset of a fixed cube.

\smallskip

The local Nash bounds, this exit estimate and continuity of the semigroup on compactly supported continuous functions show, by approximation in~$L^1(\mu)$ on an enclosing cube, that the semigroup of the diffusion without killing maps bounded Borel functions to continuous functions at each positive time. The exit estimate also gives continuity in the starting point of the transition probabilities of the diffusion killed on leaving~$\cu_0$. For each fixed positive time and Borel subset of~$\cu_0$, the killed transition probability and the integral of~$\mathsf p_{\cu_0}$ over that subset are continuous and agree~$\mu$-almost everywhere. They therefore agree at every interior starting point, since~$\mu$ has full support.

\smallskip

Coercivity gives a positive spectral gap. Together with the bound on the killed semigroup, this gives~$\mathbf P_x^\omega(\tau_{\cu_0}>t)\leq Ce^{-c_1t}$ for~$t\geq1$, uniformly in~$x\in\cu_0$. Integrate the survival probability in time to obtain a function of the starting point that is continuous at interior points; as an element of~$L^2(\mu)$, it equals the inverse of the positive killed generator applied to~$1$. This function is the continuous representative of~$u$ at interior points. With this representative, for every~$x\in\cu_0$,
\begin{equation*}
	u(x)=\int_0^\infty\int_{\cu_0}\mathsf p_{\cu_0}(t,x,y)\,\mu(dy)\,dt
	=\int_0^\infty \mathbf P_x^\omega(\tau_{\cu_0}>t)\,dt
	=\mathbf E_x^\omega\tau_{\cu_0}\,.
\end{equation*}
The bound on~$\|u\|_\infty$ obtained above therefore holds at every interior point for this representative. This proves the upper bound after rescaling, and restoring the time unit~\eqref{tight:eq-local-time} gives~\eqref{lim:eq-mean-exit}. The lower bound follows from~\eqref{tight:eq-fast-exit}: taking~$t=(2M)^{-2}$ gives~$\mathbf P_x^\omega(\tau_{\cu_0}>t)\geq\nf12$ for every~$x\in\cu_{-2}$.
\end{proof}

\subsection{Annealed and quenched path laws}
\label{lim:sec-path-laws}

In this subsection we show that the scaling limit is not Brownian motion. Its annealed and quenched path laws are mutually singular with every mixture of Brownian laws (Theorem~\ref{lim:thm-nonbrownian}), and its paths oscillate more than Brownian paths at small times (Corollaries~\ref{lim:cor-paths} and~\ref{lim:cor-exit-moments}). By Lemmas~\ref{lim:lem-mean-exit} and~\ref{lim:lem-strict-decay}, the mean time to leave a cube of side~$3^{-k}$ is at most~$C3^{-(2+\eta)k}$, whereas Brownian scaling gives order~$3^{-2k}$. This bound passes to the limit and implies that~$3^{2k}\sigma_{3^{-k}}\to0$ almost surely, which has probability zero for every Brownian law.

\smallskip

For a deterministic starting point~$x$, the annealed law~$\E[K_x^\omega]$ of the limiting process is the average over~$\omega$ of its quenched law~$K_x^\omega$.

\smallskip

Let~$\mathsf P$ be a subsequential limit of the annealed laws of the rescaled diffusions~$X^{(L)}$ or~$X$ in Proposition~\ref{tight:prop-tightness}, started at deterministic points that converge in rescaled coordinates; in particular,~$\mathsf P$ may be~$\E[K_x^\omega]$. Here~$\mathsf P$ is a probability measure on paths. We use the exit times~$\sigma_r$ of Theorem~\ref{t.scaling.limit}(iii) and write~$Z^*_t\coloneqq\sup_{0\leq s\leq t}|Z_s-Z_0|$. Here~$\sigma_r$ is the exit time from the cube of side~$r$ centered at the starting point, and~$\tau_U$ is the exit time from a fixed set~$U$.

\begin{theorem}[Annealed and quenched singularity]
\label{lim:thm-nonbrownian}
There exists~$C<\infty$ such that every law~$\mathsf P$ as above satisfies, for every~$k\in\N_0$,
\begin{equation}
	\mathbf E_{\mathsf P}\sigma_{3^{-k}}\leq C3^{-(2+\eta)k}\,,
	\label{lim:eq-limit-exits}
\end{equation}
where~$\eta$ is the exponent of Lemma~\ref{lim:lem-strict-decay}. Every such annealed law~$\mathsf P$ is mutually singular with every mixture of Brownian laws with arbitrary covariance matrices. For each deterministic~$x$, the quenched law~$K_x^\omega$ has the same singularity for~$\P$-almost every~$\omega$. More generally, every subsequential limit in distribution of the random quenched laws of either family in Proposition~\ref{tight:prop-tightness}, started at deterministic points that converge in rescaled coordinates, is almost surely mutually singular with every such mixture of Brownian laws.
\end{theorem}

This proves~\eqref{e.intro.limit.exit} and the mutual singularity in Theorem~\ref{t.scaling.limit}(iii).

\begin{proof}
\emph{Step 1.} We prove~\eqref{lim:eq-limit-exits}. For~$L\geq k$, the mean exit-time bound~\eqref{lim:eq-mean-exit}, applied at scale~$L-k$ with time divided by~$3^{2L}/\ahom_L$, and Lemma~\ref{lim:lem-strict-decay} give, for either family,
\begin{equation*}
	\E\sigma_{3^{-k}}(Z^{(L)})
	\leq C3^{-2k}\frac{\ahom_L}{\ahom_{L-k}}
	\leq C3^{-(2+\eta)k}\,.
\end{equation*}
The exit time~$\sigma_{3^{-k}}$ is a lower semicontinuous function of the path for the topology of locally uniform convergence. The portmanteau theorem, applied to this nonnegative lower semicontinuous functional, gives~\eqref{lim:eq-limit-exits}.

\smallskip

\emph{Step 2.} We prove annealed singularity. Consider the event that~$3^{2k}\sigma_{3^{-k}}\to0$ as~$k\to\infty$. By~\eqref{lim:eq-limit-exits},
\begin{equation*}
	\mathbf E_{\mathsf P}\sum_{k\geq0}3^{2k}\sigma_{3^{-k}}
	\leq C\sum_{k\geq0}3^{-\eta k}<\infty\,,
\end{equation*}
so this event has~$\mathsf P$-probability one. For a Brownian law with a nonzero finite covariance matrix, Brownian scaling gives~$3^{2k}\sigma_{3^{-k}}\overset{\mathrm{law}}=\sigma_1>0$. Apply Fatou's lemma to the events that~$3^{2k}\sigma_{3^{-k}}\leq\ep$. For every~$\ep>0$, the probability of the event~$\{3^{2k}\sigma_{3^{-k}}\to0\}$ under this law is at most the probability that~$\sigma_1\leq\ep$. Letting~$\ep\downarrow0$ shows that the event has probability zero. When the covariance matrix is zero, the law of the constant path also gives the event probability zero.

\smallskip

\emph{Step 3.} We prove the quenched assertions. For each deterministic~$x$, Step~2 with~$\mathsf P=\E[K_x^\omega]$ shows that the average over~$\omega$ of the~$K_x^\omega$-probability of the event is one. Since this probability lies in~$[0,1]$, it equals one for~$\P$-almost every~$\omega$.

\smallskip

For the general assertion, convergence in distribution of the random quenched laws implies weak convergence of their averages. The average of the random limiting law is therefore an annealed subsequential limit and gives the event mass one by Step~2. The argument of the preceding paragraph then shows that the random limiting law gives the event mass one, almost surely.
\end{proof}

For every law~$\mathsf P$ as above, the exit-time bound~\eqref{lim:eq-limit-exits} gives lower bounds on displacement at small times and an upper bound on the pointwise H\"older exponent defined in~\eqref{e.intro.limit.Holder}.

\begin{corollary}[Small-time behavior of the limiting paths]
\label{lim:cor-paths}
Let~$\mathsf P$ be a law as above. The following statements hold, with constants that do not depend on~$\mathsf P$.
\begin{enumerate}[label=(\roman*)]
\item \underline{\emph{Oscillation}}: for every~$\ep>0$, there is~$c_\ep>0$ such that,~$\mathsf P$-almost surely, for all sufficiently small~$t>0$,
\begin{equation}
	Z^*_t\geq c_\ep
	\left(\frac{t}{(\log(1/t))^{1+\ep}}\right)^{\nf1{(2+\eta)}}\,.
	\label{lim:eq-oscillation}
\end{equation}
In particular,~$\mathsf P$-almost surely,~$Z^*_t/\sqrt t\to\infty$ as~$t\downarrow0$.

\item \underline{\emph{H\"older exponent}}:~$\mathsf P$-almost surely, the pointwise H\"older exponent at zero is at most~$\nf1{(2+\eta)}$, and at most~$(2+\tau^2/\log3)^{-1}$ in dimension two.

\item \underline{\emph{Exit times}}: for every~$\ep>0$,~$\mathsf P$-almost surely, for all sufficiently large~$k$,
\begin{equation}
	\sigma_{3^{-k}}\leq k^{1+\ep}3^{-(2+\eta)k}\,.
	\label{lim:eq-endpoint-exits}
\end{equation}

\item \underline{\emph{Moments of the running maximum}}: for every~$p>0$, there is~$c_p>0$ such that, for every~$0<t\leq1$,
\begin{equation}
	\mathbf E_{\mathsf P}(Z^*_t)^p\geq c_pt^{\nf p{(2+\eta)}}\,.
	\label{lim:eq-maximal-moments}
\end{equation}
\end{enumerate}
\end{corollary}

For each deterministic~$x$, parts~(i)--(iii) also hold~$K_x^\omega$-almost surely for~$\P$-almost every~$\omega$, by Step~3 of the proof of Theorem~\ref{lim:thm-nonbrownian}. Parts~(i) and~(ii) give~\eqref{e.intro.limit.Holder} in Theorem~\ref{t.scaling.limit}.

\begin{proof}
Markov's inequality gives~$\mathsf P(\sigma_{3^{-k}}>k^{1+\ep}3^{-(2+\eta)k})\leq Ck^{-1-\ep}$, so the Borel--Cantelli lemma proves~(iii). For~$t$ sufficiently small, choose~$k$ with
\begin{equation*}
	k^{1+\ep}3^{-(2+\eta)k}\leq t
	<(k-1)^{1+\ep}3^{-(2+\eta)(k-1)}\,.
\end{equation*}
Then~\eqref{lim:eq-endpoint-exits} gives~$Z_t^*\geq3^{-k}/2$, and~$k$ is comparable to~$\log(1/t)$, proving~(i). Part~(i) implies~(ii), since~$|Z_t-Z_0|=O(t^\gamma)$ would give~$Z_t^*=O(t^\gamma)$ and contradict~\eqref{lim:eq-oscillation} for~$\gamma>\nf1{(2+\eta)}$.

\smallskip

For~(iv), enlarge~$C$ in~\eqref{lim:eq-limit-exits} to at least one and take the smallest~$k\geq0$ with~$C3^{-(2+\eta)k}\leq t/2$. Then~$\mathsf P(\sigma_{3^{-k}}\leq t)\geq1/2$ and~$3^{-k}$ is comparable to~$t^{\nf1{(2+\eta)}}$, giving the claimed moment bound.
\end{proof}

\begin{corollary}[Exit moments and small displacements]
\label{lim:cor-exit-moments}
For every~$p>0$, there exists~$\gmcdelta_0>0$ such that, if~$0<\gmcdelta\leq\gmcdelta_0$, then there are~$0<c_p\leq C_p<\infty$ with the following properties.
\begin{enumerate}[label=(\roman*)]
\item For every~$L,m\in\N_0$ and~$z\in\Rd$,
\begin{equation}
	\E\left[\sup_{x\in z+\cu_m}\mathbf E_x^\omega
	(\tau_{z+\cu_m}^{(L)})^p\right]\leq C_p3^{2mp}\ahom_{m\wedge L}^{-p}\,.
	\label{lim:eq-exit-moments}
\end{equation}
For starting points in the inner cube, the matching lower bound is
\begin{equation*}
	\E\left[\inf_{x\in z+\cu_{m-2}}\mathbf E_x^\omega
	(\tau_{z+\cu_m}^{(L)})^p\right]\geq c_p3^{2mp}\ahom_{m\wedge L}^{-p}\,.
\end{equation*}
The same bounds hold for~$X$ with~$m$ in place of~$m\wedge L$.
\item Every~$\mathsf P$ as above satisfies~$\mathbf E_{\mathsf P}\sigma_{3^{-k}}^p\leq C_p3^{-p(2+\eta)k}$ for every~$k\in\N_0$, and, for all~$t,\ep\in(0,1]$,
\begin{equation*}
	\mathsf P\left(Z_t^*\leq\ep t^{\nf1{(2+\eta)}}\right)
	\leq C_p\ep^{p(2+\eta)}\,.
\end{equation*}
\item For every~$0<q<p(2+\eta)$ there is~$C_{p,q}<\infty$ such that every~$\mathsf P$ as above satisfies, for every~$t\in(0,1]$,
\begin{equation*}
	\mathbf E_{\mathsf P}(Z_t^*)^{-q}
	\leq C_{p,q}t^{-\nf q{(2+\eta)}}\,.
\end{equation*}
\end{enumerate}
\end{corollary}

\begin{proof}
After rescaling~$z+\cu_m$ to~$\cu_0$ as in the proof of Lemma~\ref{lim:lem-mean-exit}, Markov's inequality and the strong Markov property give, for each environment and every~$n\in\N_0$, the bound~$\sup_{x\in\cu_0}\mathbf P_x^\omega(\tau_{\cu_0}>2nM)\leq2^{-n}$. Integrating yields~$\sup_x\mathbf E_x^\omega\tau_{\cu_0}^p\leq C_pM^p$. For the lower bound,~\eqref{tight:eq-fast-exit} and~$\E M\leq C$ give, for a sufficiently small deterministic~$c>0$,
\begin{equation*}
	\E\left[\inf_{x\in\cu_{-2}}\mathbf E_x^\omega\tau_{\cu_0}^p\right]
	\geq c^p\E\left[\inf_{x\in\cu_{-2}}\mathbf P_x^\omega(\tau_{\cu_0}>c)\right]
	\geq\frac{c^p}2\,.
\end{equation*}
Taking expectations in the upper bound and restoring the time unit~\eqref{tight:eq-local-time} proves both bounds. Rescaling as in Step~1 of the proof of Theorem~\ref{lim:thm-nonbrownian}, Lemma~\ref{lim:lem-strict-decay} and the portmanteau theorem for~$\sigma_{3^{-k}}^p$ give the exit-moment bound for~$\mathsf P$.

\smallskip

For~$0<r\leq1/12$, choose~$k\geq0$ with~$4r\leq3^{-k}<12r$. If~$Z_t^*\leq r$, then~$\sigma_{3^{-k}}>t$, so Markov's inequality gives
\begin{equation*}
	\mathsf P(Z_t^*\leq r)
	\leq t^{-p}\mathbf E_{\mathsf P}\sigma_{3^{-k}}^p
	\leq C_pt^{-p}r^{p(2+\eta)}\,.
\end{equation*}
Take~$r=\ep t^{\nf1{(2+\eta)}}$ to obtain the stated probability bound when~$r\leq1/12$. If~$1/12<r\leq1$, then~$\ep\geq r>1/12$, so the probability is at most~$1\leq12^{p(2+\eta)}\ep^{p(2+\eta)}$. Enlarging~$C_p$ to at least~$12^{p(2+\eta)}$ covers this interval. Integrating the probability bound against~$q\ep^{-q-1}\,d\ep$ on~$(0,1)$ gives a finite contribution because~$q<p(2+\eta)$; the range~$\ep\geq1$ contributes at most~$\int_1^\infty q\ep^{-q-1}\,d\ep=1$. Restoring the factor~$t^{-\nf q{(2+\eta)}}$ proves the negative-moment bound.
\end{proof}

\subsection{Singularity of the measure and of the transition laws}
\label{lim:sec-singularity}

In this subsection we show that the reversible measure is singular with respect to Lebesgue measure and invariant for the limiting process (Theorem~\ref{lim:thm-measure}), and that the transition laws are absolutely continuous with respect to it and singular with respect to every Gaussian law (Theorem~\ref{lim:thm-nongaussian}). Singularity comes from the multiplicative structure of the densities of~$\mu_N^0$: their logarithm is a sum of independent terms, with negative mean at Lebesgue-typical points and positive mean at points sampled from the limiting measure. Absolute continuity comes from a Nash inequality for the killed semigroups at finite cutoff, which passes to the limit.

\smallskip

The limiting measures~$\mu^0$ and~$\mu=e^H\mu^0$ of Theorem~\ref{t.scaling.limit}(ii) are constructed in Proposition~\ref{mfd:prop-chaos-growth} and its proof. For Gaussian fields,~$\tau^2=\frac12\E[\g_0(0)^2]$, so~$\mu^0$ is Gaussian multiplicative chaos, and~$\mu$ is its product with the positive, locally~$C^{1,1}$ function~$e^H$. Existence of the limit, moment bounds, absence of atoms, full support, and singularity with respect to Lebesgue measure are standard properties of Gaussian multiplicative chaos in the subcritical regime, in which the variance of the field is small enough for the limiting measure to be nonzero. Moreover, the densities of the approximating measures diverge at almost every point with respect to the limiting measure; see~\citet{Kahane} and~\citet[Sections~2 and~4.1]{RVsurvey}. The proofs in Section~\ref{mfd:sec-speed} and below also apply to the non-Gaussian fields allowed by~\ref{a.g1}--\ref{a.g4}.

\smallskip

The measures~$\a_N(3^Nx)\,dx$ have the laws of~$\mu_N^0$ and therefore converge in distribution to~$\mu^0$. The measures~$\mu^0$ and~$\mu$ are nonatomic and have full support, and Proposition~\ref{mfd:prop-chaos-growth} with~$\ep=\nf12$ gives, on each bounded set of centers, a random constant~$M$ such that, for every~$x$ in that set and every~$0<r\leq1$,
\begin{equation}
	\mu^0(B_r(x))+\mu(B_r(x))\leq Mr^{d-\nf12}\,.
	\label{lim:eq-measure-growth}
\end{equation}
For every finite moment order,~$M$ has a finite moment of that order at weak disorder, with a threshold depending on the order, by the same proposition and by Lemma~\ref{mfd:lem-infrared} for~$e^{\pm H}$.

\begin{theorem}[The limiting measures]
\label{lim:thm-measure}
The following statements hold.
\begin{enumerate}[label=(\roman*)]
\item \underline{\emph{Translation invariance}}: the law of~$\mu^0$ is invariant under translations,~$\mu^0$ is not deterministic, and it satisfies~$\E[\mu^0(dx)]=dx$. Its restrictions to two sets at distance greater than~$\sqrt d$ are independent.

\item \underline{\emph{Singularity}}: almost surely, the density of~$\mu_N^0$ at~$x$ tends to zero as~$N\to\infty$ for Lebesgue-almost every~$x$, whereas it tends to infinity for~$\mu^0$-almost every and~$\mu$-almost every~$x$. In particular, both limiting measures are singular with respect to Lebesgue measure.

\item \underline{\emph{Invariance}}: almost surely, the reversible measure~$\mu$ is invariant for the limiting process: for every nonnegative Borel function~$f$ and every~$t\geq0$,
\begin{equation}
	\int\!\!\int f(Z_t)\,dK_x^\omega\,\mu(dx)=\int f\,d\mu\,.
	\label{lim:eq-invariant-measure}
\end{equation}
\end{enumerate}
\end{theorem}

\begin{proof}
Write~$\mathsf d_N$ for the density of~$\mu_N^0$. Lemma~\ref{mfd:lem-chaos-moments} gives uniform integrability of the masses~$\mu_N^0(z+r\cu_0)$ on bounded cubes~$z+r\cu_0$, and hence
\begin{equation*}
	\E[\mu^0(dx)\mid\g_0,\ldots,\g_{-N}]=\mathsf d_N(x)\,dx\quad\text{and}\quad
	\E\mu^0(dx)=dx\,.
\end{equation*}
The strong law and Fubini's theorem give, for~$\P\otimes dx$-almost every~$(\omega,x)$,
\begin{equation}
	\frac{\log\mathsf d_N(x)}{N+1}
	=\frac1{N+1}\sum_{i=0}^N\g_{-i}(x)-\tau^2
	\longrightarrow-\tau^2<0\,.
	\label{lim:eq-lebesgue-log-density}
\end{equation}
For the behavior at a point sampled from~$\mu^0$, fix a bounded cube~$z+r\cu_0$ and consider the probability measure~$r^{-d}\P(d\omega)\mu^{0,\omega}(dx)$ on~$\Omega\times(z+r\cu_0)$. This is the change of measure used for Gaussian chaos in~\citet[Section~4.1]{RVsurvey}; here independence of the fields suffices. The conditional-expectation identity above and stationarity show that, under this measure, the variables~$\g_{-j}(x)$ for~$j\geq0$ are independent, and each has the law of~$\g_0(0)$ tilted by~$e^{\g_0(0)-\tau^2}$.

\smallskip

Their common mean is~$\E[\g_0(0)e^{\g_0(0)-\tau^2}]>\tau^2$: under the tilted law,~$e^{-\g_0(0)}$ has expectation~$e^{-\tau^2}$, and Jensen's inequality for the exponential is strict because~$\g_0(0)$ is not almost surely constant. The strong law therefore gives
\begin{equation*}
	\frac{\log\mathsf d_N(x)}{N+1}
	\longrightarrow\E[\g_0(0)e^{\g_0(0)-\tau^2}]-\tau^2>0
\end{equation*}
almost surely under this probability measure. Thus~$\mathsf d_N\to\infty$ at~$\mu^0$-almost every point, almost surely. The positive factor~$e^H$ gives the same conclusion for~$\mu$. Together with~\eqref{lim:eq-lebesgue-log-density}, this proves~(ii). Since~$\E\mu^0(dx)=dx$, singularity also shows that~$\mu^0$ is not deterministic. Translation invariance of the law and independence of restrictions to sets at distance greater than~$\sqrt d$ pass from the fields to~$\mu^0$, which is a measurable function of them.

\smallskip

Finally, symmetry and conservativeness from Proposition~\ref{mfd:prop-limit-properties} give~\eqref{lim:eq-invariant-measure}.
\end{proof}

The Sobolev estimate~\eqref{lim:eq-killed-sobolev} gives absolute continuity of the transition laws at positive times with respect to~$\mu$. The measure~$\mu$ is carried by a Lebesgue-null set and gives no mass to affine hyperplanes, so these transition laws are singular with respect to every Gaussian law.

\begin{theorem}[Quenched transition laws]
\label{lim:thm-nongaussian}
\label{lim:lem-transition-domination}
Almost surely, the following statements hold.
\begin{enumerate}[label=(\roman*)]
\item For every~$x\in\Rd$ and every~$t>0$,
\begin{equation}
	K_x^\omega(Z_t\in\cdot)\ll\mu
	\label{lim:eq-transition-domination}
\end{equation}
and the law of~$Z_t$ under~$K_x^\omega$ is singular with respect to Lebesgue measure and every Gaussian law.
\item For every bounded open set~$U$, there is~$C_U<\infty$ such that, for every~$x\in U$ and~$t>0$,
\begin{equation*}
	K_x^\omega(Z_t\in\cdot,\ t<\tau_U)
	\leq C_Ut^{-d}\mu|_U\,.
\end{equation*}
\item For every starting point, the quenched joint law of the positions at every finite set of positive times is mutually singular with every Gaussian law, and the quenched path law on every interval~$[s,t]$ with~$0\leq s<t$ is mutually singular with the law of every continuous Gaussian process on that interval.
\end{enumerate}
\end{theorem}

\begin{proof}
Fix~$j\in\N_0$ and the cube~$\cu_j$ of the deterministic countable exhaustion~$(\cu_j)_{j\in\N_0}$ of~$\Rd$, and let~$P_{N,t}$ be the symmetric killed semigroup of~$\mathscr E_N^{\cu_j}$ on~$L^2(\cu_j,\mu_N)$. By the relabelling~\eqref{eq:mfd-normalization}, the coefficient of~$\mathscr E_N^{\cu_j}$ and the measure~$\mu_N$ have the joint law of the coefficient and the reversible measure after rescaling of~$X$ at scale~$N$. Hence~\eqref{lim:eq-killed-sobolev} applies on~$\cu_j$ with constants~$M_N$ satisfying~$\sup_N\E M_N<\infty$.

\smallskip

We use Nash's argument. For nonzero~$f\in L^1(\cu_j,\mu_N)\cap L^2(\cu_j,\mu_N)$, interpolation of~\eqref{lim:eq-killed-sobolev} with~$L^1$, the energy identity and~$L^1$ contraction give
\begin{equation*}
	\frac d{dt}\|P_{N,t}f\|_{L^2(\cu_j,\mu_N)}^2
	\leq-2M_N^{-1}\|P_{N,t}f\|_{L^2(\cu_j,\mu_N)}^{2+\nf2d}
	\|f\|_{L^1(\cu_j,\mu_N)}^{-\nf2d}\,.
\end{equation*}
Integration, symmetry and the semigroup property yield, for every~$t>0$,
\begin{equation}
	\|P_{N,t}f\|_\infty
	\leq C(M_N/t)^d\|f\|_{L^1(\cu_j,\mu_N)}\,.
	\label{lim:eq-killed-domination}
\end{equation}
For~$f\in C_c(\cu_j)$, the function~$P_{N,t}f$ is continuous by the argument in the proof of Lemma~\ref{lim:lem-mean-exit}, so this bound holds at every point of~$\cu_j$.

\smallskip

By Proposition~\ref{mfd:prop-as-quenched}, the whole-space path laws converge almost surely, locally uniformly in the starting point. Fatou's lemma gives~$\liminf_N M_N<\infty$ almost surely. For each such environment, choose a subsequence on which~$M_N$ is bounded. If~$f\in C_c(\cu_j)$ is nonnegative, then the path functional~$w\mapsto f(w(t))\indc_{\{t<\tau_{\cu_j}\}}$ is lower semicontinuous: a path that stays in~$\cu_j$ up to time~$t$ has a positive distance from its complement. The portmanteau theorem,~\eqref{lim:eq-killed-domination}, and the almost-sure local weak convergence~$\mu_N\to\mu$ of Proposition~\ref{mfd:prop-chaos-growth} therefore give the same bound for the killed transition law of~$K^\omega_x$. Positive compactly supported tests give the asserted measure domination. The subsequence depends only on the environment and~$\cu_j$, so the bound holds for every~$x\in\cu_j$ and~$t>0$.

\smallskip

Enclosing an arbitrary bounded open set in one of these cubes and comparing exit times gives the bound on that set. By Proposition~\ref{mfd:prop-limit-properties}, the limiting process is conservative and has continuous paths, so each path remains in some cube of the exhaustion up to time~$t$. Monotone convergence therefore proves~\eqref{lim:eq-transition-domination} for all~$x$ and~$t>0$ on a single event of probability one.

\smallskip

By Theorem~\ref{lim:thm-measure}, a Lebesgue-null set carries~$\mu$ and hence every transition law, proving singularity with respect to Lebesgue measure and every nondegenerate Gaussian. A bounded subset of an affine hyperplane is covered by~$O(r^{1-d})$ balls of radius~$r$, so~\eqref{lim:eq-measure-growth} bounds its~$\mu$-mass by a constant times~$r^{\nf12}$, which tends to zero as~$r\downarrow0$. Thus every affine hyperplane is~$\mu$-null. Every degenerate Gaussian, including a point mass, is supported on such a hyperplane, which has zero probability under every transition law by~\eqref{lim:eq-transition-domination}. This proves mutual singularity with every Gaussian law.

\smallskip

Since the marginals of a Gaussian law are Gaussian, and two laws are mutually singular as soon as their marginals at some time are, projecting onto a single positive time gives the assertions for joint laws and path laws.
\end{proof}

\clearpage \appendix

\section{Moments of random variables with log-normal tails}
\label{s.orlicz}

In this appendix we collect the properties of the Orlicz notation used for the stochastic bounds of Sections~\ref{s.infrared.cutoffs},~\ref{s.self.similar} and~\ref{s.Holder.regularity}: the tail bound, the bounds for sums and maxima, and the moments of a random variable with Gaussian tails and of its exponential. The last of these turns Gaussian tails for the oscillation of the fields into moment bounds for their exponentials.

\smallskip

Stochastic bounds throughout the paper are stated in the Orlicz notation of~\citet[Appendix~A]{AKMBook}. All references to~\citet[Appendix~A]{AKMBook} in this appendix refer to the preprint version \href{https://arxiv.org/abs/1705.05300v2}{arXiv:1705.05300v2}, including (A.1)--(A.2), Remark~A.2 and Lemmas~A.1, A.4, A.7 and~A.9. The subscript refers to the Orlicz function~$\exp(t^s)$, and for~$s>0$ and~$A>0$ we write
\begin{equation*}
	X\leq\O_{\Gamma_s}(A)
	\quad\text{if}\quad
	\E\!\left[\exp\!\left(\left(A^{-1}\max\{X,0\}\right)^s\right)\right]\leq2\,,
\end{equation*}
and
\begin{equation}
	X=\O_{\Gamma_s}(A)
	\quad\text{if}\quad
	X\leq\O_{\Gamma_s}(A)
	\quad\text{and}\quad
	-X\leq\O_{\Gamma_s}(A)\,.
	\label{e.orlicz.two.sided.definition}
\end{equation}
If~$X\leq\O_{\Gamma_s}(A)$, then, by~\citet[Lemma~A.1]{AKMBook}, for every~$t\geq0$,
\begin{equation}
	\P[X\geq tA]\leq2\exp(-t^s)\,.
	\label{e.orlicz.tail.bound}
\end{equation}
Conversely, the bound~$\P[X\geq tA]\leq\exp(-t^s)$ for every~$t\geq0$ implies~$X\leq\O_{\Gamma_s}(2^{\nf1s}A)$.

\smallskip

By~\citet[Lemma~A.4]{AKMBook}, for every~$s\geq1$, every finite index set~$I$, and random variables~$X_k$ satisfying~$X_k\leq\O_{\Gamma_s}(a_k)$ for every~$k\in I$,
\begin{equation}
	\sum_{k\in I}X_k
	\leq
	\O_{\Gamma_s}\!\left(\sum_{k\in I}a_k\right)\,.
	\label{e.orlicz.triangle.inequality}
\end{equation}
The same conclusion holds for a countable family of nonnegative~$X_k$ with~$\sum_k a_k<\infty$, by monotone convergence. The tail bound~\eqref{e.orlicz.tail.bound}, a union bound, and the converse part of~\citet[Lemma~A.1]{AKMBook}, in the form of~\citet[Remark~A.2]{AKMBook} that requires the tail bound only beyond a threshold, give, for every~$s>0$,~$\sigma>0$ and~$k\in\N$, and all random variables~$X_1,\ldots,X_k$ satisfying~$X_i\leq\O_{\Gamma_s}(\sigma)$ for every~$i$,
\begin{equation}
	\max_{i\in\{1,\ldots,k\}}X_i
	\leq
	\O_{\Gamma_s}\bigl(C_s(\log(2k))^{\nf1s}\sigma\bigr)\,.
	\label{e.max.scaling.O.Gamma.Two}
\end{equation}
If, in addition,~$s\in(1,2]$, the variables are independent and mean zero, and~$X_i=\O_{\Gamma_s}(\sigma)$, then, by~\citet[Lemmas~A.7 and~A.9]{AKMBook},
\begin{equation}
	\sum_{i=1}^k X_i
	=\O_{\Gamma_s}\bigl(C_s k^{\nf12}\sigma\bigr)\,.
	\label{e.CLT.scaling.O.Gamma.Two}
\end{equation}

\begin{lemma}
There exists a universal constant~$C<\infty$ such that, for every~$\sigma>0$ and every random variable~$X$ satisfying~$X=\O_{\Gamma_2}(\sigma)$, the following estimates hold. For every~$k>0$,
\begin{equation}
	\E\!\left[|X|^k\right]
	\leq
	C\sigma^k\Gamma\!\left(\frac k2+1\right)\,,
	\label{e.moments.OGamma2}
\end{equation}
for every~$t>0$,
\begin{equation}
	\E\!\left[\left|\exp(tX)-(1+tX)\right|\right]
	\leq
	C\sigma^2t^2\exp(C\sigma^2t^2)\,,
	\label{e.moments.lognormal}
\end{equation}
and, for every~$m\geq1$ and~$t>0$,
\begin{equation}
	\E\!\left[|\exp(tX)-1|^m\right]
	\leq
	(Cm^{\nf12}\sigma t)^m\exp\!\left(C(m\sigma t)^2\right)\,.
	\label{e.moments.general.lognormal}
\end{equation}
\end{lemma}
\begin{proof}
By~\eqref{e.orlicz.tail.bound}, applied to~$X$ and~$-X$, we have~$\P[|X|\geq r\sigma]\leq4e^{-r^2}$ for every~$r\geq0$. The layer-cake formula gives
\begin{equation*}
	\E[|X|^k]
	\leq
	4k\sigma^k\int_0^\infty r^{k-1}e^{-r^2}\,dr
	=
	4\sigma^k\Gamma\!\left(\frac k2+1\right)\,,
\end{equation*}
which is~\eqref{e.moments.OGamma2}. Moreover, Young's inequality and~\eqref{e.orlicz.two.sided.definition} give~$\E[e^{\lambda|X|}]\leq4e^{\lambda^2\sigma^2/2}$ for every~$\lambda\geq0$. The elementary inequalities~$|e^u-1-u|\leq\tfrac12u^2e^{|u|}$ and~$|e^u-1|^m\leq|u|^me^{m|u|}$, together with the Cauchy--Schwarz inequality,~\eqref{e.moments.OGamma2} and the bound on~$\E[e^{\lambda|X|}]$, yield~$\E[|e^{tX}-1-tX|]\leq C\sigma^2t^2e^{C\sigma^2t^2}$ and
\begin{equation*}
	\E[|e^{tX}-1|^m]
	\leq
	t^m\E[|X|^{2m}]^{\nf12}\E[e^{2mt|X|}]^{\nf12}
	\leq
	(Cm^{\nf12}\sigma t)^m e^{C(m\sigma t)^2}\,.
\end{equation*}
These are~\eqref{e.moments.lognormal} and~\eqref{e.moments.general.lognormal}.
\end{proof}

\section{Concentration inequalities}
\label{app.concentration}

In this appendix we collect concentration inequalities used in the body of the paper: Rosenthal's inequality for sums of independent random variables, used in Sections~\ref{s.infrared.cutoffs} and~\ref{s.self.similar}, and bounds on the proportion of scales at which independent or finitely dependent random variables are large, used in Sections~\ref{s.Holder.regularity} and~\ref{mfd:sec-branches}.

\smallskip

The integer-exponent case of the following estimate follows by applying~\citet[Theorem~15.11]{BLM} to both signs of the sum. We give a proof which also covers every real exponent~$p\geq2$.

\begin{lemma}[Rosenthal's inequality]
\label{l.Rosenthal}
For every~$p\in [2,\infty)$, every~$N\in\N$ and every sequence~$X_1,\ldots,X_N$ of independent random variables satisfying~$\E [ X_i ] = 0$, we have the estimate
\begin{equation*}
	\E
	\Biggl[
	\biggl|
	\sum_{i=1}^N
	X_i
	\biggr|^p
	\Biggr]^{\nf 1p}
	\leq
	6p^{\nf12}  \biggl(
	\sum_{i=1}^N
	\E [ X_i^2 ]
	\biggr)^{\!\nf12}
	+
	4p
	\E \Bigl[ \max_{i\in\{1,\ldots,N\} } |X_i| ^p \Bigr]^{\nf1p}
	\,.
\end{equation*}
\end{lemma}
\begin{proof}
Write~$\|Y\|_r=\E[|Y|^r]^{1/r}$. If the right side is infinite, there is nothing to prove; otherwise all variables below have the required moments. For every real~$r\geq1$, let~$K_r\coloneqq2^{1/r}e^{-1/2}\sqrt r$. Independent uniformly distributed signs~$\varepsilon_i\in\{-1,1\}$ satisfy~$\|\sum_i\varepsilon_i a_i\|_r\leq K_r(\sum_i a_i^2)^{1/2}$. Indeed, their exponential moment is at most~$\exp(t^2\sum_i a_i^2/2)$. Use~$x^r\leq(r/t)^r e^{-r}e^{tx}$, bound~$e^{t|x|}$ by~$e^{tx}+e^{-tx}$, and take~$t=\sqrt{r/\sum_i a_i^2}$; the zero case is immediate. Symmetrizing an independent family~$Y_i$ with an independent copy, applying this sign estimate to the differences, and using the triangle inequality for the Euclidean square function gives~$\|\sum_i(Y_i-\E Y_i)\|_r\leq2K_r\|(\sum_iY_i^2)^{1/2}\|_r$. Let~$A\coloneqq\sum_iX_i^2$,~$u\coloneqq\|A^{1/2}\|_p$,~$b\coloneqq\|\max_i|X_i|\|_p$ and~$\sigma^2\coloneqq\sum_i\E X_i^2$. Apply the preceding estimate to~$Y_i=X_i^2$ with~$r=p/2\geq1$. Since~$\sum_iX_i^4\leq(\max_i|X_i|)^2A$, H\"older's inequality gives
\begin{equation*}
	u^2=\|A\|_{p/2}
	\leq\sigma^2+2K_{p/2}\|(\max_i|X_i|)A^{1/2}\|_{p/2}
	\leq\sigma^2+2K_{p/2}bu\,.
\end{equation*}
Consequently~$u\leq\sigma+2K_{p/2}b$. Applying the centered estimate once more, now to~$X_i$ at exponent~$p$, yields
\begin{equation*}
	\Big\|\sum_iX_i\Big\|_p
	\leq2K_p\sigma+4K_pK_{p/2}b
	\leq6\sqrt p\,\sigma+4p b\,.
\end{equation*}
The last inequalities follow from~$2^{1/p}\leq\sqrt2$ and~$2^{2/p}\leq2$:~$2K_p\leq2\sqrt{2/e}\sqrt p\leq6\sqrt p$ and~$4K_pK_{p/2}\leq8p/e\leq4p$.
\end{proof}
Let~$p\in[1,\infty)$,~$r \in \N$ and~$s \in [\nf1p,1]$, and let~$\{X_{k,j}:k,j\in\Z\}$ be nonnegative random variables with~$\E[X_{k,j}^p]\leq1$ for every~$k,j\in\Z$. Assume that the columns are~$r$-dependent: for every two subsets~$I_1,I_2\subseteq\Z$ with~$\dist(I_1,I_2)\geq r$,
\begin{equation}
	\{ X_{k,j}\,:\, k \in \Z\,, \, j \in I_1 \}
	\quad \text{and} \quad
	\{ X_{k,j}\,:\, k \in \Z\,, \, j \in I_2 \}
	\quad \text{are independent}\,.
	\label{e.independent.columns.twosided}
\end{equation}

\begin{proposition}[Concentration across finitely dependent columns]
\label{p.concentration.for.scales}
There exists a universal constant~$C > 0$ such that, for every~$m_0 \in \Z$,~$m \in\N_0$ and~$\theta \in (0,1]$,
\begin{equation*}
	\P \biggl[ \frac{1}{m+1} \sum_{k=m_0}^{m_0+m} \indc{\biggl\{ \sum_{j \in \Z} 3^{-s |k-j|} X_{k,j} > 6s^{-1} C^{\nf 1p} \theta^{-\nf 1p} \biggr\}} > \theta \biggr]
	\leq
	\exp \biggl(- \frac{  s p \theta}{16 r} (m+1)  \biggr)
	\,.
\end{equation*}
\end{proposition}
\begin{proof}
This is~\citet[Proposition~D.1]{ASD}.
\end{proof}

\begin{proposition}[Concentration across scales]
\label{p.concentration.for.scales.exp.sequence}
There is a universal constant~$C<\infty$ with the following property. Let~$\mathsf a_{\mathrm{sg}}\in(0,1]$, and let~$(X_i)_{i\in\Z}$ be independent nonnegative random variables satisfying~$X_i=\O_{\Gamma_2}(\mathsf a_{\mathrm{sg}})$ for every~$i\in\Z$. Then, for every~$m_0\in\Z$,~$M\in\N_0$,~$\theta\in(0,1]$, and~$s\in(0,1]$ such that~$\mathsf a_{\mathrm{sg}}\leq C^{-1}(s\theta)^{\nf12}$,
\begin{equation}
	\P \biggl[ \frac{1}{M+1}\sum_{k=m_0}^{m_0+M}
	\indc{\biggl\{\sup_{j \geq 0}3^{- s j}
	\exp\biggl(\sum_{i=0}^\infty3^{-i}X_{k+i+j}\biggr)>3\biggr\}}
	\geq \theta \biggr]
	\leq
	\exp \biggl( - \frac{s\theta}{C\mathsf a_{\mathrm{sg}}^{2}}(M+1) \biggr)\,.
	\label{e.concentration.scales.exp.sequence}
\end{equation}
\end{proposition}

\begin{proposition}[Concentration across scales for random fields]
\label{p.concentration.for.scales.exp.field}
Let~$\mathsf a_{\mathrm{sg}}\in(0,1]$, and let~$(X_i)_{i\in\Z}$ be independent real-valued random fields. Assume that~$X_i$ is~$3^i\Zd$-stationary and~$\|X_i\|_{L^\infty(\cu_i)}\leq\O_{\Gamma_2}(\mathsf a_{\mathrm{sg}})$ for every~$i\in\Z$. There exists~$C(d)<\infty$ such that, for every~$m_0\in\Z$,~$M\in\N_0$,~$h\in\N$, and~$s,\theta\in(0,1]$ satisfying~$\mathsf a_{\mathrm{sg}}\leq C^{-1}s(\theta\wedge h^{-1})^{\nf12}$,
\begin{equation}
	\P \Biggl[ \frac{1}{M+1}\sum_{k=m_0}^{m_0+M} \indc{\biggl\{ \sup_{j \geq 0}
	3^{- s j }
	\biggl\| \prod_{i=k - j}^{k+j} \exp\bigl(| X_i |\bigr) \biggr\|_{L^\infty(\cu_{k+h+j})}  > 3 \biggr\}} \geq \theta \Biggr]
	\leq
	\exp \biggl( - \frac{s^2 \theta}{C \mathsf a_{\mathrm{sg}}^2} (M+1) \biggr)
	\,.
	\label{e.concentration.scales.exp.field}
\end{equation}
\end{proposition}

Propositions~\ref{p.concentration.for.scales.exp.sequence} and~\ref{p.concentration.for.scales.exp.field} are deduced from the next lemma.

\begin{lemma}[Selection of rare intervals]
\label{l.concentration.rare.intervals}
There are universal constants~$c>0$ and~$C<\infty$ such that the following holds. Let~$(\xi_i)_{i\in\Z}$ be independent random elements and let~$\lambda>0$. For~$k\in\Z$ and~$j\in\N_0$, let~$E_{k,j}$ be an event that depends only on the variables~$\xi_i$ with indices~$i\in[k-j,k+j]\cap\Z$, and assume that~$\P[E_{k,j}]\leq e^{-\lambda(j+1)}$. If~$0<\theta\leq1$ and~$\lambda\theta\geq C$, then every integer interval~$I$ of cardinality~$m\geq1$ satisfies
\begin{equation*}
	\P\bigl[\#\{k\in I:E_{k,j}\text{ occurs for some }j\}\geq\theta m\bigr]
	\leq e^{-c\lambda\theta m}\,.
\end{equation*}
\end{lemma}

\begin{proof}
A union bound over~$k\in I$ and~$j\geq m$ shows that the probability that~$E_{k,j}$ occurs for some~$k\in I$ and~$j\geq m$ is at most~$e^{-c\lambda m}$. On the complementary event, for each~$k\in I$ such that~$E_{k,j}$ occurs for some~$j$, choose one such~$j$, necessarily~$j<m$, and take the corresponding interval of indices~$[k-j,k+j]$. Retain these intervals greedily in decreasing order of length, discarding those at distance less than two from a retained interval. Every discarded interval has length at most that of the interval which discards it, and lies in its enlargement by that length on each side. Hence the retained intervals have total cardinality at least~$\nf{\theta m}3$.

\smallskip

The retained intervals are separated by at least one missing integer and lie in an interval of cardinality at most~$3m$. Their union determines the intervals and their starting indices, so there are at most~$2^{3m}$ configurations. For each configuration, independence and the fact that~$j+1$ is at least half the cardinality of the corresponding interval give probability at most~$e^{-\nf{\lambda\theta m}6}$. Thus the required probability is at most~$e^{-c\lambda m}+e^{Cm-\nf{\lambda\theta m}6}$. Taking the constant in the assumption~$\lambda\theta\geq C$ large enough absorbs the counting factor~$e^{Cm}$ and the sum of the two terms.
\end{proof}

\begin{proof}[Proof of Proposition~\ref{p.concentration.for.scales.exp.sequence}]
Call an index~$k$ bad if the event in the indicator in~\eqref{e.concentration.scales.exp.sequence} occurs. Every bad index~$k$ admits some~$j\geq0$ such that
\begin{equation}
	X_{k+j}>\frac{\log3}3(1+sj)\,.
	\label{e.concentration.sequence.endpoint.witness}
\end{equation}
Indeed, if~$X_{k+j}\leq\frac{\log3}3(1+sj)$ for every~$j\geq0$, then, for every~$j\geq0$,
\begin{equation*}
	\sum_{i\geq0}3^{-i}X_{k+j+i}
	\leq\frac{\log3}3\left(\frac32(1+sj)+\frac34s\right)
	\leq\frac34(1+sj)\log3\,,
\end{equation*}
which excludes the event in~\eqref{e.concentration.scales.exp.sequence}. The event~\eqref{e.concentration.sequence.endpoint.witness} depends only on the index~$k+j$ and has probability at most~$2e^{-c(1+sj)^2/\mathsf a_{\mathrm{sg}}^2}\leq e^{-c_1s(j+1)/\mathsf a_{\mathrm{sg}}^2}$, since~$(1+sj)^2\geq s(j+1)$. We apply Lemma~\ref{l.concentration.rare.intervals} to these events on the~$M+1$ indices~$m_0,\ldots,m_0+M$, with~$\lambda=\nf{c_1s}{\mathsf a_{\mathrm{sg}}^2}$, which gives
\begin{equation*}
	\P[\#\{\text{bad indices}\}\geq\theta(M+1)]
	\leq\exp\left(-\frac{s\theta}{C\mathsf a_{\mathrm{sg}}^2}(M+1)\right)\,.
\end{equation*}
This is~\eqref{e.concentration.scales.exp.sequence}.
\end{proof}

\begin{proof}[Proof of Proposition~\ref{p.concentration.for.scales.exp.field}]
Let~$E_{k,j}$ be the event that the expression at index~$j$ inside the supremum in~\eqref{e.concentration.scales.exp.field} exceeds three. It depends only on~$(X_i)_{k-j\leq i\leq k+j}$. Cover~$\cu_{k+h+j}$ by its~$3^{d(h+2j)}$ triadic cubes of scale~$k-j$. On each of these cubes, the logarithm of the product is bounded above by the sum of the suprema of~$|X_i|$ on its scale-$i$ parent cubes. These~$2j+1$ random variables are independent, have means at most~$C\mathsf a_{\mathrm{sg}}$, and, after centering, have sub-Gaussian parameters at most~$C\mathsf a_{\mathrm{sg}}$. By~\eqref{e.CLT.scaling.O.Gamma.Two}, a union bound, and~$1+sj\geq s(j+1)$,
\begin{equation*}
	\P[E_{k,j}]
	\leq\exp\left(C(h+j+1)
	-\frac{c(1+sj)^2}{\mathsf a_{\mathrm{sg}}^2(j+1)}\right)
	\leq\exp\left(-\frac{c_1s^2}{\mathsf a_{\mathrm{sg}}^2}(j+1)\right)\,.
\end{equation*}
Here~$\mathsf a_{\mathrm{sg}}\leq C^{-1}s$ absorbs the means and the term proportional to~$j+1$, while~$\mathsf a_{\mathrm{sg}}^2h\leq C^{-2}s^2$ absorbs the factor from the union bound over the covering cubes at~$j=0$, and hence at every~$j$. Lemma~\ref{l.concentration.rare.intervals}, with~$\lambda=\nf{c_1s^2}{\mathsf a_{\mathrm{sg}}^2}$ and~$I=\{m_0,\ldots,m_0+M\}$, now gives~\eqref{e.concentration.scales.exp.field}.
\end{proof}

\section{Besov norm comparisons}
\label{s.Besov.norm.comparisons}

In this appendix we prove two statements about Besov norms recorded in Section~\ref{s.notation}: the duality estimate~\eqref{e.duality.splitting}, and the comparison~\eqref{e.Wsp.vs.Bspp.intro} of Besov and fractional Sobolev seminorms. We also record the fractional Sobolev--Poincar\'e inequality used in Section~\ref{s.coarse.graining}.

\smallskip

The fractional Sobolev--Poincar\'e inequality states that, for each~$s\in(0,1]$ and~$p\in[1,\nf ds)$, there exists~$C(s,p,d)<\infty$ such that, for every~$m\in\Z$ and~$f\in W^{s,p}(\cu_m)$,
\begin{equation*}
	\| f - (f)_{\cu_m} \|_{\underline{L}^{\nf{pd}{(d-sp)}} (\cu_m)}
	\leq
	C 3^{sm}  [f ]_{\underline{W}^{s,p}(\cu_m)}
	\,.
\end{equation*}
For~$0<s<1$, Jensen's inequality gives
\begin{equation*}
	\|f-(f)_{\cu_m}\|_{\underline L^p(\cu_m)}^p
	\leq\frac1{|\cu_m|^2}\int_{\cu_m}\int_{\cu_m}|f(x)-f(y)|^p\,dx\,dy
	\leq\frac{d^{(d+sp)/2}}s 3^{msp}
	[f]_{\underline W^{s,p}(\cu_m)}^p\,.
\end{equation*}
For the second inequality, use~$|x-y|_\infty\leq3^m$, hence~$|x-y|_2\leq\sqrt d\,3^m$, together with~$|\cu_m|=3^{md}$. This keeps the single volume normalization and the factor~$s$ in the fractional seminorm; taking the~$p$th root introduces the factor~$s^{-1/p}$ and a dimensional factor. On the unit cube, the extension theorem and Sobolev embedding of~\citet[Theorems~5.4 and~6.7]{DiNezzaPalatucciValdinoci} apply also at~$p=1$, provided~$sp<d$. Apply them to~$f-(f)_{\cu_m}$ and use the preceding inequality to control the~$L^p$ term. Rescaling~$f(3^m\,\cdot)$ multiplies its normalized fractional seminorm by~$3^{ms}$ and leaves its normalized~$L^{pd/(d-sp)}$ norm unchanged, proving the displayed estimate. For~$s=1$, the seminorm is the classical gradient seminorm, and the classical Sobolev--Poincar\'e inequality on the cube, including~$p=1$, applies directly.

\begin{lemma}[Besov duality]
\label{l.circ.norm.dominates}
Let~$m\in\Z$,~$s\in(0,1]$, and~$q,r\in[1,\infty]$, with~$q^{-1}+(q')^{-1}=r^{-1}+(r')^{-1}=1$, and assume that either~$s<1$ or~$(s,q)=(1,1)$. Then~\eqref{e.duality.splitting} holds for all functions~$f,g$ on~$\cu_m$ such that~$fg\in L^1(\cu_m)$ and the right side of~\eqref{e.duality.splitting} is finite.
\end{lemma}
\begin{proof}
Let~$E_n$ denote averaging on the partition into scale-$n$ triadic cubes. At finite depth~$n<m$,
\begin{equation}
	E_ng=(g)_{\cu_m}+\sum_{k=n+1}^m(E_{k-1}g-E_kg)\,.
	\label{e.general.break.down}
\end{equation}
On a scale-$k$ parent~$z+\cu_k$ and its child~$z'+\cu_{k-1}$, one has~$|(g)_{z'+\cu_{k-1}}-(g)_{z+\cu_k}|\leq3^d\|g-(g)_{z+\cu_k}\|_{\underline L^1(z+\cu_k)}$. Pairing~\eqref{e.general.break.down} with~$f$ and applying H\"older's inequality first to the cube averages and then to the sum over scales gives
\begin{equation*}
	\left|\fint_{\cu_m}fE_ng-(f)_{\cu_m}(g)_{\cu_m}\right|
	\leq Cs^{-1}\|f\|_{\Besov{-s}{r}{q}(\cu_m)}
	[g]_{\underline B^s_{1,q',r'}(\cu_m)}\,.
\end{equation*}
To see this, insert the factor~$3^{s(k-1)}3^{-s(k-1)}$ before applying H\"older's inequality to the sum over scales. Replacing the partition averages by the overlapping averages in~\eqref{e.Besov.Morrey.seminorm} introduces a factor depending only on~$d$, and the normalizations of the two sums over scales give the factor~$s^{-1}$, also at the endpoint exponents~$1$ and~$\infty$. The term of the negative norm at the top scale gives
\begin{equation*}
	|(f)_{\cu_m}(g)_{\cu_m}|
	\leq s^{-1}\|f\|_{\Besov{-s}{r}{q}(\cu_m)}
	3^{-sm}|(g)_{\cu_m}|\,.
\end{equation*}
If~$g$ is bounded, then~$E_ng\to g$ almost everywhere and~$|fE_ng|\leq\|g\|_\infty|f|$; dominated convergence proves~\eqref{e.duality.splitting} in this case.

\smallskip

For general~$g$ and~$M>0$, let
\begin{equation*}
	T_M(t)\coloneqq(-M)\vee(t\wedge M)\,,
	\qquad g_M\coloneqq(g)_{\cu_m}+T_M(g-(g)_{\cu_m})-(T_M(g-(g)_{\cu_m}))_{\cu_m}\,.
\end{equation*}
The mean of~$g_M$ is~$(g)_{\cu_m}$. Since~$T_M$ is one-Lipschitz, every cube~$z+\cu_n$ appearing in~\eqref{e.Besov.Morrey.seminorm} satisfies
\begin{align*}
	\|g_M-(g_M)_{z+\cu_n}\|_{\underline L^1(z+\cu_n)}
	&\leq2\|T_M(g-(g)_{\cu_m})-T_M((g)_{z+\cu_n}-(g)_{\cu_m})\|_{\underline L^1(z+\cu_n)}
	\\ &\leq2\|g-(g)_{z+\cu_n}\|_{\underline L^1(z+\cu_n)}\,.
\end{align*}
Thus~$\|g_M\|_{\underline B^s_{1,q',r'}} \leq2\|g\|_{\underline B^s_{1,q',r'}}$. Moreover~$g_M\to g$ almost everywhere and~$|fg_M|\leq|fg|+(2|(g)_{\cu_m}|+\fint_{\cu_m}|g-(g)_{\cu_m}|)|f|$. Applying the bounded case to~$g_M$ and then dominated convergence proves~\eqref{e.duality.splitting}.
\end{proof}

\begin{lemma}
\label{l.Wsp.vs.Bspp}
There exists a constant~$C(d)<\infty$ such that, for every~$s\in (0,1)$,~$p\in [1,\infty)$,~$m \in \Z$ and~$u \in L^p(\cu_m)$,
\begin{equation}
	C^{-1} [  u ]_{\underline{W}^{s,p}(\cu_m)}
	\leq
	[  u ]_{\underline{B}^s_{p,p,p}(\cu_m)}
	\leq
	C [  u ]_{\underline{W}^{s,p}(\cu_m)}
	\,.
	\label{e.Wsp.vs.Bspp}
\end{equation}

\end{lemma}
\begin{proof}
Consider the cubes~$z+\cu_n$ of side~$3^n$ with centers~$z\in3^{n-1}\Zd$ which are contained in~$\cu_m$. Their number is between~$|\cu_m|3^{-nd}$ and~$3^d|\cu_m|3^{-nd}$, and their overlap multiplicity is bounded by~$4^d$. For every such cube, Jensen's inequality and the triangle inequality give
\begin{equation*}
	\fint_{z+\cu_n}|u-(u)_{z+\cu_n}|^p
	\leq\fint_{z+\cu_n}\fint_{z+\cu_n}|u(x)-u(y)|^p\,dx\,dy
	\leq2^p\fint_{z+\cu_n}|u-(u)_{z+\cu_n}|^p\,.
\end{equation*}
Average over~$z\in3^{n-1}\Zd\cap\cu_m$ with~$z+\cu_n\subseteq\cu_m$, multiply by~$s3^{-nsp}$, and sum over~$n\leq m$. Tonelli's theorem shows that~$[u]_{\underline B^s_{p,p,p}(\cu_m)}^p$ is comparable, within factors~$C(d)^p$, to
\begin{equation*}
	\frac{s}{|\cu_m|}\int_{\cu_m}\int_{\cu_m} |u(x)-u(y)|^p
	\sum_{n\leq m}3^{-n(d+sp)}
	\sum_{\substack{z\in3^{n-1}\Zd\cap\cu_m\\z+\cu_n\subseteq\cu_m}}\indc_{z+\cu_n}(x)\indc_{z+\cu_n}(y)\,dx\,dy\,.
\end{equation*}
For distinct~$x,y\in\cu_m$, the points~$x$ and~$y$ lie in a common cube of scale~$n$ only if~$3^n\geq|x-y|_\infty$, so bounded overlap and a geometric sum bound the displayed kernel above by~$C(d)|x-y|_\infty^{-d-sp}$. Conversely, if~$|x-y|_\infty\leq3^{m-1}$, choose~$n\leq m$ with~$3|x-y|_\infty\leq3^n<9|x-y|_\infty$. The mesh spacing~$3^{n-1}$ ensures that~$x,y$ lie in a common~$z+\cu_n$ with~$z\in3^{n-1}\Zd\cap\cu_m$ and~$z+\cu_n\subseteq\cu_m$. For~$|x-y|_\infty>3^{m-1}$, use the top cube~$\cu_m$ instead. The kernel is therefore bounded below by~$9^{-d-sp}|x-y|_\infty^{-d-sp}$. Since~$|x-y|_\infty\leq|x-y|_2\leq\sqrt d\,|x-y|_\infty$, these kernel bounds also apply with Euclidean distance, with a dimensional change in the constants after taking the~$p$th root:~$d^{(d+sp)/(2p)}\leq d^{(d+1)/2}$ for~$p\geq1$ and~$s\leq1$. This proves~\eqref{e.Wsp.vs.Bspp}, with a constant independent of~$s$ and~$p$.
\end{proof}

\begingroup\footnotesize \bibliographystyle{plainnat} \bibliography{refs} \endgroup

@preamble{ "\newcommand{\noop}[1]{} " }

@preamble{
   "\def\cprime{$'$} "
}

@article{ASD,
  title        = {Superdiffusion and anomalous regularization in self-similar random incompressible flows},
  author       = {Armstrong, S. and Bou-Rabee, A. and Kuusi, T.},
  journal      = {arXiv preprint arXiv:2601.22142},
  year         = {2026},
  eprint       = {2601.22142},
  archivePrefix= {arXiv},
  primaryClass = {math.PR},
  note         = {Version 2, 6 September 2026},
  url          = {https://arxiv.org/abs/2601.22142v2},
  doi          = {10.48550/arXiv.2601.22142}
}

@article{Dykhne,
  title={Conductivity of a two-dimensional two-phase system},
  author={Dykhne, A. M.},
  journal={Sov. Phys. JETP},
  volume={32},
  number={1},
  pages={63--65},
  year={1971},
}

@book {BLM,
    AUTHOR = {Boucheron, S. and Lugosi, G. and Massart, P.},
     TITLE = {Concentration inequalities},
      NOTE = {A nonasymptotic theory of independence,
              With a foreword by Michel Ledoux},
 PUBLISHER = {Oxford University Press, Oxford},
      YEAR = {2013},
     PAGES = {x+481},
      ISBN = {978-0-19-953525-5},
   MRCLASS = {60E15 (60D05 60G15 60G50 60G70 62G20)},
  MRNUMBER = {3185193},
MRREVIEWER = {Sreenivasan\ Ravi},
       DOI = {10.1093/acprof:oso/9780199535255.001.0001},
       URL = {https://doi.org/10.1093/acprof:oso/9780199535255.001.0001},
}

@article {Kahane,
    AUTHOR = {Kahane, J.-P.},
     TITLE = {Sur le chaos multiplicatif},
   JOURNAL = {Ann. Sci. Math. Qu\'ebec},
  FJOURNAL = {Annales des Sciences Math\'ematiques du Qu\'ebec},
    VOLUME = {9},
      YEAR = {1985},
    NUMBER = {2},
     PAGES = {105--150},
      ISSN = {0707-9109},
   MRCLASS = {60G57 (60G42)},
  MRNUMBER = {829798},
MRREVIEWER = {S.\ D.\ Chatterji},
}

@article {RVsurvey,
    AUTHOR = {Rhodes, R. and Vargas, V.},
     TITLE = {Gaussian multiplicative chaos and applications: a review},
   JOURNAL = {Probab. Surv.},
  FJOURNAL = {Probability Surveys},
    VOLUME = {11},
      YEAR = {2014},
     PAGES = {315--392},
      ISSN = {1549-5787},
   MRCLASS = {60G57 (28A80 60G15 60G60)},
  MRNUMBER = {3274356},
MRREVIEWER = {Dora\ Sele\v si},
       DOI = {10.1214/13-PS218},
       URL = {https://doi.org/10.1214/13-PS218},
}

@article {Berestycki1,
    AUTHOR = {Berestycki, N.},
     TITLE = {An elementary approach to {G}aussian multiplicative chaos},
   JOURNAL = {Electron. Commun. Probab.},
  FJOURNAL = {Electronic Communications in Probability},
    VOLUME = {22},
      YEAR = {2017},
     PAGES = {Paper No. 27, 12},
      ISSN = {1083-589X},
   MRCLASS = {60G57 (60G15 60J65)},
  MRNUMBER = {3652040},
       DOI = {10.1214/17-ECP58},
       URL = {https://doi.org/10.1214/17-ECP58},
}

@article {GRV,
    AUTHOR = {Garban, C. and Rhodes, R. and Vargas, V.},
     TITLE = {Liouville {B}rownian motion},
   JOURNAL = {Ann. Probab.},
  FJOURNAL = {The Annals of Probability},
    VOLUME = {44},
      YEAR = {2016},
    NUMBER = {4},
     PAGES = {3076--3110},
      ISSN = {0091-1798,2168-894X},
   MRCLASS = {60D05 (28A80 60G15)},
  MRNUMBER = {3531686},
MRREVIEWER = {Artemi\ Berlinkov},
       DOI = {10.1214/15-AOP1042},
       URL = {https://doi.org/10.1214/15-AOP1042},
}

@article {Berestycki2,
    AUTHOR = {Berestycki, N.},
     TITLE = {Diffusion in planar {L}iouville quantum gravity},
   JOURNAL = {Ann. Inst. Henri Poincar\'e{} Probab. Stat.},
  FJOURNAL = {Annales de l'Institut Henri Poincar\'e{} Probabilit\'es et
              Statistiques},
    VOLUME = {51},
      YEAR = {2015},
    NUMBER = {3},
     PAGES = {947--964},
      ISSN = {0246-0203,1778-7017},
   MRCLASS = {60K37 (60J65 60J67 83C45)},
  MRNUMBER = {3365969},
MRREVIEWER = {Dimitri\ Petritis},
       DOI = {10.1214/14-AIHP605},
       URL = {https://doi.org/10.1214/14-AIHP605},
}

@article{DW,
	author = {Ding, J. and Wang, J.},
	title = {Tightness for random walks driven by the planar {G}aussian
		free field at high temperature},
	journal = {Probab. Theory Related Fields},
	year = {2026},
	doi = {10.1007/s00440-026-01467-z},
	url = {https://doi.org/10.1007/s00440-026-01467-z},
	volume = {195},
	pages = {1929--2015},
	number = {3-4},
}

@article{BCGLD,
	author = {Bouchaud, J.-P. and Comtet, A. and Georges, A. and Le Doussal, P.},
	title = {Anomalous diffusion in random media of any dimensionality},
	journal = {J. Phys. France},
	year = 1987,
	volume = 48,
	number = 9,
	pages = "1445-1450",
}

@article {BG,
    AUTHOR = {Bouchaud, J.-P. and Georges, A.},
     TITLE = {Anomalous diffusion in disordered media: statistical
              mechanisms, models and physical applications},
   JOURNAL = {Phys. Rep.},
  FJOURNAL = {Physics Reports. A Review Section of Physics Letters},
    VOLUME = {195},
      YEAR = {1990},
    NUMBER = {4-5},
     PAGES = {127--293},
      ISSN = {0370-1573,1873-6270},
   MRCLASS = {82C44 (82-02 82C31)},
  MRNUMBER = {1081295},
       DOI = {10.1016/0370-1573(90)90099-N},
       URL = {https://doi.org/10.1016/0370-1573(90)90099-N},
}

@article {GM,
    AUTHOR = {Gwynne, E. and Miller, J.},
     TITLE = {Existence and uniqueness of the {L}iouville quantum gravity
              metric for {$\gamma\in(0,2)$}},
   JOURNAL = {Invent. Math.},
  FJOURNAL = {Inventiones Mathematicae},
    VOLUME = {223},
      YEAR = {2021},
    NUMBER = {1},
     PAGES = {213--333},
      ISSN = {0020-9910,1432-1297},
   MRCLASS = {83C45 (30F99 58J90 60D05 60G15)},
  MRNUMBER = {4199443},
MRREVIEWER = {Tom\ Hutchcroft},
       DOI = {10.1007/s00222-020-00991-6},
       URL = {https://doi.org/10.1007/s00222-020-00991-6},
}

@article{DGUniqueness,
    AUTHOR = {Ding, J. and Gwynne, E.},
     TITLE = {Uniqueness of the critical and supercritical {L}iouville quantum gravity metrics},
   JOURNAL = {Proc. Lond. Math. Soc. (3)},
    VOLUME = {126},
    NUMBER = {1},
      YEAR = {2023},
     PAGES = {216--333},
       DOI = {10.1112/plms.12492},
       URL = {https://doi.org/10.1112/plms.12492},
}

@article{DGComparison,
    AUTHOR = {Ding, J. and Gwynne, E.},
     TITLE = {Up-to-constants comparison of {L}iouville first passage percolation and {L}iouville quantum gravity},
   JOURNAL = {Sci. China Math.},
    VOLUME = {66},
    NUMBER = {5},
      YEAR = {2023},
     PAGES = {1053--1072},
       DOI = {10.1007/s11425-021-1983-0},
       URL = {https://doi.org/10.1007/s11425-021-1983-0},
    EPRINT = {2108.12060},
}

@article{KLY,
	title={Random walks in media with constrained disorder},
	author={Kravtsov, V. E. and Lerner, I. V. and Yudson, V. I.},
	journal={J. Phys. A: Math. Gen.},
	volume={18},
	number={12},
	pages={L703--L707},
	year={1985},
	doi={10.1088/0305-4470/18/12/003},
}

@article{Fisher,
  title = {Random walks in random environments},
  author = {Fisher, D. S.},
  journal = {Phys. Rev. A},
  volume = {30},
  issue = {2},
  pages = {960--964},
  numpages = {0},
  year = {1984},
  month = {Aug},
  publisher = {American Physical Society},
  doi = {10.1103/PhysRevA.30.960},
  url = {https://link.aps.org/doi/10.1103/PhysRevA.30.960}
}

@article{FFQSS,
	title={Random walks in two-dimensional random environments with constrained drift forces},
	author={Fisher, D. S. and Friedan, D. and Qiu, Z. and Shenker, S. J. and Shenker, S. H.},
	journal={Phys. Rev. A},
	volume={31},
	number={6},
	pages={3841--3845},
	year={1985},
	doi={10.1103/PhysRevA.31.3841},
}

@incollection {DRSV,
    AUTHOR = {Duplantier, B. and Rhodes, R. and Sheffield, S. and
              Vargas, V.},
     TITLE = {Log-correlated {G}aussian fields: an overview},
 BOOKTITLE = {Geometry, analysis and probability},
    SERIES = {Progr. Math.},
    VOLUME = {310},
     PAGES = {191--216},
 PUBLISHER = {Birkh\"{a}user/Springer, Cham},
      YEAR = {2017},
   MRCLASS = {60D05 (60G22)},
  MRNUMBER = {3821928},
       DOI = {10.1007/978-3-319-49638-2_9},
       URL = {https://doi.org/10.1007/978-3-319-49638-2_9},
}

@article {D,
    AUTHOR = {Dario, P.},
     TITLE = {Optimal corrector estimates on percolation cluster},
   JOURNAL = {Ann. Appl. Probab.},
  FJOURNAL = {The Annals of Applied Probability},
    VOLUME = {31},
      YEAR = {2021},
    NUMBER = {1},
     PAGES = {377--431},
      ISSN = {1050-5164,2168-8737},
   MRCLASS = {60K37 (35Q70 60K35)},
  MRNUMBER = {4254484},
       DOI = {10.1214/20-aap1593},
       URL = {https://doi.org/10.1214/20-aap1593},
}

@article {BDG,
    AUTHOR = {Biskup, M. and Ding, J. and Goswami, S.},
     TITLE = {Return probability and recurrence for the random walk driven
              by two-dimensional {G}aussian free field},
   JOURNAL = {Comm. Math. Phys.},
  FJOURNAL = {Communications in Mathematical Physics},
    VOLUME = {373},
      YEAR = {2020},
    NUMBER = {1},
     PAGES = {45--106},
      ISSN = {0010-3616,1432-0916},
   MRCLASS = {60K37 (60G15 60J10)},
  MRNUMBER = {4050092},
       DOI = {10.1007/s00220-019-03589-z},
       URL = {https://doi.org/10.1007/s00220-019-03589-z},
}

@article{Fitzsimmons2004,
    AUTHOR = {Fitzsimmons, P. J.},
     TITLE = {Superposition operators on {Dirichlet} spaces},
   JOURNAL = {Tohoku Mathematical Journal. Second Series},
    VOLUME = {56},
    NUMBER = {3},
      YEAR = {2004},
     PAGES = {327--340},
       DOI = {10.2748/tmj/1113246670},
       URL = {https://doi.org/10.2748/tmj/1113246670},
}

@article{KuwaeShioya2003,
    AUTHOR = {Kuwae, K. and Shioya, T.},
     TITLE = {Convergence of spectral structures: a functional analytic theory and its applications to spectral geometry},
   JOURNAL = {Communications in Analysis and Geometry},
    VOLUME = {11},
    NUMBER = {4},
      YEAR = {2003},
     PAGES = {599--673},
       DOI = {10.4310/CAG.2003.v11.n4.a1},
       URL = {https://doi.org/10.4310/CAG.2003.v11.n4.a1},
}

@article{GyongyKrylov1996,
    AUTHOR = {Gy{\"o}ngy, I. and Krylov, N.},
     TITLE = {Existence of strong solutions for {It\^o}'s stochastic equations via approximations},
   JOURNAL = {Probability Theory and Related Fields},
    VOLUME = {105},
      YEAR = {1996},
     PAGES = {143--158},
       DOI = {10.1007/BF01203833},
       URL = {https://doi.org/10.1007/BF01203833},
}

@book {Kallenberg,
    AUTHOR = {Kallenberg, O.},
     TITLE = {Foundations of modern probability},
    SERIES = {Probability and its Applications (New York)},
   EDITION = {Second},
 PUBLISHER = {Springer-Verlag, New York},
      YEAR = {2002},
     PAGES = {xx+638},
      ISBN = {0-387-95313-2},
   MRCLASS = {60-01},
  MRNUMBER = {1876169},
MRREVIEWER = {Klaus\ D.\ Schmidt},
       DOI = {10.1007/978-1-4757-4015-8},
       URL = {https://doi.org/10.1007/978-1-4757-4015-8},
}

@misc{AK.HC,
	title={Renormalization group and elliptic homogenization in high contrast},
	Author = {Armstrong, S. and Kuusi, T.},
	note={arXiv:2405.10732, version 3},
	year={2025},
	url = {https://arxiv.org/abs/2405.10732v3},
}

@misc{ABK.SD,
Author = {Armstrong, S. and Bou-Rabee, A.  and Kuusi, T.},
Title = {Superdiffusive central limit theorem for a {B}rownian particle in a critically-correlated incompressible random drift},
Year = {2024},
note = {arXiv:2404.01115, version 2},
	url = {https://arxiv.org/abs/2404.01115v2},
}

@book {AKMBook,
    AUTHOR = {Armstrong, S. and Kuusi, T. and Mourrat,
              J.-C.},
     TITLE = {Quantitative stochastic homogenization and large-scale
              regularity},
    SERIES = {Grundlehren der mathematischen Wissenschaften [Fundamental
              Principles of Mathematical Sciences]},
    VOLUME = {352},
 PUBLISHER = {Springer, Cham},
      YEAR = {2019},
     PAGES = {xxxviii+518},
      ISBN = {978-3-030-15544-5; 978-3-030-15545-2; 978-3-030-15547-6},
   MRCLASS = {35-02 (35B27 60F17 60H15)},
  MRNUMBER = {3932093},
       DOI = {10.1007/978-3-030-15545-2},
      NOTE = {Appendix~A (pp.~416--421) is cited from arXiv:1705.05300v2 (\url{https://arxiv.org/abs/1705.05300v2})},
       URL = {https://doi.org/10.1007/978-3-030-15545-2},
}

@book{FOT,
	author = {Fukushima, M. and Oshima, Y. and Takeda, M.},
	title = {Dirichlet forms and symmetric {M}arkov processes},
	series = {De Gruyter Studies in Mathematics},
	volume = {19},
	edition = {Second revised and extended},
	publisher = {Walter de Gruyter \& Co., Berlin},
	year = {2011}
}

@article{DiNezzaPalatucciValdinoci,
	author = {Di Nezza, E. and Palatucci, G. and Valdinoci, E.},
	title = {Hitchhiker's guide to the fractional {S}obolev spaces},
	journal = {Bull. Sci. Math.},
	volume = {136},
	number = {5},
	year = {2012},
	pages = {521--573},
	doi = {10.1016/j.bulsci.2011.12.004}
}

@article{BBLM2005,
  author = {Boucheron, S. and Bousquet, O. and Lugosi, G. and Massart, P.},
  title = {Moment inequalities for functions of independent random variables},
  journal = {The Annals of Probability},
  volume = {33},
  number = {2},
  pages = {514--560},
  year = {2005},
  doi = {10.1214/009117904000000856}
}

@misc{Devlin2024,
  author = {Devlin, VI, C.},
  title = {Almost Sure Convergence of {Liouville} First Passage Percolation},
  year = {2024},
  note = {arXiv:2309.08001, version 2},
  url = {https://arxiv.org/abs/2309.08001}
}

@article{Meyers1963,
    AUTHOR = {Meyers, N. G.},
     TITLE = {An {$L^p$}-estimate for the gradient of solutions of second order elliptic divergence equations},
   JOURNAL = {Annali della Scuola Normale Superiore di Pisa. Classe di Scienze. Serie III},
    VOLUME = {17},
    NUMBER = {3},
      YEAR = {1963},
     PAGES = {189--206},
       URL = {https://www.numdam.org/item/ASNSP_1963_3_17_3_189_0/},
}

@incollection{BiskupSurvey,
  author = {Biskup, M.},
  title = {Extrema of the two-dimensional {Discrete Gaussian Free Field}},
  booktitle = {Random Graphs, Phase Transitions, and the Gaussian Free Field},
  editor = {Barlow, M. T. and Slade, G.},
  series = {Springer Proceedings in Mathematics \& Statistics},
  volume = {304},
  pages = {163--407},
  publisher = {Springer},
  address = {Cham},
  year = {2020},
  doi = {10.1007/978-3-030-32011-9_3},
  url = {https://arxiv.org/abs/1712.09972}
}

@book{Bass,
  author = {Bass, R. F.},
  title = {Diffusions and Elliptic Operators},
  series = {Probability and Its Applications},
  publisher = {Springer},
  address = {New York},
  year = {1998},
  doi = {10.1007/b97611}
}

@article{CastilloLeDoussal,
    author = {Castillo, H. and Le Doussal, P.},
    title = {Freezing of dynamical exponents in low dimensional random media},
    journal = {Phys. Rev. Lett.},
    volume = {86},
    number = {21},
    pages = {4859--4862},
    year = {2001},
    doi = {10.1103/PhysRevLett.86.4859},
    eprint = {cond-mat/0006373},
    archivePrefix = {arXiv},
    url = {https://arxiv.org/abs/cond-mat/0006373}
}

\end{document}